\documentclass{stml-l}

\usepackage{bbm}
\usepackage{amssymb}
\usepackage{epic,graphicx}
\usepackage{empheq}
\usepackage{url}
\usepackage{calrsfs}
\usepackage{ulem}

\usepackage{hyperref}
\hypersetup{
  colorlinks=true, 
  linktoc=all,     
  linkcolor=blue,  
}
\hypersetup{linktocpage}

\usepackage{backref}

\newcommand{\F}{{\mathbb{F}}}

\newcommand{\Z}{{\mathbb{Z}}}
\newcommand{\Q}{{\mathbb{Q}}}
\newcommand{\N}{{\mathbb{N}}}
\newcommand{\R}{{\mathbb{R}}}
\newcommand{\C}{{\mathbb{C}}}

\newcommand{\cH}{{\mathfrak{h}}}
\newcommand{\cL}{{\mathfrak{g}}}
\newcommand{\gl}{{\mathfrak{gl}}}
\newcommand{\fn}{{\mathfrak{n}}}
\newcommand{\slm}{{\mathfrak{sl}}}
\newcommand{\gom}{{\mathfrak{go}}}

\newcommand{\fb}{{\mathfrak{b}}}
\newcommand{\fg}{{\mathfrak{g}}}
\newcommand{\fh}{{\mathfrak{h}}}
\newcommand{\fj}{{\mathfrak{j}}}
\newcommand{\fs}{{\mathfrak{s}}}

\newcommand{\fH}{{\mathfrak{H}}}
\newcommand{\fS}{{\mathfrak{S}}}

\newcommand{\fX}{{\mathfrak{X}}}
\newcommand{\fY}{{\mathfrak{Y}}}

\newcommand{\be}{{\mathbf{e}}}
\newcommand{\bbf}{{\mathbf{f}}}
\newcommand{\bh}{{\mathbf{h}}}

\newcommand{\bB}{{\mathbf{B}}}

\newcommand{\bE}{{\mathbf{E}}}

\newcommand{\bM}{{\mathbf{M}}}

\newcommand{\cC}{{\mathcal{C}}}
\newcommand{\cE}{{\mathcal{E}}}

\newcommand{\cK}{{\mathcal{K}}}
\newcommand{\cM}{{\mathcal{M}}}

\newcommand{\GAP}{{\sf GAP}}
\newcommand{\Python}{{\sf Python}}

\newcommand{\CHEVIE}{{\sf CHEVIE}}

\newcommand{\cB}{{\mathcal{B}}}

\newcommand{\ad}{{\operatorname{ad}}}
\newcommand{\rad}{{\operatorname{rad}}}
\newcommand{\chark}{{\operatorname{char}}}
\newcommand{\id}{{\operatorname{id}}}
\newcommand{\Aut}{{\operatorname{Aut}}}
\newcommand{\Der}{{\operatorname{Der}}}
\newcommand{\Hom}{{\operatorname{Hom}}}
\newcommand{\End}{{\operatorname{End}}}
\newcommand{\GL}{{\operatorname{GL}}}
\newcommand{\SL}{{\operatorname{SL}}}
\newcommand{\hgt}{{\operatorname{ht}}}
\newcommand{\trc}{{\operatorname{Trace}}}
\renewcommand{\leq}{\leqslant}
\renewcommand{\geq}{\geqslant}

\newcommand{\nm}[1]{{\textit{#1\index{#1}}}}
\newcommand{\nmi}[2]{{\textit{#1\index{#2}}}}
\newcommand{\nms}[2]{{\textit{#1\index{#2@#1}}}}

\newcommand*\widefbox[1]{\fbox{\hspace{0em}#1\hspace{0em}}}

\def\dddots{\mathinner{\mkern1mu\raise1pt
    \vbox{\kern7pt\hbox{.}}\mkern2mu
    \raise4pt\hbox{.}\mkern2mu\raise7pt\hbox{.}\mkern1mu}}
\renewcommand{\atop}[2]{\genfrac{}{}{0pt}{}{#1}{#2}}

\newtheorem{thm}{Theorem}[section]
\newtheorem{lem}[thm]{Lemma}
\newtheorem{cor}[thm]{Corollary}
\newtheorem{prop}[thm]{Proposition}

\theoremstyle{definition}
\newtheorem{defn}[thm]{Definition}
\newtheorem{exmp}[thm]{Example}
\newtheorem{xca}[thm]{Exercise}

\theoremstyle{remark}
\newtheorem{rem}[thm]{Remark}

\numberwithin{section}{chapter}
\numberwithin{equation}{chapter}

\begin{document}

\pagenumbering{roman}
\thispagestyle{empty}

\date{}
%
%


\setcounter{page}{1}

\begin{center}
{\Large\bf A Course on Lie algebras and Chevalley groups}
\end{center}

\medskip
\begin{center}
Meinolf Geck, University of Stuttgart\\
{\tt meinolf.geck@mathematik.uni-stuttgart.de}
\end{center}

\medskip
\begin{quotation} {\small {\bf Abstract}. These are expanded notes from 
graduate courses about Lie algebras and Chevalley groups held at the 
University of Stuttgart.
In the 1950s Chevalley showed how linear groups over arbitrary fields 
could be obtained~---~by a uniform procedure~---~from the simple Lie 
algebras over $\C$ occurring in the Cartan--Killing classification. 
Together with subsequent variations, Chevalley's work had a profound and
long-lasting impact on group theory and Lie theory in general. Classical, 
and widely used references are the lectures notes by Steinberg (1967) and 
the monograph by Carter (1972). Our aim here is to present a self-contained 
introduction to the theory of Chevalley groups, based on recent 
simplifications arising from Lusztig's fundamental theory of ``canonical
bases''. A further feature of our text is that we explicitly incorporate 
algorithmic methods in our treatment, both for the handling of substantial 
examples and regarding some aspects of the general theory.}
\end{quotation}


\section*{Preface} \label{sec0}

\textit{Root systems} are highly symmetrical configurations of vectors
in Euclidean vector spaces. These were first classified by Wilhelm Killing 
around 1890. He found four infinite series, labelled $A_n$, $B_n$, $C_n$,
$D_n$ where $n=1,2,3,\ldots$, and five exceptional ones, labelled $G_2$, 
$F_4$, $E_6$, $E_7$, $E_8$. John Stembridge \cite{Stem} has a beautiful
picture of a plane projection of the $240$ vectors in the $E_8$ system 
(computer-generated from a hand-drawn picture by Peter McMullen):
\begin{center}
\includegraphics[scale=0.6]{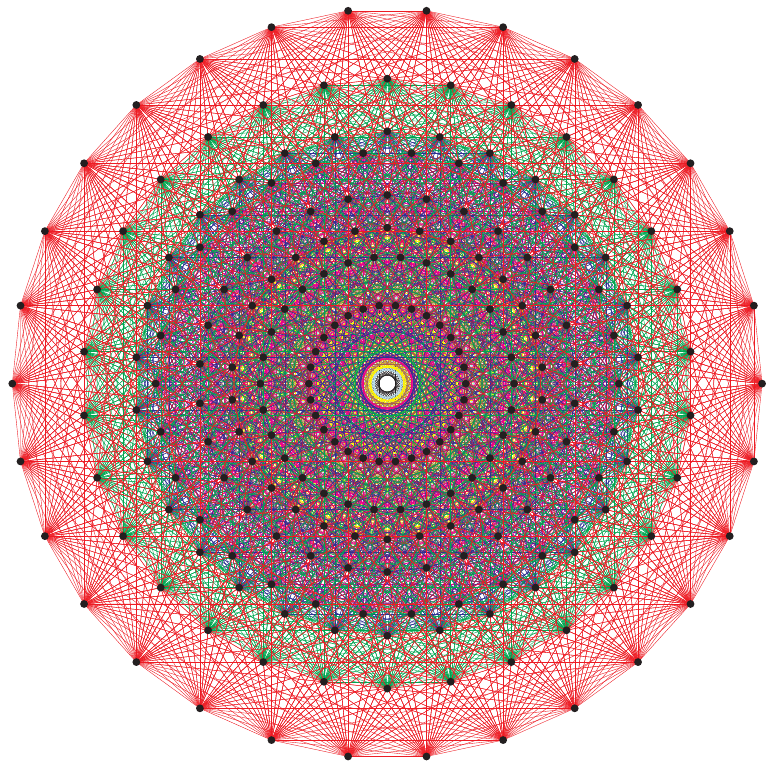}
\end{center}
The corresponding symmetry group~---~nowadays called \textit{Weyl 
group}~---~has a total of $696{,}729{,}600$ elements! 
In Killing's work, root systems arose in the context of 
classifying semisimple Lie algebras. This work has been named by
Coleman \cite{Cole} ``the greatest mathematical paper of all time''. 
For the special role that $E_8$ plays in this story, see Garibaldi's 
recent survey article \cite{gari}.

Nowadays, root systems and their Weyl groups are seen as the combinatorial 
skeleton of various Lie-theoretic structures: the original Lie groups 
and Lie algebras, Kac--Moody algebras and groups, as well as algebraic 
groups and quantum groups. Last but not least, the classification of 
finite simple groups (see Gorenstein et al. \cite{Go1}) highlights the 
importance of Chevalley groups, that is, algebraic versions of Lie groups
over arbitrary fields, in particular, finite fields. An in-depth study of 
these structures requires background material from representation theory, 
differential and algebraic geometry; there is a vast literature on all 
of these subjects. The principal aim of this text is to provide a 
self-contained introduction to the algebraic theory of Chevalley groups, 
together with the required background material about semisimple Lie 
algebras, where we take into account some new developments. 

Decades after Chevalley's seminal work \cite{Ch}, Carter's popular
monograph \cite{Ca1} and Steinberg's famous lectures \cite{St}, there have 
recently been some simplifications of the original construction of 
Chevalley groups. As far as we are aware of, this was first explicitly 
stated as a remark in a short note of Lusztig \cite{L4}:
\begin{quotation}
\textit{The Lie group $E_8$ can be obtained from the graph $E_8$ [...] 
by a method of Chevalley (1955), simplified using the 
theory of ``canonical bases'' (1990).}
\end{quotation}
This remark was further expanded in \cite{G1}, \cite{G2}, \cite{L5}; it 
also sheds some new light on the classical theory of semisimple 
Lie algebras. Our aim here is to develop this in a systematic and 
self-contained way, where we take a purely algebraic point of view and try 
to limit the prerequisites to a minimum; in fact, a good understanding of 
the material in a standard course on Linear Algebra should be sufficient,
together with some basic notions about groups and algebras.

The famous Cartan--Killing classification leads to the consideration of 
specific examples of groups and algebras of ``classical'' types $A_n$, 
$B_n$, $C_n$, $D_n$, and of ``exceptional'' types $G_2$, $F_4$, $E_6$, 
$E_7$, $E_8$. Being able to work with these examples in an efficient way 
is an essential part of the study of Lie theory. As Fulton--Harris 
write \cite[p.~vi]{FH}:
\begin{quotation}
\textit{In most subjects such an approach means one has a few out of an 
unknown infinity of examples which are useful to illuminate the general 
situation. When the subject is the representation theory of complex 
semisimple Lie groups and algebras, however, something special happens: 
once one has worked through all the examples readily at hand~---~the 
``classical'' cases of the special linear, orthogonal, and symplectic 
groups~---~one has not just a few useful examples, one has all but 
five ``exceptional'' cases.}
\end{quotation}
As far as the ``exceptional'' types are concerned, we shall use modern
computer algebra systems to deal with them, both in examples, exercises 
and in some aspects of the general theory. (Of course, the programs also 
work for the ``classical'' types, as long as the dimension is not 
too large.) We also believe that this ``computer algebra approach'' may 
help to better grasp some subtleties of the general theory, e.g., the 
explicit construction~---~in terms of matrices~---~of a spin group of 
type~$D_4$, or a simply connected group of type~$E_7$. (Note that such 
groups are not covered by Carter's book \cite{Ca1}.) And, in any case,
there is a good chance that experiments with large examples on a computer
may lead to new insights (which certainly has happened to the author, and 
probably to many other researchers as well).

Here is a rough outline of the various chapters. (For further details see
the introduction of each individual chapter.)

Chapter~\ref{chap0} introduces a few basic notions and constructions 
concerning Lie algebras. We prove a useful semisimplicity criterion 
in Section~\ref{sec04a}, as well as Lie's Theorem on solvable Lie 
algebras. The final section introduces the Lie algebras of ``classical''
type. We determine their dimensions and show that they are semisimple.
This chapter could even be a basis for a short lecture course on 
Lie algebras, or a topic in a seminar on that subject.

Chapter~\ref{chap2} establishes the main structural results about semisimple
Lie algebras, where we exclusively work over the ground field~$\C$. Our 
treatment deviates from the usual one in textbooks like those of Humphreys 
\cite{H} or Erdmann--Wildon \cite{EW}, for example. More precisely, we 
do not prove here that every semisimple Lie algebra has a Cartan 
subalgebra and a Cartan decomposition, the main ingredients that 
eventually lead to root systems and Dynkin diagrams. Instead, 
inspired by Moody--Pianzola \cite{MP}, we formulate a small 
set of axioms (see Definition~\ref{defTD}) from which the Cartan 
decomposition, root systems etc.\ can be derived without too much effort. 
A Lie algebra satisfying these axioms will be called a Lie algebra
of \nm{Cartan--Killing type}. We will show that such Lie algebras are 
semisimple.  (And, conversely, it is known that all semisimple algebras
are ``of Cartan--Killing type'' but the point is that we will not prove 
this here.) One advantage of this approach is that it allows us to focus
rather quickly on certain more modern aspects of the theory of Lie algebras.

In Section~\ref{sec1a7} we prove Lusztig's fundamental result \cite{L1}, 
\cite{L2}, \cite{L5} that every Lie algebra of Cartan--Killing type has a 
``canonical'' basis, with highly remarkable positivity properties (see 
Remark~\ref{canbas1}). This is a profound strengthening of the existence 
of ``integral'' bases due to Chevalley~\cite{Ch}. (Explicit formulae for 
the structure constants with respect to that basis have been recently 
determined in~\cite{GeLa}.)

Chapter~\ref{chap3} begins with the classification of the Cartan matrices of 
simple Lie algebras, where we use the approach taken in the theory of 
Kac--Moody algebras. Then the main theme of the chapter is to show that, 
starting from such a Cartan matrix, one can construct a corresponding Lie 
algebra of Cartan--Killing type (over $\C$) and a Chevalley group (over any 
field) by purely combinatorial procedures. More precisely, the Chevalley 
groups constructed here are those of ``adjoint type'', with a trivial
center. Our approach will be based on the existence of Lusztig's canonical 
basis. This typically leads to more rigid formulae as compared to the 
traditional approach (as, e.g., in Carter's book \cite{Ca1}), where a number 
of signs may appear which are difficult to control; see, e.g.,
Theorem~\ref{luform}. Here, in Chapter~\ref{chap3}, we will arrive 
at the point where we can show that a Chevalley group is equal to its own 
commutator subgroup (when the base field is not too small); this is one
big step towards proving that the group is simple. 

The subject of Chapter~\ref{chap4} is the construction of a more 
general type of Chevalley groups, which may have a non-trivial center. 
This relies on the existence of ``admissible'' lattices, as defined by
Chevalley \cite[\S 4]{Ch2}, in the finite-dimensional representations of 
the underlying Lie algebra. This was discovered independently by Ree 
\cite{ree2}~---~in quite a compact and elementary fashion. (But, somewhat 
surprisingly, \cite{ree2} seems to have received only very few citations.) 
Nowadays, one usually proceeds using Kostant's $\Z$-form of the universal
enveloping algebra; see Bourbaki \cite[Ch.~VIII, \S 12]{B78}, Humphreys 
\cite[\S 27]{H} or Steinberg \cite[Chap.~2]{St}. See also Lusztig 
\cite{L3} for a different approach, using extensively the theory of 
``canonical bases''. 

Here, we follow the approach in \cite{G2} which relies on elementary 
properties of ``minuscule'' weights and the corresponding representations. 
In Sections~\ref{sec5a1} and \ref{sechwm} this will be developed from first 
principles based on Stembridge \cite[\S 1]{Stem1} and Jantzen \cite[5A.1]{Ja}. 
One advantage is that this yields very explicit models of Chevalley groups 
of non-adjoint type. The further exposition in Chapter~\ref{chap4} is a
synthesis of Ree \cite{ree2} and Steinberg \cite{St}. But, generally 
speaking, we provide considerably more details on a number of arguments 
and calculations, similar in spirit to the style of Carter's book \cite{Ca1}
(``in which nothing is left to the reader'', as Solomon writes in his
impressive MathSciNet review of \cite{Ca1}). We hope that readers will 
find this a useful addition to the existing literature on Chevalley groups.

The current plan is to add a further chapter in the future, which will be 
concerned with the main structural properties of Chevalley groups, of 
adjoint or non-adjoint type. This includes the fundamental ``commutator 
relations'' due to Chevalley and the discussion of $(B,N)$-pairs in the 
sense of Tits~\cite{Ti0}. We shall also place the Chevalley groups in the 
context of the theory of (semisimple) algebraic groups over algebraically 
closed fields. 

\smallskip
\centerline{-----------------------------------------------------------}

\smallskip
My own motivation for studying Chevalley groups comes from finite
group theory, especially the classification of finite simple groups. 
As Curtis writes in \cite[p.~viii]{Cu}, the representation theory of 
finite groups has always had an abundance of challenging problems. And 
it is particularly interesting to study these problems for Chevalley 
groups over finite fields; see, e.g., Lusztig \cite{shaw} (and 
\cite{GeMa} for further references). A link between Lie's theory of 
continuous groups and the theory of finite groups is provided by 
Chevalley's famous classification \cite{Ch3} of semisimple algebraic 
groups over an arbitrary algebraically closed field~---~considered
by many as his masterpiece; see Dieudonn\'e--Tits \cite{DT} and the
postface of \cite{Ch3} (written by Cartier) for further historical context. 
In this setting, the finite Chevalley groups are seen to be groups of 
rational points of an algebraic group $\mathbf{G}$ over~$\overline{\F}_p$, 
an algebraic closure of the finite field with $p$ elements. This viewpoint,
fully developed by Steinberg \cite{St68}, is essential for many further
developments, e.g., the construction of representations after 
Deligne--Lusztig~\cite{DeLu}. Another advantage of this viewpoint is that
it provides a common framework for dealing both with ``split'' and with 
``twisted'' Chevalley groups (even the ``very twisted'' groups of Ree and 
Suzuki): all these arise as groups of fixed points
\[\textbf{G}^F:=\{g \in \mathbf{G}\mid F(g)=g\}\]
where $F\colon \mathbf{G} \rightarrow \mathbf{G}$ is a (generalized) 
Frobenius map. Our text stops once we will have arrived at this point.
For further developments (e.g., the detailed study of twisted groups, 
the general theory of reductive algebraic groups and so on), we then
refer the reader to other sources.

I taught most of the material in Chapters~1--3 in a Master course at the 
University of Stuttgart in the academic year 2019/20, and again in 2024. 
(There were 14 weeks of classes, with two lectures plus one exercise class 
per week, each 90 minutes long.) I thank the students Moritz G\"osling, 
Janik Hess, Alexander Lang, Paul Schwahn for corrections, lists of misprints 
etc. Special thanks go to Gunter Malle for a careful reading of the whole
manuscript.

Carter's book \cite{Ca1} and Steinberg's lectures \cite{St} contain much 
more about Chevalley groups than we can~---~or even want to~---~present in 
this book. (See also Borel's paper in \cite{Borel}.) For those new to the 
theory, our exposition may provide a gentler introduction, with a 
self-contained background from the theory of Lie algebras, with first steps 
towards the theory of algebraic groups, and by avoiding some of the heavier 
machinery usually required for the construction of Chevalley groups of 
non-adjoint type. In essence, we try to be as elementary and detailed as 
Carter \cite{Ca1}, but also include non-adjoint groups in the discussion 
(which are not dealt with at all in \cite{Ca1}). At the same time, our 
treatment remotely touches upon some more recent developments that go 
beyond the classical theory of semisimple Lie algebras: Kac--Moody algebras 
and ``canonical bases''. 

\medskip
\hfill \textit{Stuttgart, September 2026}

\newpage
\medskip
\noindent {\large\bf Contents}

\smallskip

{\small
\contentsline {chapter}{\tocchapter {Chapter}{1}{Introducing Lie algebras}}{1}{chapter.1}%
\contentsline {section}{\tocsection {}{1.1}{Non-associative algebras}}{1}{section.1.1}%
\contentsline {section}{\tocsection {}{1.2}{Matrix Lie algebras and derivations}}{8}{section.1.2}%
\contentsline {section}{\tocsection {}{1.3}{Solvable and semisimple Lie algebras}}{15}{section.1.3}%
\contentsline {section}{\tocsection {}{1.4}{Representations of Lie algebras}}{21}{section.1.4}%
\contentsline {section}{\tocsection {}{1.5}{Lie's Theorem}}{28}{section.1.5}%
\contentsline {section}{\tocsection {}{1.6}{The classical Lie algebras}}{35}{section.1.6}%
\contentsline {section}{\tocsection {}{}{Notes on Chapter\nonbreakingspace \ref {chap0}}}{41}{section*.9}%
\contentsline {chapter}{\tocchapter {Chapter}{2}{Semisimple Lie algebras}}{43}{chapter.2}%
\contentsline {section}{\tocsection {}{2.1}{Weights and weight spaces}}{44}{section.2.1}%
\contentsline {section}{\tocsection {}{2.2}{Lie algebras of Cartan--Killing type}}{52}{section.2.2}%
\contentsline {section}{\tocsection {}{2.3}{The Weyl group}}{62}{section.2.3}%
\contentsline {section}{\tocsection {}{2.4}{Semisimplicity}}{72}{section.2.4}%
\contentsline {section}{\tocsection {}{2.5}{Classical Lie algebras revisited}}{82}{section.2.5}%
\contentsline {section}{\tocsection {}{2.6}{The structure constants $N_{\alpha ,\beta }$}}{90}{section.2.6}%
\contentsline {section}{\tocsection {}{2.7}{Lusztig's canonical basis}}{100}{section.2.7}%
\contentsline {section}{\tocsection {}{}{Notes on Chapter\nonbreakingspace \ref {chap2}}}{113}{section*.43}%
\contentsline {chapter}{\tocchapter {Chapter}{3}{Generalised Cartan matrices}}{115}{chapter.3}%
\contentsline {section}{\tocsection {}{3.1}{Classification}}{116}{section.3.1}%
\contentsline {section}{\tocsection {}{3.2}{Finite root systems}}{125}{section.3.2}%
\contentsline {section}{\tocsection {}{3.3}{A glimpse of Kac--Moody theory}}{136}{section.3.3}%
\contentsline {section}{\tocsection {}{3.4}{Using computers: {\sf CHEVIE} and {\sf ChevLie}}}{146}{section.3.4}%
\contentsline {section}{\tocsection {}{3.5}{Introducing Chevalley groups}}{157}{section.3.5}%
\contentsline {section}{\tocsection {}{3.6}{A first example: Groups of type $A_{n-1}$}}{168}{section.3.6}%
\contentsline {section}{\tocsection {}{3.7}{The elements $\bar {n}_i(\xi )$ and $\bar {h}_i(\xi )$}}{173}{section.3.7}%
\contentsline {section}{\tocsection {}{}{Notes on Chapter\nonbreakingspace \ref {chap3}}}{183}{section*.64}%
\contentsline {chapter}{\tocchapter {Chapter}{4}{General construction of Chevalley groups}}{187}{chapter.4}%
\contentsline {section}{\tocsection {}{4.1}{The weight lattice of a representation}}{188}{section.4.1}%
\contentsline {section}{\tocsection {}{4.2}{Minuscule weights}}{199}{section.4.2}%
\contentsline {section}{\tocsection {}{4.3}{Highest weight modules}}{208}{section.4.3}%
\contentsline {section}{\tocsection {}{4.4}{Admissible lattices in ${\mathfrak {g}}$-modules}}{219}{section.4.4}%
\contentsline {section}{\tocsection {}{4.5}{The elements $\bar {x}_\alpha (\zeta ;V,{\mathrsfs {B}})$ in $G_K(V,{\mathrsfs {B}})$}}{230}{section.4.5}%
\contentsline {section}{\tocsection {}{4.6}{The diagonal and monomial subgroups}}{243}{section.4.6}%
\contentsline {section}{\tocsection {}{4.7}{Chevalley groups of type $A_1$}}{255}{section.4.7}%
\contentsline {section}{\tocsection {}{}{Notes on Chapter\nonbreakingspace \ref {chap4}}}{263}{section*.81}%
\contentsline {chapter}{\tocappendix {Appendix}{A}{Some complements and auxiliary results}}{267}{appendix.A}%
\contentsline {section}{\tocsection {}{A.1}{Generation of ${\operatorname {SL}}_n(K)$}}{267}{section.A.1}%
\contentsline {section}{\tocsection {}{A.2}{Matsumoto's Lemma}}{269}{section.A.2}%
\contentsline {chapter}{\tocchapter {\chaptername }{}{Bibliography}}{273}{appendix*.86}%
\contentsline {chapter}{\tocchapter {\chaptername }{}{Index}}{281}{appendix*.87}%
}

\mainmatter

\chapter{Introducing Lie algebras} \label{chap0}

This chapter introduces Lie algebras and describes some fundamental 
constructions related to them, e.g., representations and derivations. This 
is illustrated with a number of examples, most notably certain matrix Lie 
algebras. As far as the general theory is concerned, we will arrive at 
the point where we can single out the important class of ``semisimple'' 
Lie algebras. 

Throughout this chapter, $k$ denotes a fixed base field. All vector spaces
will be understood to be vector spaces over this field~$k$. We use 
standard notions from Linear Algebra: dimension (finite or infinite), 
linear and bilinear maps, matrices, eigenvalues. Everything else will be 
formally defined but we will assume a basic familiarity with general algebraic 
constructions, e.g., substructures and homomorphisms.

\section{Non-associative algebras} \label{sec01}

Let $A$ be a vector space (over $k$). If we are also given a bilinear map
\begin{center}
$A\times A\rightarrow A, \qquad (x,y)\mapsto  x\cdot y$,
\end{center}
then $A$ is called an \nm{algebra} (over $k$). Familiar examples 
from Linear Algebra are the algebra $A=M_n(k)$ of all $n\times n$-matrices
with entries in~$k$ (and the usual matrix product), or the algebra $A=k[T]$
of polynomials with coefficients in $k$ (where $T$ denotes an 
indeterminate). In these examples, the product in $A$ is associative;
in the second example, the product is also commutative. But for us here, 
the term ``algebra'' does not imply any further assumptions on the product 
in~$A$ (except bi-linearity). If the product in $A$ happens to be
associative (or commutative or $\ldots$), then we say explicitly that $A$ 
is an ``associative algebra'' (or ``commutative algebra'' or $\ldots$).

The usual basic algebraic constructions also apply in this general setting.
We will not completely formalize all of this, but assume that the reader
will fill in some (easy) details if required. Some examples:

$\bullet$ If $A$ is an algebra and $B\subseteq A$ is a subspace, then $B$ 
is called a \nm{subalgebra} if $x\cdot y\in B$ for all $x,y\in B$. In this
case, $B$ itself is an algebra (with product given by the restriction of
$A\times A\rightarrow A$ to $B\times B$). One easily checks that, if 
$\{B_i\}_{i\in I}$ is a family of subalgebras (where $I$ is any
indexing set), then $\bigcap_{i\in I} B_i$ is a subalgebra.

$\bullet$ If $A$ is an algebra and $B\subseteq A$ is a subspace, then $B$ 
is called an \nm{ideal} if $x\cdot y\in B$ and $y\cdot x\in B$ for all
$x\in A$ and $y\in B$. In particular, $B$ is a subalgebra in this case.
Furthermore, the quotient vector space $A/B=\{x+B\mid x\in A\}$ is an
algebra with product given by
\[ A/B \times A/B\rightarrow A/B,\qquad (x+B,y+B)\mapsto x\cdot y+B.\]
(One checks as usual that this product is well-defined and bilinear.)
Again, one easily checks that, if $\{B_i\}_{i\in I}$ is a family of ideals
(where $I$ is any indexing set), then $\bigcap_{i\in I} B_i$ is an ideal.

$\bullet$ If $A,B$ are algebras, then a linear map $\varphi\colon A
\rightarrow B$ is called an \nm{algebra homomorphism} if $\varphi(x\cdot
y)=\varphi(x)*\varphi(y)$ for all $x,y\in A$. (Here, ``$\cdot$'' is 
the product in $A$ and ``$*$'' is the product in $B$.) If, furthermore,
$\varphi$ is bijective, then we say that $\varphi$ is an \nm{algebra 
isomorphism}. In this case, the inverse map $\varphi^{-1} \colon B
\rightarrow A$ is also an algebra homomorphism and we write $A\cong B$ 
(saying that $A$ and $B$ are \textit{isomorphic}).

$\bullet$ If $A,B$ are algebras and $\varphi\colon A\rightarrow B$
is an algebra homomorphism, then the kernel $\ker(\varphi)$ is an 
ideal in $A$ and the image $\varphi(A)$ is a subalgebra of $B$. 
Furthermore, we have a canonical induced homomorphism $\bar{\varphi}
\colon A/\ker(\varphi)\rightarrow B$, $x+\ker(\varphi)\mapsto\varphi(x)$, 
which is injective and whose image equals $\varphi(A)$. Thus, we have
$A/\ker(\varphi)\cong \varphi(A)$.

Some further pieces of general notation. If $V$ is a vector space and $X
\subseteq V$ is a subset, then we denote by $\langle X\rangle_k\subseteq V$
the subspace spanned by $X$. Now let $A$ be an algebra. Given $X\subseteq 
A$, we denote by $\langle X\rangle_{\text{alg}}\subseteq A$ the subalgebra
generated by $X$, that is, the intersection of all subalgebras of $A$ that 
contain $X$. One easily checks that 
$\langle X \rangle_{\text{alg}}=\langle \hat{X}\rangle_k$ where $\hat{X}=
\bigcup_{n\geq 1} X_n$ and the subsets $X_n\subseteq A$ are inductively 
defined by $X_1:=X$ and
\[ X_n:=\{x\cdot y \mid x\in X_i, y \in X_{n-i} \mbox{ for $1\leq i 
\leq n-1$}\}\quad \mbox{for $n\geq 2$}.\]
Thus, the elements in $X_n$ are obtained by taking the iterated product, 
in any order and for any bracketing, of $n$ elements of $X$. We call the 
elements of $X_n$ \nm{monomials} in $X$ (of level $n$). For example, if 
$X=\{x,y,z\}$, then $((z\cdot (x\cdot y))\cdot z)\cdot ((z\cdot y)\cdot 
(x\cdot x))$ is a monomial of level~$8$ and, in general, we have to respect 
the parentheses in working with such products.

\begin{exmp} \label{magma} Let $M$ be a non-empty set and $\mu \colon M
\times M \rightarrow M$ be a map. Then the pair $(M,\mu)$ is called a
\nm{magma}. Now the set of all functions $f\colon M\rightarrow k$ is a vector
space over $k$, with pointwise defined addition and scalar multiplication.
Let $k[M]$ be the subspace consisting of all $f\colon M \rightarrow k$
such that $\{x \in M \mid f(x)\neq 0\}$ is finite. For $x \in M$, let 
$\varepsilon_x\in k[M]$ be defined by $\varepsilon_x(y)= 1$ if $x=y$ and 
$\varepsilon_x(y)=0$ if $x\neq y$. Then one easily sees that 
$\{\varepsilon_x\mid x\in M\}$ is a basis of $k[M]$. Furthermore, we can 
uniquely define a bilinear map
\[k[M]\times k[M]\rightarrow k[M] \quad \mbox{such that} \quad 
(\varepsilon_x, \varepsilon_y) \mapsto \varepsilon_{\mu(x,y)}.\]
Then $A=k[M]$ is an algebra, called the \nm{magma algebra} of $(M,\mu)$ 
over $k$.
\end{exmp}

We have the following useful result. Let $(M,\mu)$ be a magma as 
above. We call a function $\chi\colon M \rightarrow k$ a \nm{character} if 
$\chi$ is not identically zero and if $\chi(\mu(x,y))=\chi(x)\chi(y)$ for 
all $x,y\in M$. For example, if $A$ is an algebra, then any non-zero algebra 
homomorphism $\chi\colon A\rightarrow k$ is a character (where we regard~$A$ 
as a monoid with $\mu\colon A \times A \rightarrow A$ given by the 
multiplication in~$A$).

\begin{lem}[Linear independence of characters] \label{dedek} \nmi{}{Linear 
independence of characters} Let $n\geq 1$ and $\chi_1,\ldots,\chi_n\colon 
M\rightarrow k$ be distinct characters. Then $\chi_1,\ldots,\chi_n$ are 
linearly independent (in the vector space of all functions $f\colon 
M\rightarrow k$).
\end{lem}

\begin{proof} This is a standard result in Algebra; it is usually stated
for distinct homomorphisms of a group into the multiplicative group of~$k$.
See, for example, Milne \cite[Theorem~5.14]{milne2}. But exactly the same 
proof also works in the present, more general situation.
\end{proof}

\begin{exmp} \label{dirprodalg} Let $r\geq 1$ and $A_1,\ldots,A_r$ be
algebras (all over $k$). Then the cartesian product $A:=A_1\times\ldots
\times A_r$ is a vector space with component-wise defined addition
and scalar multiplication. But then $A$ also is an algebra with
product 
\[A\times A\rightarrow A,\quad \bigl((x_1,\ldots,x_r),(y_1,\ldots,y_r)
\bigr) \mapsto (x_1\cdot y_1,\ldots, x_r\cdot y_r),\]
where, in order to simplify the notation, we denote the product in 
each $A_i$ by the same symbol ``$\cdot$''. For a fixed $i$, we have
an injective algebra homomorphism $\iota_i \colon A_i\rightarrow A$
sending $x\in A_i$ to $(0,\ldots,0,x,0,\ldots,0)\in A$ (where $x$ appears
in the $i$-th position). If $\underline{A}_i\subseteq A$ denotes the image
of $\iota_i$, then we have a direct sum $A=\underline{A}_1 \oplus \ldots 
\oplus \underline{A}_r$ where each $\underline{A}_i$ is an ideal in $A$ and,
for $i\neq j$, we have $\underline{x}\cdot \underline{y}=0$ for all 
$\underline{x}\in \underline{A}_i$ and $\underline{y}\in \underline{A}_j$. 
The algebra $A$ is called the \nm{direct product} of $A_1,\ldots,A_r$.
\end{exmp}

\begin{rem} \label{adjoint0} Let $A$ be an algebra and $\End(A)$ be the
vector space of all linear maps of $A$ into itself; then $\End(A)$ is an
associative algebra with product given by the composition of maps. For 
$x\in A$, we have maps $L_x\colon A\rightarrow A$, $y\mapsto x\cdot y$, 
and $R_x\colon A \rightarrow A$, $y\mapsto y\cdot x$. Then note:
\begin{center}
$A$ is associative $\quad \Leftrightarrow \quad L_x\circ R_y=
R_y\circ L_x\;$ for all $x,y\in A$.
\end{center}
This simple observation is a useful ``trick'' in proving certain
identities. Here is one example. For $x\in A$, we denote $\ad_A(x):=L_x-
R_x\in \End(A)$. Thus, $\ad_A(x)(y)=x\cdot y-y\cdot x$ for all $x,y\in 
A$. The following result may be regarded as a \nm{generalized binomial 
formula}; it will turn out to be useful at a few places in the sequel. 
\end{rem}

\begin{lem} \label{superbin} Let $A$ be an associative algebra with 
identity element~$1_A$. Let $x,y\in A$, $a,b\in k$ and $n\geq 0$. Then 
\begin{align*} 
(x&+(a+b)1_A)^n\cdot y\\&=\sum_{i=0}^n \binom{n}{i} 
(\ad_A(x)+b\,\id_A)^i(y) \cdot (x+a\,1_A)^{n-i}.
\end{align*}
{\rm (Here, $\id_A\colon A\rightarrow A$ denotes the identity map.)}
\end{lem}

\begin{proof} As above, we have $\ad_A(x)=L_x-R_x$. Now $L_{x+(a+b)1_A}(y)=
x\cdot y+(a+b)y=(L_x+(a+b)\id_A)(y)$ for all $y\in A$ and so 
\[L_{x+(a+b)1_A}=L_x+(a+b)\id_A=(R_x+a\,\id_A)+(\ad_A(x)+b\,\id_A).\]
Since $A$ is associative, $L_x$ and $R_x$ commute with each other and, 
hence, $\ad_A(x)$ commutes with both $L_x$ and $R_x$. Consequently, the
maps $\ad_A(x)+b\,\id_A$ and $R_{x+a1_A}=R_x+a\,\id_A$ commute with each 
other. Hence, working in $\End(A)$, we can apply the usual binomial formula
to $L_{x+(a+b)1_A}=R_{x+a1_A}+(\ad_A(x)+b\,\id_A)$ and obtain:
\[ L_{x+(a+b)1_A}^n=\sum_{i=0}^n \binom{n}{i} R_{x+a1_A}^{n-i}
\circ (\ad_A(x)+b\,\id_A)^i.\]
Evaluating at $y$ yields the desired formula. 
\end{proof}

After these general considerations, we now introduce the particular 
(non-associative) algebras that we are interested in here. 

\begin{defn} \label{deflie} Let $A$ be an algebra (over $k$), with
product $x\cdot y$ for $x,y\in A$. We say that $A$ is a \textit{Lie 
algebra} if this product has the following two properties:
\begin{itemize}
\item (Anti-symmetry) We have $x\cdot x=0$ for all $x\in A$. Note that,
using bi-linearity, this implies $x\cdot y=-y\cdot x$ for all 
$x,y\in A$.
\item (Jacobi identity) We have $x\cdot (y\cdot z)+y\cdot (z\cdot x)
+z\cdot (x\cdot y)=0$ for all $x,y,z\in A$.
\end{itemize}
The above two rules imply the formula $x\cdot (y\cdot z)=(x\cdot y)\cdot z+
y\cdot (x\cdot z)$ which has some resemblance to the rule for 
differentiating a product. 

Usually, the product in a Lie algebra is denoted by $[x,y]$ (instead
of $x\cdot y$) and called \nm{bracket}. So the above formulae read as 
follows. 
\[ [x,x]=0 \qquad \mbox{and}\qquad [x,[y,z]]+[y,[z,x]]+[z,[x,y]]=0.\]
Usually, we will use the symbol ``$L$'' to denote a Lie algebra.
\end{defn}

\begin{exmp} \label{xcagrass} Let $L=\R^3$ (row vectors). Let $(x,y)$ be the 
usual scalar product of $x,y\in \R^3$, and $x\times y$ be the ``vector
product'' (perhaps known from a Linear Algebra course). That is, given
$x=(x_1,x_2,x_3)$ and $y=(y_1,y_2,y_3)$ in $L$, we have $x\times y=
(v_1,v_2,v_3)\in L$ where
\[v_1=x_2y_3-x_3y_2,\qquad v_2=x_3y_1-x_1y_3,\qquad v_3=x_1y_2-x_2y_1.\]
One easily verifies the ``Grassmann identity'' $x \times (y \times z)=
(x,z)\, y-(x,y)\, z$ for $x,y,z \in \mathbb{R}^3$. Setting $[x,y]:=
x \times y$ for $x,y\in L$, a straightforward computation shows 
that $L$ is a Lie algebra over $k=\R$. 
\end{exmp}

\begin{exmp} \label{lienorma} Let $L$ be a Lie algebra. If 
$V\subseteq L$ is any subspace, the \nm{normalizer} of $V$ is defined as 
\[I_L(V):=\{x\in L \mid [x,v]\in V \mbox{ for all $v\in V$}\}.\]
Clearly, $I_L(V)$ is a subspace of $L$. We claim that $I_L(V)$ is a Lie
subalgebra of $L$. Indeed, let $x,y\in I_L(V)$ and $v\in V$. By
the Jacobi identity and anti-symmetry, we have 
\[[[x,y],v]=-[v,[x,y]]=[x,\underbrace{[y,v]}_{\in V}]-[y,\underbrace{[x,
v]}_{\in V}]\in V\]
and so $[x,y]\in I_L(V)$. If $V$ is a Lie subalgebra, then $V\subseteq
I_L(V)$ and $V$ is an ideal in $I_L(V)$.
\end{exmp}

\begin{xca} \label{xcagenerator} Let $L$ be a Lie algebra and 
$X \subseteq L$ be a subset.

(a) Let $z\in L$ be such that $[x,z]=0$ for all $x\in X$. Then show 
that $[y,z]=0$ for all $y \in \langle X \rangle_{\text{alg}}$.

(b) Let $V\subseteq L$ be a subspace such that $[x,v]\in V$ for all 
$x\in X$ and $v\in V$. Then show that $[y, v]\in V$ for all $y\in 
\langle X\rangle_{\text{alg}}$ and $v\in V$. Furthermore, if 
$X\subseteq V$, then $\langle X \rangle_{\text{alg}}\subseteq V$.

(c) Let $I:=\langle X\rangle_{\text{alg}} \subseteq L$. Assume that
$[x,y]\in I$ for all $x\in X$ and $y\in L$. Then show that $I$ is an
ideal of~$L$.

(d) Let $L'$ be a further Lie algebra and $\varphi\colon L \rightarrow L'$
be a linear map. Assume that $L=\langle X \rangle_{\text{alg}}$. Then show 
that $\varphi$ is a Lie algebra homomorphism if $\varphi([x,y])=[\varphi(x),
\varphi(y)]$ for all $x\in X$ and $y\in L$.

\noindent {\footnotesize [{\it Hint}. We have $\langle X\rangle_{\text{alg}}
= \langle X_n\mid n \geq 1\rangle_k$ where $X_n$ is the set of monomials
of level $n$ in~$X$. Now use induction on $n$ and the Jacobi identity.]}
\end{xca}

\begin{exmp} \label{deflie1} (a) Let $V$ be a vector space. We define
$[x,y]:=0$ for all $x,y\in V$. Then, clearly, $V$ is a Lie algebra.
A Lie algebra in which the bracket is identically $0$ is called an 
\nm{abelian Lie algebra}.

(b) Let $A$ be an algebra that is associative. Then we define a new 
product on $A$ by $[x,y]:=x\cdot y-y\cdot x$ for all $x,y\in A$. Clearly, 
this is bilinear and we have $[x,x]=0$; furthermore, for $x,y,z\in A$, 
we have
\begin{align*} 
[x,[y,z]]+&[y,[z,x]]+[z,[x,y]]\\&=[x,y\cdot z-z\cdot y]+[y,z\cdot x-
x\cdot z]+[z,x\cdot y-y\cdot x]\\
&=x\cdot (y\cdot z-z\cdot y)-(y\cdot z-z\cdot y)\cdot x\\
&\qquad\quad + y\cdot (z\cdot x- x\cdot z)- (z\cdot x- x\cdot z)\cdot y\\
&\quad\qquad\qquad+ z\cdot (x\cdot y-y\cdot x)-(x\cdot y-y\cdot x)\cdot z.
\end{align*}
By associativity, we have $x\cdot (y\cdot z)=(x\cdot y)\cdot z$ and so 
on. We then leave it to the reader to check that the above sum collapses 
to~$0$. Thus, every associative algebra becomes a Lie algebra by this
construction.
\end{exmp}

A particular role in the general theory is played by those algebras that
do not have non-trivial ideals. This leads to: 

\begin{defn} \label{defsimple} Let $A$ be an algebra such that 
$A\neq \{0\}$ and the product of $A$ is not identically zero. Then $A$ is 
called a \nm{simple algebra} if $\{0\}$ and $A$ are the only ideals of $A$.
\end{defn} 

We shall see first examples in the following section.

\begin{xca} \label{xcafree} This exercise (which may be skipped on a 
first reading) presents a very general method for constructing algebras 
with prescribed properties. Recall from Example~\ref{magma} the 
definition of a magma. Given a non-empty set $X$, we want to define the 
``most general magma'' containing~$X$, following Bourbaki \cite[Chap.~I, 
\S 7, no.~1]{B1}. For this purpose, we define inductively sets $X_n$ 
for $n=1,2,\ldots$, as follows. We set $X_1:=X$. Now let $n\geq 2$ and 
assume that $X_i$ is already defined for $1\leq i\leq n-1$. Then define 
$X_n$ to be the disjoint union of the sets $X_i\times X_{n-i}$ for $1\leq
i \leq n-1$. Finally, we define $M(X)$ to be the disjoint union of all
the sets $X_n$, $n\geq 1$.

Now let $w,w'\in M(X)$. Since $M(X)$ is the disjoint union of all~$X_n$,
there are unique $p,p'\geq 1$ such that $w\in X_p$ and $w'\in X_{p'}$.
Let $n:=p+p'$. By the definition of $X_n$, we have $X_p\times X_{p'}
\subseteq X_n$. Then define $w*w'\in X_n$ to be the pair $(w,w')\in X_p
\times X_{p'}\subseteq X_n$. In this way, we obtain a product $M(X) 
\times M(X)\rightarrow M(X)$, $(w,w')\mapsto w*w'$. So $M(X)$ is a magma,
called the \nm{free magma} on $X$.

Thus, one may think of the elements of $M(X)$ as arbitrary 
``non-associative words'' formed using $X$. For example, if $X=
\{a,b\}$, then $(a*b)*a,(b*a)*a,a*(b*a),(a*(a*b))*b, (a*a)*(b*b)$ are 
pairwise distinct elements of~$M(X)$; and all elements of $M(X)$ are 
obtained by forming such products.

(a) Show the following \nm{universal property of the free magma}. For 
any magma~$(N,\nu)$ and any map $\varphi\colon X \rightarrow N$, there 
exists a unique map $\hat{\varphi}\colon M(X) \rightarrow N$ such that 
$\hat{\varphi}|_X=\varphi$ and $\hat{\varphi}$ is a magma homomorphism 
(meaning that $\hat{\varphi}(w*w')=\nu(\hat{\varphi}(w),\hat{\varphi}(w'))$ 
for all $w,w'\in M(X)$).

(b) As in Example~\ref{magma}, let $F_k(X):=k[M(X)]$ be the magma algebra 
over $k$ of the free magma $M(X)$. Note that, as an algebra, $F_k(X)$ is 
generated by $\{\varepsilon_x \mid x\in M(X)\}$. We denote the product of 
two elements $a,b\in F_k(X)$ by $a\cdot b$. Let $I$ be the ideal of
$F_k(X)$ which is generated by all elements of the form
\[a\cdot a \qquad \mbox{or}\qquad a\cdot (b\cdot c)+ b\cdot (c\cdot a)+ 
c\cdot (a\cdot b),\]
for $a,b,c\in F_k(X)$. (Thus, $I$ is the intersection of all ideals
of $F_k(X)$ that contain the above elements.) Let $L(X):=F_k(X)/I$ and 
$\iota \colon X\rightarrow L(X)$, $x\mapsto \varepsilon_x+I$. Show that 
$L(X)$ is a Lie algebra over $k$ which has the following \nmi{universal 
property}{universal property of the free Lie algebra}. 
For any Lie algebra $L'$ over $k$ and any map $\varphi\colon 
X\rightarrow L'$, there exists a unique Lie algebra homomorphism 
$\hat{\varphi}\colon L(X)\rightarrow L'$ such that $\varphi=
\hat{\varphi}\circ \iota$. Deduce that $\iota$ is injective.

The Lie algebra $L(X)$ is called the \nm{free Lie algebra} over~$X$.
By taking factor algebras of $L(X)$ by an ideal, we can construct Lie 
algebras in which prescribed relations hold. (See, e.g.,
Exercise~\ref{expsolv2}.)
\end{xca}

\section{Matrix Lie algebras and derivations} \label{sec02}

We have just seen that every associative algebra can be turned into
a Lie algebra. This leads to the following concrete examples.

\begin{exmp} \label{deflie2} Let $V$ be a vector space. Then 
$\End(V)$ denotes as usual the vector space of all linear maps 
$\varphi\colon V\rightarrow V$. In fact, $\End(V)$ is an associative 
algebra where the product is given by the composition of maps; the identity
map $\id_V\colon V\rightarrow V$ is the identity element for this product. 
Applying the construction in Example~\ref{deflie1}, we obtain a bracket 
on $\End(V)$ and so $\End(V)$ becomes a Lie algebra, denoted $\gl(V)$. 
Thus, $\gl(V)=\End(V)$ as vector spaces and 
\[ [\varphi,\psi]=\varphi\circ \psi -\psi\circ \varphi\qquad
\mbox{for all $\varphi,\psi\in \gl(V)$}.\]
Now assume that $\dim V<\infty$ and let $B=\{v_i\mid i\in I\}$ be a basis
of~$V$. We denote by $M_I(k)$ the algebra of all matrices with entries
in $k$ and rows and columns indexed by $I$, with the usual matrix 
product. For $\varphi\in \End(V)$, we denote by $M_B(\varphi)$ the matrix of 
$\varphi$ with respect to~$B$; thus, $M_B(\varphi)=(a_{ij})_{i,j\in I}
\in M_I(k)$ where $\varphi(v_j)=\sum_{i\in I} a_{ij}v_i$ for 
all~$j$. Now applying the construction in Example~\ref{deflie1}, we 
obtain a bracket on $M_I(k)$ and so $M_I(k)$ also becomes a Lie algebra, 
denoted $\gl_I(k)$. Thus, $\gl_I(k)=M_I(k)$ as vector spaces and 
\[ [X,Y]=X\cdot Y-Y\cdot X\qquad \mbox{for all $X,Y\in \gl_I(k)$}.\]
The map $\varphi\mapsto M_B(\varphi)$ defines an isomorphism of associative 
algebras $\End(V)\cong M_I(k)$. Consequently, this map also defines an
isomorphism of Lie algebras $\gl(V)\cong \gl_I(k)$. (Of course, if 
$I=\{1,\ldots,n\}$ where $n=\dim V$, then we write as usual $M_n(k)$ 
and $\gl_n(k)$ instead of $M_I(k)$ and $\gl_I(k)$, respectively.)
\end{exmp}

\begin{exmp} \label{deflie2a} Let $\gl(V)$ be as in the previous 
example, where $\dim V<\infty$. Then consider the map $\trc\colon \gl(V)
\rightarrow k$ which sends each $\varphi\in \gl(V)$ to the trace of 
$\varphi$ (that is, the sum of the diagonal entries of $M_B(\varphi)$, for 
some basis $B=\{v_i\mid i\in I\}$ of $V$). Since $\trc(\varphi\circ \psi)
=\trc(\psi \circ\varphi)$ for all $\varphi,\psi\in\gl(V)$, we deduce that 
\[\slm(V):=\{\varphi\in\gl(V)\mid \trc(\varphi)=0\}\]
is a Lie subalgebra of $\gl(V)$. (Note that $\slm(V)$ is not a subalgebra
with respect to the matrix product!) Considering matrices as above, we have 
analogous definitions of $\slm_I(k)$ and $\slm_n(k)$ (where $I=
\{1,\ldots,n\}$).
\end{exmp}

\begin{xca} \label{triplephi} Let $V$ be a vector space and $L=\gl(V)$, 
with Lie bracket as in Example~\ref{deflie2}. Show that, for 
$\varphi_1,\varphi_2,\varphi_3 \in L$, we have
\[ [\varphi_1,\varphi_2\circ \varphi_3]=[\varphi_1,\varphi_2]
\circ \varphi_3 +\varphi_2\circ [\varphi_1,\varphi_3].\]
\end{xca}

\begin{xca} \label{expsolv1} Let $L$ be a Lie algebra. If $\dim L=1$, 
then $L$ is clearly abelian. Now assume that $\dim L=2$ and that $L$ is
not abelian. Show that $L$ has a basis $\{x,y\}$ such that $[x,y]=y$; 
in particular, $\langle y\rangle_k$ is a non-trivial ideal of $L$ and so
$L$ is not simple. Show that $L$ is isomorphic to the following Lie 
subalgebra of $\gl_2(k)$:
\[\renewcommand{\arraystretch}{0.9} \left\{\left(
\begin{array}{c@{\hspace{6pt}}c} a & b \\ 0 & 0 \end{array}\right)
\;\Big|\; a,b\in k\right\}.\]
In particular, if $L$ is a simple Lie algebra, then $\dim L\geq 3$.
\end{xca}

\begin{xca} \label{xcanilp} This is a reminder of a basic result from Linear
Algebra. Let $V$ be a vector space and $\varphi\colon V\rightarrow V$ be a 
linear map. Let $v\in V$. We say that $\varphi$ is \nm{locally nilpotent} 
at $v$ if there exists some $d\geq 1$ (which may depend on~$v$) such that
$\varphi^d(v)=0$. We say that $\varphi$ is \textit{nilpotent} if 
$\varphi^d=0$ for some $d\geq 1$. Assume now that $\dim V<\infty$.

(a) Let $X\subseteq V$ be a subset such that $V=\langle X\rangle_k$. 
Assume that $\varphi$ is locally nilpotent at every $v\in X$.
Show that $\varphi$ is nilpotent.

(b) Show that, if $\varphi$ is nilpotent, then there is a basis $B$ of 
$V$ such that the matrix of $\varphi$ with respect to $B$ is triangular 
with $0$ on the diagonal; in particular, we have $\varphi^{\dim V}=0$ 
and the trace of $\varphi$ is~$0$.
\end{xca}

\begin{exmp} \label{defadj} Let $L$ be a Lie algebra. In analogy
to Remark~\ref{adjoint0} and Example~\ref{deflie1}(b), we define for 
$x\in L$ the linear map 
\[\ad_L(x) \colon L\rightarrow L, \qquad y\mapsto [x,y].\]
Hence, we obtain a linear map $\ad_L\colon L\rightarrow \End(L)$, 
$x\mapsto \ad_L(x)$. By the Jacobi identity and anti-symmetry, we have 
\begin{align*}
\ad_L([x,y])(z)&=[[x,y],z]=-[z,[x,y]]\\
&=[x,[y,z]]+[y,[z,x]]=[x,[y,z]]-[y,[x,z]]\\
& =(\ad_L(x)\circ \ad_L(y)-\ad_L(y)\circ \ad_L(x))(z)
\end{align*}
for all $z\in L$ and so $\ad_L([x,y])=[\ad_L(x),\ad_L(y)]$.  Thus,
we obtain a Lie algebra homomorphism $\ad_L\colon L\rightarrow \gl(L)$.
(See also Example~\ref{defmod2} below.) The kernel of $\ad_L$ is called 
the \nm{center} of $L$ and will be denoted by $Z(L)$; thus, $Z(L)$ is an
ideal of $L$ and 
\[ Z(L)=\ker(\ad_L)=\{x\in L\mid [x,y]=0 \mbox{ for all $y\in L$}\}.\]
Finally, for $x,y,z\in L$, we also have the identity
\begin{align*}
\ad_L(z)&([x,y])=[z,[x,y]]=-[x,[y,z]]-[y,[z,x]]\\&=[x,[z,y]]+[[z,x],y]
=[x,\ad_L(z)(y)]+[\ad_L(z)(x),y]
\end{align*} 
which shows that $\ad_L(z)$ is a derivation in the following sense.
\end{exmp}

\begin{defn} \label{defder} Let $A$ be an algebra. A linear map
$\delta\colon A\rightarrow A$ is called a \nm{derivation} if $\delta(x
\cdot y)=x\cdot \delta(y)+\delta(x)\cdot y$ for all $x,y\in A$. Let 
$\Der(A)$ be the set of all derivations of~$A$. One immediately checks 
that $\Der(A)$ is a subspace of $\End(A)$.
\end{defn}

\begin{xca} \label{xcader} Let $A$ be an algebra.

(a) Show that $\Der(A)$ is a Lie subalgebra of $\gl(A)$.

(b) Let $\delta\colon A \rightarrow A$ be a derivation. Show
that, for any $n\geq 0$, we have the \nm{Leibniz rule}
\[\quad\qquad \delta^n(x\cdot y)=\sum_{i=0}^n \binom{n}{i} \delta^i(x)\cdot
\delta^{n-i}(y) \quad \mbox{for all $x,y\in A$}.\]

\noindent Derivations are a source for Lie algebras which do not arise from 
associative algebras as in Example~\ref{deflie1}; see Example~\ref{wittalg} 
below. The following construction with nilpotent derivations will play a 
major role in Chapter~\ref{chap3}; see also Exercises~\ref{expocomm}
and~\ref{xcaexp}.
\end{xca}

\begin{lem} \label{exponential} Let $A$ be an algebra where the ground 
field $k$ has characteristic~$0$. If $d\colon A\rightarrow A$ is a 
derivation such that $d^n=0$ for some $n>0$ (that is, $d$ is nilpotent), 
we obtain a map 
\[ \exp(d)\colon A\rightarrow A, \qquad x\mapsto \sum_{0\leq i<n} 
\frac{d^i(x)}{i!}=\sum_{i\geq 0}\frac{d^i(x)}{i!}.\]
Then $\exp(d)$ is an algebra isomorphism, with inverse $\exp(-d)$.
\end{lem}

\begin{proof} Since $d^i$ is linear for all $i\geq 0$, it is clear that
$\exp(d)\colon A\rightarrow A$ is a linear map. For $x,y\in A$, we have
\begin{align*}
\exp(d)(x)&\cdot \exp(d)(y)=\Bigl(\sum_{i\geq 0}\frac{d^i}{i!}(x)\Bigr)
\cdot \Bigl(\sum_{j\geq 0} \frac{d^j}{j!}(y)\Bigr)\\
&=\sum_{i,j\geq 0}\frac{d^i}{i!}(x)\cdot \frac{d^j}{j!}(y)
=\sum_{m\geq 0} \Bigl(\sum_{\atop{i,j\geq 0}{i+j=m}}\frac{d^i}{i!}(x)
\cdot \frac{d^j}{j!}(y)\Bigr)\\
&=\sum_{m\geq 0} \frac{1}{m!}\Bigl(\sum_{0\leq i\leq m}\binom{m}{i}
d^i(x) \cdot d^{m-i}(y)\Bigr)=\sum_{m\geq 0} \frac{d^m}{m!}(x\cdot y),
\end{align*}
where the last equality holds by the Leibniz rule. Hence, the right side
equals $\exp(d)(x\cdot y)$. Thus, $\exp(d)$ is an algebra homomorphism. 

Now, we can also 
form $\exp(-d)$ and $\exp(0)$, where the definition immediately shows that 
$\exp(0)=\id_A$. So, for any $x\in A$, we obtain:
\[ x=\exp(0)(x)=\exp(d{+}(-d))(x)=\sum_{m\geq 0} 
\frac{(d{+}(-d))^m(x)}{m!}.\]
Since $d$ and $-d$ commute with each other, we can apply the binomial 
formula to $(d+(-d))^m$. So the right hand side evaluates to
\begin{align*}
&\sum_{m\geq 0} \frac{1}{m!}\sum_{\atop{i,j \geq 0}{i+j=m}} 
\frac{m!}{i!\,j!} (d^i{\circ}(-d)^j)(x)=\sum_{i,j\geq 0} 
\frac{(d^i{\circ} (-d)^j)(x)}{i!\,j!}\\&\qquad =\sum_{i,j\geq 0}
\frac{d^i}{i!}\Bigl( \frac{(-d)^j }{j!}(x)\Bigr)=\sum_{i\geq 0}
\frac{d^i}{i!}\Bigl( \sum_{j\geq 0}\frac{(-d)^j }{j!}(x)\Bigr)\\
&\qquad\qquad =\sum_{i\geq 0}\frac{d^i}{i!}\bigl(\exp(-d)(x)\bigr)
=\exp(d)\bigl(\exp(-d)(x)\bigr).
\end{align*}
Hence, we see that $\exp(d)\circ \exp(-d)=\id_A$; similarly, 
$\exp(-d)\circ \exp(d)=\id_A$. So $\exp(d)$ is invertible, with 
inverse $\exp(-d)$.
\end{proof}

\begin{exmp} \label{wittalg} Let $A=k[T,T^{-1}]$ be the algebra of
Laurent polynomials in the indeterminate~$T$. Let us determine $\Der(A)$. 
Since~$A=\langle T,T^{-1}\rangle_{\text{alg}}$, the product rule for 
derivations implies that every $\delta\in \Der(A)$ is uniquely 
determined by $\delta(T)$ and $\delta(T^{-1})$. Now $\delta(1)=\delta(T
\cdot T^{-1})=T\delta(T^{-1})+\delta(T)T^{-1}$. Since $\delta(1)=
\delta(1)+\delta(1)$, we have $\delta(1)=0$ and so $\delta(T^{-1})=
-T^{-2}\delta(T)$. Hence, we conclude: 
\begin{equation*}
\mbox{Every $\delta\in\Der(A)$ is uniquely determined by its 
value $\delta(T)$}.\tag{a}
\end{equation*}
For $m\in \Z$ we define a linear map $L_m\colon A\rightarrow A$ by 
\[ L_m(f)=-T^{m+1}D(f) \qquad \mbox{for all $f\in A$},\]
where $D\colon A\rightarrow A$ denotes the usual formal derivate
with respect to $T$, that is, $D$ is linear and $D(T^n)=nD(T^{n-1})$ for
all $n\in \Z$. Now $D\in \Der(A)$ (by the product rule for formal 
derivates) and so $L_m\in \Der(A)$. We have $L_m(T)=-T^{m+1}D(T)=
-T^{m+1}$. Hence, if $\delta\in \Der(A)$ and $\delta(T)=\sum_i a_iT^i$
with $a_i\in k$, then $-\delta$ and the sum $\sum_i a_iL_{i-1}$ have
the same value on $T$. So $-\delta$ must be equal to that sum by (a). 
Thus, we have shown that 
\begin{equation*}
\Der(A)=\langle L_m\mid m\in\Z\rangle_k.\tag{b}
\end{equation*}
In fact, $\{L_m\mid m\in\Z\}$ is a basis of $\Der(A)$. (Just apply
a linear combination of the $L_m$'s to $T$ and use the fact that 
$L_m(T)=-T^{m+1}$.) Now let $m,n\in\Z$. Using the bracket in $\gl(A)$, 
we obtain that 
\[[L_m,L_n](T)=(L_m\circ L_n-L_n\circ L_m)(T)= \ldots =(n-m)T^{m+n+1},\]
which is also the result of $(m-n)L_{m+n}(T)$. By Exercise~\ref{xcader}(a),
we have $[L_m,L_n]\in \Der(A)$. So (a) shows again that 
\begin{equation*}
[L_m,L_n]=(m-n)L_{m+n} \quad \mbox{for all $m,n\in\Z$}.\tag{c}
\end{equation*}
Thus, $\Der(A)$ is an infinite-dimensional Lie subalgebra of $\gl(A)$,
with basis $\{L_m\mid m\in\Z\}$ and bracket determined as above; this
Lie algebra is called a \nm{Witt algebra} (or \textit{centerless 
Virasoro algebra}; see also the notes at the end of this chapter). 
\end{exmp} 

\begin{prop} \label{witt2} Let $L=\Der(A)$ be the Witt algebra in
Example~\ref{wittalg}. If $\chark(k)=0$, then $L$ is a simple Lie algebra.
\end{prop}

\begin{proof} Let $I\subseteq L$ be a non-zero ideal and $0\neq x\in I$. 
Then we can write $x=c_1L_{m_1}+\ldots +c_rL_{m_r}$ where $r\geq 1$, 
$m_1<\ldots <m_r$ and all $c_i\in k$ are non-zero. Choose $x$ such 
that $r$ is as small as possible. We claim that $r=1$. Assume, if 
possible, that $r\geq 2$. Since $[L_0,L_m]=-mL_m$ for all $m\in \Z$, 
we obtain that $[L_0,x]=-c_1m_1L_{m_1}- \ldots -c_rm_rL_{m_r}
\in I$. Hence, 
\[ m_rx+[L_0,x]=c_1(m_r-m_1)L_{m_1}+\ldots +c_{r-1}(m_r-m_{r-1})
L_{m_{r-1}}\]
is a non-zero element of $I$, contradiction to the minimality of $r$.
Hence, $r=1$ and so $L_{m_1}\in I$. Now $[L_{m-m_1},L_{m_1}]=
(m-2m_1)L_m$ and so $L_m \in I$ for all $m\in\Z$, $m\neq 2m_1$. But
$[L_{m_1+1},L_{m_1-1}]=2L_{2m_1}$ and so we also have $L_{2m_1}\in I$. 
Hence, we do have $I=L$, as desired.
\end{proof}

\begin{xca} \label{expsolv2} Let $L=\slm_2(k)$, as in 
Example~\ref{deflie2a}. Then $\dim L=3$ and $L$ has a basis $\{e,h,f\}$ 
where 
\[\renewcommand{\arraystretch}{0.9} 
e=\left(\begin{array}{cc} 0 & 1 \\ 0 & 0\end{array}\right),\qquad 
h=\left(\begin{array}{c@{\hspace{5pt}}r} 1 & 0 \\ 0 & -1\end{array}
\right), \qquad f=\left(\begin{array}{cr} 0 & 0 \\ 1 & 0\end{array}
\right).\]
(a) Check that $[e,f]=h$, $[h,e]=2e$, $[h,f]=-2f$. Show that $L$ is 
simple if $\chark(k)\neq 2$. What happens if $\chark(k)=2$?
Consider also the Lie algebra $L'$ in Example~\ref{xcagrass}. Is 
$L'\cong \slm_2(\R)$? Is $L'$ simple? What happens if we work with $\C$ 
instead of $\R$?

\noindent (b)  Let $\hat{L}$ be the free Lie algebra over the set
$X=\{E,H,F\}$; see Exercise~\ref{xcafree}. Let $I\subseteq \hat{L}$ be the
ideal generated by $[E,F]-H$, $[H,E]-2E$, $[H,F]+2F$ (that is, the 
intersection of all ideals containing those elements). By the universal
property, there is a unique homomorphism of Lie algebras $\varphi\colon
\hat{L} \rightarrow L$ such that $\varphi(E)=e$, $\varphi(H)=h$ and 
$\varphi(F)=f$. By (a), we have $I\subseteq \ker(\varphi)$. Show that
the induced homomorphism $\bar{\varphi}\colon \hat{L}/I\rightarrow L$
is an isomorphism. 
 \end{xca}

\begin{xca} \label{xcazgl} Show that $Z(\gl_n(k))=\{a I_n 
\mid a \in k\}$ (where $I_n$ denotes the $n\times n$-identity matrix). 
What happens for $Z(\slm_n(k))$?
\end{xca}

\begin{xca} \label{xcasemidir} This exercise describes a useful method
for constructing new Lie algebras out of two given ones. So let $S,I$ be 
Lie algebras over $k$ and $\theta \colon S\rightarrow \Der(I)$, $s\mapsto
\theta_s$, be a homomorphism of Lie algebras. Consider the vector space 
$L=S\times I= \{(s,x)\mid s\in S,x\in I\}$ (with component-wise defined 
addition and scalar multiplication). For $s_1,s_2 \in S$ and $x_1,x_2\in I$ 
we define 
\[[(s_1,x_1),(s_2,x_2)]:= \bigl([s_1,s_2],[x_1,x_2]+\theta_{s_1}(x_2)-
\theta_{s_2}(x_1)\bigr).\]
Show that $L$ is a Lie algebra such that $L=\underline{S}\oplus 
\underline{I}$, where 
\begin{align*}
\underline{S}&:=\{(s,0) \mid s \in S\} \subseteq L\quad\mbox{is a 
subalgebra,}\\
\underline{I}&:=\{(0,x) \mid x\in I\} \subseteq L\quad\mbox{is an ideal}.
\end{align*} 
We also write $L=S\ltimes_\theta I$ and call $L$ the \nm{semidirect 
product} of $I$ by $S$ (via $\theta$). If $\theta(s)=0$ 
for all $s\in S$, then $[(s_1,x_1),(s_2,x_2)]=([s_1,s_2],[x_1,x_2])$ for 
all $s_1,s_2\in S$ and $x_1,x_2\in I$. Hence, in this case, $L$ is the 
\nm{direct product} of $S$ and~$I$, as in Example~\ref{dirprodalg}.
\end{xca}

\begin{xca} \label{expocomm} Let $A$ be an algebra where the ground
field $k$ has characteristic~$0$. Let $d\colon A\rightarrow A$ and 
$d'\colon A\rightarrow A$ be nilpotent derivations such that
$d\circ d'=d'\circ d$. Show that $d+d'$ also is a nilpotent derivation
and that $\exp(d+d')=\exp(d)\circ \exp(d')$.
\end{xca}

\begin{xca} \label{xcaexp} This exercise gives a first outlook to some
constructions that will be studied in much greater depth and generality
in Chapter~\ref{chap3}. Let $L\subseteq \gl(V)$ be a Lie subalgebra,
where $V$ is a finite-dimesional $\C$-vector space. Let $\Aut(L)$ be the
group of all Lie algebra automorphisms of~$L$ (that is, linear maps 
$L\rightarrow L$ which are Lie algebra isomorphisms). 


(a) Assume that $x\in L$ is nilpotent (as linear map $x\colon V\rightarrow 
V$). Then show that the linear map $\ad_L(x)\colon L\rightarrow L$ is
nilpotent. (Hint: use the ``trick'' in Remark~\ref{adjoint0}.) Is the
converse also true?

(b) Let $L=\slm_2(\C)$ with basis elements $e,h,f$ as in
Exercise~\ref{expsolv2}. Note that $e$ and $f$ are nilpotent matrices.
Hence, by (a), the derivations $\ad_L(e)\colon L \rightarrow L$ and
$\ad_L(f) \colon L \rightarrow L$ are nilpotent. Consequently, $t\,
\ad_L(e)$ and $t\, \ad_L(f)$ are nilpotent derivations for all $t \in 
\C$. By Lemma~\ref{exponential}, we obtain Lie algebra automorphisms
\[ \exp\bigl(t\,\ad_L(e)\bigr)\colon L \rightarrow L \qquad \mbox{and}
\qquad \exp\bigl(t\,\ad_L(f)\bigr)\colon L \rightarrow L; \]
we will denote these by $x(t)$ and $y(t)$, respectively. Determine the
matrices of these automorphisms with respect to the basis $\{e,h,f\}$
of~$L$. Check that $x(t+t')=x(t)x(t')$ and $y(t+t')=y(t)y(t')$ for all
$t,t'\in \C$. The subgroup $G:=\langle x(t),y(t')\mid t,t'\in \C\rangle 
\subseteq \Aut(L)$ is called the \nm{Chevalley group} associated with 
$L=\slm_2(\C)$. 
\end{xca}

\section{Solvable and semisimple Lie algebras} \label{sec03}

Let $A$ be an algebra. If $U,V\subseteq A$ are subspaces, then we denote
\[ U\cdot V:=\langle u\cdot v\mid u\in U,v\in V\rangle_k\subseteq A.\]
In general, $U\cdot V$ will only be a subspace of $A$, even
if $U$, $V$ are subalgebras or ideals. On the other hand, taking 
$U=V=A$, then 
\[A^2:=A\cdot A=\langle x\cdot y\mid x,y\in A\rangle_k \] 
clearly is an ideal of $A$, and the induced product on $A/A^2$ is identically
zero. So we can iterate this process: Let us set $A^{(0)}:=A$ and then 
\[A^{(1)}:=A^2, \quad A^{(2)}:=(A^{(1)})^2, \quad A^{(3)}:=(A^{(2)})^2, 
\quad \ldots.\]
Thus, we obtain a chain of subalgebras $A=A^{(0)}\supseteq A^{(1)}
\supseteq A^{(2)}\supseteq \ldots$ such that $A^{(i+1)}$ is an ideal in
$A^{(i)}$ for all $i$ and the induced product on $A^{(i)}/A^{(i+1)}$ is
identically zero. An easy induction on $j$ shows that 
$A^{(i+j)}=(A^{(i)})^{(j)}$ for all $i,j\geq 0$.

\begin{defn} \label{defsolv} We say that $A$ is a \nm{solvable algebra} 
if $A^{(m)}=\{0\}$ for some $m\geq 0$ (and, hence, $A^{(l)}=\{0\}$ for 
all $l\geq m$.)
\end{defn}

Note that the above definitions are only useful if $A$ does not have an
identity element which is, in particular, the case for Lie algebras by the 
anti-symmetry condition in Definition~\ref{deflie}. 

\begin{exmp} \label{defborell} (a) All Lie algebras of dimension $\leq 2$ 
are solvable; see Exercise~\ref{expsolv1}. 

(b) Let $n\geq 1$ and $\fb_n(k)\subseteq \gl_n(k)$ be the subspace 
consisting of all upper triangular matrices, that is, all $(a_{ij})_{1
\leq i,j\leq n} \in \gl_n(k)$ such that $a_{ij}=0$ for all $i>j$. Since 
the product of two upper triangular matrices is again upper 
triangular, it is clear that $\fb_n(k)$ is a Lie subalgebra of $\gl_n(k)$. 
An easy matrix calculation shows that $\fb_n(k)^{(1)}=[\fb_n(k),\fb_n(k)]$ 
consists of upper triangular matrices with $0$ on the diagonal.
More generally, $\fb_n(k)^{(r)}$ for $1\leq r\leq n$ consists of upper 
triangular matrices $(a_{ij})$ such that $a_{ij}=0$ for all $i\leq j<i+r$.
In particular, we have $\fb_n(k)^{(n)}=\{0\}$ and so $\fb_n(k)$ is solvable. 
\end{exmp}

\begin{xca} \label{excsolv} For a fixed $0\neq \delta\in k$, we 
define
\[ \renewcommand{\arraystretch}{0.8}
L_\delta:=\left\{\left(\begin{array}{ccc} a & b & 0 \\ 
0 & 0 & 0 \\ 0 & c & a\delta \end{array}\right) \;\Big|\; a,b,c\in 
k\right\}\subseteq \gl_3(k).\]
Show that $L_\delta$ is a solvable Lie subalgebra of $\gl_3(k)$, where
$[L_\delta,L_\delta]$ is abelian. Show that, if $L_\delta \cong 
L_{\delta'}$, then $\delta=\delta'$ or $\delta^{-1}=\delta'$. Hence, if 
$|k|=\infty$, then there are infinitely many pairwise non-isomorphic 
solvable Lie algebras of dimension~$3$. (See \cite[Chap.~3]{EW} for a 
further discussion of ``low-dimensional'' examples of solvable Lie 
algebras.)

\smallskip
\noindent \footnotesize{[{\it Hint.} A useful tool to check that two Lie 
algebras cannot be isomorphic is as follows. Let $L_1,L_2$ be 
finite-dimensional Lie algebras over $k$. Let $\varphi\colon L_1\rightarrow 
L_2$ be an isomorphism. Show that $\varphi\circ \ad_{L_1}(x)=\ad_{L_2}
(\varphi(x)) \circ \varphi$ for $x\in L_1$. Deduce that $\ad_{L_1}(x)\colon 
L_1 \rightarrow L_1$ and $\ad_{L_2}(\varphi(x)) \colon L_2 \rightarrow 
L_2$ must have the same characteristic polynomial. Try to apply this with 
the element $x\in L_\delta$ where $a=1$, $b=c=0$.]}
\end{xca}

\begin{xca} \label{heisenb} Let $L$ be a Lie algebra over $k$ with 
$\dim L=2n+1$, $n\geq 1$. Suppose that $L$ has a basis $\{z\}\cup
\{e_i,f_i\mid 1\leq i\leq n\}$ such that $[e_i,f_i]=z$ and $[z,e_i]=[z,f_i]
=0$ for $1\leq i \leq n$. Then $L$ is called a \nm{Heisenberg Lie algebra} 
(see \cite[\S 1.4]{MP} or \cite[\S 1.7]{FLM} for further background). Check 
that $[L,L]=Z(L)=\langle z\rangle_k$; in particular, $L$ is solvable. 
Show that, for $n=1$,
\[\renewcommand{\arraystretch}{0.8}
 L:=\left\{\left(\begin{array}{ccc} 0  & a & b \\ 
0 & 0 & c \\ 0 & 0 & 0 \end{array}\right) \;\Big|\; a,b,c\in 
k\right\}\subseteq \gl_3(k)\]
is a Heisenberg Lie algebra; find a basis $\{z\}\cup \{e_1,f_1\}$ as above.
\end{xca}
 
\begin{lem} \label{defsolv1} Let $A$ be an algebra.
\begin{itemize}
\item[(a)] Let $B$ be an algebra and $\varphi\colon A\rightarrow B$ be a
surjective algebra homomorphism. Then $\varphi(A^{(i)})=B^{(i)}$ for
all $i\geq 0$. 
\item[(b)] Let $B\subseteq A$ be a subalgebra. Then $B^{(i)}\subseteq 
A^{(i)}$ for all $i\geq 0$.
\item[(c)] Let $I\subseteq A$ be an ideal. Then $A$ is solvable if 
and only if $I$ and $A/I$ are solvable.
\end{itemize}
\end{lem}

\begin{proof} (a) Induction on $i$. If $i=0$, then this holds by
assumption. Let $i\geq 0$. Then $\varphi(A^{(i+1)})=\varphi(A^{(i)}
\cdot A^{(i)})=\langle \varphi(x)\cdot \varphi(y)\mid x,y\in A^{(i)}
\rangle_k$ which equals $B^{(i)}\cdot B^{(i)}$ since 
$\varphi(A^{(i)})=B^{(i)}$ by induction. 

(b) Induction on $i$. If $i=0$, then this is clear. Now let $i\geq 0$. 
By induction, $B^{(i)}\subseteq A^{(i)}$ and so $B^{(i+1)}=
(B^{(i)})^2 \subseteq (A^{(i)})^2=A^{(i+1)}$.

(c) If $A$ is solvable, then $I$ and $A/I$ are solvable by (a), (b).
Conversely, let $m,l\geq 0$ be such that $I^{(l)}=\{0\}$ and 
$(A/I)^{(m)}=\{0\}$. Let $\varphi\colon A\rightarrow A/I$ be the canonical 
map. Then $\varphi(A^{(m)})=(A/I)^{(m)}=\{0\}$ by (a), hence, $A^{(m)}
\subseteq \ker(\varphi)=I$. Using (b), we obtain $A^{(m+l)}=(A^{(m)})^{(l)}
\subseteq I^{(l)}=\{0\}$ and so $A$ is solvable.
\end{proof}

\begin{cor} \label{defsolv2} Let $A$ be an algebra with $\dim A<\infty$.
Then the set of all solvable ideals of $A$ is non-empty and contains a
unique maximal element (with respect to inclusion). This unique 
maximal solvable ideal will be denoted $\rad(A)$ and called the 
\nm{radical} of $A$. We have $\rad(A/\rad(A))=\{0\}$.
\end{cor}

\begin{proof} First note that $\{0\}$ is a solvable ideal of $A$. Now let 
$I\subseteq A$ be a solvable ideal such that $\dim I$ is as large as 
possible. Let $J\subseteq A$ be another solvable ideal. Clearly, 
$B:=\{x+y\mid x\in I,y\in J\} \subseteq A$ also is an ideal. We claim that 
$B$ is solvable. Indeed, we have $I\subseteq B$ and so $I$ is a solvable 
ideal of~$B$; see Lemma~\ref{defsolv1}(b). Let $\varphi\colon B\rightarrow
B/I$ be the canonical map. By restriction, we obtain an algebra 
homomorphism $\varphi'\colon J\rightarrow B/I$, $x\mapsto x+I$. By the 
definition of $B$, this map is surjective. Hence, since $J$ is solvable, 
then so is $B/I$ by Lemma~\ref{defsolv1}(a). But then $B$ itself is
solvable by Lemma~\ref{defsolv1}(c). Hence, since $\dim I$ was maximal,
we must have $B=I$ and so $J\subseteq I$. Thus, $I=\rad(A)$ is the 
unique maximal solvable ideal of $A$. 

Now consider $B:=A/\rad(A)$ and the canonical map $\varphi \colon A
\rightarrow B$. Let $J\subseteq B$ be a solvable ideal. Then
$\varphi^{-1}(J)$ is an ideal of $A$ containing $\rad(A)$. Now 
$\varphi^{-1}(J)/\rad(A)\cong J$ is solvable. Hence, $\varphi^{-1}(J)$ 
itself is solvable by Lemma~\ref{defsolv1}(c). So $\varphi^{-1}(J)=
\rad(A)$ and $J=\{0\}$. 
\end{proof}

Now let $L$ be a Lie algebra with $\dim L<\infty$. 

\begin{defn} \label{defsolv3} 
We say that $L$ is a \nm{semisimple Lie algebra} 
if $\rad(L)=\{0\}$. By Corollary~\ref{defsolv2}, $L$ itself
or $L/\rad(L)$ is semisimple. 
\end{defn}

Note that $L=\{0\}$ is considered to be semisimple. Clearly, simple Lie 
algebras are semisimple. For example, $L=\slm_2(\C)$ is semisimple. 

\begin{rem} \label{defsolv3a}  Since the center $Z(L)$ is an abelian 
ideal of $L$, we have $Z(L)\subseteq \rad(L)$. Hence, if $L$ is semisimple, 
then $Z(L)=\{0\}$ and so the homomorphism $\ad_L\colon L\rightarrow 
\gl(L)$ in Example~\ref{defadj} is injective. Thus, if $L$ is semisimple 
and $n=\dim L$, then $L$ is isomorphic to a Lie subalgebra of 
$\gl_n(k) \cong \gl(L)$.
\end{rem}

\begin{lem} \label{liesolv} Let $H\subseteq L$ be an ideal. Then
$H^{(i)}$ is an ideal of $L$ for all $i\geq 0$. In particular, if 
$H\neq \{0\}$ is solvable, then there exists a non-zero abelian 
ideal $I\subseteq L$ with $I\subseteq H$.
\end{lem}

\begin{proof} To show that $H^{(i)}$ is an ideal for all $i$, we use 
induction on~$i$. If $i=0$, then $H^{(0)}=H$ is an ideal of $L$ by 
assumption. Now let $i\geq 0$; we have $H^{(i+1)}=[H^{(i)},
H^{(i)}]$. So we must show that $[z,[x,y]]\in [H^{(i)},H^{(i)}]$ and 
$[[x,y],z] \in [H^{(i)},H^{(i)}]$, for all $x,y \in H^{(i)}$, $z \in L$. 
By anti-symmetry, it is enough to show this for $[z,[x,y]]$. 
By induction, $[z,x]\in H^{(i)}$ and $[z,y]\in H^{(i)}$.
Using anti-symmetry and the Jacobi identity, $[z,[x,y]]=-[x,[y,z]]-
[y,[z,x]]\in [H^{(i)},H^{(i)}]$, as required.

Now assume that $H=H^{(0)}\neq \{0\}$ is solvable. So there is some
$m>0$ such that $I:=H^{(m-1)}\neq \{0\}$ and $I^2=H^{(m)}=\{0\}$. We 
have just seen that $I$ is an ideal of $L$, which is abelian since 
$I^2=\{0\}$. 
\end{proof}

By Lemma~\ref{liesolv}, $L$ is semisimple if and only if $L$ has no 
non-zero abelian ideal: This is the original definition of semisimplicity 
given by Killing. A further criterion for checking if $L$ is semisimple 
is given by considering a certain bilinear form on $L$, defined as follows. 

\begin{defn} \label{defkill0} The \nm{Killing form} on $L$ is the 
map defined by 
\[ \kappa_L \colon L\times L\rightarrow k, \qquad (x,y)\mapsto
\trc\bigl(\ad_L(x)\circ \ad_L(y)\bigr).\]
It is clear that $\kappa_L$ is bilinear; it is a \nm{symmetric bilinear 
form} because $\trc(\varphi\circ\psi)=\trc(\psi\circ \varphi)$ for all
linear maps $\varphi,\psi\colon L\rightarrow L$. See also
Section~\ref{sec05} below for basic notions related to bilinear forms.
In particular, for any subset $X\subseteq L$, we define 
\[X^\perp:=\{y\in L \mid \kappa_L(x,y)=0 \mbox{ for all $x\in X$}\};\]
this is a subspace of $L$. We say that $\kappa_L$ is non-degenerate
if $L^\perp=\{0\}$.
\end{defn}

\begin{lem} \label{defkill0a} {\rm (a)} $\kappa_L(x,[y,z])=
\kappa_L([x,y],z)$ for all $x,y,z\in L$.\\
{\rm (b)} If $I\subseteq L$ is an ideal, then $I^\perp\subseteq L$ also is
an ideal.\\
{\rm (c)} If $\kappa_L$ is non-degenerate, then $L$ is 
semisimple\footnote{The converse of (b) also holds but requires more work; 
see, for example, \cite[\S 9.2]{EW}. In the setting of Chapter~\ref{chap2}, 
we will deduce directly the non-degeneracy of~$\kappa_L$.}.
\end{lem}

\begin{proof} (a) Since $\ad_L\colon L \rightarrow \gl(L)$ is a Lie
algebra homomorphism, we have $\ad_L([x,y])=\ad_L(x){\circ} \ad_L(y)-
\ad_L(y){\circ} \ad_L(x)$; similarly, $\ad_L([y,z])=\ad_L(y){\circ} \ad_L(z)-
\ad_L(z){\circ} \ad_L(y)$. \\
This yields that $\kappa_L(x,[y,z])=A-B$ where 
\begin{align*}
A&:=\trc\bigl(\ad_L(x){\circ} (\ad_L(y){\circ} \ad_L(z)\bigr), \\
B&:=\trc\bigl(\ad_L(x){\circ}(\ad_L(z){\circ} \ad_L(y)\bigr).
\end{align*}
Now note that 
\begin{align*}
B&=\trc\bigl((\ad_L(z){\circ} \ad_L(y)){\circ}\ad_L(x)\bigr)\\
&= \trc\bigl((\ad_L(y){\circ} \ad_L(x)){\circ}\ad_L(z)\bigr)
\end{align*} 
Hence, $A-B=\kappa_L([x,y],z)$, as claimed.

(b) Let $I\subseteq L$ be an ideal, $x\in I^\perp$ and $y\in L$.
Then, for any $z\in I$ we have $[y,z] \in I$. So, using (a), we obtain 
$\kappa_L([x,y],z)= \kappa_L(x,[y,z])=0$. Since this holds for
all $z\in I$, we have $[x,y]\in I^\perp$.

(c) Let $A\subseteq L$ be an abelian ideal. For $a\in A$, $x,y\in L$, we have 
\[\bigl(\ad_L(a) \circ \ad_L(x)\bigr)^2(y)=[a,[x,[a,[x,y]]]]=0,\]
since $[a,[x,y]]\in A$ and so $[x,[a,[x,y]]]\in A$. Hence, we conclude
that $(\ad_L(a)\circ \ad_L(x))^2=0$. But then $\kappa_L(a,x)=
\trc(\ad_L(a)\circ \ad_L(x))=0$; see Exercise~\ref{xcanilp}. This 
holds for all $x\in L$. Hence, $a=0$ since $\kappa_L$ is non-degenerate. So 
$\rad(L)=\{0\}$ by Lemma~\ref{liesolv}.
\end{proof}

\begin{exmp} \label{defkill0b} (a) If $L$ is abelian, then $\ad_L(x)=0$
for all $x\in L$ and so $\kappa_L(x,y)=0$ for all $x,y\in L$.

(b) Let $L=\slm_2(\C)$ with basis $\{e,h,f\}$ as in 
Exercise~\ref{expsolv2}. The corresponding matrices of $\ad_L(e)$, 
$\ad_L(h)$, $\ad_L(f)$ are given by
\[\renewcommand{\arraystretch}{0.8}  
\left(\begin{array}{c@{\hspace{3pt}}rc} 0 & -2 & 0 \\ 0 & 0 & 1 \\
0 & 0 & 0 \end{array}\right),\quad 
\left(\begin{array}{cc@{\hspace{3pt}}r} 2 & 0 & 0 \\ 0 & 0 & 0 \\
0 & 0 & -2\end{array}\right),\quad 
\left(\begin{array}{@{\hspace{3pt}}rcc} 0 & 0 & 0 \\ -1 & 0 & 0 \\
0 & 2 & 0 \end{array}\right),\]
respectively. Then $\kappa_L(h,h)=8$, $\kappa_L(e,e)=
\kappa_L(f,f)=0$; furthermore, $\kappa_L(e,f)=4$, $\kappa_L(h,e)=
\kappa_L(h,f)=0$. Hence, the Gram matrix of $\kappa_L$ has determinant 
$-128$ and so $\kappa_L$ is non-degenerate. So Lemma~\ref{defkill0a}
shows once more that $L$ is semisimple.
\end{exmp}

\begin{xca} \label{xcakill} Let $I\subseteq L$ be an ideal  and 
regard $I$ as a Lie algebra by itself. Then show that the Killing form
$\kappa_I\colon I \times I \rightarrow \C$ is equal to the restriction 
of $\kappa_L$ to $I\times I$.\\ 
{\footnotesize  [{\it Hint}. Choose a basis of $I$ and extend it to
a basis $B$ of $L$. Let $x\in L$. How does the matrix of $\ad_L(x)$
with respect to $B$ look like? Similar question for $y \in I$.]}
\end{xca}

The non-degeneracy of the Killing form $\kappa_L$ is a very powerful
tool. For example, the following result shows that, if this is the case,
then~$L$ is built up from simple Lie algebras in a very simple way.

\begin{prop} \label{killing8b} Let $L\neq \{0\}$ and $\kappa_L$ be
non-degenerate.\\
{\rm (a)} If $I\subseteq L$ is an ideal, then $I^\perp\subseteq L$ is
an ideal and $L=I\oplus I^\perp$.\\
{\rm (b)} We have $L=L_1 \oplus \ldots \oplus L_r$ ($r\geq 1$) where 
each $L_i\subseteq L$ is an ideal which is simple as a Lie algebra.
We have $[L_i,L_j]=\{0\}$ for all $i\neq j$.
\end{prop}

\begin{proof} (a) By Lemma~\ref{defkill0a}(b), $I^\perp$ is an ideal. 
Since $\kappa_L$ is non-generate, we also know that $\dim L=\dim I+
\dim I^\perp$. Hence, it will be sufficient to show that $I\cap
I^\perp=\{0\}$. Let $J:=I\cap I^\perp$. Then $\kappa_L(x,y)=0$ for all
$x\in J\subseteq I$ and $y\in J\subseteq I^\perp$. Now fix $x,y\in J$ 
and let $z\in L$. Then, by Lemma~\ref{defkill0a}(a), we have $\kappa_L
([x,y],z)=\kappa_L(x,[y,z])=0$ since $[y,z]\in J$. Since this holds for
all $z\in L$, we must have $[x,y]=0$ since $\kappa_L$ is non-degenerate. 
Hence, $J$ is an abelian ideal. But $L$ is semisimple by 
Lemma~\ref{defkill0a}(c). Hence, $J=\{0\}$.
 
(b) We use induction on $\dim L$. If $L$ itself is simple, then there is
nothing to prove. Now assume that $L$ is not simple and let $\{0\}\neq L_1
\subsetneqq L$ be an ideal such that $\dim L_1$ is as small as possible. 
By (a), we have $L=L_1\oplus L'$ where $L':=L_1^\perp\neq \{0\}$ is an ideal.
Now, if $\{0\}\neq J\subseteq L_1$ is an ideal (inside $L_1$), then $[L,J]
\subseteq [L_1,J]+[L_1^\perp,J]=[L_1,J]\subseteq J$ and so $J$ actually is 
an ideal in all of $L$. Hence, by the minimality of $\dim L_1$, we must
have $J=L_1$. Thus, $L_1$ does not have any proper ideals. Furthermore,
$L_1$ is non-abelian because $L$ is semisimple by Lemma~\ref{defkill0a}(c).
Thus, $L_1$ is simple as a Lie algebra.

By Exercise~\ref{xcakill}, we have $\kappa_L|_{L'\times L'}=\kappa_{L'}$. 
We claim that $\kappa_{L'}$ is also non-degenerate. Indeed, let $x\in L'$ 
and assume that $\kappa_{L'} (x,y)=0$ for all $y' \in L'$. Let $z\in L$. 
Then $z=x'+y'$ where $x'\in I$ and $y'\in L'$. Hence, $\kappa_L(x,z)=
\kappa_L(x,x')+\kappa_L(x,y')=\kappa_{L'}(x,y')=0$ and so $x=0$, since
 $\kappa_L$ is non-degenerate. By induction, we can write $L'=L_2\oplus 
\ldots \oplus L_r$ where each $L_i$ is an ideal in $L'$ which is simple 
as a Lie algebra. Finally, for $i\geq 2$, we have $[L_1,L_i] \subseteq 
[L_1,L']\subseteq L_1\cap L'=\{0\}$ and so $[L,L_i]\subseteq L_i$. Hence, 
each $L_i$ actually is an ideal in~$L$. This also implies that, for 
if $i\neq j$, we have $[L_i,L_j]\subseteq L_i\cap L_j=\{0\}$.
\end{proof}

%
%

\noindent This now sets the programme that we will have to pursue:

1) Obtain some idea of how solvable Lie algebras look like. 

2) Study in more detail semisimple and simple Lie algebras.  

\noindent In order to attack 1) and~2), the representation theory of Lie 
algebras will play a crucial role. This is introduced in the following 
section.

\section{Representations of Lie algebras} \label{sec04}

A fundamental tool in the theory of groups is the study of actions of
groups on sets. There is an analogous notion for the action of Lie algebras
on vector spaces, taking into account the Lie bracket. Throughout,
let $L$ be a Lie algebra over our given field~$k$.

\begin{defn} \label{defmod} Let $V$ be a vector space (also over $k$). 
Then $V$ is called an \nms{$L$-module}{L-module} if there is a bilinear map
\[ L\times V \rightarrow V, \qquad (x,v) \mapsto x.v\]
such that $[x,y].v=x.(y.v)-y.(x.v)$ for all $x,y \in L$ and $v\in V$. In this
case, we obtain for each $x\in L$ a linear map 
\[ \rho_x\colon V\rightarrow V,\qquad v\mapsto x.v,\]
and one immediately checks that $\rho\colon L\rightarrow \gl(V)$,
$x\mapsto \rho_x$, is a Lie algebra homomorphism, that is, 
$\rho_{[x,y]}=[\rho_x,\rho_y]=\rho_x{\circ} \rho_y-\rho_y{\circ} \rho_x$
for all $x,y\in L$.
This homomorphism $\rho$ will also be called the corresponding
\nm{representation} of $L$ on $V$. If $\dim V <\infty$ and $B=\{v_i\mid i
\in I\}$ is a basis of $V$, then we obtain a \nm{matrix representation} 
\[\rho_B\colon L \rightarrow \gl_I(k), \qquad x \mapsto M_B(\rho(x)),\]
where $M_B(\rho(x))$ denotes the matrix of $\rho(x)$ with respect to $B$.
Thus, we have $M_B(\rho(x))=(a_{ij})_{i,j\in I}$ where $x.v_j=
\sum_{i\in I} a_{ij} v_i$ for all~$j$. 
\end{defn}

If $V$ is an $L$-module with $\dim V<\infty$, then all the known 
techniques from Linear Algebra can be applied to the study of the maps
$\rho_x\colon V\rightarrow V$: these have a trace, a determinant, 
eigenvalues and~so~on. 

\begin{rem} \label{defmod1} Let $\rho\colon L\rightarrow \gl(V)$ be a Lie 
algebra homomorphism, where $V$ is a vector space over $k$; then $\rho$ is 
called a \textit{representation} of~$L$. One 
immediately checks that $V$ is an $L$-module via
\[ L\times V\rightarrow V,\qquad (x,v)\mapsto \rho(x)(v);\]
furthermore, $\rho$ is the homomorphism associated with this $L$-module 
structure on $V$ as in Definition~\ref{defmod}. Thus, speaking about
``$L$-modules'' or ``representations of $L$'' are just two equivalent 
ways of expressing the same mathematical concept.
\end{rem}

\begin{exmp} \label{defmod2} 
(a) If $V$ is a vector space and $L$ is a Lie subalgebra of $\gl(V)$, 
then the inclusion $L\hookrightarrow \gl(V)$ is a representation. So 
$V$ is an $L$-module in a canonical way, where $\rho_x\colon V\rightarrow V$
is given by $v\mapsto x(v)$, that is, we have $\rho_x=x$ for all $x\in L$.

(b) The map $\ad_L\colon L\rightarrow \gl(L)$ in Example~\ref{defadj} is 
a Lie algebra homomorphism, called the \nm{adjoint representation} of $L$. 
So $L$ itself is an $L$-module via this map. 
\end{exmp}

\begin{xca} \label{dualspace} Let $V$ be an $L$-module and $V^*=\Hom(V,k)$
be the dual vector space. Show that $V^*$ is an $L$-module via
$L\times V^*\rightarrow V^*$, $(x,\mu)\mapsto \mu_x$,
where $\mu_x\in V^*$ is defined by $\mu_x(v)=-\mu(x.v)$ for $v\in V$.
Assume now that $n:=\dim V<\infty$ and let $B$ be a basis of $V$.
Let $B^*$ be the dual basis of $V^*$. Let $x\in L$. Then describe 
the matrix (with respect to~$B^*$) of the linear map $V^*\rightarrow V^*$, 
$\mu\mapsto \mu_x$, in terms of the matrix (with respect to~$B$) of the
linear map $V\rightarrow V$, $v \mapsto x.v$. 
\end{xca}

\begin{exmp} \label{expsemi} Let $V$ be an $L$-module and $\rho\colon 
L\rightarrow \gl(V)$ be the corresponding representation. Now $V$ is an 
abelian Lie algebra with Lie bracket $[v,v']=0$ for all $v,v'\in V$. Hence, 
we have $\Der(V)=\gl(V)$ and we can form the \nm{semidirect product} 
$L\ltimes_\rho V$, see Exercise~\ref{xcasemidir}. We have $[(x,0),(0,v)]=
(0,x.v)$ for all $x\in L$ and $v\in V$.
\end{exmp}

\begin{defn} \label{submod} Let $V$ be an $L$-module; for $x\in L$, we
denote by $\rho_x\colon V\rightarrow V$ the linear map defined by~$x$.
Let $U\subseteq V$ be a subspace. We say that $U$ is an $L$-\nm{submodule} 
(or an $L$-\nm{invariant subspace}) if $\rho_x(U) \subseteq U$ for all
$x\in L$. If $V\neq \{0\}$ and $\{0\}, V$ are the only $L$-invariant
subspaces of $V$, then $V$ is called an \nm{irreducible module}.

Assume now that $U$ is an $L$-invariant subspace. Then $U$ itself is an 
$L$-module, via the restriction of $L\times V\rightarrow V$ to a bilinear 
map $L\times U\rightarrow U$. Furthermore, $V/U$ is an $L$-module via 
\[ L\times V/U\rightarrow V/U,\qquad (x,v+U)\mapsto x.v+U.\]
(One checks as usual that this is well-defined and bilinear.)
Finally, assume that $n=\dim V<\infty$ and let $d=\dim U$. Let $B=\{v_1,
\ldots, v_n\}$ be a basis of $V$ such that $B'=\{v_1,\ldots,v_d\}$ is a 
basis of $U$. Since $x.v_i\in U$ for $1\leq i\leq d$, the corresponding 
matrix representation has the following block triangular shape:
\[ \renewcommand{\arraystretch}{1.2} \rho_B(x)=\left(\begin{array}{c|c} 
\rho'(x) & * \\ \hline 0 & \rho''(x) \end{array}\right) \qquad 
\mbox{for all $x\in L$},\]
where $\rho'\colon L\rightarrow \gl_d(k)$ is the matrix representation
corresponding to $U$ (with respect to the basis $B'$ of $U$) and
$\rho''\colon L\rightarrow \gl_{n-d}(k)$ is the matrix representation
corresponding to $V/U$ (with respect to the basis $B''=\{v_{d+1}+U,
\ldots, v_n+U\}$ of $V/U$).
\end{defn}

\begin{exmp} \label{dirprodmod} Let $V_1$ and $V_2$ be $L$-modules.
Then the vector space direct product $V:=V_1\times V_2=\{(v_1,v_2)\mid 
v_1\in V_1, v_2\in V_2\}$ also is an $L$-module, with operation defined
by 
\[ L\times (V_1\times V_2)\rightarrow V_1\times V_2,\qquad \bigl(x,
(v_1,v_2)\bigr) \mapsto (x.v_1,x.v_2).\]
Now we may identity $V_1$ with the subspace $\{(v_1,0)\mid v_1\in V_1\} 
\subseteq V$ and~$V_2$ with the subspace $\{(0,v_2) \mid v_2\in V_2\}
\subseteq V$. Then $V=V_1\oplus V_2$ becomes the direct 
sum of $V_1$ and $V_2$. If $B_1$ is a basis of $V_1$ and $B_2$ is a 
basis of $V_2$, then $B:=B_1\cup B_2$ is a basis of $V$ and the 
corresponding matrix representation has the following block diagonal shape:
\[ \renewcommand{\arraystretch}{1.2} \rho_B(x)=\left(\begin{array}{c|c} 
\rho_{B_1}(x) & 0 \\ \hline 0 & \rho_{B_2}(x) \end{array}\right) \qquad 
\mbox{for all $x\in L$}.\]
We have the following extension of the above discussion. For $i=1,2$ let
$L_i$ be a Lie algebra and $V_i$ be an $L_i$-module; let $\rho_i\colon 
L_i \rightarrow \gl(V_i)$ be the corresponding representation. Now 
$L:=L_1\times L_2$ also is a Lie algebra (see Example~\ref{dirprodalg}). 
Since, for $i=1,2$, the two projections $\pi_i\colon L\rightarrow L_i$, 
$(x_1,x_2)\mapsto x_i$ are Lie algebra homomorphisms, we also obtain 
representations $\tilde{\rho}_i:=\rho_i \circ \pi_i \colon L\rightarrow 
\gl(V_i)$. Thus, each $V_i$ can be regarded as $L$-module, with operation 
as follows:
\[L\times V_i \rightarrow V_i,\quad ((x_1,x_2),v_i)\mapsto x_i.v_i
\quad\qquad (i=1,2).\]
Consequently, $V_1\oplus V_2$ also is a module for $L=L_1\times L_2$, 
with operation given by $L\times (V_1\oplus V_2)\rightarrow V_1\oplus V_2$, 
\[ ((x_1,x_2),v_1+v_2)\mapsto x_1.v_1+x_2.v_2.\]
\end{exmp}

\begin{exmp} \label{liedirsum}
Assume that we have a direct sum decomposition $L=L_1\oplus \ldots \oplus 
L_r$ ($r\geq 1$) as in Proposition~\ref{killing8b}(b). We regard $L$ as an 
$L$-module via the adjoint representation (see Example~\ref{defmod2}).
Since each $L_i$ is an ideal in $L$, it is clear that $L_i$ is an 
$L$-submodule. Furthermore, since $L_i$ is a simple Lie algebra, $L_i$
is irreducible as a submodule of $L$. Thus, $L=L_1\oplus \ldots \oplus L_r$
is a decomposition of the $L$-module $L$ as a direct sum of irreducible
$L$-submodules.
\end{exmp}

\begin{prop} \label{submod1} Let $V\neq \{0\}$ be an $L$-module with
$\dim V<\infty$. There is a sequence of $L$-submodules
$\{0\}=V_0\subsetneqq V_1\subsetneqq V_2 \subsetneqq \ldots \subsetneqq
V_r=V$ such that $V_i/V_{i-1}$ is irreducible for $1\leq i\leq r$. Let 
$n_i=\dim (V_i/V_{i-1})$ for all~$i$. Then there is a basis $B$ of $V$ 
such that the matrices of the representation $\rho\colon L\rightarrow 
\gl(V)$ have the following shape
\[\renewcommand{\arraystretch}{0.8}
\rho_B(x)=\left(\begin{array}{c@{\hspace{3pt}}c@{\hspace{3pt}}
c@{\hspace{3pt}}c} \rho_1(x) & * & \ldots & *\\
0 & \rho_2(x) & \ddots & \vdots \\ \vdots & \ddots  &\ddots & * \\ 
0 & \ldots & 0 &\rho_r(x) \end{array}\right) \qquad \mbox{for all $x\in L$},\]
where $\rho_i\colon L\rightarrow \gl_{n_i}(k)$ is an irreducible
matrix representation corresponding to the $L$-module $V_i/V_{i-1}$.
\end{prop}

\begin{proof} Let $U\subsetneqq V$ be an $L$-submodule with $\dim U$ 
as large as possible. If $W\subseteq V/U$ is a submodule, then one
easily checks that $\{v\in V \mid v+U\in W\}\subseteq V$ is a submodule 
containing $U$, so $W=\{0\}$ or $W=V/U$. Hence, $V/U$ is irreducible and 
we continue with $U$. 
\end{proof}

\begin{exmp} \label{slnss0a} If $V$ is an $L$-module with $\dim V=1$, 
then $V$ is obviously irreducible. Let $V=\langle v\rangle_k$
where $0\neq v\in V$. Then, for all $x\in L$, we have $x.v=\varphi(x)v$
where $\varphi(x)\in k$. It follows that $\varphi\colon L\rightarrow k$ is 
linear. Furthermore, $\varphi([x,y])v=[x,y].v=x.(y.v)-y.(x.v)=\varphi(y)x.v
-\varphi(x)y.v=0$ and so $\varphi([x,y])=0$ for all $x,y\in L$. In 
particular, if $L=[L,L]$, then $L$ acts as zero on~$V$.
\end{exmp}

\begin{exmp} \label{xcahom} Let $V$ and $W$ be $L$-modules; let $\rho
\colon L\rightarrow \gl(V)$ and $\sigma\colon L\rightarrow \gl(W)$
be the corresponding representations. A linear map $\varphi\colon V
\rightarrow W$ is called an \nms{$L$-module homomorphism}{L-module
homomorphism} if $\varphi$ commutes with the actions of $L$ on $V$
and $W$, that is, we have 
\[\varphi(x.v)=x.\varphi(v) \qquad \mbox{for all $x\in L$ and $v\in V$}\]
or, equivalently, $\varphi\circ \rho_x=\sigma_x \circ \varphi$ for all 
$x\in L$. In this case, one easily sees that the kernel $\ker(\varphi)$ 
is an $L$-submodule of $V$ and the image $\varphi(V)$ is an $L$-submodule 
of $W$. Furthermore, if $W'\subseteq W$ is an $L$-submodule, then the
preimage $\varphi^{-1}(W')$ is an $L$-submodule of~$V$.
\end{exmp}

\begin{xca} \label{liethmc}
Let $k$ be a field of characteristic~$2$ and $L$ be the Lie algebra over
$k$ with basis $\{x,y\}$ such that $[x,y]=y$ (see 
Exercise~\ref{expsolv1}). Show that the linear map defined by
\[ \renewcommand{\arraystretch}{0.8}\rho\colon L\rightarrow \gl_2(k),
\quad x\mapsto\left(\begin{array}{c@{\hspace{5pt}}c} 0 & 0 \\ 0 & 1 
\end{array}\right), \quad y\mapsto\left(\begin{array}{c@{\hspace{5pt}}c} 
0 & 1 \\ 1 & 0 \end{array}\right), \]
is a Lie algebra homomorphism and so $V=k^2$ is an $L$-module.
Show that $V$ is an irreducible $L$-module. Check that $L$ is solvable.
\end{xca}

\begin{xca} \label{xcarepwitt} Let $L=\mbox{Der}(k[T,T^{-1}])$ be the 
\nm{Witt algebra} in Example~\ref{wittalg}, with basis $\{L_m\mid m 
\in \Z\}$. Let $V$ be a vector space with a basis $\{v_i\mid i\in
\Z\}$. Let $a,b\in k$ be fixed. For $m\in\Z$ define a 
linear map $\rho_m\colon V \rightarrow V$ by
\[ \rho_m(v_i):=\bigl(i+a+b(m+1)\bigr)v_{m+i} \qquad 
\mbox{for all $i\in \Z$}.\]
Show that $V$ is an $L$-module, where $L_m.v_i=\rho_m(v_i)$ for
all $i,m\in\Z$. Determine the conditions on $a$ and $b$ under
which $V$ is irreducible. 
\end{xca}

\begin{xca} \label{xcagenerator1} Let $X\subseteq L$ be a subset such that 
$L=\langle X\rangle_{\text{alg}}$. \\
{\rm (a)} Let $V$ be an $L$-module and $U \subseteq V$ be a subspace. 
Assume that $x.u\in U$ for all $x\in X$ and $u\in U$. Then show that 
$U$ is an $L$-submodule of~$V$.\\
{\rm (b)} Let $V,W$ be $L$-modules and $\varphi \colon V\rightarrow W$
be a linear map such that $\varphi(x.v)=x.\varphi(v)$ for all $x\in X$
and $v \in V$. Then show that $\varphi$ is an $L$-module homomorphism.
\end{xca}

Up to this point, $k$ could be any field (of any characteristic). 
Stronger results will hold if $k$ is algebraically closed.

\begin{lem}[Schur's Lemma] \label{schurl} \nmi{}{Schur's Lemma} Assume 
that $k$ is algebraically closed. Let $V$ be an irreducible $L$-module, 
$\dim V<\infty$. If $\varphi\in \End(V)$ is such that $\varphi\circ\rho_x=
\rho_x\circ \varphi$ for all $x\in L$, then $\varphi=c\,\id_V$ where
$c\in k$.
\end{lem}

\begin{proof} By Example~\ref{xcahom}, $\ker(\varphi)$ is an $L$-submodule 
of $V$. Since $V$ is irreducible, $\varphi=\underline{0}$ or 
$\ker(\varphi)=\{0\}$. If $\varphi=\underline{0}$, then the desired
assertion holds with $c=0$. Now assume that $\varphi\neq \underline{0}$. 
Then $\ker(\varphi)=\{0\}$ and $\varphi$ is bijective. Since $k$ is 
algebraically closed, there is an eigenvalue $c\in k$ for 
$\varphi$. Setting $\psi:=\varphi-c\,\id_V\in \End(V)$, we also have 
$\psi(x.v)=x.(\psi(v))$ for all $x\in L$ and $v\in V$. Hence, the previous 
argument shows that either $\psi=\underline{0}$ or $\psi$ is bijective. But 
an eigenvector of $\varphi$ for the eigenvalue~$c$ lies in $\ker(\psi)$ and 
so $\psi=\underline{0}$.
\end{proof}

\begin{prop} \label{abelsimple} Assume that $k$ is algebraically closed and 
$L$ is abelian. Let $V\neq\{0\}$ be an $L$-module with $\dim V<\infty$. 
Then there exists a basis $B$ of $V$ such that, for any $x\in L$, the matrix
of the linear map $\rho_x\colon V\rightarrow V$, $v\mapsto x.v$, with
respect to $B$ has the following shape:
\[\renewcommand{\arraystretch}{0.8}
M_B(\rho_x)=\left(\begin{array}{c@{\hspace{3pt}}c@{\hspace{3pt}}
c@{\hspace{3pt}}c} \lambda_1(x) & * & 
\ldots & *\\ 0 & \lambda_2(x) & \ddots & \vdots \\ \vdots & \ddots 
&\ddots & * \\ 0 & \ldots & 0 & \lambda_n(x)\end{array}\right)\qquad
(n=\dim V),\]
where $\lambda_i\colon L\rightarrow k$ are linear maps for $1\leq i\leq n$.
In particular, if $V$ is irreducible, then $\dim V=1$.
\end{prop}

\begin{proof} Assume first that $V$ is irreducible. We show that $\dim V=1$.
Let $x\in L$ be fixed and $\varphi:=\rho_x$. Since $L$ is abelian, we have 
$0=\rho_0=\rho_{[x,y]}=\varphi \circ \rho_y-\rho_y\circ \varphi$ for all
$y\in L$. By Schur's Lemma, $\varphi=\lambda(x)\,\id_V$ where $\lambda(x)
\in k$. Hence, if $0\neq v \in V$, then $x.v=\lambda(x)v$ for all $x\in L$ 
and so $\langle v\rangle_k \subseteq V$ is an $L$-submodule. Clearly, 
$\lambda \colon L\rightarrow k$ is linear. Since $V$ is irreducible, 
$V=\langle v \rangle_k$ and so $\dim V=1$. The general case follows from 
Proposition~\ref{submod1}.
\end{proof}

\begin{exmp} \label{liethm1} Assume that $k$ is algebraically closed. Let
$V$ be a vector space over $k$ with $\dim V<\infty$. Let $\fX\subseteq 
\End(V)$ be a subset such that $\varphi \circ \psi=\psi\circ \varphi$ for 
all $\varphi,\psi\in\fX$. Then there exists a basis $B$ of $V$ such that 
the matrix of any $\varphi \in \fX$ with respect to $B$ is upper triangular. 
Indeed, just note that $L:=\langle \fX\rangle_k \subseteq \gl(V)$ is an 
abelian Lie subalgebra and $V$ is an $L$-module; then apply 
Proposition~\ref{abelsimple}. (Of course, one could also prove this more 
directly.)
\end{exmp}


\begin{xca} \label{xca000} This exercise establishes an elementary result 
from Linear Algebra that will be useful at several places. Let $k$ be an 
infinite field and $V$ be a $k$-vector space with $1\leq \dim V<\infty$. 
Let $V^*:=\Hom(V,k)$ be the dual space.

(a) Show that, if $X \subseteq V$ is a finite subset such that 
$0\not\in X$, then there exists $\mu_0\in V^*$ such that $\mu_0(x)\neq 0$ 
for all $x \in X$. 

(b) Similarly, if $\Lambda\subseteq V^*$ is a finite subset such that 
$\underline{0} \not\in \Lambda$ (where $\underline{0}\colon V\rightarrow k$ 
denotes the linear map with value $0$ for all $v\in V$), then there exists
$v_0\in V$ such that $f(v_0)\neq 0$ for all $f\in \Lambda$. 

Show that the above statements remain true if we only assume that 
$|k|\geq |X|$ in (a), or $|k|\geq |\Lambda|$ in (b). In any case,
deduce that $V$ is not the union of finitely many proper subspaces.
\end{xca}

\begin{exmp} \label{tenslie} Let $V,W$ be $L$-modules. Then the 
tensor product $V\otimes W$ also is an $L$-module, with operation 
given by the map
\[ L\times (V\otimes W)\rightarrow V\otimes W, \quad (x,v\otimes w)
\mapsto (x.v)\otimes w+v\otimes (x.w).\]
The existence of this map is shown as follows. For a fixed $x\in L$, we
have a bilinear map
\[\varphi_x\colon V\times W \rightarrow V\otimes W, \qquad (v,w)\mapsto 
(x.v)\otimes w+v\otimes (x.w).\]
So, by the defining property of $V\otimes W$, there is a unique
linear map $\tilde{\varphi}_x\in \End(V \otimes W)$ such that 
$\tilde{\varphi}_x(v\otimes w)=\varphi_x(v,w)$ for all $v\in V$, 
$w\in W$. Next, we show that $L\times (V\otimes W) \rightarrow V\otimes 
W$ is bilinear. For this purpose, let $x,y\in L$ 
and $c,c'\in k$. We can apply the above discussion to the element $z:=
cx+c'y\in L$ and obtain a unique $\tilde{\varphi}_z \in \End(V \otimes 
W)$ such that $\tilde{\varphi}_z(v\otimes w)=\varphi_z(v,w)$ for all
$v\in V$, $w\in W$. Then it is straightforward to check that 
$(c\tilde{\varphi}_x+c'\tilde{\varphi}_y)(v\otimes w)= \tilde{\varphi}_z
(v\otimes w)$ for all $v\in V$, $w\in W$; hence, we must have 
$\tilde{\varphi}_z=c\tilde{\varphi}_x+c' \tilde{\varphi}_y$.

Thus, we have a bilinear map $L\times (V\otimes W)\rightarrow V\otimes W$
as claimed. It remains to show the condition for the action of $[x,y]$
where $x,y\in L$: 
\[[x,y].(v\otimes w)=x.(y.(v\otimes w))-y.(x.(v\otimes w))\quad
\mbox{for $v\in V$, $w\in W$}.\]
We leave this as an exercise to the readrer.
\end{exmp}

\begin{rem} \label{tenslie1} We have the following extension of the
above discussion. For $i=1,2$ let $L_i$ be a Lie algebra and $V_i$
be an $L_i$-module; let $\rho_i\colon L_i \rightarrow \gl(V_i)$ be the
corresponding representation. Now $L:=L_1\times L_2$ also is a Lie algebra
(see Example~\ref{dirprodalg}). Since, for $i=1,2$, the two projections
$\pi_i\colon L\rightarrow L_i$, $(x_1,x_2)\mapsto x_i$ are Lie algebra 
homomorphisms, we also obtain representations $\tilde{\rho}_i:=\rho_i
\circ \pi_i \colon L\rightarrow \gl(V_i)$. Thus, each $V_i$ can be 
regarded as $L$-module, with operation as follows:
\[L\times V_i \rightarrow V_i,\quad ((x_1,x_2),v_i)\mapsto x_i.v_i
\quad\qquad (i=1,2).\]
Consequently, by Example~\ref{tenslie}, $V_1\otimes V_2$ also is a module 
for $L=L_1\times L_2$, with operation given by $(L_1\times L_2)\times 
(V_1\otimes V_2) \rightarrow V_1\otimes V_2$,
\[ ((x_1,x_2),v_1\otimes v_2)\mapsto (x_1.v_1)\otimes v_2+ v_1\otimes 
(x_2.v_2).\]
\end{rem}

\begin{xca}[Difficult!] \label{tenslie2} In the set-up of 
Remark~\ref{tenslie1}, assume that $V_i$ is a finite-dimensional and
irreducible $L_i$-module, for $i=1,2$. Then show that $V_1\otimes V_2$ 
is an irreducible $(L_1\times L_2)$-module. \\
(This statement will not be used in this text.)\\
{\footnotesize [{\it Hint}. See Steinberg \cite[Lemma~68 and Cor., 
p.~117]{St}. If $L$ is of ``Cartan--Killing type'' as in Chapter~2, then 
we will see this much later in an appendix.]}
\end{xca}

\section{Lie's Theorem} \label{sec04a}

The content of Lie's Theorem is that Proposition~\ref{abelsimple} (which 
was concerned with representations of abelian Lie algebras) remains true 
for the more general class of solvable Lie algebras, assuming that $k$ is 
not only algebraically closed but also has characteristic~$0$. 
(Exercice~\ref{liethmc} shows that this will definitely not work in 
positive characteristic.) So, in order to use the full power of the 
techniques developed so far, we will assume that $k=\C$.

Let $L$ be a Lie algebra over $k=\C$. If $V$ is an $L$-module, then we
denote as usual by $\rho_x\colon V \rightarrow V$ the linear map defined 
by~$x\in L$. Our approach to Lie's Theorem is based on the following 
technical result.

\begin{lem} \label{newinvlem} Let $V$ be an irreducible $L$-module (over
$k=\C$), with $\dim V<\infty$. Let $H\subseteq L$ be an abelian ideal in 
$L$ such that $\trc(\rho_x)=0$ for all $x\in H$. Then $\rho_x=0$ for all 
$x\in H$. 
\end{lem}

\begin{proof} Let $x\in H$ and consider the linear map $\rho_x\colon V
\rightarrow V$. Since $k=\C$, this map has eigenvalues. Let $c\in \C$ be 
an eigenvalue and consider the \nm{generalized eigenspace} 
\[V_c(\rho_x):=\{v\in V\mid (\rho_x-c\,\id_V)^l(v)=0\mbox{ for some
$l\geq 1$}\}\neq\{0\}.\]
We claim that $V_c(\rho_x)\subseteq V$ is an $L$-submodule. To see this,
let $v\in V_c(\rho_x)$ and $y \in L$. We must show that $y.v=\rho_y(v) \in 
V_c(\rho_x)$. Let $l\geq 1$ be such that $(\rho_x-c\,\id_V)^l(v)=0$.
We apply the \nm{generalized binomial formula} (Lemma~\ref{superbin}) to 
the associative algebra $A:=\End(V)$, the elements $\rho_x,\rho_y\in A$
and the scalars $a:=-c$, $b:=0$. This yields
\begin{equation*}
(\rho_x-c\,\id_V)^{l+1} \circ \rho_y=\sum_{i=0}^{l+1}
\binom{l{+}1}{i} \psi_i\circ (\rho_x-c\,\id_V)^{l+1-i},\tag{$*$}
\end{equation*}
where $\psi_i:=\ad_A(\rho_x)^i(\rho_y)\in A$ for $i\geq 0$. We claim 
that $\psi_i= \underline{0}$ for $i\geq 2$. Indeed, since $\ad_A(\rho_x)
(\rho_z)= \rho_x\circ \rho_z-\rho_z \circ \rho_x=\rho_{[x,z]}$ for any
$z\in L$, we obtain:
\[ \psi_i=\ad_A(\rho_x)^{i-2}\bigl(\ad_A(\rho_x)^2(\rho_y)\bigr)
=\ad_A(\rho_x)^{i-2}(\rho_{x,[x,y]]}).\]
But $[x,y] \in H$ because $H$ is an ideal, and $[x,[x,y]]=0$ because $H$
is abelian. So $\psi_i=\underline{0}$ for $i\geq 2$, as claimed. Now apply
both sides of ($*$) to~$v$. If $i=0,1$, then $l+1-i\geq l$ and so
$(\rho_x-c\,\id_V)^{l+1-i}(v)=0$. On the other hand, $\psi_i=\underline{0}$
for $i\geq 2$. Hence, the right hand side of~($*$), applied to~$v$,
equals~$0$. Consequently, we also have 
\[(\rho_x-c\,\id_V)^{l+1}(y.v)= \bigl((\rho_x-c\, \id_V)^{l+1}\circ
\rho_y\bigr)(v)=0\]
and so $y.v\in V_c(\rho_x)$, as desired.

Now, since $V$ is irreducible and $V_c(\rho_x)\neq \{0\}$, we conclude that 
$V=V_c(\rho_x)$. Let $\psi_x:=\rho_x-c\,\id_V$. Then, for $v\in V$, there 
exists some $l \geq 1$ with $\psi_x^l(v)=0$. So Exercise~\ref{xcanilp} 
shows that $\psi_x$ is nilpotent and $\trc(\psi_x)=0$. But then $\trc
(\rho_x)=\trc(\psi_x+c\,\id_V)=(\dim V)c$. So our assumptions on the 
characteristic of $k$ and on $\trc(\rho_x)$ imply that $c=0$. Thus, $c=0$ 
is the only eigenvalue of $\rho_x$, for any $x\in H$. 

Finally, regarding $V$ as an $H$-module (by restricting the action of 
$L$ on $V$ to $H$), we can apply Proposition~\ref{abelsimple}. This yields
a basis $B$ of $V$ such that, for any $x\in H$, the matrix of $\rho_x$ with 
respect to $B$ is upper triangular; by the above discussion, the entries 
along the diagonal are all~$0$. Let $v_1$ be the first vector 
in $B$. Then $x.v_1=\rho_x(v_1)= 0$ for all $x\in H$.
 Hence, the subspace 
\[ U:=\{v\in V\mid x.v=0\mbox{ for all $x\in H$}\}\]
is non-zero. Now we claim that $U$ is an $L$-submodule. Let $v\in V$ and
$y\in L$. Then, for $x\in H$, we have $x.(y.v)=[x,y].v+y.(x.v)=
[x,y].v=0$ since $v\in U$ and $[x,y]\in H$. Since $V$ is irreducible,
we conclude that $U=V$ and so $\rho_x=0$ for all $x\in H$. 
\end{proof}

\begin{prop}[Semisimplicity criterion] \label{critss} Let $k=\C$ and $V$ 
be a vector space with $\dim V<\infty$. Let $L\subseteq \slm(V)$ be a Lie 
subalgebra such that $V$ is an irreducible $L$-module. Then $L$ is 
semisimple.
\end{prop}

\begin{proof} If $\rad(L)\neq \{0\}$ then, by Lemma~\ref{liesolv}, there 
exists a non-zero abelian ideal $H\subseteq L$ such that $H\subseteq 
\rad(L)$. Since $L\subseteq \slm(V)$, Lemma~\ref{newinvlem} implies that 
$x=\rho_x=0$ for all $x\in H$, contradiction.
\end{proof}

\begin{exmp} \label{slnss0} Let $k=\C$ and $V$ be a vector space with 
$\dim V<\infty$. Clearly (!), $V$ is an irreducible $\gl(V)$-module. Next 
note that $\gl(V)=\slm(V)\oplus\C\,\id_V$. Hence, if $U\subseteq V$ is an 
$\slm(V)$-invariant subspace, then $U$ will also be $\gl(V)$-invariant. 
Consequently, $V$ is an irreducible $\slm(V)$-module. Hence, 
Proposition~\ref{critss} shows that $\slm(V)$ is semisimple. 

Note that, if $\mbox{char}(k)=p>0$ and $L=\slm_p(k)$,
then $Z:=\{aI_p \mid a \in k\}$ is an abelian ideal in $L$ and so 
$L$ is not semisimple in this case. 
\end{exmp}

\begin{thm}[Lie's Theorem] \label{liethm} Let $k=\C$. Let $L$ be solvable 
and $V\neq\{0\}$ be an $L$-module, where $\dim L<\infty$ and $\dim V
<\infty$. Then the conclusions in Proposition~\ref{abelsimple} still hold,
that is, there exists a basis $B$ of $V$ such that, for any $x\in L$, the 
matrix of the linear map $\rho_x\colon V\rightarrow V$, $v\mapsto x.v$, with
respect to $B$ has the following shape:
\[\renewcommand{\arraystretch}{0.8}
M_B(\rho_x)=\left(\begin{array}{c@{\hspace{3pt}}c@{\hspace{3pt}}
c@{\hspace{3pt}}c} \lambda_1(x) & * & 
\ldots & *\\ 0 & \lambda_2(x) & \ddots & \vdots \\ \vdots & \ddots 
&\ddots & * \\ 0 & \ldots & 0 & \lambda_n(x)\end{array}\right)\qquad
(n=\dim V),\]
where $\lambda_i\colon L\rightarrow k$ are linear maps such 
that $[L,L]\subseteq \ker(\lambda_i)$ for $1\leq i\leq n$. In particular, if
$V$ is irreducible, then $\dim V=1$.
\end{thm}

\begin{proof} First we show that, if $V$ is irreducible, then $\dim V=1$.
We use induction on $\dim L$. If $\dim L=0$, there is nothing to prove. 
Now assume that $\dim L>0$. If $L$ is abelian, then see 
Proposition~\ref{abelsimple}. Now assume that $[L,L]\neq \{0\}$. By 
Lemma~\ref{liesolv}, there exists a non-zero abelian ideal $H\subseteq 
L$ such that $H\subseteq [L,L]$. Let $x\in H$. Since $H\subseteq [L,L]$, 
we can write $x$ as a finite sum $x=\sum_i [y_i,z_i]$ where $y_i,z_i\in L$
for all~$i$. Consequently, we also have $\rho_x=\sum_i (\rho_{y_i}\circ
\rho_{z_i}- \rho_{z_i}\circ \rho_{y_i})$ and, hence, $\trc(\rho_x)=0$. By 
Lemma~\ref{newinvlem}, $\rho_x=0$ for all $x\in H$. Let $L_1:=L/H$. Then $V$ 
also is an $L_1$-module via 
\[L_1\times V\rightarrow V, \qquad (y+H,v) \mapsto y.v.\]
(This is well-defined since $x.v=0$ for $x\in H$, $v\in V$.) If $V'
\subseteq V$ is an $L_1$-invariant subspace, then $V'$ is also $L$-invariant.
Hence, $V$ is an irreducible $L_1$-module. By Lemma~\ref{defsolv1}(c), $L_1$
is solvable. So, by induction, $\dim V=1$.

The general case follows again from Proposition~\ref{submod1}. The fact
that $[L,L]\subseteq \ker(\lambda_i)$ for all~$i$ is seen as in
Example~\ref{slnss0a}. 
\end{proof}

\begin{lem} \label{liethmu} In the setting of Theorem~\ref{liethm},
the set of linear maps $\{\lambda_1,\ldots,\lambda_n\}$ does not depend
on the choice of the basis $B$ of $V$.\\
{\rm We shall call $P_L(V):=\{\lambda_1,\ldots,\lambda_n\}$ the set of 
\nmi{weights}{weight} of $L$ on~$V$}. 
\end{lem}

\begin{proof} Let $B'$ be another basis of $V$ such that, for any $x\in L$,
the matrix of $\rho_x\colon V\rightarrow V$ with respect to $B'$ has a 
triangular shape with $\lambda_1'(x)$, $\ldots$, $\lambda_n'(x)$ along
the diagonal, where $\lambda_i'\colon L\rightarrow k$ are linear maps such 
that $[L,L]\subseteq \ker(\lambda_i')$ for $1\leq i\leq n$. We must show
that $\{\lambda_1,\ldots,\lambda_n\}=\{\lambda_1',\ldots,\lambda_n'\}$. 
Assume, if possible, that there exists some $j$ such that $\lambda_j'
\neq \lambda_i$ for $1\leq i\leq n$. Let $\Lambda:=\{\lambda_i-\lambda_j'
\mid 1\leq i \leq n\}$. Then $\Lambda$ is a finite subset of $\Hom(L,\C)$
such that $\underline{0}\not\in \Lambda$. So, by Exercise~\ref{xca000}(b),
there exists some $x_0\in L$ such that $\lambda_j'(x_0)\neq \lambda_i(x_0)$
for $1\leq i\leq n$. But then $\lambda_j'(x_0)$ is an eigenvalue
of $M_{B'}(\rho_{x_0})$ that is not an eigenvalue of $M_B(\rho_{x_0})$,
contradiction since $M_B(\rho_{x_0})$ and $M_{B'}(\rho_{x_0})$ are similar 
matrices and, hence, they have the same characteristic polynomials. Thus, 
we have shown that $\{\lambda_1',\ldots,\lambda_n'\}\subseteq \{\lambda_1,
\ldots, \lambda_n\}$. Exchanging the roles of $B$, $B'$, we also have 
the reverse inclusion.
\end{proof}

\begin{exmp} \label{killing1} Let $k=\C$ and $L\neq \{0\}$ be solvable
with $\dim L<\infty$. Then, by Lemma~\ref{defkill0a}, the Killing form 
$\kappa_L\colon L\times L\rightarrow \C$ is degenerate. Actually, much 
more is true. Namely, applying Theorem~\ref{liethm} to the adjoint
representation $\ad_L\colon L\rightarrow\gl(L)$, there exists a basis $B$ 
of $L$ such that $M_B(\ad_L(x))$ is upper triangular for all $x\in L$. 
Consequently, if $x, y\in L$, then 
\[\ad_L([x,y])=\ad_L(x)\circ \ad_L(y)-\ad_L(y) \circ \ad_L(x)\]
is represented by a matrix which is upper triangular with $0$ on the
diagonal. Hence, we have 
\[\kappa_L([x,y],z)=0\qquad \mbox{for all $x,y,z\in L$},\]
which is one half of ``\nm{Cartan's First Criterion}''. The 
other half says that, if $\kappa_L([x,y],z)=0$ for all $x,y,z\in L$, then
$L$ is solvable. The proof requires much more work; see, for example, 
\cite[\S 9.2]{EW}.
\end{exmp}

\begin{xca} \label{xcaweightu} Let $k=\C$ and $L$ be solvable with 
$\dim L<\infty$. Let~$V$ be a finite-dimensional $L$-module and $U
\subseteq V$ be a non-zero, proper $L$-submodule. Show that $P_L(V)=
P_L(U) \cup P_L(V/U)$, where the set of weights of a module is defined 
by Lemma~\ref{liethmu}. In particular, if there is an $L$-submodule $U'
\subseteq V$ such that $V=U\oplus U'$, then $P_L(V)=P_L(U)\cup P_L(U')$.
\end{xca}

\begin{xca} \label{nonsplitsolv} Assume that $k\subseteq\C$. Show that 
\[\renewcommand{\arraystretch}{0.8}
L=\left\{\left(\begin{array}{@{\hspace{2pt}}rcc} 0 & t & x \\ -t & 0 & y 
\\ 0 & 0 & 0 \end{array}\right) \;\Big|\; t,x,y \in k\right\}\]
is a solvable Lie subalgebra of $\gl_3(k)$. Regard $V=k^3$ as an 
$L$-module via the inclusion $L\hookrightarrow \gl_3(k)$ 
(cf.\ Example~\ref{defmod2}). If $k=\C$, find a basis $B$ of $V$ 
such that the corresponding matrices of $L$ will be upper triangular. 
Does this also work with $k=\R$? 
\end{xca}

Finally, we develop some very basic aspects of the representation theory 
of $\slm_2(\C)$. As pointed out in \cite[\S 2.4]{MP}, this is of the 
utmost importance for the general theory of semisimple Lie algebras. (We 
shall see this in Section~\ref{sec1a2} and, again, much later in 
Chapter~\ref{chap4}.) For the remainder of this section, let 
$L=\slm_2(\C)$, with standard basis
\[\renewcommand{\arraystretch}{0.9}
 e=\left(\begin{array}{c@{\hspace{6pt}}c} 0 & 1 \\ 0 & 0\end{array}
\right),\quad h=\left(\begin{array}{c@{\hspace{5pt}}r} 1 & 0 \\ 
0 & -1\end{array}\right), \quad f=\left(\begin{array}{c@{\hspace{6pt}}r} 
0 & 0 \\ 1 & 0\end{array}\right),\]
where 
\[ [e,f]=h, \qquad [h,e]=2e, \qquad [h,f]=-2f;\]
see Exercise~\ref{expsolv2}. The following result is obtained by an easy
application of Lie's Theorem (but one can also prove it easily without 
reference to Lie's Theorem).

\begin{lem} \label{sl2higha} Let $V$ be an $\slm_2(\C)$-module with 
$\dim V<\infty$. Then there exists a non-zero vector $v^+\in V$ such that 
$e.v^+=0$ and $h.v^+=cv^+$ for some $c \in\C$. 
\end{lem}

\begin{proof} Let $S:=\langle h,e\rangle_{\C} \subseteq  \slm_2(\C)$. 
This is precisely the subalgebra of $\slm_2(\C)$ consisting of all 
upper triangular matrices with trace~$0$. Since $[h,e]=2e$, we have 
$[S,S]=\langle e\rangle_\C$ and so $S$ is solvable. By restricting the
action of $\slm_2(\C)$ on $V$ to $S$, we can regard $V$ as $S$-module.
So, by Theorem~\ref{liethm}, there exist a basis $B$ of $V$ and
$\lambda_1,\ldots,\lambda_n \in S^*$ (where $n=\dim V$) such that,
for any $x\in S$, the matrix of $\rho_x\colon V \rightarrow V$
is upper triangular with $\lambda_1(x),\ldots,\lambda_n(x)$ along
the diagonal; furthermore, $[S,S]\subseteq \ker(\lambda_i)$ for 
$1\leq i \leq n$. Let $v^+$ be the first vector in $B$. Then $\rho_x(v^+)
=\lambda_1(x)v^+$ for all $x\in S$. So $v^+$ has the required properties, 
where $c:=\lambda_1(h)\in \C$; we have $e.v^+=0$ since $e\in [S,S]$.
\end{proof}

\begin{rem} \label{sl2highb} Let $V\neq \{0\}$ be an $\slm_2(\C)$-module 
with $\dim V<\infty$. Let $v^+\in V$ be as in Lemma~\ref{sl2higha}; any 
such vector will be called a \nm{primitive vector} of~$V$. Then we define 
a sequence $(v_n)_{n\geq 0}$ in $V$ by 
\[ v_0:=v^+\qquad \mbox{and}\qquad \textstyle v_{n+1}:=
\frac{1}{n+1}f.v_n \quad \mbox{for all $n \geq 0$}.\]
Let $V':=\langle v_n\mid n \geq 0\rangle_\C\subseteq V$. 
We claim that the following relations hold for all $n\geq 0$ (where
we also set $v_{-1}:=0$):
\begin{equation*}
h.v_n=(c-2n)v_n \qquad \mbox{and}\qquad e.v_n=(c-n+1)v_{n-1}.\tag{a}
\end{equation*}
We use induction on $n$. If $n=0$, the formulae hold by definition. Now 
let $n\geq 0$. First note that $f.v_{n-1}=nv_n$. We compute:
\begin{align*}
(n+&1)e.v_{n+1} =e.(f.v_n)=[e,f].v_n+f.(e.v_n)=h.v_n+f.(e.v_n)\\
&=(c-2n)v_n+(c-n+1)f.v_{n-1}\qquad\mbox{(by induction)}\\&=(c-2n)v_n+
(c-n+1)nv_n=((n+1)c-n^2-n)v_n,
\end{align*}
and so $e.v_{n+1}=(c-n)v_n$, as required. Next, we compute:
\begin{align*}
(n+1)&h.{v_{n+1}}=h.(f.v_n)=[h,f].v_n+f.(h.v_n)\\&=-2f.v_n+(c-2n)f.v_n=
(c-2n-2)(n+1)v_{n+1},
\end{align*}
so (a) holds. Now, if $v_n\neq 0$ for all $n$, then $v_0,v_1,v_2,
\ldots$ are eigenvectors for $\rho_h\colon V\rightarrow V$ with distinct 
eigenvalues (see (a)) and so $v_0,v_1,v_2,\ldots$ are linearly independent, 
contradiction to $\dim V<\infty$. So there is some $n_0\geq 0$ such that 
$v_0,v_1,\ldots,v_{n_0}$ are linearly independent and $v_{n_0+1}=0$. But 
then, by the definition of the $v_n$, we have $v_n=0$ for all $n>n_0$ and 
so $V'=\langle v_0,v_1,\ldots,v_{n_0} \rangle_\C$. Furthermore, $0=e.0=
e.v_{n_0+1}=(c-n_0)v_{n_0}$ and so $c=n_0$. Thus, we obtain: 
\begin{equation*}
h.v^+=cv^+\qquad\mbox{where} \qquad c=\dim V'-1\in\Z_{\geq 0}.\tag{b}
\end{equation*}
So, the eigenvalue of our primitive vector $v^+$ has a very special form!

If $c\geq 1$, then the above formulae also yield an expression of $v^+=v_0$ 
in terms of $v_c=v_{n_0}$; indeed, by (a), we have $e.v_c=v_{c-1}$, 
$e.v_{c-1}= 2v_{c-2}$, $e.v_{c-2}=3v_{c-3}$ and so on. Thus, we obtain:
\begin{equation*}
\underbrace{e.e\ldots e}_{\text{$c$ times}}.v_c=
(1{\cdot} 2 {\cdot} 3 {\cdot} \ldots {\cdot} c)\, v^+.\tag{c}
\end{equation*}
We now state some useful consequences of the above discussion.
\end{rem}

\begin{cor} \label{sl2modc} In the setting of Remark~\ref{sl2highb},
assume that $V$ is irreducible. Write $\dim V=m+1$, $m \geq 0$. Then
$\rho_h$ is diagonalisable with eigenvalues $\{m-2i\mid 0\leq i \leq m\}$
(each with multiplicity~$1$). In particular, if $m\geq 1$, then $1$ or $2$ 
is an eigenvalue. The matrices of $\rho_e$ and $\rho_f$ with respect to 
the basis $\{v_0,v_1,\ldots,v_m\}$ of $V$ are given by
\[ \renewcommand{\arraystretch}{0.9} \rho_e: \left(
\begin{array}{@{\hspace{4pt}}c@{\hspace{4pt}}c@{\hspace{4pt}}
c@{\hspace{4pt}}c@{\hspace{4pt}}c} 0 & m &&&\\ &0 & m{-}1 &&\\ && \ddots & 
\ddots &  \\ & & & 0 & 1 \\ & & & & 0 \end{array}\right) \quad 
\mbox{and}\quad\; \rho_f: \left(\begin{array}{@{\hspace{5pt}}c@{\hspace{5pt}}
c@{\hspace{5pt}}c@{\hspace{5pt}}c@{\hspace{5pt}}c} 0 & &&&\\ 1 & 0 & &&\\ 
& 2 & 0 & & \\ & & \ddots & \ddots &  \\ & & & m & 0\end{array}\right)\]
(where all non-specified entries are $0$). 
\end{cor}

\begin{proof} Using the formulae in Remark~\ref{sl2highb} and an induction
on~$n$, one sees that $h.v_n\in V'$, $e.v_n\in V'$, $f.v_n\in V'$ for all
$n \geq 0$. Thus, $V'\subseteq V$ is an $\slm_2(\C)$-submodule. Since
$V'\neq \{0\}$ and $V$ is irreducible, we conclude that $V'=V$ and $m=c$.
By Remark~\ref{sl2highb}(a), we have $h.v_i=(m-2i)v_i$ for $i=0,1,\ldots,m$.
Hence, $\rho_h$ is diagonalisable, with eigenvalues as stated above. 
\end{proof}

%

\begin{prop} \label{sl2modd} Let $V$ be any finite-dimensional 
$\slm_2(\C)$-module, with $e,h,f$ as above. Then all the eigenvalues 
of $\rho_h\colon V \rightarrow V$ are integers and we have $\trc(\rho_h)=0$. 
Furthermore, if $n\in\Z$ is an eigenvalue of $\rho_h$, then so is $-n$
(with the same multiplicity as~$n$). If~$n$ is the largest eigenvalue, 
then $-n$ is the smallest eigenvalue.
\end{prop}

\begin{proof} Note that the desired statements can be read off the 
characteristic polynomial of $\rho_h\colon V \rightarrow V$. If $V$ is 
irreducible, then these hold by Corollary~\ref{sl2modc}. In general, let 
$\{0\}=V_0\subsetneqq V_1 \subsetneqq V_2 \subsetneqq \ldots \subsetneqq
V_r=V$ be a sequence of $L$-submodules as in Proposition~\ref{submod1}, 
such that $V_i/V_{i-1}$ is irreducible for $1\leq i \leq r$. It remains 
to note that the characteristic polynomial of $\rho_h\colon V 
\rightarrow V$ is the product of the characteristic polynomials of the 
actions of $h$ on $V_i/V_{i-1}$ for $1\leq i \leq r$.
\end{proof}

\section{The classical Lie algebras} \label{sec05}

Let $V$ be a vector space over $k$ and $\beta \colon V\times V
\rightarrow k$ be a bilinear map. Then we define $\gom(V,\beta)$ to be the 
set of all $\varphi\in \End(V)$ such that 
\[ \beta(\varphi(v),w)+ \beta(v,\varphi(w))=0 \qquad \mbox{for all 
$v,w\in V$}.\]
(The symbol ``$\gom$'' stands for ``general orthogonal''.) We check 
that $\gom(V,\beta)$ is a Lie subalgebra of $\gl(V)$.
Let $\varphi,\psi\in \gom(V,\beta)$. Then 
\begin{align*}
\beta([\varphi,&\psi](v),w)+\beta(v,[\varphi,\psi](w))\\&=
\beta(\varphi(\psi(v)),w)-\beta(\psi(\varphi(v)),w)\\
&\qquad +\beta(v,\varphi(\psi(w)))- \beta(v,\psi(\varphi(w))).
\end{align*}
Now $\beta(\varphi(\psi(v)),w)=-\beta(\psi(v),\varphi(w))$ since
$\varphi\in \gom(V,\beta)$; the three remaining terms can be re-written
analogously and, hence, the above sum equals zero. Thus, $\gom(V,\beta)$ 
is a Lie algebra, called a \nm{classical Lie algebra}. The further 
developments will show that these form an important class of semisimple 
Lie algebras (for certain $\beta$, over $k=\C$).

We assume throughout that $\beta$ is a \nm{reflexive bilinear form}, that 
is, for any $v,w\in V$, we have $\beta(v,w)=0\Leftrightarrow\beta(w,v)=0$. 
Thus, for any subset $X\subseteq V$, we can define 
\[X^\perp:=\{v\in V\mid \beta(v,x)=0 \mbox{ for all $x\in X$}\},\]
where it does not matter if we write ``$\beta(v,x)=0$'' or
``$\beta(x,v)=0$''. Note that $X^\perp$ is a subspace of $V$ (even if $X$
is not a subspace). We say that $\beta$ is a \nm{non-degenerate bilinear 
form} if $V^\perp=\{0\}$. 

We shall also assume throughout that $\mbox{char}(k) \neq 2$. (This 
avoids the consideration of some special cases that are not relevant to 
us here; for further details see \cite[\S 2]{Grove2}.) An elementary (but 
slightly tricky) argument shows that, since $\beta$ is reflexive, there 
exists a sign $\epsilon\in \{\pm 1\}$ such that $\beta(v,w)=\epsilon 
\beta(w,v)$ for all $v,w\in V$; see, for example, \cite[Prop.~2.7]{Grove2}. 
If $\epsilon=1$, then $\beta$ is called a \nm{symmetric bilinear form}; if 
$\epsilon=-1$, then $\beta$ is called an \nm{alternating bilinear form}. 

As in Example~\ref{defmod2}(a), the vector space $V$ is a $\gom(V,
\beta)$-module in a natural way. Again, this module turns out to be 
irreducible, if $\beta$ is non-degenerate.

\begin{prop} \label{classic1} Assume that $3\leq \dim V<\infty$. If $\beta$ 
is reflexive and non-degenerate, then $V$ is an irreducible 
$\gom(V,\beta)$-module.
\end{prop}

\begin{proof} First we describe a method for producing elements in
$\gom(V,\beta)$. For fixed $x,y\in V$ we define a linear map 
$\varphi_{x,y}\colon V\rightarrow V$ by $\varphi_{x,y}(v):=\beta(v,x)y-
\beta(y,v)x$ for all $v\in V$. We claim that $\varphi_{x,y}\in 
\gom(V,\beta)$. Indeed, for all $v,w\in V$, we have
\begin{align*}
\beta(\varphi_{x,y}(v),w)&+\beta(v,\varphi_{x,y}(w))\\
&=\bigl(\beta(v,x) \beta(y,w)-\beta(y,v)\beta(x,w)\bigr)\\
&\qquad\quad  + \bigl(\beta(w,x) \beta(v,y)-\beta(y,w)\beta(v,x)\bigr)\\
&=-\beta(y,v)\beta(x,w)+\beta(w,x) \beta(v,y),
\end{align*}
which is $0$ since $\beta(v,y)=\epsilon \beta(y,v)$ and $\beta(w,x)=
\epsilon \beta(x,w)$.  

Now let $W\subseteq V$ be a $\gom(V,\beta)$-submodule and assume, if 
possible, that $\{0\}\neq W\neq V$. Let $0\neq w\in W$. Since $\beta$ 
is non-degenerate, we have $\beta(y,w) \neq 0$ for some $y\in V$. If 
$x\in V$ is such that $\beta(x,w)=0$, then $\varphi_{x,y}(w)=\beta(w,x)
y-\beta(y,w)x=-\beta(y,w)x$. But then $\varphi_{x,y}(w)\in W$ 
(since $W$ is a submodule) and so $x\in W$. Thus, 
\[ U_w:=\{x\in V\mid \beta(x,w)=0\}\subseteq W.\]
Since $U_w$ is defined by a single, non-trivial linear equation, we have
$\dim U_w=\dim V-1$ and so $\dim W\geq \dim V-1$. Since $W\neq V$, we 
have $\dim W=\dim U_w$ and $U_w=W$. This holds for all $0\neq w\in W$ 
and so $W\subseteq W^\perp$. Since $\beta$ is non-degenerate, we have 
$\dim V=\dim W+\dim W^\perp$ (by a general result in Linear Algebra); hence,
\[\dim V=\dim W+\dim W^\perp \geq 2\dim W\geq 2(\dim V-1)\]
and so $\dim V\leq 2$, a contradiction. 
\end{proof}

In the sequel, it will be convenient to work with matrix descriptions of
$\gom(V,\beta)$; these are provided by the following exercise.

\begin{xca} \label{xcaclassic1} Let $n=\dim V<\infty$ and $B=\{v_1,
\ldots,v_n\}$ be a basis of $V$. We form the corresponding Gram matrix 
\[Q=\bigl(\beta(v_i,v_j)\bigr)_{1\leq i,j\leq n} \in M_n(k).\]
The following equivalences are well-known from Linear Algebra:
\begin{align*}
Q^{\text{tr}}=Q \;\;\;\quad&\Leftrightarrow\quad \mbox{$\beta$ symmetric},\\
Q^{\text{tr}}=-Q \quad &\Leftrightarrow \quad \mbox{$\beta$ alternating},
\\ \det(Q)\neq 0 \quad &\Leftrightarrow \quad 
\mbox{$\beta$~non-degenerate}.
\end{align*}
Recall that we are assuming $\mbox{char}(k)\neq 2$.

\noindent (a) Let $\varphi\in \End(V)$ and $A=(a_{ij})\in M_n(k)$ be the 
matrix of $\varphi$ with respect to $B$. Then show that $\varphi\in\gom(V,
\beta)\Leftrightarrow A^{\text{tr}}Q+QA=0$, where $A^{\text{tr}}$ denotes 
the transpose matrix. Hence, we obtain a Lie subalgebra
\[\gom_n(Q,k):=\{A\in M_n(k)\mid A^{\text{tr}}Q+QA=0\}\subseteq\gl_n(k).\]
Deduce that $V=k^n$ is an irreducible $\gom_n(Q,k)$-module if
$Q^{\text{tr}}=\pm Q$, $\det(Q)\neq 0$ and $n\geq 3$.

\noindent (b) Show that if $\det(Q)\neq 0$, then $\gom_n(Q,k)
\subseteq \slm_n(k)$. (In particular, for $n=1$, we have $\gom_1(Q,k)
=\{0\}$ in this case.)
\end{xca}

\begin{prop} \label{slnss} Let $n\geq 3$ and $k=\C$. If 
$Q^{\operatorname{tr}}= \pm Q$ and $\det(Q)\neq 0$, then 
$\gom_n(Q,\C)$ is semisimple.
\end{prop} 

\begin{proof} This follows from Exercise~\ref{xcaclassic1} and 
the semisimplicity criterion in Proposition~\ref{critss}.
\end{proof}

Depending on what $Q$ looks like, computations in $\gom_n(Q,k)$ can be 
more, or less complicated. Let us assume from now on that $k=\C$, $n=
\dim V<\infty$ and $Q$ is given by\footnote{If $k=\C$ and $\beta$ is 
reflexive and non-degenerate, then one can always find a basis $B$ of
$V$ such that $Q$ has this form. For $\beta$ alternating, this holds 
even over any field $k$; see \cite[Theorem~2.10]{Grove2}. For $\beta$ 
symmetric, this follows from the fact that, over $\C$, any two 
non-degenerate symmetric bilinear forms are equivalent; see
\cite[Theorem~4.4]{Grove2}.}
\[\renewcommand{\arraystretch}{0.8}
Q=Q_n:=\left(\begin{array}{c@{\hspace{5pt}}c@{\hspace{5pt}}
c@{\hspace{5pt}}c} 0 & \cdots & 0 & \delta_n \\ \vdots &\dddots & 
\dddots & 0 \\ 0 & \delta_2 &\dddots & \vdots \\ \delta_1 & 0 & \cdots 
& 0 \end{array} \right)\in M_n(\C)\qquad (\delta_i\in \{\pm 1\}),\]
where $\delta_i\delta_{n+1-i}=\epsilon$ for all~$i$ and, hence, 
$Q_n=\epsilon Q_n^{\text{tr}}$, $\det(Q_n)\neq 0$. Note that this
forces $\epsilon=1$ if $n$ is odd.

\begin{xca} \label{xcaclassic3} (a) Assume that $n=2$. Show that 
$\gom_2(Q_2,\C)=\slm_2(\C)$ if $Q_2=-Q_2^{\text{tr}}$. Hence, in this 
case, $V=\C^2$ still is an irreducible $\gom_2(Q_2,\C)$-module (see
Example~\ref{slnss0}). Also show that 
\[\renewcommand{\arraystretch}{0.9} 
\gom_2(Q_2,\C) =\left\{\left(\begin{array}{r@{\hspace{5pt}}r} c & 0 \\ 0 
& -c\end{array} \right)\;\Big|\; c\in \C\right\} \qquad \mbox{if 
$Q_2=Q_2^{\text{tr}}$},\]
and so $V=\C^2$ is not an irreducible $\gom_2(Q_2,\C)$-module in this case.

(b) Assume that $n=3$ and $Q_3=Q_3^{\text{tr}}$. Show that 
\[\renewcommand{\arraystretch}{0.9} \gom_3(Q_3,\C)=\left\{\left(
\begin{array}{c@{\hspace{5pt}}r@{\hspace{5pt}}r}  a & b & 0 \\
c & 0 & -\delta b\\ 0 & -\delta c & -a \end{array}\right)\;\Big|\;
a,b,c\in \C \right\}\quad\; (\delta:=\delta_1\delta_2)\]
is isomorphic to $\slm_2(\C)$. 

(c) Assume that $n=4$ and $Q_4=Q_4^{\text{tr}}$. Show that
\[\renewcommand{\arraystretch}{0.9} L_1:=\left\{
\left(\begin{array}{r@{\hspace{5pt}}r@{\hspace{5pt}}r@{\hspace{5pt}}r} a& 0 
& b & 0 \\ 0 & a & 0 & -b\\ c & 0 & -a & 0 \\ 0 & -c & 0 & -a
\end{array}\right)\;\Big|\;a,b,c\in \C \right\}\subseteq \gom_4(Q_4,\C)\]
is an ideal and $L_1\cong\slm_2(\C)$. Show that $\gom_4(Q_4,\C) \cong 
\slm_2(\C) \times \slm_2(\C)$ (where the direct product of two algebras is 
defined in Example~\ref{dirprodalg}). 
\end{xca}

\begin{exmp} \label{classic3a} We have the following implication:
\[ A \in \gom_n(Q_n,\C) \quad \Rightarrow \quad A^{\text{tr}}\in 
\gom_n(Q_n,\C).\]
Indeed, if $A^{\text{tr}}Q_n+Q_nA=0$, then $Q_n^{-1}A^{\text{tr}}+
AQ_n^{-1}=0$. Now note that $Q_n^{-1}=Q_n^{\text{tr}}=
\epsilon Q_n$. Hence, we also have $Q_nA^{\text{tr}}+AQ_n=0$. 
\end{exmp}

Finally, we determine a vector space basis of $\gom_n(Q_n,\C)$. We set 
\[A_{ij}:= \delta_i E_{ij} -\delta_{j} E_{n+1-j,n+1-i} \in M_n(\C)\]
for $1\leq i,j\leq n$, where $E_{ij}$ denotes the elementary matrix with 
$1$ as its $(i,j)$-entry and zeroes elsewhere. Note that, if $i=j$, then 
$A_{ii}$ is a diagonal matrix. If $i<j$, then $A_{ij}$ is a strictly upper 
triangular matrix and, if $i>j$, then $A_{ij}$ is a strictly lower 
triangular matrix. 

\begin{prop} \label{classic4} Recall that $k=\C$ and $Q=Q_n$ is as above.
We have $A_{ij}\in \gom_n(Q_n,k)$ for all $1\leq i,j\leq n$. Furthermore:
\begin{itemize}
\item[{\rm (a)}] If $Q_n^{\operatorname{tr}}=Q_n$, then 
$\{A_{ij}\mid 1\leq i,j\leq n,i+j\leq n\}$ is a basis of $\gom_n(Q_n,
\C)$ and so $\dim \gom_n(Q_n,\C)=n(n-1)/2$.
\item[{\rm (b)}] If $Q_n^{\operatorname{tr}}=-Q_n$, then $\{A_{ij}\mid 1\leq 
i,j \leq n,i+j\leq n+1\}$ is a basis of $\gom_n(Q_n,\C)$ and so $\dim 
\gom_n(Q_n,\C)=n(n+1)/2$.
\end{itemize}
\end{prop}

\begin{proof} Let $1\leq i,j\leq n$. First note that $Q_nE_{ij}=\delta_i
E_{n+1-i,j}$. Hence, 
\begin{align*}
Q_n&A_{ij}  =\delta_iQ_nE_{ij}-\delta_j Q_nE_{n+1-j,n+1-i}\\&=
\delta_i^2E_{n+1-i,j}-\delta_j\delta_{n+1-j}E_{j,n+1-i}=E_{n+1-i,j}
-\epsilon E_{j,n+1-i}.
\end{align*}
Furthermore, $A_{ij}^{\text{tr}}Q_n=\epsilon (Q_nA_{ij})^{\text{tr}}=
\epsilon(E_{n+1-i,j}^{\text{tr}}-\epsilon E_{j,n+1-i}^{\text{tr}})$
and so $A_{ij}^{\text{tr}}Q_n+Q_nA_{ij}=0$, that is, $A_{ij}\in 
\gom_n(Q_n,\C)$ as claimed. 

Now let $A\in M_n(\C)$ be arbitrary. Then $A\in\gom_n(Q_n,
\C)$ if and only if $A^{\text{tr}}Q_n=-Q_nA$. Since $A^{\text{tr}}Q_n=
\epsilon (Q_nA)^{\text{tr}}$, this is equivalent to the condition
$(Q_nA)^{\text{tr}}=-\epsilon Q_nA$. We obtain a bijective linear map
\[\gom_n(Q_n,\C)\rightarrow \{S\in M_n(\C)\mid S^{\text{tr}}=
-\epsilon S\},\quad A\mapsto Q_nA.\]
If $\epsilon=-1$, then the space on the right hand side consists precisely
of all symmetric matrices in $M_n(\C)$; hence, its dimension equals $n(n+1)
/2$. Similarly, if $\epsilon=1$, then the space on the right hand side 
consists precisely of all skew-symmetric matrices in $M_n(\C)$; hence, its 
dimension equals equals $n(n-1)/2$.  

It remains to prove the statements about bases. All we need to do now is
to find the appropriate number of linearly independent elements.
Consider the set 
\[ I:=\{(i,j)\mid 1\leq i,j\leq n,i+j\leq n\};\] 
note that $|I|=n(n-1)/2$. Now, if $(i,j)\in I$, then $(n+1-i)+
(n+1-j)\geq n+2$ and so $(n+1-j,n+1-i)\not\in I$. This implies that
the set $\{A_{ij} \mid (i,j)\in I\} \subseteq \gom_n(Q_n,\C)$ is 
linearly independent. Furthermore, for $1\leq i\leq n$, we have 
$(i,n+1-i)\not\in I$, $(n+1-i,i)\not\in I$ and 
\[A_{i}:=A_{i,n+1-i}=\delta_i(1-\epsilon)E_{i,n+1-i}.\]
Hence, if $\epsilon=-1$, then $A_i\neq 0$ and $\{A_{ij}\mid (i,j)\in I\}
\cup\{A_i\mid 1\leq i\leq n\}$ is linearly independent. Thus, (a) and
(b) are proved.
\end{proof}

\begin{rem} \label{classic3b} 
Denote by $\mbox{diag}(x_1,\ldots,x_n) \in M_n(\C)$ the diagonal matrix 
with diagonal coefficients $x_1,\ldots,x_n\in \C$. Then 
\[\mbox{diag}(x_1,\ldots,x_n)\in \gom_n(Q_n,\C)\quad\Leftrightarrow\quad
x_i+x_{n+1-i}=0 \mbox{ for all $i$}.\]
This easily follows by an explicit matrix calculation. Let $H$ be the subspace
of $\gom_n(Q_n,\C)$ consisting of all matrices in $\gom_n(Q_n,\C)$ that are 
diagonal. Let $m\geq 1$ be such that $n=2m+1$ (if $n$ is odd) or $n=2m$ (if 
$n$ is even). Then $H$ consists precisely of all diagonal matrices of the 
form
\[\left\{\begin{array}{cl} \mbox{diag}(x_1,\ldots,x_m,0,-x_m,\ldots,-x_1)
&\qquad \mbox{if $n$ is odd}, \\ \mbox{diag}(x_1,\ldots,x_m,-x_m,\ldots,
-x_1) \;\;\; &  \qquad \mbox{if $n$ is even}. \end{array}\right.\]
In particular, $\dim H=m$. With the above definition of~$m$, the dimension 
formulae in Proposition~\ref{classic4} are re-written as follows:
\[ \dim \gom_n(Q_n,\C)=\left\{\begin{array}{cl} 2m^2-m & \mbox{ if $n=2m$ and
$Q_n^{\text{tr}}=Q_n$},\\ 2m^2+m &\mbox{ otherwise}.\end{array}\right.\]
\end{rem}

\begin{cor}[Triangular decomposition] \nmi{}{triangular decomposition}
\label{triangom} Let $L=\gom_n(Q_n,\C)$, as above. Then every $x\in L$ 
has a unique expression $x=h+n^++n^-$ where $h\in L$ is a diagonal matrix, 
$n^+\in L$ is a strictly upper triangular matrix, and $n^-\in L$ is a 
strictly lower triangular matrix.
\end{cor}

\begin{proof} Note that $A_{ij}$ is diagonal if $i=j$, strictly upper
triangular if $i<j$, and strictly lower triangular if $i>j$. So the
assertion follows from Proposition~\ref{classic4}.
\end{proof}

We shall see later that the algebras $\slm_n(\C)$ and $\gom_n(Q_n,\C)$
are not only semisimple but simple (with the exceptions in 
Exercise~\ref{xcaclassic3}(a) and (c)). The following result highlights 
their importance.

\begin{thm}[Cartan--Killing Classification] \label{cartkill}
\nmi{}{Cartan--Killing classification} Let $L\neq \{0\}$ be a semisimple 
Lie algebra over $\C$ with $\dim L<\infty$. Then $L$ is a direct product 
of simple Lie algebras, each of which is isomorphic to either $\slm_n(\C)$ 
($n\geq 2$), or $\gom_n(Q_n,\C)$ ($n \geq 4$ and $Q_n$ as above), or to 
one of five ``exceptional'' algebras that are denoted by $G_2$, $F_4$, 
$E_6$, $E_7$, $E_8$ and are of dimension $14$, $52$, $78$, $133$, $248$, 
respectively.
\end{thm}

This classification result is proved in textbooks like those of Carter 
\cite{Ca3}, Erdmann--Wildon \cite{EW} or Humphreys \cite{H}, to mention 
just a few (see also Bourbaki \cite{B78} and Jacobson \cite{Jac}). It is 
achieved as the culmination of an elaborate chain of arguments. Here,
we shall take a shortcut around that proof. Following Moody--Pianzola 
\cite{MP}, we will work in a setting where the existence of something
like a ``triangular decomposition'' (as in Corollary~\ref{triangom}) is 
systematically adopted at the outset. Such a decomposition is an essential 
feature of semisimple Lie algebras, and also of the wider class of 
Kac--Moody algebras (which we will very briefly encounter in 
Section~\ref{sec3a3}). 

This approach provides a uniform framework for studying the various 
Lie algebras appearing in Theorem~\ref{cartkill}. It is still completely 
self-contained; no prior knowledge about simple Lie algebras 
is required. One advantage is that it allows us to reach more directly
the point where we can deal with certain more modern aspects of the theory 
of Lie algebras, and with the construction of Chevalley groups.

\section*{Notes on Chapter~\ref{chap0}}

For further reading about \nmi{Witt algebras}{Witt algebra} see, e.g., 
Moody--Pianzola \cite[\S 1.4]{MP}, Schottenloher \cite[\S 5.1]{schott} 
(connections with mathematical physics) and Strade \cite{Strade} 
(connections with Lie algebras over fields of positive characteristic). 
There are also quite remarkable connections of Lie algebras with 
finite group theory, especially the ``Monster'' sporadic simple group; see 
Frenkel--Lepowsky--Meurman \cite{FLM}. Another source for Lie algebras 
(over fields of characteristic $p>0$) is the theory of finite $p$-groups; 
see De Graaf \cite[\S 1.4]{graaf} for an introduction. 

Proposition~\ref{killing8b} is a standard result about semisimple Lie
algebras. Several authors (e.g., Erdmann--Wildon \cite[\S 9]{EW} or
Humphreys \cite[\S 5.2]{H}) use ``\nm{Cartan's First Criterion}'' in 
the proof (which we did not establish here), but this is actually not 
necessary; see, e.g., the argument in Bourbaki \cite[\S 3, no.~6, 
Prop.~7]{Bchap1} (which we adopt here).

Exercise~\ref{xcarepwitt} is taken from \cite[Chap.~1, Exc.~1.6]{MP}.
The usual proofs for Proposition~\ref{critss} and Theorem~\ref{liethm}
proceed somewhat differently; see, e.g., \cite[\S 6.4 and p.~245]{EW}. 
As far as $\slm_2(\C)$-modules are concerned, for any $m\geq 0$, there
exists an irreducible $\slm_2(\C)$-module of dimension $m+1$ and with
a corresponding matrix representation as in Corollary~\ref{sl2modc}.
This is, perhaps, best explained in terms of a natural action of 
$\slm_2(\C)$ on homogeneous polynomials of degree $m$ in two 
indeterminates; see, for example, \cite[Chap.~8]{EW}. 

For further properties of reflexive bilinear forms, see 
Grove \cite{Grove2} and the references there. The elegant proof of 
Proposition~\ref{classic1} is taken from Tauvel--Yu 
\cite[Theorem~20.2.3]{TY}. The study of the classical Lie algebras
will be continued in Section~\ref{sec1a5}.


\chapter{Semisimple Lie algebras} \label{chap2}

\addtocounter{footnote}{2}

Before we start, a word about notation. It is quite common in the 
literature (e.g., Bourbaki \cite{B78}, Fulton--Harris \cite{FH} or Kac 
\cite{K}) to use small Fraktur style letters to denote Lie algebras. We 
will adopt that convention from now on, and reserve capital letters for 
groups, vector spaces etc.; in this way, we can also avoid, later on, 
any conflict between the notation for groups and for Lie algebras.

The classical Cartan--Killing theory shows that a finite-di\-men\-sional
semisimple Lie algebra $\fg$ over $\C$ has a ``Cartan decomposition''. 
This is a direct sum decomposition $\fg=\fh\oplus \bigoplus_{\alpha\in 
\Phi} \fg_\alpha$, where $\fh \subseteq \fg$ is an abelian subalgebra 
and $\Phi$ is a subset of the dual space $\fh^*=\Hom(\fh,\C)$, such that 
all the $\fg_\alpha$ are one-dimensional ``weight spaces'' for $\fh$ and 
$\Phi$ is an abstract ``root system''. Separating $\Phi$ into a positive 
and a negative part, we obtain a ``triangular decomposition'' of~$\fg$.

In our prime examples $\slm_n(\C)$ and $\gom_n(Q_n,\C)$, it is relatively
straightforward to exhibit the ingredients of a Cartan decomposition 
as above (see Example~\ref{cartsln} and Section~\ref{sec1a5}.) For 
a general~$\fg$, this typically forms a core chapter of the classical
treatment in books like those mentioned at the end of Chapter~\ref{chap0}.
Here, we propose to take a shortcut by adopting an axiomatic setting from 
which the Cartan decomposition can be deduced without too much effort. 

After some preliminaries about weight spaces in Section~\ref{sec1a1}, 
the required axioms are formulated in Definition~\ref{defTD}. Readers 
familiar with the theory of Kac--Moody algebras will recognise
the influence of Kac \cite{K} and Moody--Pianzola \cite{MP} on the
formulation of those axioms. Sections~\ref{sec1a3} and \ref{sec1a4} 
introduce the all-important Weyl group and establish basic structural 
properties of Lie algebras satisfying our Definition~\ref{defTD}, e.g., 
the fact that they are indeed semisimple and that the Killing form is
non-degenerate. 

One of the most spectacular developments in the modern theory of Lie 
algebras is the introduction of ``quantum groups'' (see, e.g., 
Chari--Pressley \cite{CP}) and the subsequent discovery of 
``canonical bases'' and ``crystal bases''; see Lusztig \cite{L6}, 
\cite{Lhistory} and Kashiwara \cite{Kas2}. It would be beyond the scope 
of this text to introduce these ideas in any more detail. But, in 
Section~\ref{sec1a7}, we can at least present one highlight of this theory 
in a completely elementary fashion: Lusztig's ``canonical basis'' for the 
Lie algebra $\cL$ itself.  This constitutes a highly remarkable
strengthening of the existence of integral bases of $\cL$ due to
Chevalley \cite{Ch}.

Throughout this chapter, we work over the base field $k=\C$.

\section{Weights and weight spaces} \label{sec1a1}

Throughout this section, we let $\fh$ be a finite-dimensional 
\underline{\itshape abelian} Lie algebra. Let $\rho \colon \fh \rightarrow
\gl(V)$ be a representation of $\fh$ on a finite-dimensional vector space
$V\neq \{0\}$ (all over $k=\C$). Thus, $V$ is an $\fh$-module as in 
Section~\ref{sec04}. By Proposition~\ref{abelsimple}, there exists a 
basis $B$ of $V$ such that, for any $x\in \fh$, the matrix of the 
linear map $\rho_x\colon V \rightarrow V$, $v \mapsto x.v$, with respect 
to $B$ has an upper triangular shape as follows: 
\[\renewcommand{\arraystretch}{0.7} M_B(\rho_x)=\left(
\begin{array}{c@{\hspace{3pt}}c@{\hspace{3pt}}c@{\hspace{3pt}}c} 
\lambda_1(x) & * & \ldots & *\\ 0 & \lambda_2(x) & \ddots & \vdots \\ 
\vdots & \ddots &\ddots & * \\ 0 & \ldots & 0 & \lambda_n(x)\end{array}
\right)\qquad (n=\dim V),\]
where $\lambda_i\in \fh^*:=\Hom(\fh,\C)$ are linear maps for $1\leq i\leq n$.
By Lemma~\ref{liethmu}, the set $P_\fh(V):=\{\lambda_1,\ldots,\lambda_n\}
\subseteq \fh^*$ does not depend on the choice of the basis $B$ and is called 
the set of \nmi{weights}{weight} of $\fh$ on~$V$. A particularly favourable 
situation occurs when the matrices $M_B(\rho_x)$ are diagonal for all 
$x\in \fh$. This leads to the following definition.

\begin{defn} \label{wsddef1} In the above setting (with $\fh$ abelian), we 
say that the $\fh$-module $V$ is $\fh$-\nm{diagonalisable} if, for each
$x\in \fh$, the linear map $\rho_x\colon V \rightarrow V$ is diagonalisable,
that is, there exists a basis of $V$ such that the corresponding matrix of 
$\rho_x$ is a diagonal matrix (but, a priori, the basis may depend on the 
element $x\in \fh$).
\end{defn} 

A linear map $\rho\colon \fh \rightarrow \End(V)$ is a representation of
Lie algebras if and only if $\rho([x,x'])=\rho(x)\circ \rho(x')-
\rho(x')\circ \rho(x)$ for all $x,x'\in \fh$. Since $\fh$ is abelian, this 
just means that the maps $\{\rho(x) \mid x \in\fh\}\subseteq \End(V)$ 
commute with each other. Thus, the following results are really 
statements about commuting matrices, but it is useful to formulate them 
in terms of the abstract language of modules for Lie algebras in view 
of the later applications to ``weight space decompositions''.

\begin{lem} \label{wsdlem1} Assume that $V$ is $\fh$-diagonalisable. Let
$U \subseteq V$ be an $\fh$-submodule. Then $U$ is also 
$\fh$-diagonalisable.
\end{lem}

\begin{proof} Let $x\in \fh$ and $\lambda_1,\ldots,\lambda_r\in \C$ (where
$r \geq 1$) be the distinct eigenvalues of $\rho_x\colon V \rightarrow V$. 
Then $V=V_1+\ldots + V_r$ where $V_i$ is the $\lambda_i$-eigenspace of 
$\rho_x$. Setting $U_i:=U \cap V_i$ for $1\leq i \leq r$, we claim that 
$U=U_1+\ldots +U_r$. Indeed, let $u \in U$ and write $u=v_1+
\ldots +v_r$ where $v_i \in V_i$ for $1\leq i \leq r$. We must show that
$v_i \in U$ for all~$i$. For this purpose, we define a sequence of vectors 
$(u_j)_{j \geq 1}$ by $u_1:=u$ and $u_j:=x.u_{j-1}$ for $j \geq 2$. Then 
a simple induction on $j$ shows that 
\[ u_j=\lambda_1^{j-1}v_1+\ldots +\lambda_r^{j-1}v_r \qquad \mbox{for all 
$j\geq 1$}.\]
Since the Vandermonde matrix $\bigl(\lambda_i^{j-1}\bigr)_{1\leq i,j 
\leq r}$ is invertible, we can invert the above equations (for $j=1,
\ldots,r$) and find that each $v_i$ is a linear combination of $u_1,\ldots,
u_r$. Since $U$ is an $\fh$-submodule of $V$, we have $u_j \in U$ for all~$j$,
and so $v_i \in U$ for all~$i$, as claimed. 

Now $U_i=U\cap V_i=\{u \in U \mid x.u=\lambda_iu\}$ for all~$i$. Hence, 
all non-zero vectors in $U_i$ are eigenvectors of the restricted map 
$\rho_x|_U\colon U \rightarrow U$. Consequently, $U=U_1+\ldots +U_r$ is 
spanned by eigenvectors for $\rho_x|_U$ and, hence, $\rho_x|_U$ is 
diagonalisable.
\end{proof}

\begin{prop} \label{wsdlem2} Assume that $V$ is $\fh$-diagonalisable; let
$n=\dim V\geq 1$. Then there exist $\lambda_1,\ldots,\lambda_n\in \fh^*$ 
and one basis~$B$ of $V$ such that, for all $x\in \fh$, the matrix
of $\rho_x\colon V \rightarrow V$ with respect to $B$ is diagonal, with
$\lambda_1(x),\ldots, \lambda_n(x)$ along the diagonal. 
\end{prop}

\begin{proof} We proceed by induction on $\dim V$. If $\rho_x$ is a scalar
multiple of the identity for all $x\in \fh$ then the result is clear.
In particular, this covers the case where $\dim V=1$. Now assume that 
$\dim V>1$ and that there exists some $y \in \fh$ such that $\rho_y$ is
not a scalar multiple of the identity. Since $\rho_y$ is diagonalisable
by assumption, there are at least two distinct eigenvalues. So let 
$\lambda_1,\ldots,\lambda_r\in \C$ be the distinct eigenvalues of
$\rho_y$, where $r \geq 2$. Then $V=V_1\oplus \ldots \oplus V_r$
where $V_i$ is the $\lambda_i$-eigenspace of $\rho_y$. We claim that
each $V_i$ is an $\fh$-submodule of $V$. Indeed, let $v \in V_i$ and
$x\in \fh$. Since $\fh$ is abelian, we have $\rho_x\circ \rho_y=
\rho_y \circ \rho_x$. This yields 
\begin{align*}
\rho_y(x.v)&=\rho_y(\rho_x(v))=(\rho_y\circ \rho_x)(v)=(\rho_x\circ 
\rho_y)(v)\\&=\rho_x(y.v)=\lambda_i(y)\rho_x(v)=\lambda_i(y)(x.v)
\end{align*}
and so $x.v \in V_i$. By Lemma~\ref{wsdlem1}, each $V_i$ is 
$\fh$-diagonalisable. Now $\dim V_i<\dim V$ for all~$i$. So,
by induction, there exist bases $B_i$ of $V_i$ such that the matrices
of $\rho_x|_{V_i}\colon V_i \rightarrow V_i$ are diagonal for all
$x\in \fh$. Since $V=V_1\oplus \ldots \oplus V_r$, the set $B:=B_1\cup
\ldots \cup B_r$ is a basis of $V$ with the required property.
\end{proof}

Given $\lambda \in \fh^*$, a non-zero vector $v \in V$ is called a
\nm{weight vector} (with weight $\lambda)$ if $x.v=\lambda(x)v$ for all
$x\in \fh$. We set 
\[V_\lambda:=\{v \in V \mid x.v=\lambda(x)v \mbox{ for all $x\in \fh$}\}.\]
Clearly, $V_\lambda$ is a subspace of $V$. If $V_\lambda\neq \{0\}$, then 
$V_\lambda$ is called a \nm{weight space} for $\fh$ on $V$. With this
notation, we have the equivalence:
\begin{center}
\fbox{$\mbox{$V$ is $\fh$-diagonalisable} \quad\Leftrightarrow \quad V=
{\textstyle\sum_{\lambda \in \Lambda}} V_\lambda\quad\mbox{for some 
$\Lambda \subseteq \fh^*$};$}
\end{center}
furthermore, $P_\fh(V)\subseteq \Lambda$ in this case. This is almost 
trivial, but let us go once more through the arguments. First, if $V$ is 
$\fh$-diagonalisable, then we are in the setting of 
Proposition~\ref{wsdlem2}; hence, each vector of the basis~$B$ belongs 
to $V_{\lambda_i}$ for some~$i$ and so $V=\sum_{\lambda \in \Lambda} 
V_\lambda$ where $\Lambda=P_\fh(V)=\{\lambda_1, \ldots, \lambda_n\}$, as 
desired. Conversely, if $V=\sum_{\lambda\in \Lambda} V_{\lambda}$ for
some $\Lambda\subseteq \fh^*$, then we can extract a basis $B$ from the 
union of subspaces $\bigcup_{\lambda \in \Lambda} V_{\lambda}$. 
With respect to $B$, the map $\rho_x \colon V \rightarrow V$ is represented by 
a diagonal matrix, for every $x\in \fh$. Hence, $V$ is $\fh$-diagonalisable 
where $P_\fh(V)\subseteq \Lambda$.

\begin{exmp} \label{wsdexp1} Assume that $V$ is $\fh$-diagonalisable. 
Let $U,U'\subseteq V$ be $\fh$-submodules such that $V=U\oplus U'$.
By Lemma~\ref{wsdlem1}, both $U$ and $U'$ are also $\fh$-diagonalisable.
We claim that
\[ V_\lambda=U_\lambda \oplus U_\lambda' \qquad\mbox{for all $\lambda\in 
\fh^*$}.\]
Indeed, the inclusion ``$\supseteq$'' is clear. Conversely, let $v \in 
V_\lambda$ and write $v=u+u'$ where $u\in U$ and $u'\in U'$. For $x\in 
\fh$ we have $x.v=x.u+x.u'$ where $x.u\in U$ and $x.u'\in U'$. Since the 
left hand side equals $\lambda(x)v=\lambda(x)u+\lambda(x)u'$, and since
$V=U \oplus U'$, we conclude that $x.u=\lambda(x)u$ and $x.u'=\lambda(x)u'$. 
Hence, $v\in U_\lambda+U_\lambda'$.
\end{exmp}

\begin{prop}\label{wsdprop1} Assume that $V$ is $\fh$-diagonalisable. Recall
the definition of the set of weights $P_\fh(V)\subseteq \fh^*$ (see 
Lemma~\ref{liethmu}).
\begin{itemize}
\item[{\rm (a)}] For $\lambda \in \fh^*$, we have $\lambda \in P_\fh(V)$ if 
and only if $V_\lambda\neq \{0\}$.
\item[{\rm (b)}] We have $V=\bigoplus_{\lambda\in P_\fh(V)} V_\lambda$.
\item[{\rm (c)}] If $U\subseteq V$ is an $\fh$-submodule, then $U=
\bigoplus_{\lambda \in P_\fh(U)} U_\lambda$ where $P_\fh(U)\subseteq 
P_\fh(V)$ and $U_\lambda=U\cap V_\lambda$ for all $\lambda \in P_\fh(U)$.
\end{itemize}
\end{prop}

\begin{proof} Let $n=\dim V$. By Proposition~\ref{wsdlem2}, there exist
a basis $B=\{v_1,\ldots,v_n\}$ of $V$ and $\lambda_1,\ldots,\lambda_n
\in \fh^*$ such that $v_i \in V_{\lambda_i}$ for all~$i$. Hence, we have
$V=\sum_{1\leq i \leq n} V_{\lambda_i}$ and $P_\fh(V)=\{\lambda_1,\ldots,
\lambda_n\}$. 

(a) If $\lambda \in P_\fh(V)$, then $\lambda=\lambda_i$ and $v_i\in 
V_{\lambda_i}$ for some~$i$; hence, $V_\lambda=V_{\lambda_i}\neq \{0\}$. 
Conversely, if $V_\lambda \neq \{0\}$, then let $0\neq v\in V_\lambda$.
We also have $v \in V=\sum_{1\leq i \leq n} V_{\lambda_i}$ and so 
Exercise~\ref{xcaweights} below shows that $\lambda=\lambda_i$ for some~$i$.

(b) The $\lambda_i$ need not be distinct. So assume that $|P_\fh(V)|=r\geq 1$ 
and write $P_\fh(V)=\{\mu_1,\ldots,\mu_r\}$; then $V=\sum_{1\leq i 
\leq r} V_{\mu_i}$. We now show that the sum is direct. If $r=1$, there
is nothing to prove. So assume now that $r \geq 2$ and consider the 
finite subset
\[ \{\mu_i-\mu_j \mid 1\leq i<j\leq r\}\subseteq \fh^*.\]
By Exercice~\ref{xca000}, we can choose $x_0\in \fh$ such that all
elements of that subset have a non-zero value on $x_0$. Thus, $\mu_1(x_0),
\ldots,\mu_r(x_0)$ are all distinct. Then $V=V_1\oplus \ldots
\oplus V_r$ where $V_i$ is the $\mu_i(x_0)$-eigenspace of~$V$. Now,
we certainly have $V_{\mu_i}\subseteq V_i$ for all~$i$. Since
$V=\sum_{1\leq i \leq r} V_{\mu_i}$ and $\sum_{1\leq i \leq r}
\dim V_i=\dim V$, it follows that $V_{\mu_i}=V_i$ for all~$i$. 

(c) By Lemma~\ref{wsdlem1}, $U$ is $\fh$-diagonalisable. So, applying (b) 
to~$U$, we obtain that $U=\bigoplus_{\lambda \in P_\fh(U)} U_\lambda$.
Now, we certainly have $U_\lambda=U\cap V_\lambda$ for any $\lambda \in 
\fh^*$. Using (a), this shows that $P_\fh(U)\subseteq P_\fh(V)$. 
\end{proof}
  
\begin{xca} \label{xcaweights} Let $V$ be any $\fh$-module (with $\dim V
<\infty$). Let $r\geq 1$ and $\lambda,\lambda_1,\ldots, \lambda_r\in 
\fh^*$. Assume that $0\neq v \in V_\lambda$ and $v\in \sum_{1\leq i \leq r} 
V_{\lambda_i}$. Then show that $\lambda=\lambda_i$ for some~$i$. (This 
generalizes the familiar fact that eigenvectors corresponding to 
pairwise distinct eigenvalues are linearly independent.)
\end{xca}

Now assume that $\fh$ is a subalgebra of a larger Lie algebra $\cL$ 
with $\dim \cL<\infty$. Then $\cL$ becomes an $\fh$-module via the 
restriction of $\ad_\cL\colon \cL\rightarrow \gl(\cL)$ to $\fh$. So, 
for any $\lambda\in \fh^*$, we have 
\[\cL_\lambda=\{y\in\cL\mid [x,y]=\lambda(x)y\mbox{ for all $x\in\fh$}\}.\]
In particular, $\cL_{\underline{0}}=C_\cL(\fh):=\{y \in \cL\mid [x,y]=0 
\mbox{ for all $x\in \fh$}\}\supseteq \fh$, where $\underline{0} \in\fh^*$ 
denotes the $0$-map. If $\cL$ is $\fh$-diagonalisable, then we can apply 
the above discussion and obtain a decomposition
\[ \cL=\bigoplus_{\lambda \in P_\fh(\cL)} \cL_\lambda \quad \mbox{where 
$P_\fh(\cL)$ is the set of weights of $\fh$ on $\cL$}.\]

\begin{prop} \label{wsdprop2} We have $[\cL_\lambda,\cL_\mu]\subseteq 
\cL_{\lambda+\mu}$ for all $\lambda,\mu\in \fh^*$; furthermore, 
$\cL_{\underline{0}}$ is a subalgebra of $\cL$. If $\cL$ is 
$\fh$-diagonalisable, then we have the equivalence: $\;\cL=[\cL,\cL]\;
\Leftrightarrow\; \cL_{\underline{0}}=\sum_{\lambda \in P_\fh(\cL)}
[\cL_\lambda,\cL_{-\lambda}]$.
\end{prop}

\begin{proof} Let $v\in \cL_\lambda$ and $w\in \cL_\mu$. Thus,
$[x,v]=\lambda(x)v$ and $[x,w]=\mu(x)w$ for all $x\in \fh$. Using
anti-symmetry and the Jacobi identity, we obtain that 
\begin{align*}
[x,[v,w]]&=-[v,[w,x]]-[w,[x,v]]=[v,[x,w]]+[[x,v],w]\\&=
\mu(x)[v,w]+\lambda(x)[v,w]=(\lambda(x)+\mu(x))[v,w]
\end{align*}
for all $x\in \fh$ and so $[v,w]\in \cL_{\lambda+\mu}$. Furthermore, since 
$\fh$ is abelian, $\fh\subseteq \cL_{\underline{0}}=\{y \in \cL \mid [x,y]=0
\mbox{ for all $x\in \fh$}\}$. We have $[\cL_{\underline{0}}, 
\cL_{\underline{0}}] \subseteq \cL_{\underline{0}}$ and so 
$\cL_{\underline{0}}\subseteq \cL$ is a subalgebra. Now assume that $\cL$ 
is $\fh$-diagonalisable and write $P=P_\fh(\cL)$. Since $\cL=
\bigoplus_{\lambda \in P} \cL_\lambda$, we have 
\[[\cL,\cL]=\langle [x,y]\mid x,y\in \cL\rangle_\C=\sum_{\lambda,\mu\in P}
[\cL_\lambda,\cL_\mu].\]
Since $[\cL_\lambda,\cL_\mu]\subseteq \cL_{\lambda+\mu}$ for all 
$\lambda,\mu$, we obtain that 
\[[\cL,\cL]\subseteq \sum_{\lambda\in P}[\cL_\lambda,\cL_{-\lambda}]+
\sum_{\atop{\lambda,\mu \in P}{\lambda\neq -\mu}} \cL_{\lambda+\mu}
\subseteq \sum_{\lambda\in P} \underbrace{[\cL_\lambda,
\cL_{-\lambda}]}_{\subseteq \cL_{\underline{0}}}+\sum_{0\neq \nu\in P} 
\cL_{\nu}.\]
Hence, if $\cL=[\cL,\cL]$, then we must have $\cL_{\underline{0}}=
\sum_{\lambda\in P} [\cL_\lambda,\cL_{-\lambda}]$. Conversely, assume
that $\cL_{\underline{0}}=\sum_{\lambda\in P} [\cL_\lambda,\cL_{-\lambda}]$. 
Then $\cL_{\underline{0}}\subseteq [\cL,\cL]$. Now let $\lambda \in P$,
$\lambda \neq \underline{0}$. Then there exists some $h \in \fh$ such that
$\lambda(h)\neq 0$. For any $v \in \cL_{\lambda}$ we have $[h,v]=
\lambda(h)v$. So $v$ is a non-zero multiple of $[h,v]\in [\cL,\cL]$. It 
follows that $\cL_\lambda \subseteq [\cL,\cL]$. Consequently, we have 
$\cL=\sum_{\lambda \in P} \cL_\lambda \subseteq [\cL,\cL]$ and so $\cL=
[\cL,\cL]$.
\end{proof}

The following result will be useful to verify $\fh$-diagonalisability.

\begin{lem} \label{weightgen} Let $X\subseteq \cL$ be a non-empty subset 
and form the subalgebra $\cL':=\langle X \rangle_{\operatorname{alg}}
\subseteq \cL$. Assume that there is a subset $\{\lambda_x\mid x\in X\}
\subseteq \fh^*$ such that $x\in \cL_{\lambda_x}$ for all $x\in X$. Then 
$\cL'$ is an $\fh$-diagonalisable $\fh$-submodule of $\fg$, where every 
$\lambda\in P_\fh(\cL')$ is a $\Z_{\geq 0}$-linear combination of 
$\{\lambda_x\mid x\in X\}$.
\end{lem}

\begin{proof} Recall from Section~\ref{sec01} that $\langle X
\rangle_{\operatorname{alg}}=\langle X_n\mid n\geq 1\rangle_\C$,
where $X_n$ consists of all Lie monomials in $X$ of level~$n$. Let 
us also set 
\[\Lambda_n:=\{\lambda\in \fh^*\mid \lambda=\lambda_{x_1}+
\ldots +\lambda_{x_n}\mbox{ for some $x_i\in X$}\}.\]
We show by induction on $n$ that, for each $x\in X_n$, there exists some 
$\lambda\in\Lambda_n$ such that $x\in \cL_\lambda'$. If $n=1$, then this 
is clear by our assumptions on~$X$. Now let $n\geq 2$ and $x\in X_n$. 
So $x=[v,w]$ where $v\in X_i$, $w\in X_{n-i}$ and $1\leq i\leq n-1$. By 
induction, there are $\lambda\in \Lambda_i$ and $\mu \in \Lambda_{n-i}$ 
such that $v\in \cL_\lambda'$ and $w\in \cL_\mu'$. By a computation
analogous to that in the proof Proposition~\ref{wsdprop2}, we see that
$x=[y,z]\in \cL_{\lambda+\mu}'$, where $\lambda+\mu\in \Lambda_{i+(n-i)}
=\Lambda_n$, as desired. We conclude that 
\[ \cL'=\langle X_n\mid n\geq 1\rangle_\C=\sum_{n\geq 1} \sum_{\lambda
\in \Lambda_n} \cL_\lambda'\]
and so $\fg'$ is an $\fh$-diagonalisable $\fh$-submodule of $\fg$.
Furthermore, extracting a basis of $\fg'$ from each subspace $\fg_\lambda'$
occurring in the above sum, we see that $P_\fh(\fg')\subseteq 
\bigcup_{n\geq 1} \Lambda_n$.
\end{proof}

The following result will allow us to apply the exponential construction
in Lemma~\ref{exponential} to many elements in $\cL$.

\begin{lem} \label{wsdnil} Let $\fh\subseteq \cL$ be abelian and $\cL$ 
be $\fh$-diagonalisable. Let $\underline{0}\neq \lambda\in P_\fh(\cL)$ and
$y\in \cL_\lambda$. Then $\ad_\cL(y)\colon \cL\rightarrow \cL$ is nilpotent.
\end{lem}

\begin{proof} Let $\mu \in P_\fh(\cL)$ and $v\in \cL_\mu$. Then $\ad_\cL(y)
(v)=[y,v]\in \cL_{\lambda+\mu}$ by Proposition~\ref{wsdprop2}. A simple 
induction on $m$ shows that $\ad_\cL(y)^m(v)\in  \cL_{m\lambda+\mu}$ for all 
$m\geq 0$. Since $\{m\lambda+\mu\mid m\geq 0\} \subseteq \fh^*$ is an 
infinite subset and $P_\fh(\cL)$ is finite, there is some $m>0$ such that 
$m\lambda+\mu \not\in P_\fh(\cL)$ and so $\ad_\cL(y)^m(v)=0$. Hence, since 
$\cL=\langle \cL_\mu\mid\mu\in P_\fh(\cL)\rangle_\C$, we conclude that 
$\ad_\cL(y)$ is nilpotent (see Exercise~\ref{xcanilp}(a)).
\end{proof}


\begin{exmp} \label{cartangln} Let $\cL=\gl_n(\C)$, the Lie algebra of 
all $n\times n$-matrices over $\C$. A natural candidate for an abelian
subalgebra is 
\[ \fh:=\{x\in \cL\mid \mbox{$x$ diagonal matrix}\} \qquad (\dim \fh=n).\]
For $1\leq i\leq n$, let $\varepsilon_i\in \fh^*$ be the map that sends a 
diagonal matrix to its $i$-th diagonal entry. Then $\{\varepsilon_1,
\ldots,\varepsilon_n\}$ is a basis of $\fh^*$. If $n=1$, then $\cL=\fh$. 
Assume now that $n\geq 2$; then $\fh\subsetneqq \cL$. For $i\neq j$ let 
$e_{ij} \in \cL$ be the matrix with entry~$1$ at position $(i,j)$, and $0$ 
everywhere else. Then a simple matrix calculation shows that 
\begin{equation*}
[x,e_{ij}]=(\varepsilon_i(x)-\varepsilon_j(x))e_{ij}
\qquad \mbox{for all $x\in \fh$}.\tag{a}
\end{equation*}
Thus, $\varepsilon_i-\varepsilon_j\in P_\fh(\cL)$ and $e_{ij}\in 
\cL_{\varepsilon_i-\varepsilon_j}$. Furthermore, we have 
\begin{equation*}
\cL\;\;=\;\;\fh\;\;\oplus \; \bigoplus_{1\leq i,j\leq n\,:\,i\neq j} 
\C e_{ij}, \tag{b}
\end{equation*}
where $\fh\subseteq \fg_{\underline{0}}$ and $\C e_{ij}\subseteq 
\fg_{\varepsilon_i-\varepsilon_j}$. So $\cL$ is $\fh$-diagonalisable, 
where $P_\fh(\cL)=\{\underline{0}\}\cup \{\varepsilon_i- \varepsilon_j
\mid i\neq j\}$. Next, note that the weights $\varepsilon_i-\varepsilon_j$
for $i \neq j$ are pairwise distinct and non-zero. Since there are $n^2-n$
of them, Proposition~\ref{wsdprop1} shows that $\dim \cL=\dim 
\cL_{\underline{0}}+\sum_{i\neq j} \dim \cL_{\varepsilon_i{-} 
\varepsilon_j} \geq n+(n^2-n)=n^2=\dim \cL$. Hence, all the above 
inequalities and inclusions must be equalities. We conclude that 
\begin{equation*}
\cL_{\underline{0}}=\fh \qquad \mbox{and}\qquad \cL_{\varepsilon_i{-}
\varepsilon_j}=\langle e_{ij}\rangle_\C \quad \mbox{for all $i\neq j$}.
\tag{c} 
\end{equation*}
Finally, as in Corollary~\ref{triangom}, we have a \nm{triangular
decomposition} $\cL=\fn^+\oplus \fh\oplus \fn^-$ where $N^+$ is the 
subalgebra consisting of all strictly upper triangular matrices in 
$\gl_n(\C)$ and $\fn^-$ is the subalgebra consisting of all strictly lower 
triangular matrices in $\gl_n(\C)$. This decomposition is reflected in 
properties of $P_\fh(\cL)$ as follows. We set 
\[ \Phi^+=\{\varepsilon_i-\varepsilon_j\mid 1\leq i<j\leq n\}\quad
\mbox{and}\quad \Phi^-:=-\Phi^+.\]
Then $P_\fh(\cL)=\{\underline{0}\}\sqcup \Phi^+\sqcup \Phi^-$ (disjoint union) 
and $\fn^\pm =\bigoplus_{\alpha\in\Phi^\pm} \cL_\alpha$. Thus, the 
decomposition $\cL=\fn^+\oplus \fh\oplus \fn^-$ gives rise to a partition 
of $P_\fh(\cL)\setminus \{\underline{0}\}$ into a ``positive'' part $\Phi^+$ 
and a ``negative'' part $\Phi^-$. We also note that, for $1\leq i<j\leq n$, 
we have
\[ \varepsilon_i-\varepsilon_j=(\varepsilon_i-\varepsilon_{i+1})+
(\varepsilon_{i+1}-\varepsilon_{i+2})+\ldots+ (\varepsilon_{j-1}-
\varepsilon_j).\]
Hence, if we set $\alpha_i:=\varepsilon_{i}-\varepsilon_{i+1}$ for $1\leq 
i\leq n-1$, then 
\begin{equation*}
\Phi^\pm=\bigl\{\pm\bigl(\alpha_{i}+ \alpha_{i+1}+ \ldots +\alpha_{j-1}
\bigr) \mid 1\leq i<j\leq n\bigr\}.\tag{d}
\end{equation*}
Thus, setting $\Delta=\{\alpha_1,\ldots,\alpha_{n-1}\}$, every non-zero 
weight of $\fh$ on $\cL$ can be expressed uniquely as a sum of elements of 
$\Delta$ or of~$-\Delta$. (Readers familiar with the theory of abstract 
root systems will recognise the concept of ``simple roots'' in the above 
properties of $\Delta$; see, e.g., Bourbaki \cite[Ch.~VI, \S 1]{B}.) In 
any case, this picture is the prototype of what is also going on in the 
Lie algebras $\slm_n(\C)$ and $\gom_n(Q_n,\C)$, and this is what we will 
formalise in Definition~\ref{defTD} below. For the further discussion of 
examples, the following remark will be useful.
\end{exmp}

\begin{rem} \label{cartansubgln} Let $\cL\subseteq \gl_n(\C)$ be a 
subalgebra, and $\fh\subseteq \cL$ be the abelian subalgebra consisting of 
all diagonal matrices that are contained in $\cL$. First we claim that
\begin{itemize}
\item[(a)] $\cL$ is $\fh$-diagonalisable. 
\end{itemize}
Indeed, by the previous example, $\ad_{\gl_n(\C)}(x) \colon \gl_n(\C)
\rightarrow\gl_n(\C)$ is diagonalisable for all diagonal matrices $x\in 
\gl_n(\C)$ and, hence, also for all $x\in \fh$. Thus, $\gl_n(\C)$ is 
$\fh$-diagonalisable. Now $[\fh,\cL]\subseteq \cL$ and so $\cL$ is an 
$\fh$-submodule of $\gl_n(\C)$. So $\cL$ is $\fh$-diagonalisable by 
Lemma~\ref{wsdlem1}. Furthermore, we have the following useful criterion:
\begin{itemize}
\item[(b)] We have $\fh=C_\cL(\fh)$ if there exists some $x_0\in \fh$ 
with distinct diagonal entries.
\end{itemize}
Indeed, let $x_0=\mbox{diag}(x_1,\ldots,x_n)\in \fh$ with distinct entries 
$x_i\in \C$ and $y=(y_{ij})\in \cL$ be such that $[x_0,y]=x_0\cdot y-y\cdot 
x_0=0$. Then $x_iy_{ij}=y_{ij}x_j$ for all $i,j$ and so $y_{ij}=0$ for 
$i\neq j$.  Thus, $y$ is a diagonal matrix. Since $y\in \cL$, we have 
$y \in \fh$, as required.

For example, let $\cL=\slm_n(\C)$ where $\fh$ consists of all diagonal
matrices with trace~$0$. In this case, we have $C_\fg(\fh)=\fh$ since the
condition in (b) holds for the diagonal matrix $x_0\in \cL$ with diagonal
entries $1,2, \ldots,n-1,-n(n-1)/2$.

Now let $\cL=\gom_n(Q_n,\C)$ as in Section~\ref{sec05}. Then the matrices
in $\fh\subseteq \cL$ are explicitly described in Remark~\ref{classic3b}.
Writing $n=2m+1$ (if $m$ is odd) or $n=2m$ (if $n$ is even), the condition
in (b) holds for 
\[x_0=\left\{\begin{array}{ll} \mbox{diag}(1,\ldots,m,0,-m,\ldots,-1)
&\quad \mbox{if $n$ is odd}, \\ \mbox{diag}(1,\ldots,m,-m,\ldots,-1) 
& \quad \mbox{if $n$ is even}. \end{array}\right.\]
Hence, again, we have $C_\fg(\fh)=\fh$ in this case.
\end{rem}

\begin{exmp} \label{tdnon} Consider the subalgebra $\cL_\delta\subseteq 
\gl_3(\C)$ in Exercise~\ref{excsolv}, where $0\neq \delta\in\C$; also
assume that $\delta\neq 1$. Now 
\[ \renewcommand{\arraystretch}{0.8}
e=\left(\begin{array}{ccc} 0 & 1 & 0 \\ 0 & 0 & 0 \\ 0 & 0 & 0 
\end{array}\right),\quad 
h:=\left(\begin{array}{ccc} 1 & 0 & 0 \\ 0 & 0 & 0 \\ 0 & 0 & \delta 
\end{array}\right),\quad
f=\left(\begin{array}{ccc} 0 & 0 & 0 \\ 0 & 0 & 0 \\ 0 & 1 & 0 
\end{array}\right)\]
form a basis of $\cL_\delta$ and one checks by an explicit computation that 
\[ [h,e]=e, \qquad [h,f]=\delta f, \qquad [e,f]=0.\]
Hence, we have a triangular decomposition $\cL_\delta=\fn^+\oplus \fh\oplus 
\fn^-$, where 
\[\fn^+=\langle e\rangle_\C, \qquad \fn^-=\langle f\rangle_\C, \qquad
\fh:=\langle h \rangle_\C.\]
We have $C_{\cL_\delta}(\fh)=\fh$ since the condition in 
Remark~\ref{cartansubgln}(b) holds for~$h$. The corresponding weights are 
given by $P_\fh(\cL_\delta)=\{\underline{0},\alpha,\delta \alpha\}$, 
where $\alpha \in \fh^*$ is defined by $\alpha(h)=1$. Thus, if $\delta=-1$, 
then we have a partition of $P_\fh(\cL_\delta)\setminus \{\underline{0}\}$ 
into a ``positive'' and a ``negative'' part (symmetrical to each other). 
On the other hand, if $\delta=2$ (for example), then we only have a 
``positive'' part but no ``negative'' part at all. 
So this example appears to differ from that of $\gl_n(\C)$ in a 
crucial way. We shall see that this difference has to do with the fact 
that $[e,f]=0$, that is, $[\fn^+,\fn^-]=\{0\}$. We also know from 
Exercise~\ref{excsolv} that $\cL_\delta$ is solvable, while $\gl_n(\C)$ 
(for $n\geq 2$) is not. 
\end{exmp}

\section{Lie algebras of Cartan--Killing type} \label{sec1a2}

Let $\cL$ be a finite-dimensional Lie algebra over $k=\C$, and $\fh\subseteq 
\cL$ be an abelian subalgebra. Then we regard $\cL$ as an $\fh$-module via 
the restriction of $\ad_\cL\colon \cL \rightarrow \gl(\cL)$ to~$\fh$. Let 
$P_\fh(\cL) \subseteq \fh^*$ be the corresponding set of weights. Motivated 
by the examples and the discussion in the previous section, we introduce the 
following definition.

\begin{defn}[Cf.\ Kac \protect{\cite[Chap.~1]{K}} and Moody--Pianzola 
\protect{\cite[\S 2.1 and \S 4.1]{MP}}] \label{defTD} 
We say that $(\cL,\fh)$ is of \nm{Cartan--Killing type} if there exists a 
linearly independent subset $\Delta=\{\alpha_i \mid i \in I\}\subseteq\fh^*$ 
(where $I$ is a finite index set) such that the following conditions are 
satisfied.
\begin{itemize}
\item[(CK1)] We have $\Delta\subseteq P_\fh(\fg)$ and $\cL$ is 
$\fh$-diagonalisable, with $\cL_{\underline{0}}=\fh$.
\item[(CK2)] Each $\lambda \in P_\fh(\cL)$ is a $\Z$-linear combination 
of $\Delta=\{\alpha_i \mid i \in I\}$ where the coefficients are either
all $\geq 0$ or all $\leq 0$.
\item[(CK3)] We have $\cL_{\underline{0}}=\sum_{i \in I} [\cL_{\alpha_i},
\cL_{-\alpha_i}]$.
\end{itemize}
We set $\Phi:=\{\alpha \in P_\fh(\cL)\mid \alpha \neq \underline{0}\}$. Thus, 
$\cL=\fh\oplus \bigoplus_{\alpha \in \Phi} \cL_\alpha$, which is called the
\nm{Cartan decomposition} of $\cL$ (but we do not yet know $\dim 
\cL_\alpha$). Then $\fh$ is called a \nm{Cartan subalgebra} and $\Phi$ the 
set of \nm{roots} of $\cL$ with respect to~$\fh$. We may also speak of 
$(\Phi, \Delta)$ as a \nm{based root system}.

We say that $\alpha \in \Phi$ is a \nm{positive root} if 
$\alpha=\sum_{i \in I} n_i \alpha_i$ where $n_i \geq 0$ for all $i \in I$; 
similarly, $\alpha \in \Phi$ is a \nm{negative root} if $\alpha=\sum_{i 
\in I} n_i \alpha_i$ where $n_i \leq 0$ for all $i \in I$. Let $\Phi^+$ 
be the set of all positive roots and $\Phi^-$ be the set of all negative 
roots. Thus, $\Phi=\Phi^+\sqcup \Phi^-$ (disjoint union). 

As far as extreme cases are concerned, we note the following.
If $I=\varnothing$, then $P_\fh(\cL)=\{\underline{0}\}$ by (CK2); furthermore, 
$\cL=\cL_{\underline{0}}=\fh$ by (CK1) and so $\cL=\fh=\{0\}$ by (CK3). On 
the other hand, if $I\neq \varnothing$, then $\fh=\cL_{\underline{0}}
\subsetneqq \cL$ and so $\cL$ is not abelian. 
\end{defn}

\begin{rem} \label{remsemis} By Theorem~\ref{mainideal} below, a Lie 
algebra $\cL$ as in Definition~\ref{defTD} is semisimple; so all of 
the above notions (``Cartan subalgebra'', ``roots'' etc.) are 
consistent with the common usage in the general theory of semisimple 
Lie algebras. Conversely, any semisimple Lie algebra is of 
Cartan--Killing type. This result is in fact proved along with the 
classification result in Theorem~\ref{cartkill} (and it will not
be proved here). 
\end{rem}

\begin{rem} \label{defTDa} In the above setting, let $\fn^{\pm}:=
\sum_{\alpha \in \Phi^{\pm}} \cL_\alpha\subseteq \cL$. First note: If 
$\alpha, \beta\in \Phi^\pm$ are such that $\alpha+\beta\in \Phi$, 
then we automatically have $\alpha+\beta\in \Phi^\pm$. So 
Proposition~\ref{wsdprop2} immediately implies that $\fn^\pm$ are 
subalgebras of $\cL$, such that $[\fh,\fn^{\pm}] \subseteq \fn^{\pm}$. Since 
$\Phi=\Phi^+ \sqcup \Phi^-$ and $\fh=\cL_{\underline{0}}$, we have $\cL=\fn^+
\oplus \fh \oplus \fn^-$. So we are close to having a ``\nm{triangular 
decomposition}'', but it is not yet clear that $\Phi^-=-\Phi^+$; we will 
settle this point in Theorem~\ref{mainthm1} below. Finally, condition 
(CK3) ensures that the equality $\cL=[\cL,\cL]$ holds in the most economical 
way (see again Proposition~\ref{wsdprop2}).
\end{rem}

\begin{rem} \label{defheight} Since $\{\alpha_i\mid i\in I\}$ is linearly
independent, we can define a function $\hgt\colon P_\fh(\cL)\rightarrow \Z$ 
as follows. Let $\lambda\in P_\fh(\cL)$ and write $\lambda= \sum_{i\in I} 
n_i \alpha_i$ where $n_i\in\Z$ for all~$i$. Then set $\hgt(\lambda):=
\sum_{i \in I} n_i\in\Z$; this is called the \nm{height} of~$\lambda$.
Since $\Phi=\Phi^+\sqcup \Phi^-$, we have
\[ \hgt(\lambda)=1 \Leftrightarrow \lambda \in \Delta;\;
\hgt(\lambda)\geq 1 \Leftrightarrow \lambda \in \Phi^+;\;
\hgt(\lambda)\leq -1 \Leftrightarrow \lambda \in \Phi^-.\]
Clearly, if $-\lambda\in P_\fh(\cL)$, then $\hgt(-\lambda)=-\hgt(\lambda)$. 
Also note that, if $\lambda,\mu\in P_\fh(\cL)$ are such that $\lambda+\mu 
\in P_\fh(\cL)$, then $\hgt(\lambda+\mu)=\hgt(\lambda)+\hgt(\mu)$. This 
function is a useful technical tool. Here is one example:
\end{rem}

\begin{prop} \label{borelsub} Consider the subalgebra $\fn^\pm\subseteq \cL$ 
defined in Remark~\ref{defTDa}. Then $\fb^\pm:=\fh+\fn^\pm\subseteq \cL$ is a
solvable subalgebra, and we have $[\fb^\pm,\fb^\pm]=\fn^\pm$. 
\end{prop}

\begin{proof} We only consider $\fb^+$. (The proof for $\fb^-$ is completely
analogous.) Since $[\fh,\fn^+]\subseteq \fn^+$, it is clear that $\fb^+$ is a
subalgebra; furthermore, $[\fb^+,\fb^+]\subseteq [\fh,\fh]+[\fh,\fn^+]+
[\fn^+,\fh]+[\fn^+,\fn^+]\subseteq \fn^+$. Conversely, let $\alpha\in\Phi^+$ 
and $x\in \cL_\alpha$. Since $\alpha\neq \underline{0}$, there exists some 
$h\in \fh$ such that $\alpha(h)\neq 0$. Then $\alpha(h)x=[h,x]\in [\fh,
\cL_\alpha] \subseteq [\fb^+,\fb^+]$ and so $x\in [\fb^+,\fb^+]$. 
Thus, $\cL_\alpha\subseteq [\fb^+,\fb^+]$ for all 
$\alpha\in \Phi^+$ and so $\fn^+\subseteq [\fb^+,\fb^+]$. 

It remains to show that $\fn^+$ is solvable. For this purpose, we write 
$P_\fh(\cL)=\{\lambda_1,\ldots,\lambda_d\}$ where the numbering is chosen
such that $\hgt(\lambda_1)\geq \ldots\geq\hgt(\lambda_d)$. For 
each $i$ let $B_i$ be a basis of $\cL_{\lambda_i}$. Then $B:=B_1 \cup \ldots 
\cup B_d$ is a basis of~$\cL$. Let $x\in \fn^+$. We claim that the matrix of 
$\mbox{ad}_\cL(x) \in \gl(\cL)$ with respect to~$B$ is upper triangular. Since 
$\fn^+=\sum_{\alpha \in \Phi^+} \cL_{\alpha}$, it is enough to consider the 
case where $x\in \cL_\alpha$ for some $\alpha\in \Phi^+$; note that $\hgt
(\alpha)>0$. Now let $b \in B_i$ for some~$i$. Then $[x,b]
\in [\cL_\alpha,\cL_{\lambda_i}]\subseteq \cL_{\alpha+\lambda_i}$. If $\alpha+
\lambda_i\not\in P_\fh(\cL)$, then $[x,b]=0$. Otherwise, $\alpha+\lambda_i=
\lambda_k$ for some $k \in \{1,\ldots,d\}$. Since $\hgt(\lambda_k)=\hgt
(\alpha+\lambda_i)=\hgt(\alpha)+\hgt(\lambda_i)>\hgt(\lambda_i)$, we must 
have $k<i$. Hence, $\mbox{ad}_\cL(x)(b)=[x,b]$ is a linear combination of 
basis elements in $B_1\cup \ldots \cup B_{i-1}$. This means that the matrix 
of $\ad_\cL(x)$ with respect to $B$ is upper triangular, as claimed. 
Now, via the basis $B$, we have an isomorphism $\gl(\cL)\cong \gl_n(\C)$
where $n=\dim \cL$. Thus, $\mbox{ad}_\cL(\fn^+)$ is isomorphic to a 
subalgebra of the solvable subalgebra $\fb_n(\C)\subseteq \gl_n(\C)$ in
Definition~\ref{defborell}(b); hence, $\mbox{ad}_\cL(\fn^+)$ is solvable. 
On the other hand, the kernel of $\mbox{ad}_\cL\colon \fn^+ \rightarrow 
\gl(\cL)$ is contained in $\ker(\mbox{ad}_\cL)=Z(\cL)$, which is abelian. 
Hence, $\fn^+$ itself is solvable (see Lemma~\ref{defsolv1}). 
\end{proof}

The further theory will now be developed from the axioms in
Definition~\ref{defTD}. We begin with the following two basic results.

\begin{lem} \label{wsdlem4} Assume that $\cL$ is $\fh$-diagonalisable. Let 
$\lambda\in \fh^*$ be such that $[\cL_\lambda,\cL_{-\lambda}] \subseteq \fh$. 
If the restriction of $\lambda$ to $[\cL_\lambda,\cL_{-\lambda}]$ is zero, 
then $\ad_\cL(x)=0$ for all $x\in [\cL_\lambda,\cL_{-\lambda}]$.
\end{lem}

\begin{proof} Let $y \in \cL_\lambda$, $z \in \cL_{-\lambda}$, and set 
$x:=[y,z] \in [\cL_\lambda,\cL_{-\lambda}]\subseteq \fh$. Consider the 
subspace $\fs:=\langle x,y,z \rangle_\C \subseteq \cL$. Since 
$\lambda(x)=0$, we have $[x,y]=\lambda(x)y=0$, $[x,z]=-\lambda(x)z=0$
and $[y,z]=x$. Thus, $\fs$ is a subalgebra of $\cL$; furthermore, $[\fs,
\fs]=\langle x \rangle_\C$ and so $\fs$ is solvable. We regard $\cL$ as an 
$\fs$-module via the restriction of $\ad_\cL\colon \cL\rightarrow \gl(\cL)$ 
to~$\fs$. Since $\fs$ is solvable, Lie's Theorem~\ref{liethm} shows that
there is a basis $B$ of $\cL$ such that, for any $s\in \fs$, the matrix of 
$\ad_\cL(s)$ with respect to $B$ is upper triangular. Now $x=[y,z]$ and so 
\[\ad_\cL(x)=\ad_\cL(y)\circ \ad_\cL(z)-\ad_\cL(z)\circ \ad_\cL(y).\] 
Hence, the matrix of $\ad_\cL(x)$ is upper triangular with $0$ along the
diagonal.  But $\ad_\cL(x)$ is diagonalisable and so $\ad_\cL(x)=0$. 
\end{proof}

\begin{lem} \label{wsdlem3} Assume that $\cL$ is $\fh$-diagonalisable. Let 
$\lambda\in \fh^*$ be such that $[\cL_\lambda,\cL_{-\lambda}] \subseteq \fh$ 
and the restriction of $\lambda$ to $[\cL_\lambda,\cL_{-\lambda}]$ is 
non-zero; in particular, $\lambda \neq \underline{0}$ and $\cL_\lambda
\neq\{0\}$. Then we have $\dim \cL_{\pm \lambda}=1$ and $P_\fh(\cL)\cap 
\{n\lambda \mid n\in\Z\}=\{\underline{0}, \pm \lambda\}$.
\end{lem}

\begin{proof} By assumption, there exist elements $e\in \cL_\lambda$ and 
$f\in \cL_{-\lambda}$ such that $h:=[e,f]\in [\cL_\lambda,\cL_{-\lambda}] 
\subseteq \fh$ and $\lambda(h)\neq 0$. Note that $e\neq 0$, $f\neq 0$, 
$h\neq 0$. Replacing $f$ by a scalar multiple if necessary, we may 
assume that $\lambda(h)=2$. Then we have the relations
\[ [e,f]=h,\quad [h,e]=\lambda(h)e=2e, \quad [h,f]=-\lambda(h)f=-2f.\]
Thus, $\fs:=\langle e,h,f\rangle_\C$ is a $3$-dimensional subalgebra of $\cL$ 
that is isomorphic to $\slm_2(\C)$ (see Exercise~\ref{expsolv2}). Let 
$p:=\max\{n\geq 1 \mid \cL_{n\lambda} \neq \{0\}\}$ and consider the subspace
\[M:=\C f\oplus \fh \oplus \cL_\lambda \oplus \cL_{2\lambda}\oplus \ldots 
\oplus \cL_{p\lambda}\subseteq \cL,\]
where $\C f \subseteq \cL_{-\lambda}$, $\fh\subseteq \cL_{\underline{0}}$ and 
some terms $\cL_{n\lambda}$ may be $\{0\}$ for $2\leq n<p$. 
By Proposition~\ref{wsdprop2}, we have $[\cL_{n\lambda},\cL_{m\lambda}]
\subseteq \cL_{(n+m)\lambda}$ for all $n,m\in \Z$. Furthermore, $[f,y]
\in \fh$ for all $y\in \cL_\lambda$ (by assumption), $[x,f]=-\lambda(x)f\in
\C f$ for all $x\in \fh$, and $[\fh,\cL_{n\lambda}] \subseteq \cL_{n\lambda}$ 
for all $n \in\Z$. It follows that $[\fs,M] \subseteq M$ and so $M$ may be 
regarded as an $\fs$-module via the restriction of $\ad_\cL\colon \cL
\rightarrow \gl(\cL)$ to $\fs$. The set of eigenvalues of $h$ on $M$ is 
contained in $\{-2,0,2,4, \ldots,2p\}$, where $-2$ has multiplicity~$1$ as
an eigenvalue and $0,2,2p$ have multiplicity at least~$1$. Now, if we had 
$p\geq 2$, then $-2p$ should also be an eigenvalue by 
Proposition~\ref{sl2modd}, contradiction. So we have $p=1$. But then the 
trace of $h$ on $M$ is $-2+2m$ where $m\geq 1$ is the multiplicity of $2$ 
as an eigenvalue. By Proposition~\ref{sl2modd}, that trace is $0$ and so 
$m=1$. Thus, we have shown that $\dim \cL_{\lambda}=1$ and 
$n\lambda\not\in P_\fh(\cL)$ for all $n\geq 2$.

Finally, note that the assumptions also hold for $-\lambda$ instead of 
$\lambda$; so we also obtain $\dim \cL_{-\lambda}=1$ and $\cL_{-n\lambda}=
\{0\}$ for all $n \geq 2$.
\end{proof}

\begin{prop} \label{wsdprop3} Assume that the conditions in 
Definition~\ref{defTD} hold for $(\cL,\fh)$. Then, for each $i\in I$, we
have $\pm \alpha_i \in \Phi$ and 
\[\dim \cL_{\alpha_i}=\dim \cL_{-\alpha_i}=\dim [\cL_{\alpha_i},
\cL_{-\alpha_i}]=1.\]
There is a unique $h_i \in [\cL_{\alpha_i},\cL_{-\alpha_i}]$ 
with $\alpha_i(h_i)=2$. Furthermore, the set $\Delta=\{\alpha_i\mid i  
\in I\}$ is a basis of $\fh^*$ and $\{h_i\mid i \in I\}$ is a basis of~$\fh$.
\end{prop}

\begin{proof} By (CK1) we have $\alpha_i\in P_\fh(\fg)$; so there is 
some $0\neq e_i\in \fg_{\alpha_i}$, for each $i\in I$. Let $I'$ be the set
of all $i\in I$ such that the restriction of $\alpha_i$ to $[\cL_{\alpha_i},
\cL_{-\alpha_i}]$ is non-zero. Then $\{0\}\neq [\cL_{\alpha_i},
\cL_{-\alpha_i}]\subseteq \cL_{\underline{0}}=\fh$ and $\cL_{-\alpha_i}
\neq \{0\}$ for $i\in I'$. Now let us fix $i\in I'$. By Lemma~\ref{wsdlem3},
we have $\dim \cL_{\alpha_i}=\dim \cL_{-\alpha_i}=1$. So $\cL_{\alpha_i}=
\langle e_i \rangle_\C$ and $\cL_{-\alpha_i}=\langle f_i\rangle_\C$ for
some $f_i\neq 0$. Consequently, $[\cL_{\alpha_i},\cL_{-\alpha_i}]=\langle 
h_i\rangle_\C$ where $0\neq h_i:=[e_i,f_i]$ and $\alpha_i(h_i) \neq 0$.
Replacing $f_i$ by a scalar multiple if necessary, we can assume that
$\alpha_i(h_i)=2$; then $h_i$ is uniquely determined (since $\dim 
[\cL_{\alpha_i},\cL_{-\alpha_i}]= 1$). So (CK1) and (CK3) show that
\[ \fh=\fh'+\langle h_i\mid i \in I'\rangle_\C\quad\mbox{where} \quad
\fh':=\sum_{j\in I\setminus I'} [\cL_{\alpha_j},\cL_{-\alpha_j}]
\subseteq \fh.\]
Next we claim that $\alpha_i|_{\fh'}=\underline{0}$ for all $i \in I$.
Indeed, if $j \in I\setminus I'$ and $x\in [\cL_{\alpha_j},\cL_{-\alpha_j}]$,
then $\ad_\cL(x)=\underline{0}$ by Lemma~\ref{wsdlem4}. Consequently, we
have $\ad_\fg(x)=\underline{0}$ for all $x\in \fh'$. For $i\in I$ and
$x\in \fh'$, this yields that $0=\ad_\fg(x)(e_i)=[x,e_i]=\alpha_i(x)e_i$
and so $\alpha_i(x)=0$. Hence, the restriction of $\alpha_i$ to
$\fh'$ is zero for any $i\in I$, as claimed.

Assume, if possible, that $I'\subsetneqq I$. Then the restrictions of
the linear maps $\{\alpha_i\mid i\in I\}$ to the subspace $\langle h_j 
\mid j\in I'\rangle_\C\subsetneqq\fh$ are linearly dependent. So there
are scalars $c_i\in \C$, not all $0$, such that $\sum_{i\in I}c_i 
\alpha_i(h_j)=0$ for all $j\in I'$. Now, we have just seen that
$\alpha_i|_{\fh'}=\underline{0}$. Hence, $\sum_{i\in I}c_i\alpha_i(x)=0$
for all $x \in \fh$, contradiction to $\{\alpha_i\mid i \in I\}$ being
linearly independent. So we must have $I'=I$, which shows that $\fh=
\langle h_i\mid i \in I \rangle_\C$ and so $\dim \fh\leq |I|$. On the
other hand, since $\{\alpha_i \mid i \in I\}$ is linearly independent,
we have $\dim \fh=\dim \fh^* \geq |I|$. Hence, $\{h_i\mid i \in I\}$ is
a basis of $\fh$ and $\{\alpha_i\mid i \in I\}$ is a basis of~$\fh^*$.
\end{proof}

\begin{defn} \label{defTD2} Assume that the conditions in 
Definition~\ref{defTD} hold. Let $h_i\in \fh$ ($i\in I$) be as
in Proposition~\ref{wsdprop3}. Then 
\[A=\bigl(\alpha_j(h_i)\bigr)_{i,j\in I}\]
is called the \nm{structure matrix} of $\cL$ (with respect to $\Delta$).

Note that, since $\{h_i\mid i \in I\}$ is a basis of $\fh$ and $\{\alpha_i
\mid i \in I\}$ is a basis of $\fh^*$, we certainly have $\det(A)\neq 0$.
\end{defn}

\begin{exmp} \label{CKgo3} Let $\cL=\gom_3(Q_3,\C)$ as in 
Exercise~\ref{xcaclassic3}(b). The following three matrices form
a basis of $\cL$:
\[\renewcommand{\arraystretch}{0.9} 
h:=\left( \begin{array}{c@{\hspace{4pt}}r@{\hspace{4pt}}r}  1 & 0 & 0 \\
0 & 0 & 0\\ 0 & 0 & -1 \end{array}\right),\quad
e:=\left( \begin{array}{c@{\hspace{4pt}}r@{\hspace{4pt}}r}  0 & 1 & 0 \\
0 & 0 & -\delta \\ 0 & 0 & 0 \end{array}\right),\quad
f:=\left( \begin{array}{c@{\hspace{4pt}}r@{\hspace{4pt}}r}  0 & 0 & 0 \\
1 & 0 & 0\\ 0 & -1 & 0 \end{array}\right).\]
We have $[e,f]=h$, $[h,e]=e$ and $[h,f]=-f$. Let $\fh:=\langle h\rangle_\C$.
Then $\fh$ is an abelian subalgebra such that $\cL$ is $\fh$-diagonalisable
and $\cL_{\underline{0}}=\fh$. Define $\alpha_1\in \fh^*$ by $\alpha_1(h):=1$.
Then $\cL=\fh\oplus \cL_{\alpha_1}\oplus \cL_{-\alpha_1}$. Thus $(\cL,\fh)$ is
of Cartan--Killing type with respect to $\Delta=\{\alpha_1\}\subseteq \fh^*$.
But $h$ is not yet the ``correct'' element in $\fh$ according to 
Proposition~\ref{wsdprop3}. We can fix this as follows. Setting $e_1:=e$,
$f_1:=2f$ and $h_1:=2h$, we obtain the required relations $[e_1,f_1]=h_1$, 
$[h_1,e_1]=2e_1$, $[h_1,f_1]=-2f_1$ and $\alpha_1(h_1)=2$. The 
corresponding structure matrix is $A=(2)$. Since the Lie brackets among 
the basis elements $h_1,e_1,f_1$ are exactly the same as in $\slm_2(\C)$,
we now see that $\cL\cong \slm_2(\C)$.
\end{exmp}

\begin{xca} \label{xcagensl2} Show that, if $|I|=1$, then 
$\cL\cong\slm_2(\C)$.
\end{xca}

\begin{exmp} \label{cartsln} 
Let $\cL=\slm_n(\C)$ ($n\geq 2$) and $\fh\subseteq \cL$ be the abelian 
subalgebra of all diagonal matrices in~$\cL$; we have $\dim \fh=\dim \fh^*
=n-1$. By Remark~\ref{cartansubgln}, $\cL$ is $\fh$-diagonalisable 
and $C_\cL(\fh)=\fh$. Thus, (CK1) holds. For $1\leq i\leq n$, let 
$\varepsilon_i\in \fh^*$ be the map which sends a diagonal matrix to its 
$i$-th diagonal entry. (Note that, now, we have the linear relation 
$\varepsilon_1+ \ldots +\varepsilon_n= \underline{0}$.) For $i\neq j$ let 
$e_{ij} \in \cL$ be the matrix with entry~$1$ at position $(i,j)$, and 
$0$ everywhere else. Then we have again $\cL=\fh\oplus \bigoplus_{i\neq j} 
\C e_{ij}$. By the same computations as in Example~\ref{cartangln}, we see 
that $P_\fh(\cL)= \{\underline{0}\}\cup \Phi$, where 
\[\Phi:=\{ \pm (\alpha_i+\alpha_{i+1}+\ldots + \alpha_{j-1}) \mid 
1\leq i<j\leq n\},\] 
with $\alpha_i=\varepsilon_i-\varepsilon_{i+1}$ for $1\leq i \leq n-1$.
Thus, (CK2) holds, but we still need to check that $\{\alpha_1,\ldots, 
\alpha_{n-1} \}\subseteq \fh^*$ is linearly independent. If not, then
there would exist some $0\neq x\in \fh$ such that $\alpha_i(x)=0$
for $1\leq i \leq n-1$. But then $\alpha_1(x)=0$ and so the first two 
diagonal entries of $x$ are equal. Next, $\alpha_2(x)=0$ and so the second 
and third diagonal entries are equal. We conclude that all diagonal entries 
are equal and so $\trc(x)\neq 0$, contradiction. Hence, since $\dim \fh^*=
n-1$, the set $\{\alpha_1, \ldots,\alpha_{n-1}\}$ is a basis of $\fh^*$. 
Given the above description of $\Phi$, this now shows that $|\Phi|=n^2-n$,
and so a dimension argument as in Example~\ref{cartangln} yields that
$\cL_{\underline{0}}=\fh$ and $\dim \cL_\alpha=1$ for all $\alpha \in \Phi$.

Finally, we set $e_i:=e_{i,i+1}\in \cL_{\alpha_i}$ and $f_i:=e_{i+1,i}
\in \cL_{-\alpha_i}$ for $1\leq i \leq n-1$. Then $h_i:=[e_i,f_i]\in \fh$
is the diagonal matrix with entries $1,-1$ at positions $i,i+1$ (and $0$ 
otherwise). We see that $\{h_1, \ldots,h_{n-1}\}$ is a basis of $\fh$ and,
hence, that $\fh=\sum_{1\leq i \leq n-1} [\cL_{\alpha_i},\cL_{-\alpha_i}]$.
Thus, (CK3) also holds and so $(\cL,\fh)$ is of Cartan--Killing type with
respect to $\Delta=\{\alpha_1,\ldots,\alpha_{n-1}\}$. We compute that
\[ \renewcommand{\arraystretch}{0.9} A=\bigl(\alpha_j(h_i)\bigr)=
\left(\begin{array}{@{\hspace{4pt}}r@{\hspace{4pt}}r@{\hspace{4pt}}
r@{\hspace{4pt}}r@{\hspace{4pt}}r@{\hspace{4pt}}r} 2 & -1 &&&&\\ -1 & 2 & 
-1 &&&\\ & -1 & 2 & -1 & & \\ & & \ddots & \ddots & \ddots & \\ & & & -1 
& 2 & -1 \\ &&&& -1 & 2 \end{array}\right)\in M_{n-1}(\Z)\]
where all non-specified entries are $0$. Note that 
$h_i\in [\cL_{\alpha_i},\cL_{-\alpha_i}]$ and $\alpha_i(h_i)=2$. Hence, the 
above elements $\{h_1,\ldots,h_{n-1}\}$ are indeed the elements whose 
existence and uniqueness is proved in Proposition~\ref{wsdprop3}. We 
know that $\det(A)\neq 0$ but we leave it as an exercise to compute that 
$\det(A)=n$. 
\end{exmp}

Assume from now on that $(\cL,\fh)$ is of Cartan--Killing type with respect
to $\Delta=\{\alpha_i\mid i\in I\}$, as in Definition~\ref{defTD}.

\begin{lem} \label{wsdprop4} Let $\alpha\in\Phi^+$ and $i\in I$. 
If $\alpha+m\alpha_i\in\Phi$ for some $m\in\Z$, then $\alpha=\alpha_i$ 
or $\alpha+m\alpha_i\in \Phi^+$.
\end{lem}

\begin{proof} Write $\alpha=\sum_{j \in I} n_j \alpha_j$ where $n_j\in
\Z_{\geq 0}$ for all~$j$. Assume that $\alpha\neq \alpha_i$; since $\alpha
\in\Phi^+$, we also have $\alpha\neq -\alpha_i$. By 
Proposition~\ref{wsdprop3}, the restriction of $\alpha_i$ to
$[\cL_{\alpha_i},\cL_{-\alpha_i}]$ is non-zero and so Lemma~\ref{wsdlem3} 
implies that $\alpha\not\in \Z \alpha_i$.
Hence, we must have $n_{i_0}>0$ for some $i_0\neq i$. But then $n_{i_0}>0$ 
is also the coefficient of $\alpha_{i_0}$ in $\alpha+m\alpha_i$. Since
every root is either in $\Phi^+$ or in $\Phi^-$, we conclude that 
$\alpha+ m\alpha_i\in \Phi^+$.
\end{proof}

\begin{rem} \label{astring0} Let $i\in I$ and $h_i\in [\cL_{\alpha_i},
\cL_{-\alpha_i}]$ be as in Proposition~\ref{wsdprop3}. Let $e_i\in 
\cL_{\alpha_i}$ and $f_i\in \cL_{-\alpha_i}$ be such that $h_i=[e_i,f_i]$. 
Since $\dim \cL_{\pm \alpha_i}=1$, we have $\cL_{\alpha_i}=\langle e_i
\rangle_\C$ and $\cL_{-\alpha_i}=\langle f_i \rangle_\C$. Furthermore, 
since $\alpha_i(h_i)=2$, we have $[h_i,e_i]=2e_i$ and $[h_i,f_i]=-2f_i$. 
Thus, 
\[ \fs_i:=\langle e_i,h_i,f_i \rangle_\C\subseteq \cL\]
is a $3$-dimensional subalgebra isomorphic to $\slm_2(\C)$. We call 
$\{e_i,h_i,f_i\}$ an \nms{$\slm_2$-triple}{sl$_2$-triple} in $\cL$.
This will provide a powerful tool in the study of $\cL$. The elements 
$\{e_i,f_i\mid i \in I\}$ are called \nm{Chevalley generators} of~$\cL$.
Note that the $f_i$ are determined once the $e_i$ are chosen (via the
relations $h_i=[e_i,f_i]$); the $e_i$ are only unique up to non-zero 
scalar multiples. We also have the following relations for all $i,j\in I$ 
such that $i\neq j$:
\[ [h_i,h_j]=0, \;\; [h_i,e_j]=a_{ij}e_j, \;\; [h_i,f_j]=-a_{ij}f_j,\;\;
[e_i,f_j]=0.\]
The first relation holds since $\fh$ is abelian; the second and third
relations hold since $e_j\in \fg_{\alpha_j}$, $f_j\in \fg_{-\alpha_j}$ and
by the definition of $A$. Finally, we have $[e_i,f_j]\in [\cL_{\alpha_i},
\cL_{-\alpha_j}] \subseteq \cL_{\alpha_i-\alpha_j}$ by 
Proposition~\ref{wsdprop2}. But, for $i\neq j$, we have 
$\alpha_i-\alpha_j\not\in P_\fh(\fg)$ by (CK2) and so $[e_i,f_j]=0$.
\end{rem}

\begin{rem} \label{astring} In the proof of Lemma~\ref{wsdlem3}, we used 
the results on representations of $\slm_2(\C)$ that we obtained in 
Section~\ref{sec04a}. We can now push this argument much further. So let 
us fix $i\in I$ and let $\{e_i,h_i,f_i\}$ be a corresponding 
$\slm_2$-triple, as above. Then $\slm_2(\C)\cong \fs_i:=\langle e_i,h_i,
f_i \rangle_\C \subseteq \cL$. Let us also fix $\beta\in\Phi$ such that
$\beta \neq \pm\alpha_i$. Since $\Phi$ is finite, there are well-defined 
integers $p,q\geq 0$ such that 
\[\beta-q\alpha_i,\quad\ldots,\quad\beta-\alpha_i,\quad\beta,\quad
\beta+\alpha_i,\quad\ldots,\quad \beta+p\alpha_i\]
are all contained in $\Phi$, but $\beta+(p+1)\alpha_i\not\in\Phi$
and $\beta-(q+1)\alpha_i\not\in \Phi$. (It could be that $p=0$ or $q=0$.) 
The above sequence of roots is called the \nms{$\alpha_i$-string through 
$\beta$}{alphai-string through beta}. Now consider the subspace
\[ M:=\cL_{\beta-q\alpha_i}\oplus \ldots\oplus \cL_{\beta-\alpha_i}\oplus
\cL_{\beta}\oplus \cL_{\beta+\alpha_i}\oplus\ldots\oplus \cL_{\beta+p\alpha_i}
\subseteq \cL.\]
We claim that $M$ is an $\fs_i$-submodule of $\cL$. Now, we certainly have
$[\fh,M]\subseteq M$ and so $M$ is invariant under~$h_i$. By 
Proposition~\ref{wsdprop2}, we have $[\cL_{\pm \alpha_i},\cL_{\beta+
n\alpha_i}] \subseteq \cL_{\beta+(n\pm 1)\alpha_i}$ for all $n\in \Z$. This 
shows that all subspaces $\cL_{\beta+n\alpha_i}$ with $-q<n<p$ are invariant 
under $e_i$ and~$f_i$. Finally, by Lemma~\ref{wsdlem3} (applied to $\lambda=
\alpha_i$), we have $\beta\neq n\alpha_i$ for all $n\in \Z$. Hence, 
$\underline{0}\neq \beta+(p+1)\alpha_i\not\in \Phi$ and so $[\cL_{\alpha_i},
\cL_{\beta+p\alpha_i}]\subseteq \cL_{\beta+(p+1)\alpha_i}=\{0\}$. Similarly, 
we have $[\cL_{-\alpha_i},\cL_{\beta-q\alpha_i}]\subseteq \cL_{\beta-(q+1)
\alpha_i}=\{0\}$. Thus, $M$ is an $\fs_i$-submodule of $\cL$, as claimed. 
Now recall that the module action is given by $\ad_\cL\colon \cL\rightarrow 
\gl(\cL)$. Since $\cL$ is $\fh$-diagonalisable, the eigenvalues of $x\in\fh$ 
are given by $\lambda(x)$ for $\lambda \in P_\fh(\cL)$ (each with 
multiplicity $\dim \cL_\lambda\geq 1$). So the eigenvalues of $h_i$ on $M$ 
are given by $(\beta+n\alpha_i)(h_i)$ for $-q\leq n\leq p$, each with 
multiplicity $\dim \cL_{\beta+n\alpha_i}\geq 1$. Explicitly, the list of 
eigenvalues (not counting multiplicities) is 
\[ \beta(h_i)-2q,\;\ldots,\;\beta(h_i)-2,\;\beta(h_i),\;\beta(h_i)+2,
\;\ldots,\;\beta(h_i)+2p.\]
By Proposition~\ref{sl2modd}, all eigenvalues of $h_i$ are integers,
and if $m\in\Z$ is an eigenvalue, then so is $-m$. In particular, the
largest eigenvalue is the negative of the smallest eigenvalue. First of
all, this implies that $\beta(h_i)+2p=-(\beta(h_i)-2q)$ and so 
\begin{equation*}
\beta(h_i)=q-p\in\Z.\tag{a}
\end{equation*}
Furthermore, $-q\leq p-q=-\beta(h_i)\leq p$. Thus, we conclude that 
\begin{equation*}
\beta-\beta(h_i)\alpha_i \in\Phi\; \mbox{ belongs to the
$\alpha_i$-string through $\beta$}. \tag{b}
\end{equation*}
We can go even one step further. Let $0\neq v^+\in \cL_{\beta+p\alpha_i}$
be fixed. Then $h_i.v^+=cv^+$ where $c=\beta(h_i)+2p=(q-p)+2p=p+q$. Since
$[e_i,v^+]\in \cL_{\beta+(p+1)\alpha_i}=\{0\}$, we have $e_i.v^+=\{0\}$ and 
so $v^+\in M$ is a \nm{primitive vector}, as in Remark~\ref{sl2highb}.
Correspondingly, we have a subspace $E:=\langle v_n\mid n\geq 0
\rangle_\C\subseteq M$, where 
\[ v_0:=v^+\qquad \mbox{and}\qquad \textstyle v_{n+1}:=
\frac{1}{n+1}[f_i,v_n] \quad \mbox{for all $n \geq 0$}.\]
(We also set $v_{-1}:=0$.) As shown in Remark~\ref{sl2highb}, we have
\[\dim E=c+1=p+q+1\quad \mbox{and}\quad E=\langle v_0,v_1,\ldots,
v_{p+q}\rangle_\C.\]
In particular, $v_0,v_1,\ldots,v_{p+q}$ are all non-zero. We can exploit 
this as follows. First, $v_0=v^+\in 
\cL_{\beta+p\alpha_i}$. Hence, if $p\geq 1$, then $v_1=[f_i,v_0]\in 
[\cL_{-\alpha_i},\cL_{\beta+p\alpha_i}]\subseteq \cL_{\beta+(p-1)\alpha_i}$; 
furthermore, if $p\geq 2$, then $v_2=\frac{1}{2}[f_i,v_1] \in 
[\cL_{-\alpha_i},\cL_{\beta+(p-1)\alpha_i}] \subseteq \cL_{\beta+(p-2)
\alpha_i}$. Going on in this way, we find that $0\neq v_p \in \cL_{\beta}$. 
Since $[e_i,v_p]=(c-p+1)=(q+1) v_{p-1}$ (see Remark~\ref{sl2highb}), we 
conclude that 
\begin{equation*}
\begin{array}{l@{\hspace{3pt}}c@{\hspace{3pt}}l}
\left[f_i,\left[e_i,v_p\right]\right]&=&(q+1)\left[f_i,v_{p-1}
\right]=p(q+1)v_p,\\
\left[e_i,\left[f_i,v_p\right]\right]&=&(p+1)\left[e_i,v_{p+1}
\right]=q(p+1)v_p.  \end{array} \qquad\qquad\tag{c}
\end{equation*}
In particular, since $0\neq v_p\in \cL_\beta$, this implies that 
\begin{equation*}
\begin{array}{ll}
\{0\}\neq \left[\cL_{\alpha_i},\cL_\beta\right] \subseteq \cL_{\beta+\alpha_i}
&\quad \mbox{if $p>0$, that is, $\beta+\alpha_i\in \Phi$},\\ 
\{0\}\neq \left[\cL_{-\alpha_i},\cL_\beta\right]\subseteq \cL_{\beta-\alpha_i}
&\quad \mbox{if $q>0$, that is, $\beta-\alpha_i\in\Phi$}.
\end{array} \tag{c'}
\end{equation*}
These relations will be very helpful for inductive arguments (see, e.g.,
Proposition~\ref{genlie} or Theorem~\ref{canbas} below).
\end{rem}

\begin{rem} \label{poshi} For future reference, we note that
$\beta(h_i)\in\Z$ for all $\beta\in\Phi$ and all $i\in I$. Indeed,
if $\beta\neq \pm \alpha_i$, then this holds by Remark~\ref{astring}(a).
But if $\beta=\pm\alpha_i$, then $\beta(h_i)=\pm \alpha_i(h_i)=\pm 2$.
\end{rem}

\begin{cor} \label{gencart0} Consider the matrix
$A=(a_{ij})_{i,j\in I}$ in Definition~\ref{defTD2}, where
$a_{ij}=\alpha_j(h_i)$ for $i,j\in I$. Then the following hold.
\begin{itemize}
\item[{\rm (a)}] $a_{ij}\in \Z$ and $a_{ii}=2$ for all $i,j\in I$.
\item[{\rm (b)}] $a_{ij}\leq 0$ for all $i,j\in I$, $i\neq j$.
\item[{\rm (c)}] $a_{ij}\neq 0\Leftrightarrow a_{ji}\neq 0$ for all
$i,j\in I$.
\end{itemize}
\end{cor}

\begin{proof} (a) See Proposition~\ref{wsdprop3} and 
Remark~\ref{poshi}.

(b) Assume, if possible, that $a_{ij}>0$. Then, by Remark~\ref{astring}(b),
we have $\alpha_j-n\alpha_i\in\Phi$, where $n=\alpha_j(h_i)>0$, contradiction
to (CK2).

(c) This is clear for $i=j$. Now assume that $i\neq j$ and $a_{ji}\neq 0$; 
then $a_{ji}<0$ by (b). By Remark~\ref{astring}(b), we have $\alpha_i+n
\alpha_j\in\Phi$, where $n=-\alpha_i(h_j)=-a_{ji}>0$; furthermore, 
$\alpha_i+n\alpha_j$ belongs to the $\alpha_j$-string through $\alpha_i$. 
Hence, since $n>0$, we also have $\alpha_i+\alpha_j \in\Phi$. Now we 
reverse the roles of $\alpha_i$ and $\alpha_j$ and consider the 
$\alpha_i$-string through $\alpha_j$. Let $p,q \geq 0$ in 
Remark~\ref{astring} be defined with respect to $\alpha_i$ and $\alpha
:=\alpha_j$. Since $\alpha_j+\alpha_i \in\Phi$, we have $p\geq 1$. By 
(CK2), we have $\alpha_j-\alpha_i\not\in \Phi$ and so $q=0$. Hence, 
Remark~\ref{astring}(a) shows that $a_{ij}= \alpha_j(h_i)=-p<0$.
\end{proof}

\begin{xca} \label{xcastring} In the setting of Remark~\ref{astring},
show that 
\begin{align*}
p&=\max\{n\geq 0\mid \beta+n\alpha_i\in\Phi\},\\
q&=\max\{n\geq 0\mid \beta-n\alpha_i\in\Phi\}.
\end{align*}
Deduce that, if $\beta\pm n\alpha_i\in \Phi$ for some $n>0$, then 
$\beta\pm \alpha_i\in\Phi$.
\end{xca}

\begin{xca} \label{xcastringm} Let $i,j\in I$ be such that $i\neq j$.
Show that $a_{ij}\neq 0$ if and only if $\alpha_i+\alpha_j\in \Phi$.
\end{xca}

\begin{xca} \label{xcadirck} Let $\cL_1$ and $\cL_2$ be finite-dimensional
Lie algebras. Let $\fh_1\subseteq \cL_1$ and $\fh_2\subseteq \cL_2$ be
abelian subalgebras such that $\cL_1$ is of Cartan--Killing  type
with respect to $\Delta_1=\{\alpha_i\mid i\in I_1\}\subseteq \fh_1^*$ and
$\cL_2$ is of Cartan--Killing  type with respect to $\Delta_2=\{\beta_j
\mid j\in I_2\}\subseteq \fh_2^*$. 
Now consider the direct product $\cL:=\cL_1\times \cL_2$ (see 
Example~\ref{dirprodalg}). Then $\fh:=\fh_1\times\fh_2 \subseteq \cL$ is 
an abelian subalgebra. For $\lambda\in \fh_1^*$ we define 
$\dot{\lambda}\in \fh^*$ by $\dot{\lambda}(h_1,h_2):=\lambda(h_1)$ 
for $(h_1,h_2)\in \fh$; similarly, for $\mu\in \fh_2^*$ 
we define $\dot{\mu}\in \fh^*$ by $\dot{\mu}(h_1,h_2):=\mu(h_2)$ for
$(h_1,h_2)\in \fh$. 
\begin{itemize}
\item[(a)] Show that $\cL$ is $\fh$-diagonalisable and that 
$P_\fh(\cL)=\{\dot{\alpha}\mid \alpha\in P_\fh(\cL_1)\} \cup\{\dot{\beta}
\mid \beta \in P_\fh(\cL_2)\}$.
\item[(b)] Show that $\cL$ is of Cartan--Killing type with respect to 
$\Delta:=\{\dot{\alpha}_i\mid i\in I_1\}\cup \{\dot{\beta}_j\mid j\in 
I_2\}\subseteq \fh^*$.
\item[(c)] Let $A_1$ and $A_2$ be the structure matrices of $(\cL_1,\fh_1)$ 
and $(\cL_2,\fh_2)$, respectively. Show that the structure matrix 
of $(\cL,\fh)$ is block diagonal with diagonal blocks $A_1$ and~$A_2$. 
\end{itemize}
\end{xca}

\section{The Weyl group} \label{sec1a3}

We keep the basic setting of the previous section, where $(\cL,\fh)$ is of 
Cartan--Killing type with respect to $\Delta=\{\alpha_i\mid i \in I\}
\subseteq \fh^*$. The formula in Remark~\ref{astring}(b) suggests the 
following definition.

\begin{defn} \label{weyl1} For $i \in I$, let $h_i\in [\cL_{\alpha_i},
\cL_{-\alpha_i}]$ be as in Propositon~\ref{wsdprop3}. We define a linear map
$s_i\colon \fh^*\rightarrow \fh^*$ by 
\[ s_i(\lambda):=\lambda-\lambda(h_i)\alpha_i \qquad\mbox{for
$\lambda \in \fh^*$}.\]
Note that $s_i(\alpha_i)=\alpha_i-2\alpha_i=-\alpha_i$ and 
$s_i(\lambda)=\lambda$ for all $\lambda \in \fh^*$ with $\lambda(h_i)=0$. 
Since $\fh^*=\langle \alpha_i \rangle_\C \oplus \{\lambda\in \fh^*\mid
\lambda(h_i)=0\}$,
we conclude that $s_i$ is diagonalisable, with one eigenvalue equal to~$-1$ 
and $|I|-1$ eigenvalues equal to~$1$. In particular, $s_i^2=\id_{\fh^*}$,
$\det(s_i)=-1$ and $s_i \in \GL(\fh^*)$. The subgroup 
\[W:=\langle s_i \mid i \in I\rangle \subseteq \GL(\fh^*)\]
is called the \nm{Weyl group} of $\cL$ (with respect to $\Delta$).
Note that, since $s_i^{-1}=s_i$ for all $i\in I$, every element $w\in W$
can be written as a product $w=s_{i_1}\cdots s_{i_r}$ where $r\geq 0$ and
$i_1,\ldots,i_r\in I$. (Such an expression for $w$ is by no means unique; 
we have $w=\id$ if $r=0$.) 
\end{defn}

\begin{rem} \label{weyl0} By Remark~\ref{astring}, we have $s_i(\alpha)
\in\Phi$ for all $\alpha\in\Phi$ with $\alpha\neq\pm \alpha_i$. But we 
also have $s_i(\alpha_i)=-\alpha_i$ and so $s_i(\Phi)=\Phi$. Consequently, 
we have $w(\Phi)=\Phi$ for all $w\in W$. So we have an action of the 
group $W$ on the finite set $\Phi$ via 
\[W \times \Phi \rightarrow \Phi, \qquad (w,\alpha)\mapsto w(\alpha).\]
Let $\mbox{Sym}(\Phi)$ denote the symmetric group on $\Phi$. Then we obtain 
a group homomorphism $\pi \colon W \rightarrow \mbox{Sym}(\Phi)$, $w 
\mapsto \pi_w$, where $\pi_w(\alpha):=w(\alpha)$ for all $\alpha\in \Phi$. 
If $\pi_w=\id_\Phi$, then $w(\alpha)=\alpha$ for all $\alpha\in \Phi$.
In particular, $w(\alpha_i)=\alpha_i$ for all $i\in I$. Since $\{\alpha_i
\mid i\in I\}$ is a basis of~$\fh^*$, it follows that $w=\id_{\fh^*}$. Thus, 
$\pi$ is injective and $W$ is isomorphic to a subgroup of 
$\mbox{Sym}(\Phi)$; in particular, $W$ is a finite group.
\end{rem}

In order to prove the ``Key Lemma'' below, we shall use a construction that 
essentially relies on the fact that $W$ is a finite group. For this purpose,
let $E:=\langle \alpha_i\mid i\in I\rangle_\R \subseteq \fh^*$. Then $E$ is 
an $\R$-vector space, and $\{\alpha_i\mid i\in I\}$ still is a basis of 
$E$. By (CK2), we have $\Phi\subseteq E$. Since $\alpha(h_i)\in\Z$ for all 
$\alpha\in \Phi$ and $i\in I$ (see Remark~\ref{poshi}), we also have 
$s_i(E)\subseteq E$ for all $i\in I$ and so $w(E)\subseteq E$ for all 
$w\in W$. Thus, we may regard $W$ as a subgroup of $\GL(E)$ (but we will 
not introduce a separate notation for this). Let $\langle\;,\;\rangle_0
\colon E\times E \rightarrow \R$ be the standard scalar product for which 
$\{\alpha_i\mid i\in I\}$ is an orthonormal basis. Thus, for $v,v'\in E$
we have $\langle v,v'\rangle_0=\sum_{i,j\in I} x_ix_j'$ where
$v=\sum_{i\in I}x_i\alpha_i$ and $v'=\sum_{j\in I} x_j'\alpha_j$,
with $x_i,x_j'\in\R$ for all $i,j\in I$. Then we define a new map 
$\langle\;,\;\rangle \colon E\times E \rightarrow\R$ by
\[ \langle v,v'\rangle:=\sum_{w\in W} \langle w(v),w(v')\rangle_0 \qquad
\mbox{for $v,v'\in E$}.\]  
Since $E\rightarrow E$, $v\mapsto w(v)$, is linear for each $w\in W$,
it is clear that $\langle \;,\;\rangle$ is a symmetric bilinear form.
For $v\in E$, we have 
\[\langle v,v\rangle=\sum_{w\in W} \underbrace{\langle w(v), w(v)
\rangle_0}_{\geq 0} \geq 0.\]
If $\langle v,v\rangle=0$, then $\langle w(v),w(v)\rangle_0=0$ for all 
$w\in W$. In particular, this holds for $w=\id_E$ and so $\langle v,v
\rangle_0=0$. But $\langle \;,\;\rangle_0$ is positive-definite and so $v=0$.
Thus, $\langle\;,\;\rangle$ is also positive-definite. Finally, taking the
sum over all $w\in W$ implies the following invariance property:
\[ \langle w(v),w(v')\rangle=\langle v,v'\rangle\qquad \mbox{for all
$w\in W$ and $v,v'\in E$}.\]
Indeed, for a fixed $w\in W$, we have 
\[\langle w(v),w(v')\rangle=\sum_{w'\in W}\langle w'w(v),w'w(v')\rangle_0.\]
Now, since $W$ is a group, the map $W\rightarrow W$, $w'\mapsto w'w$,
is a bijection. Hence, up to reordering terms, the sum on the right hand 
side is the same as the sum in the definition of $\langle v,v'\rangle$.

\begin{rem} \label{keylem0} Let $i\in I$ and $\lambda \in E$; recall
that $E=\langle\alpha_i\mid i\in I\rangle_\R\subseteq \fh^*$. Using 
the relation $s_i(\alpha_i)=-\alpha_i$, the defining formula for
$s_i(\lambda)$, and the above invariance property, we obtain the
following identities:
\begin{align*}
-\langle \alpha_i,\lambda\rangle&=\langle s_i(\alpha_i),\lambda\rangle=
\langle s_i^2(\alpha_i),s_i(\lambda)\rangle=\langle \alpha_i,s_i(\lambda)
\rangle\\ & =\langle \alpha_i,\lambda-\lambda(h_i)\alpha_i\rangle
=\langle \alpha_i,\lambda\rangle -\lambda(h_i)\langle \alpha_i,
\alpha_i\rangle.
\end{align*}
Since $\langle\alpha_i,\alpha_i\rangle\in\R_{>0}$, this yields the fomula
\[ \lambda(h_i)=2\frac{\langle \alpha_i,\lambda\rangle}{\langle \alpha_i,
\alpha_i\rangle}\in \R\qquad\mbox{for all $\lambda\in E$ and $i\in I$}.\]
This formula shows that each $s_i\colon E\rightarrow E$ is an 
\nm{orthogonal reflection} with root $\alpha_i$ (and with respect to
$\langle\;,\;\rangle$). 
\end{rem}

\begin{lem}[Key Lemma] \label{keylem} \nmi{}{Key Lemma} Let $\alpha\in 
\Phi^+$ but $\alpha\not\in\Delta$. Write $\alpha=\sum_{i \in I} n_i
\alpha_i$ where $n_i\in \Z_{\geq 0}$ for all~$i$. Then there exists some 
$i\in I$ such that $n_i>0$ and $\alpha(h_i)\in\Z_{>0}$. Furthermore, we 
have $s_i(\alpha)=\alpha-\alpha(h_i) \alpha_i \in\Phi^+$ and 
$\alpha-\alpha_i \in \Phi^+$.  
\end{lem}

\begin{proof} Since $\underline{0}\neq \alpha\in E$, the above discussion 
shows that 
\[ \sum_{i \in I} n_i\underbrace{\langle \alpha_i,\alpha\rangle}_{\in \R}=
\langle \alpha,\alpha\rangle>0.\]
Since $n_i\geq 0$ for all $i$, there must be some $i\in I$ such 
that $n_i>0$ and $\langle \alpha_i,\alpha\rangle >0$. Furthermore, since 
$\langle \alpha_i,\alpha_i\rangle >0$, the formula in Remark~\ref{keylem0} 
shows that we also have $\alpha(h_i)>0$. By Remark~\ref{poshi}, $\alpha(h_i)
\in \Z$ and so $\alpha(h_i)\in \Z_{>0}$, as desired. Now, since $\alpha\in 
\Phi^+ \setminus \Delta$, we have $\alpha\neq \pm \alpha_i$. Hence, 
Remark~\ref{astring}(b) shows that $\alpha-\alpha(h_i)\alpha_i\in \Phi$ 
belongs to the $\alpha_i$-string through $\alpha$. Since $\alpha(h_i)\in 
\Z_{>0}$, we conclude that $\alpha-\alpha_i$ also belongs to that 
$\alpha_i$-string and so $\alpha-\alpha_i\in\Phi$. It remains to show that 
$\alpha-\alpha_i \in \Phi^+$ and $\alpha-\alpha(h_i)\alpha_i\in \Phi^+$. 
But this follows from Lemma~\ref{wsdprop4}, since $\alpha\neq \alpha_i$. 
\end{proof} 

\begin{rem} \label{defbasephi1} Recall from Remark~\ref{defheight} the
definition of the height function $\hgt\colon P_\fh(\fg)\rightarrow \Z$. 
Since $\{\alpha_i\mid i\in I\}$ is a basis of $\fh^*$, we can actually 
extend it linearly to a function $\hgt\colon \fh^* \rightarrow \C$.
The ``Key Lemma'' often allows us to argue by induction on the height of
roots; here is a first example. 

Let $\alpha \in \Phi^+$ and 
$n=\hgt(\alpha) \geq 1$. Claim: We can write $\alpha=\alpha_{i_1}+
\ldots+\alpha_{i_n}$ where $i_j \in I$ for all~$j$ and, for each 
$j\in\{1,\ldots,n\}$, we also have $\alpha_{i_j}+ \ldots +\alpha_{i_n}
\in \Phi^+$. 

We argue by induction on $n:=\hgt(\alpha)\geq 1$. If $n=1$, then 
$\alpha=\alpha_i$ for some $i\in I$ and there is nothing to prove. Now 
let $n\geq 2$. Then $\alpha\not\in \Delta$ and so, by Lemma~\ref{keylem}, 
we have $\beta:=\alpha-\alpha_{i_1}\in \Phi^+$ for some $i_1\in I$. Now 
$\hgt(\beta)=n-1$. By induction, there exist $i_2,\ldots,i_n\in I$ 
such that the required conditions hold for $\beta$. But then $\alpha=
\alpha_{i_1}+\alpha_{i_2}+\ldots+\alpha_{i_n}$ and the required 
conditions hold for~$\alpha$.
\end{rem}

\begin{thm} \label{mainthm1} Recall that $(\cL,\fh)$ is of Cartan--Killing 
type with respect to $\Delta=\{\alpha_i\mid i \in I\}$. Then the 
following hold.
\begin{itemize}
\item[{\rm (a)}] $\Phi=\{w(\alpha_i) \mid w\in W,i\in I\}$ and 
$\Phi^-=-\Phi^+$. 
\item[{\rm (b)}] If $\alpha\in\Phi$ and $0\neq c\in\C$ are such that 
$c\alpha\in\Phi$, then $c\in \{\pm 1\}$.
\end{itemize}
\end{thm}

\begin{proof} (a) Let $\Phi_0:=\{w(\alpha_i) \mid w\in W,i\in I\}$. By
Remark~\ref{weyl0}, $\Phi_0\subseteq \Phi$. Next, let $\alpha\in
\Phi^+$. We show by induction on $n:=\hgt(\alpha)\geq 1$ that $\alpha\in
\Phi_0$. If $n=1$, then $\alpha=\alpha_i$ for some $i\in I$ and so 
$\alpha=\id(\alpha_i)\in\Phi_0$. Now let $n\geq 2$. By 
Lemma~\ref{keylem}, there is some $j\in I$ such that $\alpha(h_j)
\in \Z_{>0}$ and $\beta:=s_j(\alpha)=\alpha-\alpha(h_j)\alpha_j\in\Phi^+$. 
We have $\hgt(\beta)=n-\alpha(h_j)<n$. By induction, $\beta\in\Phi_0$ and
so $\beta=w'(\alpha_i)$ for some $w'\in W$ and $i\in I$. But then 
$\alpha=s_j^2(\alpha)=s_j\bigl(s_j(\alpha)\bigr)=s_j(\beta)=s_j
w'(\alpha_i)\in\Phi_0$, as required. Thus, we have shown that 
$\Phi^+\subseteq\Phi_0$. 

Next, let $\alpha\in\Phi^+$. Since $\alpha\in\Phi_0$, we can write $\alpha=
w(\alpha_i)$, where $w\in W$ and $i\in I$, as above. Since $s_i(\alpha_i)=
-\alpha_i$, we obtain $-\alpha=w(-\alpha_i)=ws_i(\alpha_i)\in\Phi_0
\subseteq \Phi$. Furthermore, since $\alpha\in\Phi^+$, we have 
$-\alpha\in\Phi^-$. Thus, we have shown that $-\Phi^+\subseteq \Phi^-
\cap \Phi_0$. 

Now, there is a symmetry in Definition~\ref{defTD}. If we set $\alpha_i':=
-\alpha_i$ for all $i\in I$, then $(\cL,\fh)$ also is of Cartan--Killing 
type with respect to $\Delta':=\{\alpha_i'\mid i \in I\}$. (We have 
$-\alpha_i\in \Phi$ for all $i\in I$ by Proposition~\ref{wsdprop3}.) Then,
clearly, $\Phi^-$ is the corresponding set of positive roots and 
$\Phi^+$ is the set of negative roots. Now, the previous argument 
applied to $\Delta'$ instead of $\Delta$ shows that $-\Phi^-\subseteq 
\Phi^+$ and, hence, $|\Phi^-|\leq |\Phi^+|$. Since we also have $-\Phi^+
\subseteq \Phi^-\cap \Phi_0 \subseteq \Phi^-$, it now follows that 
$|\Phi^+|=|\Phi^-|$ and $\Phi^-=-\Phi^+\subseteq \Phi_0$. Hence, 
$\Phi=\Phi^+\cup \Phi^- \subseteq \Phi_0$ and so $\Phi=\Phi_0$.

(b) Assume that $\alpha\in\Phi$ and $c\alpha\in\Phi$, where $0\neq c\in\C$.
By (a) we can write $\alpha=w(\alpha_i)$ for some $w\in W$ and $i\in I$.
Then $c\alpha_i=cw^{-1}(\alpha)=w^{-1}(c\alpha)\in\Phi$ and so $c\alpha_i
(h_i)\in\Z$ by Remark~\ref{poshi}. But $\alpha_i(h_i)=2$ and so $2c\in\Z$; 
thus, $c\alpha_i\in\Phi$, where $c=n/2$ with $n\in\Z$. On the other hand, 
we can run the same argument with $\beta:=c\alpha\in \Phi$ and 
$c^{-1}\beta=\alpha\in\Phi$. So we also obtain that $c^{-1}\alpha_j
\in\Phi$ for some $j\in I$, where $c^{-1}=m/2$ for some $m\in\Z$. Thus, 
we have $nm=4$. If $m=\pm 1$, then $n=\pm 4$ and so $c=\pm 2$; hence,
$\pm 2\alpha_i\in\Phi$, contradiction to Lemma~\ref{wsdlem3} (applied 
to $\lambda=\alpha_i$). Similarly, if $n=\pm 1$, then $m=\pm 4$ and so 
$c^{-1}=\pm 2$; hence, $\pm 2\alpha_j\in \Phi$, contradiction to
Lemma~\ref{wsdlem3} (applied to $\lambda=\alpha_j$). Thus, we must 
have $n=\pm 2$ and so $c=\pm 1$. 
\end{proof}

We would like to make it completely explicit that $W$ and $\Phi$ are
determined by the single knowledge of the structure matrix~$A$ of~$\cL$.

\begin{rem} \label{explicit} Recall that $A=(a_{ij})_{i,j\in I}$, where
$a_{ij}=\alpha_j(h_i)\in \Z$ for all $i,j\in I$. Thus, the defining
equation of $s_i$ yields that 
\[ s_i(\alpha_j)=\alpha_j-a_{ij}\alpha_i \quad \mbox{for all $i,j\in I$}.\]
Hence, if $\lambda\in \fh^*$ and $\lambda=\sum_{i \in I} \lambda_i 
\alpha_i$ with $\lambda_i\in \C$, then we have 
\begin{equation*}
s_i(\lambda)=\sum_{j \in I} \lambda_j\bigl(\alpha_j-a_{ij}\alpha_i\bigr)=
\lambda-\Bigl(\sum_{j \in I} a_{ij}\lambda_j\Bigr)\alpha_i.\tag{$\clubsuit$}
\end{equation*}
This shows that the action of $s_i$ on $\fh^*$ is completely determined 
by~$A$. For each $w\in W$, let $M_w\in \GL_I(\C)$ be the matrix of
$w$ with respect to the basis $\{\alpha_i\mid i \in I\}$ of $\fh^*$. 
We have $w=s_{i_1} \cdots s_{i_l}$ for some $i_1,\ldots,i_l\in I$ and,
hence, also $M_w=M_{s_{i_1}}\cdot \ldots \cdot M_{s_{i_l}}$. The above
formulae show that each $M_{s_i}$ is completely determined by $A$, and has 
entries in $\Z$. Hence, the set of matrices $\{M_w\mid w\in W\} \subseteq 
\GL_I(\Z)$ is also determined by~$A$. Finally, by Theorem~\ref{mainthm1}(a), 
every $\alpha\in \Phi$ can be written as $\alpha:=w(\alpha_i)$ where 
$w\in W$ and $i\in I$. Then $\alpha=\sum_{i \in I} n_i\alpha_i$ where 
$(n_i)_{i \in I}\in \Z^I$ is the $i$-th column of $M_w$. Thus, 
\[ \cC(A):=\Bigl\{(n_i)_{i \in I} \in \Z^I\,\big|\, 
\sum_{i \in I} n_i \alpha_i\in \Phi\Bigr\}\subseteq \Z^I\]
is completely determined by $A$. More concretely, every $\alpha\in\Phi$ is 
obtained by repeatedly applying the generators $s_j$ of $W$ to the 
various~$\alpha_i$, using formula ($\clubsuit$). If, in the process, we
avoid the relation $s_i(\alpha_i)=-\alpha_i$, then we just obtain the set 
\[ \cC^+(A):=\Bigl\{(n_i)_{i \in I} \in \Z_{\geq 0}^I\,\big|\, 
\sum_{i \in I} n_i \alpha_i\in \Phi^+\Bigr\}\subseteq \Z^I.\]
(See the proof of Theorem~\ref{mainthm1}.) Here are a few examples. 
\end{rem}

\begin{exmp} \label{cartanA2} Let $\cL=\slm_3(\C)$, where $\Delta=\{\alpha_1,
\alpha_2\}$ and 
\[\renewcommand{\arraystretch}{0.8}
A=\left(\begin{array}{rr} 2 & -1 \\ -1 & 2\end{array}\right); \qquad
\mbox{see Example~\ref{cartsln}}.\]
The matrices of $s_1,s_2\in W$ with respect to the basis $\Delta$ are 
given by:
\[\renewcommand{\arraystretch}{0.8}
s_1:\left(\begin{array}{rr} -1 & 1 \\ 0 & 1 \end{array}\right),\qquad 
s_2:\left(\begin{array}{rr} 1 & 0 \\ 1 & -1 \end{array}\right);\]
see ($\clubsuit$). A direct computation shows that the product $s_1s_2
\in W$ has order $3$ and so $W\cong \fS_3$. Applying $s_1,s_2$ repeatedly 
to $\alpha_1,\alpha_2$ (avoiding $s_i(\alpha_i)=-\alpha_i$ for $i=1,2$), 
we obtain that
\[\cC^+(A)=\{(1,0),(0,1),(1,1)\} \qquad \mbox{or}\qquad
\Phi^+=\{\alpha_1, \alpha_2, \alpha_1+\alpha_2\}\]
which is, of course, consistent with the general description of the set 
of roots $\Phi$ for $\slm_n(\C)$, $n\geq 2$, in Example~\ref{cartsln}.
\end{exmp}

\begin{exmp} \label{cartanB2} Let $\cL=\gom_4(Q_4,\C)$ where 
$Q_4^{\text{tr}}=-Q_4$, as in Section~\ref{sec05}. We will see in 
Proposition~\ref{CKbcd} below that $\cL$ is of Cartan--Killing type with 
respect to a set $\Delta=\{\alpha_1,\alpha_2\}$ and structure matrix 
\[\renewcommand{\arraystretch}{0.8}
A=\left(\begin{array}{rr} 2 & -1 \\ -2 & 2\end{array}\right).\]
Using ($\clubsuit$), the matrices of $s_1,s_2\in W$ with respect to $\Delta$ 
are:
\[\renewcommand{\arraystretch}{0.8}
s_1:\left(\begin{array}{rr} -1 & 1 \\ 0 & 1 \end{array}\right),\qquad 
s_2:\left(\begin{array}{rr} 1 & 0 \\ 2 & -1 \end{array}\right).\]
Now $s_1s_2\in W$ has order $4$ and so $W$ is dihedral of order~$8$, 
consisting of the elements:
\[\renewcommand{\arraystretch}{0.8} 
\pm \left(\begin{array}{rr} 1 & 0 \\0 & 1\end{array}\right),
\pm \left(\begin{array}{rr} -1 & 1 \\ 0 & 1\end{array}\right),
\pm \left(\begin{array}{rr} 1 & 0 \\2 & -1\end{array}\right),
\pm \left(\begin{array}{rr} 1 & -1 \\ 2& -1\end{array}\right).\]
As above, we obtain that $\cC^+(A)=\{\,(1,0),\,(0,1),\,(1,1),\,(1,2)\,\}$.
Of course, this will turn out to be consistent with the general description 
of the set of roots $\Phi$ for $\gom_n(Q_n,\C)$ in Remark~\ref{sec164} below.
\end{exmp}

\begin{exmp} \label{cartanG2} Consider the matrix
$\renewcommand{\arraystretch}{0.8} A=\left(\begin{array}{rr} 2 & -1 
\\ -3 & 2\end{array}\right)$.\\
We have not yet seen a corresponding Lie algebra but we can just 
formally apply the above procedure, where $\{\alpha_1,\alpha_2\}$ 
denotes the standard basis of~$\C^2$. Using ($\clubsuit$), the 
matrices of $s_1,s_2 \in\GL_2(\C)$ are:
\[\renewcommand{\arraystretch}{0.8}
s_1:\left(\begin{array}{rr} -1 & 1 \\ 0 & 1 \end{array}\right),\qquad 
s_2:\left(\begin{array}{rr} 1 & 0 \\ 3 & -1 \end{array}\right).\]
The product $s_1s_2$ has order $6$ and so $\langle s_1,s_2\rangle
\subseteq \GL_2(\C)$ is a dihedral group of order~$12$. Applying 
$s_1,s_2$ repeatedly to $\alpha_1,\alpha_2$ (avoiding $s_i(\alpha_i)=
-\alpha_i$ for $i=1,2$), we find the following set $\cC^+(A)$:
\[\{\;(1,0),\quad(0,1),\quad(1,1),\quad(1,2),\quad(1,3),\quad(2,3)\;\}\]
(or $\{\alpha_1,\alpha_2, \alpha_1+\alpha_2,\alpha_1+2\alpha_2,\alpha_1+3
\alpha_2, 2\alpha_1+ 3\alpha_2\}\subseteq \C^2$). This discussion will be 
continued in Example~\ref{hypog2} below. 
\end{exmp}

\begin{table}[htbp] \caption{A {\sf Python} program for computing $\cC^+(A)$}
\label{pythontab} 
{\small \begin{verbatim}
>>> def refl(A,n,r,i):          # apply s_i to root r
...   nr=r[:]                   # make a copy of the root r
...   nr[i]-=sum(A[i][j]*nr[j] for j in range(n))
...   return nr
>>> def rootsystem(A):             # A=structure matrix
...   n=len(A)
...   R=[[0]*n for i in range(n)]  # initialise R with
...   for i in range(n):           # unit basis vectors
...     R[i][i]=1
...   for r in R:                 
...     for i in range(n):         
...       if R[i]!=r:          # avoid s_i(alpha_i)=-alpha_i
...         nr=refl(A,n,r,i)   # apply s_i to r
...         if not nr in R:    # check if we get something new
...           R.append(nr)       
...   R.sort(reverse=True)     # sort list nicely
...   R.sort(key=sum)          
...   return R
>>> rootsystem([[2, -1], [-3, 2]])  # see Example 2.3.10 
[[1, 0], [0, 1], [1, 1], [1, 2], [1, 3], [2, 3]]
\end{verbatim}}
\end{table}

The above examples illustrate how $\Phi=\Phi^+\cup (-\Phi^+)$ can be 
computed by a purely mechanical procedure from the structure matrix~$A$. 
In fact, we do not have to do this by hand, but we can simply write a 
computer program for this purpose. Table~\ref{pythontab} contains such 
a program written in the \Python\ language; see 
\url{http://www.python.org}. (It is a version of the basic \nm{orbit 
algorithm}; see, e.g., Holt et al.\cite[\S 4.1]{handbook}.) 
The function {\tt refl(A,|I|,r,i)} implements the formula ($\clubsuit$) 
in Remark~\ref{explicit}.) It outputs the set $\cC^+(A)$, where the 
ordering of the roots is exactly the same as in \CHEVIE\ \cite{chv},
\cite{jchv}. If we apply the program to an arbitrary matrix~$A$, then 
it will either return some nonsense or run into an infinite loop.

\begin{xca} \label{xcacartanAA} Of course, the above procedure will 
not work with any integer matrix $A$, even if the entries of $A$ 
satisfy the various conditions that we have seen so far. For example, 
let $A$ be 
\[\renewcommand{\arraystretch}{0.8} \left(\begin{array}{rrr} 2 & -1 & 
-1 \\ -1 & 2 & -1\\ -1 & -1 & 2 \end{array}\right) \quad\mbox{or}
\quad \left(\begin{array}{rrr} 2 & -1 & 0 \\ -2 & 2 & -1\\ 0 & -3 & 2 
\end{array}\right).\]
Define $s_1,s_2,s_3\in\GL_3(\C)$ using ($\clubsuit$); show that 
$|\langle s_1,s_2,s_3\rangle|=\infty$. 
\end{xca}

\begin{rem} \label{ggtroots} Let $\alpha\in\Phi$ and write $\alpha=
\sum_{i\in I} n_i\alpha_i$, with $n_i\in\Z$ for $i\in I$. We claim that
there is no prime number $p$ such that $p\mid n_i$ for all $i\in I$. 
Indeed, we can write $\alpha=w(\alpha_j)$ for some $w\in W$ and $j\in I$. 
Furthermore, $w=s_{i_1}\cdots s_{i_r}$ where $i_1,\ldots,i_r\in I$. By 
Remark~\ref{explicit}, each $s_{i_j}$ is represented by a matrix with 
entries in $\Z$ with respect to the basis $\Delta$ of~$\fh^*$. Hence, 
the same is also true for~$w$. Since $\alpha=w(\alpha_j)$, the entries 
in the $j$-th column of the matrix of~$w$ are precisely the coefficients 
$(n_i)_{i\in I}$. If there was a prime number $p$ such that $p\mid n_i$ 
for all $i\in I$, then we could conclude that $p\mid \det(w)$. But this is 
a contradiction since $s_i^2=\id_{\fh^*}$ for all~$i$, and so $\det(w)=\pm 1$.
\end{rem}

\begin{rem} \label{keylem1} Consider the structure matrix $A=(a_{ij})_{i,j
\in I}$. The formula in Remark~\ref{keylem0} shows that
\begin{equation*}
a_{ij}=\alpha_j(h_i)=2\frac{\langle \alpha_i,\alpha_j\rangle}{\langle
\alpha_i,\alpha_i\rangle} \qquad \mbox{for all $i,j\in I$}.\tag{$*$}
\end{equation*}
This has the following implication on $A$. Let us set $d_i:=\langle 
\alpha_i,\alpha_i\rangle$ for $i \in I$. Since all elements $w\in W$
are represented by integer matrices with respect to the basis $\Delta$ 
of $\fh^*$ (see Remark~\ref{explicit}), we see from the above definition 
of $\langle\;,\;\rangle$ that $d_i\in\Z_{>0}$. 
Then ($*$)
implies that 
\[ d_ia_{ij}=2\langle \alpha_i,\alpha_j\rangle=2\langle \alpha_j,
\alpha_i\rangle=a_{ji}d_j \quad \mbox{for all $i,j \in I$}.\]
Hence, if we denote by $D\in M_I(\Z)$ the diagonal matrix with diagonal 
entries $d_i$ ($i \in I$), then $D\cdot A\in M_I(\Z)$ is a symmetric matrix. 
In fact, $D\cdot A$ is (up to the factor $2$) the Gram matrix of $\langle\;,\;
\rangle$ with respect to the basis $\Delta$ of $E$. Since 
$\langle \;,\;\rangle$ is positive-definite, a well-known result from 
Linear Algebra shows that $\det(D\cdot A)>0$; since also $\det(D)>0$,
we have $\det(A)>0$. 
\end{rem}


\nmi{}{graph of $A$} 
The above remarks have the following consequence on the \nm{combinatorial
graph of $A$}, which is defined as follows\footnote{Here, we only 
use very basic notions from graph theory, as in Bourbaki 
\cite[Ch.~IV, Annexe]{B}. There are no loops, that is, only distinct edges 
can be joined by an edge; there are no multiple edges and no 
orientations on the edges between vertices.}. 
The set of vertices is~$I$; two vertices $i,j\in I$, $i\neq j$, are 
joined by an edge if $a_{ij}\neq 0$. (Recall that $a_{ij} \neq 0 
\Leftrightarrow a_{ji}\neq 0$.) A \nmi{cycle}{cycle (graph theory)}
(or ``circuit'') in this graph is a sequence of {\itshape distinct} 
indices $i_1,\ldots,i_r \in I$, where $r\geq 3$, such that 
\[ a_{i_1i_2}a_{i_2i_3}\cdots a_{i_{r-1}i_r}a_{i_ri_1}\neq 0.\]
Thus, $i_1,i_2$ are joined, then $i_2,i_3$ are joined and so on until 
$i_{r-1},i_r$ are joined; finally, $i_r$ is joined back to~$i_1$. The graph 
is called a \nmi{forest}{forest (graph theory)} if there are no cycles.


\begin{prop} \label{forestg} The graph of the structure matrix $A$ is
a forest. 
\end{prop}


\begin{proof} Assume that the graph of $A$ has a cycle as above; let 
$I':=\{i_1,\ldots,i_r\}$. Define $x=(x_i)_{i\in I} \in \R^I$ by 
$x_i:=1$ if $i \in I'$, and $x_i:=0$ otherwise. Regard $x$ as a column 
vector and consider the product $y:=A\cdot x$. For $i\in I$, the
$i$-th component of $y$ is given by
\[ y_i=\sum_{j \in I} a_{ij}x_j=\sum_{j\in I'}a_{ij}.\]
Now, if $i\in I'$, there are at least two distinct vertices
$j_1,j_2\in I'\setminus\{i\}$ that are joined to~$i$; then 
$a_{ij_1}\leq -1$ and $a_{ij_2}\leq -1$. Since $a_{ii}=2$ and $a_{ij}
\leq 0$ for all $j\in I'\setminus\{i\}$, we conclude that $y_i\leq 
a_{ii}+a_{ij_1}+a_{ij_2}\leq 2-1-1=0$. Now let $D$ be the diagonal matrix
in Remark~\ref{keylem1}; then $S:=D\cdot A$ is a symmetric matrix which 
defines a positive-definite symmetric bilinear form on $\R^I$. Hence, 
$x^{\text{tr}}\cdot S \cdot x>0$ since $x\neq 0$. But we also have
$d_i>0$ and $y_i\leq 0$ for all $i\in I'$; this implies that 
\[ x^{\text{tr}} \cdot S\cdot x=x^{\text{tr}}\cdot D\cdot y=
\sum_{i \in I} x_id_iy_i=\sum_{i \in I'} d_iy_i\leq 0,\]
contradiction. Thus, there are no cycles in the graph of $A$. 
\end{proof}

\begin{xca} \label{xcaforest} Consider a graph as above with a finite 
(non-empty) vertex set $I$. Assume that the graph is a forest. \\
(a) Show that there is a partition $I=I_1\sqcup I_2$ such that two distinct 
vertices that both belong to $I_1$ or to $I_2$ are never joined by an edge.\\
(b) Deduce that there is a function $f\colon I \rightarrow \{\pm 1\}$ such 
that $f(i)=-f(j)$ whenever $i,j\in I$ are joined by an edge in the graph.\\
{\footnotesize [{\it Hint}. For (a) see \cite[Ch.~IV, Annexe, no.~3]{B}. 
For (b) define the function $f\colon I \rightarrow \{\pm 1\}$
by $f(i)=1$ if $i\in I_1$ and $f(i)=-1$ if $i\in I_2$.]}
\end{xca}

\section{Semisimplicity} \label{sec1a4}

We continue to assume that $(\cL,\fh)$ is of Cartan--Killing type with 
respect to $\Delta=\{\alpha_i\mid i\in I\}$. In this section, we establish
the main structural properties of $\cL$. For each $i\in I$ let $\{e_i,h_i,
f_i\}$ be a corresponding $\slm_2$-triple, as in Remark~\ref{astring0}.
Let $W=\langle s_i\mid i \in I\rangle\subseteq \GL(E)$ be the Weyl group of
$(\fg,\fh)$.

The first step consists of ``lifting'' the generators $s_i$ of $W$ to 
Lie algebra automorphisms of $\cL$. By Lemma~\ref{wsdnil}, the derivations
$\ad_\cL(e_i) \colon \cL\rightarrow \cL$ and $\ad_\cL(f_i)\colon \cL
\rightarrow \cL$ are nilpotent. Hence, $t\,\ad_\cL(e_i)$ and $t\,
\ad_\cL(f_i)$ are nilpotent derivations for all $t\in\C$. So we can apply 
the exponential construction in Lemma~\ref{exponential}, and set
\begin{align*}
x_i(t)&:=\exp\bigl(t\,\ad_\cL(e_i)\bigr) \in\mbox{Aut}(\cL)\quad
\mbox{for all $t\in \C$},\\
y_i(t)&:=\exp\bigl(t\,\ad_\cL(f_i)\bigr) \in\mbox{Aut}(\cL)\quad
\mbox{for all $t\in \C$}.
\end{align*}

\begin{lem} \label{weyl5} With the above notation, we define
\[n_i(t):=x_i(t)\circ y_i(-t^{-1})\circ x_i(t) \in\operatorname{Aut}(\cL)
\quad\mbox{for $0\neq t\in \C$}.\]
Then the following hold.
\begin{align*}
n_i(t)(h)&=h-\alpha_i(h)h_i\in \fh \qquad\mbox{for all $h\in \fh$}, \tag{a}\\
\lambda\bigl(n_i(t)(h)\bigr)&=s_i(\lambda)(h) \qquad\mbox{for all 
$\lambda\in \fh^*$ and $h\in \fh$},\tag{b}\\ 
n_i(t)(\cL_\alpha)&=\cL_{s_i(\alpha)} \qquad \mbox{for all $\alpha \in\Phi$}.
\tag{c}
\end{align*}
\end{lem}

\begin{proof} (a) Let $h\in \fh$. Let us first determine $x_i(t)(h)$. 
For this purpose, we need to work out $\ad_\cL(e_i)^m(h)$ for all $m\geq 1$.
Now, we have $\ad_\cL(e_i)(h)=[e_i,h]=-[h,e_i]=-\alpha_i(h)e_i$ and, 
consequently, $\ad_\cL(e_i)^m(h)=0$ for all $m\geq 2$. This already shows 
that 
\[x_i(t)(h)=\sum_{m\geq 0} \frac{(t\,\ad_\cL(e_i))^m(h)}{m!}=
h-\alpha_i(h)te_i.\]
Similarly, we have $\ad_\cL(f_i)(h)=[f_i,h]=-[h,f_i]=\alpha_i(h)f_i$ and, 
consequently, $\ad_\cL(f_i)^m(h)=0$ for all $m\geq 2$. This shows that
\[y_i(t)(h)=\sum_{m\geq 0} \frac{(t\,\ad_\cL(e_i))^m(h)}{m!}=
h+\alpha_i(h)tf_i.\]
Next, we determine $y_i(t)(e_i)$. We have $\ad_\cL(f_i)(e_i)=
-[e_i,f_i]=-h_i$, $\ad_\cL^2(f_i)(e_i)=-[f_i,h_i]=-2f_i$ and, 
consequently, $\ad_\cL(f_i)^m(e_i)=0$ for all $m\geq 3$. This yields that 
\[y_i(t)(e_i)=\sum_{m\geq 0} \frac{(t\,\ad_\cL(f_i))^m(e_i)}{m!}
=e_i-th_i-t^2f_i.\]
(Similarly, one has $x_i(t)(f_i)=f_i+th_i-t^2e_i$.)
Combining the above formulae, we obtain that
\begin{align*}
\bigl(y_i(-t^{-1})&\circ x_i(t)\bigr)(h)=y_i(-t^{-1})\bigl(h- \alpha_i(h)
te_i \bigr)\\&=\bigl(h-\alpha_i(h)t^{-1}f_i\bigr)-\alpha_i(h)t
\bigl(e_i+t^{-1}h_i-t^{-2}f_i\bigr)\\
& =h-\alpha_i(h)h_i -\alpha_i(h)te_i.
\end{align*}
Finally, $\ad_\cL(e_i)^m(e_i)=0$ for all $m\geq 1$ and so
$x_i(t)(e_i)=e_i$. Hence, 
\begin{align*}
n_i(t)(h)&=x_i(t)\bigl(h-\alpha_i(h)h_i-\alpha_i(h)te_i\bigr)\\&=
\bigl(h-\alpha_i(h)te_i\bigr)-\alpha_i(h)\bigl(h_i-2te_i\bigr)-
\alpha_i(h)te_i\\&=h-\alpha_i(h)h_i.
\end{align*}

(b) Recall that $s_i(\lambda)=\lambda-\lambda(h_i) \alpha_i$. Using (a), 
this yields:
\begin{align*}
\lambda\bigl(n_i(t)(h)\bigr)&=\lambda\bigl(h-\alpha_i(h)h_i\bigr)
=\lambda(h)-\alpha_i(h)\lambda(h_i)\\&=(\lambda-\lambda(h_i)\alpha_i)(h)=
s_i(\lambda)(h)
\end{align*}
for all $h\in \fh$, as desired.

(c) Let $h\in \fh$ and set $h':=n_i(t)(h)\in \fh$. Since $\alpha_i(h_i)=2$, we 
see using (a) that $n_i(t)(h_i)=-h_i$; furthermore,
\[ n_i(t)(h')=n_i(t)\bigl(h-\alpha_i(h)h_i\bigr)=n_i(t)(h)+\alpha_i(h)h_i
=h.\]
Now let $y \in \cL_{\alpha}$ and set $y':=n_i(t)(y)\in \cL$. Then 
\begin{align*}
[h,y']&=[n_i(t)(h'),n_i(t)(y)]=n_i(t)\bigl([h',y]
\bigr)\\ &=n_i(t)\bigl(\alpha(h')y\bigr)
=\alpha(h')n_i(t)(y)=\alpha(h')y',
\end{align*}
where the second equality holds since $n_i(t)$ is a Lie algebra
automorphism. Now, by (b), we have $\alpha(h')=s_i(\alpha)(h)$ and so 
$y'=n_i(t)(y)\in \cL_{s_i(\alpha)}$. Hence, $n_i(t)(\cL_\alpha) \subseteq 
\cL_{s_i(\alpha)}$ and $\dim \cL_\alpha\leq \dim \cL_{s_i(\alpha)}$. Since 
$s_i^2=\id_{\fh^*}$, we also obtain that $n_i(t)(\cL_{s_i(\alpha)})
\subseteq \cL_{s_i^2(\alpha)}=\cL_\alpha$ and so $\dim \cL_{s_i(\alpha)}\leq 
\dim \cL_\alpha$. Hence, we must have $n_i(t)(\cL_\alpha)=\cL_{s_i(\alpha)}$.
\end{proof}

\begin{xca} \label{weyl5a} Let $i\in I$ and $0\neq t\in \C$. Using the
formulae obtained in the above proof, deduce that 
\[ n_i(t)(e_i)=-t^{-2}f_i\qquad\mbox{and}\qquad n_i(t)(f_i)=-t^2e_i.\]
(Much later, in Theorem~\ref{luform3}, we will even obtain explicit
formulae for $n_i(t)(y)$ where $y \in \cL_\alpha$ and $\alpha\in \Phi$ is
arbitrary.)
\end{xca}

\begin{prop}  \label{weyl5b} We have $\dim \cL_\alpha=1$ and 
$\dim [\cL_\alpha,\cL_{-\alpha}]=1$ for all $\alpha\in \Phi$.  In
particular, $\dim \cL=|I|+|\Phi|$.
\end{prop}

\begin{proof} Let $\alpha\in\Phi$. By Theorem~\ref{mainthm1}(a) we can
write $\alpha=w(\alpha_i)$ for some $w\in W$ and $i\in I$. Furthermore,
we can write $w=s_{i_1}\cdots s_{i_r}$, where $r\geq 0$ and $i_1,\ldots,
i_r \in I$. Let us set $\eta:=n_{i_1}(1)\circ \ldots \circ n_{i_r}(1)
\in \mbox{Aut}(\cL)$. Now Lemma~\ref{weyl5}(c) and a simple induction 
on~$r$ show that 
\[ \cL_\alpha=\cL_{(s_{i_1}\cdots s_{i_r})(\alpha_i)}=
\bigl(n_{i_1}(1)\circ \ldots \circ n_{i_r}(1)\bigr)(\cL_{\alpha_i})=
\eta(\cL_{\alpha_i}).\]
Since $\eta\in\mbox{Aut}(\cL)$, we conclude that $\dim \cL_\alpha=
\dim \cL_{\alpha_i}=1$, where the last equality holds by 
Proposition~\ref{wsdprop3}. Furthermore, since $-\alpha=-w(\alpha_i)=
w(-\alpha_i)$, the same argument shows that 
$\cL_{-\alpha}=\eta(\cL_{-\alpha_i})$. Again, since $\eta\in
\mbox{Aut}(\cL)$, we also have 
\[[\cL_\alpha,\cL_{-\alpha}]=[\eta(\cL_{\alpha_i}),\eta(\cL_{-\alpha_i})]
=\eta\bigl([\cL_{\alpha_i},\cL_{-\alpha_i}]\bigr), \]
and this is $1$-dimensional by Proposition~\ref{wsdprop3}. Finally, the
formula for $\dim \cL$ follows from the direct sum decomposition 
$\cL=\fh\oplus \bigoplus_{\alpha\in \Phi} \cL_\alpha$ and the fact 
that $\{h_i\mid i \in I\}$ is a basis of $\fh$.
\end{proof}

\begin{prop} \label{weyl5c} For each $\alpha\in\Phi$, there exists
a unique element $h_\alpha\in [\cL_\alpha,\cL_{-\alpha}]$ such that 
$\alpha(h_\alpha)=2$. (We have $h_{\alpha_i}=h_i$ for $i \in I$.)
Furthermore, $h_{-\alpha}=-h_\alpha$ and 
\[ h_{s_i(\alpha)}=n_i(t)(h_\alpha)=h_\alpha-\alpha_i(h_\alpha)h_i
\qquad \mbox{for $i\in I$ and $t\in \C^\times$}.\] 
\end{prop}

\begin{proof} By Proposition~\ref{weyl5b}, we have $[\cL_\alpha,
\cL_{-\alpha}]=\langle h\rangle_\C$ for some $0\neq h\in \fh$. If 
$\alpha(h)=0$, then Lemma~\ref{wsdlem4} would imply that $\ad_\cL(h)=0$. 
In particular, all eigenvalues of $\ad_\cL(h)$ are zero and so 
$\alpha_i(h)=0$ for all $i\in I$, contradiction since $\{\alpha_i\mid i
\in I\}$ is a basis of $\fh^*$. Thus, $\alpha(h)\neq 0$ and so there is a 
unique scalar multiple of $h$ on which $\alpha$ takes value~$2$. This 
defines the required element $h_\alpha$.

Since $-\alpha\in \Phi$ and $[\cL_{-\alpha},\cL_{\alpha}]=
[\cL_{\alpha},\cL_{-\alpha}]$ is $1$-dimensional, we have $h_{-\alpha}=
\xi h_\alpha$ for some $0\neq \xi\in \C$. But then we conclude that 
$2=(-\alpha)(h_{-\alpha})=-\xi\alpha(h_\alpha)=-2\xi$ and so $\xi=-1$. 

Now let $i \in I$. By Lemma~\ref{weyl5}(c), we have $n_i(t)(\cL_\alpha)=
\cL_{s_i(\alpha)}$ and $n_i(t)(\cL_{-\alpha})=\cL_{-s_i(\alpha)}$. Hence,
we obtain 
\[ \langle n_i(t)(h_\alpha)\rangle_\C=n_i(t)([\fg_\alpha,\fg_{-\alpha}])=
[\fg_{s_i(\alpha)},\fg_{-s_i(\alpha})]=\langle h_{s_i(\alpha)}\rangle_\C\]
and so $h_{s_i(\alpha)}=\xi n_i(t)(h_\alpha)$ for some $0\neq \xi\in \C$. 
Now, applying Lemma~\ref{weyl5}(b) with $\lambda=s_i(\alpha)$, we obtain 
\[ s_i(\alpha)\bigl(n_i(t)(h_\alpha)\bigr)=\lambda(n_i(t)(h_\alpha))=
s_i(\lambda)(h_\alpha)=\alpha(h_\alpha)=2.\]
Since also $s_i(\alpha)(h_{s_i(\alpha)})=2$, we conclude that $\xi=1$.
Finally, we have $n_i(t)(h_\alpha)=h_\alpha-\alpha_i(h_\alpha)h_i$
by Lemma~\ref{weyl5}(a).
\end{proof}

\begin{xca} \label{xcaalgoh} (a) By Lemma~\ref{weyl5}, we have
$n_i(t)(\fh)\subseteq \fh$ for all $i\in I$ and $0\neq t\in \C$. Show that
$n_i(t)^2(h)=h$ for all $h\in \fh$. Furthermore, show that the matrix 
of $n_i(t)|_\fh\colon \fh\rightarrow \fh$ with respect to the basis $\{h_i
\mid i \in I\}$ of $\fh$ has integer coefficients and determinant~$-1$.

(b) Let $\alpha\in\Phi$ and write 
$\alpha=w(\alpha_i)$ where $w\in W$ and $i\in I$; further write 
$w=s_{i_1}\cdots s_{i_r}$ where $i_1,\ldots,i_r\in I$. Show that 
\[ h_\alpha=\bigl(n_{i_1}(1) \circ \ldots \circ n_{i_r}(1)\bigr)(h_i)
\in \langle h_j \mid j \in I\rangle_\Z.\]
Using also Remark~\ref{astring}(a), deduce that $\beta(h_\alpha)\in \Z$ 
for all $\beta\in \Phi$.
\end{xca}

The following result shows that the ``\nm{Chevalley generators}'' in 
Remark~\ref{astring0} are indeed generators for $\cL$ as a Lie algebra.

\begin{prop} \label{genlie} We have $\cL=\langle e_i,f_i\mid i\in I
\rangle_{\operatorname{alg}}$.
\end{prop}

\begin{proof} Let $\cL_0:=\langle e_i, f_i \mid i \in I\rangle_{\text{alg}}
\subseteq \cL$. Since $h_i=[e_i,f_i]\in \cL_0$ for all $i$, we have 
$\fh\subseteq \cL_0$. So it remains to show that $\cL_{\pm \alpha}\subseteq
\cL_0$ for all $\alpha\in\Phi^+$. We proceed by induction on $\hgt(\alpha)$. 

If $\hgt(\alpha)=1$, then $\alpha=\alpha_i$ for some $i\in I$. Since 
$\cL_{\alpha_i}= \langle e_i\rangle_\C$ and $\cL_{-\alpha_i}=\langle f_i
\rangle_\C$, we have $\cL_{\pm\alpha_i}\subseteq \cL_0$ by the definition of 
$\cL_0$. Now let $\hgt(\alpha)>1$. By the Key Lemma~\ref{keylem}, there exists 
some $j \in I$ such that $\beta:=\alpha-\alpha_j\in \Phi^+$. We have 
$\hgt(\beta)=\hgt(\alpha)-1$ and so, by induction, $\cL_{\pm\beta}\subseteq 
\cL_0$. By Remark~\ref{astring}(c'), since $\alpha_j+\beta=\alpha\in\Phi$, we 
have $\{0\}\neq [\cL_{\alpha_j},\cL_{\beta}]\subseteq \cL_{\alpha_j+\beta}=
\cL_{\alpha}$. Since $\dim \cL_{\alpha}=1$ (see Proposition~\ref{weyl5b}), we 
conclude that $\cL_\alpha=[\cL_{\alpha_j},\cL_{\beta}]$, and this is contained 
in $\cL_0$ because $\cL_0$ is a subalgebra and $\cL_{\alpha_j},\cL_{\beta}
\subseteq \cL_0$. Similarly, $-\alpha=-\alpha_j- \beta$ and 
$\cL_{-\alpha}=[\cL_{-\alpha_j}, \cL_{-\beta}]\subseteq \cL_0$. 
\end{proof}

\begin{lem} \label{mainideal1} Let $\fj\subseteq \cL$ be an ideal. 
If $\fj\neq \{0\}$, then $\fj\cap \fh\neq \{0\}$. In fact,
in this case, we have $h_i\in \fj$ for some $i \in I$.
\end{lem}

\begin{proof} Since $\fj$ is an ideal, we have $[\fh,\fj]\subseteq \fj$ 
which means that $\fj$ is an $\fh$-submodule of $\cL$. So 
Proposition~\ref{wsdprop1}(b) implies that 
\[ \fj=(\fj \cap \fh)\oplus\bigoplus_{\alpha\in \Phi}(\fj\cap\cL_\alpha).\]
Hence, if $\fj\neq \{0\}$, then $\fj\cap \fh \neq \{0\}$ or $\fj\cap
\cL_\alpha\neq \{0\}$ for some $\alpha \in \Phi$. Assume that we are in
the second case. By Proposition~\ref{weyl5b}, we have $\dim \cL_\alpha=1$ 
and so $\cL_\alpha \subseteq \fj$. So, by Proposition~\ref{weyl5c}
(and since $\fj$ is an ideal), we also have $h_\alpha\in [\cL_\alpha,
\cL_{-\alpha}] \subseteq \fj$. Hence, in any case, we have $\fj\cap \fh
\neq \{0\}$. Let $0\neq h \in \fj\cap \fh$. Since $\{\alpha_i \mid i 
\in I\}$ is a basis of $\fh^*$, we have $\alpha_i(h) \neq 0$ for some
$i \in I$. We deduce that $\alpha_i(h) e_i=[h,e_i]\in \fj$ and so
$e_i \in \fj$. But then we also have $h_i=[e_i,f_i]\in \fj$, as desired.
\end{proof}

\begin{thm} \label{mainideal} The Killing form $\kappa_\cL\colon \cL \times 
\cL \rightarrow\C$ (as in De\-fini\-tion~\ref{defkill0}) is non-degenerate. 
Consequently, $\cL$ is semisimple.
\end{thm}

\begin{proof} Consider $\fg^\perp=\{x\in \cL\mid \kappa_\cL(x,y)=0 
\mbox{ for all $y\in \cL$}\}$. We want to show that $\fg^\perp=\{0\}$. 
Assume that $\fg^\perp\neq \{0\}$. By Lemma~\ref{defkill0a}(b), 
$\fg^\perp$ is an ideal. Hence, by Lemma~\ref{mainideal1}, there exists
some $i\in I$ such that $h_i \in \fg^\perp$ and so $\kappa_\cL(h_i,h_i)=0$. 
For each $\alpha\in \Phi$, let $0\neq e_\alpha\in \cL_\alpha$. Then 
$\{h_j\mid j \in I\} \cup \{e_\alpha\mid \alpha\in \Phi\}$ is a basis of 
$\cL$; see Proposition~\ref{weyl5b}. With respect to this basis, the matrix 
of $\ad_\cL(h_i)$ is diagonal with eigenvalues $0$ (with multiplicity~$|I|$) 
and $\alpha(h_i)$ (each with multiplicity~$1$). Hence, we obtain that
\[ 0=\kappa_\cL(h_i,h_i)=\trc\bigl(\ad_\cL(h_i)\circ \ad_\cL(h_i)\bigr)=
\sum_{\alpha\in \Phi} \alpha(h_i)^2.\]
By Remark~\ref{poshi}, we have $\alpha(h_i)\in \Z$ for all $\alpha\in 
\Phi$. So we must have $\alpha(h_i)=0$ for all $\alpha\in \Phi$,
contradiction since $h_i\neq 0$ and $\fh^*=\langle \Phi\rangle_\C$. Thus,
our assumption was wrong and so $\kappa_\fg$ is non-degenerate. Finally, 
by Lemma~\ref{defkill0a}(c), this implies that $\cL$ is semisimple.
\end{proof}

\begin{defn} \label{cindec1} Consider the structure matrix $A=(a_{ij})_{i,
j\in I}$ of~$\cL$ or, somewhat more generally, any matrix $A=(a_{ij})_{i,j
\in I}$ such that the entries~$a_{ij}$ satisfy the conditions (a), (b), 
(c) in Corollary~\ref{gencart0}. Assume that $I\neq \varnothing$. We say 
that $A$ is \textit{decomposable} if there is a partition $I=I_1\sqcup 
I_2$ (where $I_1,I_2\subsetneqq I$ and $I_1\cap I_2=\varnothing$) such 
that $a_{ij}= a_{ji}=0$ for all $i\in I_1$ and $j\in I_2$. In this case
we can arrange~$I$ such that $A$ has a block diagonal shape
\begin{center}
$A=\left(\begin{array}{c|c} A_1 & 0 \\\hline 0 & A_2\end{array}\right)$
\end{center}
where $A_1$ has rows and columns labelled by $I_1$, and $A_2$ has rows 
and columns labelled by $I_2$. If no such partition of $I$ exists, then 
we say that $A$ is \nmi{indecomposable}{indecomposable matrix}. Note that
the indexing set $I$ can always be arranged such that $A$ has a block 
diagonal matrix where the diagonal blocks are indecomposable. 
\end{defn}

\begin{rem} \label{indecgraph} Consider the \nm{combinatorial graph of $A$} 
introduced at the end of Section~\ref{sec1a3}. A standard argument in graph 
theory shows that this graph is connected if and only if $A$ is 
indecomposable (see, e.g., \cite[Ch.~IV, Annexe, Cor.~1]{B}). Hence, the 
indecomposability of $A$ can be alternatively expressed as follows.
For any $i,j\in I$ such that $i\neq j$, there exists a sequence of 
(distinct) indices $i=i_0,i_1,\ldots,i_r=j$ in $I$, where $r\geq 1$ and 
$a_{i_li_{l+1}} \neq 0$ for $0\leq l\leq r-1$. 
\end{rem} 

\begin{rem} \label{mainideal3} Let $I=I_1\sqcup I_2$ be a partition
as in Definition~\ref{cindec1}. Then we define the following Lie 
subalgebras of $\fg$:
\[ \fg_1:=\langle e_i,f_i\mid i \in I_1\rangle_{\text{alg}} \qquad\mbox{and}
\qquad \fg_2:=\langle e_j,f_j\mid j \in I_2\rangle_{\text{alg}}.\]
We claim that $[\fg_1,\fg_2]=\{0\}$. Let $i\in I_1$ and $j\in I_2$.
Since $a_{ij}=0$, we have $\alpha_i+\alpha_j\not\in \Phi$ and so 
$[e_i,e_j]=0$; see Exercise~\ref{xcastringm}. Since we also have 
$[e_i,f_j]=0$, it follows that $[e_i,\fg_2]=\{0\}$; see 
Exercise~\ref{xcagenerator}(a). Similarly, we see that $[f_i,\fg_2]=\{0\}$. 
But then Exercise~\ref{xcagenerator}(a) also implies that $[\fg_1,
\fg_2]=\{0\}$, as claimed. In particular, $\fg_1$ and $\fg_2$ are ideals
in $\fg$. Now set $\fg':=\fg_1+\fg_2\subseteq \fg$. Since $[\fg_1,\fg_2]=
\{0\}$, this is a Lie subalgebra of $\fg$. But then 
Proposition~\ref{genlie} implies that $\fg'=\fg$, since $e_i,f_i\in 
\fg'$ for all $i\in I$. Furthermore, for any $x,y\in \fg_1\cap \fg_2$, 
we have $[x,y]\in [\fg_1,\fg_2]=\{0\}$. Hence, $\fg_1\cap \fg_2$ is 
an abelian ideal of $\fg$ and so $\fg_1\cap \fg_2=\{0\}$. Thus, we
obtain a direct sum decomposition
\[\fg=\fg_1\oplus \fg_2\qquad\mbox{where}\qquad [\fg_1,\fg_2]=\{0\}.\]
\end{rem}

\begin{rem} \label{mainideal2} Let $\fg_1\subseteq \fg$ be an ideal. By 
Theorem~\ref{mainideal}, the Killing form $\kappa_\cL$ is non-degenerate. 
Hence, $\fg_2:=\fg^\perp$ also is an ideal and $\fg=\fg_1\oplus \fg_2$;
see Proposition~\ref{killing8b}(a). Now set 
\[ I_1:=\{i \in I \mid h_i\in \fg_1\}\qquad \mbox{and}\qquad 
I_2:=\{i \in I \mid h_i\in \fg_2\}.\] 
We show that $I_1,I_2$ define a partition of $I$ as in 
Definition~\ref{cindec1}. Clearly, we have $I_1\cap I_2=\varnothing$. We 
claim that $I=I_1\cup I_2$. Indeed, let $i\in I$. Since $\fg_1,\fg_2
\subseteq \fg$ are $\fh$-submodules, Example~\ref{wsdexp1} shows that
\[ \fg_{\alpha_i}=(\fg_1\cap \fg_{\alpha_i})\oplus (\fg_2\cap 
\fg_{\alpha_i}).\]
Hence, since $\fg_{\alpha_i}=\langle e_i\rangle_\C$, either $e_i \in 
\fg_1$ or $e_i\in \fg_2$. In the first case, it follows that $h_i=
[e_i,f_i]\in \fg_1$ and so $i \in I_1$. Similarly, in the second case, 
it follows that $i\in I_2$. Thus, $I=I_1\sqcup I_2$. Finally, let 
$i \in I_1$ and $j\in I_2$. Then $a_{ij}e_j=\alpha_j(h_i)e_j=[h_i,e_j] 
\in \fg_1$; furthermore, $2e_j=[h_j,e_j]\in \fg_2$. Hence, $2a_{ij}e_j
\in \fg_1\cap \fg_2=\{0\}$ and so $a_{ij}=0$.
\end{rem}

\begin{thm} \label{cindec4} Assume that $\cL\neq\{0\}$. Then $\cL$ is
simple if and only if $A$ is indecomposable.
\end{thm}

\begin{proof} Assume first that $A$ indecomposable. Let $\fg_1\subseteq 
\cL$ be an ideal. Let $\fg_2:=\fg_1^\perp$ and define $I_1,I_2\subseteq I$ 
as in Remark~\ref{mainideal2}. So $I=I_1\sqcup I_2$ is a partition as
in Definition~\ref{cindec1}. Since $A$ is indecomposable, we must have 
$I=I_1$ or $I=I_2$. Assume that $I=I_1$. Then $h_i\in \fg_1$ for all $i 
\in I$. Since $2e_i=[h_i,e_i]$ and $-2f_i=[h_i,f_i]$, it follows that  
$e_i,f_i\in \fg_1$ for all $i\in I$. So Proposition~\ref{genlie} implies
that $\fg=\fg_1$. If $I=I_2$, then an analogous argument shows that $\fg
=\fg_2$ and so $\fg_1=\{0\}$. Hence, $\fg$ does not have any non-trivial 
ideals. Since $\cL\neq \{0\}$, we have $I\neq \varnothing$ and $\cL$ is 
not abelian (see Definition~\ref{defTD}). Hence, $\fg$ is simple.

Conversely, assume hat $\fg$ is simple. Let $I=I_1\sqcup I_2$ be
a partition as in Definition~\ref{cindec1}. As in Remark~\ref{mainideal3},
we obtain a direct sum decomposition $\fg=\fg_1\oplus \fg_2$ where
$[\fg_1,\fg_2]=\{0\}$. Then $\fg_1,\fg_2$ are ideals in $\fg$. These
are not abelian since $[e_i,f_i]=h_i\neq 0$ for $i\in I$. Hence, since
$\fg$ is simple, we must have $\fg_1=\{0\}$ or $\fg_2=\{0\}$. But then 
$I_1=\varnothing$ or $I_2=\varnothing$. So $A$ is indecomposable. 
\end{proof}

\begin{rem} \label{cindec1a} The above result leads to a simple (!) method
for testing if $\cL$ is a simple Lie algebra: we just have to see if $A$ 
indecomposable. For example, let $\cL=\slm_n(\C)$, where $n\geq 2$. By
Example~\ref{cartsln}, we know that $\cL$ is of Cartan--Killing type.
Now we just note that the structure matrix $A$ in Example~\ref{cartsln} 
is indecomposable. Indeed, the graph of $A$ (as introduced at the end 
of Section~\ref{sec1a3}) is given by 
\begin{center}
\begin{picture}(160,12)
\put( 20,1){\circle*{5}}
\put( 18,6){$1$}
\put( 20,1){\line(1,0){20}}
\put( 40,1){\circle*{5}}
\put( 38,6){$2$}
\put( 40,1){\line(1,0){30}}
\put( 60,1){\circle*{5}}
\put( 58,6){$3$}
\put( 80,1){\circle*{1}}
\put( 90,1){\circle*{1}}
\put(100,1){\circle*{1}}
\put(110,1){\line(1,0){10}}
\put(120,1){\circle*{5}}
\put(111,6){$n{-}1$}
\end{picture}
\end{center}
and this is obviously connected. Hence, $\cL$ is simple.~---~In the next 
section, we will employ a similar argument to show that the Lie algebras 
$\gom_n(Q_n,\C)$ are simple (with the exceptions in 
Exercise~\ref{xcaclassic3}.)
\end{rem}

Finally, we would like to understand a bit better what happens in the 
situation when $A$ is decomposable.  By Proposition~\ref{killing8b}(b),
we can write $\cL$ as a direct sum of ideals which are all simple as 
Lie algebras. Our aim is to show that each of these simple ideals is also 
of Cartan--Killing type. To see this, let us begin with a partition $I=
I_1\sqcup I_2$ as in Definition~\ref{cindec1}. By Remark~\ref{mainideal3} 
we have a corresponding direct sum decomposition $\fg=\fg_1\oplus \fg_2$
where 
\[ \fg_1:=\langle e_i,f_i\mid i \in I_1\rangle_{\text{alg}} \qquad
\mbox{and} \qquad \fg_2:=\langle e_j,f_j\mid j \in I_2
\rangle_{\text{alg}};\]
here, $\fg_1,\fg_2\subseteq \fg$ are ideals such that $[\fg_1,\fg_2]=
\{0\}$. We set 
\begin{alignat*}{2}
\Delta_1&:=\{\alpha_i\mid i\in I_1\} &\qquad \mbox{and}\qquad
&\Delta_2:=\{\alpha_j\mid j \in I_1\};\\
\fh_1&:=\langle h_i\mid i \in I_1\rangle_\C &\qquad\mbox{and} 
\qquad &\;\fh_2:=\langle h_j\mid j \in I_2\rangle_\C;\\
W_1&:=\langle s_i \mid i\in I_1\rangle &\qquad \mbox{and}\qquad 
&W_2:=\langle s_j \mid j\in I_2\rangle.
\end{alignat*}
Since $h_i=[e_i,f_i]$ for all $i\in I$, it is clear that $\fh_1\subseteq
\fg_1$ and $\fh_2\subseteq \fg_2$; note that $\fh=\fh_1\oplus \fh_2$. 
We also set 
\[\Phi_s:=\{w(\alpha_i)\mid w\in W_s, i\in I_s\}\subseteq \Phi 
\qquad \mbox{for $s=1,2$}.\]

\begin{prop} \label{levisub2} In the above setting, the following hold.\\
{\rm (a)} We have $W=W_1\cdot W_2$, $W_1\cap W_2=\{\id\}$ and $w_1w_2=
w_2w_1$ for $w_1\in W_1$ and $w_2\in W_2$. Thus, $W$ is the direct 
product of $W_1$ and~$W_2$. \\
{\rm (b)} We have $\Phi=\Phi_1\sqcup 
\Phi_2$ (disjoint union) and $\Phi_s \subseteq \langle \Delta_s\rangle_\Z$ 
for $s=1,2$. For any $\alpha\in \Phi_1$ and $\beta \in \Phi_2$, we have 
$\alpha\pm \beta \not\in \Phi\cup \{\underline{0}\}$.
\end{prop}

\begin{proof} For $s=1,2$ let $E_s:=\langle \Delta_s\rangle_\Z$. Let 
$i\in I_1$ and $j\in I$. Then $s_i(\alpha_j)=\alpha_j-\alpha_j(h_i)
\alpha_i=\alpha_j-a_{ij}\alpha_i$. Hence, if $j\in I_1$, then 
$s_i(\alpha_j)\in E_1$; if $j\in I_2$, then $s_i(\alpha_j)=\alpha_j$, 
since $a_{ij}=0$.  Consequently, we have:
\begin{equation*}
i\in I_1 \quad\Rightarrow \quad s_i(E_1)\subseteq E_1 \;\mbox{ and }\; 
s_i(v)=v \;\mbox{ for all $v\in E_2$}.\tag{1}
\end{equation*}
Similarly, we see that 
\begin{equation*}
j\in I_2 \quad\Rightarrow \quad s_j(E_2)\subseteq E_2 \;\mbox{ and }\; 
s_j(v)=v \;\mbox{ for all $v\in E_1$}.\tag{2}
\end{equation*}
This immediately implies that 
\begin{equation*}
s_is_j=s_js_i \qquad \mbox{for $i \in I_1$ and $j \in I_2$}.\tag{3}
\end{equation*}
First of all, this shows that $w(E_1)\subseteq E_1$ and $w(E_2)
\subseteq E_2$ for all $w\in W$. (Indeed, by (1) and (2), this property 
holds for all $s_i$ and, hence, it holds for all elements of $W$.) By a 
similar argument, (3) implies that $w_1w_2=w_2 w_1$ for all $w_1\in W_1$ 
and $w_2\in W_2$. But then $W_1 \cdot W_2 \subseteq W$ is a subgroup which 
contains all $s_i$ ($i\in I$). Hence, $W=W_1\cdot W_2$. If $w\in W_1\cap 
W_2$, then (1) implies that $w(v)=v$ for all $v \in E_2$ and (2) implies 
that $w(v)=v$ for all $v \in E_1$. Hence, $w=\id$. 

It remains to show the assertions about $\Phi_s$. Let $\alpha\in \Phi$.
By Theorem~\ref{mainthm1}(a), we have $\alpha=w(\alpha_i)$ for some 
$i\in I$ and $w\in W$. Write $w=w_1w_2=w_2w_1$ where $w_1\in W_1$ and
$w_2\in W_2$. If $i\in I_1$, then (2) implies that $w_2(\alpha_i)=
\alpha_i$ and so $\alpha=w(\alpha_i)=w_1(\alpha_i)\in \Phi_1$. Similarly, 
if $i\in I_2$, then $\alpha\in \Phi_2$. Hence, $\Phi=\Phi_1\cup \Phi_2$. 
Furthermore, by (1), we have $w_1(E_1)\subseteq E_1$ for all $w_1\in W_1$;
hence, $\Phi_1\subseteq E_1$. Similarly, using (2), we obtain $\Phi_2
\subseteq E_2$. So we must have $\Phi_1\cap \Phi_2=\varnothing$.

Finally, let $\alpha\in \Phi_1$ and $\beta\in \Phi_2$. If $\alpha=\pm 
\beta$ then $\alpha\in E_1\cap E_2=\{\underline{0}\}$, contradiction. 
Hence, $\alpha\pm \beta \neq \underline{0}$. Now assume that $\gamma:=
\alpha \pm \beta \in \Phi$. Since $\Phi=\Phi_1\sqcup \Phi_2$, we have 
either $\gamma \in \Phi_1$ or $\gamma \in \Phi_2$. In the first case, 
$\pm \beta =\gamma-\alpha\in E_1\cap E_2=\{\underline{0}\}$, 
contradiction. The second case leads to a similar contradiction.
Thus, $\alpha\pm \beta \not\in \Phi$.
\end{proof}

\begin{prop} \label{levisub} In the above setting, let $s\in \{1,2\}$. 
Then $(\fg_s,\fh_s)$ is of Cartan--Killing type with respect to the 
subset $\Delta_s':=\{\alpha_i|_{\fh_s}\mid i \in I_s\}\subseteq 
\fh_s^*$; the corresponding structure matrix is $A_s$. 
\end{prop}

\begin{proof} First we show that $\Delta_s'\subseteq \fh_s^*$ is linearly 
independent. Now, we have $A_s=(\alpha_j(h_i))_{i,j\in I_s}$ and $\alpha_j
(h_i)=\alpha_j|_{\fh_s}(h_i)$ for $i,j\in I_s$. Hence, $\Delta_s'$ will be 
linearly independent if $\det(A_s)\neq 0$. But $A$ is a block diagonal
matrix with diagonal blocks $A_1$ and $A_2$. Hence, since $\det(A)\neq 0$, 
we also have $\det(A_s)\neq 0$, as required. 

Let us prove (CK1) for $(\fg_1,\fh_1)$. Let $x \in \fg_1$ be such that 
$[h,x]=0$ for all $h \in \fh_1$. We must show that $x\in \fh_1$. Now, 
since $[\fg_1,\fg_2]=\{0\}$, we also have $[h',x]=0$ for all $h'\in\fh_2$. 
Since $\fh=\fh_1+\fh_2$, we conclude that $x\in \fg_{\underline{0}}=\fh$, 
where the last equality holds by (CK1) for $(\fg,\fh)$. Now write $x=x_1+
x_2$ where $x_1\in \fh_1\subseteq \fg_1$ and $x_2\in \fh_2\subseteq \fg_2$.
Then $x-x_1=x_2\in \fg_1\cap \fg_2=\{0\}$ and so $x=x_1\in \fh_1$. Hence, 
(CK1) holds for $(\fg_1,\fh_1)$. The argument is completely analogous 
for $(\fg_2,\fh_2)$. 

Now consider (CK2). We will do this for $s=1$; the proof for $s=2$ is 
completely analogous. Let $\lambda\in P_{\fh_1}(\fg_1)$. So there is some 
$0\neq x\in \fg_1$ such that $[h,x]=\lambda(h)x$ for all $h\in \fh_1$. Since 
$[\fg_1,\fg_2]=\{0\}$ we have $[h,x]=0$ for all $h \in \fh_2$. Hence, we 
have $x\in \fg_{\tilde{\lambda}}$ where $\tilde{\lambda} \in \fh^*$ is 
defined by $\tilde{\lambda}|_{\fh_1}:=\lambda$ and $\tilde{\lambda}|_{\fh_2}
:=\underline{0}$. Since $x\neq 0$, this means that $\tilde{\lambda}\in 
P_\fh(\fg)$.  By (CK2) for $(\fg,\fh)$ we can write $\tilde{\lambda}=
\sum_{i \in I} n_i\alpha_i$ where either $n_i\in\Z_{\geq 0}$ for all 
$i\in I$, or $n_i\in \Z_{\leq 0}$ for all $i\in I$. Now, if $i\in I_2$, 
then $\alpha_i|_{\fh_1}=\underline{0}$, since $\alpha_i(h_j)=a_{ji}=0$ 
for all $j\in I_1$. Hence, we have 
\[ \lambda=\tilde{\lambda}|_{\fh_1}=\sum_{i \in I} n_i \alpha_i|_{\fh_1}
=\sum_{i \in I_1} n_i \alpha_i|_{\fh_1},\]
as required. Finally, consider (CK3). We have $[h,e_i]=\alpha_i|_{\fh_s}
(h)e_i$ for all $h\in \fh_s$. So $e_i\in \fg_s$ belongs to the 
$\fh_s$-weight space corresponding to $\alpha_i|_{\fh_s}\in \fh_s^*$. 
Similarly, $f_i\in \fg_s$ belongs to the $\fh_s$-weight space corresponding 
to $-\alpha_i|_{\fh_s}\in \fh_s^*$. Since $\fh_s=\langle h_i\mid i \in I_s
\rangle$ and $[e_i,f_i]=h_i$, it follows that (CK3) holds for 
$(\fg_s,\fh_s)$.
\end{proof}

\begin{xca} \label{levisub3} In the above setting, show that 
$\{\alpha|_{\fh_s}\mid \alpha\in \Phi_s\}$ is the root system of 
$\fg_s$; furthermore, $\alpha|_{\fh_s}\neq \beta|_{\fh_s}$ for 
$\alpha\neq \beta$ in $\Phi_s$.
\end{xca}

\begin{rem} \label{cindec5} Assume that $A$ is decomposable. As mentioned 
in Definition~\ref{cindec1}, there is a finite index set $S$ such that
$I=\bigsqcup_{s\in S} I_s$ (disjoint union), where $I_s\neq \varnothing$ 
for all $s\in S$, and $A$ is a block diagonal matrix with indecomposable 
diagonal blocks $A_s=(a_{ij})_{i,j\in I_s}$ for all $s\in S$. For each
$s\in S$, let $\Phi_s$ be the set of all roots $\alpha\in \Phi$ that can 
be expressed as linear combinations of $\{\alpha_i\mid i \in I_s\}$. Then
an iterated application of Proposition~\ref{levisub2} shows that $\Phi$
is the disjoint union of the sets $\Phi_s$ ($s\in S$). Also iterating
Remark~\ref{mainideal3}, we obtain a direct sum decomposition $\cL=
\bigoplus_{s\in S} \cL_s$, where each $\cL_s$ is an ideal and $[\cL_s,
\cL_{s'}]=\{0\}$ for all $s\neq s'$. We have abelian subalgebras 
$\fh_s:=\langle h_i \mid i \in I_s \rangle_\C \subseteq \cL_s$ for all 
$s\in S$. Finally, $(\cL_s,\fh_s)$ is of Cartan--Killing type with 
structure matrix~$A_s$ and, hence, $\cL_s$ is a simple Lie algebra; 
see Proposition~\ref{levisub} and Theorem~\ref{cindec4}. In this way, 
the study of Lie algebras of Cartan--Killing type is reduced to the 
case where the structure matrix $A$ is indecomposable. 
\end{rem}



\begin{defn} \label{indecA} In the above setting, the various matrices 
$A_s$ (for $s\in S$) will be called the \nm{indecomposable blocks} 
of~$A$. Furthermore, the subsets $\Phi_s\subseteq \Phi$ (for $s \in S$) 
will be called the \nm{indecomposable subsystems} of~$\Phi$. Thus, every 
root $\alpha\in \Phi$ belongs to a unique indecomposable subsystem of 
$\Phi$, and the ``type'' of that subsystem (or of~$A_s$) will be one of 
the ``\nmi{Dynkin diagrams}{Dynkin diagram}'' in Table~\ref{Mdynkintbl} 
(p.~\pageref{Mdynkintbl}).
\end{defn}

\section{Classical Lie algebras revisited} \label{sec1a5}

We return to the \nmi{classical Lie algebras}{classical Lie algebra} 
in Section~\ref{sec05}. Recall that 
\[\cL:=\gom_n(Q_n,\C):=\{A\in M_n(\C)\mid A^{\text{tr}} Q_n+Q_nA=0\} 
\subseteq \gl_n(\C)\]
where $Q_n\in M_n(\C)$ is fixed such that $\det(Q_n)\neq 0$ and 
$Q_n^{\text{tr}}=\epsilon Q_n$, $\epsilon\in \{\pm 1\}$.
We assume throughout that $n\geq 3$. Then we have already seen in 
Proposition~\ref{slnss} that $\fg=\gom_n(Q_n,\C)$ is semisimple. Our aim 
is to show that $\fg$ is simple (with the exception in 
Exercise~\ref{xcaclassic3}(c)). For this purpose, we make a specific
choice of~$Q_n$, as in Section~\ref{sec05}:
\[\renewcommand{\arraystretch}{0.8} Q_n=\left(
\begin{array}{c@{\hspace{5pt}}c@{\hspace{5pt}}c@{\hspace{5pt}}c} 0 
& \cdots & 0 & \delta_n \\ \vdots &\dddots & \dddots & 0 \\ 0 & 
\delta_2 &\dddots & \vdots \\ \delta_1 & 0 & \cdots & 0 \end{array} 
\right)\in M_n(\C),\]
where $\delta_1,\ldots,\delta_n\in \{\pm 1\}$ are such that 
$\delta_i\delta_{n+1-i}=\epsilon$ for all~$i$. 

Let $\fh$ be the subspace of diagonal matrices in $\cL$. Let $m\geq 1$ be 
such that $n=2m+1$ (if $n$ is odd) or $n=2m$ (if $n$ is 
even). By the explicit description of $\fh$ in Remark~\ref{classic3b}, 
we have $\dim \fh=m$ and $\fh=\{h(x_1,\ldots,x_m)\mid x_i \in \C\}$, where 
\[h(x_1,\ldots,x_m){:=}\left\{\begin{array}{c@{\hspace{3pt}}l@{\hspace{1pt}}} 
\mbox{diag}(x_1,\ldots,x_m, 0,-x_m,\ldots,-x_1) &\mbox{if $n=2m{+}1$}, \\ 
\mbox{diag}(x_1,\ldots,x_m,-x_m,\ldots, -x_1) \;\;\; &  \mbox{if $n=2m$}. 
\end{array}\right.\]
Furthermore, by Remark~\ref{cartansubgln}, we have $C_\cL(\fh)=\fh$ and 
$\cL$ is $\fh$-diagonalisable. Thus, we have a weight space decomposition
\[ \cL=\fh\oplus \bigoplus_{\alpha\in\Phi} \cL_\alpha\qquad\mbox{where}\qquad
\fh=\cL_{\underline{0}} \;\mbox{ and } \;\Phi\subseteq \fh^*\setminus
\{\underline{0}\}.\]
In order to determine $\Phi$, we use the elements 
\[A_{ij}:= \delta_i E_{ij} -\delta_{j} E_{n+1-j,n+1-i} \in \gom_n(Q_n,\C)\]
for all $1\leq i,j\leq n$, where $E_{ij}$ denotes the matrix 
with~$1$ as its $(i,j)$-entry and zeroes elsewhere. (See 
Proposition~\ref{classic4}.) If $x=\mbox{diag}(x_1,\ldots,x_n)\in \fh$, 
we write $\varepsilon_l(x)=x_l$ for $1\leq l\leq n$; this defines a 
linear map $\varepsilon_l\colon \fh\rightarrow \C$. Note that
$\varepsilon_l+\varepsilon_{n+1-l}=0$ for $1\leq l\leq n$.

\begin{lem} \label{sec161} We have $[x,A_{ij}]=\bigl(\varepsilon_i(x)-
\varepsilon_j(x)\bigr)A_{ij}$ for all $x\in \fh$.
\end{lem}

\begin{proof} If $x=\mbox{diag}(x_1,\ldots,x_n)$, then   
$[x,E_{ij}]=(x_i-x_j) E_{ij}$ and so 
\begin{align*}
[x,A_{ij}]&=\delta_i [x,E_{ij}]-\delta_j [x,E_{n+1-j,n+1-i}]\\&=
\delta_i(x_i-x_j)E_{ij}-\delta_j(x_{n+1-j}-x_{n+1-i})E_{n+1-j,n+1-i}.
\end{align*}
But, since $x\in \fh$, we have $x_{n+1-l}=-x_l$ for $1\leq l \leq n$ and so 
$[x,A_{ij}]=(x_i-x_j)(\delta_iE_{ij}-\delta_jE_{n+1-j,n+1-i})=(x_i-x_j)
A_{ij}$.
\end{proof}

\begin{rem} \label{commaij} Later on, we shall also need to know at least 
some Lie brackets among the elements $A_{ij}$. A straightforward 
computation yields the following formulae. If $i+j\neq n+1$, then 
\[ [A_{ij},A_{ji}]=\delta_i\delta_j(E_{ii}-E_{jj})+\delta_j\delta_i
(E_{n+1-j,n+1-j}-E_{n+1-i,n+1-i});\]
note that this is a diagonal matrix in $\fh$. Furthermore, a particular 
situation occurs when $i+j=n+1$ and $\epsilon=-1$. Then
\[A_{ij}=2\delta_iE_{ij}\qquad\mbox{and}\qquad
[A_{ij},A_{ji}]=4(E_{jj}-E_{ii})\in \fh.\]
\end{rem}

\begin{lem} \label{sec162} Recall that $m\geq 1$ is such that $n=2m+1$ 
or~$n=2m$. 
\begin{itemize}
\item[{\rm (a)}] In all cases, $\{\pm \varepsilon_i\pm \varepsilon_j \mid 
1\leq i,j \leq m,i\neq j\}\subseteq \Phi$. This subset contains 
precisely $2m(m-1)$ distinct elements.
\item[{\rm (b)}] $\{\pm \varepsilon_i\mid 1\leq i\leq m\} \subseteq
\Phi$ if $n=2m+1$ is odd and $Q_n^{\operatorname{tr}}=Q_n$.
\item[{\rm (c)}] $\{\pm 2\varepsilon_i\mid 1\leq i\leq m\}\subseteq 
\Phi$ if $n=2m$ is even and $Q_n^{\operatorname{tr}}=-Q_n$.
\end{itemize}
\end{lem}

\begin{proof} (a) Let $1\leq i,j\leq m$, $i\neq j$. Then
Lemma~\ref{sec161} shows that $\varepsilon_i-\varepsilon_j\in\Phi$, 
with $A_{ij}$ as a corresponding eigenvector. (We have $A_{ij}\neq 0$ in 
this case.) Now set $l:=n+1-j$. Then $l\neq i$ and so Lemma~\ref{sec161} 
also shows that $\varepsilon_i - \varepsilon_l\in\Phi$. (Note that, 
again, $A_{il}\neq 0$.) But $\varepsilon_l=\varepsilon_{n+1-j}=
-\varepsilon_j$ and so $\varepsilon_i+\varepsilon_j\in\Phi$. 
Similarly, let $k:=n+1-i$; then $k\neq j$ and so $\varepsilon_k-
\varepsilon_j \in \Phi$. But $\varepsilon_k=\varepsilon_{n+1-i}=
-\varepsilon_i$ and so $-\varepsilon_i- \varepsilon_j\in\Phi$. Since 
$\{\varepsilon_1,\ldots, \varepsilon_m\}\subseteq \fh^*$ are linearly 
independent, the functions $\pm \varepsilon_i\pm \varepsilon_j\in \fh^*$ 
($1\leq i<j\leq m$) are all distinct. So we have precisely $2m(m-1)$ 
such functions.

(b) Let $1\leq i\leq m$. Then $[x,A_{i,m+1}]=(x_i-x_{m+1})A_{i,m+1}$ for
all $x\in \fh$. But $x_{m+1}=-x_{n+1-(m+1)}=-x_{m+1}$ and so $x_{m+1}=0$.
Hence, we have $[x,A_{i,m+1}]=x_iA_{i,m+1}=\varepsilon_i(x)A_{i,m+1}$ for all
$x\in \fh$. So $\varepsilon_i\in\Phi$ (since $A_{i,m+1}\neq 0$). Similarly,
we see that $[x,A_{m+1,i}]=-\varepsilon_i(x)A_{m+1,i}$ for all $x\in \fh$.
Hence, $-\varepsilon_i\in\Phi$.

(c) Let $1\leq i\leq m$ and $x\in \fh$. Since $x_{2m+1-i}=-x_{i}$, we have 
$[x,A_{i,2m+1-i}]=(x_i-x_{2m+1-i})A_{i,2m+1-i}=2\varepsilon_i(x)A_{i,2m+1-
i}$. Since $Q_n^{\text{tr}}=-Q_n$, we have $\delta_i=-\delta_{2m+1-i}$ 
and so $A_{i,2m+1-i}\neq 0$. This shows that $2\varepsilon_i\in\Phi$. 
Similarly, we see that $[x,A_{2m+1-i,i}]=-2\varepsilon_i(x)A_{2m+1-i,i}$
for all $x\in \fh$. Hence, $-2\varepsilon_i\in\Phi$.
\end{proof}

\begin{prop} \label{sec163} Let $\fh\subseteq \cL=\gom_n(Q_n,\C)$ as above. 
\begin{itemize}
\item[{\rm (a)}] If $Q_n^{\operatorname{tr}}=Q_n$ and $n=2m$ is even, then 
we have $|\Phi|=2(m^2-m)$ and $\Phi=\{\pm \varepsilon_i\pm \varepsilon_j
\mid 1\leq i,j\leq m, i\neq j\}$. 
\item[{\rm (b)}] If $Q_n^{\operatorname{tr}}=Q_n$ and $n=2m+1$ is odd,
then we have $|\Phi|=2m^2$ and $\Phi=\{\pm \varepsilon_i\pm 
\varepsilon_j,\pm\varepsilon_i\mid 1\leq i,j\leq m, i\neq j\}$.
\item[{\rm (c)}] If $Q_n^{\operatorname{tr}}=-Q_n$, then $n=2m$ is 
necessarily even, we have $|\Phi|=2m^2$ and $\Phi=\{\pm \varepsilon_i\pm 
\varepsilon_j,\pm 2\varepsilon_i \mid 1\leq i,j \leq m, i\neq j\}$.
\end{itemize}
\end{prop}

\begin{proof} By Lemma~\ref{sec162}, $|\Phi|\geq 2m^2-2m$ (if $n=2m$ and 
$Q_n^{\operatorname{tr}}=Q_n$) and $|\Phi|\geq 2m^2$ (otherwise). Since
$\dim \fh=m$, this implies that $\dim \cL\geq \dim \fh+|\Phi|\geq 2m^2-m$ (if 
$n=2m$ and $Q_n^{\operatorname{tr}}=Q_n$) and $\dim \cL\geq 2m^2+m$ 
(otherwise). Combining this with the formulae in Remark~\ref{classic3b}, we 
conclude that equality holds everywhere. In particular, $\Phi$ is given by
the functions described in Lemma~\ref{sec162}. In (c), note that 
$Q_n^{\text{tr}}=-Q_n$ implies that $n$ must be even.
\end{proof}

\begin{rem} \label{sec164} In all three cases in Proposition~\ref{sec163}, 
we have $\Phi':=\{\varepsilon_i- \varepsilon_j \mid 1\leq i,j\leq m, 
i\neq j\} \subseteq \Phi$, which is like the set of roots of $\slm_m(\C)$ 
in Example~\ref{cartsln}. We reverse the notation there\footnote{The reason
for this notational reversion is to maintain consistence with the labelling 
of the Dynkin diagrams in Table~\ref{Mdynkintbl} (see 
p.~\pageref{Mdynkintbl}); see also Remark~\ref{trianBCD} below.} and set 
\[\alpha_i:=\varepsilon_{m+1-i}-\varepsilon_{m+2-i}\qquad \mbox{for $2\leq 
i\leq m$}.\]
Thus, $\alpha_m=\varepsilon_1-\varepsilon_2,\;\alpha_{m-1}=\varepsilon_2-
\varepsilon_3,\;\ldots,\;\alpha_2=\varepsilon_{m-1}-\varepsilon_m$; or
$\alpha_{m+2-i}=\varepsilon_{i-1}-\varepsilon_{i}$. For $1 \leq i<j\leq m$, 
we obtain: 
\[\alpha_{i+1}+\alpha_{i+2}+\ldots +\alpha_j=\varepsilon_{m+1-j}-
\varepsilon_{m+1-i}\]
and so $\Phi'=\{\pm (\alpha_{i+1}+\alpha_{i+2}+\ldots +\alpha_{j})\mid 
1\leq i<j \leq m\}$. Furthermore, in all three cases, we have 
$\Phi'':=\{\pm (\varepsilon_i+\varepsilon_j) \mid 1\leq i<j\leq m\} 
\subseteq \Phi$. 
We will now try to obtain convenient descriptions for $\Phi''$. 

$\bullet$ In case (a), $\Phi=\Phi'\cup \Phi''$. If we also set 
$\alpha_1:=\varepsilon_{m-1}+\varepsilon_m$, then $\alpha_1,\alpha_2,\ldots,
\alpha_{m}$ are linearly independent. For $1\leq i<j\leq m$, we have 
\[\alpha_2+\ldots+\alpha_i=\varepsilon_{m+1-i}-\varepsilon_m,\quad
\alpha_3+\ldots+\alpha_j=\varepsilon_{m+1-j}-\varepsilon_{m-1},\]
and so $(\alpha_1+\alpha_2+\ldots+\alpha_i)+(\alpha_3+\alpha_4 +\ldots +
\alpha_j)=\varepsilon_{m+1-i}+\varepsilon_{m+1-j}$. (Note that $m\geq 2$ 
since $n\geq 3$ is even.) Hence, these expressions (and their negatives) 
describe all elements of $\Phi''$.

$\bullet$ In case (b), $\Phi=\Phi'\cup \Phi''\cup 
\{\pm \varepsilon_i\mid 1\leq i \leq m\}$. If we also set $\alpha_1:=
\varepsilon_m$, then $\alpha_1,\alpha_2, \ldots,\alpha_m$ are
linearly independent. We have 
\[\alpha_1+(\alpha_2+\ldots+\alpha_i)=\varepsilon_m+(\varepsilon_{m+1-i}-
\varepsilon_m)=\varepsilon_{m+1-i} \]
for $1\leq i \leq m$. Furthermore, for $1\leq i<j\leq m$, we obtain
\begin{align*}
2(\alpha_1+&\alpha_2+\ldots+\alpha_i)+ \alpha_{i+1}+ \alpha_{i+2} +\ldots +
\alpha_j\\ &=2\varepsilon_{m+1-i}+(\varepsilon_{m+1-j}-\varepsilon_{m+1-i})
=\varepsilon_{m+1-i}+\varepsilon_{m+1-j}.
\end{align*}
Hence, the above expressions describe all elements of~$\Phi''$.

$\bullet$ In case (c), $\Phi=\Phi'\cup \Phi''\cup \{\pm 2\varepsilon_i\mid 
1 \leq i\leq m\}$. If we also set $\alpha_1:=2\varepsilon_m$, then 
$\alpha_1, \alpha_2, \ldots,\alpha_m$ are linearly independent. We have 
\[\alpha_1+2(\alpha_2+\ldots +\alpha_i)=2\varepsilon_m+2(\varepsilon_{m+1-i}
-\varepsilon_m)=2\varepsilon_{m+1-i}\]
for $1\leq i\leq m$. Furthermore, for $1 \leq i <j \leq m$, we obtain
\begin{align*}
\alpha_1+&2(\alpha_2+ \ldots+ \alpha_i)+ \alpha_{i+1}+\alpha_{i+2} +
\ldots +\alpha_j\\ &=2\varepsilon_{m+1-i}+(\varepsilon_{m+1-j}-
\varepsilon_{m+1-i})=\varepsilon_{m+1-i}+\varepsilon_{m+1-j}.
\end{align*}
Hence, again, the above expressions describe all elements of $\Phi''$.
\end{rem}

\begin{cor} \label{sec165} Let $\cL=\gom_n(Q_n,\C)$. Then, with notation 
as in Remark~\ref{sec164}, $\Delta:=\{\alpha_1,\ldots,\alpha_m\}$ is a
basis of $\fh^*$ and each $\alpha\in\Phi$ can be written as $\alpha=
\pm \sum_{1\leq i\leq m} n_i\alpha_i$ with $n_i\in\{0,1,2\}$ for all~$i$. 
\end{cor}

\begin{proof} We already noted that $\{\alpha_1,\ldots,\alpha_m\}$ is
linearly independent. The required expressions of $\alpha$ are explicitly 
given above.
\end{proof}

\begin{rem} \label{trianBCD} Let $x\in \cL=\gom_n(Q_n,\C)$ and write
$x=h+n^++n^-$ as in Corollary~\ref{triangom}. Then one easily checks that
our choice of $\alpha_1,\ldots,\alpha_m$ in Remark~\ref{sec164} is such 
that $n^{\pm}\in\sum_\alpha \cL_{\pm \alpha}$ where the sum runs over all
$\alpha \in \Phi$ such that $\alpha=\sum_{1\leq i\leq m} n_i \alpha_i$ 
with $n_i\geq 0$.
\end{rem}

\begin{table}[htbp] \caption{Structure matrices $A$ 
for the Lie algebras $\cL=\gom_n(Q_n,\C)$} \label{cartanBCD} 
\begin{center} {$\renewcommand{\arraystretch}{0.8} 
\left(\begin{array}{r@{\hspace{4pt}}r@{\hspace{4pt}}r@{\hspace{4pt}}
r@{\hspace{4pt}}r@{\hspace{4pt}}r@{\hspace{4pt}}r@{\hspace{4pt}}r} 2 & 0 
& -1 &&&&&\\ 0 & 2 & -1 & &&&&\\ -1 & -1 & 2 &-1&&&\\ && -1 & 2 & -1 &&& \\ 
&&&\ddots &\ddots &\ddots && \\ && & & -1 & 2 & -1& \\ &&&&& -1 & 2 &
\end{array}\right)\qquad \mbox{($Q_n^{\text{tr}}=Q_n$ and $n=2m$)},
\quad\;$}\end{center}
\begin{center} {$\renewcommand{\arraystretch}{0.9} \begin{array}{c}
\left(\begin{array}{r@{\hspace{4pt}}r@{\hspace{4pt}}r@{\hspace{4pt}}
r@{\hspace{4pt}}r@{\hspace{4pt}}r} 2 & -2 &&&&\\ -1 & 2 & -1 &&&\\ & -1 & 
2 & -1 & & \\ & & \ddots & \ddots & \ddots & \\ & & & -1 & 2 & -1 \\ 
&&&& -1 & 2 \end{array}\right) \\\\\text{($Q_n^{\text{tr}}=Q_n$ 
and $n=2m+1$)}\end{array} \;\mbox{ and }\;\begin{array}{c}
\left(\begin{array}{r@{\hspace{4pt}}r@{\hspace{4pt}}r@{\hspace{4pt}}
r@{\hspace{4pt}}r@{\hspace{4pt}}r} 2 & -1 &&&&\\ -2 & 2 & -1 &&&\\ & -1 & 
2 & -1 & & \\ & & \ddots & \ddots & \ddots & \\ & & & -1 & 2 & -1 \\ 
&&&& -1 & 2 \end{array}\right)\\\\\text{($Q_n^{\text{tr}}=-Q_n$ 
and $n=2m$)}\end{array}$.}
\end{center}
\end{table}

\begin{prop} \label{CKbcd} Let $\cL=\gom_n(Q_n,\C)$ and $\fh\subseteq \cL$ 
be as above; assume that $n\geq 3$ and write $n=2m+1$ or $n=2m$, where $m
\geq 1$. Then $(\cL,\fh)$ is of Cartan--Killing type with respect to 
$\Delta=\{\alpha_1, \ldots, \alpha_m\}\subseteq \fh^*$, as defined in 
Remark~\ref{sec164}; the structure matrix $A$ is given in 
Table~\ref{cartanBCD}. (Each of those matrices has size $m\times m$.) 
\end{prop}

\begin{proof} We already noted that $\cL$ is $\fh$-diagonalisable and 
$C_\cL(\fh)=\fh$; we also have $\alpha_i\in P_\fh(\fg)$ for $1\leq i \leq m$.
Hence, (CK1) in Definition~\ref{defTD} holds. Furthermore, 
(CK2) holds by Corollary~\ref{sec165}. 

It remains to verify (CK3) and to specify $e_i\in \cL_{\alpha_i}$ and 
$f_i \in \cL_{-\alpha_i}$ such that $\alpha_i(h_i)=2$, where 
$h_i:=[e_i,f_i]\in \fh$. For $2\leq i \leq m$, we 
have $\alpha_i=\varepsilon_{m+1-i}-\varepsilon_{m+2-i}$, or $\alpha_{m+2-i}
=\varepsilon_{i-1}-\varepsilon_i$. So Lemma~\ref{sec161} shows that
\[ e_{m+2-i}:=\delta_{i-1}A_{i-1,i}\in \cL_{\alpha_{m+2-i}}, \quad
f_{m+2-i}:=\delta_iA_{i,i-1}\in \cL_{-\alpha_{m+2-i}}.\]
Using the formulae in Remark~\ref{commaij}, we find that 
\[ h_{m+2-i}:=[e_{m+2-i},f_{m+2-i}]=h(0,\ldots,0,1,-1,0,\ldots,0)\in \fh,\]
where the entry $1$ is at the $(i-1)$-th position and $-1$ is at the 
$i$-th position. Hence, $\alpha_i(h_i)=2$ for $2\leq i\leq m$, as required.

If $Q_n^{\text{tr}}=Q_n$ and $n=2m$, then we have $\alpha_1=\varepsilon_{m-1}
+\varepsilon_m$. As in the proof of Lemma~\ref{sec162}(a), we see that
\[ e_1:=\delta_{m-1}A_{m-1,m+1} \in \cL_{\alpha_1}\quad \mbox{and}\quad 
f_1:=\delta_{m+1}A_{m+1,m-1} \in \cL_{-\alpha_1}.\]
Using Remark~\ref{commaij}, we find that $h_1:=[e_1,f_1]=
h(0,\ldots, 0,1,1)\in \fh$ and $\alpha_1(h_1)=2$, as required.
If $Q_n^{\text{tr}}=Q_n$ and $n=2m+1$, then we have $\alpha_1=\varepsilon_m$.
As in the proof of Lemma~\ref{sec162}(b), we see that 
\[ e_1:=\delta_mA_{m,m+1}\in \cL_{\alpha_1} \quad \mbox{and}\quad f_1:=
2\delta_{m+1}A_{m+1,m}\in \cL_{-\alpha_1}.\]
Now $h_1:=[e_1,f_1]=h(0,\ldots,0,2)\in \fh$ and $\alpha_1(h_1)=2$, as 
required. Finally, if $Q_n^{\text{tr}}=-Q_n$ and $n=2m$, then we have 
$\alpha_1=2\varepsilon_m$. As in the proof of Lemma~\ref{sec162}(c), we 
see that 
\[\textstyle e_1:=\frac{1}{2}\delta_mA_{m,m+1} \in \cL_{\alpha_1}, \qquad 
f_1:=\frac{1}{2}\delta_{m+1}A_{m+1,m}\in \cL_{-\alpha_1}.\]
By Remark~\ref{commaij}, we actually have $e_1=E_{m,m+1}$ and $f_1=
E_{m+1,m}$ in this case; furthermore, $h_1:=[e_1,f_1]=h(0,\ldots,
0,1)\in \fh$ and $\alpha_1(h_1)=2$, as required.

In all cases, we see that $\fh=\langle h_1,\ldots,h_m\rangle_\C$ and so 
(CK3) holds. Finally, $A$ is obtained by evaluating $\alpha_j(h_i)$ for all
$i,j$.
\end{proof}

\begin{rem} \label{CKbcda} The above proof shows that, in each case, 
there is a system of Chevalley generators $\{e_i,f_i\mid 1\leq i 
\leq m\}$ for $\cL$ such that all $e_i,f_i$ are given by matrices with 
entries in~$\Z$. 
\end{rem}

\begin{thm} \label{sec169} Recall that $n\geq 3$. If 
$Q_n^{\operatorname{tr}}=Q_n$ and $n$ is even, also assume that 
$n\geq 6$. Then $\cL=\gom_n(Q_n,\C)$ is a simple Lie algebra. 
{\rm (Note that, by Exercise~\ref{xcaclassic3}(c), we really do have to 
exclude the case where $n=4$ and $Q_4=Q_4^{\text{tr}}$.) }
\end{thm}

\begin{proof} By Proposition~\ref{CKbcd}, $(\cL,\fh)$ is of Cartan--Killing
type with respect to $\Delta=\{\alpha_1,\ldots,\alpha_m\}$. We can now
use Remark~\ref{cindec1a} to show that $\cL$ is simple (exactly as
for $\cL=\slm_n(\C)$). Just note that, for the specified conditions on~$n$, 
each of the structure matrices in Table~\ref{cartanBCD} is indecomposable. 
(Draw the corresponding graph as in Remark~\ref{cindec1a}.)
\end{proof}


Finally, we determine the Weyl group $W$ of $\cL=\gom_n(Q_n,\C)$. With the 
above description of $\Delta=\{\alpha_1,\ldots,\alpha_m\}$, we have 
$W=\langle s_1,\ldots,s_m\rangle\subseteq \GL(\fh^*)$, where
$s_i(\lambda)=\lambda-\lambda(h_i)\alpha_i$ for all $\lambda\in \fh^*$.

First we consider the cases (b) and (c) in Proposition~\ref{sec163}.
If $Q_n^{\text{tr}}=Q_n$ and $n=2m+1$,  we set $d=1$; if $Q_n^{\text{tr}}
=-Q_n$ and $n=2m$, we set $d=2$. It will be convenient to define
$v_1:=\varepsilon_m$ and $v_i:=\varepsilon_{m+1-i}$ for $i\geq 2$. Thus,
$\{v_1,\ldots,v_m\}$ is a basis of $\fh^*$ such that $\alpha_1=dv_1$ and 
$\alpha_i=v_{i}-v_{i-1}$ for $2\leq i\leq m$ (see Remark~\ref{sec164}). 
We have
\begin{center}
$v_1=\frac{1}{d}\alpha_1 \quad\mbox{and}\quad v_i=\alpha_i+\alpha_{i-1}+
\ldots+\alpha_2+\frac{1}{d}\alpha_1 \quad\mbox{for $2\leq i\leq m$}$.
\end{center}
We compute the effect of $s_i\in W$ on these basis vectors, using
the formulae in Remark~\ref{explicit} and the knowledge of the structure
matrix~$A$. For $i\in\{2,\ldots,m\}$, we obtain 
\[s_i(v_i)=v_{i-1},\quad s_i(v_{i-1})=v_i\quad\mbox{and}\quad s_i(v_j)=
v_j \mbox{ if $j\not\in\{i-1,i\}$}.\]
So $s_i$ permutes the basis vectors $v_j$. We also find that $s_1(v_1)=
-v_1$ and $s_1(v_j)=v_j$ for $j\geq 2$. (Details of these computations
are left to the reader.) Thus, the effect of an arbitrary element of $W$ 
on the list of basis vectors $v_1,\ldots,v_m$ will be a ``permutation 
with signs''.

\begin{defn} \label{hyperoct} Consider the symmetric group $\mbox{Sym}(X)$
where $X=\{\pm 1,\ldots,\pm m\}$. An element $\pi\in \mbox{Sym}(X)$ is
called a \nm{signed permutation} if $\pi(-i)=-\pi(i)$ for $1\leq i\leq m$. 
Such a $\pi$ is uniquely determined by its values on $1,\ldots,m$; so we 
can simply write $\pi$ as 
\[ \renewcommand{\arraystretch}{0.8}
\pi=\left(\begin{array}{cccc} 1 & 2 & \ldots & m\\
\pi(1) & \pi(2) & \ldots & \pi(m)\end{array}\right)\quad\mbox{where 
$\pi(i) \in \{\pm 1,\ldots,\pm m\}$}.\]
Let $\fH_m\subseteq \mbox{Sym}(X)$ be the subset consisting of all signed
permutations. One immediately checks that $\fH_m$ is a subgroup of 
$\mbox{Sym}(X)$, called the \nm{hyperoctahedral group} (of degree $m$). 
\end{defn}

\begin{xca} \label{xcahyperoct} (a) Show that $|\fH_m|=2^mm!$ and that 
$\fH_m$ is generated by the following signed permutations $\tau$ and 
$\sigma_1,\ldots,\sigma_{m-1}$:
\begin{align*} \renewcommand{\arraystretch}{0.8}
\tau & :=\left(\begin{array}{c@{\hspace{8pt}}c@{\hspace{8pt}}
c@{\hspace{8pt}}c@{\hspace{8pt}}c} 1 & 2 & 3 & \ldots & m\\ -1 & 2 & 3 & 
\ldots & m\end{array}\right)\\
\sigma_i&:=\left(\begin{array}{c@{\hspace{8pt}}c@{\hspace{8pt}}
c@{\hspace{8pt}}c@{\hspace{8pt}}c@{\hspace{8pt}}c@{\hspace{8pt}}cc} 1 & 
\ldots & i{-}1 & i & i{+}1 & i{+}2 & \ldots & m \\ 1 & \ldots & i{-}1 
& i{+}1 & i & i{+}2 & \ldots & m \end{array} \right)
\end{align*}
for $1\leq i\leq m-1$. Show that $\langle \sigma_1,\ldots,\sigma_{m-1}
\rangle \cong \fS_m$.

(b) Let $m\geq 2$. For $\pi\in\fH_m$ we denote by $\nu_\pi$ the number 
of ``sign changes'', that is, the number of $i\in\{1,\ldots,m\}$ such that 
$\pi(i)<0$. Show that $\fH_m':=\{\pi\in\fH_m\mid \nu_\pi \mbox{ even}\}$ 
is a subgroup of $\fH_m$ of index~$2$. Show that $\fH_m'$ is generated 
by $\tau\circ \sigma_1\circ \tau$ and $\sigma_1,\ldots,\sigma_{m-1}$.
\end{xca}

\begin{prop} \label{sec168} Assume that either $Q_n^{\operatorname{tr}}=Q_n$
and $n=2m+1$, or $Q_n^{\operatorname{tr}}=-Q_n$ and $n=2m$, where $m\geq 1$.
Then $W\cong \fH_m$; in particular, $|W|=2^mm!$.
\end{prop}

\begin{proof} Let $X':=\{\pm v_1,\ldots,\pm v_m\}\subseteq \fh^*$. We have 
seen above that $s_i(X')=X'$ for $1\leq i\leq m$ and so $w(X')=X'$ for all 
$w\in W$. Hence, there is an action of the group $W$ on the set $X'$ via 
\[W\times X'\rightarrow X',\qquad (w,\pm v_i) \mapsto \pm w(v_i).\]
Since $X'$ contains a basis of $\fh^*$, and since $W$ is faithful on 
$\fh^*$, the corresponding group homomorphism $W\rightarrow \mbox{Sym}(X')$ 
is injective. Since $w(-v_i)=-w(v_i)$ for all~$i$, the image of that 
homomorphism is contained in the set of signed permutations of $X'$. Thus, 
identifying $\pm v_i\leftrightarrow \pm i$, we obtain an injective
homomorphism $W\rightarrow \fH_m$. We have also seen that $s_{i+1}$
induces the signed permutation $\sigma_i\in \fH_m$ for $1\leq i\leq m-1$. 
Hence, $\sigma_1, \ldots,\sigma_{m-1}$ belong to the image of 
$W\rightarrow \fH_m$. Furthermore, $\tau \in\fH_m$ is induced by $s_1$. 
So, by Exercise~\ref{xcahyperoct}(a), the map $W \rightarrow\fH_m$ is 
surjective. 
\end{proof}

\begin{rem} \label{weylD} Now consider the case where 
$Q_n^{\operatorname{tr}}=Q_n$ and $n=2m$ ($m\geq 2$), as in
Proposition~\ref{sec163}(a). Then we claim that $W\cong \fH_m'$. To see 
this, it will now be convenient to consider the basis $\{v_1,\ldots,v_m\}$
of $\fh^*$ such that $\alpha_1=v_1+v_2$ and $\alpha_i=v_{i}-v_{i-1}$ for
$2\leq i\leq m$. Thus,

$v_1=\frac{1}{2}(\alpha_1-\alpha_2), \qquad v_2=\frac{1}{2}(\alpha_1+
\alpha_2)$, 

$v_i=\alpha_i+\alpha_{i-1}+ \ldots+\alpha_3+\frac{1}{2}
(\alpha_1+\alpha_2)\;$ for $3\leq i\leq m$.

\noindent We find again that $s_{i+1}$ induces the signed permutation 
$\sigma_i$ on the basis vectors $v_1,\ldots,v_m$. Furthermore, $s_1$ 
induces the signed permutation $\tau\circ \sigma_1\circ \tau$. Hence, by 
Exercise~\ref{xcahyperoct}(b), we conclude that $W\cong\fH_m'$.
\end{rem}

\begin{rem} \label{weylan} Finally, to complete the picture, we also determine
the Weyl group of the Lie algebra $\cL=\slm_n(\C)$, where $n\geq 2$. For
this purpose, we use the inclusion $\cL\subseteq \hat{\cL}=\gl_n(\C)$. Let 
$\hat{\fh}:=\{\mbox{diag}(x_1,\ldots,x_n)\mid x_i\in\C\}\subseteq \hat{\cL}$ 
be the subspace of all diagonal matrices in $\hat{\cL}$. For $1\leq i\leq n$,
let $\hat{\varepsilon}_i\in \hat{\fh}^*$ be the map that sends a diagonal
matrix to its $i$-th diagonal entry. Then $\{\hat{\varepsilon}_1,\ldots,
\hat{\varepsilon}_n\}$ is a basis of $\hat{\fh}^*$. Another basis is given by
$\{\delta,\hat{\alpha}_1,\ldots,\hat{\alpha}_{n-1}\}$ where
\[\delta:=\hat{\varepsilon}_1+\ldots+\hat{\varepsilon}_n\quad\mbox{and}\quad
\hat{\alpha}_i:=\hat{\varepsilon}_i-\hat{\varepsilon}_{i+1} \quad
\mbox{for $1\leq i\leq n-1$}.\]
Now consider the Weyl group $W=\langle s_1,\ldots,s_{n-1}\rangle\subseteq 
\fh^*$ of $\cL$, where $\fh=\ker(\delta)\subseteq \hat{\fh}$. We define a map 
$\pi\colon W\rightarrow \GL(\hat{\fh}^*)$ as follows. Let $w\in W$ and write 
$w(\alpha_j)= \sum_{i} m_{ij}(w)\alpha_i$ with $m_{ij}(w)\in\Z$ for $1\leq 
i,j\leq n-1$. Thus, $M_w=\bigl(m_{ij}(w) \bigr) \in\GL_{n-1}(\C)$ is the 
matrix of $w$ with respect to the basis $\Delta=\{\alpha_1,\ldots, 
\alpha_{n-1}\}\subseteq \fh^*$. Then we define $\hat{w}\in
\GL(\hat{\fh}^*)$ by setting
\[ \hat{w}(\delta):=\delta \quad \mbox{and}\quad \hat{w}(\hat{\alpha}_j)
:=\sum_{1\leq i\leq n-1} m_{ij}(w)\hat{\alpha}_i\quad \mbox{for $1\leq j
\leq n-1$}.\]
Thus, the matrix of $\hat{w}$ with respect to the basis $\{\delta,
\hat{\alpha}_1,\ldots,\hat{\alpha}_{n-1}\}$ of $\hat{\fh}^*$ is a block
diagonal matrix of the following shape:
\[\renewcommand{\arraystretch}{1.1} \left(\begin{array}{c|c} \,\;1\,\; & 
0 \\\hline 0 & M_w \end{array}\right).\]
Now $\pi\colon W\rightarrow \GL(\hat{\fh}^*)$, $w\mapsto \hat{w}$, is an 
injective group homomorphism, and we have $\pi(W)=\langle \hat{s}_1,
\ldots,\hat{s}_{n-1}\rangle$. Since $\delta(h_i)=0$ for all~$i$, we 
see that $\hat{s}_i\colon \hat{\fh}^*\rightarrow\hat{\fh}^*$ is given by 
the formula 
\[\hat{s}_i(\mu)=\mu-\mu(h_i)\hat{\alpha}_i\qquad\mbox{for all $\mu\in 
\hat{\fh}^*$}.\] 
A straightforward computation shows that  
\[\hat{s}_i(\hat{\varepsilon}_i)=\hat{\varepsilon}_{i+1},\quad \hat{s}_i
(\hat{\varepsilon}_{i+1})=\hat{\varepsilon}_i\quad\mbox{and}\quad \hat{s}_i
(\hat{\varepsilon}_j)=\hat{\varepsilon}_j \mbox{ if $j\not\in\{i,i+1\}$}.\]
Thus, the matrix of $\hat{s}_i$ with respect to the basis 
$\{\hat{\varepsilon}_1,\ldots, \hat{\varepsilon}_n\}$ of $\hat{\fh}^*$ is
the permutation matrix corresponding to the transposition in $\fS_n$ that
exchanges $i$ and $i+1$. Since $\fS_n$ is generated by these transpositions,
we conclude that $W\cong \pi(W) \cong \fS_n$. 
\end{rem}

\section{The structure constants $N_{\alpha,\beta}$} \label{sec1a6}

Returning to the general situation, let again $(\cL,\fh)$ be of 
Cartan--Killing type with respect to $\Delta=\{\alpha_i\mid i \in I\}$. 
Let $\Phi\subseteq \fh^*$ be the set of roots of $\cL$ and fix a collection 
of elements
\[\{0\neq e_\alpha \in \cL_\alpha \mid \alpha \in \Phi\}.\]
Then, since $\dim \cL_\alpha=1$ for all $\alpha\in\Phi$, the set 
\[ \{h_i\mid i \in I\}\cup \{e_\alpha\mid \alpha\in\Phi\} \quad
\mbox{is a basis of $\cL$}.\]
If $\alpha,\beta\in\Phi$ and $\alpha+\beta\in\Phi$, then $[\cL_\alpha,
\cL_\beta]\subseteq \cL_{\alpha+\beta}$ and 
\[ [e_\alpha,e_\beta]=N_{\alpha,\beta}e_{\alpha+\beta},\qquad
\mbox{where} \qquad N_{\alpha,\beta}\in\C.\]
The knowledge of the structure constants $N_{\alpha,\beta}$ is, of
course, crucial for doing explicit computations inside $\cL$. Eventually,
one would hope to find purely combinatorial formulae  for $N_{\alpha,
\beta}$ in terms of properties of $\Phi$. In this section, we establish 
some basic properties of the $N_{\alpha,\beta}$. 

It will be convenient to set $N_{\alpha,\beta}:=0$ if 
$\alpha+\beta\not\in\Phi\cup \{\underline{0}\}$.

\begin{rem} \label{propnrs0} Let $\alpha\in\Phi$. By 
Proposition~\ref{weyl5c}, there is a unique $h_\alpha\in [\cL_\alpha,
\cL_{-\alpha}]$ such that $\alpha(h_\alpha)=2$. Now recall that $\Phi=-\Phi$.
We claim that the elements $\{e_\alpha\mid \alpha\in\Phi\}$ can be adjusted 
such that 
\begin{equation*}
[e_\alpha,e_{-\alpha}]=h_\alpha\qquad \mbox{for all 
$\alpha\in\Phi$}.\tag{a}
\end{equation*}
Indeed, we have $\Phi=\Phi^+\sqcup \Phi^-$ (disjoint union), where $\Phi^-=
-\Phi^+$. Let $\alpha\in\Phi^+$. Then $[e_\alpha,e_{-\alpha}]=\xi h_\alpha$ 
for some $0\neq \xi\in \C$. Hence, replacing $e_{-\alpha}$ by a suitable 
scalar multiple if necessary, we can achieve that $[e_\alpha,e_{-\alpha}]=
h_\alpha$. Thus, the desired relation holds for all $\alpha\in \Phi^+$.
Now let $\beta\in\Phi^-$; then $\alpha=-\beta\in\Phi^+$. So 
$[e_\beta,e_{-\beta}]=-[e_\alpha,e_{-\alpha}]=-h_\alpha=h_\beta$, 
where the last equality holds by Proposition~\ref{weyl5c}. So (a) holds 
in general. Now, writing $f_\alpha:=e_{-\alpha}$ we have $[e_\alpha,
f_\alpha]=h_\alpha$, $[h_\alpha,e_\alpha]=\alpha(h_\alpha)e_\alpha=
2e_\alpha$ and $[h_\alpha,f_\alpha]=-2f_\alpha$. Hence, as in
Remark~\ref{astring}, we obtain a $3$-dimensional subalgebra 
\begin{equation*}
\fs_\alpha:=\langle e_\alpha,h_\alpha,f_\alpha\rangle_\C \subseteq \cL
\quad \mbox{such that} \quad \fs_\alpha\cong \slm_2(\C).\tag{b}
\end{equation*} 
Regarding $\cL$ as an $\fs_\alpha$-module, we obtain results completely
analogous to those in Remark~\ref{astring}. Here is a first example.
As in Section~\ref{sec1a3}, let $E:=\langle\alpha_i\mid i\in I\rangle_\R
\subseteq \fh^*$ and $\langle \;,\;\rangle\colon E\times E\rightarrow\R$
be a $W$-invariant scalar product, where $W$ is the Weyl group of $(\cL,
\fh)$. 
\end{rem}

\begin{lem} \label{astring2} Let $\alpha \in \Phi$. Then we have 
\[\lambda(h_\alpha)=2\frac{\langle \alpha,\lambda \rangle}{\langle 
\alpha,\alpha\rangle} \qquad \mbox{for all $\lambda \in E$}.\]
Furthermore, if $\beta\in \Phi$ is such that $\beta\neq \pm \alpha$, then 
$\beta(h_\alpha)=q-p\in\Z$, where $p,q\geq 0$ are defined by the 
condition that
\[ \beta-q\alpha,\quad\ldots,\quad\beta-\alpha,\quad\beta,\quad
\beta+\alpha,\quad\ldots,\quad \beta+p\alpha\]
all belong to $\Phi$, but $\beta+(p+1)\alpha\not\in\Phi$ and 
$\beta-(q+1)\alpha\not\in \Phi$. 
\end{lem}
In analogy to Remark~\ref{astring}, the above sequence of roots is 
called the \nms{$\alpha$-string through $\beta$}{alpha-string through 
beta}. The element $h_\alpha$ is also called a \nm{co-root} of $\cL$.

\begin{proof} We write $\alpha=w(\alpha_i)$, where $w\in W$ and $i
\in I$. Applying $w^{-1}$ to the above sequence of roots and setting
$\beta':=w^{-1}(\beta)$, we see that  
\[ \beta'-q\alpha_i,\quad\ldots,\quad\beta'-\alpha_i,\quad\beta',
\quad\beta'+\alpha_i, \quad\ldots,\quad\beta'+p\alpha_i\]
all belong to $\Phi$. If we had $\beta'+(p+1)\alpha_i\in \Phi$, then 
also $\beta+(p+1)\alpha=w(\beta'+(p+1)\alpha_i)\in \Phi$, contradiction. 
Similarly, we have $\beta'-(q+1)\alpha_i\not\in \Phi$. Hence, the above
sequence is the $\alpha_i$-string through $\beta'$ and so $\beta'(h_i)=q-p$;
see Remark~\ref{astring}(a). Using the $W$-invariance of $\langle\;,\;
\rangle$ and the formula in Remark~\ref{keylem0}, we obtain that 
\[ 2\frac{\langle \alpha,\beta\rangle}{\langle \alpha,\alpha\rangle}=
2\frac{\langle w(\alpha_i),w(\beta')\rangle}{\langle w(\alpha_i),
w(\alpha_i)\rangle}=
2\frac{\langle \alpha_i,\beta'\rangle}{\langle \alpha_i,
\alpha_i\rangle}=\beta'(h_i)=q-p.\]
Furthermore, using $\fs_\alpha=\langle h_\alpha,e_\alpha,f_\alpha
\rangle_\C \subseteq \cL$ as above, one sees that $\beta(h_\alpha)=q-p$, 
exactly as in Remark~\ref{astring}(a) (where $e_\alpha,h_\alpha,
f_\alpha$ play the role of $e_i,h_i,f_i$, respectively). Hence, the 
formula $\lambda(h_\alpha)=2\frac{\langle \alpha,\lambda \rangle}{\langle 
\alpha,\alpha\rangle}$ holds for all $\lambda\in \Phi$ such that 
$\lambda\neq \pm \alpha$. By the definition of $h_\alpha$, it also holds 
for $\lambda=\pm \alpha$. Finally, since $E=\langle \Phi\rangle_\R$, 
it holds in general.
\end{proof}

\begin{lem} \label{astring3} Let $\alpha \in \Phi$ and write 
$\alpha=\sum_{i \in I} n_i\alpha_i$ with $n_i\in \Z$. Then 
$h_\alpha=\sum_{i \in I} n_i^\vee h_i$, where
\[ n_i^\vee=\frac{\langle \alpha_i,\alpha_i \rangle}{\langle \alpha,\alpha
\rangle}n_i\in \Z \qquad \mbox{for all $i\in I$}.\]
\end{lem}

\begin{proof} Given the expression $\alpha=\sum_{i \in I} n_i\alpha_i$,
we obtain 
\[ \frac{2\alpha}{\langle \alpha,\alpha\rangle}=\sum_{i \in I } n_i
\frac{2}{\langle \alpha,\alpha \rangle} \frac{\langle \alpha_i,\alpha_i
\rangle}{\langle \alpha_i,\alpha_i\rangle}\alpha_i=\sum_{i\in I} n_i
\frac{\langle \alpha_i,\alpha_i \rangle}{\langle \alpha,\alpha\rangle}
\frac{2\alpha_i}{\langle \alpha_i,\alpha_i\rangle}.\]
Now let $\lambda \in E$. Using the formula in Lemma~\ref{astring2}, 
we obtain:
\[ \lambda(h_\alpha)=\sum_{i \in I} n_i\frac{\langle \alpha_i,\alpha_i
\rangle}{\langle \alpha,\alpha\rangle}\lambda(h_i)=\lambda\Bigl(\sum_{i 
\in I} n_i\frac{\langle \alpha_i,\alpha_i \rangle}{\langle \alpha,\alpha
\rangle} h_i\Bigr).\]
Since this holds for all $\lambda$, we obtain the desired formula.
The fact that the coefficients $n_i^\vee$ are integers follows from
Exercise~\ref{xcaalgoh}.
\end{proof}

\begin{rem} \label{intstruct} In the following discussion, we assume 
throughout that (a) in Remark~\ref{propnrs0} holds, that is, we have 
$[e_\alpha,e_{-\alpha}]=h_\alpha$ for all $\alpha\in \Phi$. 
This assumption leads to the following summary about the Lie brackets 
in $\cL$. We have:
\begin{alignat*}{2}
[h_i,h_j]&=0,\quad\qquad&&\mbox{for all $i,j\in I$},\\ [h_i,e_\alpha]&=
\alpha(h_i)e_\alpha,\quad &&\mbox{where $\alpha(h_i)\in\Z$},\\
[e_\alpha,e_{-\alpha}]&=h_\alpha \in\langle h_i\mid i\in I\rangle_\Z 
\quad&& \mbox{(see Lemma~\ref{astring3})},\\
[e_\alpha,e_\beta]&=0 \quad &&\mbox{if $\alpha+\beta\not\in\Phi
\cup\{\underline{0}\}$},\\
[e_\alpha,e_\beta]&=N_{\alpha,\beta}e_{\alpha+\beta} \quad &&\mbox{if 
$\alpha+\beta\in\Phi$}.
\end{alignat*}
Since $\{h_i\mid i \in I\}\cup \{e_\alpha\mid \alpha\in\Phi\}$ is a basis 
of $\cL$, the above formulae completely determine the multiplication
in~$\cL$. At this point, the only unknown quantities in those formulae are 
the constants~$N_{\alpha,\beta}$. 
\end{rem}

\begin{lem} \label{propnrs1} If $\gamma_1,\gamma_2,\gamma_3\in \Phi$ 
are such that $\gamma_1+\gamma_2+\gamma_3=\underline{0}$, then
\[ N_{\gamma_1,\gamma_2}=-N_{\gamma_2,\gamma_1} \qquad \mbox{and}\qquad
\frac{N_{\gamma_1,\gamma_2}}{\langle \gamma_3,\gamma_3\rangle}=
\frac{N_{\gamma_2,\gamma_3}}{\langle \gamma_1,\gamma_1\rangle}=
\frac{N_{\gamma_3,\gamma_1}}{\langle \gamma_2,\gamma_2\rangle}.\]
\end{lem}

\begin{proof} Since $\gamma_1+\gamma_2=-\gamma_3\in \Phi$, the 
anti-symmetry of $[\;,\;]$ immediately yields $N_{\gamma_1,\gamma_2}
=-N_{\gamma_2,\gamma_1}$. Now, since also $\gamma_2+\gamma_3=-\gamma_1
\in \Phi$, we have $[e_{\gamma_2},e_{\gamma_3}]=N_{\gamma_2,\gamma_3}
e_{\gamma_2+\gamma_3}=N_{\gamma_2,\gamma_3}e_{-\gamma_1}$ and so 
\[ [e_{\gamma_1},[e_{\gamma_2},e_{\gamma_3}]]=N_{\gamma_2,\gamma_3}
[e_{\gamma_1},e_{-\gamma_1}]=N_{\gamma_2,\gamma_3}h_{\gamma_1},\]
where we used Remark~\ref{propnrs0}(a). Since the assumption is symmetric
in $\gamma_1$, $\gamma_2$, $\gamma_3$, we also obtain that 
\[[e_{\gamma_2},[e_{\gamma_3},e_{\gamma_1}]]=N_{\gamma_3,\gamma_1}
h_{\gamma_2} \quad\mbox{and} \quad [e_{\gamma_3},[e_{\gamma_1},
e_{\gamma_2}]]=N_{\gamma_1,\gamma_2}h_{\gamma_3}.\]
So the Jacobi identity $[e_{\gamma_1},[e_{\gamma_2},e_{\gamma_3}]]+
[e_{\gamma_2},[e_{\gamma_3}, e_{\gamma_1}]]+[e_{\gamma_3},[e_{\gamma_1},
e_{\gamma_2}]]=0$ yields the identity $N_{\gamma_2,\gamma_3}h_{\gamma_1}+
N_{\gamma_3,\gamma_1}h_{\gamma_2}+N_{\gamma_1,\gamma_2}h_{\gamma_3}=0$. 
Now apply any $\beta\in \Phi$ to the above relation. Using 
Lemma~\ref{astring2}, we obtain 
\begin{align*}
2\Big\langle \beta, & \;\frac{N_{\gamma_2,\gamma_3}}{\langle \gamma_1,
\gamma_1\rangle} \gamma_1+\frac{N_{\gamma_3,\gamma_1}}{\langle 
\gamma_2,\gamma_2\rangle}\gamma_2+ \frac{N_{\gamma_1,\gamma_2}}{\langle 
\gamma_3,\gamma_3\rangle}\gamma_3\Big\rangle \\ &=\frac{2N_{\gamma_2,
\gamma_3}\langle \beta,\gamma_1\rangle}{\langle \gamma_1,\gamma_1\rangle}
+\frac{2N_{\gamma_3,\gamma_1}\langle \beta,\gamma_2\rangle}{\langle 
\gamma_2,\gamma_2\rangle} + \frac{2N_{\gamma_1,\gamma_2}\langle\beta,
\gamma_3 \rangle}{\langle \gamma_3,\gamma_3\rangle}\\
&=\beta\bigl(N_{\gamma_2,\gamma_3}h_{\gamma_1}+N_{\gamma_3,\gamma_1}
h_{\gamma_2}+N_{\gamma_1,\gamma_2}h_{\gamma_3}\bigr)=0.  
\end{align*} 
Since this holds for all $\beta \in\Phi$ and since $E=\langle \Phi
\rangle_\R$, we deduce that 
\[ \frac{N_{\gamma_2,\gamma_3}}{\langle \gamma_1,\gamma_1\rangle}\gamma_1
+\frac{N_{\gamma_3,\gamma_1}}{\langle \gamma_2,\gamma_2\rangle}\gamma_2
+\frac{N_{\gamma_1,\gamma_2}}{\langle \gamma_3,\gamma_3\rangle}\gamma_3=
\underline{0}.\]
Since $\gamma_3=-\gamma_1-\gamma_2$, we obtain
\[ \Bigl(\frac{N_{\gamma_2,\gamma_3}}{\langle \gamma_1,\gamma_1\rangle}
-\frac{N_{\gamma_1,\gamma_2}}{\langle \gamma_3,\gamma_3\rangle}\Bigr)
\gamma_1+\Bigl(\frac{N_{\gamma_3,\gamma_1}}{\langle \gamma_2,\gamma_2
\rangle}-\frac{N_{\gamma_1,\gamma_2}}{\langle \gamma_3,\gamma_3\rangle}
\Bigr)\gamma_2=\underline{0}.\]
Now $\{\gamma_1,\gamma_2\}$ are linearly independent. For otherwise,
we would have $\gamma_2=\pm \gamma_1$ and so $\gamma_3=-2\gamma_1$ or
$\gamma_3=\underline{0}$, contradiction. Hence, the coefficients of
$\gamma_1$, $\gamma_2$ in the above equation must be zero.
\end{proof}

\begin{lem} \label{propnrs2} Let $\alpha,\beta\in\Phi$ be such 
that $\alpha+\beta\in\Phi$. Then 
\[ N_{\alpha,\beta}N_{-\alpha,-\beta}=-p(q+1)\frac{\langle \alpha+\beta,
\alpha+\beta\rangle}{\langle \beta,\beta\rangle},\]
where $\beta-q\alpha,\ldots,\beta-\alpha,\beta, \beta+\alpha,\ldots,
\beta+p\alpha$ is the $\alpha$-string through $\beta$. In particular, 
this shows that $N_{\alpha,\beta} \neq 0$ (since $p\geq 1$ by assumption).
\end{lem}

\begin{proof} We have $[e_{-\alpha},[e_\alpha,e_\beta]]=N_{\alpha,\beta}
[e_{-\alpha},e_{\alpha+\beta}]=N_{\alpha,\beta}N_{-\alpha,\alpha+\beta}
e_\beta$. Applying Lemma~\ref{propnrs1} with $\gamma_1=-\alpha$,
$\gamma_2=\alpha+\beta$, $\gamma_3=-\beta$, we obtain 
\[\frac{N_{-\alpha,\alpha+\beta}}{\langle \beta,\beta\rangle}=
-\frac{N_{-\alpha, -\beta}}{\langle \alpha+\beta,\alpha+\beta\rangle}.\]
On the other hand, let $\slm_2(\C)\cong \fs_\alpha=\langle e_\alpha,
h_\alpha,f_\alpha\rangle \subseteq \cL$ as in Remark~\ref{propnrs0}(b). 
Then, arguing as in Remark~\ref{astring} (where $e_\alpha,h_\alpha,
f_\alpha$ play the role of $e_i,h_i,f_i$, respectively), we find that 
\[ [e_{-\alpha},[e_\alpha,e_\beta]]=[f_\alpha,[e_\alpha,e_\beta]]=
p(q+1) e_\beta.\]
This yields the desired formula. 
\end{proof}

There is also the following result involving four roots.

\begin{lem} \label{propnrs4} Let $\beta_1,\beta_2,\gamma_1,\gamma_2
\in \Phi$ be such that $\beta_1+\beta_2=\gamma_1+\gamma_2\in \Phi$.
Assume that $\beta_1-\gamma_1\not\in \Phi\cup\{\underline{0}\}$ and
that $\beta_2\neq \gamma_1$. Then $\gamma':=\beta_2-\gamma_1= \gamma_2-
\beta_1\in \Phi$ and
\[ N_{\beta_1,\beta_2}N_{-\gamma_1,-\gamma_2}=N_{\beta_1,\gamma'}
N_{-\gamma_1,-\gamma'}\frac{\langle \gamma_2,\gamma_2\rangle}{
\langle \beta_2,\beta_2\rangle}\frac{\langle \gamma',\gamma'
\rangle}{\langle \beta_1+\beta_2,\beta_1+\beta_2 \rangle}.\]
\end{lem}

\begin{proof} By the Jacobi identity we have
\[ [e_{\beta_2},[e_{\beta_1},e_{-\gamma_1}]]+[e_{\beta_1},[e_{-\gamma_1},
e_{\beta_2}]]+[e_{-\gamma_1},[e_{\beta_2},e_{\beta_1}]]=0.\]
Now $[e_{\beta_1},e_{-\gamma_1}]\in \cL_{\beta_1-\gamma_1}$ and, hence,
$[e_{\beta_1},e_{-\gamma_1}]=0$ since $\beta_1-\gamma_1\not\in\Phi\cup
\{\underline{0}\}$. So the first of the above summands is zero and we obtain:
\begin{equation*}
[e_{-\gamma_1},[e_{\beta_1},e_{\beta_2}]]=-[e_{-\gamma_1},[e_{\beta_2},
e_{\beta_1}]]=[e_{\beta_1},[e_{-\gamma_1},e_{\beta_2}]].\tag{$\dagger$}
\end{equation*}
The left hand side of ($\dagger$) evaluates to
\begin{align*}
[e_{-\gamma_1},[e_{\beta_1},e_{\beta_2}]]&=N_{\beta_1,\beta_2}[e_{-\gamma_1},
e_{\beta_1+\beta_2}]\\
&=N_{\beta_1,\beta_2}[e_{-\gamma_1}, e_{\gamma_1+\gamma_2}]=
N_{\beta_1,\beta_2}N_{-\gamma_1,\gamma_1+\gamma_2}e_{\gamma_2}.
\end{align*}
Now $N_{\beta_1,\beta_2}\neq 0$ and $N_{-\gamma_1,\gamma_1+\gamma_2}\neq 0$
by Lemma~\ref{propnrs2}. Hence, the left hand side of ($\dagger$) is
non-zero. So we must have $[e_{-\gamma_1},e_{\beta_2}]\neq 0$. Since
$\beta_2\neq \gamma_1$, this means that $-\gamma_1+\beta_2 \in\Phi$. Then, 
similarly, we find that
\begin{align*}
[e_{\beta_1},[e_{-\gamma_1},e_{\beta_2}]]&=N_{-\gamma_1,\beta_2}
[e_{\beta_1},e_{-\gamma_1+\beta_2}]\\
&=N_{-\gamma_1,\beta_2} [e_{\beta_1},e_{\gamma_2-\beta_1}]
=N_{-\gamma_1,\beta_2}N_{\beta_1,\gamma_2-\beta_1}e_{\gamma_2}.
\end{align*}
This yields $N_{\beta_1,\beta_2}N_{-\gamma_1,\gamma_1+\gamma_2}=
N_{-\gamma_1,\beta_2}N_{\beta_1,\gamma_2-\beta_1}=
N_{-\gamma_1,\beta_2}N_{\beta_1,\gamma'}$. Finally, we have
\[ \frac{N_{-\gamma_1,\beta_2}}{\langle \gamma',\gamma'\rangle}=
\frac{N_{-\gamma',-\gamma_1}}{\langle \beta_2,\beta_2
\rangle}=-\frac{N_{-\gamma_1,-\gamma'}}{\langle \beta_2,\beta_2
\rangle},\]
using Lemma~\ref{propnrs1} with $(-\gamma_1)+\beta_2
+(-\gamma')=\underline{0}$. Furthermore,
\[ \frac{N_{-\gamma_1,\gamma_1+\gamma_2}}{\langle \gamma_2,\gamma_2
\rangle}=\frac{N_{-\gamma_2,-\gamma_1}}{\langle \gamma_1+\gamma_2,
\gamma_1+\gamma_2\rangle}=-\frac{N_{-\gamma_1,-\gamma_2}}{\langle 
\gamma_1+\gamma_2,\gamma_1+\gamma_2\rangle},\]
using Lemma~\ref{propnrs1} with $(-\gamma_1)+(\gamma_1+\gamma_2)+
(-\gamma_2)=\underline{0}$.
\end{proof}

As observed by Chevalley \cite[p.~23]{Ch}, the right hand side of the 
formula in Lemma~\ref{propnrs2} can be simplified, as follows. Let 
$\alpha,\beta\in\Phi$ be such that $\beta\neq \pm \alpha$.  Define
$p,q\geq 0$ as in Lemma~\ref{astring2}. Then 
\[ 2\frac{\langle \alpha,\beta\rangle}{\langle \alpha,\alpha\rangle}
=\beta(h_\alpha)=q-p\in\Z.\]
To simplify the notation, let us denote $\lambda^\vee:=2\lambda/
\langle \lambda,\lambda\rangle\in E$ for any $0 \neq \lambda \in E$. 
Thus, $\langle \alpha^\vee,\beta\rangle=q-p$.
Now, by the \nm{Cauchy--Schwarz inequality}, we have $0\leq \langle
\alpha,\beta\rangle^2<\langle\alpha,\alpha\rangle\cdot \langle \beta,
\beta\rangle$. This yields that 
\[ 0\leq \langle \alpha^\vee,\beta\rangle\cdot \langle\alpha,\beta^\vee
\rangle=2\frac{\langle\alpha,\beta\rangle}{\langle \alpha,\alpha\rangle}
\cdot 2\frac{\langle\alpha,\beta\rangle}{\langle \beta,\beta\rangle}<4.\]
Since $\langle \alpha^\vee,\beta\rangle$ and $\langle \alpha,\beta^\vee
\rangle$ are integers, we conclude that 
\begin{align*} \label{spade12}
\langle &\alpha^\vee,\beta \rangle=q-p\in\{0,\pm 1,\pm 2, \pm 3\}, 
\tag{$\spadesuit_1$} \\
\langle&\alpha^\vee,\beta\rangle=\pm 2 \mbox{ or} \pm 3 
\quad \Rightarrow \quad \langle\alpha,\beta^\vee\rangle=\pm 1.
\tag{$\spadesuit_2$} 
\end{align*}
Now let $\gamma:=\beta-q\alpha\in\Phi$; note that also $\gamma\neq \pm 
\alpha$. Then one immediately sees that the $\alpha$-string through 
$\gamma$ is given by
\[ \gamma,\quad \gamma+\alpha,\quad, \ldots,\quad \gamma+(p+q)\alpha.\]
Applying ($\spadesuit_1$) to $\alpha,\gamma$ yields $\langle 
\alpha^\vee,\gamma\rangle=-(p+q) \in\{0,\pm 1,\pm 2, \pm 3\}$. So 
\begin{align*}
 p&+q=-\langle\alpha^\vee,\gamma\rangle\in\{0,1,2,3\}.
\qquad\qquad\tag{$\spadesuit_3$}
\end{align*}
Now assume that $\alpha+\beta \in \Phi$, as in Lemma~\ref{propnrs2}. 
Then we claim that
\begin{equation*}
r=r(\alpha,\beta):=\frac{\langle\alpha+\beta,\alpha+\beta
\rangle}{\langle \beta,\beta\rangle}=\frac{q+1}{p}.
\;\qquad\qquad\tag{$\spadesuit_4$} 
\end{equation*}
This can now be proved as follows. By ($\spadesuit_3$), we have $0 \leq 
p+q\leq 3$. Since $\alpha+\beta \in \Phi$, we have $p\geq 1$. This leads
to the following cases.

\smallskip
\noindent \fbox{$p=1$, $q=0$ or $p=2$, $q=1$.} Then $\langle
\alpha^\vee,\beta\rangle=q-p=-1$, which means that $2\langle \alpha,
\beta\rangle=-\langle \alpha,\alpha\rangle$. So $\langle \alpha+ \beta,
\alpha+\beta \rangle=\langle\alpha,\alpha \rangle+2\langle\alpha,\beta
\rangle+\langle \beta,\beta\rangle=\langle \beta,\beta\rangle$. Hence, 
$r=1$; we also have $(q+1)/p=1$, as required.

\smallskip
\noindent \fbox{$p=1$, $q=1$.} Then $\langle\alpha^\vee,\beta\rangle=
q-p=0$ and so $\langle \alpha^\vee,\gamma\rangle=-2$, where $\gamma:=
\beta-\alpha$. By ($\spadesuit_2$), we must have $\langle\alpha,
\gamma^\vee\rangle=-1$ and so $2\langle \alpha,\gamma\rangle=-\langle 
\gamma,\gamma\rangle$. Since $\gamma=\beta-\alpha$, this yields 
$\langle\alpha,\alpha\rangle=\langle \beta, \beta\rangle$. Now 
$\langle\alpha^\vee,\beta\rangle =0$ and so $\langle \alpha, \beta 
\rangle=0$. Hence, we obtain $\langle \alpha+\beta,\alpha+\beta\rangle =
\langle \alpha,\alpha\rangle+ \langle \beta,\beta\rangle=2\langle\beta, 
\beta\rangle$. Thus, we have $r=2$ which equals $(q+1)/p=2$ as required. 

\smallskip
\noindent \fbox{$p=1$, $q=2$.} Then $\langle\alpha^\vee,\beta\rangle=
q-p=1$ and so $\langle\alpha^\vee,\gamma\rangle=-3$, where $\gamma:=
\beta-2\alpha$. By ($\spadesuit_2$), we must have $\langle\alpha,
\gamma^\vee\rangle=-1$ and so $2\langle \alpha,\gamma\rangle=-\langle 
\gamma,\gamma\rangle$. Since $\gamma= \beta-2\alpha$, this yields that 
$2\langle\alpha,\beta \rangle=\langle \beta, \beta\rangle$. Now 
$\langle \alpha^\vee,\beta \rangle=1$ also implies that $2\langle 
\alpha,\beta \rangle=\langle \alpha,\alpha\rangle$ and so $\langle
\alpha,\alpha\rangle=\langle \beta,\beta\rangle$. Hence, we obtain 
$\langle \alpha+\beta,\alpha+\beta\rangle=\langle \alpha,\alpha\rangle+
2\langle \alpha,\beta\rangle+\langle \beta,\beta\rangle=3\langle \beta,
\beta\rangle$ and so $r=3$, which equals $(q+1)/p=3$, as required.  

\smallskip
\noindent \fbox{$p\geq 2$, $q=0$.} Then $\langle\alpha^\vee,\beta\rangle
=-p\leq -2$ and so $\langle\alpha,\beta^\vee\rangle=-1$, by 
($\spadesuit_2$). This yields $-p\langle \alpha,\alpha\rangle=2\langle 
\alpha,\beta\rangle=-\langle \beta,\beta\rangle$ and so $\langle\alpha+
\beta,\alpha+ \beta \rangle=\langle \alpha,\alpha\rangle +2\langle\alpha, 
\beta \rangle+\langle \beta,\beta\rangle=\frac{1}{p}\langle \beta,\beta 
\rangle$. Hence, $r=\frac{1}{p}=\frac{q+1}{p}$, as required. 

Thus, the identity in ($\spadesuit_4$) holds in all cases and we obtain:

\begin{prop}[Chevalley] \label{propnrs3} Let $\alpha,\beta\in\Phi$ be such 
that $\alpha+\beta\in\Phi$. Using the notation in Lemma~\ref{propnrs2}, we
have  
\[ N_{\alpha,\beta}N_{-\alpha,-\beta}=-(q+1)^2.\]
\end{prop}

\begin{proof} Since $\alpha+\beta\in \Phi$, we have $\beta \neq \pm 
\alpha$. We have seen above that then ($\spadesuit_4$) holds. It 
remains to use the formula in Lemma~\ref{propnrs2}.
\end{proof}

The above formula suggests that there might be a clever choice of the
elements $\{e_\alpha\mid \alpha\in \Phi\}$ such that $N_{\alpha,\beta}=
\pm (q+1)$ whenever $\alpha+\beta\in \Phi$. We will pursue this issue 
further in the following section.

\begin{exmp} \label{propnrs1a} Suppose we know all $N_{\alpha_j,\beta}$, 
where $j\in I$ and $\beta\in\Phi^+$. We claim that then all structure 
constants $N_{\pm\alpha_i,\alpha}$ for $i\in I$ and $\alpha\in\Phi$ can be 
determined, using only manipulations with roots in $\Phi$. 

(1) First, let $i\in I$ and $\alpha\in\Phi^-$. Then 
Proposition~\ref{propnrs3} shows how to express $N_{-\alpha_i,\alpha}$
in terms of $N_{\alpha_i,-\alpha}$ (which is known by assumption).

(2) Next, we determine $N_{-\alpha_i,\alpha}$ for $i\in I$ and 
$\alpha\in\Phi^+$. If $\alpha-\alpha_i\not\in\Phi$, then
$N_{-\alpha_i,\alpha}=0$. Now assume that $\alpha-\alpha_i\in\Phi$.
Then $(-\alpha_i)+\alpha-(\alpha-\alpha_i)=\underline{0}$ and so 
Lemma~\ref{propnrs1} yields that
\[ \frac{N_{-\alpha_i,\alpha}}{\langle \alpha-\alpha_i,\alpha-\alpha_i
\rangle}=\frac{N_{-(\alpha-\alpha_i),-\alpha_i}}{\langle \alpha,
\alpha \rangle}=-\frac{N_{-\alpha_i,-(\alpha-\alpha_i)}}{\langle 
\alpha,\alpha \rangle}.\] 
Since $-(\alpha-\alpha_i)\in\Phi^-$, the right hand side
can be handled by~(1).

(3) Finally, if $i\in I$ and $\alpha\in\Phi^-$, then 
Proposition~\ref{propnrs3} expresses $N_{\alpha_i,\alpha}$ in terms of 
$N_{-\alpha_i,-\alpha}$, which is handled by (2) since
$-\alpha\in\Phi^+$.  

Of course, if we want to do this in a concrete example, then we need to 
be able to perform computations with roots in $\Phi$: check if the sum 
of roots is again a root, or calculate the scalar product of a
root with itself. More precisely, we do not need to know the actual values 
of those scalar products, but rather the values of fractions 
$r(\alpha,\beta)=\langle\alpha+\beta,\alpha+\beta \rangle/\langle \beta,
\beta \rangle$ as above; we have seen in ($\spadesuit_4$) how such 
fractions are determined.
\end{exmp}

To illustrate the above results, let us consider the matrix 
\[ \renewcommand{\arraystretch}{0.9} A=\left(\begin{array}{rr} 2 & -1 
\\ -3 & 2 \end{array}\right).\]
In Example~\ref{cartanG2}, we have computed corresponding ``roots'',
although we do not know (yet) if there is a Lie algebra with $A$  as 
structure matrix. We can now push this discussion a bit further. 



\begin{table}[htbp] \caption{Structure constants for type $G_2$} 
\label{nrsG2}
\begin{center}
{\small $\renewcommand{\arraystretch}{1.0}
\renewcommand{\arraycolsep}{3.5pt}
\begin{array}{r|cccccc|cccccc} \hline N_{\alpha,\beta} & \;\;10 & \;\;01 & 
\;\;11 &\;\;12 & \;\;13 & \;\;23 &-10&-01&-11&-12&-13&-23\\ \hline
 10 &  . & 1 & . & . & 1 & . & * & . & 1 & . & . &-1\\
 01 & -1 & . &-2 &-3 & . & . & . & * &-3 & 2 &-1 & .\\
 11 &  . & 2 & . &-3 & . & . &-1 & 3 & * & 2 & . &-1\\
 12 &  . & 3 & 3 & . & . & . & . & 2 &-2 & * & 1 &-1\\
 13 & -1 & . & . & . & . & . & . & 1 & . &-1 & * &-1\\
 23 &  . & . & . & . & . & . &-1 & . & 1 &-1 & 1 & *\\\hline
-10 &  * & . & 1 & . & . & 1 & . &-1 & . & . &-1 & .\\
-01 &  . & * &-3 &-2 &-1 & . & 1 & . & 2 & 3 & . & .\\
-11 & -1 & 3 & * & 2 & . &-1 & . &-2 & . & 3 & . & .\\
-12 &  . &-2 &-2 & * & 1 & 1 & . &-3 &-3 & . & . & .\\
-13 &  . & 1 & . &-1 & * &-1 & 1 & . & . & . & . & .\\
-23 &  1 & . & 1 & 1 & 1 & * & . & . & . & . & . & .\\\hline
\multicolumn{13}{c}{\text{(Here, e.g., $-12$ stands for $-(\alpha_1+
2\alpha_2)\in \Phi$, ``$*$'' for $h_\alpha$ and ``$.$'' for $0$.)}}
\end{array}$}
\end{center}
\end{table}

\begin{exmp} \label{hypog2} Assume that there exists a Lie algebra $\cL$ 
with subalgebra $\fh\subseteq \cL$ such that $(\cL,\fh)$ is of Cartan--Killing
type with respect to $\Delta=\{\alpha_1,\alpha_2\}$ and corresponding
structure matrix as above:
\[\renewcommand{\arraystretch}{0.9} A=\left(\begin{array}{rr} 2 & -1 
\\ -3 & 2 \end{array}\right) \qquad \mbox{(called of ``type $G_2$'')}.\]
Then, as in Example~\ref{cartanG2}, $W$ is dihedral of order $12$ and 
\[ \Phi^+=\{\alpha_1,\;\alpha_2,\;\alpha_1+\alpha_2,\;\alpha_1+2\alpha_2,\;
\alpha_1+3\alpha_2,\;2\alpha_1+3\alpha_2\}.\]
We have $-\langle\alpha_1,\alpha_1\rangle=2\langle \alpha_1,
\alpha_2\rangle=-3\langle \alpha_2,\alpha_2\rangle$ and so 
$\langle \alpha_1,\alpha_1\rangle=3\langle\alpha_2,\alpha_2\rangle$.
From the computation in Example~\ref{cartanG2}, we also see that
\begin{align*}
\Phi_1:=\{w(\alpha_1)\mid w\in W\}&=\{\alpha_1,\alpha_1+3\alpha_2,2\alpha_1+
3\alpha_2\},\\
\Phi_2:=\{w(\alpha_2)\mid w\in W\}&=\{\alpha_2,\alpha_1+\alpha_2,\alpha_1+
2\alpha_2\}.
\end{align*}
Thus, $\langle\alpha,\alpha\rangle/\langle \beta,\beta\rangle$ is known
for all $\alpha,\beta\in\Phi$. Let $\{e_1,e_2,f_1,f_2\}$ be Chevalley 
generators for $\cL$. Let us try to determine a collection of elements 
$\{\be_\alpha \mid \alpha \in \Phi\}$ and the corresponding 
structure constants. Anticipating what we will do in the following 
section, let us set 
\[ \be_{\alpha_1}=e_1, \qquad \be_{\alpha_2}=-e_2, \qquad 
\be_{-\alpha_1}=f_1, \qquad \be_{-\alpha_2}=-f_2.\]
For $i\in I$ and $\alpha\in \Phi$, let $q_{i,\alpha}:=\max\{m\geq 0
\mid \alpha-m\alpha_i\in \Phi\}$. In view of the formula in 
Proposition~\ref{propnrs3}, we define successively:
\begin{alignat*}{2}
\be_{\alpha_1+\alpha_2}&:=\;\;\,[e_1,\be_{\alpha_2}]
& \in \cL_{\alpha_1+\alpha_2}\quad & \qquad (q_{1,\alpha_2}=0),\\ 
\be_{\alpha_1+2\alpha_2}&:=\textstyle{\frac{1}{2}}[e_2,\be_{\alpha_1
+\alpha_2}] & \in \cL_{\alpha_1+2\alpha_2}\;\;& \qquad (q_{2,\alpha_1+
\alpha_2}=1),\\ \be_{\alpha_1+3\alpha_2}&:=\textstyle{\frac{1}{3}}
[e_2,\be_{\alpha_1+2\alpha_2}]&\in \cL_{\alpha_1+3\alpha_2}\;\;& 
\qquad (q_{2,\alpha_1+2\alpha_2}=2),\\\be_{2\alpha_1+3\alpha_2}&:=
\;\;\,[e_1,\be_{\alpha_1+3\alpha_2}]\;\;& \in 
\cL_{2\alpha_1+3\alpha_2}\,& \qquad (q_{1,\alpha_1+3\alpha_2}=0).
\end{alignat*}
All these are non-zero by Lemma~\ref{propnrs2}. Hence, for $\alpha \in 
\Phi^+$, there is a unique $\be_{-\alpha}\in \cL_{-\alpha}$ such that 
$[\be_\alpha,\be_{-\alpha}]=h_\alpha$. Thus, we have defined elements 
$\be_\alpha\in \cL_\alpha$ for all $\alpha \in \Phi$, such that
Remark~\ref{propnrs0}(a) holds. Let $N_{\alpha,\beta}$ be the 
corresponding structure constants; we leave it as an exercise for the 
reader to check that these are given by Table~\ref{nrsG2}. (In order 
to compute that table, one only needs arguments like those in 
Example~\ref{propnrs1a}.) Thus, without knowing that $\cL$ exists at 
all, we are able to compute all the structure constants $N_{\alpha,
\beta}$~---~and we see that they are all integers!  Furthermore, using 
Lemma~\ref{astring3}, we obtain
\begin{gather*}
h_{\alpha_1+\alpha_2}=3h_1+h_2,\qquad
h_{\alpha_1+2\alpha_2}=3h_1+2h_2,\\
h_{\alpha_1+3\alpha_2}=h_1+h_2,\qquad
h_{2\alpha_1+3\alpha_2}=2h_1+h_2.
\end{gather*}
Thus, all the Lie brackets in $\cL$ are known and the whole
situation is completely rigid. One could try to construct a Lie algebra 
with these properties using a suitable factor algebra of the free Lie
 algebra over $X=\{e_1,e_2,f_1,f_2\}$ (as in Example~\ref{expsolv2})
but, still, one has to show that such a factor algebra has the correct
dimension (namely, $14$). 
\end{exmp}

Here is a further illustration of the power of the relation 
($\spadesuit_1$).

\begin{exmp} \label{xcastdbase22} Let $\alpha,\beta\in \Phi$ be such
that $\beta\neq \pm\alpha$. We claim that, if $r,s\geq 1$ are integers
such that $r\alpha+ s\beta\in\Phi$, then $\alpha+\beta\in \Phi$. \\ (Roots
of the form $r\alpha+s\beta$ where $r,s\geq 1$ will occur in Chevalley's
commutator relations, to be discussed in a later section.)\\
This is seen as follows. Let $p,q\geq 0$ be as in Lemma~\ref{astring2};
then $\langle \alpha^\vee,\beta\rangle=q-p$. If $\langle \alpha^\vee,
\beta\rangle<0$, then $p>0$ and so $\alpha+\beta\in \Phi$, as desired. Now 
let $\langle \alpha^\vee,\beta\rangle\geq 0$. Then $\langle 
\alpha^\vee,r\alpha+s\beta\rangle=2r+s\langle \alpha^\vee, \beta\rangle 
\geq 2r$. Since, by ($\spadesuit_1$), the left hand side has absolute 
value $\leq 3$, we conclude that $r=1$. We also have $\langle \beta^\vee, 
\alpha\rangle\geq 0$ and so $\langle \beta^\vee,r\alpha+s\beta\rangle=
r\langle \beta^\vee, \alpha\rangle +2s\geq 2s$. Again, we conclude that 
$s=1$. Thus, $\alpha+\beta\in \Phi$. 
\end{exmp}


\section{Lusztig's canonical basis} \label{sec1a7}

We keep the general setting of the previous section. Recall that 
$\dim \fg_\alpha=1$ for each root $\alpha\in \Phi$. The aim of this 
section is to show the remarkable fact that one can single out a 
``canonical'' collection of elements in the various weight 
spaces~$\cL_\alpha$.
 
\begin{rem} \label{canbas0a} Let $i\in I$ and $\beta\in\Phi$ be such
that $\beta\neq \pm\alpha_i$. As in Remark~\ref{astring}, let 
$\beta-q\alpha_i,\ldots,\beta-\alpha_i,\beta,\beta+\alpha_i,\ldots,
\beta+p \alpha_i$ be the \nms{$\alpha_i$-string through 
$\beta$}{alphai-string through beta}. By Exercise~\ref{xcastring}, we 
have 
\begin{align*}
p=p_{i,\beta}&:=\max\{m\geq 0\mid \beta+m\alpha_i\in \Phi\},\\
q=q_{i,\beta}&:=\max\{m\geq 0\mid \beta-m\alpha_i\in \Phi\}.
\end{align*}
Also note that, for any $m\geq 0$, we have $\beta-m\alpha_i\in\Phi$ if 
and only if $-\beta+m\alpha_i=-(\beta-m\alpha_i)\in\Phi$. 
Thus, we have $q_{i,\beta}=p_{i,-\beta}$.
\end{rem}

\begin{thm}[Lusztig \protect{\cite[\S 1]{L1}, \cite[\S 2]{L2},
\cite[Theorem~0.6]{L5}}\footnote{The result, as stated here, is just the
shadow of a much more sophisticated and powerful result about quantized
enveloping algebras.}] \label{canbas} Given Chevalley generators $\{e_i,
f_i\mid i \in I\}$ of $\cL$, there is a collection of elements $\{0\neq 
\be_\alpha^+ \in \cL_\alpha \mid \alpha \in \Phi\}$ with the following 
properties:
\begin{itemize}
\item[{\rm (L1)}] $[f_i,\be_{\alpha_i}^+]=[e_i,\be_{-\alpha_i}^+]$ for 
all $i\in I$.
\item[{\rm (L2)}] $[e_i,\be_\alpha^+]=(q_{i,\alpha}+1) \be_{\alpha+
\alpha_i}^+$ if $i\in I$, $\alpha\in \Phi$ and $\alpha+\alpha_i \in \Phi$.
\item[{\rm (L3)}] $[f_i,\be_\alpha^+]=(p_{i,\alpha}+1) \be_{\alpha-
\alpha_i}^+$ if $i\in I$, $\alpha\in \Phi$ and $\alpha-\alpha_i \in \Phi$.
\end{itemize}
If $A$ is indecomposable, then this collection $\{\be_\alpha^+\mid \alpha
\in \Phi\}$ is unique up to a global constant, that is, if $\{0\neq
\be_\alpha'\in \cL_\alpha\mid\alpha\in \Phi\}$ is another collection 
satisfying {\rm (L1)--(L3)}, then there exists some $0\neq \xi\in\C$ such 
that $\be_{\alpha}'=\xi\be_{\alpha}^+$ for all $\alpha\in \Phi$.
\end{thm}

The proof will be given later in this section (from 
Definition~\ref{startproof} on), after the following remarks. First note 
that, even for $\cL=\slm_2(\C)$, we have to modify the standard 
elements $e,h,f$ in order to obtain the above formulae. Indeed, setting 
$\be^+:=e$ and $\bbf^+:=-f$, we have
\[[e,\bbf^+]=-[e,f]=-h=[f,e]=[f,\be^+].\]
Hence, $\{\be^+,\bbf^+\}$ is a collection satisfying (L1); the conditions
in (L2) and (L3) are empty in this case. (See also Exercise~\ref{lucanAn}
below.)

\begin{rem} \label{canbas0} Assume that a collection $\{\be_\alpha^+
\mid \alpha \in \Phi\}$ as in Theorem~\ref{canbas} exists. Since 
$\be_{\alpha_i}^+\in \cL_{\alpha_i}$ for $i\in I$, we have $\be_{\alpha_i}^+
=c_ie_i$, where $0\neq c_i\in\C$. Similarly, we have $\be_{-\alpha_i}^+
\in \cL_{-\alpha_i}$ and so $\be_{-\alpha_i}^+=d_if_i$, where $0\neq d_i
\in\C$. Hence, we obtain 
\begin{align*}
[f_i,\be_{\alpha_i}^+]&=c_i[f_i,e_i]=-c_i[e_i,f_i]=-c_ih_i,\\ 
[e_i,\be_{-\alpha_i}^+]&=d_i[e_i,f_i]=d_ih_i,
\end{align*}
and so (L1) implies that $d_i=-c_i$ for all $i\in I$. This also shows
that $[\be_{\alpha_i}^+,\be_{-\alpha_i}^+]=c_id_i[e_i,f_i]=-c_i^2h_i$ for 
$i\in I$. Thus, Remark~\ref{propnrs0}(a) does not seem to hold for the
collection $\{\be_{\alpha}^+\mid \alpha\in \Phi\}$. (This issue will 
be resolved later in Corollary~\ref{canbash} below.)

Now, the possibilities for the constants $c_i$ are severely restricted, 
as follows. Let $i,j\in I$ be such that $i\neq j$ and $a_{ij}\neq 0$. 
Then $\beta=\alpha_i+\alpha_j\in\Phi$; see Exercise~\ref{xcastringm}. 
Applying (L2) twice, we obtain:
\begin{align*}
[e_i,e_j]&=[e_i,c_j^{-1}\be_{\alpha_j}^+]=(q_{i,\alpha_j}+1)c_j^{-1}
\be_\beta^+= c_j^{-1}\be_\beta^+,\\
[e_j,e_i]&=[e_j,c_i^{-1}\be_{\alpha_i}^+]= (q_{j,\alpha_i}+1)c_i^{-1}
\be_\beta^+=c_i^{-1}\be_\beta^+.
\end{align*}
Note that $\pm(\alpha_i-\alpha_j)\not\in \Phi$ and so $q_{j,\alpha_i}=
q_{i,\alpha_j}=0$. Since $[e_i,e_j]=-[e_j,e_i]$, we conclude that 
$c_j=-c_i$. Thus
\begin{equation*}
c_j=-c_i \quad \mbox{whenever $i,j\in I$ are such that $a_{ij}<0$}.\tag{$*$}
\end{equation*}
Thus, the function $i \mapsto c_i$ has the property in 
Exercise~\ref{xcaforest}(c).

If $A$ is indecomposable, then ($*$) implies that $\{c_i\mid i\in I\}$ is 
completely determined by $c_{i_0}$, for one particular choice of $i_0\in I$.
Indeed, let $i\in I$, $i\neq i_0$. By Remark~\ref{indecgraph}, there is 
a sequence of (distinct) indices $i_0,i_1,\ldots,i_r=i$ ($r\geq 1$) 
such that $a_{i_li_{l+1}}\neq 0$ for $0\leq l\leq r-1$. Hence, by ($*$),
we find that $c_i=(-1)^rc_{i_0}$. Consequently, if $\{c_i'\mid i\in I\}$ is 
another collection of non-zero constants satisfying~($*$), then $c_i'=
\xi c_i$ for all $i\in I$, where $\xi=c_{i_0}' c_{i_0}^{-1}\in\C^\times$ 
is a constant.
\end{rem}

\begin{rem} \label{canbas1} Assume that a collection $\{\be_\alpha^+
\mid \alpha \in \Phi\}$ as in Theorem~\ref{canbas} exists. 
Using (L1), we can define
\[ h_j^+:=[e_j,\be_{-\alpha_j}^+]=[f_j,\be_{\alpha_j}^+] \in \fh
\qquad \mbox{for all $j\in I$}.\]
Writing $\be_{\alpha_j}^+=c_je_j$ as in Remark~\ref{canbas0}, we see
that $h_j^+=-c_jh_j$. So
\[ \bB:=\{h_j^+\mid j\in I\}\cup\{\be_\alpha^+\mid \alpha \in\Phi\}
\quad \mbox{is a basis of $\cL$}.\]
We claim that the action of the Chevalley generators $\{e_i,f_i\mid 
i \in I\}$ on this basis is given as follows, where $j\in I$ and 
$\alpha\in\Phi$:
\begin{empheq}[box=\widefbox]{align*}
[e_i,h_j^+] &= |a_{ji}|\be_{\alpha_i}^+, \qquad [f_i,h_j^+] = |a_{ji}| 
\be_{-\alpha_i}^+, \\ [e_i,\be_\alpha^+] &= \left\{
\begin{array}{c@{\hspace{5pt}}l} (q_{i,\alpha}+1) \be_{\alpha+\alpha_i}^+ 
& \quad \mbox{if $\alpha+\alpha_i \in\Phi$},\\ h_i^+ & \quad \mbox{if 
$\alpha=-\alpha_i$},\\ 0  & \quad \mbox{otherwise},\end{array}\right.\\
[f_i,\be_\alpha^+] &=\left\{\begin{array}{c@{\hspace{5pt}}l} 
(p_{i,\alpha}+1) \be_{\alpha -\alpha_i}^+ & \quad \mbox{if $\alpha-
\alpha_i \in\Phi$},\\ h_i^+ & \quad \mbox{if $\alpha=\alpha_i$},\\ 0 & 
\quad \mbox{otherwise}.  \end{array} \right.
\end{empheq}
Indeed, let $\alpha\in \Phi$. If $\alpha=-\alpha_i$, then
$\be_\alpha^+=-c_if_i$ and so $[e_i,\be_\alpha^+]=-c_i[e_i,f_i]=-c_ih_i=
h_j^+$. Now let $\alpha\neq -\alpha_i$; if $\alpha+ \alpha_i\not
\in \Phi$, then $[e_i,\be_\alpha^+]=0$; otherwise, $[e_i,\be_{\alpha_i}^+]$ 
is given by (L2). Similarly, if $\alpha=\alpha_i$, then $[f_i,\be_\alpha^+]
=h_i^+$. Now let $\alpha\neq \alpha_i$. If $\alpha-\alpha_i\not
\in\Phi$, then $[f_i,\be_\alpha^+]=0$; otherwise, $[f_i,\be_\alpha^+]$ is 
given by (L3). Now let $j\in I$. Then
\[[e_i,h_j^+]=-[h_j^+,e_i]=c_j[h_j,e_i]=c_j\alpha_i(h_j)e_i=
c_ja_{ji}e_i.\]
If $i=j$, then $a_{ji}=2$ and $c_je_i=c_ie_i=\be_{\alpha_i}^+$; 
thus, $[e_i,h_i^+]=2\be_{\alpha_i}^+$. Now let $i\neq j$. If $a_{ji}=0$, 
then $[e_i,h_j^+]=0$. If $a_{ji}\neq 0$, then $c_i=-c_j$ by
Remark~\ref{canbas0}. So $[e_i,h_j^+]=-c_ia_{ji}e_i=-a_{ji}
\be_{\alpha_i}^+$, where $a_{ji}<0$. This yields the above formula for 
$[e_i,h_j^+]$. Finally, consider $f_i$. We have seen in Remark~\ref{canbas0}
that $\be_{-\alpha_i}^+=-c_if_i$. This yields that 
\[[f_i,h_j^+]=-[h_j^+,f_i]=c_j[h_j,f_i]=-c_j\alpha_i(h_j)f_i=
-c_ja_{ji}f_i.\]
Now we argue as before to obtain the formula for $[f_i,h_j^+]$. 

Thus, all the entries of the matrices of $\ad_\cL(e_i)$ and $\ad_\cL(f_i)$ 
with respect to the basis $\bB$ are non-negative integers! This is one of 
the remarkable features of Lusztig's theory of ``canonical bases'' (see 
\cite{L6}, \cite{L5} and further references there).
\end{rem}

\begin{rem} \label{canbas2} Assume that $A$ is indecomposable and that a 
collection $\{\be_\alpha^+ \mid \alpha \in \Phi\}$ as in
Theorem~\ref{canbas} exists. First note that, if $0\neq \xi\in \C$ is 
fixed and we set $\be_{\alpha}':= \xi\be_{\alpha}^+$ for all $\alpha\in
\Phi$, then the new collection $\{\be_\alpha'\mid \alpha\in\Phi\}$ also 
satisfies (L1)--(L3). Conversely, we show that any two collections 
satisfying (L1)--(L3) are related by such a global constant~$\xi$. 

Now, as above, for $i\in I$ we have $\be_{\alpha_i}^+=c_ie_i$, where 
$0\neq c_i\in\C$. Then (L2) combined with the Key Lemma~\ref{keylem} 
determines $\be_\alpha^+$ for all $\alpha \in\Phi^+$. Furthermore, as
above, we have $\be_{-\alpha_i}^+=-c_if_i$ for $i\in I$. But then (L3) also 
determines $\be_{-\alpha}^+$ for all $\alpha \in\Phi^+$. Thus, the whole 
collection $\{\be_\alpha^+\mid \alpha\in \Phi\}$ is completely determined 
by $\{c_i\mid i\in I\}$ and properties of $\Phi$ (e.g., the numbers 
$p_{i,\alpha}$, $q_{i,\alpha}$). 

Now assume that $\{\be_\alpha'\mid \alpha\in\Phi\}$ is any other collection
that satisfies (L1)--(L3). For $i\in I$, we have again $\be_{\alpha_i}'=
c_i'e_i$, where $0\neq c_i'\in \C$. Now both collections of constants 
$\{c_i\mid i \in I\}$ and $\{c_i'\mid i \in I\}$ satisfy ($*$) in
Remark~\ref{canbas0}. So there is some $0\neq \xi \in \C$ such that
$c_i'=\xi c_i$ for all $i\in I$. Hence, we have $\be_{\alpha_i}'=
\xi\be_{\alpha_i}^+$ for all $i\in I$. But then the previous discussion
shows that $\be_\alpha'=\xi \be_\alpha^+$ for all $\alpha\in \Phi$.
This proves the uniqueness part of Theorem~\ref{canbas}.
\end{rem}

We now turn to the existence part of Theorem~\ref{canbas}. The crucial 
step will be the construction of $\be_\alpha^+\in \fg_\alpha$ for 
$\alpha\in \Phi^+$. Lusztig's argument in \cite[Lemma~1.4]{L1} assumes 
that $A=(a_{ij})_{i,j \in I}$ is indecomposable and proceeds by a downward
induction on $\hgt(\alpha)$, starting with a root of maximal height; it 
is also assumed that $a_{ij} \in \{0,-1\}$ for all $i\neq j$ in~$I$. (In 
\cite[\S 2]{L2}, the latter assumption is removed, but there are no 
details about the proof; in \cite{L5}, the proof is based on general
results on canonical bases in \cite{L6}.) We shall proceed here by an 
{\itshape upward} induction on $\hgt(\alpha)$ for $\alpha\in \Phi^+$, 
one side effect of which is that we do not need to assume that $A$ is 
indecomposable.

\begin{defn}  \label{startproof} We fix a collection of signs $c_i\in
\{\pm 1\}$ ($i\in I$) satisfying condition~($*$) in Remark~\ref{canbas0}. 
(This is easily seen to exist since the graph of $A$ is a 
\nm{forest}\footnote{In Section~\ref{sec3a1} below, we shall show that
each connected component of the graph of $A$ is described by one of the 
diagrams in Table~\ref{Mdynkintbl} (p.~\pageref{Mdynkintbl}); a specific 
choice of the signs $c_i$ for each of those diagrams will be displayed 
in Table~\ref{Mdynkineps} (p.~\pageref{Mdynkineps}).}; see 
Proposition~\ref{forestg} and Exercise~\ref{xcaforest}.) Let $|I|=n$ 
and fix a numbering $I=\{1,\ldots,n\}$. 

For $\alpha\in \Phi^+$ we now define a particular element $\be_\alpha
\in \cL_\alpha$, where we proceed by induction on $\hgt(\alpha)$. If
$\hgt(\alpha)=1$, then $\alpha=\alpha_i$ for $i\in I$; in this case, we 
set $\be_{\alpha_i}:=c_ie_i$. Now assume that $\hgt(\alpha)>1$ and that 
$\be_\beta\in \cL_\beta$ has been already defined for all $\beta \in 
\Phi^+$ with $\hgt(\beta)<\hgt(\alpha)$. By Lemma~\ref{keylem} there is 
some $i\in I$ such that $\beta:=\alpha-\alpha_i\in \Phi^+$. There may be 
several $i$ with this property; in order to make a specific choice, we 
take the smallest $i\in I=\{1,\ldots,n\}$ such that $\alpha-\alpha_i\in
\Phi^+$. Then $0\neq [e_i,\be_\beta]\in \cL_\alpha$ and we define 
$\be_\alpha \in \cL_\alpha$ be the condition that 
\[ [e_i,\be_\beta]=(q_{i,\beta}+1)\be_\alpha.\]
Once $\be_\alpha$ is defined for each $\alpha \in \Phi^+$, there is a 
unique $\be_{-\alpha}\in \cL_{-\alpha}$ such that $[\be_\alpha,
\be_{-\alpha}]=h_\alpha$. Thus, we obtain a complete collection
\begin{center}
\fbox{$\{\be_\alpha\mid \alpha\in\Phi\} \quad$ such that
Remark~\ref{propnrs0}(a) holds.}
\end{center}
Let $N_{\alpha,\beta}$ be the structure constants with respect to the
above collection; since Remark~\ref{propnrs0}(a) holds (by construction), 
all the results in Section~\ref{sec1a6} can be used. 
\end{defn}

The following result is the crucial step in the proof of 
Theorem~\ref{canbas}. It shows that the  collection of elements
$\{\be_\alpha\mid\alpha\in\Phi\}$ does not depend at all on the choice
of the numbering $I=\{1,\ldots,n\}$.

\begin{lem} \label{canbas5} Let $\alpha\in\Phi^+$ and $j\in I$ 
be arbitrary such that $\gamma:=\alpha-\alpha_j\in\Phi^+$. Then we also
have $[e_j,\be_\gamma]=(q_{j,\gamma}+1)\be_\alpha$.
\end{lem}

\begin{proof} We proceed by induction on $\hgt(\alpha)$. If $\hgt(\alpha)=1$,
then $\alpha=\alpha_i$ for some $i\in I$. In that case, there is no
$j\in I$ such that $\alpha-\alpha_j\in \Phi^+$ and so there is nothing 
to prove. Now assume that $\hgt(\alpha)>1$ and let $i\in I$ be minimal 
such that $\beta:=\alpha-\alpha_i \in \Phi^+$, as in 
Definition~\ref{startproof}; thus, $[e_i, \be_\beta]=(q_{i,\beta}+1)
\be_\alpha$. Let also $j\in I$ be such that $\gamma:=\alpha-\alpha_j\in 
\Phi^+$. We must show that $[e_j,\be_\gamma]=(q_{j,\gamma}+1) \be_\alpha$. 
If $i=j$, then this is trivially true. Now assume that $i\neq j$. Then
we have two expressions 
\[ \alpha_i+\beta=\alpha=\alpha_j+\gamma \qquad \mbox{where} \qquad
\alpha_i- \alpha_j\not\in \Phi\cup\{\underline{0}\}.\]
If $\hgt(\alpha)=2$, then $\alpha=\alpha_i+\alpha_j$ where $\beta=
\alpha_j$ and $\gamma=\alpha_i$; furthermore, $q_{i,\alpha_j}=q_{j,
\alpha_i}=0$. Now, we have $\be_{\alpha_i}=c_ie_i$ and $\be_{\alpha_j}=
c_je_j$. Hence,
\[ \be_\alpha=[e_i,\be_\beta]=c_j[e_i,e_j]=-c_j[e_j,e_i]=-c_ic_j[e_j,
\be_\gamma].\]
Since $\alpha_i+\alpha_j=\alpha\in \Phi$, we have $a_{ij}\neq 0$ (see 
Exercise~\ref{xcastringm}) and so $c_i=-c_j$. Hence, the desired identity 
holds in this case. Now assume that $\hgt(\alpha)>2$. Then we still have 
$\hgt(\beta)>1$ and so $\beta \neq \alpha_j$. Hence, we can apply 
Lemma~\ref{propnrs4} which yields that 
\begin{equation*}
N_{\alpha_i,\beta}N_{-\alpha_j,-\gamma}=N_{\alpha_i,\gamma'}
N_{-\alpha_j,-\gamma'}\frac{\langle \gamma,\gamma\rangle}{\langle 
\beta,\beta\rangle}\frac{\langle \gamma',\gamma'\rangle}{\langle 
\alpha, \alpha\rangle},\tag{$\dagger_1$}
\end{equation*}
where $\gamma':=\beta-\alpha_j=\gamma-\alpha_i \in \Phi$; note that 
$\gamma'\in \Phi^+$. Now, one could try to simplify the right hand side 
using the formulae in the previous section. But there is a simple trick
(taken from \cite[\S 2.9, Lemma~E]{Sam}) to avoid such calculations.
Namely, we can also apply Lemma~\ref{propnrs4} to the two expressions 
$-\alpha_i-\beta=-\alpha= -\alpha_j-\gamma$, where $\alpha_j-\alpha_i
\not\in \Phi\cup \{\underline{0}\}$. This yields the identity:
\begin{equation*}
N_{-\alpha_i,-\beta}N_{\alpha_j,\gamma}=N_{-\alpha_i,-\gamma'}
N_{\alpha_j,\gamma'}\frac{\langle \gamma,\gamma\rangle}{\langle 
\beta,\beta\rangle}\frac{\langle \gamma',\gamma'\rangle}{\langle 
\alpha, \alpha\rangle},\tag{$\dagger_2$}
\end{equation*}
Now, we have $\alpha_i+\gamma'=\gamma$ and $\hgt(\gamma)=\hgt(\alpha)-1$;
similarly, $\alpha_j+\gamma'=\beta$ and $\hgt(\beta)=\hgt(\alpha)-1$. 
So we can apply induction and obtain that 
\[ [e_i,\be_{\gamma'}]=(q_{i,\gamma'}+1)\be_{\gamma}
\qquad \mbox{and}\qquad [e_j,\be_{\gamma'}]=(q_{j,\gamma'}+1)
\be_{\beta}.\]
Since $\be_{\alpha_i}=c_ie_i$ and $\be_{\alpha_j}=c_je_j$, the above 
formulae mean that 
\[N_{\alpha_i,\gamma'}=c_i(q_{i,\gamma'}+1)\qquad\mbox{and}\qquad
N_{\alpha_j,\gamma'}=c_j(q_{j,\gamma'}+1).\]
But then the formula in Proposition~\ref{propnrs3} shows that 
\begin{align*}
N_{-\alpha_i,-\gamma'}&=-c_i^{-1}(q_{i,\gamma'}+1)=-c_i(q_{i,\gamma'}+1)=
-N_{\alpha_i,\gamma'},\\ N_{-\alpha_j,-\gamma'}&=-c_j^{-1}(q_{j,\gamma'}+1)
=-c_j(q_{j,\gamma'}+1)=-N_{\alpha_j,\gamma'}.
\end{align*}
(Recall that $c_i,c_j\in \{\pm 1\}$.) Hence, the right hand side of
($\dagger_1$) is equal to the right hand side of ($\dagger_2$). 
Consequently, the two left hand sides are also equal and we obtain
\[ N_{\alpha_i,\beta}N_{-\alpha_j,-\gamma}=N_{-\alpha_i,-\beta}
N_{\alpha_j,\gamma}.\]
Now, we have $[e_i,\be_\beta]=(q_{i,\beta}+1)\be_\alpha$ and so
$N_{\alpha_i,\beta}=c_i(q_{i,\beta}+1)$. By Proposition~\ref{propnrs3}, 
this also yields that $N_{-\alpha_i-\beta}=-c_i(q_{i,\beta}+1)$. Inserting 
this into the above identity, we deduce that $N_{-\alpha_j,-\gamma}=
-N_{\alpha_j,\gamma}$. Hence, a further application of 
Proposition~\ref{propnrs3} shows that 
\[ N_{\alpha_j,\gamma}=\pm (q_{j,\gamma}+1).\]
It remains to determine the sign. But this can be done using ($\dagger_2$) 
and the formulae obtained above. Indeed, we have seen that
\begin{align*} 
N_{-\alpha_i,-\beta}&=-c_i(q_{i,\beta}+1),\\
N_{-\alpha_i,-\gamma'}&=-c_i(q_{i,\gamma'}+1),\\
N_{\alpha_j,\gamma'}&=c_j(q_{j,\gamma'}+1).
\end{align*}
Inserting this into ($\dagger_2$), we obtain that 
\[-c_i(q_{i,\beta}+1)N_{\alpha_j,\gamma}=-c_ic_j(q_{i,\gamma'}+1)
(q_{j,\gamma'}+1)\frac{\langle \beta,\beta\rangle}{\langle \alpha,
\alpha\rangle}\frac{\langle \gamma',\gamma'\rangle}{\langle \gamma, 
\gamma\rangle}\]
and, hence,
\[c_jN_{\alpha_j,\gamma}=(q_{i,\beta}+1)^{-1}(q_{i,\gamma'}+1)
(q_{j,\gamma'}+1)\frac{\langle \beta,\beta\rangle}{\langle \alpha,
\alpha\rangle}\frac{\langle \gamma',\gamma'\rangle}{\langle \gamma, 
\gamma\rangle}.\]
Now all terms on the right hand side are positive real numbers and so 
$c_jN_{\alpha_j,\gamma}$ must also be positive. Since $N_{\alpha_j,
\gamma}=\pm (q_{j,\gamma}+1)$, we conclude that $c_jN_{\alpha_j,\gamma}=
q_{j,\gamma}+1$. Since $\be_{\alpha_j}=c_je_j$, this finally yields 
that $[e_j,\be_\gamma]=(q_{j,\gamma}+1)\be_\alpha$, as desired.
\end{proof}

By the discussion in Example~\ref{propnrs1a}, the above result should 
now determine all $N_{\pm\alpha_i,\alpha}$ for $i\in I$ and $\alpha\in\Phi$. 

\begin{rem} \label{canbas5a} Let $i\in I$ and $\alpha\in \Phi^+$
be such that $\alpha+\alpha_i\in \Phi$. Then $(\alpha+\alpha_i)-\alpha_i
=\alpha\in \Phi^+$ and so Lemma~\ref{canbas5} yields that 
$[e_i,\be_\alpha]=(q_{i,\alpha}+1)\be_{\alpha+\alpha_i}$. Thus,
(L2) holds for positive roots.
\end{rem}

\begin{lem} \label{canbas6} Let $\alpha\in\Phi^+$ and $i\in I$ be such that
$\alpha-\alpha_i\in\Phi$. Then $[f_i,\be_\alpha]=(p_{i,\alpha}+1)
\be_{\alpha-\alpha_i}$. Thus, {\rm (L3)} holds for positve roots.
\end{lem}

\begin{proof} Set $\beta:=\alpha-\alpha_i\in\Phi^+$ and write $[f_i,
\be_\alpha]=c\,\be_\beta$, where $c\in \C$. By Lemma~\ref{canbas5}, we 
have $[e_i,\be_{\beta}]=(q_{i,\beta}+1)\be_\alpha$. Next note that 
\begin{align*}
q_{i,\alpha}&=\max\{m\geq 0\mid \alpha-m\alpha_i\in\Phi\}\\&=
\max\{m\geq 0\mid \beta-(m-1)\alpha_i\in\Phi\}\\&=
\max\{m'\geq 0\mid \beta-m'\alpha_i\in\Phi\}+1=q_{i,\beta}+1.
\end{align*}
Hence, we have $[e_i,\be_{\beta}]=q_{i,\alpha}\be_{\alpha}$. 
Consequently, we obtain the identity $[e_i,[f_i,\be_\alpha]]=c[e_i,
\be_\beta]=cq_{i,\alpha}\be_{\alpha}$. Since $\alpha\neq \pm
\alpha_i$, we can apply Remark~\ref{astring}(c). This shows that the 
left hand side of the identity equals $q_{i,\alpha} (p_{i,\alpha}+1) 
\be_\alpha$. Hence, we have $c=p_{i,\alpha}+1$, as desired.
\end{proof}

\begin{lem} \label{canbas7} Let $i\in I$ and $\alpha\in\Phi^-$ be
negative. 
\begin{itemize}
\item[{\rm (a)}] If $\alpha+\alpha_i\in\Phi$, then  $[e_i,\be_\alpha]=
-(q_{i,\alpha}+1)\be_{\alpha+\alpha_i}$.
\item[{\rm (b)}] If $\alpha-\alpha_i\in\Phi$, then $[f_i,\be_\alpha]=
-(p_{i,\alpha}+1)\be_{\alpha-\alpha_i}$.
\end{itemize}
\end{lem}

\begin{proof} (a) Set $\beta:=-\alpha\in\Phi^+$. Then $\beta-\alpha_i
=-(\alpha+\alpha_i)\in\Phi$. Since $\hgt(\beta)\geq 1$, we have
$\hgt(\beta-\alpha_i)\geq 0$ and so $\beta-\alpha_i\in\Phi^+$.
So we can apply Lemma~\ref{canbas6}, which yields that 
\[[f_i,\be_{-\alpha}]=[f_i,\be_\beta]=(p_{i,\beta}+1)\be_{\beta-\alpha_i}
=(q_{i,\alpha}+1)\be_{-(\alpha+\alpha_i)},\]
where the last equality holds by Remark~\ref{canbas0a}. Write 
$\be_{-\alpha_i}=c_i'f_i$ where $0\neq c_i'\in \C$. Since
\[ h_i=h_{\alpha_i}=[\be_{\alpha_i},\be_{-\alpha_i}]=c_ic_i'[e_i,f_i]
=c_ic_i'h_i,\]
we conclude that $c_i'=c_i^{-1}=c_i$. Hence, $\be_{-\alpha_i}=c_if_i$
and so $N_{-\alpha_i,-\alpha}=c_i(q_{i,\alpha}+1)$. Then 
Proposition~\ref{propnrs3} implies that $N_{\alpha_i,\alpha}=
-c_i(q_{i,\alpha}+1)$ and, hence, $[e_i,\be_\alpha]=-(q_{i,\alpha}+1)
\be_{\alpha+\alpha_i}$, as claimed. 

(b) Set again $\beta:=-\alpha\in\Phi^+$. Then $\beta+\alpha_i
=-(\alpha-\alpha_i)\in\Phi$ and so Remark~\ref{canbas5a} yields that
\[ [e_i,\be_{-\alpha}]=[e_i,\be_{\beta}]=(q_{i,\beta}+1)\be_{\beta+
\alpha_i}.\]
Since $\be_{\alpha_i}=c_ie_i$ this yields $N_{\alpha_i,\beta}=
c_i(q_{i,\beta}+1)$, and Proposition~\ref{propnrs3} shows again that 
$N_{-\alpha_i,\alpha}=N_{-\alpha_i,-\beta}=-c_i(q_{i,\beta}+1)$. By 
Remark~\ref{canbas0a}, we have $q_{i,\beta}=p_{i,\alpha}$. 
\end{proof}

Thus, we have found explicit formulae for the structure constants
$N_{\pm \alpha_i,\alpha}$, for all $i\in I$ and $\alpha\in\Phi$, 
summarized as follows:
\begin{alignat*}{2}
[e_i,\be_\alpha] &= +(q_{i,\alpha}+1)\be_{\alpha+\alpha_i}\qquad&&
\mbox{if $\alpha\in\Phi^+$ and $\alpha+\alpha_i\in\Phi$},\\
[e_i,\be_\alpha] &= -(q_{i,\alpha}+1)\be_{\alpha+\alpha_i}\qquad&&
\mbox{if $\alpha\in\Phi^-$ and $\alpha+\alpha_i\in\Phi$},\\
[f_i,\be_\alpha] &= +(p_{i,\alpha}+1)\be_{\alpha-\alpha_i}\qquad&&
\mbox{if $\alpha\in\Phi^+$ and $\alpha-\alpha_i\in\Phi$},\\
[f_i,\be_\alpha] &= -(p_{i,\alpha}+1)\be_{\alpha-\alpha_i}\qquad&&
\mbox{if $\alpha\in\Phi^-$ and $\alpha-\alpha_i\in\Phi$}.
\end{alignat*}
Hence, the signs are not yet right as compared to the desired formulae 
in Theorem~\ref{canbas}. To fix this, we define for $\alpha \in\Phi$:
\begin{center}
\fbox{$\displaystyle \be_\alpha^+:=\left\{\begin{array}{cl} \be_\alpha 
& \quad \mbox{if $\alpha\in\Phi^+$},\\ (-1)^{\hgt(\alpha)}\be_\alpha & 
\quad \mbox{if $\alpha\in\Phi^-$}.\end{array}\right.$}
\end{center}
We claim that (L1), (L2), (L3) in Theorem~\ref{canbas} hold. 
First consider (L2). Let $i\in I$ and $\alpha\in\Phi$ be such that
$\alpha+ \alpha_i\in\Phi$. If $\alpha\in\Phi^+$, then $\be_\alpha^+=
\be_\alpha$ and the required formula holds. If $\alpha\in\Phi^-$, then 
$[e_i,\be_\alpha^+]=(-1)^{\hgt(\alpha)}[e_i,\be_\alpha]=
-(-1)^{\hgt(\alpha)}(q_{i,\alpha}+1)\be_{\alpha+\alpha_i}$; so the desired 
formula holds again, since $\be_{\alpha+\alpha_i}^+=(-1)^{\hgt(\alpha+
\alpha_i)}\be_{\alpha+\alpha_i}$. The argument for (L3) is
analogous. Now consider (L1). In the proof of Lemma~\ref{canbas7}(a),
we already saw that $\be_{-\alpha_i}=c_if_i$ for $i\in I$. Hence,
we obtain 
\begin{align*}
[e_i,\be_{-\alpha_i}^+]&=-[e_i,\be_{-\alpha_i}]=-c_i[e_i,f_i]=
-c_ih_i, \\ [f_i,\be_{\alpha_i}^+]\,&=+c_i[f_i,e_i]\;\;\,=
-c_i[e_i,f_i]=-c_ih_i.
\end{align*}
Thus, (L3) holds and the proof of Theorem~\ref{canbas} is complete. As
a by-product, we also obtain:

\begin{cor} \label{canbash} There is a collection of elements
$\{\be_\alpha^+\mid \alpha\in\Phi\}$ satisfying {\rm (L1)--(L3)} in 
Theorem~\ref{canbas} and such that 
\[[\be_\alpha^+,\be_{-\alpha}^+]=(-1)^{\hgt(\alpha)}h_\alpha\qquad
\mbox{for all $\alpha\in\Phi$}.\] 
We have $\be_{\alpha_i}^+=c_ie_i$ and $\be_{-\alpha_i}^+=-c_i f_i$,
with $c_i\in \{\pm 1\}$ for all $i\in I$. If~$A$ is indecomposable, then 
such a collection $\{\be_\alpha^+\mid \alpha\in \Phi\}$ is unique up to 
a global sign, that is, if $\{\be_\alpha' \mid\alpha\in \Phi\}$ is 
another collection satisfying {\rm (L1)--(L3)} and the above identity, 
then there exists $\xi\in \{\pm 1\}$ such that $\be_{\alpha}'=\xi
\be_{\alpha}^+$ for all $\alpha\in \Phi$.
\end{cor}

\begin{proof} Since $[\be_\alpha, \be_{-\alpha}]=h_\alpha$, the formula 
for $[\be_\alpha^+,\be_{-\alpha}^+]$ is clear by the definition of 
$\be_\alpha^+$ and the fact that $h_{-\alpha}=-h_\alpha$ for all $\alpha
\in \Phi$. Now assume that $A$ is indecomposable and let $\{\be_\alpha' 
\mid \alpha\in \Phi\}$ be another collection satisyfing (L1)--(L3) 
and the above identity. As discussed in Remark~\ref{canbas2}, there 
exists $0\neq \xi\in \C$ such that $\be_\alpha'=\xi\be_\alpha^+$ for 
all $\alpha\in\Phi$.  But then $(-1)^{\hgt(\alpha)}h_\alpha=
[\be_\alpha',\be_{-\alpha}']= \xi^2 [\be_\alpha^+,\be_{-\alpha}^+]=
\xi^2(-1)^{\hgt(\alpha)} h_\alpha$ and so $\xi=\pm 1$, as desired.
\end{proof}

\begin{xca} \label{lucanAn} Let $\cL=\slm_n(\C)$ and $\fh\subseteq \cL$ be 
the usual abelian subalgebra of diagonal matrices. Let $\{h_i,e_i,f_i\mid 
1\leq i\leq n-1\}$ be as in Example~\ref{cartsln}; also recall that
\[ \Phi=\{\varepsilon_i-\varepsilon_j \mid 1\leq i,j\leq n,i\neq j\},
\qquad \cL_{\varepsilon_i-\varepsilon_j}=\langle e_{ij}\rangle_\C.\]
We set $\be_{\alpha}^+:= (-1)^je_{ij}$ for $\alpha=\varepsilon_i-
\varepsilon_j$, $i\neq j$. Show that the collection $\{\be_\alpha^+\mid 
\alpha\in\Phi\}$ satisfies the conditions in Corollary~\ref{canbash}.
In particular, we have $\be_{\alpha_i}^+=-(-1)^ie_i$ and 
$\be_{-\alpha_i}^+=(-1)^i f_i$ for $1\leq i\leq n-1$; furthermore,
$h_i^+=[e_i,\be_{-\alpha_i}^+]=(-1)^ih_i$.

\smallskip
\noindent {\footnotesize [{\it Hint}. Just verify (L1), (L2), (L3) 
for $\be_\alpha^+$ as defined above.]}
\end{xca}

\begin{exmp} \label{lucanG2} Let $\cL$ be a Lie algebra of
Cartan--Killing type as in Example~\ref{hypog2}, with structure matrix 
$A$ of type $G_2$. (We do not know yet that such an algebra exists.)
We already constructed elements $\be_\alpha\in \cL_\alpha$ for the
twelve roots in $\Phi$. If we now define $\be_\alpha^+$ as
above, then $\{\be_\alpha^+\mid \alpha\in \Phi\}$ is a collection
of elements as in Corollary~\ref{canbash}, where $\be_{\alpha_1}^+=e_1$,
$\be_{\alpha_2}^+=-e_2$, $\be_{-\alpha_1}^+=-f_1$ and $\be_{-\alpha_2}^+=
f_2$. We leave it to the reader to check that, indeed, the formulae in 
Remark~\ref{canbas1} hold.
\end{exmp}

We now establish an important consequence of Theorem~\ref{canbas}. 
Let also $\tilde{\cL}$ be a Lie algebra of Cartan--Killing type, that is,
there is a subalgebra $\tilde{\fh}\subseteq \tilde{\cL}$ and a subset 
$\tilde{\Delta}=\{\tilde{\alpha}_i\mid i\in \tilde{I}\}$ (for some 
finite index set $\tilde{I}$) such that the conditions in 
Definition~\ref{defTD} hold. Let $\tilde{A}=(\tilde{a}_{ij})_{i,j
\in \tilde{I}}$ be the corresponding structure matrix. 

\begin{thm}[Isomorphism Theorem] \label{isothm} \nmi{}{Isomorphism Theorem} 
With the above notation, assume that $I=\tilde{I}$ and $A=\tilde{A}$. Then 
there is a unique isomorphism of Lie algebras $\varphi\colon \cL\rightarrow 
\tilde{\cL}$ such that $\varphi(e_i)=\tilde{e}_i$ and $\varphi(f_i)=
\tilde{f}_i$ for all $i\in I$, where $\{e_i,f_i\mid i\in I\}$ and 
$\{\tilde{e_i}, \tilde{f_i}\mid i\in I\}$ are Chevalley generators 
for $\cL$ and $\tilde{\cL}$, respectively (as in Remark~\ref{astring0}).
\end{thm}

\begin{proof} The uniqueness of $\varphi$ is clear since $\cL=\langle e_i, 
f_i\mid i\in I\rangle_{\text{alg}}$; see Proposition~\ref{genlie}. The 
problem is to prove the existence of $\varphi$. Let $\Phi\subseteq \fh^*$
be the set of roots of $\cL$ and $\tilde{\Phi}\subseteq \tilde{\fh}^*$ be 
the set of roots of $\tilde{\cL}$. Since $A=\tilde{A}$, the discussion 
in Remark~\ref{explicit} shows that we have a canonical bijection $\Phi
\stackrel{\sim}{\longrightarrow}\tilde{\Phi}$, $\alpha\mapsto 
\tilde{\alpha}$, given as follows. If $\alpha=\sum_{i\in I}n_i\alpha_i
\in\Phi$ (with $n_i\in\Z$), then $\tilde{\alpha}=\sum_{i\in I}n_i
\tilde{\alpha_i}\in\tilde{\Phi}$. Then this bijection has
the following property: for any $\alpha,\beta \in\Phi$, we have 
\begin{equation*}
\alpha+\beta\in \Phi \qquad \Leftrightarrow\qquad \tilde{\alpha}+
\tilde{\beta} \in\tilde{\Phi}.\tag{$\heartsuit$}
\end{equation*}
Now let us fix a collection of signs $c_i=\pm 1$ ($i\in I$) as in
Definition~\ref{startproof}. This yields the basis 
\[ \bB =\{h_i^+\mid i \in I\} \cup\{\be_\alpha^+\mid \alpha\in \Phi\}\]
for $\cL$, as in Corollary~\ref{canbash}. Similarly, we obtain the basis 
\[ \tilde{\bB}=\{\tilde{h}_i^+\mid i \in I\} \cup\{\tilde{\be}_\alpha^+ 
\mid \tilde{\alpha}\in \tilde{\Phi}\}\]
for $\tilde{\cL}$. Now define a (bijective) linear map $\varphi\colon 
\cL\rightarrow\tilde{\cL}$ by 
\[\varphi(h_i^+):=\tilde{h}_i^+ \quad (i\in I)\qquad\mbox{and}\qquad 
\varphi(\be_\alpha^+):=\tilde{\be}_{\tilde{\alpha}}^+\quad 
(\alpha\in\Phi).\] 
We have $\be_{\alpha_i}^+=c_ie_i$ and $\be_{-\alpha_i}^+=-c_if_i$ 
for all $i\in I$; similarly, $\tilde{\be}_{\tilde{\alpha}_i}^+=c_ie_i$ 
and $\tilde{\be}_{-\tilde{\alpha}_i}^+=-c_if_i$ for all $i\in I$. 
Consequently, we have 
\[\varphi(e_i)=\tilde{e}_i \qquad\mbox{and}\qquad \varphi(f_i)=
\tilde{f}_i \qquad \mbox{for all $i\in I$}.\]
Furthermore, let $i\in I$ and $\alpha\in\Phi$ be such that $\alpha+
\alpha_i\in\Phi$. By ($\heartsuit$), we also have $\tilde{\alpha}+
\tilde{\alpha}_i\in \tilde{\Phi}$ and 
\begin{align*}
q_{i,\alpha}&=\max\{m \geq 0 \mid \alpha-m\alpha_i\in \Phi\}\\&=
\max\{m \geq 0\mid \tilde{\alpha}-m\tilde{\alpha}_i\in \tilde{\Phi}\}
=q_{i,\tilde{\alpha}}.
\end{align*} 
Similarly, if $\alpha-\alpha_i\in\Phi$, then $\tilde{\alpha}-
\tilde{\alpha}_i\in\tilde{\Phi}$ and $p_{i,\alpha}=p_{i,\tilde{\alpha}}$.
Hence, (L2) shows that the matrix of $\ad_\cL(e_i)\colon \cL\rightarrow \cL$ 
with respect to the basis $B$ is equal to the matrix of $\ad_{\tilde{\cL}}
(\tilde{e}_i)\colon \tilde{\cL} \rightarrow \tilde{\cL}$ with respect to 
the basis $\tilde{B}$; by (L3), similar statements also hold for 
$\ad_\cL(f_i)$ and $\ad_{\tilde{\cL}}(\tilde{f}_i)$. Since $\varphi$ is 
linear, this implies that 
\begin{align*}
\varphi([e_i,y])&=[\tilde{e}_i,\varphi(y)]=[\varphi(e_i),\varphi(y)],\\
\varphi([f_i,y])&=[\tilde{f}_i,\varphi(y)]=[\varphi(f_i),\varphi(y)]
\end{align*}
for all $i\in I$, $y\in \cL$. Since $\cL=\langle e_i,f_i\mid i\in I
\rangle_{\text{alg}}$, it follows that $\varphi([x,y])=[\varphi(x),
\varphi(y)]$ for all $x,y\in \cL$; see Exercise~\ref{xcagenerator}(d). 
So $\varphi$ is an isomorphism of Lie algebras.
\end{proof}

\begin{exmp} \label{stdbase3} Let $\tilde{\alpha}_i:=-\alpha_i$ for all 
$i \in I$. Then $(\cL,\fh)$ also is of Cartan--Killing type with
respect to $\tilde{\Delta}:=\{\tilde{\alpha}_i\mid i \in I\}$. (This was 
already used in the proof of Theorem~\ref{mainthm1}(a).) We have 
$\tilde{h}_i=h_{\tilde{\alpha}_i}=h_{-\alpha_i}=-h_i$ for $i\in I$. 
Hence, the structure matrix $\tilde{A}$ of $\cL$ with respect to 
$\tilde{\Delta}$ is the same as the original structure matrix $A$ of 
$\cL$ with respect to $\Delta$. Finally, if we set $\tilde{e}_i:=f_i$ 
and $\tilde{f}_i:=e_i$ for $i\in I$, then $\{\tilde{e}_i,\tilde{f}_i
\mid i \in I\}$ are Chevalley generators for $\cL$ with respect to 
$\tilde{\Delta}$. So Theorem~\ref{isothm} shows that there is a unique 
automorphism of Lie algebras $\omega \colon \cL \rightarrow \cL$ such that:
\[ \omega(e_i)=f_i, \qquad \omega(f_i)=e_i, \qquad 
\omega(h_i)=-h_i \qquad (i\in I).\]
This is called the \nm{Chevalley involution} of $\cL$; we have $\omega^2=
\id_\cL$. (Note that, alternatively, one can also apply the
whole argument with $\tilde{e}_i:=-f_i$ and $\tilde{f}_i:=-e_i$ for 
$i\in I$; we still have $\tilde{h}_i=-h_i=[\tilde{e}_i,\tilde{f}_i]$.)
\end{exmp}

\begin{thm}[Cf.\ Chevalley \protect{\cite[\S I]{Ch}}] \label{xcaomega} 
Let $\{\be_\alpha^+\mid \alpha \in \Phi\}$ be a collection as in 
Corollary~\ref{canbash}. Then the following hold.
\begin{itemize}
\item[{\rm (a)}] We have $\omega (\be_\alpha^+)=-\be_{-\alpha}^+$ for all 
$\alpha\in \Phi$. 
\item[{\rm (b)}] Let $\alpha, \beta\in\Phi$ be such that $\alpha+\beta
\in\Phi$. Then $[\be_\alpha^+,\be_\beta^+]=\pm (q+1)\be_{\alpha+\beta}^+$, 
where $q\geq 0$ is defined as in Lemma~\ref{astring2}.
\end{itemize}
\end{thm}

\begin{proof} (a) Let $\alpha\in\Phi^+$. We show the assertion
by induction on $\hgt(\alpha)$. If $\hgt(\alpha)=1$, then $\alpha=
\alpha_i$ for some $i\in I$. We have $\be_{\alpha_i}^+=c_ie_i$ and 
$\be_{-\alpha_i}^+=-c_if_i$, where $c_i\in\{\pm 1\}$ for all $i\in I$.
Hence, using Example~\ref{stdbase3}, we obtain $\omega(\be_{\alpha_i}^+)
=c_i\omega(e_i)=c_if_i=-\be_{-\alpha_i}^+$, as required. Now let 
$\hgt(\alpha)>1$. By the Key Lemma~\ref{keylem}, there exists some 
$i\in I$ such that $\beta:=\alpha-\alpha_i\in\Phi^+$. We have 
$\hgt(\beta)=\hgt(\alpha)-1$ and so $\omega(\be_{\beta}^+)=
-\be_{-\beta}^+$, by induction. By condition (L1) in 
Theorem~\ref{canbas}, we have $[e_i,\be_\beta^+]=(q_{i,\beta}+1)
\be_\alpha^+$. Applying $\omega$ yields that
\[ (q_{i,\beta}+1)\omega(\be_\alpha^+)=\omega\bigl([e_i,
\be_\beta^+]\bigr)=[\omega(e_i),\omega(\be_\beta^+)]=-[f_i,
\be_{-\beta}^+].\]
Now, we have $-\beta-\alpha_i=-\alpha\in\Phi$ and so condition (L2)
in Theorem~\ref{canbas} yields that $[f_i,\be_{-\beta}^+]=
(p_{i,-\beta}+1)\be_{-\alpha}^+$. Hence, we deduce that
$\omega(\be_\alpha^+)=-\be_{-\alpha}^+$, since $p_{i,-\beta}=
q_{i,\beta}$ as pointed out in Remark~\ref{canbas0a}. Thus, the
assertion holds for all $\alpha\in\Phi^+$. But, since $\omega^2=\id_\cL$,
we then also have $\omega(\be_{-\alpha}^+)=\omega\bigl(-\omega
(\be_{\alpha}^+)\bigr)=-\omega^2(\be_\alpha^+)=-\be_\alpha^+$, as required.

(b) We would like to use Proposition~\ref{propnrs3}, but we can not do
that directly because the condition in Remark~\ref{propnrs0}(a) does not 
hold for the collection $\{\be_\alpha^+\mid \alpha\in\Phi\}$. So we go
back to the collection $\{0\neq \be_\alpha\in \cL_\alpha\mid\alpha\in\Phi\}$ 
in Definition~\ref{startproof}; thus, 
\[ \be_\alpha:=\left\{\begin{array}{cl} \be_{\alpha}^+ & \quad \mbox{if 
$\alpha\in\Phi^+$},\\ (-1)^{\hgt(\alpha)}\be_{\alpha}^+ & \quad 
\mbox{if $\alpha\in \Phi^-$}.  \end{array}\right.\]
Then $[\be_\alpha,\be_{-\alpha}]=h_\alpha$ for all $\alpha\in\Phi$.
By (a), we also have the formula:
\[ \omega(\be_\alpha)=-(-1)^{\hgt(\alpha)}\be_{-\alpha} \qquad\mbox{for 
all $\alpha\in\Phi$}.\]
Let again $N_{\alpha,\beta}$ be the structure constants with respect to 
$\{\be_\alpha\mid \alpha\in\Phi\}$, as in Section~\ref{sec1a6}. Writing 
$[\be_\alpha,\be_\beta]=N_{\alpha,\beta} \be_{\alpha+\beta}$, we certainly 
have $[\be_\alpha^+, \be_\beta^+]=\pm N_{\alpha,\beta}\be_{\alpha+
\beta}^+$. So it suffices to show that $N_{\alpha,\beta}=\pm (q+1)$.
This is seen as follows. Using the above formula for $\omega$, we obtain
\[\omega\bigl([\be_\alpha,\be_\beta])=N_{\alpha,\beta} \omega
(\be_{\alpha+\beta})=-(-1)^{\hgt(\alpha+\beta)} N_{\alpha,\beta}
\be_{-(\alpha+\beta)}.\]
On the other hand, we can also evaluate the left hand side as follows.
\begin{align*} 
\omega\bigl([\be_\alpha,\be_\beta])&=[\omega(\be_\alpha),\omega
(\be_\beta)]=(-1)^{\hgt(\alpha)+\hgt(\beta)}[\be_{-\alpha},\be_{-\beta}]\\
&=(-1)^{\hgt(\alpha)+\hgt(\beta)}N_{-\alpha,-\beta}\be_{-(\alpha-\beta)}.
\end{align*}
Hence, we conclude that $N_{-\alpha,-\beta}=-N_{\alpha,\beta}$ and so 
Proposition~\ref{propnrs3} implies that $N_{\alpha,\beta}^2=(q+1)^2$.  
Thus, $N_{\alpha,\beta}=\pm (q+1)$, as claimed.
\end{proof}

\begin{xca} \label{xcaothero} This exercise provides a more direct 
construction of the Chevalley involution in Example~\ref{stdbase3}, 
without reference to the Isomorphism Theorem. For this purpose, consider 
the basis $\bB$ of~$\cL$ in Remark~\ref{canbas1} and define a linear map 
$\tilde{\omega}\colon \cL\rightarrow \cL$ by 
\[ \tilde{\omega}(h_j^+):=-h_j^+ \quad (j\in I)\qquad\mbox{and}\qquad
\tilde{\omega}(\be_\alpha^+):=-\be_{-\alpha}^+ \quad (\alpha\in \Phi).\]
Use (L1), (L2), (L3) in Theorem~\ref{canbas} to verify that 
\[\tilde{\omega}\circ \ad_\cL(e_i)=\ad_\cL(f_i)\circ \tilde{\omega}\quad
\mbox{and}\quad \tilde{\omega}\circ \ad_\cL(h_i^+)=-\ad_\cL(h_i^+)\circ 
\tilde{\omega}\]
for all $i\in I$. Then use Exercise~\ref{xcagenerator}(d) to deduce that 
$\tilde{\omega}$ is a Lie algebra automorphism and that $\tilde{\omega}$
equals $\omega$ in Example~\ref{stdbase3}; note that $\tilde{\omega}(e_i)
=f_i$ for $i\in I$. 
\end{xca}

\begin{xca}[Difficult!] \label{isothmexp2} Let us fix $j\in I$ and
consider the simple reflection $s_j\in W$. Let $\Delta'=\{\alpha_i'\mid 
i\in I\}$ where $\alpha_i'=s_j(\alpha_i)$ for all $i\in I$. We set 
$h_i':=h_i-a_{ij}h_j$ and
\begin{align*}
e_j'&:=-f_j, \qquad\qquad f_j':=-e_j,\\
e_i'&:=(-a_{ji})!^{-1}\ad_\cL(e_j)^{-a_{ji}}(e_i)\;\qquad\qquad
\mbox{if $i\neq j$},\\ f_i'&:=(-1)^{a_{ji}}(-a_{ji})!^{-1}
\ad_\cL(f_j)^{-a_{ji}}(f_i)\quad\, \mbox{if $i\neq j$}.
\end{align*}
(These formulae can be found in \cite[Appendix A6]{CP}.) Recall that 
$a_{ji}\leq 0$ if $i\neq j$. For example, if $i\neq j$ and $a_{ji}=-1$, 
then $s_j(\alpha_i)=\alpha_i-a_{ji}\alpha_j=\alpha_i+\alpha_j$ and
\[ e_i'=[e_j,e_i]\in \cL_{s_j(\alpha_i)}, \qquad f_i'=-[f_j,f_i]\in 
\cL_{-s_j(\alpha_i)};\]
if $i\neq j$ and $a_{ji}=-2$, then $s_j(\alpha_i)=\alpha_i-a_{ji}\alpha_j=
\alpha_i+2\alpha_j$ and 
\[ e_i'=\textstyle{\frac{1}{2}}[e_j,[e_j,e_i]]\in \cL_{s_j(\alpha_i)},
\qquad f_i'=\textstyle{\frac{1}{2}}[f_j,[f_j,f_i]]\in\cL_{-s_j(\alpha_i)}.\]
Then show that $\cE'=\{e_i',f_i'\mid i \in I\}$ is a system of
Chevalley generators for $\cL$ with respect to $(\fh,\Delta')$. Hence, there 
exists a unique automorphism of Lie algebras $\eta_j\colon \cL\rightarrow 
\cL$ such that $\eta_j(e_i)=e_i'$, $\eta_j(f_i)=f_i'$ and $\eta_j(h_i)=h_i'$ 
for all $i\in I$. Thus, $\eta_j$ may be regarded as an automorphism which 
``lifts'' the simple reflection $s_j\in W$ to $\cL$. (Far-reaching 
generalizations of these maps and formulae are contained in Lusztig 
\cite[Part~VI]{L6}; see also Jantzen \cite[Chap.~8]{Ja}.)
\end{xca}

\section*{Notes on Chapter~\ref{chap2}}

Lemmas~\ref{wsdlem4} and \ref{wsdlem3} (which lead to the definition
of the structure matrix of $\cL$) contain suitably adapted standard 
arguments from the theory of semisimple Lie algebras; see, e.g., the 
proofs of \cite[8.3(e)]{H} and \cite[4.20]{Ca3}. As already mentioned, 
the idea of developing the theory from a set of axioms in the spirit of 
Definition~\ref{defTD} is taken from Moody--Pianzola \cite{MP}. In 
Section~\ref{sec1a3}, we point to an algorithmic approach around root 
systems, Weyl groups and Lie algebras, which has a long tradition in Lie 
theory and has proved to be extremely powerful in various situations; 
see, e.g., De Graaf \cite{graaf} and further references there. (We will 
say more about this in Section~\ref{sec3comp} below.) 
Proposition~\ref{forestg} appears in Erdmann--Wildon \cite[Prop.~13.5]{EW},
Humphreys \cite[\S 11.4(3)]{H}; see also Bourbaki \cite[Ch.~V, \S 4, 
no.~8]{B} for a similar statement in the more general setting of finite 
reflection groups.

The early introduction of the automorphisms $x_i(t)$ and $y_i(t)$ in 
Section~\ref{sec1a4} is convenient because these will later be used to
construct Chevalley groups. The statement in Lemma~\ref{mainideal1} (that 
every non-trivial ideal of $\fg$ intersects $\fh$ non-trivially) is a 
crucial property in the theory Kac--Moody Lie algebras; see 
Kac \cite[Prop.~1.4]{K}. The results on structure constants in 
Section~\ref{sec1a6} are standard and can be found, for example, in 
Bourbaki \cite[Ch.~VIII, \S 2, no.~4]{B78} and Carter \cite[\S 4.1]{Ca1}. 
For further properties of these constants, see Casselman \cite{Cass}, 
\cite{Cass1}, \cite{Cass2} and Tits \cite{Ti}. 
 
The proof of Lusztig's Theorem~\ref{canbas} would have been somewhat 
easier if the Isomorphism Theorem~\ref{isothm} and the existence of the
Chevalley involution $\omega\colon \cL\rightarrow \cL$ (see 
Example~\ref{stdbase3}) were known in advance. Here, we first work a 
little harder to get Theorem~\ref{canbas}, but then Theorem~\ref{isothm} 
is a relatively easy consequence. We also remark that a function 
$i\mapsto c_i$ satisfying the condition in Remark~\ref{canbas0}($*$) 
already appeared in the work of Rietsch \cite[\S 4]{KR}.

Exercise~\ref{lucanAn} contains an explicit description of the 
canonical basis for $\cL=\slm_n(\C)$. Similar explicit descriptions 
are determined by Lang \cite[Chap.~2]{Lang} for all the classical Lie 
algebras $\cL=\gom_n(Q_n,\C)$.

The Isomorphism Theorem is a classical result which is covered in most 
textbooks on Lie algebras (independently of Lusztig's theorem); see, for 
example, Jacobson \cite[Chap.~IV, \S 4]{Jac}, Serre \cite[Chap.~V, 
\S 5]{S}, Humphreys \cite[\S 14.2]{H}, Samelson \cite[\S 2.9]{Sam}, or 
DeGraaf \cite[\S 5.11]{graaf}. In Carter \cite[\S 7.2]{Ca3} (see also 
Stewart \cite[Chap.~8]{ist}), the proof is based on the consideration of 
``special'' and ``extraspecial'' pairs of roots, which are also often 
used for algorithmic purposes (see, e.g., \cite[\S 3]{Coet}), but which 
we do not need here at all. 

There is also a different proof of Lusztig's Theorem~\ref{canbas}: 
in \cite{G1}, we explicitly construct a specific Lie algebra of 
Cartan--Killing type with structure matrix $A$ and with a basis such 
that (L1)--(L3) in Theorem~\ref{canbas} hold. Then the Isomorphism 
Theorem (which would thus have to be proved differently, as indicated 
above) shows that Theorem~\ref{canbas} holds in general. The proof in 
\cite{G1} has the advantage that it gives a construction of a Lie 
algebra with a given~$A$. Here, we will obtain this existence
result in Section~\ref{sec3a3}.

The involution $\omega\colon \cL\rightarrow \cL$ in Example~\ref{stdbase3}
is used to construct the compact real form of $\cL$; see, e.g., 
\cite[\S 2.10]{Sam} for further details. The proof of 
Theorem~\ref{xcaomega}(b), based on the existence of $\omega\colon 
\fg\rightarrow \fg$ and the identity in Proposition~\ref{propnrs3}, 
essentially follows the original argument of Chevalley 
\cite[Th\'eor\`eme~1 (p.~24)]{Ch}; see also Carter \cite[\S 4.2]{Ca1}.
Somewhat different arguments can be found in Bourbaki \cite[Ch.~VIII, 
\S 12, no.~6, Cor.~4]{B78} and Samelson \cite[\S 2.9]{Sam}. 


\chapter{Generalised Cartan matrices} \label{chap3}

\addtocounter{footnote}{6}

In the previous chapter we have seen that a Lie algebra $\cL$ of 
Cartan--Killing type is determined (up to isomorphism) by its structure 
matrix $A=(a_{ij})_{i,j\in I}$. The entries of $A$ are integers, we have
$a_{ii}=2$ and $a_{ij}\leq 0$ for $i\neq j$; furthermore, $a_{ij}<0
\Leftrightarrow a_{ji}<0$. In Section~\ref{sec3a1} we show that 
every indecomposable matrix satisfying those conditions has one of 
three possible types: (FIN), (AFF) or (IND). There is a complete 
classification of all such matrices of types (FIN) and (AFF). The 
structure matrix $A$ of $\cL$ does turn out to be of type (FIN) and, 
hence, it is encoded by one of the graphs in the famous list of Dynkin 
diagrams of type $A_n$, $B_n$, $C_n$, $D_n$, $G_2$, $F_4$, $E_6$, 
$E_7$ and $E_8$. 

Once the results in Section~\ref{sec3a1} are established, the central 
theme of this chapter is as follows. We start with an arbitrary matrix 
$A$ as above, of type (FIN). Then we can construct the following objects:

1) An abstract root system $\Phi$. In Section~\ref{sec1a3} we already
made first steps in that direction, and presented a {\sf Python} program 
to determine~$\Phi$ from~$A$. This will be further developed in
Section~\ref{sec3a2}.

2) A Lie algebra $\cL$ of Cartan--Killing type with structure matrix~$A$
and root system $\Phi$. This will be done by a process that reverses
the construction of Lusztig's canonical basis; see Section~\ref{sec3a3}.

3) A Chevalley group $G_K(\fg)$ ``of type $\cL$'', first over $\C$ and then
over any field $K$. Here we follow Lusztig's simplified construction using 
the canonical basis of $\cL$; see Section~\ref{sec3a7}.

We emphasise the fact that the constructions are by means of purely 
combinatorial procedures, which do not involve any other ingredients (or 
choices) and, hence, can also be implemented on a computer: the single 
input datum for the computer programs is the matrix~$A$ (plus the field~$K$
for the Chevalley groups). We present a specific computer algebra package 
with these features in Section~\ref{sec3comp}.

\section{Classification} \label{sec3a1}

Let $I$ be a finite, non-empty index set. We consider matrices
$A=(a_{ij})_{i,j\in I}$ with entries in $\R$ satisfying the following 
two conditions:
\begin{itemize}
\item[(C1)] $a_{ij}\leq 0$ for all $i\neq j$ in $I$;
\item[(C2)] $a_{ij}\neq 0\Leftrightarrow a_{ji}\neq 0$, for all $i,j\in I$.
\end{itemize}
Examples of such matrices are the structure matrices of Lie algebras of
Cartan--Killing type; see Corollary~\ref{gencart0}. One of our aims will 
be to find the complete list of all possible such structure matrices.
For this purpose, it will be convenient to first work in a more
general setting, where we only assume that (C1) and (C2) hold.

In analogy to Definition~\ref{cindec1}, we say that $A$ is 
\textit{indecomposable} if there is no partition $I=I_1\sqcup I_2$ (where 
$I_1,I_2\subsetneqq I$ and $I_1\cap I_2=\varnothing$) such that $a_{ij}=
a_{ji}=0$ for all $i\in I_1$ and $j\in I_2$. 

Some further notation. Let $u=(u_i)_{i\in I}\in \R^I$. We write $u\geq 0$ 
if $u_i\geq 0$ for all $i\in I$; we write $u>0$ if $u_i>0$ for all 
$i\in I$. Finally, $Au\in \R^I$ is the vector with $i$-th component given 
by $\sum_{j\in I}a_{ij}u_j$ (usual product of $A$ with $u$ regarded as a 
column vector). 

\begin{lem} \label{kac1} Assume that $A$ satisfies {\rm (C1)}, {\rm (C2)} 
and is indecomposable. If $u\in \R^I$ is such that $u\geq 0$, $Au\geq 0$, 
then $u=0$ or $u>0$.
\end{lem}

\begin{proof} Let $I_1:=\{i\in I\mid u_i=0\}$ and $I_2:=\{i\in I\mid u_i>
0\}$. Then $I=I_1\cup I_2$, $I_1\cap I_2=\varnothing$. Let $i\in I_1$ 
and $v_i$ be the $i$-th component of $Au$; by assumption, $v_i\geq 0$. 
On the other hand, $v_i=\sum_{j\in I} a_{ij}u_j=\sum_{j\in I_2} a_{ij}u_j$
where all terms in the sum on the right hand side are $\leq 0$ since 
$A$ satisfies (C1) and $u_j>0$ for all $j\in I_2$; furthermore, if $a_{ij}
<0$ for some $j\in I_2$, then $v_i<0$, contradiction to $v_i\geq 0$. So 
we must have $a_{ij}=0$ for all $i\in I_1,j\in I_2$. Since $A$ satisfies 
(C2), we also have $a_{ji}=0$ for all $i\in I_1,j\in I_2$. Since $A$ is 
indecomposable, either $I_1=I$ (and so $u=0$) or $I_2=I$ (and so $u>0$). 
\end{proof}

\begin{thm}[Vinberg \protect{\cite{Vi}}] \label{kac2} Assume that $A$ 
satisfies {\rm (C1)}, {\rm (C2)} and is indecomposable. Let $\cK_A:=
\{u\in \R^I\mid Au\geq 0\}$. Then exactly one of the following three 
conditions holds.
\begin{itemize}
\item[{\rm (FIN)}] $\{0\}\neq \cK_A\subseteq \{u\in \R^I\mid u>0\}\cup
\{0\}$.
\item[{\rm (AFF)}] $\cK_A=\{u\in \R^I\mid Au=0\}=\langle u_0\rangle_\R$
where $u_0>0$. 
\item[{\rm (IND)}] $\cK_A\cap\{u\in \R^I\mid u\geq 0\}=\{0\}$.
\end{itemize}
Accordingly, we say that $A$ is of \nmi{finite}{finite type}, 
\nmi{affine}{affine type} or \nmi{indefinite}{indefinite type} type.
\end{thm}

\begin{proof} First we show that the three conditions are disjoint.
If (FIN) or (AFF) holds, then there exists some $u\in \R^I$ such that $u>0$ 
and $Au\geq 0$. Hence, (IND) does not hold. If (AFF) holds, then 
there exists some $u\in \R^I$ such that $u>0$ and $Au=0\geq 0$. But
then also $A(-u)\geq 0$ and so (FIN) does not hold. Hence, the conditions
are indeed disjoint. It remains to show that we are always in one of the
three cases. Assume that (IND) does not hold. Then there exists some $0\neq 
v\in \cK_A$ such that $v\geq 0$. By Lemma~\ref{kac1}, we have $v>0$. We 
want to show that (FIN) or (AFF) holds. Assume that (FIN) does not hold. 
Since $\cK_A\neq\{0\}$, this means that there exists $0\neq u\in \cK_A$ 
such that $u_l\leq 0$ for some $l\in I$. We have $v>0$ and so we can consider 
the ratios $u_i/v_i$ for $i\in I$. Let $j\in I$ be such that $u_j/v_j\leq 
u_i/v_i$ for all $i\in I$. If $u_j\geq 0$, then $u_i\geq 0$ for all 
$i\in I$ and so $u\geq 0$. But then Lemma~\ref{kac1} would imply that $u>0$,
contradiction to our choice of~$u$. Hence, $u_j<0$ and so $s:=-u_j/v_j>0$. 
Now let us look at the vector $u+sv$; its $i$-th component is 
\[(u+sv)_i=u_i+sv_i=v_i(u_i/v_i-u_j/v_j)\;\left\{\begin{array}{rl} =0 &
\mbox{if $i=j$},\\ \geq 0 & \mbox{if $i\neq j$}.
\end{array}\right.\] 
Hence, we have $u+sv\geq 0$ and $A(u+sv)=Au+sAv\geq 0$. By
Lemma~\ref{kac1}, either $u+sv=0$ or $u+sv>0$. But $(u+sv)_j=0$
and so we must have $u+sv=0$, that is, $u=-sv$. But then $0\leq Au =(-s)Av
\leq 0$ (since $s>0$ and $Av\geq 0$) and so $Av=Au=0$. 

Finally, consider any $0\neq w\in \cK_A$. Again, let $j\in J$ be such 
that $w_j/v_j\leq w_i/v_i$ for all $i\in I$, and set $t:=-w_j/v_j$. As 
above, we see that $w+tv\geq 0$ and $(w+tv)_j=0$. Furthermore, $A(w+tv)=
Aw+tAv=Aw\geq 0$ (since $Av=0$). So Lemma~\ref{kac1} implies that either 
$w+tv>0$ (which is not the case) or $w+tv=0$; hence, $w=-tv\in \langle v 
\rangle_\R$. So $\cK_A\subseteq \langle v\rangle_\R \subseteq \{z\in 
\R^I\mid Az=0\}$ and the right hand side is contained in $\cK_A$. Hence, 
(AFF) holds with $u_0=v$. 
\end{proof}

\begin{cor} \label{kac3} Let $A$ be as in Theorem~\ref{kac2}. If $A$
is of finite or affine type, then there exists some $u\in \R^I$ such that
$u>0$ and $Au\geq 0$. Furthermore, we have the following equivalences.
\begin{itemize}
\item[{\rm (a)}] $A$ is of finite type if and only if there exists $u\in
\R^I$ such that $u\geq 0$, $Au\geq 0$ and $Au\neq 0$. In this case, 
$\det(A)\neq 0$.
\item[{\rm (b)}] $A$ is of affine type if and only if there exists 
$0\neq u\in\R^I$ such that $u\geq 0$ and $Au=0$. In this case, $A$ has 
rank $|I|-1$. 
\end{itemize}
\end{cor}

\begin{proof} The first statement is clear by the characterisations of
(FIN) and (AFF) in Theorem~\ref{kac2}.

(a) Assume that there exists $u\in\R^I$ such that $u\geq 0$, $Au\geq 0$ 
and $Au\neq 0$; in particular, $u\neq 0$ and so (IND) does not hold. 
Furthermore, $Au\neq 0$ and so (AFF) does not hold. Hence, the only 
remaining possibility is that (FIN) holds. Conversely, assume that
(FIN) holds. Then, by Theorem~\ref{kac2}, there is some $u\in\R^I$ 
such that $u>0$ and $Au\geq 0$. If we had $Au=0$, then also 
$A(-u)=0$ and so $-u\in \cK_A$, contradiction to $\cK_A\subseteq 
\{u\in \R^I\mid u>0\}\cup \{0\}$. Finally assume, if possible, that
$\det(A)=0$. Then there exists some $0\neq v\in\R^I$ such that $Av=0$. 
But then $v,-v\in \cK_A\subseteq \{u \in \R^I\mid u>0\}\cup\{0\}$ 
and so $v=0$, contradiction.

(b) If (AFF) holds, then Theorem~\ref{kac2} shows that there is 
some $u\in\R^I$ such that $u>0$ and $Au=0$, as required. Conversely, assume 
that there exists $0\neq u\in\R^I$ such that $u\geq 0$ and $Au=0$; in 
particular, $u\in\cK_A$ and $\det(A)=0$. But then neither (FIN) nor (IND) 
holds, so (AFF) must hold. The statement about the rank of $A$ is clear by 
condition (AFF).
\end{proof}

\begin{rem} \label{kaccart} Let $A=(a_{ij})_{i,j \in I}$ be the structure
matrix of a Lie algebra $\cL$ of Cartan--Killing type, as in
Chapter~\ref{chap2}. Assume that $A$ is indecomposable and, hence,
$\cL$ is simple. As already remarked above, $A$ satisfies (C1) and (C2). 
So we can now ask whether $A$ is of finite, affine or indefinite type. 
We claim that $A$ is of finite type. To see this, let $\alpha\in\Phi^+$ 
be such that $\hgt(\alpha)$ is as large as possible. Write $\alpha=
\sum_{j\in I} n_j \alpha_j$ where $n_j\in \Z_{\geq 0}$. Let $i\in I$. 
Using the formula in Remark~\ref{explicit}, we obtain 
\[\alpha-\Bigl(\sum_{j\in I} a_{ij}n_j\Bigr)\alpha_i=\sum_{j\in I} 
n_j\bigl(\alpha_j-a_{ij}\alpha_i\bigr)=s_i(\alpha) \in\Phi.\]
Now $\hgt(s_i(\alpha))\leq \hgt(\alpha)$ and so $\sum_{j\in I} a_{ij}
n_j\geq 0$ for all $i\in I$. Hence, we have $A u\geq 0$ where $u:=
(n_i)_{i\in I}\geq 0$. Furthermore, $\det(A)\neq 0$ and so $Au\neq 0$. 
So $A$ is of finite type by Corollary~\ref{kac3}(a).
\end{rem}

\begin{defn}[Kac \protect{\cite[\S 1.1]{K}}, Moody--Pianzola
\protect{\cite[\S 3.4]{MP}}] \label{kac5} 
Assume that $A=(a_{ij})_{i,j\in I}$ satisfies {\rm (C1)}, {\rm (C2)}. 
We say that $A$ is a \nm{generalized Cartan matrix} if $a_{ij}\in\Z$ 
and $a_{ii}=2$ for all $i,j\in I$. 
\end{defn}

Our aim is to classify the indecomposable generalized Cartan matrices 
of finite and affine type. We begin with some preparations.

\begin{lem} \label{kac4} Assume that $A$ satisfies {\rm (C1)}, {\rm (C2)} 
and is indecomposable. Let $A_J:=(a_{ij})_{i,j\in J}$ where $\varnothing 
\neq J\subsetneqq I$. Then, clearly, $A_J$ also satisfies {\rm (C1)}, 
{\rm (C2)}. If $A$ is of finite or affine type and if $A_J$ is 
indecomposable, then $A_J$ is of finite type.
\end{lem}

\begin{proof} Since $A$ is of finite or affine type, there exists
$u\in\R^I$ such that $u>0$ and $Au\geq 0$. Define $u':=(u_i)_{i\in J} 
\in \R^J$. For $i\in J$ we have
\[ 0\leq (Au)_i=\sum_{j\in I}  a_{ij}u_j=\sum_{j\in J} a_{ij} u_j+
\sum_{j\in I\setminus J} \underbrace{a_{ij}u_j}_{\leq 0}\leq (A_Ju')_i.\]
Hence, $u'>0$ and $u'\in \cK_{A_J}$ which means that $A_J$ is of finite or 
affine type (see Theorem~\ref{kac2}). By Corollary~\ref{kac3}, it remains 
to show that $A_Ju'\neq 0$. Assume, if possible, that $(A_Ju')_i=0$ for
all $i\in J$. Then the above inequality shows that $a_{ij}u_j=0$ and,
hence, $a_{ij}=0$ for all $i\in J$ and $j\in I\setminus J$. But then $A$ 
is decomposable, contradiction.
\end{proof}

\begin{lem} \label{kac6} Let $A:=(a_{ij})_{i,j\in I}$ be an indecomposable 
generalized Cartan matrix of finite or affine type. Then $0\leq a_{ij}
a_{ji}\leq 4$ for all $i,j\in I$. If $|I|\geq 3$, then 
$0\leq a_{ij}a_{ji}\leq 3$ for all $i\neq j$ in $I$.
\end{lem}

\begin{proof} If $i=j$, then $a_{ii}=2$ and so the assertion
is clear. Now let $|I|\geq 2$ and $J=\{i,j\}$, where $i\neq j$ 
in $I$ are such that $a_{ij} \neq 0$. Then 
$\renewcommand{\arraystretch}{0.8} A_J=\left(\begin{array}{rr} 
2 & -a \\ -b & 2\end{array} \right)$ where $a=-a_{ij}$, $b=-a_{ji}$, 
$a,b>0$. If $|I|=2$, then $A_J=A$; otherwise, $A_J$ is of finite type
by Lemma~\ref{kac4}. In any case, there exists some $u\in \R^J$ such 
that $u>0$ and $A_Ju\geq 0$. Now 
\[ \renewcommand{\arraystretch}{0.8} 0\leq A_Ju=\left(\begin{array}{rr} 
2 & -a \\ -b & 2\end{array}\right)\left(\begin{array}{c} u_i \\ u_j
\end{array}\right)= \left(\begin{array}{c} 2u_i-au_j\\ -bu_1+2u_j
\end{array}\right).\]
This shows that $b/2\leq u_j/u_i\leq 2/a$ and so $ab\leq 4$, as desired. 
Finally, if $|I|\geq 3$, then $A_J$ is of finite type (as already noted) 
and so $\det(A_J)\neq 0$ by Corollary~\ref{kac3}(a). This implies that 
$ab\neq 4$, as claimed.
\end{proof}

\addtocounter{table}{3}
\begin{table}[htbp] 
\caption{Dynkin diagrams of finite type} \label{Mdynkintbl} 
\begin{center} {\small 
\begin{picture}(290,155)
\put( 00, 25){$E_7$}
\put( 20, 25){\circle*{5}}
\put( 18, 30){$1$}
\put( 20, 25){\line(1,0){20}}
\put( 40, 25){\circle*{5}}
\put( 38, 30){$3$}
\put( 40, 25){\line(1,0){20}}
\put( 60, 25){\circle*{5}}
\put( 58, 30){$4$}
\put( 60, 25){\line(0,-1){20}}
\put( 60, 05){\circle*{5}}
\put( 65, 03){$2$}
\put( 60, 25){\line(1,0){20}}
\put( 80, 25){\circle*{5}}
\put( 78, 30){$5$}
\put( 80, 25){\line(1,0){20}}
\put(100, 25){\circle*{5}}
\put( 98, 30){$6$}
\put(100, 25){\line(1,0){20}}
\put(120, 25){\circle*{5}}
\put(118, 30){$7$}

\put(150, 25){$E_8$}
\put(170, 25){\circle*{5}}
\put(168, 30){$1$}
\put(170, 25){\line(1,0){20}}
\put(190, 25){\circle*{5}}
\put(188, 30){$3$}
\put(190, 25){\line(1,0){20}}
\put(210, 25){\circle*{5}}
\put(208, 30){$4$}
\put(210, 25){\line(0,-1){20}}
\put(210, 05){\circle*{5}}
\put(215, 03){$2$}
\put(210, 25){\line(1,0){20}}
\put(230, 25){\circle*{5}}
\put(228, 30){$5$}
\put(230, 25){\line(1,0){20}}
\put(250, 25){\circle*{5}}
\put(248, 30){$6$}
\put(250, 25){\line(1,0){20}}
\put(270, 25){\circle*{5}}
\put(268, 30){$7$}
\put(270, 25){\line(1,0){20}}
\put(290, 25){\circle*{5}}
\put(288, 30){$8$}

\put( 00, 59){$G_2$}
\put( 20, 60){\circle*{6}}
\put( 18, 66){$1$}
\put( 20, 58){\line(1,0){20}}
\put( 20, 60){\line(1,0){20}}
\put( 20, 62){\line(1,0){20}}
\put( 26, 57.5){$>$}
\put( 40, 60){\circle*{6}}
\put( 38, 66){$2$}

\put( 80, 60){$F_4$}
\put(100, 60){\circle*{5}}
\put( 98, 65){$1$}
\put(100, 60){\line(1,0){20}}
\put(120, 60){\circle*{5}}
\put(118, 65){$2$}
\put(120, 58){\line(1,0){20}}
\put(120, 62){\line(1,0){20}}
\put(126, 57.5){$>$}
\put(140, 60){\circle*{5}}
\put(138, 65){$3$}
\put(140, 60){\line(1,0){20}}
\put(160, 60){\circle*{5}}
\put(158, 65){$4$}

\put(190, 75){$E_6$}
\put(210, 75){\circle*{5}}
\put(208, 80){$1$}
\put(210, 75){\line(1,0){20}}
\put(230, 75){\circle*{5}}
\put(228, 80){$3$}
\put(230, 75){\line(1,0){20}}
\put(250, 75){\circle*{5}}
\put(248, 80){$4$}
\put(250, 75){\line(0,-1){20}}
\put(250, 55){\circle*{5}}
\put(255, 53){$2$}
\put(250, 75){\line(1,0){20}}
\put(270, 75){\circle*{5}}
\put(268, 80){$5$}
\put(270, 75){\line(1,0){20}}
\put(290, 75){\circle*{5}}
\put(288, 80){$6$}

\put( 00,110){$D_n$}
\put( 00,100){$\scriptstyle{n \geq 3}$}
\put( 20,130){\circle*{5}}
\put( 25,130){$1$}
\put( 20,130){\line(1,-1){21}}
\put( 20, 90){\circle*{5}}
\put( 26, 85){$2$}
\put( 20, 90){\line(1,1){21}}
\put( 40,110){\circle*{5}}
\put( 38,115){$3$}
\put( 40,110){\line(1,0){30}}
\put( 60,110){\circle*{5}}
\put( 58,115){$4$}
\put( 80,110){\circle*{1}}
\put( 90,110){\circle*{1}}
\put(100,110){\circle*{1}}
\put(110,110){\line(1,0){10}}
\put(120,110){\circle*{5}}
\put(117,115){$n$}

\put(170,110){$C_n$}
\put(170,100){$\scriptstyle{n \geq 2}$}
\put(190,110){\circle*{5}}
\put(188,115){$1$}
\put(190,108){\line(1,0){20}}
\put(190,112){\line(1,0){20}}
\put(196,107.5){$>$}
\put(210,110){\circle*{5}}
\put(208,115){$2$}
\put(210,110){\line(1,0){30}}
\put(230,110){\circle*{5}}
\put(228,115){$3$}
\put(250,110){\circle*{1}}
\put(260,110){\circle*{1}}
\put(270,110){\circle*{1}}
\put(280,110){\line(1,0){10}}
\put(290,110){\circle*{5}}
\put(287,115){$n$}

\put( 00,150){$A_n$}
\put( 00,140){$\scriptstyle{n \geq 1}$}
\put( 20,150){\circle*{5}}
\put( 18,155){$1$}
\put( 20,150){\line(1,0){20}}
\put( 40,150){\circle*{5}}
\put( 38,155){$2$}
\put( 40,150){\line(1,0){30}}
\put( 60,150){\circle*{5}}
\put( 58,155){$3$}
\put( 80,150){\circle*{1}}
\put( 90,150){\circle*{1}}
\put(100,150){\circle*{1}}
\put(110,150){\line(1,0){10}}
\put(120,150){\circle*{5}}
\put(117,155){$n$}

\put(170,150){$B_n$}
\put(170,140){$\scriptstyle{n \geq 2}$}
\put(190,150){\circle*{5}}
\put(188,155){$1$}
\put(190,148){\line(1,0){20}}
\put(190,152){\line(1,0){20}}
\put(196,147.5){$<$}
\put(210,150){\circle*{5}}
\put(208,155){$2$}
\put(210,150){\line(1,0){30}}
\put(230,150){\circle*{5}}
\put(228,155){$3$}
\put(250,150){\circle*{1}}
\put(260,150){\circle*{1}}
\put(270,150){\circle*{1}}
\put(280,150){\line(1,0){10}}
\put(290,150){\circle*{5}}
\put(287,155){$n$}
\end{picture}}
{\footnotesize (The numbers attached to the vertices define a standard 
labelling of the graph.)}
\end{center}
\end{table}

\begin{table}[htbp] \caption{Dynkin diagrams of affine type} 
\label{Mdynkintbl1} 
\begin{center} {\small 
\begin{picture}(290,210)
\put(  0, 45){$\tilde{E}_6$}
\put( 20, 45){\circle*{5}}
\put( 18, 50){$1$}
\put( 20, 45){\line(1,0){20}}
\put( 40, 45){\circle*{5}}
\put( 38, 50){$2$}
\put( 40, 45){\line(1,0){20}}
\put( 60, 45){\circle*{5}}
\put( 58, 50){$3$}
\put( 60, 45){\line(0,-1){20}}
\put( 60, 25){\circle*{5}}
\put( 65, 23){$2$}
\put( 60, 25){\line(0,-1){17.5}}
\put( 60,  5){\circle{5}}
\put( 65,  3){$1$}
\put( 60, 45){\line(1,0){20}}
\put( 80, 45){\circle*{5}}
\put( 78, 50){$2$}
\put( 80, 45){\line(1,0){20}}
\put(100, 45){\circle*{5}}
\put( 98, 50){$1$}

\put(130, 25){$\tilde{E}_8$}
\put(150, 25){\circle*{5}}
\put(148, 30){$2$}
\put(150, 25){\line(1,0){20}}
\put(170, 25){\circle*{5}}
\put(168, 30){$4$}
\put(170, 25){\line(1,0){20}}
\put(190, 25){\circle*{5}}
\put(188, 30){$6$}
\put(190, 25){\line(0,-1){20}}
\put(190, 05){\circle*{5}}
\put(195, 03){$3$}
\put(190, 25){\line(1,0){20}}
\put(210, 25){\circle*{5}}
\put(208, 30){$5$}
\put(210, 25){\line(1,0){20}}
\put(230, 25){\circle*{5}}
\put(228, 30){$4$}
\put(230, 25){\line(1,0){20}}
\put(250, 25){\circle*{5}}
\put(248, 30){$3$}
\put(250, 25){\line(1,0){20}}
\put(270, 25){\circle*{5}}
\put(268, 30){$2$}
\put(270, 25){\line(1,0){17.5}}
\put(290, 25){\circle{5}}
\put(288, 30){$1$}

\put(  0, 74){$\tilde{F}_4$}
\put( 20, 75){\circle{5}}
\put( 18, 81){$1$}
\put( 22.5, 75){\line(1,0){18}}
\put( 40, 75){\circle*{5}}
\put( 38, 81){$2$}
\put( 40, 75){\line(1,0){20}}
\put( 60, 75){\circle*{5}}
\put( 58, 81){$3$}
\put( 60, 73){\line(1,0){20}}
\put( 60, 77){\line(1,0){20}}
\put( 66, 72.5){$>$}
\put( 80, 75){\circle*{5}}
\put( 78, 81){$4$}
\put( 80, 75){\line(1,0){20}}
\put(100, 75){\circle*{5}}
\put( 98, 81){$2$}

\put(150, 75){$\tilde{E}_7$}
\put(170, 75){\circle{5}}
\put(168, 80){$1$}
\put(172.5, 75){\line(1,0){18}}
\put(190, 75){\circle*{5}}
\put(188, 80){$2$}
\put(190, 75){\line(1,0){20}}
\put(210, 75){\circle*{5}}
\put(208, 80){$3$}
\put(230, 75){\line(0,-1){20}}
\put(230, 55){\circle*{5}}
\put(235, 53){$2$}
\put(210, 75){\line(1,0){20}}
\put(230, 75){\circle*{5}}
\put(228, 80){$4$}
\put(230, 75){\line(1,0){20}}
\put(250, 75){\circle*{5}}
\put(248, 80){$3$}
\put(250, 75){\line(1,0){20}}
\put(270, 75){\circle*{5}}
\put(268, 80){$2$}
\put(270, 75){\line(1,0){20}}
\put(290, 75){\circle*{5}}
\put(288, 80){$1$}

\put(225,105){$\tilde{G}_2$}
\put(250,105){\circle{5}}
\put(248,110){$1$}
\put(252.5,105){\line(1,0){18}}
\put(270,105){\circle*{5}}
\put(268,110){$2$}
\put(270,103){\line(1,0){20}}
\put(270,105){\line(1,0){20}}
\put(270,107){\line(1,0){20}}
\put(276,102.5){$>$}
\put(290,105){\circle*{5}}
\put(288,110){$3$}

\put(  0,120){$\tilde{B}_n$}
\put(  0,110){$\scriptstyle{n \geq 3}$}
\put( 20,120){\circle*{5}}
\put( 17,125){$2$}
\put( 26,117.5){$<$}
\put( 20,122){\line(1,0){20}}
\put( 20,118){\line(1,0){20}}
\put( 40,120){\circle*{5}}
\put( 38,125){$2$}
\put( 40,120){\line(1,0){30}}
\put( 60,120){\circle*{5}}
\put( 58,125){$2$}
\put( 80,120){\circle*{1}}
\put( 90,120){\circle*{1}}
\put(100,120){\circle*{1}}
\put(110,120){\line(1,0){10}}
\put(120,120){\circle*{5}}
\put(117,125){$2$}
\put(120,120){\line(1,1){13}}
\put(120,120){\line(1,-1){14}}
\put(135,135){\circle{5}}
\put(135,105){\circle*{5}}
\put(138, 98){$1$}
\put(137,139){$1$}

\put(145,170){$\tilde{D}_n$}
\put(145,160){$\scriptstyle{n \geq 4}$}
\put(175,170){\circle*{5}}
\put(166,172){$1$}
\put(175,170){\line(1,-1){14}}
\put(175,140){\circle*{5}}
\put(166,132){$1$}
\put(175,140){\line(1,1){14}}
\put(190,155){\circle*{5}}
\put(188,160){$2$}
\put(190,155){\line(1,0){30}}
\put(210,155){\circle*{5}}
\put(208,160){$2$}
\put(230,155){\circle*{1}}
\put(240,155){\circle*{1}}
\put(250,155){\circle*{1}}
\put(260,155){\line(1,0){10}}
\put(270,155){\circle*{5}}
\put(267,160){$2$}
\put(270,155){\line(1,1){13}}
\put(270,155){\line(1,-1){14}}
\put(285,170){\circle{5}}
\put(285,140){\circle*{5}}
\put(288,132){$1$}
\put(287,174){$1$}

\put(150,205){$\tilde{C}_n$}
\put(150,195){$\scriptstyle{n \geq 2}$}
\put(170,205){\circle*{5}}
\put(168,210){$1$}
\put(170,203){\line(1,0){20}}
\put(170,207){\line(1,0){20}}
\put(176,202.5){$>$}
\put(190,205){\circle*{5}}
\put(188,210){$2$}
\put(190,205){\line(1,0){30}}
\put(210,205){\circle*{5}}
\put(208,210){$2$}
\put(230,205){\circle*{1}}
\put(240,205){\circle*{1}}
\put(250,205){\circle*{1}}
\put(260,205){\line(1,0){10}}
\put(270,205){\circle*{5}}
\put(267,210){$2$}
\put(270,203){\line(1,0){18.5}}
\put(270,207){\line(1,0){18.5}}
\put(276,202.5){$<$}
\put(290,205){\circle{5}}
\put(287,210){$1$}

\put(  0,205){$\tilde{A}_n$}
\put(  0,195){$\scriptstyle{n \geq 2}$}
\put( 20,205){\circle*{5}}
\put( 18,210){$1$}
\put( 20,205){\line(1,0){20}}
\put( 40,205){\circle*{5}}
\put( 38,210){$1$}
\put( 40,205){\line(1,0){30}}
\put( 60,205){\circle*{5}}
\put( 58,210){$1$}
\put( 78,205){\circle*{1}}
\put( 88,205){\circle*{1}}
\put( 98,205){\circle*{1}}
\put(105,205){\line(1,0){15}}
\put(120,205){\circle*{5}}
\put(117,210){$1$}
\put( 20,205){\line(2,-1){48}}
\put( 70,179){\circle{5}}
\put(120,205){\line(-2,-1){48}}
\put( 75,174){$1$}

\put( 0,155){$\tilde{A}_1$}
\put(20,155){\circle*{5}}
\put(17,160){$1$}
\put(20,157){\line(1,0){18.5}}
\put(20,153){\line(1,0){18.5}}
\put(40,155){\circle{5}}
\put(38,160){$1$}
\end{picture}}
\end{center}
\begin{center} {\small 
\begin{picture}(290,110)
\put( 0, 92){$A_{2n}^{(2)}$}
\put( 1, 82){$\scriptstyle{n \geq 2}$}
\put(30, 90){\circle{5}}
\put(28, 95){$2$}
\put(31.5, 88){\line(1,0){18}}
\put(31.5, 92){\line(1,0){18}}
\put(36, 87.5){$<$}
\put(50, 90){\circle*{5}}
\put(48, 95){$2$}
\put(50, 90){\line(1,0){30}}
\put(70, 90){\circle*{5}}
\put(68, 95){$2$}
\put(90, 90){\circle*{1}}
\put(100,90){\circle*{1}}
\put(110,90){\circle*{1}}
\put(120,90){\line(1,0){10}}
\put(130,90){\circle*{5}}
\put(128,95){$2$}
\put(130,88){\line(1,0){19}}
\put(130,92){\line(1,0){19}}
\put(136,87.5){$<$}
\put(150,90){\circle*{5}}
\put(147,95){$1$}

\put( 0, 12){$D_{n{+}1}^{(2)}$}
\put( 3,  2){$\scriptstyle{n \geq 2}$}
\put(30, 10){\circle{5}}
\put(28, 15){$1$}
\put(31.5,  8){\line(1,0){18}}
\put(31.5, 12){\line(1,0){18}}
\put(36,  7.5){$<$}
\put(50, 10){\circle*{5}}
\put(48, 15){$1$}
\put(50, 10){\line(1,0){30}}
\put(70, 10){\circle*{5}}
\put(68, 15){$1$}
\put(90, 10){\circle*{1}}
\put(100,10){\circle*{1}}
\put(110,10){\circle*{1}}
\put(120,10){\line(1,0){10}}
\put(130,10){\circle*{5}}
\put(128,15){$1$}
\put(130, 8){\line(1,0){19}}
\put(130,12){\line(1,0){19}}
\put(136, 7.5){$>$}
\put(150,10){\circle*{5}}
\put(147,15){$1$}

\put(180, 12){$E_6^{(2)}$}
\put(210, 10){\circle{5}}
\put(212.5, 10){\line(1,0){18}}
\put(208, 15){$1$}
\put(230, 10){\circle*{5}}
\put(230, 10){\line(1,0){20}}
\put(228, 15){$2$}
\put(250, 10){\circle*{5}}
\put(250, 12){\line(1,0){20}}
\put(250,  8){\line(1,0){20}}
\put(256, 7.5){$<$}
\put(248, 15){$3$}
\put(270, 10){\circle*{5}}
\put(270, 10){\line(1,0){20}}
\put(268, 15){$2$}
\put(290, 10){\circle*{5}}
\put(288, 15){$1$}

\put(220, 52){$D_4^{(3)}$}
\put(250, 50){\circle{5}}
\put(252.5, 50){\line(1,0){18}}
\put(248, 56){$1$}
\put(270, 50){\circle*{5}}
\put(268, 55){$2$}
\put(270, 52){\line(1,0){20}}
\put(270, 50){\line(1,0){20}}
\put(270, 48){\line(1,0){20}}
\put(276, 47.5){$<$}
\put(288, 55){$1$}
\put(290, 50){\circle*{5}}

\put(240, 92){$A_2^{(2)}$}
\put(270, 90){\circle{6}}
\put(268, 95){$2$}
\put(270, 93){\line(1,0){20}}
\put(272.5, 91){\line(1,0){18}}
\put(272.5, 89){\line(1,0){18}}
\put(270, 87){\line(1,0){20}}
\put(276, 87.5){$<$}
\put(288, 95){$1$}
\put(290, 90){\circle*{6}}

\put( 0,52){$A_{2n{-}1}^{(2)}$}
\put( 5,42){$\scriptstyle{n \geq 3}$}
\put(35,65){\circle{5}}
\put(25,65){$1$}
\put(37,63){\line(1,-1){14}}
\put(35,35){\circle*{5}}
\put(25,31){$1$}
\put(35,35){\line(1,1){14}}
\put(50,50){\circle*{5}}
\put(48,55){$2$}
\put(50,50){\line(1,0){30}}
\put(70,50){\circle*{5}}
\put(68,55){$2$}
\put(90,50){\circle*{1}}
\put(100,50){\circle*{1}}
\put(110,50){\circle*{1}}
\put(120,50){\line(1,0){10}}
\put(130,50){\circle*{5}}
\put(128,55){$2$}
\put(130,48){\line(1,0){20}}
\put(130,52){\line(1,0){20}}
\put(136,47.5){$<$}
\put(150,50){\circle*{5}}
\put(147,55){$1$}
\end{picture}}
\end{center}
{\footnotesize (Each diagram denoted $\tilde{X}_n$ has $n+1$ vertices;
$A_{2n}^{(2)}$, $A_{2n-1}^{(2)}$, $D_{n+1}^{(2)}$ have $n+1$ vertices;
the numbers attached to the vertices define a vector 
$u=(u_i)_{i\in I}$ such that $Au=0$.)} 
\end{table}

\begin{defn} \label{kac7} Let $A=(a_{ij})_{i,j\in I}$ be an indecomposable 
generalized Cartan matrix of finite or affine type. Then we encode $A$ in a 
diagram, called \nm{Dynkin diagram} and denoted by $\Gamma(A)$, as follows.

The vertices of $\Gamma(A)$ are in bijection to $I$. Now let $i,j\in I$,
$i\neq j$. If $a_{ij}=a_{ji}=0$, then there is no edge between the vertices
labelled by~$i$ and~$j$. Now assume that $a_{ij}\neq 0$. By
Lemma~\ref{kac6}, we have $1\leq a_{ij}a_{ji}\leq 4$. If $a_{ij}=a_{ji}=-2$,
then the vertices labelled by $i$, $j$ will be joined by a double edge.
Otherwise, we can choose the notation such that $a_{ij}=-1$; let $m:=
-a_{ji}\in\{1,2,3,4\}$. Then the vertices labelled by $i$, $j$ will be 
joined by $m$ edges; if $m\geq 2$, then we put an additional arrow 
pointing towards~$j$. 

Note that $A$ and $\Gamma(A)$ determine each other completely; the fact 
that $A$ is indecomposable means that $\Gamma(A)$ is connected. Examples:

\noindent If $\renewcommand{\arraystretch}{0.8} A=\left(\begin{array}{rr} 
2 & -2 \\ -2 & 2 \end{array} \right)$, then $\Gamma(A)$ is the graph 
$\tilde{A}_1$ in Table~\ref{Mdynkintbl1}. 

\noindent If $\renewcommand{\arraystretch}{0.8} A=\left(\begin{array}{rr} 
2 & -4 \\ -1 & 2 \end{array} \right)$, then $\Gamma(A)$ is the graph 
$\tilde{A}_2^{(2)}$ in Table~\ref{Mdynkintbl1}. 

\noindent If $A$ corresponds to the Lie algebra $\slm_n(\C)$ ($n\geq 2$), 
then $\Gamma(A)$ is the graph $A_{n-1}$ in Table~\ref{Mdynkintbl}; see 
Example~\ref{cartsln}. If $A$ corresponds to a classical 
Lie algebra $\gom_n(Q_n,\C)$, then Table~\ref{cartanBCD} 
(p.~\pageref{cartanBCD}) shows that
\[\mbox{$\Gamma(A)$ is the graph } \left\{\begin{array}{cl} 
B_m &\quad \mbox{if $Q_n^{\text{tr}}=Q_n$ and $n=2m+1\geq 5$},\\ 
C_m &\quad \mbox{if $Q_n^{\text{tr}}=-Q_n$ and $n=2m\geq 4$},\\
D_m &\quad \mbox{if $Q_n^{\text{tr}}=Q_n$ and $n=2m\geq 6$}.
\end{array}\right.\]
\end{defn}

\begin{lem} \label{kac8} The graphs in Table~\ref{Mdynkintbl} correspond to 
indecomposable generalized Cartan matrices of finite type; those in
Table~\ref{Mdynkintbl1} to indecomposable generalized 
Cartan matrices of affine type.
\end{lem}

\begin{proof} Let $\Gamma$ be one of the diagrams in 
Table~\ref{Mdynkintbl1}. Let $|I|=n+1$ and write $I=\{0,1,\ldots,n\}$ 
where $1,\ldots,n$ correspond to the vertices ``$\bullet$'' and $0$ 
corresponds to the vertex~``$\circ$''. Using the scheme in 
Definition~\ref{kac7}, we obtain an indecomposable generalized Cartan 
matrix~$A$ such that $\Gamma=\Gamma(A)$. Let $u=(u_i)_{i\in I}$ be the 
vector defined by the numbers attached to the vertices in 
Table~\ref{Mdynkintbl1}. One checks in each case that $u>0$, $Au=0$ and 
so $A$ is of affine type by Corollary~\ref{kac3}(b). For example, the graph
$D_4^{(3)}$ leads to:
\[\renewcommand{\arraystretch}{0.8}
A=\left(\begin{array}{rrr} 2 & -1 & 0 \\ -1 & 2 & -3 \\ 0 & -1 & 2 
\end{array}\right), \qquad u=\left(\begin{array}{c} 1 \\ 2 \\ 1 
\end{array}\right), \qquad Au=\left(\begin{array}{c} 0 \\ 0 \\ 0 
\end{array}\right).\] 
Finally, all graphs in Table~\ref{Mdynkintbl} are obtained as
$\Gamma(A_J)$ where $J=I\setminus\{0\}$. Now Lemma~\ref{kac4} shows, 
without any further calculations, that $A_J$ is of finite type. 
\end{proof}

\begin{lem} \label{kac9} Let $A=(a_{ij})_{i,j\in I}$ and $A'=(a_{ij}')_{i,
j\in I}$ be indecomposable generalized Cartan matrices such that $A\neq A'$
and $a_{ij} \leq a_{ij}'$ for all $i,j\in I$. If $A$ is of finite or
affine type, then $A'$ is of finite type.
\end{lem}

\begin{proof} Let $A$ be of finite or affine type. There exists some 
$u\in \R^I$ such that $u>0$ and $Au\geq 0$. Let $i\in I$. Then 
\[(A'u)_i=\sum_{j\in I} a_{ij}'u_j \geq \sum_{j\in I} a_{ij}u_j=
(Au)_i \geq 0.\]
So $A'u\geq 0$ and $A'$ is of finite or affine type, by 
Corollary~\ref{kac3}. Since $A\neq A'$, there exist $i,j\in I$ 
such that $a_{ij}<a_{ij}'$. Then the above computation shows that 
$(A'u)_i>(Au)_i\geq 0$. Hence, $A'u\neq 0$ and so $A'$ is of 
finite type (again, by Corollary~\ref{kac3}).
\end{proof}

\begin{lem} \label{kac9a} Let $A=(a_{ij})_{i,j\in I}$ be an 
indecomposable generalized Cartan matrix of finite or affine type. 
Assume that there is a cycle in $\Gamma(A)$, that is, there exist 
indices $i_1, i_2,\ldots,i_r$ in $I$ ($r\geq 3$) such that 
\begin{equation*}
a_{i_1i_2}a_{i_2i_3}\cdots a_{i_{r-1}i_r}a_{i_ri_1}\neq 0 \quad\mbox{and} 
\quad i_1,i_2,\ldots,i_r \mbox{ are distinct}.\tag{$\circlearrowleft$}
\end{equation*}
Then $A$ is of affine type, $|I|=r$ and $\Gamma(A)=\tilde{A}_{r-1}$ in 
Table~\ref{Mdynkintbl1}.
\end{lem}

\begin{proof} Let $J:=\{i_1,\ldots,i_r\}$. By ($\circlearrowleft$) and 
Remark~\ref{indecgraph}, the submatrix $A_J$ is indecomposable. By 
Lemma~\ref{kac4}, $A_J$ is of finite or affine type. Now define $A_J'
=(a_{ij}')_{i,j\in J}$ by 
\begin{align*}
a_{i_1i_2}'&=a_{i_2i_3}'=\ldots =a_{i_{r-1}i_r}'=
a_{i_ri_1}'=-1, \qquad a_{jj}'=2,\\
a_{i_2i_1}'&=a_{i_3i_2}'=\ldots =a_{i_ri_{r-1}}'=
a_{i_1i_r}'=-1, \qquad a_{jj}'=2,
\end{align*}
and $a_{jj'}'=0$ for all other indices $j\neq j'$ in~$J$. Then 
$\Gamma(A_J')$ is the graph $\tilde{A}_{r-1}$ and so $A_J'$ is of 
affine type; see Lemma~\ref{kac8}. We claim that $A_J=A_J'$. Indeed,
by ($\circlearrowleft$), we have $a_{ij}\leq a_{ij}'$ for all $i,j\in J$. 
Hence, if we had $A_J\neq A_J'$, then Lemma~\ref{kac9} would imply that 
$A_J'$ is of finite type, contradiction. So we must have $A_J=A_J'$. 
Next we claim that $J=I$. Indeed, if we had $J\subsetneqq I$, then 
Lemma~\ref{kac4} would imply that $A_J=A_J'$ is of finite type,
contradiction. So we must have $J=I$ and, hence, $A=A_J=A_J'$ is
of affine type.
\end{proof}

\begin{thm} \label{classcart} The Dynkin diagrams of indecomposable
generalized Cartan matrices of finite type are precisely those in
Table~\ref{Mdynkintbl}.
\end{thm}

\begin{proof} By Lemma~\ref{kac8}, we already know that all diagrams
in Table~\ref{Mdynkintbl} satisfy this condition. Conversely, let 
$A=(a_{ij})_{i,j\in I}$ be an arbitrary indecomposable generalized 
Cartan matrix of finite type. We must show that the corresponding diagram 
$\Gamma(A)$ appears in Table~\ref{Mdynkintbl}. If $|I|=1$, then $A=(2)$ 
and $\Gamma(A)=A_1$. Now let $|I|=2$. Then 
\[ \renewcommand{\arraystretch}{0.8} A=\left(\begin{array}{rr} 2 & -a \\ 
-b & 2 \end{array}\right) \quad \mbox{where}\quad a,b\in \Z_{>0}
\;\mbox{and}\; 1\leq ab\leq 4;\]
see Lemma~\ref{kac6}. Since $\det(A)\neq 0$, we also have $ab\neq 4$.
So $\Gamma(A)$ is one of the graphs $A_2$, $B_2$, $C_2$ or $G_2$. Assume
from now on that $|I|\geq 3$. By Lemma~\ref{kac6}, there are only single, 
double or triple edges in $\Gamma(A)$ (and an arrow is attached to a 
double or triple edge). Using Lemmas~\ref{kac4} and~\ref{kac9}, 
one obtains further restrictions on $\Gamma(A)$ which eventually lead to
the list of graphs in Table~\ref{Mdynkintbl}. We give full details for 
one example.

\underline{Claim:} $\Gamma(A)$ does not have a triple edge. This is seen as
follows. Assume, if possible, that there are $i\neq j$ in $I$ which are 
connected by a triple edge. Since $|I|\geq 3$ and $A$ is indecomposable, 
there is a further $k\in I$ connected to $i$ or $j$; we choose the notation 
such that $k$ is connected to $i$. By Lemma~\ref{kac9a}, there are no cyles
in $\Gamma(A)$ and so there is no edge between $j,k$. Let $J:=\{k,i,j\}$ 
and consider the submatrix $A_J$. We have
\[\renewcommand{\arraystretch}{0.8} 
A_J=\left(\begin{array}{rrr} 2 & -a & 0 \\ -b & 2 & -c \\ 0 & -d & 
2\end{array}\right)\qquad \mbox{where $a,b,c,d\in \Z_{>0}$ and $cd=3$}.\]
Then $A_J$ must also be of finite type; see Lemma~\ref{kac4}. Let 
\[ \renewcommand{\arraystretch}{0.8}
A_J'=\left(\begin{array}{rrr} 2 & -1 & 0 \\ -1 & 2 & -c \\ 0 & -d & 
2\end{array}\right).\]
Then $A_J'$ is still of finite type by Lemma~\ref{kac9}. But 
$\Gamma(A_J')$ is the graph $\tilde{G}_2$ or the graph $D_4^{(3)}$, 
contradiction to Lemma~\ref{kac8}. 

By similar arguments one shows that, if $\Gamma(A)$ has a double edge,
then there is only one double edge and no branch point (that is,
a vertex connected to at least three other vertices). Furthermore, 
if $|I|\geq 5$, then $\Gamma(A)$ must be one of the graphs $B_n$ or
$C_n$. (For otherwise, we would have a subset $J\subseteq I$ such 
that $|J|=5$ and $\Gamma(A_J)$ is the graph $\tilde{F}_4$ or 
$\tilde{E}_6^{(2)}$.) If $|I|=4$, then $\Gamma(A)$ must be the graph 
$F_4$. Finally, if $\Gamma(A)$ has only single edges, then one shows that 
there is at most one branch point, and that the remaining possibilities are 
$A_n$, $D_n$, $E_6$, $E_7$ and $E_8$.
\end{proof}

\begin{rem} \label{kac10} By similar arguments, one can also show that the 
Dynkin diagrams of indecomposable generalized Cartan matrices of affine type 
are precisely those in Table~\ref{Mdynkintbl1}; see 
Kac \cite[Chap.~4]{K}.
\end{rem}

\begin{xca} \label{posinv} Let $A$ be an indecomposable generalized 
Cartan matrix of type (FIN). Then $\det(A)\neq 0$ and we can form 
$A^{-1}$. Use condition (FIN) to show that all entries of $A^{-1}$ are 
strictly positive rational numbers. Work out some examples 
explicitly. Explicit formulae for the entries of $A^{-1}$ are found
in Lusztig--Tits \cite{LT}.
\end{xca}

\begin{rem} \label{kac3a} By Vinberg \cite[p.~1099]{Vi}, the type of 
$A$ can also be characterised in terms of the eigenvalues of $A$, as 
follows. Choose any $c\in \R$ such that all diagonal entries (and, hence,
all entries) of $B:=c\,\id_I-A$ are $\geq 0$. Then, by a weak form of 
the \nm{Frobenius--Perron Theorem} (see, e.g., \cite[\S 8.2]{dserre}), 
$B$ has at least one real eigenvalue; furthermore, if $\mu_0$ is the 
largest real eigenvalue, then $\mu_0\geq 0$ and there exists a corresponding 
eigenvector $v\in \R^I$ such that $v\geq 0$. It follows that $A=c\,\id_I 
-B$ also has at least one real eigenvalue. Let $\lambda_0\in\R$ be 
the smallest real eigenvalue of~$A$. Then $\lambda_0=c-\mu_0$ and we still 
have $Av=(c-\mu_0)v= \lambda_0 v$. Then we have:
\[ \mbox{(FIN)} \;\Leftrightarrow \;\lambda_0>0, \qquad
\mbox{(AFF)}\; \Leftrightarrow\; \lambda_0=0, \qquad
\mbox{(IND)}\; \Leftrightarrow \;\lambda_0<0.\]
This is seen as follows. Assume that $\lambda_0\geq 0$. Then 
$v\geq 0$, $Av=\lambda_0v\geq 0$ and so $A$ is of finite type or
affine type by Corollary~\ref{kac3}. Furthermore, if $\lambda_0=0$,
then $Av=0$ and so $A$ is of affine type; if $\lambda_0>0$, then $Av\neq 0$
and so $A$ is of finite type. Conversely, assume that $A$ is of finite 
or affine type. There exists some $u\in \R^I$ such that $u>0$ and 
$Au\geq 0$. For $\lambda\in \R_{>0}$, we have $(A+\lambda\id_I)u=Au+
\lambda u\geq 0$ and $(A+\lambda\id_I)u\neq 0$. Hence, by 
Corollary~\ref{kac3}(a), $A+\lambda \id_I$ is of finite type and 
$\det(A+\lambda\id_I)\neq 0$. Thus, all real eigenvalues of $A$ are 
non-negative and so $\lambda_0\geq 0$. If $A$ is of finite type, then 
$\det(A)\neq 0$ and so $\lambda_0>0$; if $A$ is of affine type, then 
$\det(A)=0$ and so $\lambda_0=0$. Thus, the first two equivalences
are proved; but then the third equivalence follows from 
Theorem~\ref{kac2}. (See also Moody--Pianzola \cite[\S 3.6]{MP}.)
\end{rem}

\begin{rem} \label{proctor}
The diagrams of type $A_n$, $D_n$, $E_n$ arise in a number of situations
and can actually be characterised in a very simple way; see, e.g., 
the short note by Proctor \cite{proct}. (According to Lusztig 
\cite[\S 2]{Lu23}, this is originally due to Coxeter.)
\end{rem}

\section{Finite root systems} \label{sec3a2}

Consider a generalized Cartan matrix $A=(a_{ij})_{i,j\in I}$, where $I$ 
is a non-empty finite index set. Eventually, we would like to construct 
a Lie algebra of Cartan--Killing type with structure matrix~$A$ (at least
for certain~$A$). As a first step, we need to construct the underlying 
root system, directly from~$A$. Let $E$ be an $\R$-vector space with a 
basis $\Delta=\{\alpha_i\mid i \in I\}$. For each $i\in I$, we define
a linear map $s_i\colon E\rightarrow E$ by the formula
\[ s_i(\alpha_j):=\alpha_j-a_{ij}\alpha_i \quad \mbox{for 
$j \in I\quad$ (cf.\ Remark~\ref{explicit})}.\]
Since $a_{ii}=2$, we have $s_i(\alpha_i)=-\alpha_i$. Furthermore, we 
compute $s_i^2(\alpha_j)=s_i(\alpha_j-a_{ij}\alpha_i)=
s_i(\alpha_j)+a_{ij}\alpha_i=\alpha_j$ for all $j\in I$. Hence, 
we have $s_i^2=\id_E$ and so $s_i\in\GL(E)$. The subgroup 
\[ W=W(A):=\langle s_i \mid i \in I\rangle\subseteq \GL(E)\]
is called the \nm{Weyl group} associated with $A$. In analogy to
Theorem~\ref{mainthm1}(a), the corresponding \nm{abstract root system} 
is defined by
\[ \Phi=\Phi(A):=\{w(\alpha_i)\mid w\in W,i \in I\}; \]
the roots $\{\alpha_i\mid i\in I\}$ are also called \nm{simple roots}.
Clearly, if $W$ is finite, then so is $\Phi$. Conversely, assume that
$\Phi$ is finite. By definition, it is clear that $w(\alpha)\in \Phi$ for
all $w\in W$ and $\alpha\in\Phi$. So there is an action of $W$ on $\Phi$.
Since $\Phi$ contains a basis of $E$, we have a corresponding 
\textit{injective} group homomorphism
\[ \pi\colon W=W(A)\hookrightarrow \mbox{Sym}(\Phi).\]
By exactly the same argument as in Remark~\ref{weyl0}, it follows
that $W$ is finite. Hence, we have:
\[ |W(A)|<\infty \qquad \Leftrightarrow \qquad |\Phi(A)|<\infty.\]
In Example~\ref{cartanG2}, we have computed $W(A)$ and $\Phi(A)$
for the matrix $A$ with Dynkin diagram $G_2$ in Table~\ref{Mdynkintbl};
in this case, $|W(A)|=12<\infty$. In Exercise~\ref{xcacartanAA}, there are 
two examples where $|W(A)|=\infty$. (The first of those matrices has
affine type with Dynkin diagram $\tilde{A}_2$ in Table~\ref{Mdynkintbl1};
the second matrix is of indefinite type.)

\begin{rem} \label{decabst} Assume that $A$ is decomposable. So there is 
a partition $I=I_1\sqcup I_2$ such that $A$ has a block diagonal shape
\begin{center}
$A=\left(\begin{array}{c|c} A_1 & 0 \\\hline 0 & A_2\end{array}\right)$
\end{center}
where $A_1$ has rows and columns labelled by $I_1$, and $A_2$ has rows 
and columns labelled by $I_2$. Let 
\[W_1:=\langle s_i \mid i \in I_1 \rangle\subseteq W\qquad\mbox{and}
\qquad W_2:=\langle s_i \mid i \in I_2\rangle \subseteq W.\]
As in Proposition~\ref{levisub2}, one sees that $W=W_1\cdot W_2$, $W_1\cap 
W_2=\{1\}$ and $w_1w_2=w_2w_1$ for all $w_i\in W_i$; furthermore, 
$\Phi=\Phi_1 \sqcup\Phi_2$ where 
\begin{align*}
\Phi_1&:=\{w(\alpha_i)\mid i \in I_1,w\in W_1\} \subseteq \langle
\alpha_i\mid i \in I_1\rangle_\Z,\\ \Phi_2&:=\{w(\alpha_i)\mid i 
\in I_2,w\in W_2\} \subseteq \langle \alpha_i \mid i \in I_2\rangle_\Z.
\end{align*}
Since $W_1 \cong W(A_1)$ and $W_2\cong W(A_2)$, we obtain the equivalence:
\[ |W(A)|<\infty \qquad \Leftrightarrow \qquad |W(A_1)|<\infty
\;\; \mbox{and} \;\; |W(A_2)|<\infty.\]
Thus, in order to characterise those $A$ for which $W(A)$ is finite,
we may assume without loss of generality that $A$ is indecomposable.
\end{rem}

\begin{rem} \label{decabst1a} Assume that $|W(A)|<\infty$. Then we 
can construct a $W(A)$-invariant scalar product $\langle \;,\;\rangle 
\colon E\times E\rightarrow \R$ by the same method as in 
Section~\ref{sec1a3}. (In the sequel, it will not be important how 
exactly $\langle\;,\;\rangle$ is defined; it just needs to be symmetric, 
positive-definite and $W(A)$-invariant.) This yields the formula
\[ a_{ij}=2\frac{\langle\alpha_i,\alpha_j\rangle}{\langle\alpha_i,
\alpha_i\rangle}\qquad\mbox{for all $i,j\in I$};\] 
see the argument in Remark~\ref{keylem0}. Consequently, we have
\[ s_i(v)=v-\langle \alpha_i^\vee,v\rangle \alpha_i \qquad 
\mbox{for all $v\in E$}.\]
Here, we write $\alpha^\vee:=2\alpha/\langle \alpha,\alpha\rangle\in E$ 
for any $\alpha\in \Phi(A)$. As in Remark~\ref{keylem1}, it follows
that $\det(A)>0$.
\end{rem}

\begin{lem} \label{decabst1} Assume that $A$ is indecomposable
and $|W(A)|<\infty$. Then $A$ is of type {\rm (FIN)}.
\end{lem}

\begin{proof} Let $X$ be the set of all $\alpha\in\Phi$ such that 
$\alpha$ can be written as a $\Z$-linear combination of
$\Delta$, where all coefficients are $\geq 0$. Then $X$ is
non-empty; for example, $\Delta\subseteq X$. Let $\alpha_0\in X$
be such that the sum of the coefficients is as large as possible. (This
exists since $|\Phi|<\infty$.) Write $\alpha_0=\sum_{j \in I} n_j\alpha_j$ 
where $n_j\geq 0$ for all $j\in I$. If $m:=\langle \alpha_i^\vee,
\alpha_0\rangle<0$ for some $i\in I$, then
\[ s_i(\alpha_0)=\alpha_0-\langle\alpha_i^\vee,\alpha_0\rangle\alpha_i
=\bigl(\underbrace{n_i-m}_{>n_i}\bigr)\alpha_i+
\sum_{\atop{j \in I}{j\neq i}} n_j \alpha_j \in \Phi,\]
where all coefficients are still non-negative but the sum of the 
coefficients is strictly larger than that of $\alpha_0$, contradiction.
So we must have $\langle \alpha_i^\vee,\alpha_0\rangle\geq 0$ for 
all $i\in I$. But this means $\sum_{j \in I} a_{ij}n_j=\sum_{j \in I} 
n_j\langle \alpha_i^\vee,\alpha_j\rangle\geq 0$. So, if 
$u:=(n_j)_{j\in I} \in \R^I$, then $u\geq 0$, $u\neq 0$, and $Au\geq 0$.
Since $\det(A)\neq 0$, we also have $Au\neq 0$. So $A$ is
of type (FIN) by Corollary~\ref{kac3}(a).
\end{proof}

\begin{table}[htbp] \caption{Positive roots for exceptional types $F_4$,
$E_6$, $E_7$, $E_8$} 
\label{plusF4E678} {\small $\renewcommand{\arraystretch}{1.04}
\begin{array}{l}\begin{array}{c@{\hspace{13pt}}
c@{\hspace{13pt}}c@{\hspace{13pt}}c@{\hspace{13pt}}c@{\hspace{13pt}}
c@{\hspace{13pt}}c@{\hspace{13pt}}c@{\hspace{13pt}}c@{\hspace{13pt}}c}
\hline \multicolumn{3}{l}{\text{Type $F_4$, $\;|\Phi^+|=24$:}} &
   1 0 0 0 &  0 1 0 0 &  0 0 1 0 &  0 0 0 1 & 
   1 1 0 0 &  0 1 1 0 \\  0 0 1 1 &  1 1 1 0 &
   0 1 2 0 &  0 1 1 1 &  1 1 2 0 &  1 1 1 1 &
   0 1 2 1 &  1 2 2 0 &  1 1 2 1 \\  0 1 2 2 & 
   1 2 2 1 &  1 1 2 2 &  1 2 3 1 &  1 2 2 2 & 
   1 2 3 2 &  1 2 4 2 &  1 3 4 2 &  2 3 4 2 
\end{array}\\ \begin{array}{c@{\hspace{8pt}}c@{\hspace{8pt}}
c@{\hspace{8pt}}c@{\hspace{8pt}}c@{\hspace{8pt}}c@{\hspace{8pt}}
c@{\hspace{8pt}}c} \hline \multicolumn{4}{l}{\text{Type $E_6$, 
$\;|\Phi^+|=36$:}} &
   1 0 0 0 0 0 &  0 1 0 0 0 0 &  0 0 1 0 0 0 & 
   0 0 0 1 0 0 \\  0 0 0 0 1 0 &  0 0 0 0 0 1 &
   1 0 1 0 0 0 &  0 1 0 1 0 0 &  0 0 1 1 0 0 & 
   0 0 0 1 1 0 &  0 0 0 0 1 1 &  1 0 1 1 0 0 \\ 
   0 1 1 1 0 0 &  0 1 0 1 1 0 & 0 0 1 1 1 0 & 
   0 0 0 1 1 1 &  1 1 1 1 0 0 &  1 0 1 1 1 0 & 
   0 1 1 1 1 0 &  0 1 0 1 1 1 \\  0 0 1 1 1 1 & 
   1 1 1 1 1 0 &  1 0 1 1 1 1 &  0 1 1 2 1 0 &
   0 1 1 1 1 1 &  1 1 1 2 1 0 &  1 1 1 1 1 1 & 
   0 1 1 2 1 1 \\  1 1 2 2 1 0 &  1 1 1 2 1 1 &
   0 1 1 2 2 1 &  1 1 2 2 1 1 &  1 1 1 2 2 1 & 
   1 1 2 2 2 1 &  1 1 2 3 2 1 &  1 2 2 3 2 1  
\end{array}\\ \begin{array}{c@{\hspace{8pt}}c@{\hspace{8pt}}c@{\hspace{8pt}}
c@{\hspace{8pt}}c@{\hspace{8pt}}c@{\hspace{8pt}}c} \hline 
\multicolumn{7}{l}{\text{Type $E_7$, $\;|\Phi^+|=63$:}}\\ 
   1 0 0 0 0 0 0 &  0 1 0 0 0 0 0 &  0 0 1 0 0 0 0 & 
   0 0 0 1 0 0 0 &  0 0 0 0 1 0 0 &  0 0 0 0 0 1 0 & 
   0 0 0 0 0 0 1 \\  1 0 1 0 0 0 0 &  0 1 0 1 0 0 0 &
   0 0 1 1 0 0 0 &  0 0 0 1 1 0 0 &  0 0 0 0 1 1 0 &
   0 0 0 0 0 1 1 &  1 0 1 1 0 0 0 \\  0 1 1 1 0 0 0 & 
   0 1 0 1 1 0 0 &  0 0 1 1 1 0 0 &  0 0 0 1 1 1 0 &
   0 0 0 0 1 1 1 &  1 1 1 1 0 0 0 &  1 0 1 1 1 0 0 \\ 
   0 1 1 1 1 0 0 &  0 1 0 1 1 1 0 &  0 0 1 1 1 1 0 & 
   0 0 0 1 1 1 1 &  1 1 1 1 1 0 0 &  1 0 1 1 1 1 0 &
   0 1 1 2 1 0 0 \\  0 1 1 1 1 1 0 &  0 1 0 1 1 1 1 & 
   0 0 1 1 1 1 1 &  1 1 1 2 1 0 0 &  1 1 1 1 1 1 0 & 
   1 0 1 1 1 1 1 &  0 1 1 2 1 1 0 \\  0 1 1 1 1 1 1 & 
   1 1 2 2 1 0 0 &  1 1 1 2 1 1 0 &  1 1 1 1 1 1 1 & 
   0 1 1 2 2 1 0 &  0 1 1 2 1 1 1 &  1 1 2 2 1 1 0 \\ 
   1 1 1 2 2 1 0 &  1 1 1 2 1 1 1 &  0 1 1 2 2 1 1 & 
   1 1 2 2 2 1 0 &  1 1 2 2 1 1 1 &  1 1 1 2 2 1 1 &
   0 1 1 2 2 2 1 \\  1 1 2 3 2 1 0 &  1 1 2 2 2 1 1 & 
   1 1 1 2 2 2 1 &  1 2 2 3 2 1 0 &  1 1 2 3 2 1 1 & 
   1 1 2 2 2 2 1 &  1 2 2 3 2 1 1 \\  1 1 2 3 2 2 1 & 
   1 2 2 3 2 2 1 &  1 1 2 3 3 2 1 &  1 2 2 3 3 2 1 & 
   1 2 2 4 3 2 1 &  1 2 3 4 3 2 1 & 2 2 3 4 3 2 1  
\end{array}\\
\begin{array}{cccccc} \hline \multicolumn{6}{l}{\text{Type $E_8$, 
$\;|\Phi^+|=120$:}}\\ 1 0 0 0 0 0 0 0 &  0 1 0 0 0 0 0 0 & 
   0 0 1 0 0 0 0 0 &  0 0 0 1 0 0 0 0 & 
   0 0 0 0 1 0 0 0 &  0 0 0 0 0 1 0 0 \\
   0 0 0 0 0 0 1 0 &  0 0 0 0 0 0 0 1 & 
   1 0 1 0 0 0 0 0 &  0 1 0 1 0 0 0 0 & 
   0 0 1 1 0 0 0 0 &  0 0 0 1 1 0 0 0 \\ 
   0 0 0 0 1 1 0 0 &  0 0 0 0 0 1 1 0 & 
   0 0 0 0 0 0 1 1 &  1 0 1 1 0 0 0 0 & 
   0 1 1 1 0 0 0 0 &  0 1 0 1 1 0 0 0 \\ 
   0 0 1 1 1 0 0 0 &  0 0 0 1 1 1 0 0 & 
   0 0 0 0 1 1 1 0 &  0 0 0 0 0 1 1 1 & 
   1 1 1 1 0 0 0 0 &  1 0 1 1 1 0 0 0 \\ 
   0 1 1 1 1 0 0 0 &  0 1 0 1 1 1 0 0 & 
   0 0 1 1 1 1 0 0 &  0 0 0 1 1 1 1 0 & 
   0 0 0 0 1 1 1 1 &  1 1 1 1 1 0 0 0 \\
   1 0 1 1 1 1 0 0 &  0 1 1 2 1 0 0 0 & 
   0 1 1 1 1 1 0 0 &  0 1 0 1 1 1 1 0 & 
   0 0 1 1 1 1 1 0 &  0 0 0 1 1 1 1 1 \\
   1 1 1 2 1 0 0 0 &  1 1 1 1 1 1 0 0 & 
   1 0 1 1 1 1 1 0 &  0 1 1 2 1 1 0 0 & 
   0 1 1 1 1 1 1 0 &  0 1 0 1 1 1 1 1 \\ 
   0 0 1 1 1 1 1 1 &  1 1 2 2 1 0 0 0 & 
   1 1 1 2 1 1 0 0 &  1 1 1 1 1 1 1 0 & 
   1 0 1 1 1 1 1 1 &  0 1 1 2 2 1 0 0 \\ 
   0 1 1 2 1 1 1 0 &  0 1 1 1 1 1 1 1 & 
   1 1 2 2 1 1 0 0 &  1 1 1 2 2 1 0 0 & 
   1 1 1 2 1 1 1 0 &  1 1 1 1 1 1 1 1 \\ 
   0 1 1 2 2 1 1 0 &  0 1 1 2 1 1 1 1 & 
   1 1 2 2 2 1 0 0 &  1 1 2 2 1 1 1 0 & 
   1 1 1 2 2 1 1 0 &  1 1 1 2 1 1 1 1 \\ 
   0 1 1 2 2 2 1 0 &  0 1 1 2 2 1 1 1 & 
   1 1 2 3 2 1 0 0 &  1 1 2 2 2 1 1 0 & 
   1 1 2 2 1 1 1 1 &  1 1 1 2 2 2 1 0 \\ 
   1 1 1 2 2 1 1 1 &  0 1 1 2 2 2 1 1 & 
   1 2 2 3 2 1 0 0 &  1 1 2 3 2 1 1 0 & 
   1 1 2 2 2 2 1 0 &  1 1 2 2 2 1 1 1 \\ 
   1 1 1 2 2 2 1 1 &  0 1 1 2 2 2 2 1 & 
   1 2 2 3 2 1 1 0 &  1 1 2 3 2 2 1 0 & 
   1 1 2 3 2 1 1 1 &  1 1 2 2 2 2 1 1 \\ 
   1 1 1 2 2 2 2 1 &  1 2 2 3 2 2 1 0 & 
   1 2 2 3 2 1 1 1 &  1 1 2 3 3 2 1 0 & 
   1 1 2 3 2 2 1 1 &  1 1 2 2 2 2 2 1 \\ 
   1 2 2 3 3 2 1 0 &  1 2 2 3 2 2 1 1 & 
   1 1 2 3 3 2 1 1 &  1 1 2 3 2 2 2 1 & 
   1 2 2 4 3 2 1 0 &  1 2 2 3 3 2 1 1 \\ 
   1 2 2 3 2 2 2 1 &  1 1 2 3 3 2 2 1 & 
   1 2 3 4 3 2 1 0 &  1 2 2 4 3 2 1 1 & 
   1 2 2 3 3 2 2 1 &  1 1 2 3 3 3 2 1 \\ 
   2 2 3 4 3 2 1 0 &  1 2 3 4 3 2 1 1 & 
   1 2 2 4 3 2 2 1 &  1 2 2 3 3 3 2 1 & 
   2 2 3 4 3 2 1 1 &  1 2 3 4 3 2 2 1 \\ 
   1 2 2 4 3 3 2 1 &  2 2 3 4 3 2 2 1 & 
   1 2 3 4 3 3 2 1 &  1 2 2 4 4 3 2 1 & 
   2 2 3 4 3 3 2 1 &  1 2 3 4 4 3 2 1 \\ 
   2 2 3 4 4 3 2 1 &  1 2 3 5 4 3 2 1 & 
   2 2 3 5 4 3 2 1 &  1 3 3 5 4 3 2 1 & 
   2 3 3 5 4 3 2 1 &  2 2 4 5 4 3 2 1 \\
   2 3 4 5 4 3 2 1 &  2 3 4 6 4 3 2 1 & 
   2 3 4 6 5 3 2 1 &  2 3 4 6 5 4 2 1 & 
   2 3 4 6 5 4 3 1 &  2 3 4 6 5 4 3 2 \\ \hline
\multicolumn{6}{c}{\text{For example, $2342$ stands for 
$2\alpha_1{+}3 \alpha_2{+}4\alpha_3{+}2\alpha_4$, $I=\{1,2,3,4\}$.}}
\end{array}\end{array}$}
\end{table}

\begin{prop} \label{decabst2} Assume that $A$ is indecomposable and
of type {\rm (FIN)}. Then $|W(A)|<\infty$ and $|\Phi(A)|<\infty$. 
Furthermore, $(\Phi(A),\Delta)$ is a \nm{based root system}, that is,
every $\alpha\in\Phi(A)$ can be written as a $\Z$-linear combination 
of $\Delta=\{\alpha_i\mid i\in I\}$, where the coefficients are
either all $\geq 0$ or all $\leq 0$ (as in condition 
{\rm (CK2)} of Definition~\ref{defTD}). Finally, $\Phi(A)$ is reduced,
that is, $\Phi(A)\cap \R \alpha=\{\pm\alpha\}$ for all
$\alpha\in\Phi(A)$.
\end{prop}

\begin{proof} We use the classification in Section~\ref{sec3a1} and
go through the list of Dynkin diagrams in Table~\ref{Mdynkintbl}. If 
$A$ has a diagram of type $A_n$, $B_n$, $C_n$ or $D_n$, then
$\Phi(A)$ has been explicitly described in Chapter~\ref{chap2}; the 
desired properties hold by Example~\ref{cartsln} and 
Corollary~\ref{sec165}. By inspection, one sees that $\Phi(A)$ is reduced.

Now assume that $A$ has a diagram of type $G_2$, $F_4$, $E_6$, $E_7$, 
or $E_8$. Then we take a ``\nm{computer algebra approach}'', based on 
our \Python\ programs in Table~\ref{pythontab} (p.~\pageref{pythontab}). 
We apply the program {\tt rootsystem} to~$A$; the program actually 
terminates and outputs a finite list\footnote{As shown in Jacobson 
\cite[Chap.~VII, \S 5]{Jac}, it is actually possible to produce such
lists without too much effort ``by hand'', even for type $E_8$.} of 
tuples $\cC^+(A)\subseteq \N_0^I$. For example, for type $G_2$, we obtain:
\[\{(1, 0), (0, 1), (1, 1), (1, 2), (1, 3), (2, 3)\} \quad \mbox{(see 
also Example~\ref{cartanG2})}.\]
For the types $F_4$, $E_6$, $E_7$, $E_8$, these vectors are explicitly 
listed in Table~\ref{plusF4E678}. Now we set $\Phi:=\Phi^+\cup (-\Phi^+)$, 
where 
\[ \Phi^+:=\Bigl\{\alpha:=\sum_{i\in I} n_i\alpha_i\;\big|\; 
(n_i)_{i\in I} \in \cC^+(A)\Bigr\}\subseteq E.\]
By construction, it is clear that $\Phi^+\subseteq \Phi(A)$. Since
$s_i(\alpha_i)=-\alpha_i$ for $i\in I$, it also follows that $
-\Phi^+\subseteq \Phi(A)$. 
Now we apply our program {\tt refl} to all tuples in $\cC^+(A)\cup
(-\cC^+(A))$. By inspection, we find that $\cC^+(A)\cup (-\cC^+(A))$ 
remains invariant under these operations. In other words,
we have $s_i(\Phi)\subseteq \Phi$ for all $i\in I$ (recall
that {\tt refl} corresponds to applying $s_i$ to an element of~$E$).
Since $\Delta\subseteq \Phi$, we conclude that $\Phi(A)\subseteq \Phi$ 
and, hence, that $\Phi(A)=\Phi$; in particular, $|\Phi(A)|<\infty$. The 
fact that $(\Phi(A),\Delta)$ is a based root system is clear because all 
tuples in $\cC^+(A)$ have non-negative entries. The fact that $\Phi(A)$ 
is reduced is seen by inspection of Table~\ref{plusF4E678}: one just
has to check that the coefficients $(n_i)_{i \in I}$ are always coprime.
\end{proof}

\begin{rem} \label{decabst2a} Of course, one can avoid the classification
and the use of computer algebra methods in order to obtain the above
result. The finiteness of $W(A)$ follows from a topological argument, 
based on the fact that $W(A)$ is a discrete, bounded subset of~$\GL(E)$; 
see, e.g., \cite[Ch.~V, \S 4. no.~8]{B}. The fact that $(\Phi(A),
\Delta)$ is based requires a more elaborate argument; see, 
e.g., \cite[(64.28)]{CR2} or \cite[1.1.10]{GePf}.
\end{rem}

\begin{rem} \label{myth} As Lusztig \cite[\S 22]{shaw} writes, $E_8$ has 
an almost mythical status in mathematics. By various measurements, it is 
the largest, most symmetrical and, perhaps, the most interesting root 
system (see also Garibaldi's survey \cite{gari}). As noted in \cite{L4},
the quantity $\frac{\dim \cL}{(\dim \fh)^2}$ is bounded above, where 
$(\cL,\fh)$ is of Cartan--Killing type and $\cL$ is simple; it reaches 
its maximum ($\frac{248}{8^2}\approx 4$) for $\cL$ of type $E_8$ (whose 
existence we still have to prove). See also, for example, Ebeling 
\cite[\S 1.3]{Eb}, for an interesting
connection of $E_8$ with \nm{coding theory}, which yields a construction 
of the root system out of the classical Hamming code.
Further properties and results can be found in 
{\tt https://en.wikipedia.org/wiki/E8$\_$lattice}.
\end{rem}

\begin{xca} \label{xcadaniel} Let $A$ be an indecomposable generalized 
Cartan matrix of type (FIN). Let $\bar{A}\in M_I(\Z)$ be the matrix
with $(i,j)$-entry $|a_{ij}|$ for $i,j\in I$. Show that $\det(A)=
\det(\bar{A})$.

\noindent \footnotesize{[{\it Hint} (thanks to Daniel Juteau).
We have $a_{ij}=2\langle \alpha_i,\alpha_j\rangle/\langle \alpha_i,
\alpha_i\rangle$. Note that there is a partition $I=I^+ \sqcup I^-$ 
such that $a_{ij}=0$ for all $i\neq j$ in $I^+$ and all $i\neq j$ in 
$I^-$. Then define $\alpha_i':=\alpha_i$ if $i \in I^+$, and $\alpha_i'
:=-\alpha_i$ if $i\in I^-$. Consider the matrix $A'=(a_{ij}')_{i,j 
\in I}$ where $a_{ij}':=2\langle \alpha_i',\alpha_j'\rangle/\langle 
\alpha_i', \alpha_i'\rangle$.]}
\end{xca}

Let us fix a generalized Cartan matrix $A=(a_{ij})_{i,j\in I}$. Let 
$W=W(A)$, $\Phi=\Phi(A)$ and assume that $W(A)$ is finite. We now turn to 
the discussion of some specific properties of $W$ and $\Phi$, which can 
be derived from the classification in Section~\ref{sec3a1}. Let us fix a 
$W$-invariant scalar product $\langle \;,\;\rangle \colon E\times E 
\rightarrow \R$ as in Remark~\ref{decabst1a}. For $\alpha\in\Phi$, the 
number $\sqrt{\langle\alpha,\alpha\rangle} \in\R_{>0}$ will be called 
the \nm{length of $\alpha$}. As before, we write $\alpha^\vee:=
2\alpha/\langle \alpha,\alpha\rangle\in E$ for any $\alpha\in \Phi$.
Note that, exactly as in Section~\ref{sec1a6}, the \nm{Cauchy--Schwarz 
inequality} implies that 
\[ 0 \leq \langle \alpha^\vee,\beta\rangle\cdot \langle \alpha,
\beta^\vee\rangle <4 \qquad\mbox{where} \qquad \alpha,\beta\in \Phi, 
\beta\neq \pm \alpha.\]

\begin{rem} \label{rellen} Assume that $A$ is indecomposable. First we 
note that the arrows in the Dynkin diagrams in Table~\ref{Mdynkintbl} 
indicate the relative lengths of the roots $\alpha_i$ ($i\in I$). More 
precisely, let $i\neq j$ in $I$ be joined by a possibly multiple edge; 
then $a_{ij}<0$ and $a_{ji}<0$. We choose the notation such that 
$a_{ij}=\langle \alpha_i^\vee,\alpha_j\rangle=-1$ and $a_{ji}=
\langle\alpha_j^\vee,\alpha_i\rangle=-r$, where $r\geq 1$. Then 
\[2\frac{\langle \alpha_i,\alpha_j\rangle}{\langle\alpha_j,\alpha_j
\rangle}=a_{ji}=-r=a_{ij}r=2\frac{\langle\alpha_i,\alpha_j 
\rangle}{\langle \alpha_i,\alpha_i\rangle}r\] 
and so $\langle\alpha_i,\alpha_i\rangle=r\langle\alpha_j,\alpha_j\rangle$.
Now set $m:=\min\{ \langle\alpha_i,\alpha_i\rangle\mid i\in I\}$ and
$e:=\max\{-a_{ij}\mid i,j\in I,i\neq j, a_{ij}\neq 0\}$. By inspection 
of Table~\ref{Mdynkintbl}, we conclude that we are in one of the 
following two cases.
\begin{itemize}
\item[(a)] $e=1$ (the \nm{simply laced} case). This is the case for
$A$ of type $A_n$, $D_n$, $E_6$, $E_7$, $E_8$. Then $\langle \alpha_i,
\alpha_i\rangle=m$ for all~$i\in I$. 
\item[(b)] $e\in\{2,3\}$. This is the case for $A$ of type $B_n$, $C_n$,
$F_4$ ($e=2$) or $G_2$ ($e=3$). Then $\langle \alpha_i,\alpha_i\rangle
\in \{m,em\}$ for all~$i\in I$. 
\end{itemize}
Now consider any $\alpha \in \Phi$. By definition, we can write 
$\alpha=w(\alpha_i)$ where $i\in I$ and $w\in W$. So $\langle \alpha,
\alpha\rangle=\langle w(\alpha_i),w(\alpha_i)\rangle=\langle \alpha_i,
\alpha_i\rangle$, by the $W$-invariance of $\langle\;,\;\rangle$. Hence, 
we conclude that 
\begin{itemize}
\item[(c)] $\langle \alpha,\alpha\rangle\in \{m,em\} \qquad \mbox{for all
$\alpha \in \Phi$}.$
\end{itemize}
Thus, in case (a), all roots in $\Phi$ have the same length; in case (b), 
there are precisely two root lengths in $\Phi$ and so we may speak of 
\nm{short roots} and \nm{long roots}. In case (a), we declare all roots 
to be long roots.
\end{rem}

\begin{lem} \label{rellen1} Assume that $A$ is indecomposable. Let $e\geq 1$
be as in Remark~\ref{rellen}. Then $\langle \alpha^\vee,\beta\rangle\in\{0,
\pm 1,\pm e\}$ for all $\alpha,\beta\in \Phi$, $\beta\neq \pm\alpha$. 
\end{lem}
 
\begin{proof} Let $\alpha,\beta\in \Phi$. We can write $\alpha=
w(\alpha_i)$ for some $w\in W$ and $i\in I$. Setting $\beta':=
w^{-1}(\beta)\in\Phi$, we obtain
\[ \langle\alpha^\vee,\beta\rangle=2\frac{\langle \alpha, 
\beta\rangle}{\langle\alpha,\alpha\rangle}=2\frac{\langle w(\alpha_i), 
w(\beta')\rangle}{\langle w(\alpha_i),w(\alpha_i)\rangle}=
2\frac{\langle \alpha_i, \beta'\rangle}{\langle \alpha_i,\alpha_i 
\rangle}=\langle\alpha_i^\vee,\beta'\rangle,\]
where we used the $W$-invariance property of $\langle \;,\;\rangle$.
Writing $\beta'=\sum_{j\in I} n_j\alpha_j$ with $n_j\in\Z$, the
right hand side evaluates to $\sum_{j\in I} n_ja_{ij}\in\Z$. 
Thus, $\langle\alpha^\vee,\beta\rangle\in\Z$ for all $\alpha,
\beta\in\Phi$. Now let $\beta\neq \pm\alpha$ and assume that
$|\langle\alpha^\vee,\beta \rangle|\geq 2$. Then the fact 
that $0 \leq \langle \alpha^\vee,\beta\rangle\cdot \langle \alpha,
\beta^\vee\rangle <4$ implies that $\langle \alpha,\beta^\vee 
\rangle=\pm 1$. We conclude that 
\[\langle\alpha^\vee,\beta\rangle=2\frac{\langle\alpha,\beta
\rangle}{\langle\alpha,\alpha\rangle}=2\frac{\langle \alpha,\beta
\rangle}{\langle \beta,\beta\rangle}\frac{\langle\beta,\beta
\rangle}{\langle\alpha,\alpha\rangle}=\frac{\langle\beta,\beta
\rangle}{\langle \alpha,\alpha\rangle}\langle\alpha,\beta^\vee\rangle
=\pm \frac{\langle\beta,\beta \rangle}{\langle \alpha,\alpha\rangle}.\]
The left hand side is an integer and the right side equals $\pm e$ 
or $\pm e^{-1}$; see Remark~\ref{rellen}(c). Hence, we must have 
$\langle\alpha^\vee,\beta\rangle=\pm e$.
\end{proof}

%

\begin{exmp} \label{notG2} Assume that $A$ is indecomposable and let
$\alpha,\beta\in \Phi$ be such that $\alpha+\beta\in \Phi$. We claim 
that, if $A$ is not of type $G_2$, then $2\alpha+\beta\not\in \Phi$ or 
$\alpha+2\beta\not\in \Phi$. This is seen as follows. Assume that 
$2\alpha+\beta\in \Phi$. Then $\langle \alpha^\vee,2\alpha+\beta
\rangle=4+\langle \alpha^\vee,\beta\rangle$. If $A$ is not of type~$G_2$,
then the values of $\langle \alpha^\vee, 2\alpha+\beta \rangle$ and
$\langle \alpha^\vee,\beta\rangle$ are in $\{0,\pm 1,\pm 2\}$ by
Lemma~\ref{rellen1}. So the only possibility is that $\langle \alpha^\vee,
\beta\rangle=-2$. Now, if we also had $\alpha+2\beta\in \Phi$, then a
similar argument would show that $\langle \beta^\vee, \alpha\rangle=-2$,
contradiction to the inequality $0 \leq \langle \alpha^\vee,\beta
\rangle\cdot \langle \alpha, \beta^\vee\rangle <4$.

Note that, if $A$ is of type $G_2$, then there are examples of
roots $\alpha,\beta\in \Phi$ such that $2\alpha+\beta\in \Phi$ and
$\alpha+2\beta\in \Phi$.
\end{exmp}

\begin{xca} \label{xcastrsleq5} Assume that $A$ is indecomposable.
Let $\alpha,\beta\in \Phi^+$ and $r,s\geq 1$ be such that $r\alpha+
s\beta\in\Phi$. Check that $r,s\leq 3$ and $r+s\leq 5$. Furthermore, 
if $r+s\geq 4$, then $r\neq s$ and $A$ is of type $G_2$. 
\end{xca}

\begin{xca} \label{xcastdbase2} Assume that $A$ is indecomposable and $\Phi$ 
is simply laced. Let $\alpha, \beta\in\Phi$ be such that $\beta\neq\pm 
\alpha$. By Lemma~\ref{rellen1}, we have $\langle\alpha^\vee,
\beta\rangle\in\{0,\pm1\}$. Then show the following implications:
\begin{align*}
\langle\alpha^\vee,\beta\rangle=\;\;\,0 \quad &\Rightarrow \quad 
\beta-\alpha\not\in\Phi \;\mbox{and}\;\beta+\alpha\not\in\Phi,\\
\langle\alpha^\vee, \beta\rangle=+1 \quad &\Rightarrow \quad 
\beta-\alpha \in\Phi, \;\beta-2\alpha\not\in\Phi \;\mbox{and}\;
\beta+\alpha\not\in\Phi,\\
\langle\alpha^\vee, \beta\rangle=-1 \quad &\Rightarrow \quad 
\beta+\alpha\in\Phi,\;\beta+2\alpha\not\in\Phi \;\mbox{and}\;
\beta-\alpha\not\in\Phi.
\end{align*}
Show that, if $\alpha\in \Phi$ is written as $\alpha=\sum_{i \in I}
n_i\alpha_i$ with $n_i\in \Z$, then $\alpha^\vee=\sum_{i \in I} n_i
\alpha_i^\vee$ (see also Lemma~\ref{astring3}).
\end{xca}


\begin{xca} (Taken from \cite[p.~85]{St}) \label{xcasteinb} Assume that $A$ 
is indecomposable and that $e>1$. Let $\alpha\in \Phi$ and write $\alpha=
\sum_{i \in I} n_i\alpha_i$ where $n_i\in \Z$ for all~$i$. Show that 
$\alpha$ is a long root if and only if $e\mid n_i$ for all $i \in I$ such 
that $\alpha_i$ is a short root.
\end{xca}
 
\begin{xca} \label{xcabraid} Let $i,j\in I$, $i\neq j$. In this exercise,
we determine a formula for the order of the element $s_is_j\in W$. Show
the following.

\noindent (a) Let $E=U\oplus U'$, where $U:=\R\alpha_i+\R\alpha_j
\subseteq E$ and $U':=U^\perp\subseteq E$. Then $s_i(U) \subseteq U$, 
$s_j(U)\subseteq U$ and $s_i(u')=s_j(u')=u'$ for $u'\in U'$.

\noindent (b) For any $w\in W$, denote by $o(w)$ the order of $w$. By (a),
the order $o(s_is_j)$ is the order of $(s_is_j)|_U\colon U \rightarrow U$. 
Explicitly, we have:
\begin{alignat*}{2}
s_is_j &= s_js_i, & \quad o(s_is_j)=2 & \qquad \mbox{if $a_{ij}=0$},\\
s_is_js_i &=s_js_is_j, & \quad o(s_is_j)=3 & \qquad \mbox{if $a_{ij}
a_{ji}=1$},\\
s_is_js_is_j &=s_js_is_js_i, & \quad o(s_is_j)=4 & \qquad \mbox{if $a_{ij}
a_{ji}=2$},\\
s_is_js_is_js_is_j &=s_js_is_js_is_js_i, & \quad o(s_is_j)=6 & 
\qquad \mbox{if $a_{ij} a_{ji}=3$}.
\end{alignat*}
The above relations $s_is_j \cdots=s_js_i\cdots$ are called \nm{braid
relations}.
\end{xca}

\begin{rem} \label{weylFE} Assume that $A$ is indecomposable. In 
Section~\ref{sec1a5}, we have given explicit descriptions of the Weyl
groups $W(A)$ for $A$ of type $A_n$, $B_n$, $C_n$ or $D_n$. Now assume 
that $A$ is of type $G_2$, $F_4$, $E_6$, $E_7$ or $E_8$. For $G_2$, the 
computation in Example~\ref{cartanG2} shows that $W(A)$ is a dihedral 
group of order~$12$. For the remaining types, we use again a 
``\nm{computer algebra approach}'' to determine the order $|W(A)|$. Let 
us write $\Phi^+=\{\alpha_1,\ldots,\alpha_N\}$, where the roots are 
ordered in the same way as in Table~\ref{plusF4E678}. Then 
\[\Phi=\Phi^+\cup (-\Phi^+)=\{\alpha_1,\ldots,\alpha_N,\alpha_{N+1},
\ldots,\alpha_{2N}\} \subseteq E,\]
where $\alpha_{N+l}=-\alpha_l$ for $1\leq l\leq N$. As discussed above, 
we can identify $W(A)$ with a subgroup 
of the symmetric group $\fS_{2N} \cong \mbox{Sym}(\Phi)$. The permutation 
$\sigma_i\in \fS_{2N}$ corresponding to $s_i\in W(A)$ is obtained by 
applying $s_i$ to a root $\alpha_l$ and identifying $l'\in \{1,\ldots,
2N\}$ such that $s_i(\alpha_l)=\alpha_{l'}$; then $\sigma_i(l)=l'$. Now, 
a computer algebra system like \GAP\ \cite{gap4} contains built-in 
algorithms to work with permutation groups; in particular, there are 
efficient algorithms to determine the order of such a group\footnote{See, 
e.g., Holt et al. \cite[Chapter~4]{handbook} for the theoretical
foundations; note that, here, we certainly do not need the most 
sophisticated versions of those algorithms, since the groups in question,
and the sets on which they act, are still of moderate size.}.
In this way, we find the numbers in Table~\ref{highroot}. For example, for 
$F_4$, we obtain the following permutations in $\fS_{48}$: 
\begin{align*}
\sigma_1 & = (25,5,3,4,\ldots), \qquad\quad \sigma_2 = (5,26,6,4,\ldots),\\
\sigma_3 & = (1,9,27,7,\ldots), \qquad\quad\sigma_4  = (1,2,7,28,\ldots),
\end{align*}
where we only list $\sigma_i(l)$ for $l=1,2,3,4$; the remaining images 
are uniquely determined by these. (See Section~\ref{sec3comp} for further 
details.)
\end{rem}

\begin{table}[htbp] \caption{Highest roots and $|W(A)|$ (labelling as 
in Table~\ref{Mdynkintbl}, p.~\pageref{Mdynkintbl})} \label{highroot} 
{\small $\renewcommand{\arraystretch}{1.2}
\begin{array}{l@{\hspace{0pt}}c@{\hspace{7pt}}c} \hline \mbox{Type} & 
\mbox{Highest root $\alpha_0$} & |W(A)|\\\hline A_{n}\, (n{\geq} 1)&
\alpha_1+\alpha_2+ \ldots+\alpha_n & (n{+}1)!\\B_{n}\, (n{\geq}2)&2(\alpha_1
+\alpha_2+\ldots+ \alpha_{n-1})+\alpha_n & 2^nn!\\ C_{n}\, (n{\geq}2)&
\alpha_1+2(\alpha_2+ \ldots+\alpha_{n-1}+\alpha_n) & 2^nn!\\ D_{n}\, 
(n{\geq}3)&\alpha_1+ \alpha_2+2(\alpha_3+\ldots+\alpha_{n-1})+ \alpha_n & 
2^{n-1}n!\\ G_2 & 2\alpha_1+3\alpha_2 & 12 \\ F_4 & 2\alpha_1+3\alpha_2+
4\alpha_3+2\alpha_4 & 1152\\ E_6 & \alpha_1+2\alpha_2+2\alpha_3+3\alpha_4+
2\alpha_5+\alpha_6 & 51840\\ E_7 & 2\alpha_1+2\alpha_2+3\alpha_3+4\alpha_4+
3\alpha_5+2\alpha_6+ \alpha_7 & 2903040\\ E_8 & 2\alpha_1+3\alpha_2+
4\alpha_3+6\alpha_4+5\alpha_5 +4\alpha_6+3\alpha_7+ 2\alpha_8 & 696729600\\
\hline \end{array}$

\smallskip
$\renewcommand{\arraystretch}{1.2}\begin{array}{ccc} \hline \mbox{Type} 
& \mbox{Highest short root $\alpha_0'$} & \mbox{expression for 
$\alpha_0'^\vee$} \\\hline B_{n}& \alpha_1+\alpha_2+\ldots+\alpha_{n-1}+
\alpha_n & \alpha_1^\vee+2(\alpha_2^\vee+\ldots+\alpha_{n-1}^\vee+
\alpha_n^\vee) \\ C_{n}&\alpha_1+2(\alpha_2+\ldots+\alpha_{n-1})+\alpha_n 
&2(\alpha_1^\vee+\alpha_2^\vee+\ldots+\alpha_{n-1}^\vee)+\alpha_n^\vee \\ 
G_2 & \alpha_1+2\alpha_2 & 3\alpha_1^\vee+2\alpha_2^\vee\\ F_4 & 
\alpha_1+2\alpha_2+3\alpha_3+2\alpha_4 & 2\alpha_1^\vee+4\alpha_2^\vee+
3\alpha_3^\vee+2\alpha_4^\vee\\\hline
\end{array}$}
\end{table}

\begin{rem} \label{highestr} As in Remark~\ref{defbasephi1},
we can define a linear map $\hgt\colon E\rightarrow\R$ such that
$\hgt(\alpha_i)=1$ for all $i\in I$. If $\alpha\in \Phi$ and
$\alpha=\sum_{i\in I} n_i\alpha_i$ with $n_i\in\Z$, then $\hgt(\alpha)
=\sum_{i\in I} n_i\in\Z$ is called the \nm{height} of $\alpha$. Assuming
that $A$ is indecomposable, there is a unique root $\alpha_0\in \Phi$ 
such that $\hgt(\alpha_0)$ takes its maximum value; this root $\alpha_0$ 
is called the \nm{highest root} of $\Phi$. One can prove this by a general 
argument (see, e.g., \cite[\S 10.4]{H}), but here we can simply extract 
this from our knowledge of all root systems, using Example~\ref{cartsln} 
($A_n$), Remark~\ref{sec164} ($B_n$, $C_n$, $D_n$), 
Example~\ref{cartanG2} ($G_2$) and Table~\ref{plusF4E678} ($F_4$,
$E_6$, $E_7$, $E_8$). See Table~\ref{highroot} for explicit expressions 
of $\alpha_0$ in terms of~$\Delta$. Note the following: Let $X_n$ be one
of the diagrams in Table~\ref{Mdynkintbl}, and $\tilde{X}_n$ be the 
corresponding diagram of affine type in Table~\ref{Mdynkintbl1}.
Then $\alpha_0=\sum_{i \in I} n_i\alpha_i$ where the $n_i$ are the
numbers attached to those nodes of $\tilde{X}_n$ that are marked by
a full circle ``$\bullet$''. 
\end{rem}

\begin{xca} \label{highshort} Assume that $A$ is indecomposable and that
$e>1$. Check that $\alpha_0$ always is a long root, and that there is 
also a unique \nm{highest short root} $\alpha_0'\in \Phi$; expressions for 
$\alpha_0'$ and $\alpha_0'^\vee$ are given in Table~\ref{highroot},
where $\alpha^\vee=2\alpha/\langle \alpha,\alpha\rangle$ for 
any $\alpha\in\Phi$. 
\end{xca}

\begin{xca} \label{excdim89} Show that there are no Lie algebras of
Cartan--Killing type of dimensions $4$, $5$ or $7$. 
\end{xca}

\begin{table}[htbp] \caption{Diagram automorphisms (labelling as in 
Table~\ref{Mdynkintbl}, p.~\pageref{Mdynkintbl})} \label{Mautos} 
\begin{center} {\small
$\renewcommand{\arraystretch}{1.2} \begin{array}{ccc} \hline \mbox{Type
of $A$} & d & \mbox{orbits of $i\mapsto i^\prime$}  
\\ \hline A_{2n-1} \;(n\geq 2) & 2 & \{n\},\{n{-}1,n{+}1\},\{n{-}2,n{+}2\},
\ldots, \{1,2n{-}1\} \\ D_{n} \; (n\geq 3) & 2 & \{1,2\},\{3\},\{4\},
\ldots, \{n\} \\ D_4 & 3 & \{3\},\{1,2,4\} \\ E_6 & 2 & \{2\},\{4\},
\{3,5\},\{1,6\} \\\hline \end{array}$}
\end{center}
\end{table}

\begin{xca} \label{xcafold} Assume that $A$ is indecomposable of
\nm{simply laced} type. Thus, $A$ is of type $A_n$ ($n\geq 1$), $D_n$ 
($n\geq 3$) or $E_n$ ($n=6,7,8$). Furthermore, let $I\rightarrow I$, 
$i \mapsto i^\prime$, be a bijection such that
\begin{alignat*}{2}
a_{ij}&=a_{i^\prime j^\prime} &&\qquad \mbox{for all $i,j\in I$}.\\
a_{ii'}& =0 &&\qquad \mbox{for all $i\in I$ such that $i'\neq i$}.
\end{alignat*}
The first condition means that $i\mapsto i^\prime$ corresponds to a
symmetry of the Dynkin diagram of $A$; the second condition means that,
if $i\neq i'$, then the nodes labelled by $i$ and $i'$ are not connected
in the Dynkin diagram. Let $d\geq 1$ denote the order of the bijection
$i\mapsto i^\prime$ (as an element of the symmetric group on~$I$). The
non-trivial possibilities for $i\mapsto i^\prime$ are listed in 
Table~\ref{Mautos}. Note that there is also a non-trivial symmetry of 
order $2$ for $A$ of type $A_{2n}$ ($n \geq 1$), but the second of the
above two conditions is not satisfied in this case. 

\noindent (a) Let $\tau \colon E \rightarrow E$ be the linear map defined
by $\tau(\alpha_i):=\alpha_{i^\prime}$ for $i \in I$. Verify that
$s_{i^\prime}=\tau \circ s_i \circ \tau^{-1}$ for $i\in I$. Deduce 
that $\tau(\Phi)=\Phi$. Thus, the bijection $i\mapsto i^\prime$
induces a permutation of $\Phi$. 

\noindent (b) Let $\cL$ be a Lie algebra of Cartan--Killing type 
with structure matrix $A$. Let $\{e_i,f_i,h_i\mid i \in I\}$ be Chevalley 
generators of~$\cL$. Use the Isomorphism Theorem~\ref{isothm} to show 
that there is a Lie algebra automorphism $\tilde{\tau}\colon \cL\rightarrow 
\cL$ such that $\tilde{\tau}(e_i)=e_{i^\prime}$, $\tilde{\tau}(f_i)=
f_{i^\prime}$ and $\tilde{\tau}(h_i)=h_{i^\prime}$ for all $i\in I$. 
(Alternatively, argue as in Exercise~\ref{xcaothero}.) The map 
$\tilde{\tau}\colon \cL\rightarrow \cL$ is called a \nm{graph 
automorphism} of~$\cL$; we have $\tilde{\tau}^d=\mbox{id}_\cL$. 

\noindent (c) Let $\bB=\{h_j^+\mid j \in I\}\cup \{\be_\alpha^+\mid 
\alpha\in \Phi\}$ be Lusztig's canonical basis of $\cL$; see 
Remark~\ref{canbas1}. Show that 
\[ \tilde{\tau}(h_j^+)= h_{j^\prime}^+ \quad (j\in I) \qquad \mbox{and}
\qquad \tilde{\tau}(\be_\alpha^+)=\be_{\tau(\alpha)}^+ \quad (\alpha\in 
\Phi).\]
\noindent {\footnotesize [{\it Hints}. (a) Let $\alpha\in \Phi$ and
write $\alpha=w(\alpha_j)$ where $j \in I$ and $w\in W$. Further
write $w=s_{i_1}\cdots s_{i_r}$ where $i_1,\ldots,i_r\in I$. Then 
evaluate $\tau(\alpha)=(\tau\circ w)(\alpha_j)$ by repeatedly applying
the rule $\tau\circ s_i=s_{i^\prime} \circ \tau$, which is verified by
a straightforward computation. (b) The argument is analogous to that in
Example~\ref{stdbase3}. (c) If $\alpha\in \Phi^+$, then proceed by 
induction on $\hgt(\alpha)$; similarly for $\alpha\in \Phi^-$.]}
\end{xca}

\section{A glimpse of Kac--Moody theory} \label{sec3a3}

Let $I$ be a finite, non-empty index set and $A=(a_{ij})_{i,j\in I}\in
M_I(\C)$ be arbitrary with entries in $\C$. We would like to study Lie 
algebras for which $A$ should play the role as a ``structure matrix''. In 
order to find out how this could possibly work, let us first return to the 
case where $A$ is the true structure matrix of a Lie algebra $\cL$ of 
Cartan--Killing type with respect to an abelian subalgebra $\fh\subseteq 
\cL$ and a subset $\Delta=\{\alpha_i\mid i \in I\}$, as in
Section~\ref{sec1a2}. Then we have 
\begin{align*}
&\cL=\langle e_i,h_i,f_i \mid i\in I\rangle_{\text{alg}}
\tag{Ch0}\\
\intertext{for a suitable collection of elements $\{e_i,h_i,f_i\mid 
i\in I\} \subseteq \cL$ such that the following ``\nm{Chevalley 
relations}'' hold:}
&[e_i,f_i]=h_i \;\mbox{ and } \; [e_i,f_j]=0 \;\mbox{ for 
$i,j\in I$ such that $i\neq j$},\tag{Ch1}\\
&[h_i,h_j]=0, \; [h_i,e_j]=a_{ij}e_j, \; [h_i,f_j]=-a_{ij}f_j
\; \mbox{ for $i,j\in I$}.\tag{Ch2}
\end{align*}
(Indeed, (Ch0) holds by Proposition~\ref{genlie}; the relations in
(Ch1), (Ch2) hold by Remark~\ref{astring0}.) 


We notice that (Ch0), (Ch1), (Ch2) only refer to the collection of 
elements $\{e_i,h_i,f_i\mid i \in I\}\subseteq \cL$ and the entries of $A$, 
but not to any further structural properties of $\cL$ (e.g., finite dimension 
or $\fh$-diagonalisability). Presenting things in this way, it seems obvious 
how to proceed (but note that it is obvious only now, with hindsight): given 
any $A\in M_I(\C)$, we try to consider a Lie algebra~$\cL$ for which there 
exist elements $\{e_i,h_i,f_i\mid i \in I\}$ such that (Ch0), (Ch1), (Ch2) 
hold. Two basic questions present themselves: 
\begin{itemize}
\item Do such Lie algebras $\cL$ exist at all?
\item If yes, then does $\cL$ have interesting structural properties?
\end{itemize}
The first question is handled by the construction in 
Exercise~\ref{xcafree} (using free Lie algebras and taking suitable 
quotients). As Kac and Moody (independently) discovered in the 1960s, 
the second question also has an affirmative answer, and this has led to 
a new area of research with many interesting applications and connections, 
for example, to mathematical physics, especially when $A$ is a generalized 
Cartan matrix of type (AFF); see the monographs \cite{K}, \cite{MP}. What 
we will do in this section is the following:
\begin{itemize}
\item exhibit the ingredients of a ``triangular decomposition'' in 
any Lie algebra $\cL$ satisfying (Ch0), (Ch1), (Ch2);
\item apply these ideas to prove the Existence Theorem~\ref{myG3}.
\end{itemize}
So let us assume now that we are given any $A\in M_I(\C)$ and a Lie
algebra~$\cL$, together with elements $\{e_i,h_i,f_i\mid i\in I\}$ such 
that the conditions (Ch0), (Ch1), (Ch2) hold. In order to avoid
the discussion of trivial cases, we assume throughout that 
\[ e_j\neq 0\qquad\mbox{or}\qquad f_j\neq 0\qquad\mbox{for each 
$j\in I$}.\]
(Note that, if $e_j=f_j=0$ for some $j$, then also $h_j=0$ by (Ch1) 
and $e_j,h_j,f_j$ can simply be omitted from the collection $\{e_i,
h_i, f_i\mid i\in I\}$.)

\begin{lem} \label{deltakm} In the above setting, let $\fh:=\langle h_i\mid i
\in I \rangle_\C \subseteq \cL$. Then $\fh$ is abelian and there is a 
well-defined collection of linear maps
\[ \Delta:=\{\alpha_j\mid j\in I\}\subseteq \fh^*,\quad\mbox{where
$\alpha_j(h_i)=a_{ij}$ for all $i,j\in I$}.\]
The set $\Delta\subseteq \fh^*$ is linearly independent if and only
if $\det(A)\neq 0$. Note that, in this case, $\{h_i\mid i \in I\}$ 
is a basis of $\fh$.
\end{lem}

\begin{proof} By (Ch2), $\fh$ is an abelian 
subalgebra of~$\cL$. Next we want to define $\alpha_j\in \fh^*$ for
$j\in I$. Let $h\in \fh$ and write $h=\sum_{i\in I} x_ih_i$ 
where $x_i\in \C$. Then set $\alpha_j(h):=\sum_{i\in I} x_ia_{ij}$. 
We must show that this is well-defined. So assume that we also have 
$h=\sum_{i\in I} y_i h_i$ where $y_i\in\C$. Then $\sum_{i\in I} 
(x_i-y_i)h_i=0$; using (Ch2), we obtain: 
\[ 0=\sum_{i\in I} (x_i-y_i)[h_i,e_j]=\Bigl(\sum_{i\in I} (x_i-y_i)a_{ij}
\Bigr)e_j \quad \mbox{for any $j\in I$}.\]
If $e_j\neq 0$, then this implies that $\sum_{i\in I} x_ia_{ij}=
\sum_{i\in I}y_ia_{ij}$, as desired. If $f_j\neq 0$, then an analogous 
argument using the relation $[h_i,f_j]=-a_{ij}f_j$ yields the same 
conclusion. Thus, we obtain a well-defined subset $\Delta=\{\alpha_j
\mid j \in I\}\subseteq \fh^*$ as above. Now let $x_j\in \C$ ($j\in I$)
be such that $\sum_{j\in I} x_j\alpha_j=\underline{0}$. Then 
\[ 0=\sum_{j\in I} x_j\alpha_j(h_i)=\sum_{j\in I} a_{ij}x_j\qquad
\mbox{for all $i\in I$}.\]
If $\det(A)\neq 0$, then this implies $x_j=0$ for all $j$ and so
$\Delta$ is linearly independent. Conversely, if $\det(A)=0$, then
there exist $x_j\in\C$ ($j\in I$), not all equal to zero, such that 
$\sum_{i\in I} a_{ij}x_j=0$ for all $i\in I$. Then we also have
$\sum_{j\in I} x_j\alpha_j=\underline{0}$ and so $\Delta$ is linearly
dependent.
\end{proof}

\begin{exmp} \label{a1tilde} Let $R=\C[T,T^{-1}]$ be the ring of
Laurent polynomials over $\C$ in the indeterminate~$T$. We consider
the Lie algebra 
\[\renewcommand{\arraystretch}{0.8} 
\cL=\Bigl\{ \left(\begin{array}{cr} a & b \\ c & -a \end{array}\right)
\;\big|\; a,b,c\in R\Bigr\} \quad \bigl(=\slm_2(R)\bigr),\]
with the usual Lie bracket for matrices. A vector space basis of $\cL$
is given by $\{ T^ke_1,T^lh_1,T^mf_1\mid k,l,m\in\Z\}$, where we set 
as usual:
\[\renewcommand{\arraystretch}{0.8}
e_1:=\left(\begin{array}{cc} 0 & 1 \\ 0 & 0\end{array}\right),\qquad 
h_1:=\left(\begin{array}{cr} 1 & 0 \\ 0 & -1\end{array}\right), \qquad 
f_1:=\left(\begin{array}{cr} 0 & 0 \\ 1 & 0\end{array}\right),\]
with relations $[e_1,f_1]=h_1$, $[h_1,e_1]=2e_1$, $[h_1,f_1]=-2f_1$. 
Now set 
\[e_2:=Tf_1, \qquad h_2:=-h_1, \qquad f_2:=T^{-1}e_1.\]
Then it is straightforward to verify that the Chevalley relations 
(Ch1), (Ch2) hold with respect to the matrix
\[ \renewcommand{\arraystretch}{0.8} A=\left(\begin{array}{rr} 
2 & -2 \\ -2 & 2 \end{array}\right) \qquad\mbox{(affine type 
$\tilde{A}_1$ in Table~\ref{Mdynkintbl1})}.\]
(For example, $[h_1,h_2]=-[h_1,h_1]=0$, $[e_2,f_2]=-[f_1,e_1]=-h_1=h_2$;
furthermore, $[h_1,e_2]=T[h_1,f_1]=-2Tf_1=-2e_2$, $[e_1,f_2]=
T^{-1}[e_1,e_1]=0$ and so on.) We also note that $[e_1,e_2]=Th_1$
and $[f_1,f_2]=T^{-1}h_2$. Starting from these relations, one 
also sees that (Ch0) holds. (Details are left as an exercise for the
reader.)
\end{exmp}

Returning to the general setting, let $\fh\subseteq \cL$ be as in 
Lemma~\ref{deltakm}. Then $\dim \fh < \infty$ but we have no information 
at all about $\dim \cL$. We can still adopt a large portion of the 
definitions and results concerning weights and weight spaces from
Section~\ref{sec1a1}. For any $\lambda\in \fh^*$, we set 
\[ \cL_\lambda:=\{x\in \cL\mid [h,x]=\lambda(h)x\mbox{ for all 
$h\in \fh$}\};\]
this is a subspace of $\cL$. If $\cL_\lambda\neq \{0\}$, then $\lambda$ is 
called a \nm{weight} and $\cL_\lambda$ the corresponding 
\nm{weight space}. Since $\fh$ is abelian, we have $\fh\subseteq 
\cL_{\underline{0}}$, where $\underline{0}\in \fh^*$ is the $0$-map. The 
same argument as in Proposition~\ref{wsdprop2} shows that 
$[\cL_\lambda,\cL_\mu] \subseteq \cL_{\lambda+\mu}$ for all $\lambda,\mu 
\in \fh^*$. Let us set 
\begin{align*}
Q_{\geq 0} &:=\bigl\{\lambda\in \fh^*\mid \textstyle \lambda=\sum_{i\in I}
n_i\alpha_i \mbox{ where $n_i\in\Z_{\geq 0}$ for all~$i$}\bigr\},\\
Q_{\leq 0} &:=\bigl\{\lambda\in \fh^*\mid \textstyle \lambda=\sum_{i\in I}
n_i\alpha_i \mbox{ where $n_i\in\Z_{\leq 0}$ for all~$i$}\bigr\}.
\end{align*}
In the following discussion, some care is needed because $\Delta$ may be 
linearly dependent, and so it might happen that $Q_{\geq 0}\cap 
Q_{\leq 0}\neq \{\underline{0}\}$. 

\begin{lem} \label{triang1} In the above setting, we have 
\begin{align*}
\fn^+&:=\langle e_i\mid i\in I \rangle_{\operatorname{alg}}\subseteq 
\textstyle \sum_{\lambda\in Q_{\geq 0}} \cL_\lambda,\\
\fn^-&:=\langle f_i\mid i\in I\rangle_{\operatorname{alg}} \subseteq 
\textstyle \sum_{\lambda\in Q_{\leq 0}} \cL_\lambda.
\end{align*}
In particular, we have $[\fh,\fn^+]\subseteq \fn^+$ and $[\fh,\fn^-]
\subseteq \fn^-$.
\end{lem}

\begin{proof} Recall from Section~\ref{sec01} that $\fn^+=\langle 
X_n\mid n\geq 1\rangle_\C$, where $X_n$ consists of all Lie monomials 
in $\{e_i\mid i \in I\}$ of level~$n$. By (Ch2) and the definition of 
$\alpha_i$, we have $e_i\in \cL_{\alpha_i}$ for all $i\in I$. Hence, exactly
as in Lemma~\ref{weightgen}, one sees that $X_n\subseteq 
\bigcup_{\lambda} \cL_\lambda$, where the union runs over all 
$\lambda\in Q_{\geq 0}$ that can be expressed as $\lambda=\sum_{i\in I} 
n_i\alpha_i$ with $\sum_{i\in I} n_i=n\geq 1$. This yields that 
\[\fn^+\subseteq \sum_{\lambda\in Q_{\geq 0}}\cL_\lambda
\qquad\mbox{and}\qquad [\fh,\fn^+]\subseteq \fn^+.\]
The argument for $\fn^-$ is completely analogous, starting with the fact
that $f_i\in \cL_{-\alpha_i}$ for all $i\in I$. 
\end{proof}

\begin{lem} \label{triang2}  We have $\cL=\fn^++\fh+\fn^-$.
\end{lem}

\begin{proof} The crucial property to show is that $[f_j,\fn^+] \subseteq 
\fn^++\fh$ for all $j\in I$. This is done as follows. As in the above
proof, $\fn^+$ is spanned by Lie monomials in $\{e_i\mid i\in I\}$.
So it is sufficient to show that $[f_j,x]\in \fn^++\fh$, where $x\in \fn^+$
is a Lie monomial of level, say $n\geq 1$. We proceed by induction on~$n$.
If $n=1$, then $x=e_i$ for some~$i$ and so $[f_j,x]=-[e_i,f_j]$ is 
either zero or equal to $h_i\in \fh$. So the assertion holds in this case.
Now let $n\geq 2$. Then $x=[y,z]$ where $y,z\in \fn^+$ are Lie monomials of 
level $k$ and $n-k$, respectively; here, $1\leq k\leq n-1$. Using the 
Jacobi identity, we obtain 
\[[f_j,x]=[f_j,[y,z]]=-[y,[z,f_j]]-[z,[f_j,y]]=[y,[f_j,z]]+[[f_j,y],z].\]
By induction, we can write $[f_j,z]=z'+h$, where $z'\in \fn^+$ and 
$h\in \fh$. This yields $[y,[f_j,z]]=[y,z']+[y,h]=[y,z']-[h,y]\in \fn^++\fh$.
(We have $[y,z']\in \fn^+$ by the definition of $\fn^+$, and $[h,y]\in\fn^+$
by Lemma~\ref{triang1}.) Similarly, one sees that $[[f_j,y],z]\in 
\fn^++\fh$. 

Thus, we have shown that $[f_j,\fn^+] \subseteq \fn^++\fh$ for all $j\in I$. 
By an analogous argument, one also shows that $[e_j,\fn^-] \subseteq
\fn^-+\fh$ for all $j\in I$. Furthermore, $[e_j,\fh]\subseteq \fn^+$ and 
$[f_j,\fh]\subseteq \fn^-$ for all $j\in I$. Hence, setting $V:=\fn^++\fh+
\fn^- \subseteq \cL$, we conclude that 
\[[e_j,V]\subseteq V\qquad\mbox{and}\qquad 
[f_j,V] \subseteq V \quad\mbox{for all $j\in I$}.\] 
By Lemma~\ref{triang1}, we also have $[h_j,V]\subseteq V$. By (Ch0), we
have $\cL=\langle e_j,h_j,f_j\mid j\in I\rangle_{\text{alg}}$ and so 
Exercise~\ref{xcagenerator}(b) implies that $[\cL,V] \subseteq V$. 
In particular, $V$ is a subalgebra. Since $V$ contains all generators 
of~$\cL$, we must have $\cL=V$.
\end{proof}

\begin{xca} \label{xa1tilde} In the setting of Example~\ref{a1tilde},
we certainly have $\fh=\langle h_1,h_2\rangle_\C=\langle h_1\rangle_\C$.
Show that $T^{n+1}h_1, T^ne_1,T^ne_2\in \fn^+$ for all integers $n\geq 0$ 
(and that an analogous result holds for $\fn^-$). Use this to determine 
explicitly the subalgebras $\fn^+\subseteq \cL$ and $\fn^-\subseteq \cL$. 
Show that $\cL=\fn^+\oplus \fh \oplus \fn^-$.\\
{\footnotesize [{\it Hint}. Note that $[e_1,e_2]=Th_1$; so 
$Th_1\in \fn^+$. Now compute $(\ad_\cL(e_1)\circ \ad_\cL(e_2))^n
(h_1)$ for any $n\geq 0$. Similarly, $[f_1,f_2]=-T^{-1}h_1$; then
consider $(\ad_\cL(f_1)\circ \ad_\cL(f_2))^n(h_1)$.]}
\end{xca}

\begin{lem} \label{triang3} If $\det(A)\neq 0$, then the sum in 
Lemma~\ref{triang2} is direct; furthermore, we have $\fh=
\cL_{\underline{0}}$ and 
\[ \fn^+=\sum_{\lambda\in Q_{\geq 0}\setminus \{\underline{0}\}} 
\cL_\lambda,\qquad \fn^-=\sum_{\lambda\in Q_{\leq 0}\setminus
\{\underline{0}\}} \cL_\lambda.\]
\end{lem}

\begin{proof} By Lemma~\ref{deltakm}, the assumption that $\det(A)\neq 0$ 
implies that $\Delta=\{\alpha_i\mid i\in I\}\subseteq \fh^*$ is linearly 
independent. This has the following consequence. In the proof of
Lemma~\ref{triang1}, we have seen that $\fn^+\subseteq \sum_\lambda 
\cL_\lambda$, where the sum runs over all $\lambda\in Q_{\geq 0}$ 
that can be expressed as $\lambda=\sum_{i\in I} n_i\alpha_i$ with 
$\sum_{i\in I} n_i\geq 1$; in particular, $n_i>0$ for at least some~$i$, 
and so $\lambda\neq \underline{0}$. This shows that 
\[ \fn^+\subseteq \sum_{\lambda \in Q_+} \cL_\lambda\qquad\mbox{where}
\qquad Q_+:=\{\lambda\in Q_{\geq 0} \mid \lambda\neq \underline{0}\}.\]
Similary, we have $\fn^-\subseteq \sum_{\lambda \in Q_-} \cL_\lambda$,
where $Q_-:=\{\lambda\in Q_{\leq 0} \mid \lambda\neq \underline{0}\}$.
Combined with Lemma~\ref{triang2}, we obtain:  
\[\cL=\fn^++\fh+\fn^-\subseteq \Bigl( \sum_{\lambda\in Q_+} \cL_\lambda
\Bigr)+ \cL_{\underline{0}}+ \Bigl(\sum_{\mu\in Q_-} \cL_\mu\Bigr).\]
So it is sufficient to show that the sum on the right hand side is 
direct. Let $x\in \cL_{\underline{0}}$, $y\in \sum_{\lambda\in Q_+}
\cL_\lambda$ and $z\in \sum_{\mu\in Q_-} \cL_\mu$ be such that  
$y+x+z=0$. We must show that $x=y=z=0$. Assume, if possible, 
that $x\neq 0$. Then $x\in \cL_{\underline{0}}$ and $x=-y-z\in 
\cL_{\lambda_1}+ \ldots +\cL_{\lambda_r}$, where $r\geq 1$ and 
$\underline{0}\neq \lambda_i\in Q_+\cup Q_-$ for all~$i$. But then 
Exercise~\ref{xcaweights} (which also holds without any assumption on 
dimensions) shows that $\lambda_i=\underline{0}$ for some $i$, 
contradiction.
\end{proof}

\begin{rem} \label{triang4} Even if $\det(A)=0$, the conclusion of
Lemma~\ref{triang3} remains true, but the proof requires a more subtle 
argument; see Kac \cite[Theorem~1.2]{K} or Moody--Pianzola 
\cite[\S 4.2, Prop.~5]{MP}. If we accept this result, then the above 
discussion shows that an arbitrary Lie algebra $\cL$ satisfying 
(Ch0), (Ch1), (Ch2) is a sum of weight spaces, where $\fh$ is just the
 $\underline{0}$-weight space; furthermore, each non-zero weight is a 
$\Z$-linear combination of the set $\Delta=\{\alpha_i\mid i\in I\}
\subseteq \fh^*$, where the coefficients are either all $\geq 0$ or all 
$\leq 0$. Thus, we have a ``triangular decomposition'' $\cL=\fn^+\oplus \fh
\oplus \fn^-$, which is analogous to that in Remark~\ref{defTDa}~---~but 
now $\fn^\pm$ may be infinite-dimensional. Finally, $\cL$ is called 
\nm{integrable} if, for all $i\in I$, the linear maps 
\[\ad_\cL(e_i)\colon \cL\rightarrow \cL\quad\mbox{and}\quad \ad_\cL(f_i)
\colon \cL \rightarrow \cL\quad\mbox{are \nm{locally nilpotent}}\]
(at each $v\in \cL$; see Exercise~\ref{xcanilp}). In this case, $\cL$ is 
also called a \nm{Kac--Moody algebra}; see \cite[\S 1.3, \S 3.6]{K} or 
\cite[\S 4.1]{MP}. The integrability condition is equivalent to $A$ 
being a \nm{generalized Cartan matrix}. 

For example, a Lie algebra of Cartan--Killing type as in 
Definition~\ref{defTD} is a finite-dimensional Kac--Moody algebra. 
(Integrability holds by Lemma~\ref{wsdnil} applied to $e_i$ and 
$f_i$.) Conversely, we have:
\end{rem}

\begin{prop} \label{triang3b} Let $A=(a_{ij})_{i,j\in I}\in M_I(\C)$ 
and $\cL$ be a Lie algebra over $\C$ for which there exist elements 
$\{e_i,h_i, f_i \mid i \in I\}\subseteq \cL$ such that {\rm (Ch0)} and the 
\nm{Chevalley relations} {\rm (Ch1), (Ch2)} hold (and, for each 
$j\in I$, we have $e_j\neq 0$ or $f_j\neq 0$). Let 
\[ \fh:=\langle h_i\mid i \in I \rangle_{\C} \subseteq \cL\qquad \mbox{and}
\qquad \Delta:=\{\alpha_j\mid j \in I\}\subseteq \fh^*\]
be defined as in Lemma~\ref{deltakm}. Assume that $\dim \cL<\infty$ and 
$\det(A)\neq 0$. Then $(\cL,\fh)$ is of Cartan--Killing type with respect 
to $\Delta$; if $a_{ii}=2$ for all $i\in I$, then $A$ is the corresponding
structure matrix.
\end{prop}

\begin{proof} By Lemma~\ref{deltakm}, the set $\Delta\subseteq \fh^*$ is
linearly independent and $\{h_i\mid i \in I\}$ is a basis of $\fh$. Since
$\det(A)\neq 0$, the elements $\{h_i\mid i\in I\}$ form a basis of $\fh$
(see Lemma~\ref{deltakm}); in particular, $0\neq h_i=[e_i,f_i]$ and so
$e_i\neq 0$ for all $i\in I$. Since $e_i\in\cL_{\alpha_i}$ for all~$i\in I$,
it follows that $\Delta \subseteq P_\fh(\fg)$, as required. By
Lemma~\ref{triang3}, $\cL$ is $\fh$-diagonalisable and $\cL_{\underline{0}}
=\fh$. Furthermore, every weight $\underline{0} \neq\lambda\in P_\fh(\cL)$
belongs to $Q_+$ or $Q_-$. Thus, (CK1) and (CK2) in Definition~\ref{defTD}
hold. Finally, we have $h_i=[e_i,f_i]\in [\cL_{\alpha_i},\cL_{-\alpha_i}]$
by (Ch1). Since $\fh=\langle h_i \mid i \in I\rangle_\C$, we conclude that
(CK3) also holds. Now assume that $a_{ii}=2$ for all $i\in I$. Then
$\alpha_i(h_i)=2$ and so the elements $\{h_i\mid i \in I\}$ are the
elements required in Definition~\ref{defTD2}.
\end{proof}

\begin{lem} \label{triang3c} Assume that we are in the set-up of 
Proposition~\ref{triang3b}, where $\dim \cL<\infty$, $\det(A)\neq 0$ and
$a_{ii}=2$ for all $i\in I$. Then $A$ is a generalized Cartan matrix. 
We have the following ``\nm{Serre relations}'': 
\[\ad_\cL(e_i)^{1-a_{ij}}(e_j)=0 \qquad \mbox{and}\qquad 
\ad_\cL (f_i)^{1-a_{ij}}(f_j)=0\]
for any $i,j\in I$, $i\neq j$. (Note that $a_{ij}\in \Z_{\leq 0}$ 
for $i\neq j$.)
\end{lem}

\begin{proof} Since $(\cL,\fh)$ is of Cartan--Killing type and $a_{ii}=2$, 
the matrix $A$ is a generalized Cartan matrix by Corollary~\ref{gencart0}. 
Let $\Phi\subseteq \fh^*$ be the root system of $\cL$. Now, for $i\in I$, 
the elements $\{e_i,h_i,f_i\}$ form an $\slm_2$-triple as in 
Remark~\ref{astring0}. Let $j\in I$, $j \neq i$, and consider the 
$\alpha_i$-string through $\alpha_j$. Let $p\geq 0$ be such that 
$\alpha_j,\alpha_j+\alpha_i, \ldots, \alpha_j+p\alpha_i\in\Phi$ and 
$\alpha_j+(p+1)\alpha_i \not\in\Phi$. Using Proposition~\ref{wsdprop2}, 
we obtain
\[ \ad_\cL(e_i)^{p+1}(e_j)=\underbrace{[e_i,[e_i,[\ldots, 
[e_i}_{\text{$p+1$ times}},e_j]\ldots ]]]\in \cL_{\alpha_j+
(p+1)\alpha_i}=\{0\}\] 
and so $\ad_\cL(e_i)^{p+1}(e_j)=0$. Since $\alpha_j-\alpha_i\not\in
\Phi$, we have $a_{ij}=\alpha_j(h_i)=-p$ by Remark~\ref{astring}(a); 
this yields the desired relation. In order to obtain the analogous 
relation with $f_i,f_j$ instead of $e_i,e_j$, one can simply use 
the Chevalley involution $\omega\colon \cL \rightarrow \cL$ in 
Example~\ref{stdbase3}.
\end{proof}

\begin{rem} \label{serre} Let $A$ be an indecomposable Cartan 
matrix of type (FIN). An important theorem of Serre shows that, if 
$\cL\neq \{0\}$ is a Lie algebra such that (Ch0), (Ch1), (Ch2) and 
the ``Serre relations'' in Lemma~\ref{triang3c} hold, then $\dim \cL<
\infty$ and so we can apply Proposition~\ref{triang3b}. See Serre 
\cite[Chap.~VI, Appendix]{S} or Humphreys \cite[\S 18]{H} for further 
details; we will not need this here. In our context, it will usually be 
sufficient to apply a combination of Proposition~\ref{triang3b} and 
the Isomorphism Theorem~\ref{isothm}, without passing through the 
Serre relations (see, e.g., Proposition~\ref{corj} below).
\end{rem}

We now use the above ideas to solve a question that was left open
in Chapter~\ref{chap2}. Let $A$ be an indecomposable generalized 
Cartan matrix of type (FIN). We have seen that, if $A$ is of type 
$A_n$, $B_n$, $C_n$ or $D_n$, then $A$ arises as the structure matrix 
of a Lie algebra of Cartan--Killing type (namely, from $\cL=\slm_{n+1}(\C)$ 
or a suitable Lie algebra of classical type). But what about 
$A$ of type $G_2$, $F_4$, $E_6$, $E_7$, or $E_8$~?  For example, at 
the end of Section~\ref{sec1a6}, we saw that all the Lie brackets 
inside a Lie algebra of type $G_2$ are easily determined~---~although 
we did not know if such an algebra exists at all. (In principle, the 
same could be done for the types $F_4$, $E_6$, $E_7$ and $E_8$.) 
We now present a general solution of the existence problem.

\begin{defn}[Cf.\ \protect{\cite{G1}, \cite{L5}}] \label{myG0} Let 
$A=(a_{ij})_{i,j \in I}$ be an indecomposable generalized Cartan matrix of 
type (FIN). As in Section~\ref{sec3a2}, consider an $\R$-vector space 
$E$ with a basis $\{\alpha_i \mid i \in I\}$, and let $\Phi=\Phi(A)
\subseteq E$ be the abstract root system determined by~$A$. (We have 
$|\Phi|< \infty$ by Proposition~\ref{decabst2}.) Having obtained
the set $\Phi$, let $\bM$ be a $\C$-vector space with a basis 
\[ \bB=\{u_j \mid j\in I\}\cup \{v_\alpha\mid \alpha \in\Phi\}; 
\qquad \dim \bM=|I|+|\Phi|.\] 
Taking the formulae in Lusztig's Theorem~\ref{canbas} as a model, we 
define for $i\in I$ linear maps $\be_i\colon \bM\rightarrow \bM$
and $\bbf_i\colon \bM\rightarrow \bM$ as follows. 
\begin{align*}
\be_i(u_j) &:= |a_{ji}|v_{\alpha_i}, \qquad \bbf_i(u_j) 
:= |a_{ji}| v_{-\alpha_i},
\end{align*}
\begin{align*}
\be_i(v_\alpha) &:=\left\{
\begin{array}{c@{\hspace{5pt}}l} (q_{i,\alpha}+1) v_{\alpha+\alpha_i} & 
\quad \mbox{if $\alpha+\alpha_i \in\Phi$},\\ u_i & \quad \mbox{if 
$\alpha=-\alpha_i$},\\ 0  & \quad \mbox{otherwise},\end{array}\right.
\end{align*}
\begin{align*}
\bbf_i(v_\alpha) &:= \left\{\begin{array}{c@{\hspace{5pt}}l} (p_{i,
\alpha}+1) v_{\alpha -\alpha_i} & \quad \mbox{if $\alpha-\alpha_i 
\in\Phi$},\\ u_i & \quad \mbox{if $\alpha=\alpha_i$},\\ 0 & \quad 
\mbox{otherwise}. \end{array} \right.
\end{align*}
It is obvious that the maps $\be_i,\bbf_i$ are all non-zero. 
Now consider the Lie algebra $\gl(\bM)$, with the usual Lie bracket 
$[\varphi,\psi]=\varphi \circ \psi-\psi\circ \varphi$ for $\varphi,\psi 
\in\gl(\bM)$. We obtain a subalgebra by setting
\[\cL(A):=\langle \be_i,\bbf_i\mid i \in I
\rangle_{\text{alg}} \subseteq \gl(\bM).\]
Since $\dim \gl(\bM)<\infty$, it is clear that $\dim \cL(A)<\infty$. Our 
aim is to show that $\cL(A)$ is of Cartan--Killing type, with $A$ as
structure matrix.
\end{defn}

\begin{lem}[Cf.\ \protect{\cite[\S 3]{G1}}] \label{myG1} In the setting
of Definition~\ref{myG0}, let us also define $\bh_i:=
[\be_i,\bbf_i]\in \gl(\bM)$ for $i\in I$. Then the linear maps 
$\be_i$, $\bbf_i$, $\bh_i\in \gl(\bM)$ satisfy the 
\nm{Chevalley relations} {\rm (Ch1), (Ch2)}:
\begin{align*}
[\be_i,\bbf_j]&=0 \quad\mbox{ for all $i,j\in I$ such 
that $i\neq j$};\\
[\bh_i,\bh_j]&=0, \quad [\bh_i,\be_j]=a_{ij}
\be_j, \quad [\bh_i,\bbf_j]=-a_{ij}\bbf_j
\quad \mbox{for all $i,j\in I$}.
\end{align*}
\end{lem}

\begin{proof} Assume first that $A$ arises as the structure matrix
of a Lie algebra $\cL$ of Cartan--Killing type with respect to an
abelian subalgebra $\fh\subseteq \cL$ and a subset $\Delta=\{\alpha_i\mid
i \in I\}\subseteq \fh^*$. Thus, $A=(a_{ij})_{i,j\in I}$, where $a_{ij}=
\alpha_j(h_i)$ and $h_i\in \fh$ is defined by Proposition~\ref{wsdprop3}. 
We already discussed at the beginning of this section that then (Ch0), 
(Ch1), (Ch2) hold for $\{e_i,h_i,f_i\mid i \in I\}\subseteq \cL$, where
$e_i,f_i$ are Chevalley generators as in Remark~\ref{astring0}. 
Since $\ad_\cL\colon \cL\rightarrow \gl(\cL)$ is a 
homomorphism of Lie algebras, it follows that (Ch1), (Ch2) also hold
for the maps $\ad_\cL(e_i),\ad_\cL(f_i),\ad_\cL(h_i)\in \gl(\cL)$. Now let 
$\{\be_\alpha^+\mid \alpha\in \Phi\}$ be a collection of elements 
as in Lusztig's Theorem~\ref{canbas}. We consider the vector space $\bM:=
\cL$ and set
\[u_i:=[e_i,\be_{-\alpha_i}^+]=[f_i,\be_{\alpha_i}^+] \quad 
(i \in I), \qquad v_\alpha:=\be_{\alpha}^+\quad (\alpha\in \Phi).\] 
Then the above formulae defining $\be_i\colon \bM\rightarrow \bM$ and
$\bbf_i\colon \bM\rightarrow \bM$ correspond exactly to the formulae in 
Remark~\ref{canbas1}; in other words, we have $\be_i=\ad_\cL(e_i)$
and $\bbf_i=\ad_\cL(f_i)$ for all $i\in I$. Hence, (Ch1), (Ch2)
also hold for $\be_i,\bbf_i,\bh_i\in\gl(\bM)$.

This argument works for $A$ of type $A_n$, $B_n$, $C_n$ or $D_n$, using 
the fact, already mentioned, that then $A$ arises as the structure
matrix of $\cL=\slm_{n+1}(\C)$ or a suitable Lie algebra of classical type.
It remains to consider $A$ of type $G_2$, $F_4$, $E_6$, $E_7$ or
$E_8$. In these cases, we use again a \nm{computer algebra approach}:
we simply write down the matrices of all the $\be_i$ and $\bbf_i$ with 
respect to the above basis $\bB$ of~$\bM$, and explicitly verify (Ch1),
(Ch2) using a computer. Note that this is a finite computation since
there are only five matrices $A$ to consider and, in each case, there are
only finitely many relations to verify; see Section~\ref{sec3comp} for
further details and examples. Readers who are not happy with this argument
may consult \cite[\S 3]{G1}, where a purely theoretical, computer-free
argument is presented.
\end{proof}

Let $\cL(A)=\langle \be_i,\bbf_i\mid i \in I\rangle_{\text{alg}}
\subseteq \gl(\bM)$ be as in Definition~\ref{myG0} and set $\bh_i:=[\be_i,
\bbf_i]$ for $i\in I$. By Lemma~\ref{myG1}, the Chevalley relations (Ch1),
(Ch2) hold. Let $\fh=\langle \bh_i\mid i\in I\rangle_\C \subseteq \cL(A)$;
then $\fh$ is an abelian subalgebra. For each $j \in I$ we define
$\dot{\alpha}_j\in \fh^*$ as in Lemma~\ref{deltakm}, that is,
$\dot{\alpha}_j(h_i):=a_{ij}$ for $i \in I$. (We write $\dot{\alpha}_j$ 
in order to have a notation that is separate from $\alpha_j\in 
\Phi=\Phi(A)$.)  More generally, if $\alpha\in \Phi$, we write
$\alpha=\sum_{i \in I} n_i\alpha_i$ with $n_i\in \Z$ and set
$\dot{\alpha}:=\sum_{i \in I} n_i\dot{\alpha}_i$. Thus, we obtain a 
subset $\dot{\Phi}:=\{\dot{\alpha}\mid \alpha \in \Phi\}\subseteq \fh^*$.

\begin{thm}[Existence Theorem] \nmi{}{Existence Theorem} \label{myG3} 
With the above notation, the Lie algebra $\cL(A)\subseteq \gl(\bM)$ is 
of Cartan--Killing type with respect to $\fh\subseteq \cL(A)$ and 
$\dot{\Delta}=\{\dot{\alpha}_j\mid j \in I\}\subseteq \fh^*$, such that
$A$ is the corresponding structure matrix and $\dot{\Phi}$ is the 
set of roots with respect to~$\fh$. In particular, $\dim \cL(A)=
|I|+|\Phi|$; furthermore, since $A$ is indecomposable, $\cL(A)$ is a simple 
Lie algebra (see Theorem~\ref{cindec4}).
\end{thm}

\begin{proof} We noted in Definition~\ref{myG0} that $\be_i\neq 0$
and $\bbf_i\neq 0$ for all $i\in I$; furthermore, $\dim \cL(A)<\infty$. 
Since $\bh_i=[\be_i,\bbf_i]\in \cL(A)$, it is clear that (Ch0) holds. 
We already noted that (Ch1), (Ch2) hold. Since $A$ is of type (FIN), 
we have $\det(A) \neq 0$; furthermore, $a_{ii}=2$ for $i\in I$. 
Hence, all the assumptions of Proposition~\ref{triang3b} are satisfied 
and so $(\cL(A),\fh)$ is of Cartan--Killing type with respect to
$\dot{\Delta}=\{\dot{\alpha}_j\mid j\in I\}$ and with structure 
matrix~$A$. The fact that $\dot{\Phi}$ is the set of roots with
respect to $\fh$ follows from Remark~\ref{explicit}.
\end{proof}

\begin{cor}[Universal property of $\cL(A)$] \label{unila} 
\nmi{}{universal property of $\cL(A)$} Let $\tilde{\cL}$ be any Lie 
algebra with $\dim \tilde{\cL}< \infty$ and $\{\tilde{e}_i,\tilde{h}_i,
\tilde{f}_i\mid i \in I\}\subseteq \tilde{\cL}$ be a collection of elements
such that {\rm (Ch0), (Ch1), (Ch2)} hold (with respect to the given 
indecomposable $A$ as in Definition~\ref{myG0}) and, for all $i\in I$, 
we have $\tilde{e}_i \neq 0$ or $\tilde{f}_i\neq 0$. Then there is a
unique isomorphism of Lie algebras $\varphi\colon \cL(A) \rightarrow 
\tilde{\cL}$ such that $\varphi(\be_i)=\tilde{e}_i$, $\varphi(\bbf_i)=
\tilde{f}_i$ for
all $i\in I$.
\end{cor}

\begin{proof} Let $\tilde{\fh}=\langle \tilde{h}_i \mid i \in I\rangle_\C 
\subseteq \tilde{\cL}$ and $\tilde{\Delta}:=\{\tilde{\alpha}_j\mid j\in I\}
\subseteq \tilde{\fh}^*$ be defined as in Lemma~\ref{deltakm}, where 
$\tilde{\alpha}_j(\tilde{h}_i)=a_{ij}$ for all $i,j\in I$. Since $A$ is 
an indecomposable Cartan matrix of type (FIN), we have $\det(A)\neq 0$
and $a_{ii}=2$ for all $i\in I$. So, since also $\dim \tilde{\cL}<\infty$, 
we can apply Proposition~\ref{triang3b} which shows that $(\tilde{\cL},
\tilde{\fh})$ is of Cartan--Killing type with respect to $\tilde{\Delta}$ 
and with structure matrix $A$. So the assertion is a direct consequence 
of the Isomorphism Theorem~\ref{isothm}.
\end{proof}

A further application of Proposition~\ref{triang3b} will be given
much later, when we consider certain subalgebras of $\fg$.

\section{Using computers: {\sf CHEVIE} and {\sf ChevLie}} \label{sec3comp}

Let $A=(a_{ij})_{i,j\in I}$ be a generalized Cartan matrix with $|W(A)|
<\infty$. In this section, we explain how one can systematically deal 
with the various constructions arising from $A$ in an algorithmic fashion, 
and effectively using a computer. Several general purpose computer algebra 
systems contain built-in functions for dealing with root systems, Weyl 
groups, Lie algebras, and so on; see the online manuals of {\sf GAP} 
\cite{gap4} and {\sf Magma} \cite{magma}, for example. We introduce the 
basic features of the package {\sf ChevLie} \cite{gj}, written in the 
{\sf Julia} language (see \url{https://julialang.org}). This builds on 
the design and the conventions of the older {\sf GAP} package 
{\sf CHEVIE} \cite{chv}, \cite{jchv}. These packages are freely available 
and particularly well suited to the topics discussed here\footnote{And, 
as of 2025, {\sf ChevLie} appears to be the only package that uses 
Lusztig's fundamental Theorem~\ref{canbas} for calculations inside simple 
Lie algebras.}. Suppose you have downloaded the file {\tt chevlie1r2.jl};
then start {\sf Julia} and load {\sf ChevLie} into your current {\sf Julia} 
session:
\begin{verbatim}
  julia> include("chevlie1r2.jl"); using .ChevLie
\end{verbatim}
The central function in {\sf ChevLie} is the {\sf Julia} constructor
{\tt LieAlg}, with holds various fields with information about a Lie 
algebra of a given type (a {\sf Julia} symbol like {\tt :g}) and rank 
(a positive integer). Let us go through an example and add further 
explanations as we go along (or just type {\tt ?LieAlg} for further 
details and examples). 
\begin{verbatim}
  julia> l=LieAlg(:g,2)        # Lie algebra of type G_2
  #I dim = 14
  LieAlg('G2')
\end{verbatim}
In the background, the following happens. Firstly, given the type and
rank, there is a corresponding Cartan matrix $A$, where the labelling 
in Table~\ref{Mdynkintbl} is used. (If you wish to use a different 
labelling, then follow the instructions in the online help of 
{\tt LieAlg}.) Then {\tt LieAlg} computes some basic data related to~$A$. 
A version of the program in Table~\ref{pythontab} (p.~\pageref{pythontab})
yields the root system $\Phi$. This is stored in the component {\tt roots} 
of {\tt LieAlg}; the Cartan matrix $A$ is also stored:
\begin{verbatim}
  julia> l.cartan
    2  -1
   -3   2
  julia> l.roots
   [1, 0]  [0, 1]  [1, 1]  [1, 2]  [1, 3]  [2, 3]  
   [-1, 0] [0, -1] [-1, -1] [-1, -2] [-1, -3] [-2, -3]
\end{verbatim}
The roots are stored in terms of the list of tuples 
\[\cC(A)=\bigl\{(n_i)_{i \in I} \in \Z^I\,\big|\, \sum_{i \in I} n_i 
\alpha_i \in \Phi\bigr\}\subseteq \Z^I,\]
exactly as in Remark~\ref{explicit}. Let $N=|\Phi^+|$; this is stored as
{\tt l.N}. Then we use an enumeration of the $2N$ elements of $\Phi$ as 
follows:
\[\underbrace{\beta_1,\ldots,\beta_{|I|}}_{\text{simple roots}},
\underbrace{\beta_{|I|+1},\ldots,\beta_{N}}_{\text{further positive roots}}, 
\underbrace{-\beta_1,\ldots,-\beta_{|I|},-\beta_{|I|+1},\ldots,
-\beta_N}_{\text{negative roots}},\]
where the simple roots are those of height $1$, followed by the remaining
positive roots ordered by increasing height, followed by the negative roots.
In particular, if $A$ is indecomposable, then {\tt l.roots[l.N]} is the 
unique highest root (see Remark~\ref{highestr}). Once all roots 
are available, the permutations induced by the generators $s_i\in W$ 
($i\in I$) of the Weyl group are computed (as explained in
Remark~\ref{weylFE}) and stored. In our example: 
\begin{verbatim}
  julia> l.perms
   (7, 3, 2, 4, 6, 5, 1, 9, 8, 10, 12, 11)
   (5, 8, 4, 3, 1, 6, 11, 2, 10, 9, 7, 12)
\end{verbatim}
Here, the permutation induced by any $w\in W$ is specified by the tuple 
of integers $(j_1,\ldots,j_{2N})$ such that $w(\beta_{j_l})=\beta_l$ for 
$1\leq l \leq 2N$. (We use that convention, and not $w(\beta_l)=
\beta_{j_l}$, in order to maintain consistency with {\sf GAP} and
{\sf CHEVIE}, where permutations act from the right; for a generator
$s_i$, both conventions yield the same tuple, because $s_i$ has 
order~$2$.) Working with the permutations induced by $W$ on $\Phi$ 
immediately yields a test for equality of two elements (which would 
otherwise be difficult by working with words in the generators). 
Multiplication inside $W$ is extremely efficient: if we also have an 
element $w'\in W$ represented by $(j_1',\ldots, j_{2N}')$, then the 
product $ww' \in W$ is represented by $(j_{j_1}',\ldots,
j_{j_{2N}}')$. Thus, in our example, the permutation induced by the
element $w=s_2s_1\in W$ is obtained as follows.
\begin{verbatim}
  julia> p1=l.perms[1]; p2=l.perms[2];
  julia> ([p1[i] for i in p2]...,)   # create a tuple
  (6, 9, 4, 2, 7, 5, 12, 3, 10, 8, 1, 11)
\end{verbatim}
We will see below in Remark~\ref{basiclength2} how a permutation can
be converted back into a word in the generators of $W$. 

\begin{table}[htbp] \caption{Constructing $G_2$ using {\sf Julia} 
and {\sf ChevLie}} \label{chvlietab1} 
{\small \begin{verbatim}
  julia> l=LieAlg(:g,2)  
  julia> mats=[l.e_i[1],l.e_i[2],l.f_i[1],l.f_i[2]];
  julia> [Array(m) for m in mats]
  [...]
  # written out as 14 x 14 - matrices 
  # e_1:           e_2:           f_1:           f_2:  
  # 01000000000000 00000000000000 00000000000000 00000000000000
  # 00000000000000 00300000000000 10000000000000 00000000000000
  # 00000000000000 00020000000000 00000000000000 01000000000000
  # 00001000000000 00000100000000 00000000000000 00200000000000
  # 00000000000000 00000012000000 00010000000000 00000000000000
  # 00000023000000 00000000000000 00000000000000 00030000000000
  # 00000000100000 00000000000000 00000100000000 00000000000000
  # 00000000000000 00000000010000 00000000000000 00001000000000
  # 00000000000000 00000000003000 00000023000000 00000000000000
  # 00000000001000 00000000000000 00000000000000 00000012000000
  # 00000000000000 00000000000200 00000000010000 00000000100000
  # 00000000000000 00000000000010 00000000000000 00000000002000
  # 00000000000001 00000000000000 00000000000000 00000000000300
  # 00000000000000 00000000000000 00000000000010 00000000000000
  julia> checkrels(l,l.e_i,l.f_i,l.h_i)
  Relations OK  
  true                                 # Chevalley relations OK
\end{verbatim}}
\end{table}

Assume now that $A$ is indecomposable. Given the roots and the permutations
induced by the generators of $W$, one can then simply compute the orbits of 
the simple roots $\alpha_i$ ($i\in I$) under the action of $W$, and check 
which ones are short and which ones are long (see Remark~\ref{rellen}). 
If $A$ is simply laced, then all roots have the same length. Otherwise,
there are short roots and long roots:
\begin{verbatim}
  julia> l.short
   2 3 4 8 9 10
\end{verbatim}
Thus, $\{\beta_i\mid i=2,3,4,8,9,10\}$ are the short roots. 

Once $\Phi$ is available, it is an almost trivial matter to set up 
the matrices of the linear maps $\be_i\colon \bM\rightarrow \bM$ and 
$\bbf_i\colon\bM\rightarrow\bM$ with respect to the basis $\bB$ in 
Definition~\ref{myG0}. These are contained in the components {\tt l.e$\_$i} 
and {\tt l.f$\_$i}; there is also a component {\tt l.h$\_$i} containing 
the matrices of $\bh_i=[\be_i,\bbf_i]$ for $i\in I$. In our example, these 
matrices are printed in Table~\ref{chvlietab1}. Here, the following 
conventions are used.

\begin{itemize}
\item The basis $\bB$ is always ordered as follows:
\[ v_{\beta_N},\;\ldots,\;v_{\beta_1},\;u_1,\;\ldots,\;u_l,\;
v_{-\beta_1}, \;\ldots,\;v_{-\beta_N},\]
where $I=\{1,\ldots,l\}$. Thus, each $\be_i$ is upper triangular and 
each $\bbf_i$ is lower triangular; each $\bh_i$ is a diagonal matrix.
\item Since the matrices representing $\be_i$, $\bbf_i$, $\bh_i$ 
are extremely sparse, they are stored as {\sf Julia} {\tt SparseArrays}. 
In order to see them in full, one has to apply the {\sf Julia} function 
{\tt Array}.
\end{itemize}
Given the matrices of $\be_i$, $\bbf_i$, $\bh_i$ for all $i\in I$, one 
can then check if the Chevalley relations (Ch1), (Ch2) hold; this is 
done by the function {\tt checkrels}. We rely on these programs in the 
proof of Lemma~\ref{myG1} for Lie algebras of type $G_2$, $F_4$, $E_6$, 
$E_7$ and $E_8$. (Even for type $E_8$, this just takes a few milliseconds.)

\begin{table}[htbp] \caption{Matrix generators for the Lie 
algebra of type $G_2$} \label{matLieG2}
\begin{center}
{$\renewcommand{\arraystretch}{0.7} \renewcommand{\arraycolsep}{3pt}
e_1:\left(\begin{array}{ccccccc}
.&.&.&.&.&.&. \\  .&.&1&.&.&.&. \\  .&.&.&.&.&.&. \\ .&.&.&.&.&.&. \\  
.&.&.&.&.&1&. \\ .&.&.&.&.&.&. \\  .&.&.&.&.&.&.  \end{array}\right),\quad
e_2:\left(\begin{array}{ccccccc}
.&1&.&.&.&.&. \\  .&.&.&.&.&.&. \\  .&.&.&2&.&.&. \\ .&.&.&.&1&.&. \\  
.&.&.&.&.&.&. \\ .&.&.&.&.&.&1 \\  .&.&.&.&.&.&.  \end{array}\right)$,}\\
{$\renewcommand{\arraystretch}{0.7} \renewcommand{\arraycolsep}{3pt}
f_1:\left(\begin{array}{ccccccc}
.&.&.&.&.&.&. \\  .&.&.&.&.&.&. \\  .&1&.&.&.&.&. \\ .&.&.&.&.&.&. \\  
.&.&.&.&.&.&. \\ .&.&.&.&1&.&. \\  .&.&.&.&.&.&. \end{array}\right),\quad
f_2:\left(\begin{array}{ccccccc}
.&.&.&.&.&.&. \\  1&.&.&.&.&.&. \\  .&.&.&.&.&.&. \\ .&.&1&.&.&.&. \\  
.&.&.&2&.&.&. \\ .&.&.&.&.&.&. \\  .&.&.&.&.&1&. \end{array}\right).$}
\end{center}
{\footnotesize (Here, a dot ``.'' stands for $0$.)}
\end{table}

\begin{xca} \label{otherG2} Define matrices $e_1,e_2,f_1,f_2\in 
\gl_7(\C)$ as in Table~\ref{matLieG2}. (Note again that these matrices
only have non-negative entries.) Verify (for example, using a 
computer) that
\begin{align*}
h_1&:=[e_1,f_1]=\mbox{diag}(0, 1,-1, 0,1,-1,0),\\
h_2&:=[e_2,f_2]=\mbox{diag}(1,-1,2,0,-2,1,-1),
\end{align*}
and that the Chevalley relations (Ch1), (Ch2) hold with respect to 
the generalized Cartan matrix $A$ of type $G_2$ (as in 
Table~\ref{Mdynkintbl}). Deduce that $\cL=\langle e_1,e_2,f_1,f_2
\rangle_{\text{alg}}\subseteq \gl_7(\C)$ is a simple Lie algebra of 
type~$G_2$. How are those matrices obtained? Verify that they 
arise from the general procedure described by Jantzen \cite[\S 5A.2]{Ja};
or see Wildberger \cite{wild03b}. Similarly, realise the Lie algebra of 
type $F_4$ as a subalgebra of $\gl_{26}(\C)$.
\end{xca}

\begin{table}[htbp] 
\caption{Dynkin diagrams with $\epsilon$-function} 
\label{Mdynkineps} 
\begin{center} {\small 
\begin{picture}(290,155)
\put( 00, 25){$E_7$}
\put( 20, 25){\circle*{5}}
\put( 18, 30){$1^+$}
\put( 20, 25){\line(1,0){20}}
\put( 40, 25){\circle*{5}}
\put( 38, 30){$3^-$}
\put( 40, 25){\line(1,0){20}}
\put( 60, 25){\circle*{5}}
\put( 58, 30){$4^+$}
\put( 60, 25){\line(0,-1){20}}
\put( 60, 05){\circle*{5}}
\put( 65, 03){$2^-$}
\put( 60, 25){\line(1,0){20}}
\put( 80, 25){\circle*{5}}
\put( 78, 30){$5^-$}
\put( 80, 25){\line(1,0){20}}
\put(100, 25){\circle*{5}}
\put( 98, 30){$6^+$}
\put(100, 25){\line(1,0){20}}
\put(120, 25){\circle*{5}}
\put(118, 30){$7^-$}

\put(150, 25){$E_8$}
\put(170, 25){\circle*{5}}
\put(168, 30){$1^+$}
\put(170, 25){\line(1,0){20}}
\put(190, 25){\circle*{5}}
\put(188, 30){$3^-$}
\put(190, 25){\line(1,0){20}}
\put(210, 25){\circle*{5}}
\put(208, 30){$4^+$}
\put(210, 25){\line(0,-1){20}}
\put(210, 05){\circle*{5}}
\put(215, 03){$2^-$}
\put(210, 25){\line(1,0){20}}
\put(230, 25){\circle*{5}}
\put(228, 30){$5^-$}
\put(230, 25){\line(1,0){20}}
\put(250, 25){\circle*{5}}
\put(248, 30){$6^+$}
\put(250, 25){\line(1,0){20}}
\put(270, 25){\circle*{5}}
\put(268, 30){$7^-$}
\put(270, 25){\line(1,0){20}}
\put(290, 25){\circle*{5}}
\put(288, 30){$8^+$}

\put( 00, 59){$G_2$}
\put( 20, 60){\circle*{6}}
\put( 18, 66){$1^+$}
\put( 20, 58){\line(1,0){20}}
\put( 20, 60){\line(1,0){20}}
\put( 20, 62){\line(1,0){20}}
\put( 26, 57.5){$>$}
\put( 40, 60){\circle*{6}}
\put( 38, 66){$2^-$}

\put( 80, 60){$F_4$}
\put(100, 60){\circle*{5}}
\put( 98, 65){$1^+$}
\put(100, 60){\line(1,0){20}}
\put(120, 60){\circle*{5}}
\put(118, 65){$2^-$}
\put(120, 58){\line(1,0){20}}
\put(120, 62){\line(1,0){20}}
\put(126, 57.5){$>$}
\put(140, 60){\circle*{5}}
\put(138, 65){$3^+$}
\put(140, 60){\line(1,0){20}}
\put(160, 60){\circle*{5}}
\put(158, 65){$4^-$}

\put(190, 75){$E_6$}
\put(210, 75){\circle*{5}}
\put(208, 80){$1^+$}
\put(210, 75){\line(1,0){20}}
\put(230, 75){\circle*{5}}
\put(228, 80){$3^-$}
\put(230, 75){\line(1,0){20}}
\put(250, 75){\circle*{5}}
\put(248, 80){$4^+$}
\put(250, 75){\line(0,-1){20}}
\put(250, 55){\circle*{5}}
\put(255, 53){$2^-$}
\put(250, 75){\line(1,0){20}}
\put(270, 75){\circle*{5}}
\put(268, 80){$5^-$}
\put(270, 75){\line(1,0){20}}
\put(290, 75){\circle*{5}}
\put(288, 80){$6^+$}

\put( 00,110){$D_n$}
\put( 00,100){$\scriptstyle{n \geq 3}$}
\put( 20,130){\circle*{5}}
\put( 25,130){$1^+$}
\put( 20,130){\line(1,-1){21}}
\put( 20, 90){\circle*{5}}
\put( 26, 85){$2^+$}
\put( 20, 90){\line(1,1){21}}
\put( 40,110){\circle*{5}}
\put( 38,115){$3^-$}
\put( 40,110){\line(1,0){30}}
\put( 60,110){\circle*{5}}
\put( 58,115){$4^+$}
\put( 80,110){\circle*{1}}
\put( 90,110){\circle*{1}}
\put(100,110){\circle*{1}}
\put(110,110){\line(1,0){10}}
\put(120,110){\circle*{5}}
\put(117,115){$n^{\pm}$}

\put(170,110){$C_n$}
\put(170,100){$\scriptstyle{n \geq 2}$}
\put(190,110){\circle*{5}}
\put(188,115){$1^+$}
\put(190,108){\line(1,0){20}}
\put(190,112){\line(1,0){20}}
\put(196,107.5){$>$}
\put(210,110){\circle*{5}}
\put(208,115){$2^-$}
\put(210,110){\line(1,0){30}}
\put(230,110){\circle*{5}}
\put(228,115){$3^+$}
\put(250,110){\circle*{1}}
\put(260,110){\circle*{1}}
\put(270,110){\circle*{1}}
\put(280,110){\line(1,0){10}}
\put(290,110){\circle*{5}}
\put(287,115){$n^{\pm}$}

\put( 00,150){$A_n$}
\put( 00,140){$\scriptstyle{n \geq 1}$}
\put( 20,150){\circle*{5}}
\put( 18,155){$1^+$}
\put( 20,150){\line(1,0){20}}
\put( 40,150){\circle*{5}}
\put( 38,155){$2^-$}
\put( 40,150){\line(1,0){30}}
\put( 60,150){\circle*{5}}
\put( 58,155){$3^+$}
\put( 80,150){\circle*{1}}
\put( 90,150){\circle*{1}}
\put(100,150){\circle*{1}}
\put(110,150){\line(1,0){10}}
\put(120,150){\circle*{5}}
\put(117,155){$n^{\pm}$}

\put(170,150){$B_n$}
\put(170,140){$\scriptstyle{n \geq 2}$}
\put(190,150){\circle*{5}}
\put(188,155){$1^+$}
\put(190,148){\line(1,0){20}}
\put(190,152){\line(1,0){20}}
\put(196,147.5){$<$}
\put(210,150){\circle*{5}}
\put(208,155){$2^-$}
\put(210,150){\line(1,0){30}}
\put(230,150){\circle*{5}}
\put(228,155){$3^+$}
\put(250,150){\circle*{1}}
\put(260,150){\circle*{1}}
\put(270,150){\circle*{1}}
\put(280,150){\line(1,0){10}}
\put(290,150){\circle*{5}}
\put(287,155){$n^{\pm}$}
\end{picture}}
\end{center}
\end{table}

\begin{rem} \label{epscanr} Let $\epsilon\colon I \rightarrow\{\pm 1\}$ be
a function such that $\epsilon(j)=-\epsilon(i)$ whenever $i,j\in I$ are 
such that $a_{ij}<0$. If $A$ is indecomposable, there are precisely
two such functions: if $\epsilon$ is one of them, then the other one is
$-\epsilon$. In Table~\ref{Mdynkineps}, we have specified a particular
$\epsilon$ for each type of~$A$. This is contained in the component
{\tt epsilon} of {\tt LieAlg}:
\begin{verbatim}
  julia> l=LieAlg(:g,2); l.epsilon
    1 -1
\end{verbatim}
Once $\epsilon$ is fixed, we obtain a collection of elements 
$\{\be_\alpha^+\mid \alpha \in \Phi\}$ as in Corollary~\ref{canbash},
which we call the \nms{$\epsilon$-canonical Chevalley 
system}{epsilon-canonical Chevalley system}\footnote{Chevalley
systems in general are defined in Bourbaki \cite[Ch.~VIII, \S 2, 
D\'ef.~3]{B78}.} of~$\cL$. We shall also write $\be_\alpha^\epsilon=
\be_\alpha^+$ in order to indicate the dependence on~$\epsilon$; note 
that, if we replace $\epsilon$ by $-\epsilon$, then 
$\be_\alpha^{-\epsilon}=-\be_\alpha^\epsilon$ for all $\alpha\in\Phi$.
\end{rem}

The matrices of all $\be_\alpha^\epsilon$ ($\alpha\in \Phi$) with respect 
to $\bB$ are obtained using the function {\tt canchevbasis}. (They are 
computed using the inductive procedure in
Definition~\ref{startproof}.) For example, for type $E_8$, the matrices
have size $248 \times 248$ but they are extremely sparse; so neither 
computer memory nor computing time is an issue here. (In {\sf ChevLie}, 
they are stored as {\tt SparseArrays}, 
with signed $8$-bit integers as entries.) Once those matrices are 
available, the function {\tt structconst} computes the corresponding 
structure constants $N_{\alpha,\beta}^\epsilon$ such that 
\[[\be_\alpha^\epsilon, \be_\beta^\epsilon]=N_{\alpha,\beta}^\epsilon
\be_{\alpha+\beta}^\epsilon\quad \mbox{for $\alpha,\beta, \alpha+ 
\beta\in \Phi$}.\] 
(Again, this is very efficient since one only needs to identify one 
non-zero entry in the matrix of $\be_{\alpha+ \beta}^\epsilon$ and 
then work out only that entry in the matrix of the Lie bracket 
$[\be_\alpha^\epsilon,\be_\beta^\epsilon]$; in Example~\ref{exmcass}
below, we will see another method for computing 
those structure constants.) In our above example (where $\cL$ has type
$G_2$), we have:
\begin{verbatim}
  julia> l.roots[1:6]
   [1, 0]  [0, 1]  [1, 1]  [1, 2]  [1, 3]  [2, 3]  
  julia> structconst(l,2,4)
  (2, 4, -3, 5)
  julia> structconst(l,1,3)
  (1, 3, 0, 0)
\end{verbatim}
Here, $(2,4,-3,5)$ means that {\tt l.roots[2]+l.roots[4]=l.roots[5]}
is a root and that $N_{\alpha,\beta}^\epsilon=-3$; the output $(1,3,0,0)$
means that {\tt l.roots[1]+l.roots[3]} is not a root (and, hence,
$N_{\alpha,\beta}^\epsilon=0$).

Finally, we briefly discuss how one can work efficiently with the elements
of the Weyl group~$W$. Recall that $W=\langle s_i \mid i\in I\rangle$ and 
that $s_i^2=\id$ for all $i\in I$. Thus, every element of $W$ can be 
written as a product of various $s_i$ (but inverses of the $s_i$ are 
not required). Similarly to the height of roots, the length function
on $W$ is a crucial tool for inductive arguments. 

\begin{defn} \label{deflength} Let $w\in W$. We define the \nm{length}
of $w$, denoted $\ell(w)$, as follows. We set $\ell(\id):=0$. Now let
$w\in W$, $w\neq \id$. Then
\[ \ell(w):=\min\{r\geq 1\mid w=s_{i_1}\cdots s_{i_r}\mbox{ for some
$i_1,\ldots,i_r\in I$}\}.\]
In particular, $\ell(s_i)=1$ for all $i\in I$. If $r=\ell(w)$ and 
$i_1,\ldots,i_r\in I$ are such that $w=s_{i_1} \cdots s_{i_r}$, then 
we call this a \nm{reduced expression} for~$w$.  In general, there may 
be several reduced expressions for~$w$.
\end{defn}

\begin{rem} \label{deflength1} The formula in Remark~\ref{decabst1a} 
shows that each $s_i\in W$ ($i\in I$) is a reflection and 
so $\det(s_i)=-1$. Hence, we obtain 
\[ \det(w)=(-1)^{\ell(w)}\qquad\mbox{for any $w\in W$}.\]
Now let $w\neq \id$ and $w=s_{i_1}\cdots s_{i_r}$ be a reduced expression 
for $w$, where $r=\ell(w)$ and $i_1,\ldots,i_r\in I$. Since $s_i^{-1}=s_i$ 
for all $i\in I$, we have $w^{-1}=s_{i_r} \cdots s_{i_1}$ and so 
$\ell(w^{-1}) \leq \ell(w)$. But then also $\ell(w)=\ell((w^{-1})^{-1})
\leq \ell(w^{-1})$ and so $\ell(w)=\ell(w^{-1})$.

Now let $i\in I$. Then, clearly, $\ell(ws_i)\leq \ell(w)+1$. 
Setting $w':=ws_i\in W$, we also have $w=w's_i$ and so $\ell(w)=
\ell(w's_i)\leq \ell(w')+1=\ell(ws_i)+1$. Hence, $\ell(ws_i)\geq
\ell(w)-1$. But, since $\det(w)=(-1)^{\ell(w)}$, we can not have
$\ell(ws_i)=\ell(w)$. So we always have 
\[ \ell(ws_i)=\ell(w)\pm 1\qquad \mbox{and}\qquad \ell(s_iw)=\ell(w)
\pm 1,\]
where the second relation follows from the first by taking inverses.
\end{rem}

\begin{rem} \label{genrefl} Let $E=\langle \alpha_i\mid i \in I\rangle_\R$
and $\langle\;,\;\rangle$ be a $W$-invariant scalar product on $E$. Let
$\alpha\in\Phi$ and define $s_\alpha\colon E\rightarrow E$ by the formula
\[ s_\alpha(v):=v-2\langle \alpha^\vee,v\rangle \alpha \qquad
\mbox{for all $v\in E$}.\]
One easily sees that $s_\alpha^2=\id_E$ and that $s_\alpha$ is the 
\nm{reflection} with root~$\alpha$. (If $\alpha=\alpha_i$ for some
$i\in I$, then $s_{\alpha_i}=s_i$, as defined earlier.) We can write
$\alpha=w(\alpha_i)$ for some $i\in I$ and $w\in W$. Then we claim that
\[ s_\alpha=ws_i w^{-1} \in W.\]
Indeed, let $v\in E$ and write $v':=w^{-1}(v)\in E$. Using the 
$W$-invariance of $\langle\;,\; \rangle$, we obtain 
\[\langle \alpha_i^\vee,v'\rangle=2\frac{\langle \alpha_i,
v'\rangle}{\langle \alpha_i,\alpha_i\rangle}=
2\frac{\langle w(\alpha_i),w(v')\rangle}{\langle w(\alpha_i),
w(\alpha_i)\rangle}=2\frac{\langle \alpha, v\rangle}{\langle 
\alpha,\alpha\rangle} =\langle \alpha^\vee,v\rangle\]
and so $(ws_iw^{-1})(v)=w\bigl(s_i(v')\bigr)=v-\langle \alpha_i^\vee,v'
\rangle\alpha=s_\alpha(v)$, as claimed.  Conversely, if $s_\alpha=
ws_iw^{-1}$ for some $i\in I$ and $w\in W$, then we must have 
$w(\alpha_i)=\pm \alpha$.  Indeed, by the previous argument we have
$ws_iw^{-1}=s_{w(\alpha_i)}$ and, hence, $s_\alpha=s_{w(\alpha_i)}$. But
then $\alpha=\pm s_i(\alpha)$, since the $(-1)$-eigenspace of a reflection 
has dimension~$1$.
\end{rem}

\begin{lem}[Exchange Condition] \label{exccond} \nmi{}{Exchange Condition}
Let $w\in W$ and $i \in I$.\\
{\rm (a)} We have $\ell(s_iw)=\ell(w)-1$ if and only if $w^{-1}(\alpha_i)
\in \Phi^-$. \\
{\rm (b)} We have $\ell(ws_i)=\ell(w)-1$ if and only if $w(\alpha_i)
\in \Phi^-$. \\
{\rm (c)} Assume that $\ell(s_iw)=\ell(w)-1$ and let $w=s_{i_1}\cdots 
s_{i_r}$ be a reduced expression where $r=\ell(w)\geq 1$ and $i_1,\ldots,
i_r \in I$. Then there exists an index $j\in \{1,\ldots,r\}$ such 
that $s_is_{i_1} \cdots s_{i_{j-1}}=s_{i_1}\cdots s_{i_{j-1}}s_{i_j}$.
\end{lem}

\begin{proof} First assume that $w^{-1}(\alpha_i)\in \Phi^-$; in particular,
$w\neq \id$. We will show that $\ell(s_iw)=\ell(w)-1$ and that the 
statement in (c) holds. Let $w=s_{i_1}\cdots s_{i_r}$ be a reduced 
epxression where $r=\ell(w) \geq 1$ and $i_1,\ldots,i_r \in I$. Then 
consider the following sequence of $r+1$ roots:
\[\alpha_i,\quad s_{i_1}(\alpha_i),\quad s_{i_2}s_{i_1}(\alpha_i),
\quad \ldots,\quad s_{i_r}\cdots s_{i_1}(\alpha_i).\]
Denote them by $\beta_0,\beta_1,\ldots,\beta_r$ (from left to right).
By assumption, we have $\beta_r=s_{i_r}\cdots s_{i_1}(\alpha_i)=w^{-1}
(\alpha_i)\in \Phi^-$. On the other hand, since $\beta_0 =\alpha_i\in 
\Phi^+$, there must be some $j\in \{1,2,\ldots,r\}$ such that $\beta_0,
\beta_1,\ldots,\beta_{j-1} \in \Phi^+$ but $\beta_j \in\Phi^-$. Now 
$\beta_j=s_{i_j}(\beta_{j-1})$ and so 
\[ \beta_j=\beta_{j-1}-m\alpha_{i_j}\in \Phi^-\quad\mbox{ where }\quad
m:=\langle \alpha_{i_j}^\vee, \beta_{j-1}\rangle \in\Z.\]
Since $\beta_{j-1}\in \Phi^+$, this forces that $\beta_{j-1}=
\alpha_{i_j}$; see Lemma~\ref{wsdprop4}.  Hence, we have
$\alpha_{i_j}=\beta_{j-1}=y(\alpha_i)$ where $y:=s_{i_{j-1}}s_{i_{j-2}}
\cdots s_{i_1}\in W$. By Remark~\ref{genrefl}, this implies that $ys_i
y^{-1}=s_{i_j}$. But then
\[ s_is_{i_1}\cdots s_{i_{j-1}}=s_iy^{-1}=y^{-1}s_{i_j}=
s_{i_1}\cdots s_{i_{j-1}}s_{i_j},\]
which is the statement in (c). Inserting this into the given reduced 
expression for $w$, we obtain
\[ w=(s_{i_1}\cdots s_{i_{j-1}}s_{i_j})s_{i_{j+1}}\cdots s_{i_r}=
(s_is_{i_1}\cdots s_{i_{j-1}})s_{i_{j+1}}\cdots s_{i_r}.\]
So then $s_iw=(s_{i_1}\cdots s_{i_{j-1}})(s_{i_{j+1}}\cdots s_{i_r})$
is a product with $r-1$ factors; hence, $\ell(s_iw)=\ell(w)-1$, as desired. 

Conversely, let $w^{-1}(\alpha_i) \in \Phi^+$. Setting $w':=s_iw$, we 
have $w'^{-1}(\alpha_i)=w^{-1}(s_i(\alpha_i))=-w^{-1}(\alpha_i)\in \Phi^-$. 
Hence, by the above argument, we must have $\ell(s_iw')=\ell(w')-1$. 
Since $w=s_iw'$, this implies $\ell(s_iw)=\ell(w')=\ell(s_iw')+1=
\ell(w)+1$, as desired. Thus, (a) and (c) are proved. The equivalence
in (b) is a simple consequence of~(a); just note that $\ell(w)=
\ell(w^{-1})$ for all $w\in W$.
\end{proof}

\begin{cor} \label{basiclength1} Let $w\in W$, $w\neq \id$. Then there
exists some $i\in I$ such that $w^{-1}(\alpha_i)\in \Phi^-$. For any such $i$,
we have $w=s_iw'$, where $w'\in W$ is such that $\ell(w')=\ell(w)-1$.
\end{cor}

\begin{proof} Let $r:=\ell(w)\geq 1$ and write $w=s_{i_1}\cdots s_{i_r}$, 
where $i_1,\ldots,i_r \in I$. Then $s_{i_1}w=s_{i_2}\cdots s_{i_r}$ and so
$\ell(s_{i_1}w)<\ell(w)$. Hence, $w^{-1}(\alpha_{i_1})\in \Phi^-$ by 
Lemma~\ref{exccond}. Now take any $i\in I$ such that 
$w^{-1}(\alpha_i)\in \Phi^-$, and set $w':=s_iw$. Again, by 
Lemma~\ref{exccond}, we have $\ell(w')=\ell(w)-1$.
\end{proof}

\begin{rem} \label{basiclength2} We now obtain an efficient algorithm
for computing a reduced expression of an element $w\in W$, given as
a permutation on the roots as above. Let $(j_1,\ldots,j_{2N})$ be
the tuple representing that permutation. If $j_l=l$ for $1\leq l
\leq 2N$, then $w=\id$. Otherwise, by Corollary~\ref{basiclength1}, 
there exists some $i\in I$ such that $w^{-1}(\alpha_i)\in \Phi^-$. Using
the above conventions about the tuple $(j_1,\ldots,j_{2N})$, this
means that $j_i>N$. In order to make a definite choice, we take the 
smallest $i\in \{1,\ldots,|I|\}$ such that $j_i>N$. Then $\ell(s_iw)=
\ell(w)-1$ and we can proceed with $w':=s_iw$. In {\sf ChevLie}, this is 
implemented in the function {\tt permword}. 
\begin{verbatim}
  julia> l=LieAlg(:g,2)    # G2 with generators s1,s2
  julia> permword(l,(6,9,4,2,7,5,12,3,10,8,1,11))
  2 1                      # reduced expression s2*s1
\end{verbatim}
Conversion from a word (reduced or not), like $[2,1,2,1]$, to a 
permutation is done by the function {\tt wordperm}. 
Corollary~\ref{basiclength1} also shows how to produce all elements 
of $W$ systematically, up to a given length. Indeed, if $W(n)$ denotes 
the set of all $w\in W$ such that $\ell(w)=n$, then the set of all 
elements of length $n+1$ is obtained by taking the set of all products 
$s_iw$, where $w\in W(n)$ and $i \in I$ are such that $\ell(s_iw)=
\ell(w)+1$. This procedure is implemented in the function 
{\tt allwords}. In our above example:
\begin{verbatim}
  julia> allwords(l,3)    # elements up to length 3
  [] [1] [2] [1, 2] [2, 1] [1, 2, 1] [2, 1, 2]
\end{verbatim}
(All elements are obtained by {\tt allwords(l)}.)
\end{rem}


Next we will establish a geometric interpretation of the length function.
For $w\in W$ we set $\Phi_w^-:=\{\alpha\in\Phi^+\mid w(\alpha)\in\Phi^-\}$. 

\begin{lem} \label{phiwplus1} Let $w\in W$ and $i\in I$ be such that 
$\ell(ws_i)=\ell(w)+1$. Then $\Phi_{ws_i}^-=\{\alpha_i\}
\cup s_i(\Phi_w^-)\;$ (disjoint union).
\end{lem}

\begin{proof} First note that, since $\ell(ws_i)=\ell(w)+1$, we have 
$w(\alpha_i)\in \Phi^+$ by Lemma~\ref{exccond}. Now let $\alpha\in
\Phi_{ws_i}^-$. If $\alpha\neq \alpha_i$, then set $\beta:=s_i(\alpha)=
\alpha-m\alpha_i\in \Phi$ where $m:=-\alpha(h_i)\in \Z$. By 
Remark~\ref{wsdprop4}, we must have $\beta \in \Phi^+$. Since also 
$w(\beta)=ws_i(\alpha) \in \Phi^-$, we have $\beta \in \Phi_w^-$ and so 
$\alpha=s_i(\beta)\in s_i(\Phi_w^-)$, as desired. Conversely, if $\alpha=
\alpha_i$, then $ws_i(\alpha)=ws_i(\alpha_i)=-w(\alpha_i)\in \Phi^-$. On 
the other hand, if $\alpha=s_i(\beta)$ where $\beta \in \Phi_w^-$, then 
$ws_i(\alpha)=w(\beta)\in \Phi^-$. Furthermore, since $w(\alpha_i)\in 
\Phi^+$, we must have $\beta \neq \alpha_i$ and so $\alpha=s_i(\beta)\in 
\Phi^+$. Finally, the union $\{\alpha_i\}\cup s_i(\Phi_w^-)$ is disjoint.
Indeed, if $\alpha_i=s_i(\beta)$ for some $\beta \in \Phi_w^-$, then
$\beta=s_i(\alpha_i)=-\alpha_i\in \Phi^-$, contradiction.
\end{proof}

\begin{prop} \label{phiwplus} For $w\in W$ we have $\ell(w)=|\Phi_w^-|$. 
In particular, for $i\in I$, we have $\Phi_{s_i}^-=\{\alpha_i\}$.
\end{prop}

\begin{proof} We use induction on $\ell(w)$. If $\ell(w)=0$, then $w=\id$
and $\Phi_w^-=\varnothing$. So the assertion is clear in this case. Now 
let $\ell(w)\geq 1$ and write $w=w's_i$ where $w'\in W$ and $i\in I$ are 
such that $\ell(w)=\ell(w')+1$. By Lemma~\ref{phiwplus1}, we have
$|\Phi_w^-|=|\Phi_{w's_i}^-|=1+|s_i(\Phi_{w'}^-)|=1+|\Phi_{w'}^-|$. By 
induction, the right hand side equals $1+\ell(w')=\ell(w)$, as required.
\end{proof}

The above result provides a highly efficient method for computing
$\ell(w)$ for any $w\in W$, given as a permutation on the roots
as above: if $(j_1,\ldots,j_{2N})$ is the tuple representing that
permutation, then $\ell(w)=\ell(w^{-1})$ is just the number of 
$i\in \{1,\ldots,N\}$ such that $j_i>N$.


\begin{xca} \label{longest1} Show that there is a unique $w_0\in W$ 
such $\ell(w_0)=\max\{\ell(w)\mid w\in W\}$. Also show that $w_0(\Phi^+)=
\Phi^-$, $\ell(w_0)=|\Phi^+|$ and $w_0^2=\mbox{id}$. Furthermore, 
$\ell(w_0s_i) <\ell(w_0)$ for all $i\in I$, and this property 
characterises~$w_0$. Morevover, $\ell(ww_0)=\ell(w_0)-\ell(w)$ for all
$w\in W$. The element $w_0$ is called the \nm{longest element} of $W$.\\
{\footnotesize [{\it Hint}. Since $|W|<\infty$, there is some $w_0
\in W$ such that $\ell(w)\leq \ell(w_0)$ for all $w\in W$. So $\ell
(s_iw_0)\leq \ell(w_0)$ for all $i\in I$. Then use Lemma~\ref{exccond} 
and Proposition~\ref{phiwplus}.]}
\end{xca}

\begin{xca} \label{xcaspringer} Let $\alpha,\beta\in \Phi$ be such 
that $\beta\neq \pm\alpha$. Show that there exists some $w\in W$ such
that $w(\alpha)\in \Phi^+$ and $w(\beta)\in \Phi^+$. (A statement of
this kind can already be found in \cite[\S~I, Lemme~1]{Ch}; the current
formulation is taken from \cite[Lemma~7.5.2]{Spr}.)\\
{\footnotesize [{\it Hint}. If $\alpha,\beta\in \Phi^-$, then use the 
longest element $w_0\in W$; see Exercise~\ref{longest1}. Now assume that 
$\alpha\in \Phi^+$ and $\beta \in \Phi^-$. Then proceed by induction on 
$\hgt(\alpha)$.]}
\end{xca}


\section{Introducing Chevalley groups} \label{sec3a7}

Let $\cL$ be a Lie algebra (over $\C$, and with $\dim \cL<\infty$, as usual) 
and $\cH \subseteq \cL$ be an abelian subalgebra such that $(\cL,\cH)$ is of
Cartan--Killing type with respect to $\Delta=\{\alpha_i\mid i \in I\}
\subseteq \cH^*$. For each $i\in I$ let $\{e_i,h_i,f_i \mid i \in I\}$ be 
a corresponding $\slm_2$-triple in~$\cL$, as in Remark~\ref{astring0}.
Already in Section~\ref{sec1a4} we introduced the following 
automorphisms of $\cL$:
\begin{align*}
x_i(t)&:=\exp\bigl(t\,\ad_\cL(e_i)\bigr) \in\mbox{Aut}(\cL)\quad
\mbox{for all $i\in I$ and  $t\in \C$},\\
y_i(t)&:=\exp\bigl(t\,\ad_\cL(f_i)\bigr) \in\mbox{Aut}(\cL)\quad
\mbox{for all $i\in I$ and $t\in \C$}.
\end{align*}
Hence, we can form the subgroup 
\[G_\C(\fg):=\big\langle x_i(t),y_i(t)\mid i\in I,t \in\C\big\rangle 
\subseteq \mbox{Aut}(\cL).\]
In Definition~\ref{deflu} below we will see that one can define a similar 
group $G_K(\fg)$ over {\it any} field $K$ instead of~$\C$. As a first step,
we now use Lusztig's canonical basis 
\[ \bB=\{h_j^+\mid j \in I\}\cup \{\be_\alpha^+\mid \alpha\in \Phi\}
\subseteq \fg\] 
to explicitly write down matrices for the generators $x_i(t)$ and $y_i(t)$. 
We also assume that the additional conditions in Corollary~\ref{canbash} 
hold. Thus, there is a certain function $\epsilon \colon I \rightarrow 
\{\pm 1\}$ such that 
\[ \be_{\alpha_i}^+= \epsilon(i)e_i, \quad \be_{-\alpha_i}^+=
-\epsilon(i)f_i,\quad h_j^+=-\epsilon(i)h_i\quad\mbox{ for $i\in I$}.\]
A specific choice of $\epsilon$ for each indecomposable component of 
$I$ is defined by Table~\ref{Mdynkineps} (p.~\pageref{Mdynkineps}).
Note that the formulae in the following theorem are independent of 
those choices.

\begin{thm}[Lusztig \protect{\cite[\S 2]{L5}}] \label{luform} 
For $i\in I$ and $t\in\C$, the action of $x_i(t)$ and of $y_i(t)$ on 
$\bB$ are given by the following formulae.
\begin{empheq}[box=\widefbox]{gather*}
x_i(t)(h_j^+) = h_j^++|a_{ji}|t\be_{\alpha_i}^+,\qquad
x_i(t)(\be_{-\alpha_i}^+) =\be_{-\alpha_i}^++th_i^++t^2\be_{\alpha_i}^+,\\
x_i(t)(\be_{\alpha_i}^+) =\be_{\alpha_i}^+, \qquad x_i(t)(\be_{\alpha}^+)
=\sum_{0\leq r\leq p_{i,\alpha}} \binom{q_{i,\alpha}{+}r}{r}
t^r\be_{\alpha+ r\alpha_i}^+,\\ y_i(t)(h_j^+)=h_j^++|a_{ji}|t
\be_{-\alpha_i}^+, \qquad y_i(t)(\be_{\alpha_i}^+)=\be_{\alpha_i}^++
th_i^++t^2\be_{-\alpha_i}^+,\\ y_i(t)(\be_{-\alpha_i}^+)=\be_{-\alpha_i}^+,
\qquad y_i(t)(\be_{\alpha}^+)=\sum_{0\leq r\leq q_{i,\alpha}} 
\binom{p_{i,\alpha}{+}r}{r} t^r\be_{\alpha-r\alpha_i}^+,
\end{empheq}
where $j\in I$ and $\alpha\in \Phi$, $\alpha\neq \pm \alpha_i$. 
Here, $p_{i,\alpha},q_{i,\alpha}$ are the non-negative integers defining
the $\alpha_i$-string through $\alpha$ (see Remark~\ref{canbas0a}).
\end{thm}

\begin{proof} In the proof of Lemma~\ref{weyl5}, we already established the
following formulae, where $i\in I$, $t\in\C$ and $h\in\cH$:
\begin{align*}
x_i(t)(h)& =h-\alpha_i(h)te_i,\tag{a}\\ 
y_i(t)(h)&=h+\alpha_i(h)tf_i,\tag{b}\\ 
x_i(t)(e_i)&=e_i,\tag{c}\\ y_i(t)(e_i)&=e_i-th_i-t^2f_i.\tag{d}
\end{align*}
Now, since $h_j^+=-\epsilon(j)h_j$, we obtain using (a) that 
\[x_i(t)(h_j^+)=-\epsilon(j)h_j+\epsilon(j)\alpha_i(h_j)te_i=
h_j^+ +\epsilon(j)a_{ji}te_i.\]
In Remark~\ref{canbas1}, we saw that $[e_i,h_j^+]=\epsilon(j)a_{ji}e_i=
|a_{ji}|\be_{\alpha_i}^+$. This yields the desired formula for
$x_i(t)(h_j^+)$. Similarly, using (b), we obtain the desired formula
for $y_i(t)(h_j^+)$. The formula for $x_i(t)(\be_{\alpha_i}^+)$ 
immediately follows from (c). Analogously to (c), we have $y_i(t)(f_i)=f_i$
and this yields the formula for $y_i(t)(\be_{-\alpha_i}^+)$. Next, using
(d), we obtain:
\[ y_i(t)(\be_{\alpha_i}^+)=\epsilon(i)e_i-\epsilon(i)th_i-
\epsilon(i)t^2f_i= \be_{\alpha_i}^++th_i^++t^2\be_{-\alpha_i}^+,\]
as required. Analogously to (d), we have $x_i(t)(f_i)=f_i+th_i-t^2e_i$
and this yields the formula for $x_i(t)(\be_{-\alpha_i}^+)$. It
remains to prove the formulae for $x_i(t)(\be_\alpha^+)$ and
$y_i(t)(\be_{\alpha}^+)$, where $\alpha\neq \pm \alpha_i$. We only
do this here in detail for $x_i(t)(\be_\alpha^+)$; the argument for
$y_i(t)(\be_{\alpha}^+)$ is completely analogous. Now, by definition, we 
have
\[ x_i(t)(\be_\alpha^+)=\be_\alpha^++\sum_{r\geq 1} \frac{t^r
\ad_\cL(e_i)^r(\be_\alpha^+)}{r!}.\]
Note that $\ad_\cL(e_i)^r(\be_{\alpha}^+)\in \cL_{\alpha+r\alpha_i}=
\{0\}$ if $r>p_{i,\alpha}$. So now assume that $1\leq r\leq p_{i,\alpha}$.
Then $\alpha+\alpha_i\in\Phi$ and $\ad_\cL(e_i)(\be_\alpha^+)=[e_i,
\be_\alpha^+]=(q_{i,\alpha}+1)\be_{\alpha+\alpha_i}^+$; see (L2) in
Theorem~\ref{canbas}. Furthermore, 
\[ \ad_\cL(e_i)^2(\be_\alpha^+)=[e_i,[e_i,\be_\alpha^+]]=(q_{i,\alpha}+1)
[e_i,\be_{\alpha+\alpha_i}^+].\]
If $p_{i,\alpha}\geq 2$, then $\alpha+2\alpha_i\in\Phi$ and so the
right hand side equals $(q_{i,\alpha}+1)(q_{i,\alpha+\alpha_i}+1)
\be_{\alpha+2\alpha_i}$, again by Theorem~\ref{canbas}. Continuing in 
this way, we find that 
\[ \ad_\cL(e_i)^r(\be_\alpha^+)=(q_{i,\alpha}+1)(q_{i,\alpha+\alpha_i}+1)
\cdots (q_{i,\alpha+(r-1)\alpha_i}+1)\be_{\alpha+r\alpha_i}^+\]
for $1\leq r\leq p_{i,\alpha}$. Now note that 
\[ q_{i,\alpha+\alpha_i}=\max \{m\geq 0\mid \alpha+\alpha_i-m\alpha_i
\in\Phi\}=q_{i,\alpha}+1.\]
Similarly, $q_{i,\alpha+r\alpha_i}=q_{i,\alpha}+r$ for $1\leq r
\leq p_{i,\alpha}$. Hence, we obtain that 
\begin{align*} 
(q_{i,\alpha}+&1)(q_{i,\alpha+\alpha_i}+1)
\cdots (q_{i,\alpha+(r-1)\alpha_i}+1)\\
&=(q_{i,\alpha}+1)(q_{i,\alpha}+2)
\cdots (q_{i,\alpha}+r)=(q_{i,\alpha}+r)!/q_{i,\alpha}!
\end{align*}
Inserting this into the formula for $x_i(t)(\be_\alpha^+)$, we obtain
\[x_i(t)(\be_\alpha^+)=\sum_{r\geq 0} \frac{t^r
\ad_\cL(e_i)^r(\be_\alpha^+)}{r!}=\sum_{0\leq r \leq p_{i,\alpha}}
\frac{(q_{i,\alpha}{+}r)!}{r!\,q_{i,\alpha}!}t^r\be_{\alpha+r\alpha_i}^+,\]
and it remains to use the formula for binomial coefficients. 
\end{proof}

The above result shows that the actions of $x_i(t)$ and $y_i(t)$ on $\cL$ 
are completely determined by the structure matrix $A$ and the (abstract) 
root system~$\Phi=\Phi(A)$. As pointed out by Lusztig \cite[0.1]{L5}, 
this seems to simplify the original setting of Chevalley \cite{Ch}, 
where a number of signs appear in the formulae which depend on certain 
choices.

\begin{exmp} \label{luform1} Let $i\in I$ and $\alpha\in \Phi$ be
such that $\alpha\neq \pm \alpha_i$. If $\alpha+\alpha_i\not\in\Phi$,
then the above formulae show that $x_i(t)(\be_{\alpha}^+)=
\be_{\alpha}^+$. Similarly, if $\alpha-\alpha_i\not\in\Phi$, 
then $y_i(t)(\be_{\alpha}^+)=\be_{\alpha}^+$. Now assume that 
$\alpha+\alpha_i\in\Phi$ and that $p_{i,\alpha}=1$. Then
\[ x_i(t)(\be_\alpha^+)=\be_{\alpha}^++\binom{q_{i,\alpha}{+}1}{1}
t\be_{\alpha+\alpha_i}^+=\be_{\alpha}^++(q_{i,\alpha}+1)t
\be_{\alpha+\alpha_i}^+.\] 
Similarly, if $\alpha-\alpha_i\in\Phi$ and $q_{i,\alpha}=1$, then 
\[ y_i(t)(\be_\alpha^+)=\be_{\alpha}^++\binom{p_{i,\alpha}{+}1}{1}
t\be_{\alpha-\alpha_i}^+=\be_{\alpha}^++(p_{i,\alpha}+1)t
\be_{\alpha-\alpha_i}^+.\] 
Note that these formulae cover all cases where $A$ is of \nm{simply laced}
type, that is, all roots in $\Phi$ have the same length; see 
Exercise~\ref{xcastdbase2}. Recall from ($\spadesuit_3$) 
(p.~\pageref{spade12}) that, in general, we have $p_{i,\alpha}+
q_{i,\alpha}\leq 3$. 
\end{exmp}

\begin{rem} \label{luform2} Let $N=|\Phi^+|$ and write 
$\Phi^+=\{\beta_1,\ldots,\beta_N\}$ where the numbering is such
that $\hgt(\beta_1)\leq \hgt(\beta_2) \leq \ldots  \leq \hgt(\beta_N)$.
Let also $l=|I|$ and simply write $I=\{1,\ldots,l\}$. Then, as in
Section~\ref{sec3comp}, we order the basis $\bB$ as follows:
\[ \be_{\beta_N}^+,\;\ldots,\;\be_{\beta_1}^+,\;h_1^+,\;\ldots,\;h_l^+,
\;\be_{-\beta_1}^+, \;\ldots,\;\be_{-\beta_N}^+.\]
Let $N':=2N+l=|\bB|$ and denote the above basis elements by $v_1,
\ldots,v_{N'}$, from left to right. For $i\in I$ and $t\in \C$,
let $X_i(t)\in M_{N'}(\C)$ be the matrix of $x_i(t)$ with respect 
to the basis $\{v_1,\ldots,v_{N'}\}$; also let $Y_i(t)\in M_{N'}(\C)$ 
be the matrix of $y_i(t)$ with respect to that basis. Then the formulae in
Theorem~\ref{luform} show that 

$X_i(t)$ is an upper triangular matrix with $1$ along the diagonal,

$Y_i(t)$ is a lower triangular matrix with $1$ along the diagonal.

\noindent In particular, we have $\det(x_i(t))=\det(y_i(t))=1$. 
We also notice that each entry in $X_i(t)$ or $Y_i(t)$ is of the form
$at^r$, where the coefficient $a\in\Z$ and the exponent $r\in\Z_{\geq 0}$ 
do not depend on $t\in\C$. Now let $\Z[T]$ be the polynomial ring over
$\Z$ in an indeterminate~$T$. Replacing each entry of the form $at^r$ 
by $a T^r$, we obtain matrices 
\[X_i(T)\in M_{N'}(\Z[T]) \qquad \mbox{and}\qquad 
Y_i(T)\in M_{N'}(\Z[T]).\]
Upon substituting $T\mapsto t$ for any $t\in\C$, we get back the original 
matrices $X_i(t)$ and $Y_i(t)$. The possibility of working at a ``polynomial 
level'' will turn out to be crucial later on.
\end{rem}

\begin{exmp} \label{luform2a} Let $\cL=\slm_2(\C)$ with $I=\{1\}$ and
structure matrix $A=(2)$. We have the standard basis $\{e_1,h_1,f_1\}$, 
such that $[e_1,f_1]=h_1$, $[h_1,e_1]=2e_1$ and $[h_1,f_1]=-2f_1$.
In Exercise~\ref{xcaexp}, we already considered the automorphisms 
\[ x_1(t)=\exp\bigl(t\,\ad_\cL(e_1)\bigr) \quad \mbox{and}\quad
y_1(t)=\exp\bigl(t\,\ad_\cL(f_1)\bigr)\qquad (t\in\C)\]
and worked out the corresponding matrices. Now note that $\bB=
\{e_1, -h_1,-f_1\}$ (see the remark just after Theorem~\ref{canbas}). Hence,
\[\renewcommand{\arraystretch}{0.8} 
X_1(t)=\left(\begin{array}{ccc} 1 & 2t & t^2 \\ 0 & 1 & t \\ 0 & 0 & 1 
\end{array}\right) \quad\mbox{and}\quad 
Y_1(t)=\left(\begin{array}{ccc} 1 & 0 & 0 \\ t & 1 & 0 \\ t^2 & 2t & 1
\end{array}\right).\]
So, obviously, we have the following matrices over $\Z[T]$:
\[\renewcommand{\arraystretch}{0.8} 
X_1(T)=\left(\begin{array}{ccc} 1 & 2T & T^2 \\ 0 & 1 & T \\ 0 & 0 & 1 
\end{array}\right) \quad\mbox{and}\quad 
Y_1(T)=\left(\begin{array}{ccc} 1 & 0 & 0 \\ T & 1 & 0 \\ T^2 & 2T & 1
\end{array}\right).\]
\end{exmp}

We now show how the definition of $G_\C(\fg)$ can be extended to an 
arbitrary field~$K$. We usually attach a bar to objects defined over~$K$. 
So let $\bar{\cL}$ be a vector space\footnote{This vector space
$\bar{\cL}$ also inherits a Lie algebra structure from $\cL$; see 
Carter \cite[\S 4.4]{Ca1}. But we will not need this here.} over 
$K$ with a basis
\[ \bar{\bB}=\{\bar{h}_j^+\mid j\in I\}\cup\{\bar{\be}_\alpha^+
\mid\alpha\in\Phi\}.\]
For $i\in I$ and $\zeta\in K$ we use the formulae in Theorem~\ref{luform}
to define linear maps $\bar{x}_i(\zeta)\colon \bar{\cL}
\rightarrow\bar{\cL}$ and $\bar{y}_i(\zeta)\colon \bar{\cL}
\rightarrow\bar{\cL}$. Explicitly, we set:
\begin{gather*}
\bar{x}_i(\zeta)(\bar{h}_j^+):=\bar{h}_j^++|a_{ji}|\zeta
\bar{\be}_{\alpha_i}^+, \quad\; \bar{x}_i(\zeta)(\bar{\be}_{-\alpha_i}^+)
:=\bar{\be}_{-\alpha_i}^+ +\zeta\bar{h}_i^++\zeta^2\bar{\be}_{\alpha_i}^+,
\\ \bar{x}_i(\zeta)(\bar{\be}_{\alpha_i}^+) :=\bar{\be}_{\alpha_i}^+, 
\qquad \bar{x}_i(\zeta)(\bar{\be}_{\alpha}^+) :=\sum_{0\leq r\leq p_{i,
\alpha}} \binom{q_{i,\alpha}{+}r}{r} \zeta^r\bar{\be}_{\alpha+ 
r\alpha_i}^+,
\end{gather*} 
\begin{gather*} \bar{y}_i(\zeta)(\bar{h}_j^+):=\bar{h}_j^++|a_{ji}|
\zeta \bar{\be}_{-\alpha_i}^+, \quad\; \bar{y}_i(\zeta)
(\bar{\be}_{\alpha_i}^+):=\bar{\be}_{\alpha_i}^++ \zeta\bar{h}_i^++
\zeta^2\bar{\be}_{-\alpha_i}^+,\\ \bar{y}_i(\zeta)
(\bar{\be}_{-\alpha_i}^+):=\bar{\be}_{-\alpha_i}^+, \qquad \bar{y}_i
(\zeta)(\bar{\be}_{\alpha}^+):=\sum_{0\leq r\leq q_{i,\alpha}}
\binom{p_{i,\alpha}{+}r}{r} \zeta^r\bar{\be}_{\alpha-r\alpha_i}^+,
\end{gather*}
where $j\in I$ and $\alpha\in \Phi$, $\alpha\neq \pm \alpha_i$.
(Here, the product of an integer in $\Z$ and an element of $K$ is
defined in the obvious way.) Let $\bar{X}_i(\zeta)$ and $\bar{Y}_i(\zeta)$
be the matrices of $\bar{x}_i(\zeta)$ and $\bar{y}_i(\zeta)$, respectively,
with respect to~$\bar{\bB}$, where the elements of $\bar{\bB}$ are 
arranged as in Remark~\ref{luform2}. Then the above formulae show again 
that

$\bar{X}_i(\zeta)$ is upper triangular with $1$ along the diagonal,

$\,\bar{Y}_i(\zeta)\,$ is lower triangular with $1$ along the 
diagonal.

\noindent In particular, we have $\det(\bar{x}_i(\zeta))=\det
(\bar{y}_i(\zeta))=1$. Note that, if $K=\C$, then $\bar{x}_i(\zeta)=
x_i(\zeta)$ and $\bar{y}_i(\zeta)=y_i(\zeta)$ for all $\zeta\in\C$. 

\begin{defn} \label{deflu} Following Chevalley \cite{Ch} and Lusztig 
\cite[\S 2]{L5}, the \nm{Chevalley group}\footnote{More precisely, 
$G_K(\cL)$ is a Chevalley group of ``\nmi{adjoint type}{Chevalley group of 
adjoint type}''. (More general types of groups will be constructed in
Chapter~\ref{chap4}.) Chevalley \cite[\S IV]{Ch} denotes this group by 
$G_K^\prime(\fg)$ because there is a slightly larger group containing 
additional ``diagonal elements''; see Carter \cite[\S 7.1]{Ca1} and 
Chevalley \cite[p.~37]{Ch}. Following Steinberg \cite{St}, we will just 
consider $G_K(\fg)$ as defined above. If $K$ is an algebraically closed 
field, then the difference between $G_K(\cL)$ and the slightly larger group 
disappears. See further comments in Remark~\ref{remfoot} and 
Proposition~\ref{steinb35} in Chapter~\ref{chap4} below.} of type $\cL$ 
over the field~$K$ is 
defined by 
\[ G_K(\fg):=\big\langle \bar{x}_i(\zeta),\bar{y}_i(\zeta)\mid 
i\in I,\zeta\in K\big\rangle\subseteq \GL(\bar{\cL}).\]
If there is no danger of confusion, then we just write $\bar{G}$ instead
of $G_K(\fg)$ (where the bar is meant to indicate that we are working 
over~$K$ and not over~$\C$). Note that $G_K(\fg)$ is completely determined 
by the structure matrix~$A$ of~$\fg$, the corresponding (abstract) root 
system $\Phi$, and the field~$K$. If $K=\C$, then $G_\C(\fg)$ is the
group defined at the beginning of this section. Also note that, if $K$ 
is a finite field, then $G_K(\fg)$ is a finite group. 
\end{defn}

Chevalley \cite{Ch} showed that, if $A$ is indecomposable, then  $G_K(\fg)$ 
is a simple group in almost all cases; the finitely many exceptions only 
occur when $|I|\leq 2$ and $K$ is a field with $2$ or $3$ elements. 
As already mentioned, this discovery had a profound influence on the 
further development of group theory and Lie theory in general. 

\begin{exmp} \label{luform2aa} Let $\cL=\slm_2(\C)$. In 
Example~\ref{luform2a}, we determined the matrices of $x_1(t)$ and 
$y_1(t)$ for $t\in \C$. Now let $K$ be any field and $\zeta\in K$. Then 
the matrices of $\bar{x}_1(\zeta)$ and $\bar{y}_1(\zeta)$ are given by 
\[\renewcommand{\arraystretch}{0.8} 
\bar{X}_1(\zeta)=\left(\begin{array}{ccc} 1 & 2\zeta & \zeta^2 \\ 
0 & 1 & \zeta \\ 0 & 0 & 1 \end{array}\right) \quad\mbox{and}\quad 
\bar{Y}_1(\zeta)=\left(\begin{array}{ccc} 1 & 0 & 0 \\ 
\zeta & 1 & 0 \\ \zeta^2 & 2\zeta & 1 \end{array}\right).\]
In Section~\ref{sec3fin} we will see that $G_K(\fg)=\langle \bar{x}_1
(\zeta),\bar{y}_1(\zeta)\mid \zeta\in K\rangle$ is isomorphic to 
$\SL_2(K)/\{\pm I_2\}$.
\end{exmp}

\begin{rem} \label{xinontriv} The definition immediately shows
that $\bar{x}_i(0)=\id_{\bar{\cL}}$ and $\bar{y}_i(0)=\id_{\bar{\cL}}$. 
Now let $0\neq \zeta\in K$. Then 
\[\bar{x}_i(\zeta)(\bar{\be}_{-\alpha_i}^+)=\bar{\be}_{-\alpha_i}^+ 
+\zeta \bar{h}_i^++\zeta^2\bar{\be}_{\alpha_i}^+\neq 
\bar{\be}_{-\alpha_i}^+\]
and so $\bar{x}_i(\zeta)\neq \id_{\bar{\cL}}$. In fact, this shows
that the map $\zeta \mapsto \bar{x}_i(\zeta)$ is injective. Similarly, 
one sees that the map $\zeta \mapsto \bar{y}_i(\zeta)$ is injective. 
\end{rem}

Of course, one would hope that the elements $\bar{x}_i(\zeta)$ and 
$\bar{y}_i(\zeta)$ (over~$K$) have further properties analogous to those 
of $x_i(t)$ and $y_i(t)$ (over~$\C$). In order to justify this in concrete
cases, some extra argument is usually required because the definition of 
$\bar{x}_i(\zeta)$ or $\bar{y}_i(\zeta)$ in terms of an exponential 
construction is not available over~$K$ (at least not if $K$ has positive
characteristic). For this purpose, we make crucial use of the possibility 
of working at a ``polynomial level'', as already mentioned in 
Remark~\ref{luform2}. Here is a simple first example.

\begin{lem} \label{fromctok1} Let $i\in I$. Then $\bar{x}_i(\zeta)^{-1}=
\bar{x}_i(-\zeta)$ and $\bar{y}_i(\zeta)^{-1}=\bar{y}_i(-\zeta)$ 
for all $\zeta\in K$. Furthermore, $\bar{x}_i(\zeta+\zeta')=
\bar{x}_i(\zeta)\bar{x}_i(\zeta')$ and $\bar{y}_i(\zeta+\zeta')=
\bar{y}_i(\zeta)\bar{y}_i(\zeta')$ for all $\zeta,\zeta'\in K$. 
\end{lem}

\begin{proof} First we prove the assertion about $\bar{x}_i(\zeta)^{-1}$.
(This would also follow from the assertion about $\bar{x}_i(\zeta+
\zeta')$ and the fact that $\bar{x}_i(0)=\id_{\bar{\cL}}$, but it
may be useful to run the two arguments separately, since they involve
different ingredients.) Let $\Z[T]$ be the polynomial ring over $\Z$ 
with indeterminate~$T$. Let $X_i(T)\in M_{N'}(\Z[T])$ be the matrix 
defined in Remark~\ref{luform2}; upon substituting $T\mapsto t$ for any 
$t\in\C$, we obtain the matrix of the element $x_i(t)\in G_\C(\cL)$. 
We claim that 
\[ X_i(T)\cdot X_i(-T)=I_{N'} \qquad \mbox{(equality in 
$M_{N'}(\Z[T])$)},\]
where $I_{N'}$ denotes the $N'\times N'$-times identity matrix.
This is seen as follows. Let $f_{rs}\in \Z[T]$ be the $(r,s)$-entry
of $X_i(T)$. Writing out the matrix product $X_i(T)\cdot X_i(-T)$, we 
must show that the following identities of polynomials in 
$\Z[T]$ hold for all $r,s\in\{1,\ldots,N'\}$:
\[ \textstyle{\sum_{r'}} f_{rr'}(T)f_{r's}(-T)=\left\{\begin{array}{cl}
1 & \mbox{ if $r=s$},\\ 0 & \mbox{ if $r\neq s$}.\end{array}\right.\]
Since $x_i(t)x_i(-t)=\id_\cL$ (see Lemma~\ref{exponential}), we have 
$X_i(t)\cdot X_i(-t)=I_{N'}$ for all $t\in \C$, which means that 
\[ \textstyle{\sum_{r'}} f_{rr'}(t)f_{r's}(-t)=\left\{\begin{array}{cl}
1 & \mbox{ if $r=s$},\\ 0 & \mbox{ if $r\neq s$}.\end{array}\right.\]
So the assertion follows from the general fact that, if $g,h\in\Z[T]$ 
are such that $g(t)=h(t)$ for infinitely many $t \in\C$, then $g=h$ 
in $\Z[T]$. 

Now fix $\zeta\in K$. By the universal property of $\Z[T]$, we have 
a canonical ring homomorphism $\varphi_\zeta\colon \Z[T] \rightarrow K$ 
such that $\varphi_\zeta(T)=\zeta$ and $\varphi_\zeta(m)=m\cdot 1_K$ 
for $m\in\Z$. Applying $\varphi_\zeta$ to the entries of $X_i(T)$,
we obtain the matrix $\bar{X}_i(\zeta)\in M_{N'}(K)$, by the above 
definition of $\bar{x}_i(\zeta)$. Similarly, applying $\varphi_\zeta$ 
to the entries of $X_i(-T)$, we obtain the matrix $\bar{X}_i(-\zeta)
\in M_{N'}(K)$. Since $\varphi_\zeta$ is a ring homomorphism, the 
identity $X_i(T) \cdot X_i(-T)=I_{N'}$ over $\Z[T]$ implies the 
identity $\bar{X}_i(\zeta) \cdot \bar{X}_i(-\zeta)=\bar{I}_{N'}$ 
over~$K$. Consequently, we have $\bar{x}_i(\zeta) \bar{x}_i(-\zeta)=
\id_{\bar{\cL}}$, as desired. The argument for $\bar{y}_i(\zeta)$ is 
completely analogous.

Now consider the assertion about $\bar{x}_i(\zeta+\zeta')$. First 
we work over~$\C$. For $t,t'\in \C$, the derivations $t\,\ad_\cL
(e_i)$ and $t'\,\ad_\cL(e_i)$ of $\cL$ certainly commute with each 
other. Hence, Exercise~\ref{expocomm} shows that 
\begin{align*}
x_i(t+t') &=\exp\bigl(t\,\ad_\cL(e_i)+t'\,\ad_\cL(e_i)\bigr)\\ 
&=\exp\bigl(t\,\ad_\cL(e_i)\bigr)\circ \exp\bigl(t'\ad_\cL(e_i)
\bigr)=x_i(t)x_i(t'),
\end{align*}
where we omit the symbol ``$\circ$'' for the multiplication inside $G_\C
(\cL)$. We ``lift'' again the above identity to a ``polynomial level'', 
where now we work over $\Z[T,T']$, the polynomial ring in two 
commuting indeterminates $T,T'$ over $\Z$. Regarding $X_i(T)$ and 
$X_i(T')$ as matrices in $M_{N'}(\Z[T,T'])$, we claim that 
\[X_i(T+T')=X_i(T)\cdot X_i(T') \qquad \mbox{(equality in 
$M_{N'}(\Z[T,T'])$)}.\]
This is seen as follows. Let again $f_{rs}\in \Z[T]$ be the $(r,s)$-entry 
of $X_i(T)$. Writing out the above matrix product, we must show that the
following identities in $\Z[T,T']$ hold for all $r,s\in \{1,\ldots,N'\}$:
\[ f_{rs}(T+T')=\textstyle{\sum_{r'}} f_{rr'}(T)f_{r's}(T').\]
We have just seen that these identities do hold upon substituting 
$T\mapsto t$ and $T'\mapsto t'$ for any $t,t'\in\C$. Hence, the assertion
now follows from the general fact that, if $g,h\in\Z[T,T']$ are any 
polynomials such that $g(t,t')=h(t,t')$ for all $t,t'\in\C$, then $g=h$ 
in $\Z[T,T']$. (Proof left as an exercise; the analogous statement is 
also true for polynomials in several commuting variables.) Now fix
$\zeta,\zeta'\in K$. Then we have a canonical ring homomorphism 
$\varphi_{\zeta,\zeta'}\colon \Z[T,T'] \rightarrow K$ such that 
$\varphi_{\zeta,\zeta'}(T)=\zeta$, $\varphi_{\zeta,\zeta'}(T')=\zeta'$ 
and $\varphi_{\zeta,\zeta'}(m)=m\cdot 1_K$ for $m\in\Z$. Applying 
$\varphi_{\zeta,\zeta'}$ to the entries of $X_i(T)$, $X_i(T')$ and
$X_i(T+T')$, we obtain the matrices $\bar{X}_i(\zeta)$, 
$\bar{X}_i(\zeta')$ and $\bar{X}_i(\zeta+\zeta')$. Consequently,
the identity $X_i(T+T')=X_i(T) \cdot X_i(T')$ over $\Z[T,T']$ implies the 
identity $\bar{X}_i(\zeta+\zeta')=\bar{X}_i(\zeta)\cdot \bar{X}_i(\zeta')$
over~$K$. Hence, we have $\bar{x}_i(\zeta+\zeta')=\bar{x}_i(\zeta)
\bar{x}_i(\zeta')$, as desired. The argument for $\bar{y}_i(\zeta+
\zeta')$ is analogous.
\end{proof}

We will see similar arguments, or variations thereof, frequently in the 
development to follow. The following result will be very useful. 

\begin{lem} \label{genchev0} Let $x\in \cL$ be such that $\ad_\cL(x)\colon 
\cL\rightarrow \cL$ is nilpotent. Let $\theta\colon \cL\rightarrow \cL$ be 
any Lie algebra automorphism. Then $\ad_\cL(\theta(x))$ is nilpotent 
and $\exp\bigl(\ad_\cL(\theta(x))\bigr)=\theta \circ
\exp(\ad_\cL(x)) \circ \theta^{-1}$.
\end{lem}

\begin{proof} Let $y\in \cL$. Since $\theta$ is an automorphism, we have 
for $m\geq 0$:
\begin{align*}
\ad_\cL(\theta&(x))^m(y)=[\underbrace{\theta(x),[\theta(x),\ldots,
[\theta(x)}_{\text{$m$ terms}},\theta\bigl(\theta^{-1}(y)\bigr)]
\ldots]]\\&=\theta\bigl([\underbrace{x,[x,\ldots,[x}_{\text{$m$ terms}},
\theta^{-1}(y)]\ldots]]\bigr)=\theta\bigl(\ad_\cL(x)^m(\theta^{-1}(y))\bigr).
\end{align*}
Hence, since $\ad_\cL(x)^d=0$ for some $d\geq 1$, we also have 
$\ad_\cL(\theta(x))^d=0$, that is, $\ad_\cL(\theta(x))$ is nilpotent. 
The above identity also yields:
\begin{align*}
\bigl(\theta &\circ \exp(\ad_\cL(x))\circ \theta^{-1}\bigr)(y)=\theta
\Bigl(\sum_{m\geq 0} \frac{1}{m!}\ad_\cL(x)^m\bigl(\theta^{-1}(y)\bigr) 
\Bigr)
\end{align*}
\begin{align*}
&\quad =\sum_{m\geq 0} \frac{1}{m!}\theta\bigl(\ad_\cL(x)^m
(\theta^{-1}(y))\bigr)=\sum_{m\geq 0} \frac{1}{m!}\ad_\cL(\theta(x))^m(y),
\end{align*}
which equals $\exp\bigl(\ad_\cL(\theta(x))\bigr)(y)$, as required.
\end{proof}

\begin{exmp} \label{xcaomega1} Consider the Chevalley involution
$\omega \colon \cL\rightarrow \cL$ in Example~\ref{stdbase3}; we have
$\omega(e_i)=f_i$, $\omega(f_i)=e_i$ and $\omega(h_i)=-h_i$ for
$i\in I$. Applying Lemma~\ref{genchev0} with $\theta=\omega$, we obtain
\begin{align*}
\omega\,\circ &\, x_i(t)\circ \omega^{-1}=\omega\circ \exp\bigl(t\,
\ad_\cL(e_i) \bigr)\circ \omega^{-1}\\ &=\exp\bigl(t\,
\ad_\cL(\omega(e_i))\bigr)=\exp\bigl(t\,\ad_\cL(f_i)\bigr)=y_i(t)
\end{align*}
for all $t\in \C$. We wish to extend this formula to any field $K$. For 
this purpose, we first consider the action of $\omega$ on $\bB$. Since 
$h_j^+=-\epsilon(j)h_j$ for $j\in I$, we have $\omega(h_j^+)=-h_j^+$. 
By Theorem~\ref{xcaomega}, we also have $\omega(\be_\alpha^+)= 
-\be_{-\alpha}^+$ for $\alpha\in \Phi$. We use these formulae to 
define a linear map $\bar{\omega} \colon \bar{\cL}\rightarrow 
\bar{\cL}$; explicitly, we set:
\[\bar{\omega}(\bar{h}_j^+):=-\bar{h}_j^+ \quad (j \in I) \quad
\mbox{and}\quad \bar{\omega}(\bar{\be}_\alpha^+):=-\bar{\be}_{-\alpha}^+
\quad (\alpha\in\Phi).\]
Note that we still have $\bar{\omega}^2=\id_{\bar{\fg}}$.
With this definition, we claim that 
\[ \bar{\omega}\circ \bar{x}_i(\zeta)\circ \bar{\omega}^{-1}=
\bar{y}_i(\zeta)\qquad \mbox{for all $\zeta\in K$}.\]
To prove this, we follow the argument in Lemma~\ref{fromctok1}.
Let $\Omega\in M_{N'}(\C)$ be the matrix of $\omega$ with respect to~$\bB$. 
The above formulae show that $\Omega$ only has entries $0$ and $-1$; 
we can simply regard $\Omega$ as a matrix in $M_{N'}(\Z[T])$. Then the
above formula over $\C$ implies that
\[ \Omega\cdot X_i(T)=Y_i(T)\cdot \Omega \qquad\mbox{(equality in 
$M_{N'}(\Z[T])$).}\]
Let $\bar{\Omega}\in M_{N'}(K)$ be the matrix of~$\bar{\omega}$. Now 
fix $\zeta\in K$ and consider the canonical ring homomorphism 
$\varphi_\zeta\colon \Z[T]\rightarrow K$ with $\varphi_\zeta(T)=\zeta$. 
Applying $\varphi_\zeta$ to the entries of $\Omega$, we obtain 
$\bar{\Omega}$. Hence, the above identity over 
$\Z[T]$ implies the identity $\bar{\Omega}\cdot \bar{X}_i(\zeta)=
\bar{Y}_i(\zeta)\cdot  \bar{\Omega}$ over~$K$, which means that
$\bar{\omega}\circ \bar{x}_i(\zeta)\circ \bar{\omega}^{-1}=
\bar{y}_i(\zeta)$, as desired. Hence, conjugation by~$\bar{\omega}$
inside $\GL(\bar{\cL})$ defines a group isomorphism
\[ G_K(\cL)\rightarrow G_K(\cL), \qquad 
\bar{x}_i(\zeta)\mapsto \bar{y}_i(\zeta),\quad 
\bar{y}_i(\zeta)\mapsto \bar{x}_i(\zeta).\]
\end{exmp}



\begin{rem} \label{dirprodg} Assume that the structure matrix $A=
(a_{ij})_{i,j\in I}$ of $(\cL,\cH)$ is decomposable. So there is
a partition $I=I_1\sqcup I_2$ such that $A$ has a block diagonal shape
\begin{center}
$A=\left(\begin{array}{c|c} A_1 & 0 \\\hline 0 & A_2\end{array}\right)$
\end{center}
where $A_1$ has rows and columns labelled by $I_1$, and $A_2$ has rows
and columns labelled by $I_2$. As discussed in Remark~\ref{cindec5}, we 
have $\cL=\cL_1\oplus \cL_2$, where $\cL_1$ and $\cL_2$ are subalgebras 
of Cartan--Killing type with structure matrices $A_1$ and $A_2$, 
respectively, and such that $[\cL_1,\cL_2]=\{0\}$. One immediately sees 
that $\bB=\bB_1 \cup\bB_2$, where $\bB_1$ is the canonical basis of $\cL_1$ 
(with respect to $\epsilon|_{I_1}$) and $\bB_2$ is the canonical basis of 
$\cL_2$ (with respect to $\epsilon|_{I_2}$). Let $N_1'=|\bB_1|$ and 
\[X_i^{(1)}(T),\; Y_i^{(1)}(T)\in M_{N_1'}(\Z[T]), \qquad i \in I_1,\]
be the matrices defined in Remark~\ref{luform2} with respect to $\cL_1$ 
and the basis $\bB_1$; similarly, let $N_2':=|\bB_2|$ and 
\[X_j^{(2)}(T),\; Y_j^{(2)}(T)\in M_{N_2'}(\Z[T]), \qquad j \in I_2,\]
be the matrices defined with respect to $\cL_2$ and the basis $\bB_2$. We 
also have matrices $X_i(T),Y_i(T)\in M_{N'}(\Z[T])$ and $X_j(T),Y_j(T)\in 
M_{N'}(\Z[T])$ defined with respect to $\cL$ and $\bB$. Then the formulae 
in Theorem~\ref{luform} show that
\begin{align*} 
X_i(T)&=\left(\begin{array}{c|c} X_i^{(1)}(T) & 0 \\\hline 0 & I_{N_2'}
\end{array}\right),  \qquad X_j(T)=\left(\begin{array}{c|c} 
I_{N_1'} & 0 \\\hline 0 & X_i^{(2)}(T) \end{array}\right),\\
Y_i(T)&=\left(\begin{array}{c|c} Y_i^{(1)}(T) & 0 \\\hline 0 & I_{N_2'}
\end{array}\right),  \qquad \;Y_j(T)=\left(\begin{array}{c|c} 
I_{N_1'} & 0 \\\hline 0 & Y_i^{(2)}(T) \end{array}\right).
\end{align*}
Since this holds at the polynomial level, we obtain analogous block
diagonal shapes for $\bar{X}_i(\zeta)$, $\bar{X}_j(\zeta)$, $\bar{Y}_i
(\zeta)$, $\bar{Y}_j(\zeta)$, where $i\in I_1$, $j\in I_2$ and $\zeta 
\in K$ for any field~$K$. Consequently, if we set 
\begin{align*}
G_{K,1}(\fg)&:=\langle \bar{x}_i(\zeta),\bar{y}_i(\zeta)\mid i\in I_1,
\zeta\in K \rangle\subseteq G_K(\fg),\\
G_{K,2}(\fg)&:=\langle \bar{x}_i(\zeta),\bar{y}_i(\zeta)\mid i\in I_2,
\zeta\in K \rangle\subseteq G_K(\fg),
\end{align*}
then $G_K(\fg)=G_{K,1}(\fg)\times G_{K,2}(\fg)$ (direct product of groups). 
Furthermore, we have group isomorphisms
\[G_{K,1}(\fg)\cong G_K(\fg_1)\qquad \mbox{and}\qquad G_{K,2}(\fg)\cong 
G_K(\fg_2).\]
The first one is obtained by sending the matrix of an element of the 
group $G_K(\fg_1)$ to a block diagonal matrix as above, where the second 
diagonal block is the identity matrix; analogously for~$G_K(\fg_2)$.
\end{rem}

\begin{rem} \label{limits} In order to establish properties of 
$G_K(\fg)$ we used (and will use) several times the technique of first 
establishing an analogue of that property over $K=\C$, then to lift this 
to a ``polynomial level'' and finally to pass from there to an arbitrary 
field~$K$. We just want to mention here that there are properties of 
$G_K(\fg)$ which do not seem to be accessible via that technique. For 
example, we shall see later that the center $Z(G_K(\fg))$ is always trivial. 
But it is not at all clear how to prove this using the above technique.
\end{rem}


%

As already noted, if $K$ is a finite field, then $G_K(\fg)$ is finite. 
Even if $K$ is very small, then these groups may simply become 
enormous. For example, if $|K|=2$ and $\cL$ is of type~$E_8$, 
then $G_K(\fg)$ has
\[2^{120}{\cdot}3^{13}{\cdot}5^{5}{\cdot}7^{4}{\cdot}11^{2}{\cdot} 
13^{2}{\cdot}17^{2}{\cdot}19{\cdot}31^{2}{\cdot} 41{\cdot} 43{\cdot}
 73{\cdot}  127{\cdot} 151{\cdot}  241{\cdot}  331\]
($\approx 3,38 \times 10^{74}$) elements (see the general order formula
that we will prove later.) Nevertheless, we 
shall see that $G_K(\fg)$ has a very user-friendly internal structure, and
there are highly convenient ways how to work with the elements. Many 
manipulations with $G_K(\fg)$ can be performed in a uniform way, almost 
regardless of the specific base field~$K$. 

\section{A first example: Groups of type $A_{n-1}$} \label{sec3fin}

In this short section we look in more detail at the example where 
$\cL=\slm_n(\C)$, $n\geq 2$. We would like to identify the corresponding 
Chevalley group $G_K(\fg)$ (over a fixed field~$K$) with a ``known'' group. 

We recall some notation from Example~\ref{cartsln}. Let $\cH\subseteq \cL$ 
be the abelian subalgebra of diagonal matrices. For $1\leq i,j\leq n$ let 
$E_{ij}$ be the $n \times n$-matrix with~$1$ as its $(i,j)$-entry and zeroes 
elsewhere. Let $e_i:=E_{i,i+1}$ and $f_i:=E_{i+1,i}$ for $1\leq i \leq n-1$. 
Then $\{e_i,f_i\mid 1\leq i\leq n-1\}$ are Chevalley generators of~$\cL$; 
furthermore, $h_i=[e_i,f_i]=E_{ii}-E_{i+1,i+1}$. Also recall from 
Example~\ref{cartsln} that
\[ \Phi=\{\varepsilon_i-\varepsilon_j \mid 1\leq i,j\leq n,i\neq j\},
\qquad \cL_{\varepsilon_i-\varepsilon_j}=\langle E_{ij}\rangle_\C.\]
We set $\be_{\alpha}^+:=(-1)^jE_{ij}$ for $\alpha=\varepsilon_i-
\varepsilon_j$, $i\neq j$. By Exercise~\ref{lucanAn}, the collection 
$\{\be_\alpha^+\mid \alpha\in\Phi\}$ satisfies the conditions in 
Corollary~\ref{canbash}. In particular, $\be_{\alpha_i}^+=-(-1)^ie_i$ and
$\be_{-\alpha_i}^+=(-1)^i f_i$ for $1\leq i\leq n{-}1$; furthermore, $h_i^+
=[e_i,\be_{-\alpha_i}^+]=(-1)^ih_i$. Thus, all elements in Lusztig's 
canonical basis $\bB$ of $\fg$ are matrices with entries in~$\Z$. 

Note that all matrices $E_{ij}$ with $i\neq j$ are nilpotent. We now 
require the following result which will also be extremely useful later 
on. It is called ``\nm{Transfer Lemma}'' because it provides a tool to 
``transfer'' results about the adjoint representation of a Lie algebra 
to an arbitrary representation.

\begin{lem}[Transfer Lemma] \label{superexp} Let $\fg$ be an arbitrary 
Lie algebra over~$\C$ and $V$ be a $\cL$-module. Let $\rho \colon \fg
\rightarrow \gl(V)$ be the corresponding representation. Let $x\in\cL$ be 
such that the linear maps $\ad_\cL(x)\colon \cL\rightarrow \cL$ and $\rho(x) 
\colon V\rightarrow V$ are nilpotent. Then, for any $y\in \fg$, we have 
\[ \rho\bigl(\exp\bigl(\ad_\cL(x)\bigr)(y)\bigr)=\exp(\rho(x))\circ
\rho(y)\circ \exp(\rho(x))^{-1}.\]
\end{lem}

\begin{proof} Consider the associative algebra $A=\End(V)$ (with product
given by the composition ``$\circ$'' of maps). We write $\tilde{y}=
\rho(y)\in A$ for any $y\in \cL$. A simple induction on $n$ shows that 
\begin{equation*}
\rho\bigl(\ad_\cL(x)^n(y)\bigr)=\ad_A(\tilde{x})^n(\tilde{y})\qquad
\mbox{for all~$n\geq 0$}.\tag{$*_1$}
\end{equation*}
Now, as in Remark~\ref{adjoint0}, we have $\ad_A(\tilde{x})=L_{\tilde{x}} 
-R_{\tilde{x}}$, where $L_{\tilde{x}}$ and $R_{\tilde{x}}$ are the 
endomorphisms of $A$ given by left and right composition with~$\tilde{x}$, 
respectively. These endomorphisms commute with each other, since $A$ is 
associative. Hence, we obtain that 
\[\frac{1}{n!}\ad_A(\tilde{x})^n= \frac{1}{n!}
\bigl(L_{\tilde{x}} -R_{\tilde{x}}\bigr)^n=\sum_{\atop{i,j\geq 0}{i+j=n}}
(-1)^j\; \frac{L_{\tilde{x}}^i}{i!} \circ \frac{R_{\tilde{x}}^j}{j!}\]
(where, here, ``$\circ$'' is the composition in $\End(A)$). We apply the
above endomorphism of $A$ to $\tilde{y}$; using also ($*_1$) yields that 
\begin{equation*}
\frac{1}{n!}\rho\bigl(\ad_\cL(x)^n(y)\bigr)=
\frac{1}{n!}\ad_A(\tilde{x})^n(\tilde{y})=\sum_{\atop{i,j\geq 0}{i+j=n}}
(-1)^j\; \frac{\tilde{x}^i}{i!} \circ \tilde{y} \circ 
\frac{\tilde{x}^j}{j!}\tag{$*_2$}
\end{equation*}
(where, now, ``$\circ$'' is the composition in $A$). By assumption,
both $\ad_\cL(x)$ and $\tilde{x}$ are nilpotent. So we can now sum ($*_2$) 
over all $n\geq 0$ (there will only be finitely many non-zero terms). This 
yields that
\begin{align*}
\rho\Bigl(\sum_{n\geq 0} &\frac{1}{n!}\ad_\cL(x)^n(y)\Bigr)=
\sum_{n\geq 0} \sum_{\atop{i,j\geq 0}{i+j=n}} (-1)^j\; 
\frac{\tilde{x}^i}{i!} \circ \tilde{y} \circ \frac{\tilde{x}^j}{j!}\\
& =\sum_{i,j\geq 0} (-1)^j\; \frac{\tilde{x}^i}{i!} \circ \tilde{y}
\circ \frac{\tilde{x}^j}{j!}=\Bigl(\sum_{i\geq 0} \frac{\tilde{x}^i}{i!} 
\Bigr)\circ \tilde{y} \circ\Bigl(\sum_{j\geq 0} \frac{(-\tilde{x})^j}{j!}
\Bigr).
\end{align*}
The right hand side equals $\exp(\tilde{x})\circ \tilde{y}\circ 
\exp(\tilde{x})^{-1}$, and the left hand side equals 
$\rho\bigl(\exp\bigl(\ad_\cL(x)\bigr)(y)\bigr)$, as desired.
\end{proof}

Let us return to $\fg=\slm_n(\C)$. Then $\C^n$ is naturally a $\fg$-module; 
the corresponding representation is given by the inclusion $\fg
\hookrightarrow \gl_n(\C)$, where we identify $\gl_n(\C)=\gl(\C^n)$. Let 
$i\in \{1,\ldots,n-1\}$ and consider $e_i=E_{i,i+1}$. (The following 
argument will be similar for $f_i$.) We already noted that $e_i$ is a 
nilpotent matrix; in fact, $e_i^2=0_{n\times n}$ and so
\[ \exp(te_i)=I_n+te_i\qquad \mbox{for $t\in \C$}.\] 
Note that then we also have $\exp(te_i)^{-1}=I_n-te_i$ for $t\in 
\C$. Hence, the above Transfer Lemma yields the following identity:
\[x_i(t)(y)=(I_n+te_i)\cdot y\cdot (I_n-te_i) \quad \mbox{for 
$t\in \C$ and $y\in \slm_n(\C)$}.\]
Now, as in the previous section, there are well-defined polynomials 
$f_{b,b'}^i \in \Z[T]$ (where $T$ is an indeterminate and $b,b'\in \bB$) 
such that
\[ x_i(t)(b)=\sum_{b'\in \bB} f_{b,b'}^i(t)b' \qquad\mbox{for 
all $b\in \bB$ and $t\in \C$}.\] 
Hence, for each fixed $b\in \bB\subseteq M_n(\C)$, we obtain the 
following identity of matrices in $M_n(\C)$:
\[\sum_{b'\in \bB} f_{b,b'}^i(t)b'=(I_n+te_i)\cdot b\cdot (I_n-te_i) 
\qquad\mbox{where $t\in \C$}.\]
Since this holds for all $t\in \C$, we also obtain an identity of matrices 
in $M_n(\C[T])$ where $T$ is an indeterminate:
\begin{equation*}
\sum_{b'\in \bB} f_{b,b'}^i \,b'=(I_n+Te_i)\cdot b\cdot (I_n-Te_i).
\tag{$\dagger$}
\end{equation*}
Actually, by the above description of $\bB$, we have $b,b' \in M_n(\Z)$ 
for all $b,b'\in \bB$. So ($\dagger$) is an identity of matrices in 
$M_n(\Z[T])$.

Now let $K$ be any field. Following the construction in the previous 
section, we need to consider a vector space $\bar{\cL}$ over $K$ with a
basis indexed by the canonical basis~$\bB$ of $\cL$. Concretely, we may 
take $\bar{\cL}:=\slm_n(K)$ with basis 
\[ \bar{\bB}=\{\bar{h}_i^+\mid 1\leq j \leq n-1\}\cup\{\bar{\be}_\alpha^+
\mid \alpha\in \Phi\},\]
where $\bar{e}_i,\bar{f}_i,\bar{h}_i^+\in \slm_n(K)$ and 
$\bar{\be}_\alpha^+\in \slm_n(K)$ are defined exactly as above, using 
the matrices $E_{ij}\in M_n(K)$. We now define
\[ x_i^*(\zeta):=I_n+\zeta \bar{e}_i \qquad \mbox{for $\zeta\in K$}.\]
We still have $\bar{e}_i^2=0_{n\times n}$ and so $x_i^*(\zeta)^{-1}=
I_n-\zeta\bar{e}_i$. Applying the ring homomorphism $\Z[T]\rightarrow K$, 
$T\mapsto \zeta$, to the identity ($\dagger$), we obtain an analogous
identity over~$K$ for any $b\in \bB$:
\[\sum_{b'\in \bB} f_{b,b'}^i(\zeta)\,\bar{b}'=(I_n+\zeta\bar{e}_i) 
\cdot \bar{b}\cdot (I_n-\zeta\bar{e}_i) \qquad \mbox{for any $\zeta\in K$}.\]
Now note that the left hand side just equals $\bar{x}_i(\zeta)
(\bar{b})$. Using also the above definition of $x_i^*(\zeta)$, we 
finally obtain the following identity:
\begin{equation*}
\bar{x}_i(\zeta)(\bar{b})=x_i^*(\zeta)\cdot \bar{b}\cdot x_i^*(\zeta)^{-1} 
\qquad\mbox{for all $\zeta \in K$}.
\tag{$\dagger_K$}
\end{equation*}
A completely analogous argument (using $f_i$ instead of $e_i$) 
shows that 
\[ \bar{y}_i(\zeta)(\bar{b})=y_i^*(\zeta)\cdot \bar{b}\cdot y_i^*
(\zeta)^{-1} \qquad\mbox{for all $\zeta \in K$},\]
where we set $y_i^*(\zeta):=I_n+\zeta \bar{f}_i$. After these 
preparations, we can now prove the following identification result.

\begin{prop}[Ree \protect{\cite{ree1}}] \label{propree1} If $\cL=\slm_n
(\C)$ and $K$ is any field, then the Chevalley group $G_K(\fg)\subseteq\GL
(\bar{\cL})$ (as in Definition~\ref{deflu}) is isomorphic to $\SL_n(K)/Z$, 
where $Z=\{\zeta \bar{I}_n\mid \zeta\in K^\times,\zeta^n=1\}$. 
\end{prop}

\begin{proof} As above, let $\bar{\cL}=\slm_n(K)$. We also set 
\[ G^*:=\bigl\langle x_i^*(\zeta),y_i^*(\zeta)\mid 1\leq i \leq n-1, 
\zeta \in K\bigr\rangle\subseteq \SL_n(K).\]
Then $G^*$ acts on $\bar{\cL}$ by conjugation. Thus, for $g\in G^*$ we 
obtain a map 
\[ \gamma_g\colon \bar{\cL}\rightarrow \bar{\cL}, \qquad y \mapsto 
g\cdot y\cdot g^{-1}.\]
Then $\gamma_g\in \GL(\bar{\cL})$ and $\gamma\colon G^*\rightarrow
\GL(\bar{\cL})$, $g\mapsto\gamma_g$, is a group homomorphism. By 
($\dagger_K)$, we have $\gamma_g=\bar{x}_i(\zeta)$ for $g=x_i^*(\zeta)$ 
and, analogously, $\gamma_g=\bar{y}_i(\zeta)$ for $g=y_i^*(\zeta)$. Hence, 
the image of $\gamma$ equals $G_K(\fg)$. By Appendix~\ref{asecSLn}, we 
have in fact $\SL_n(K)=G^*$. It remains to show that $\ker(\gamma)=Z$. 
So let $g\in G^*$ be such that $\gamma_g=\id_{\bar{\cL}}$. Then $g\cdot 
y=y\cdot g$ for all $y\in \bar{\cL}$; it is a standard fact from Linear 
Algebra that then $g=\zeta \bar{I}_n$ for some $\zeta\in K$. Since 
$\det(g)=1$, we must have $\zeta^n=1$ and so $g\in Z$. Conversely, 
it is clear that $Z\subseteq \ker(\gamma)$. 
\end{proof}

\begin{rem} \label{remree1} (a) Let $n\geq 2$. It is known that 
$\SL_n(K)/Z$ is simple, unless $n=2$ and $K$ has 
$2$ or $3$ elements; see, e.g., \cite[Theorem~1.13]{Grove2}. We will 
also see this later as a special case of a more general result.

(b) The Chevalley groups associated with the classical Lie algebras
$\gom_n(Q_n,\C)$ can be identified with symplectic or orthogonal groups
in a similar way; see Carter \cite[Chap.~11]{Ca1} and Ree \cite{ree1} 
for further details. We will come back to this later, once we have 
introduced Chevalley groups of non-adjoint type.
\end{rem}

\begin{xca} \label{reesp4} The purpose of this exercise is to give at 
least one example showing that the above procedure also works 
for the classical Lie algebras introduced in Section~\ref{sec05}.
Let $\cL=\gom_4(Q_4,\C)$, where 
\[\renewcommand{\arraystretch}{0.8} Q_4=\left(
\begin{array}{r@{\hspace{5pt}}r@{\hspace{10pt}}r@{\hspace{10pt}}r} 
0 & 0 & 0 & 1 \\ 0 & 0 & 1 & 0 \\ 0 & -1 & 0 & 0 \\ -1 & 0 & 0 & 0
\end{array}\right), \qquad Q_4^{\text{tr}}=-Q_4.\]
Let $I=\{1,2\}$. We have $\Phi=\{\pm \alpha_1,\pm\alpha_2,
\pm(\alpha_1+\alpha_2),\pm (\alpha_1+2\alpha_2)\}$. 
Chevalley generators for~$\cL$ are given as follows:
\begin{alignat*}{2}
e_1&=-\textstyle{\frac{1}{2}}A_{2,3}, \quad f_1 & 
=\textstyle{\frac{1}{2}}A_{3,2}, \quad h_1&=[e_1,f_1]= 
\mbox{diag}(0,1,-1,0);\\e_2&=-A_{1,2}, \quad f_2&=-A_{2,1},\quad  
h_2&=[e_2,f_2]=\mbox{diag}(1,-1,1,-1).
\end{alignat*}
(See the proof of Proposition~\ref{CKbcd}.) We have the relations $[h_1,
e_2]=-e_2$ and $[h_2,e_1]=-2e_1$; see the structure matrix in 
Table~\ref{cartanBCD} (p.~\pageref{cartanBCD}). 

Let $\mbox{Sp}_4(K):=\{T\in M_4(K)\mid T^{\text{tr}}Q_4T =Q_4\}$. One
easily sees that $\mbox{Sp}_4(K)$ is a subgroup of $\GL_4(K)$; it is 
called the $4$-dimensional \nm{symplectic group} over $K$. Now proceed 
as follows.

\noindent (a) Let $\epsilon\colon I\rightarrow \{\pm 1\}$ be given by
$\epsilon(1)=1$ and $\epsilon(2)=-1$, as in Table~\ref{Mdynkineps} 
(p.~\pageref{Mdynkineps}). Starting with $\be_{\alpha_i}^+=\epsilon(i)e_i$
and $\be_{-\alpha_i}^+=-\epsilon(i)f_i$ for $i=1,2$, determine all 
the elements of the canonical basis $\bB$, explicitly as matrices in 
$\cL$; observe that all those matrices have entries in~$\Z$.

\noindent (b) Let $K$ be any field with $\mbox{char}(K)\neq 2$, and set
$\bar{\cL}:=\gom_4(Q_4,K)$. Check that Proposition~\ref{classic4}(b) 
also holds over~$K$ instead of~$\C$. 

\noindent (c) Define $\bar{\bB}\subseteq \bar{\cL}$ by taking analogues 
of the matrices in (a) over~$K$; check that $\bar{\bB}$ is a basis of
$\bar{\cL}$. For $i\in I$ and $\zeta\in K$, define $x_i^*(\zeta):=
\exp(\zeta \bar{e}_i)$ and $y_i^*(\zeta):=\exp(\zeta \bar{f}_i)$, as 
above. Check that $e_i^2=f_i^2=0_{4\times 4}$ for $i=1,2$ and that the 
analogous versions of ($\dagger_K$) also hold here. 

\noindent (d) Show that $x_i^*(\zeta),y_i^*(\zeta)\in \mbox{Sp}_4(K)$
for $i=1,2$ and $\zeta\in K$. Let
\[G^*:=\big\langle x_i^*(\zeta), y_i^*(\zeta)\mid \zeta\in K,i=1,2
\big\rangle \subseteq \mbox{Sp}_4(K).\] 
Analogously to Proposition~\ref{propree1}, show that $G_K(\fg)\cong 
G^*/Z$, where $Z=\{\pm \bar{I}_4\}$. Finally (and this is probably the 
difficult part) show that $G^*=\mbox{Sp}_4(K)$. (For help and further 
references, see Ree \cite{ree1}.)
\end{xca}

\begin{rem} \label{remgenclassic} Both in the proof of 
Proposition~\ref{propree1} and in Example~\ref{reesp4}, a crucial step 
consists of showing that the subgroup~$G^*$ defined there equals $\SL_n(K)$ 
or $\mbox{Sp}_4(K)$, respectively. If $K$ is algebraically closed, then 
this can be alternatively deduced from general results about algebraic 
groups. Indeed, in Chapter~\ref{chap4}, $G^*$ will be seen to be a Chevalley 
group of ``non-adjoint'' type. If $K$ is algebraically closed, then any 
such group will be shown to be a connected linear algebraic group with a 
$(B,N)$-pair, such that $\dim G^*= |I|+|\Phi|=\dim \fg$. (For all this, 
see a later chapter.) On the other hand, $\SL_n(K)$ and 
$\mbox{Sp}_4(K)$ are known to be connected linear algebraic groups with 
dimension equal to $\dim \fg$; see, e.g., \cite[Example~1.3.10, 
Corollary~1.5.14 and Summary~1.7.9]{my0}. Hence, $G^*$ must be equal to 
$\SL_n(K)$ or $\mbox{Sp}_4(K)$, respectively. Similar arguments apply to 
the generation of other classical groups.
\end{rem}

\section{The elements $\bar{n}_i(\xi)$ and $\bar{h}_i(\xi)$} 
\label{secweylG}

Now let us return to the general situation, where $G_K(\fg)$ is the 
Chevalley group (over a field $K$) associated with a Lie algebra $\cL$ 
of Cartan--Killing type. Since $\fg$ and $K$ will be fixed throughout 
the following discussion, we shall just write $\bar{G}=G_K(\fg)$. Now 
$\bar{G}$ is defined in terms of the generators $\bar{x}_i(\zeta)$ 
and $\bar{y}_i(\zeta)$. However, just knowing generators of a group does 
not tell us much about the structure of that group. (For example, it is
known that every finite simple group is generated by two of its elements; 
see Aschbacher--Guralnick~\cite{AschG}.) So one usually needs to produce 
further, more ``helpful'' elements where ``helpful'' means, for example, 
that they generate subgroups with interesting properties. We now make some
first steps in this direction, which will at least be sufficient to show 
that, if $|K| \geq 4$, then $\bar{G}$ is equal to its own commutator 
subgroup\footnote{If $\Gamma$ is any group, then 
$[g,h]:=g^{-1}h^{-1}gh$ is called the commutator of $g,h\in \Gamma$; if 
$\Gamma_1,\Gamma_2 \subseteq \Gamma$ are subgroups, we set $[\Gamma_1, 
\Gamma_2]:=\langle [g_1,g_2] \mid g_1\in \Gamma_1,g_2\in \Gamma_2
\rangle \subseteq\Gamma$.  Then $[\Gamma,\Gamma]$ is called the 
\nm{commutator subgroup} of $\Gamma$; this is the smallest normal 
subgroup $\Gamma' \subseteq \Gamma$ such that $\Gamma/\Gamma'$ is abelian. 
These are standard notions in the theory of groups.}. It will turn out later 
that this is a big step towards proving that $\bar{G}$ is even simple (if 
$|K|\geq 4$). Furthermore, we will exhibit a diagonalisable abelian
subgroup $\bar{H}\subseteq \bar{G}$ which, in the context of the
theory of algebraic groups (to be discussed in a later chapter) will
play the role of a ``maximal torus''.

As far as new, ``helpful'' elements are concerned, the next candidates to
consider are the Lie algebra automorphisms
\[ n_i(t):=x_i(t)\circ y_i(-t^{-1})\circ x_i(t) \in\mbox{Aut}(\cL)
\qquad (i \in I,\;t \in\C^\times)\]
introduced in Lemma~\ref{weyl5} (over $\C$). This leads us to define 
\[ \bar{n}_i(\xi):=\bar{x}_i(\xi)\bar{y}_i(-\xi^{-1})\bar{x}_i(\xi) 
\in\bar{G} \qquad (i \in I,\;\xi\in K^\times).\]
Here are first properties of these elements.

\begin{lem} \label{barni1} Let $i\in I$ and $\xi\in K^\times$. Then 
$\bar{n}_i(\xi)^{-1}=\bar{n}_i(-\xi)$. Furthermore, for $\zeta\in K$, we 
have 
\begin{align*}
\bar{n}_i(\xi)\bar{x}_i(\zeta)\bar{n}_i(\xi)^{-1}&=\bar{y}_i(-\zeta
\xi^{-2}),\\ \bar{n}_i(\xi)\bar{y}_i(\zeta)\bar{n}_i(\xi)^{-1}& =
\bar{x}_i (-\zeta\xi^2).
\end{align*}
Consequently, we also have $\bar{n}_i(\xi)=\bar{y}_i(-\xi^{-1}) 
\bar{x}_i(\xi)\bar{y}_i(-\xi^{-1})$.
\end{lem}

\begin{proof} By Lemma~\ref{fromctok1}, we have $\bar{x}_i(\xi)^{-1}=
\bar{x}_i(-\xi)$ and $\bar{y}_i(\xi)^{-1}=\bar{y}_i(-\xi)$. Hence,
we obtain
\[\bar{n}_i(\xi)^{-1}=\bigl(\bar{x}_i(\xi)\bar{y}_i(-\xi^{-1})
\bar{x}_i(\xi) \bigr)^{-1}=\bar{x}_i(-\xi)\bar{y}_i(\xi^{-1})
\bar{x}_i(-\xi)\]
where the right hand side equals $\bar{n}_i(-\xi)$, as claimed. In order 
to prove the remaining statements, we first work over $K=\C$. Let $t\in \C$ 
and $u\in \C^\times$. Let $\theta=n_i(u)\in \Aut(\fg)$. Then we obtain
\begin{align*}
n_i&(u)x_i(t)n_i(u)^{-1}=n_i(u)\circ \exp\bigl(t\,\ad_\fg(e_i)\bigr)
\circ n_i(u)^{-1}\\ &=\exp\bigl(t\,\ad_\fg(n_i(u)(e_i))\bigr)= 
\exp\bigl(-tu^{-2}\,\ad_\fg(f_i)\bigr)=y_i(-tu^{-2}),
\end{align*}
where the second equality holds by Lemma~\ref{genchev0} and the third
equality holds by Exercise~\ref{weyl5a}. Now let again $X_i(t)$ and 
$Y_i(t)$ be the matrices of $x_i(t)$ and $y_i(t)$ with respect to~$\bB$, 
respectively. Then $N_i(u):=X_i(u)\cdot Y_i(-u^{-1}) \cdot X_i(u)$ is the 
matrix of $n_i(u)$ with respect to $\bB$. Thus, we have 
\[ N_i(u)\cdot X_i(t)=Y_i(-tu^{-2})\cdot N_i(u) \quad \mbox{for all $t\in
\C$ and $u\in \C^\times$}.\]
We will now work over the ring $\Z[T,U^{\pm 1}]$, where $T,U$ are two 
commuting indeterminates. Let $X_i(T)$ and $Y_i(T)$ be the matrices 
defined in Remark~\ref{luform2} (with entries in $\Z[T]$). Let 
\[ N_i(U):=X_i(U)\cdot Y_i(-U^{-1}) \cdot X_i(U).\]
Since the above identity holds for all $t\in \C$ and $u \in \C^\times$, we 
also have an identity of matrices at the ``polynomial level'':
\[ N_i(U)\cdot X_i(T)=Y_i(-TU^{-2})\cdot N_i(U). \]
Given $\zeta\in K$ and $\xi\in K^\times$, we consider the 
ring homomorphism  $\Z[T,U^{\pm 1}]\rightarrow K$ which sends $T$ to
$\zeta$ and $U$ to~$\xi$. Applying that homomorphism to the above 
identity, we obtain 
\[ \bar{N}_i(\xi)\cdot \bar{X}_i(\zeta)=\bar{Y}_i(-\zeta\xi^{-2})\cdot 
\bar{N}_i(\xi),\]
where $\bar{N}_i(\xi)$ is the matrix of $\bar{n}_i(\xi)$ with respect to
$\bar{\bB}$; furthermore, $\bar{X}_i(\zeta)$ and $\bar{Y}_i(-\zeta\xi^{-2})$ 
are the matrices of $\bar{x}_i(\zeta)$ and $\bar{y}_i(-\zeta\xi^{-2})$
with respect to $\bar{\bB}$, respectively. This implies the identity
$\bar{n}_i(\xi)\bar{x}_i(\zeta)\bar{n}_i(\xi)^{-1}=\bar{y}_i(-\xi^{-2}
\zeta)$. The proof of the second identity is entirely analogous.
To prove the last identity, we write 
\[\bar{n}_i(\xi)=\bar{n}_i(\xi)\, \bar{n}_i(\xi)\,\bar{n}_i(\xi)^{-1}=
\bar{n}_i(\xi)\,\bigl(\bar{x}_i(\xi) \bar{y}_i(-\xi^{-1})\bar{x}_i(\xi)
\bigr)\, \bar{n}_i(\xi)^{-1}\]
and then use the above formulae. 
\end{proof}

\begin{cor} \label{barni2} We set $\bar{h}_i(\xi):=\bar{n}_i(\xi)
\bar{n}_i(-1) \in\bar{G}$ for $i\in I$ and $\xi\in K^\times$.
Then, for any $\zeta\in K$, we have 
\begin{align*}
\bar{h}_i(\xi)\bar{x}_i(\zeta)\bar{h}_i(\xi)^{-1}&=\bar{x}_i(\zeta\xi^{2}),
\\ \bar{h}_i(\xi)\bar{y}_i(\zeta)\bar{h}_i(\xi)^{-1}&=\bar{y}_i(\zeta
\xi^{-2}).
\end{align*}
\end{cor}

\begin{proof} A first application of Lemma~\ref{barni1} yields that
\[\bar{n}_i(-1)\bar{x}_i(\zeta) \bar{n}_i(-1)^{-1}=\bar{y}_i(-\zeta).\] 
Then a second application yields the desired formula. The argument
for $\bar{y}_i(\zeta)$ is completely analogous. 
\end{proof}



\begin{exmp} \label{expHNsl2} Let $\cL=\slm_2(\C)$ and $K$ be any field.
In Example~\ref{luform2aa}, we determined the matrices of $\bar{x}_1(\zeta)$ 
and $\bar{y}_1(\zeta)$ for $\zeta\in K$:
\[\renewcommand{\arraystretch}{0.8} 
\bar{X}_1(\zeta)=\left(\begin{array}{ccc} 1 & 2\zeta & \zeta^2 \\ 
0 & 1 & \zeta \\ 0 & 0 & 1 \end{array}\right) \quad\mbox{and}\quad 
\bar{Y}_1(\zeta)=\left(\begin{array}{ccc} 1 & 0 & 0 \\ 
\zeta & 1 & 0 \\ \zeta^2 & 2\zeta & 1 \end{array}\right).\]
Now consider $\bar{n}_1(\xi)$ and $\bar{h}_1(\xi)$ for $\xi\in K^\times$.
A straightfoward computation shows that the matrices of these elements are
given by
\[\renewcommand{\arraystretch}{0.8} 
\bar{N}_1(\xi)=\left(\begin{array}{c@{\hspace{7pt}}c@{\hspace{9pt}}c} 
0 & 0 & \xi^2 \\ 0 & -1 & 0 \\ \xi^{-2} & 0 & 0 \end{array}\right)
\quad \mbox{and}\quad \bar{H}_1(\xi)=
\left(\begin{array}{c@{\hspace{9pt}}c@{\hspace{9pt}}c} \xi^2 & 0 & 0 
\\ 0 & 1 & 0 \\ 0 & 0 & \xi^{-2} \end{array}\right).\]
(We just need to compute, first the product $\bar{X}_1(\xi)\cdot \bar{Y}_1
(-\xi^{-1}) \cdot\bar{X}_1(\xi)$, and then the product $\bar{N}_1(\xi)
\cdot \bar{N}_1(-1)$.) In particular, this shows that the maps
$\xi \mapsto \bar{n}_i(\xi)$ and $\xi\mapsto \bar{h}_i(\xi)$ need not
be injective, contrary to the maps $\zeta \mapsto \bar{x}_i(\zeta)$
and $\zeta \mapsto \bar{y}_i(\zeta)$; see Remark~\ref{xinontriv}.
\end{exmp}

In the above example, the element $\bar{n}_i(\xi)$ is represented by
a monomial matrix with respect to $\bar{\bB}$ and the element 
$\bar{h}_i(\xi)$ is represented by a diagonal matrix. We will now 
show that this holds in general. Let us agree to set, for {\itshape any} 
$i\in I$ and $\alpha\in \Phi$: 
\begin{align*}
p_{i,\alpha}&:=\max\{m\geq 0 \mid \alpha+m\alpha_i\in\Phi\},\\
q_{i,\alpha}&:=\max\{m\geq 0 \mid \alpha-m\alpha_i\in\Phi\}.
\end{align*}
For $\alpha\neq \pm \alpha_i$, this coincides with the earlier definitions. 
Furthermore, we obtain $q_{i,\alpha_i}=p_{i,-\alpha_i}=2$ and 
$p_{i,\alpha_i}=q_{i,-\alpha_i}=0$. Note that $q_{i,\alpha}-p_{i,\alpha}=
\langle \alpha_i^\vee,\alpha\rangle$ in all cases; see Lemma~\ref{astring2}.
Now we can state:
  
\begin{thm} \label{luform3} Let $i\in I$ and $\xi\in K^\times$. Then
we have 
\begin{alignat*}{2}
\bar{n}_i(\xi)(\bar{h}_j^+)&=\bar{h}_j^+-|a_{ji}|\bar{h}_i^+
&&\qquad \mbox{for all $j\in I$},\\
\bar{n}_i(\xi)(\bar{\be}_\alpha^+)&=(-1)^{q_{i,\alpha}}
\xi^{-\langle \alpha_i^\vee,\alpha\rangle}\bar{\be}_{s_i(\alpha)}^+
&&\qquad \mbox{for all $\alpha\in \Phi$}.
\end{alignat*}
\end{thm}

\begin{proof} First we work over $K=\C$. Let $t\in \C^\times$ and 
consider $n_i(t)\in\Aut(\fg)$. For $j\in I$ we have $h_j^+ =-\epsilon(j)
h_j$. Lemma~\ref{weyl5}(a) shows that 
\[ n_i(t)(h_j^+)=-\epsilon(j)h_j+\epsilon(j)\alpha_i(h_j)h_i=
h_j^+ +\epsilon(j)a_{ji}h_i.\]
If $i=j$, then $\epsilon(j)a_{ji}h_i=\epsilon(i)a_{ii}h_i=-2h_i^+=
-|a_{ii}|h_i^+$. Now let $i\neq j$ and $a_{ij}\neq 0$. Then $a_{ij}<0$ 
and $\epsilon(i)= -\epsilon(j)$ (see Remark~\ref{epscanr}). Hence, 
$\epsilon(j)a_{ij}h_i=-\epsilon(i)a_{ji}h_i=-|a_{ji}|h_i^+$. This 
yields the desired formula for $n_i(t)(h_j^+)$. 

Next, consider $n_i(t)(\be_{\alpha}^+)$ for $\alpha\in \Phi$. If
$\alpha=\pm \alpha_i$, then $\be_{\alpha_i}^+=\epsilon(i)e_i$ and 
$\be_{-\alpha_i}^+=-\epsilon(i)f_i$. Hence, using Exercise~\ref{weyl5a},
we obtain 
\begin{align*}
n_i(t)(\be_{\alpha_i}^+)&=\epsilon(i)n_i(t)(e_i)=-\epsilon(i)t^{-2}
f_i=t^{-2}\be_{-\alpha_i}^+,\\
n_i(t)(\be_{-\alpha_i}^+)&=-\epsilon(i)n_i(t)(f_i)=\epsilon(i)t^2e_i=
t^2\be_{\alpha_i}^+,
\end{align*}
as required. Now let $\alpha\neq \pm \alpha_i$; we set $p:=p_{i,\alpha}$ 
and $q:=q_{i,\alpha}$. 

First assume that $p=0$. By Lemma~\ref{weyl5}(c), we have 
$n_i(t) (\cL_\alpha)=\cL_{s_i(\alpha)}$ and so $n_i(t)(\be_\alpha^+)=z
\be_{s_i(\alpha)}^+$ for some $0\neq z\in \C$. We must determine the 
scalar~$z$. Since $p=0$, we have $x_i(t)(\be_\alpha^+)=\be_\alpha^+$ 
and so Theorem~\ref{luform} yields that 
\begin{align*}
n_i(t)(\be_\alpha^+)&=x_i(t)\bigl(y_i(-t^{-1})(\be_\alpha^+)\bigr)=
\sum_{0\leq l\leq q} (-1)^l t^{-l}x_i(t)(\be_{\alpha-l\alpha_i}^+)\\ 
&=\sum_{0\leq l\leq q}\sum_{\;0\leq k\leq p_{i,\alpha-l\alpha_i}}
\binom{q_{i,\alpha-l\alpha_i} {+}k}{k}(-1)^lt^{k-l} 
\be_{\alpha-(l-k)\alpha_i}^+.
\end{align*}
By Remark~\ref{astring}(a), we have $\alpha(h_i)=q-p=q$ and so 
$s_i(\alpha)=\alpha-\alpha(h_i)\alpha_i=\alpha-q\alpha_i$. Hence, in 
order to determine the scalar~$z$, we must consider all terms in the 
above sums that correspond to indices $l,k$ such that $\alpha-(l-k)\alpha_i=
\alpha-q\alpha_i$, that is, $l-k=q$. Now note that $p_{i,\alpha-l 
\alpha_i}=p_{i,\alpha}+l=p+l=l$ and so $0\leq l-k\leq q$ for all terms
in the above sum. Hence, the condition $l-k=q$ only holds for the indices 
$l=q$ and $k=0$. Noting that $q_{i,\alpha-q\alpha_i}=0$, we obtain
\[ z=\binom{q_{i,\alpha-q\alpha_i}+0}{0}(-1)^qt^{-q}=
(-1)^qt^{-q}=(-1)^qt^{p-q},\]
as desired. Finally, assume that $p>0$ and set $\beta:=\alpha+p\alpha_i
\in \Phi$; then $p_{i,\beta}=0$ and $q_{i,\beta}=p+q$. We have
$s_i(\beta)=\beta-\beta(h_i)\alpha_i$ where $\beta(h_i)=q_{i,\beta}-
p_{i,\beta}=p+q$.
Hence, $s_i(\beta)=(\alpha+p\alpha_i)-(p+q)\alpha_i=\alpha-q\alpha_i$.
So, by the previous argument, we already know that 
\[ n_i(t)(\be_{\beta}^+)=(-1)^{p+q}t^{-p-q}\be_{s_i(\beta)}^+
=(-1)^{p+q}t^{-p-q}\be_{\alpha-q\alpha_i}^+.\]
Now, by Theorem~\ref{canbas}, we have 
$[f_i,\be_{\beta}^+]=(p_{i,\beta}+1)\be_{\beta-\alpha_i}^+=
\be_{\beta-\alpha_i}^+$. Since $p_{i,\beta-\alpha_i}=1$, we also have
\[[f_i,[f_i,\be_{\beta}^+]]=[f_i,\be_{\beta-\alpha_i}^+]
=(p_{i,\beta-\alpha_i}+1)\be_{\beta-2\alpha_i}^+=2\be_{\beta-2\alpha_i}^+\]
and so on. Thus, after $p$ repetitions of this computation, we obtain
\[ \underbrace{[f_i,[f_i,\ldots,[f_i}_{\text{$p$ times}},\be_\beta^+]
\ldots]] =p!\,\be_{\beta-p\alpha_i}^+=p!\,\be_\alpha^+.\]
Now $n_i(t)$ is a Lie algebra automorphism; furthermore, we already know
that $n_i(t)(\be_\beta^+)=(-1)^{p+q}t^{-p-q}\be_{\alpha-q\alpha_i}^+$ and 
that $n_i(t)(f_i)=-t^2e_i$ (see once more Exercise~\ref{weyl5a}). Hence,
applying $n_i(t)$ to the left hand side of the above identity yields that 
\begin{align*}
&\underbrace{[n_i(t)(f_i),[n_i(t)(f_i),\ldots,
[n_i(t)(f_i)}_{\text{$p$ times}},n_i(t)(\be_\beta^+)]\ldots]]\\
&\qquad =[-t^2e_i,[-t^2e_i,\ldots, [-t^2e_i,(-1)^{p+q}t^{-p-q}
\be_{\alpha-q\alpha_i}^+)]\ldots]]\\
&\qquad\qquad =(-1)^qt^{p-q}[e_i,[e_i,\ldots, [e_i,\be_{\alpha-
q\alpha_i}^+)]\ldots]].
\end{align*}
Now $q_{i,\alpha-q\alpha_i}=0$, $q_{i,\alpha-(q-1)\alpha_i}=1$,
$q_{i,\alpha-(q-2)\alpha_i}=2$, and so on. Hence, by Theorem~\ref{canbas}, 
we obtain that 
\[ \underbrace{[e_i,[e_i,\ldots,[e_i}_{\text{$p$ times}},\be_{\alpha-
q\alpha_i}^+] \ldots]] =p!\,\be_{(\alpha-q\alpha_i)+p\alpha_i}^+=
p!\,\be_{s_i(\alpha)}^+.\]
So we conclude that $n_i(t)(\be_\alpha^+)=(-1)^qt^{p-q} \be_{s_i(\alpha)}^+$,
as required.

Now let again $X_i(t)$ and $Y_i(t)$ be the matrices of $x_i(t)$ and
$y_i(t)$ with respect to $\bB$, respectively. Then $N_i(t):=X_i(t)
\cdot Y_i(-t^{-1})\cdot X_i(t)$ is the matrix of $n_i(t)$ with respect
to $\bB$, where $t\in \C^\times$. The rows and columns of $N_i(t)$ are 
indexed by~$\bB$. By the above computation, the $(b,b')$-entry of 
$N_i(t)$ (where $b,b'\in \bB$) is given as follows. 
\[ N_i(t)_{bb'}=\left\{\begin{array}{cl} -1 & \mbox{if $b=b'=h_i^+$},\\
1 & \mbox{if $b=b'=h_j^+$, $i\neq j$},\\ -|a_{ji}| & \mbox{if $b=h_i^+$, 
$b'=h_j^+$, $i\neq j$}\\ 
(-1)^{q_{i,\alpha}}t^{-\langle \alpha_i^\vee,\alpha\rangle} & \mbox{if 
$b=\be_{\alpha}^+$ and $b'=\be_{s_i(\alpha)}^+$},\\
0 & \mbox{otherwise}.\end{array}\right.\]
We will now work over the ring of Laurent polynomials $\Z[T,T^{-1}]$ in
an indeterminate~$T$. Let again $X_i(T)$ and $Y_i(T)$ be the matrices 
defined in Remark~\ref{luform2} (with entries in $\Z[T]$). Let 
\[ N_i(T):=X_i(T)\cdot Y_i (-T^{-1})\cdot X_i(T).\]
Then, upon substituting $T\mapsto t$ for any $t\in\C^\times$, we obtain 
$N_i(t)$. Since the above formulae for $N_i(t)_{bb'}$ hold for all $t\in 
\C^\times$, the matrix $N_i(T)$ will have analogous entries where~$t$ is 
replaced by~$T$. Now let $K$ be arbitrary, fix $\zeta\in K^\times$ and 
consider the canonical ring homomorphism $\Z[T,T^{-1}] \rightarrow K$ such 
that $T\mapsto \zeta$. Applying that homomorphism to $N_i(T)$, we obtain 
the matrix of $\bar{n}_i(\zeta)$ with respect to the basis~$\bar{\bB}$ 
of~$\bar{\fg}$. This yields the required formulae over~$K$.
\end{proof}


The exact formulae in the above theorem (especially for the signs occurring
there) will also play a role in the discussion of Chevalley groups of
non-adjoint type in Chapter~\ref{chap4} (see Proposition~\ref{Cree0}).
Working over $K=\C$, we have the following immediate application to the 
problem of computing the structure constants $N_{\alpha,\beta}^+$ of~$\fg$.

\begin{exmp} \label{exmcass} Let $K=\C$. By Theorem~\ref{luform3}, we 
have\footnote{There is a similar formula in \cite[Prop.~6.4.2]{Ca1}, 
with signs that are not explicitly determined there, but for which 
there are a certain number of rules in \cite[Prop.~6.4.3]{Ca1}.}
\[n_i(1)(\be_\alpha^+)=(-1)^{q_{i,\alpha}} \be_{s_i(\alpha)}^+ \qquad 
\mbox{for all $i\in I$ and $\alpha\in\Phi$}.\]
This can be used to compute the structure constants $N_{\alpha,\beta}^+$
recursively. Indeed, let $\alpha,\beta\in\Phi$ be such that $\alpha+
\beta\in\Phi$. If $\alpha\in\Phi^-$, then we apply the Chevalley 
involution $\omega\colon \cL\rightarrow \cL$ to $[\be_\alpha^+,
\be_\beta^+]=N_{\alpha,\beta}^+\be_{\alpha+\beta}^+$. Using 
Theorem~\ref{xcaomega}(a), we obtain $N_{\alpha,\beta}^+=
-N_{-\alpha,-\beta}^+$. Thus, it is sufficient to compute $N_{\alpha,
\beta}^+$ for $\alpha\in \Phi^+$. We will do this by induction on 
$\hgt(\alpha)$. If $\hgt(\alpha)=1$, then $\alpha=\alpha_i$ where 
$i\in I$. In this case, $N_{\alpha_i,\beta}^+=\epsilon(i)(q_{i,\beta}
+1)$; see Theorem~\ref{canbas} and note that $\be_{\alpha_i}^+=
\epsilon(i)e_i$. Now let $\hgt(\alpha)>1$. By the Key Lemma~\ref{keylem}, 
there exists some $i\in I$ such that $\langle \alpha_i^\vee,\alpha
\rangle >0$ and $s_i(\alpha)\in \Phi^+$. Applying $n_i(1)$ to 
$[\be_\alpha^+,\be_\beta^+]=N_{\alpha,\beta}^+\be_{\alpha+\beta}^+$ and 
using the above formula, we obtain that
\[ N_{\alpha,\beta}^+=(-1)^{q_{i,\alpha}+q_{i,\beta}+q_{i,\alpha+\beta}}
N_{s_i(\alpha),s_i(\beta)}^+.\]
Now $1\leq \hgt(s_i(\alpha))=\hgt(\alpha)-\langle \alpha_i^\vee,\alpha
\rangle <\hgt(\alpha)$ and so the structure constant on the right hand 
side is known by induction. 
\end{exmp}

Having obtained the explicit formulae in Theorem~\ref{luform3}, we 
also obtain formulae for the elements $\bar{h}_i(\xi)$ introduced 
in Corollary~\ref{barni2}.

\begin{prop} \label{luform3c} Let $i\in I$ and $\xi\in K^\times$. Then 
\[ \bar{h}_i(\xi)(\bar{h}_j^+)=\bar{h}_j^+ \qquad \mbox{and}\qquad
\bar{h}_i(\xi)(\bar{\be}_\alpha^+)=\xi^{\langle \alpha_i^\vee,\alpha
\rangle}\bar{\be}_\alpha^+\]
for all $j\in I$ and $\alpha\in \Phi$. In particular, $\bar{h}_i(\xi)$ 
is represented by a diagonal matrix with respect to the basis~$\bar{\bB}$. 
Furthermore, $\bar{h}_i(1)=\id_{\bar{\fg}}$ and $\bar{h}_i(\xi\xi')=
\bar{h}_i(\xi)\bar{h}_i(\xi')$ for all $\xi,\xi'\in K^\times$. 
\end{prop}

\begin{proof} Let $j\in I$. By Theorem~\ref{luform3}, $v:=\bar{n}_i
(\xi)(\bar{h}_j^+)\in \tilde{\fg}$ does not depend on~$\xi$.
Hence, we also have $\bar{n}_i(-1)(\bar{h}_j^+)=v$ and 
$\bar{n}_i(-\xi)(\bar{h}_j^+)=v$. Since $\bar{n}_i(\xi)^{-1}=
\bar{n}_i(-\xi)$, we deduce that 
\[\bar{n}_i(\xi)(v)=\bar{n}_i(\xi)\bigl(\bar{n}_i(-\xi)(\bar{h}_j^+)
\bigr)=\bar{n}_i(\xi)\bigl(\bar{n}_i(\xi)^{-1}(\bar{h}_j^+)\bigr)=
\bar{h}_j^+\]
and so $\bar{h}_i(\xi)(\bar{h}_j^+)=\bar{n}_i(\xi)\bigl(\bar{n}_i(-1)
(\bar{h}_j^+)\bigr)=\bar{n}_i(\xi)(v)=\bar{h}_j^+$, as claimed.

Next, let $\alpha\in \Phi$ and set $m:=\langle \alpha_i^\vee,\alpha\rangle
\in \Z$. Then the formula in Theorem~\ref{luform3} shows that 
\[\bar{n}_i(\xi)(\bar{\be}_\alpha^+)=\xi^{-m}v\quad \mbox{where 
$v:=\pm \bar{\be}_{s_i(\alpha)}^+$ does not depend on~$\xi$}.\]
So we also have $\bar{n}_i(-1)(\bar{\be}_\alpha^+)=(-1)^{-m}v$ and 
$\bar{n}_i(-\xi) (\bar{\be}_\alpha^+)=(-\xi)^{-m}v$. Again, since 
$\bar{n}_i(\xi)^{-1}=\bar{n}_i(-\xi)$, we deduce that 
\[\bar{n}_i(\xi)(v)=(-\xi)^m\bar{n}_i(\xi)\bigl(\bar{n}_i(-\xi)
(\bar{\be}_\alpha^+)\bigr)=(-\xi)^m\bar{\be}_\alpha^+\]
and so 
\begin{align*}
\bar{h}_i(\xi)(\bar{\be}_\alpha^+)& =\bar{n}_i(\xi)\bigl(
\bar{n}_i(-1)(\bar{\be}_\alpha^+)\bigr) =(-1)^{-m}\bar{n}_i(\xi)(v)\\
&=(-1)^{-m}(-\xi)^m\bar{\be}_\alpha^+ =\xi^m\bar{\be}_\alpha^+,
\end{align*}
as claimed.
Once these formulae are established, it immediately follows that
$\bar{h}_i(1)=\id_{\bar{\fg}}$ and $\bar{h}_i(\xi\xi')=
\bar{h}_i(\xi)\bar{h}_i(\xi')$ for $\xi,\xi'\in K^\times$ . 
\end{proof}

\begin{defn} \label{luform3d} The formulae in Proposition~\ref{luform3c} 
show that 
\[\bar{H}:=\Bigl\{\,\prod_{i \in I} \bar{h}_i(\xi_i)\,\big|\,\xi_i
\in K^\times\mbox{ for all $i\in I$}\Bigr\}\]
is an abelian subgroup of $\bar{G}$, where all elements of $\bar{H}$ are 
represented by diagonal matrices with respect to the basis~$\bar{\bB}$.
We call $\bar{H}$ the \nm{diagonal subgroup}\footnote{The importance of 
this subgroup will be fully revealed when we work over an algebraically
closed field $K$ and show that $\bar{G}$ is a linear algebraic group
over~$K$. It will then turn out that~$\bar{H}$ is a \nm{maximal torus} 
of~$\bar{G}$. Maximal tori are a crucial
ingredient in the general structure theory of linear algebraic groups.}
of $\bar{G}$. 
\end{defn}


Finally, we consider the commutator subgroup $[\bar{G},\bar{G}]$ of
$\bar{G}$. 

\begin{lem} \label{lemkchar4} For $i\in I$ we define the subgroup
\[\bar{G}_i:=\langle \bar{x}_i (\zeta),\bar{y}_i(\zeta)\mid \zeta \in K
\rangle\; \subseteq \;\bar{G}.\]
If $|K|\geq 4$, then $[\bar{G}_i,\bar{G}_i]=\bar{G}_i$.
\end{lem}

\begin{proof} Since $K$ is a field, the set $\{\zeta\in K\mid \zeta^2=1\}$
contains at most two elements. Hence, if $|K|\geq 4$, then there is 
some $\xi\in K^\times$ such that $\xi^2\neq 1$. Now let $\zeta\in K$
be arbitrary and set $\zeta':=\zeta(1-\xi^2)^{-1}\in K$. Using
Lemma~\ref{fromctok1} and Corollary~\ref{barni2}, we obtain:
\begin{align*}
[\bar{x}_i(\zeta'\xi^2),&\bar{h}_i(\xi)] =\bar{x}_i(\zeta'\xi^2)^{-1}
\bigl(\bar{h}_i(\xi)^{-1}\bar{x}_i(\zeta'\xi^2)\bar{h}_i(\xi)\bigr)\\
&=\bar{x}_i(-\zeta'\xi^2)\bar{x}_i(\zeta')=\bar{x}_i\bigl(\zeta'(1-
\xi^2)\bigr) =\bar{x}_i(\zeta).
\end{align*}
Hence, we have $\bar{x}_i(\zeta)\in [\bar{G}_i,\bar{G}_i]$. Similarly, 
one sees that $\bar{y}_i(\zeta)\in [\bar{G}_i, \bar{G}_i]$. Consequently,
$\bar{G}_i=[\bar{G}_i,\bar{G}_i]$.
\end{proof}

\begin{cor} \label{kchar4} If $|K|\geq 4$, then $\bar{G}$ is equal 
to its own commutator subgroup.
\end{cor}

\begin{proof} By Lemma~\ref{lemkchar4}, we have $\bar{G}_i=[\bar{G}_i,
\bar{G}_i]\subseteq [\bar{G},\bar{G}]$ for all $i\in I$. 
Hence, $\bar{G}=\langle \bar{G}_i\mid i \in I \rangle\subseteq 
[\bar{G}, \bar{G}]$. 
\end{proof}

\begin{rem} \label{kchar4a} If $K$ has only $2$ or $3$ elements, then it 
can actually happen that $[\bar{G},\bar{G}]\subsetneqq \bar{G}$. The 
situation is discussed in detail in Chevalley \cite[Th\'eor\`eme~3 
(p.~63)]{Ch}, Carter \cite[\S 11.1]{Ca1} or Steinberg \cite[Chapter~4]{St}. 
Altogether, if $A$ is indecomposable, then there are only the following 
four cases where $[\bar{G}, \bar{G}] \subsetneqq 
\bar{G}$. Suppose first that $|K|=2$. If $\cL$ is of type~$A_1$, 
then $\bar{G}$ has order $6$ and is isomorphic to the symmetric 
group~$\mathfrak{S}_3$; if $\cL$ is of type $B_2$, then $\bar{G}$ 
has order $720$ and is isomorphic to the symmetric group~$\mathfrak{S}_6$; 
if $\cL$ is of type $G_2$, then $\bar{G}$ has order $12096$ and 
there is a simple normal subgroup of index~$2$. The last exception occurs 
when $|K|=3$ and $\cL$ is of type $A_1$, in which case $\bar{G}$ 
has order $12$ and is isomorphic to the alternating group~$\mathfrak{A}_4$.
\end{rem}


The subgroups $\bar{G}_i\subseteq \bar{G}$ in Lemma~\ref{lemkchar4} are 
like {\itshape $\slm_2$-triples} in the Lie algebra~$\fg$. We will
encounter them again later on. Eventually, in Section~\ref{secA1}, we 
will see that there is always a surjective homomorphism $\SL_2(K)
\rightarrow \bar{G}_i$, with kernel contained in $\{\pm I_2\}$. At this 
point, we only show a kind of ``normal form'' for the elements 
in~$\bar{G}_i$. There will be no uniqueness of expressions but, as we 
shall see, the ``normal form'' is sufficient for certain purposes; see, 
e.g., Exercise~\ref{xcasl2tripGi} below. The following proof also provides
a good illustration of an efficient use of the various relations in 
$\bar{G}$ obtained so far.


\begin{prop}[Chevalley] \label{sl2tripGi} Let $i\in I$. Then we have
\[ \bar{G}_i=\big\{\bar{x}_i(\zeta_1)\bar{y}_i(\zeta_2)\bar{h}_i(\xi)
\bar{x}_i(\zeta_3)\mid\zeta_1,\zeta_2,\zeta_3\in K,\xi\in K^\times\big\}.\]
\end{prop}


\begin{proof} We must show $\bar{G}_i=\fX_i\fY_i\fH_i\fX_i$, where we set
\begin{align*}
\fX_i:&=\{\bar{x}_i(\zeta)\mid \zeta \in K\},\\
\fY_i:&=\{\bar{y}_i(\zeta)\mid \zeta \in K\},\\
\fH_i:&=\{\bar{h}_i(\xi)\mid \xi \in K^\times\}.
\end{align*}
By Lemma~\ref{fromctok1} and Propositon~\ref{luform3c}, these are all 
subgroups of~$\bar{G}$; furthermore, the maps $\bar{x}_i\colon K^+
\rightarrow \bar{G}$ and $\bar{y}_i\colon K^+\rightarrow \bar{G}$ are
group homomorphisms. Clearly, we have $\fX_i,\fY_i\subseteq\bar{G}_i$. 
Since $\bar{n}_i(\xi)\in \bar{G}_i$, we also have $\bar{h}_i(\xi) \in
\bar{G}_i$ for all $\xi\in K^\times$. Hence, $\fH_i\subseteq\bar{G}_i$. 
By Lemma~\ref{barni1} and Corollary~\ref{barni2}, we have the following 
relations:
\begin{itemize}
\item[(1)] $\bar{n}_i(\xi)\fX_i\bar{n}_i(\xi)^{-1}=\fY_i$
and $\bar{n}_i(\xi)\fY_i\bar{n}_i(\xi)^{-1}=\fX_i$.
\item[(2)] $\fH_i\fX_i=\fX_i\fH_i$ and $\fH_i\fY_i=\fY_i\fH_i$.
\end{itemize}
Let us set $D:=\fX_i\fY_i \fH_i \fX_i$. Since $\fX_i,\fY_i,\fH_i\subseteq D$, 
we have $D\subseteq \bar{G}_i$ and it will be sufficient to show that~$D$ is 
closed under inversion and multiplication. Firstly, we have
\[ D^{-1}=\fX_i^{-1}\fH_i^{-1}\fY_i^{-1}\fX_i^{-1}=
\fX_i\fH_i\fY_i\fX_i=D,\]
where the second equality holds since $\fX_i,\fY_i,\fH_i$ are subgroups 
and the third equality holds by~(2). Thus, $D$ is closed under inversion. 
To show that $D$ is closed under multiplication, it is sufficient to show 
that $\fX_iD\subseteq D$, $\fH_iD\subseteq D$ and $\fY_iD \subseteq D$. 
This is clear for the first two inclusions, by (2) and since $\fX_i$ is a 
subgroup. So it remains to show that $\bar{y}_i(\zeta)D\subseteq D$ for 
all $\zeta\in K$. This is clear for $\bar{y}_i(0)=\id$, so now assume
that $\zeta\neq 0$. By definition, we have
\[ \bar{n}_i(-\zeta^{-1})=\bar{x}_i(-\zeta^{-1})\bar{y}_i(\zeta)
\bar{x}_i(-\zeta^{-1})\]
and so $\bar{y}_i(\zeta)\in \fX_i \bar{n}_i(-\zeta^{-1})\fX_i$.
Hence, since $\fX_iD\subseteq D$, it will be sufficient to show that 
$\bar{n}_i(-\zeta^{-1})D \subseteq D$ for all $\zeta\in K^\times$. 
In other words (and changing variables), we must show that 
\begin{equation*}
\bar{n}_i(\xi)\bar{x}_i(\zeta)\fY_i\fH_i\fX_i\subseteq D
\qquad \mbox{for all $\zeta\in K$, $\xi \in K^\times$}.\tag{$\dagger$}
\end{equation*}
If $\zeta=0$, then $\bar{x}_i(0)=\id$ and $\bar{n}_i(\xi)\fY_i
\subseteq\fX_i\bar{n}_i(\xi)$, by (1); furthermore, $\bar{n}_i(\xi)
\in\fX_i\fY_i\fX_i$ by the definition of $\bar{n}_i(\xi)$ and so
\begin{align*}
\bar{n}_i(\xi)\fY_i\fH_i\fX_i &\subseteq (\fX_i
\bar{n}_i(\xi)) \fH_i\fX_i \subseteq \fX_i(\fX_i
\fY_i\fX_i)\fH_i\fX_i\\
&\subseteq \fX_i\fY_i\fX_i\fH_i\fX_i=
\fX_i\fY_i \fH_i\fX_i\fX_i\subseteq D,
\end{align*}
as required, where we used that $\fX_i$ is a subgroup. Now let $\zeta
\neq 0$. By definition, we have $\bar{n}_i(-\zeta)=\bar{h}_i(-\zeta)
\bar{n}_i(-1)^{-1}$ and $\bar{n}_i(\xi)=\bar{h}_i(\xi)\bar{n}_i(-1)^{-1}$. 
Hence, $\bar{n}_i(\xi)=h\bar{n}_i(-\zeta)$ where $h:=\bar{h}_i(\xi)
\bar{h}_i(-\zeta)^{-1}\in \fH_i$ (since $\fH_i$ is a subgroup). It follows 
that 
\begin{align*}
\bar{n}_i(\xi)\bar{x}_i(\zeta) &=h\bar{n}_i(-\zeta)\bar{x}_i(\zeta)=h
\bar{x}_i(-\zeta)\bar{y}_i(\zeta^{-1})\bar{x}_i(-\zeta)\bar{x}_i(\zeta) 
\\&=h\bar{x}_i(-\zeta) \bar{y}_i(\zeta^{-1})\in \fH_i\fX_i
\fY_i=\fX_i\fY_i\fH_i.
\end{align*}
where we used (2) and the fact that $\bar{x}_i\colon K^+\rightarrow 
\fX_i$ is a group homomorphism. So, finally, we obtain 
\[\bar{n}_i(\xi)\bar{x}_i(\zeta)\fY_i\fH_i\fX_i\subseteq (\fX_i\fY_i\fH_i)
(\fY_i\fH_i\fX_i)\subseteq \fX_i\fY_i\fY_i\fH_i\fH_i\fX_i=D,\]
as required, where we used again (2) and that $\fH_i$, $\fY_i$ are 
subgroups. Thus, ($\dagger$) is proved and so $D$ is a subgroup; hence, 
$\bar{G}_i=D$.
\end{proof}

Much later, we shall establish a generalisation
of the above normal form for $\bar{G}_i$ to something analogous for 
{\itshape all} elements of $\bar{G}$. This will allow us, for example, to 
deduce that the diagonal subgroup $\bar{H}$ in Definition~\ref{luform3d} 
consists precisely of all elements of the {\itshape whole} group $\bar{G}$ 
that are represented by diagonal matrices with respect to $\bar{\bB}$.
This seems to be difficult to prove directly at this stage. We invite the 
reader to try to prove this as far as $\bar{G}_i$ is concerned; see 
the following exercise. 

\begin{xca} \label{xcasl2tripGi} In the setting of 
Proposition~\ref{sl2tripGi}, show that every element of $\bar{G}_i$ that 
is represented by a diagonal matrix is of the form $\bar{h}_i(\xi)$ for 
some $\xi \in K^\times$. (Use the fact that elements in $\fX_i$ are 
represented by upper triangular matrices with~$1$ along the diagonal, and 
similar properties of the matrices of the elements in $\fY_i$ and $\fH_i$.)
\end{xca}

\section*{Notes on Chapter~\ref{chap3}}

For the classification of generalized Cartan matrices we follow Kac 
\cite[Chap.~4]{K}; see also Moody--Pianzola \cite[\S 3.5, \S 3.6]{MP}. 
A somewhat different, and arguably simpler, approach is given by Lusztig 
\cite[14.1.4--14.1.6]{L6} but, in any case, we shall need the 
characterisations of (FIN) in Theorem~\ref{kac2} and Corollary~\ref{kac3}
at some places later on. Systematic descriptions of the irreducible root 
systems of the various (finite) types can be found in Bourbaki 
\cite[Ch.~VI, \S 4, no.~4.4--4.13]{B}; see also Benson--Grove 
\cite[\S 5.3]{BG} for explicit constructions.


See Kac \cite[\S 1.9]{K} for some notes about the historical development
of the study of Kac--Moody Lie algebras. The appendix of Moody--Pianzola 
\cite{MP} contains a much more thorough discussion of 
Example~\ref{a1tilde}. The idea of replacing $\C$ by a ring of Laurent 
polynomials can be generalized to all Lie algebras of Cartan--Killing 
type; see, e.g., Carter \cite[Chap.~18]{Ca3}
for a detailed exposition. Lemmas~\ref{triang1}--\ref{triang3} are 
analogous to certain steps in the proof of \cite[Theorem~18.2]{H} (which 
prepare Serre's theorem mentioned in Remark~\ref{serre}). There are 
several other proofs of the important Existence Theorem~\ref{myG3}: 
\begin{itemize} \itemsep2pt  
\item Via free Lie algebras and definitions in terms of generators 
and relations. See Jacobson \cite[Chap.~VII, \S 4]{Jac}, Serre 
\cite[Chap.~VI, Appendix]{S} (and also \cite[\S 18]{H} for further details).
\item Via explicit descriptions of structure constants. There is 
an elegant way to do this for $A$ of simply laced type; the remaining cases 
are obtained by  a ``folding'' procedure. See Kac \cite[\S 7.8 and 
\S 7.9]{K}, with further details in De Graaf \cite[\S 5.13--\S 
5.15]{graaf}. For a general approach see Tits \cite{Ti}.
\item Via explicit constructions in all cases. Historically, this is
the original method. For the classical types $A_n$, $B_n$, $C_n$, $D_n$, 
we have seen this already. For the exceptional types, see Fulton--Harris 
\cite[\S 22.4]{FH} for further details and references.
\end{itemize}
The approach via Definition~\ref{myG0} works uniformly for all types of $A$
and is completely elementary and self-contained; it does not use free Lie 
algebras or any other further theory, except for the verification of the 
relations in Lemma~\ref{myG1}.
 
The {\sf ChevLie} package presented in Section~\ref{sec3comp} is one 
example of a whole variety of software packages for Lie theory. The 
computer algebra systems {\GAP} \cite{gap4} and {\sf Magma} \cite{magma} 
contain large packages for Lie theory; we also mention the Lie Atlas 
project \cite{adams} here. Some research articles heavily relying on 
computational methods are Holt--Spaltenstein \cite{HoSp}, Gilkey--Seitz 
\cite{GiSe}, Cohen--Murray--Taylor \cite{Coet}.

See Casselman \cite{Cass1} for a slightly different version of the 
recursive algorithm for computing structure constants in 
Example~\ref{exmcass}. The appendix of \cite{GiSe} contains explicit 
tables with the structure constants $N_{\alpha,\beta}$ for types $G_2$, 
$F_4$, $E_6$ and $E_7$; for type $E_8$ see Mizuno \cite{Miz2}. These 
tables rely on \textit{some} choice of elements $0\neq e_\alpha\in \cL$ 
for $\alpha\in \Phi$; hence, in general, they will be different from our 
$N_{\alpha,\beta}^+$. See Ringel \cite{Ri} and Li--Xiao \cite{lixiao}
for a different approach to those structure constants. 

Once the explicit formulae in Theorem~\ref{luform} are available, the
identities concerning the elements $\bar{n}_i(\xi)$ and $\bar{h}_i(\xi)$ 
follow from rather straightforward computations. Theorem~\ref{luform3} 
appears in \cite[\S 5]{G1} (with a somewhat different proof); the analogous 
formulae in \cite[p.~36]{Ch} or \cite[Prop.~6.4.2]{Ca1} involve some 
unspecified signs. The Transfer Lemma~\ref{superexp} appears in 
\cite[Lemma~4.5.1]{Ca1}. The discussion of Chevalley groups associated 
with the Lie algebra $\slm_2(\C)$ will be continued in Section~\ref{secA1}.


For a completely different approach to Chevalley groups, using
extensively the general theory of ``canonical bases'', see Lusztig 
\cite{L3}. In this context, one also obtains reductive algebraic groups 
(and not just semisimple ones). We also mention that there are groups 
associated with Kac--Moody Lie algebras; see, for example, Marquis \cite{Mar}.

Finally, note that Definition~\ref{deflu} actually works with $K$ replaced 
by a commutative ring $R$ with~$1$. In this context, $G_R(\fg)$ would be 
called an \nm{elementary Chevalley group}; see, e.g., Plotkin--Vavilov 
\cite{VP} and further references there. We just note that there are subtle 
problems arising from working over a ring instead of a field, but we will
not elaborate this here in any further detail. This is also briefly
discussed in Steinberg lectures \cite[Chap.~7, pp.~56--59]{St}. 

\chapter{General construction of Chevalley groups} \label{chap4}

\addtocounter{footnote}{15}

Consider a Chevalley group $G_K(\fg)$ as defined in the previous
chapter, where $\fg$ is a Lie algebra of Cartan--Killing type and
$K$ is a field. Let us assume, for a moment, that $K$ is algebraically
closed. Then $G_K(\fg)$ is a (connected) \nm{semisimple algebraic group}, 
with a trivial center. Now Chevalley's famous seminar \cite{Ch3} provides
a complete classification of \textit{all} semisimple algebraic groups 
over~$K$. It turns out that these are still related to Lie algebras~$\cL$
as above but there may be several groups~---~in fact, finitely many up to 
isomorphism~---~corresponding to the same~$\cL$ (and these groups may have 
a non-trivial center). It is one purpose of this chapter to construct 
explicitly \textit{all} the semisimple groups occurring in Chevalley's 
classification. The construction will work for arbitrary fields $K$, not 
just algebraically closed ones.

The Lie algebra $\fg=\slm_n(\C)$ already provides a good illustration.
We have seen that $G_K(\fg)$ is isomorphic to the quotient of $\SL_n(K)$
by the normal subgroup consisting of scalar matrices with determinant~$1$.
It is quite natural to ask if there is a modified construction which 
would produce the whole group $\SL_n(K)$. Analogous questions arise for 
the possible groups associated with the Lie algebras $\gom_n(Q_n,\C)$ of 
classical type; for example, there are the symplectic, orthogonal and 
spin groups (which may have a non-trivial center). 

Recall from Section~\ref{sec3a7} the construction of $G_K(\fg)$: We start 
with a subgroup $G_\C(\fg)\subseteq \Aut(\cL)$ over $\C$, generated by 
elements that are obtained by taking the exponential of the nilpotent
maps $\ad_\cL(e_i)$ and $\ad_\cL(f_i)$ for $i \in I$. Then it is made 
sure that certain integrality conditions hold, which eventually allow 
the passage to an arbitrary field~$K$. Now the basic idea is to replace 
the homomorphism $\ad_\cL\colon \cL\rightarrow \gl(\cL)$ by suitable 
representations $\rho\colon\cL \rightarrow \gl(V)$ where $V$ are 
finite-dimensional $\C$-vector spaces. Again, one has to make sure that 
certain integrality conditions hold which eventually should allow the 
passage to an arbitrary~$K$. Already from this very informal description, 
it becomes clear that we need to know a bit more about the representation 
theory of~$\cL$. This is a vast area of research; see, for example, the 
relevant chapters of Fulton--Harris \cite{FH} or Humphreys \cite{H}. Here, 
we will give a self-contained introduction where we focus on those aspects 
that are particularly relevant for the construction of (non-adjoint) 
Chevalley groups later on. This includes some topics which are often not 
covered in textbooks (like \cite{FH} or \cite{H}, for example), most 
notably the representations corresponding to ``minuscule'' weights.
After some general preparations in Section~\ref{sechighw}, we develop 
this from first principles in Sections~\ref{sec5a1} and \ref{sechwm} 
(which may be of independent interest). We shall see that, by taking 
direct sums of the adjoint representation and, possibly, various  
``minuscule'' representations, we obtain a sufficiently rich family of 
representations which yield all the groups in Chevalley's classification.

Sections~\ref{sec5a2}--\ref{secdiagII} lay the groundwork for defining 
and studying these groups. Our exposition here (and in the following 
chapter) is a synthesis of Ree's article \cite{ree2} and Steinberg's 
lectures~\cite{St}, but we shall provide significantly more details 
on various arguments and calculations~---~similar in style to Carter's
book \cite{Ca1}. We hope that readers will find this a useful addition
to the existing literature. (For further general comments and a comparison
between the approaches of \cite{ree2} and \cite{St}, see also the notes 
at the end of this chapter.) The final Section~\ref{secA1} explicitly 
describes all possible Chevalley groups associated with the Lie 
algebra $\fg=\slm_2(\C)$; this is a model case for the general
classification to be discussed in a later chapter.

\section{The weight lattice of a representation} \label{sechighw}

We return to the setting of Chapter~\ref{chap2}. Let $\cL$ be a 
finite-dimensional Lie algebra over $\C$ and $\fh\subseteq \cL$ be an
abelian subalgebra such that $(\cL,\fh)$ is of \nm{Cartan--Killing type} 
with respect to a linearly independent subset $\Delta= \{\alpha_i \mid 
i \in I\}\subseteq \fh^*$. We have the Cartan decomposition $\cL=
\cH\oplus \bigoplus_{\alpha\in \Phi}\cL_\alpha$ where $\Phi\subseteq\cH$ 
is the root system of $\cL$. For $i\in I$ let $\{e_i,f_i,h_i\} \subseteq 
\cL$ be a corresponding \nms{$\slm_2$-triple}{sl$_2$-triple} (as in 
Remark~\ref{astring0}). Thus, $\cL_i:=\langle e_i,f_i,h_i \rangle_\C 
\subseteq \cL$ is a subalgebra isomorphic to $\slm_2(\C)$. The aim of 
this section is to establish some basic results on $\cL$-modules, which 
may be seen as generalisations of the results on $\slm_2(\C)$-modules 
in Section~\ref{sec04a}. Throughout, we will tacitly assume that 
\begin{center}
\fbox{all $\cL$-modules under consideration are finite-dimensional.}
\end{center}
We recall some further notation. We have a basis of $\fh$ consisting of the
``co-roots'' $\{h_i\mid i \in I\}$. Let $A=(a_{ij})_{i,j \in I}$ be the
corresponding structure matrix, where $a_{ij}=\alpha_j(h_i)$ for all 
$i,j\in I$. As in Section~\ref{sec1a3}, let $E\subseteq \fh^*$ be the 
$\R$-span of $\{\alpha_i \mid i \in I\}$ and $\Phi\subseteq E$ be the 
corresponding root system, with Weyl group $W=W(A)\subseteq \GL(E)$. We
have $W=\langle s_i \mid i \in I\rangle$, where $s_i$ acts on $E$ via 
\[ s_i(\lambda)=\lambda-\lambda(h_i)\alpha_i \qquad (\lambda \in E).\]
Now let $V$ be a $\cL$-module. Thus, $V$ is a $\C$-vector space where the
$\cL$-module structure is given by a bilinear map $\cL\times V \rightarrow 
V$, $(x,v)\mapsto x.v$, satisfying the conditions in Section~\ref{sec04}. 
As before, we denote by $\rho_x\colon V \rightarrow V$ the linear map
defined by $x \in \cL$. Thus, $\rho\colon \cL\rightarrow \gl(V)$, $x\mapsto 
\rho_x$, is a representation. By restricting $\rho$ to $\fh\subseteq \cL$, 
we may regard $V$ as an $\fh$-module. We set 
\[ V_\lambda:=\{v\in V\mid h.v=\lambda(h)v \mbox{ for all $h\in \fh$}\}
\quad\mbox{for any $\lambda \in \fh^*$}.\]
If $V_\lambda\neq\{0\}$, then $\lambda$ is called a \nm{weight} of $\fh$
on $V$ and $V_\lambda$ the corresponding \nm{weight space}. 

Recall from Section~\ref{sec1a1} the notion of ``$\fh$-diagonalisability''
of~$V$; this means that $V$ is a direct sum $V=\bigoplus_{\lambda \in 
P_\fh(V)} V_\lambda$, where $P_\fh(V) \subseteq \fh^*$ denotes the set 
of weights of~$\fh$ on~$V$. (In an appendix, we 
shall see that every $V$ is automatically $\fh$-diagonalizable.)


\begin{rem} \label{fundcal} We note the following analogue of 
Proposition~\ref{wsdprop2}. Let $V$ be a $\cL$-module. Let $\lambda\in 
\fh^*$ and $\alpha\in \Phi$. Then the action of $\cL_\alpha$ on $V$ maps 
$V_\lambda$ into $V_{\lambda+\alpha}$. Indeed, let $v\in V_\lambda$ and 
set $v':=e_\alpha.v$, where $0\neq e_\alpha \in \cL_\alpha$. Then we 
compute for $h\in \cH$:
\[h.v'=[h,e_\alpha].v+e_\alpha.(h.v)=\alpha(h)e_\alpha.v+\lambda(h)
e_\alpha.v=(\lambda+\alpha)(h)v'.\]
Hence, we have $e_\alpha.v=v'\in V_{\lambda+\alpha}$. Fulton--Harris
\cite[p.~148]{FH} call this the ``\nm{fundamental calculation}''.
\end{rem}

Here are the first consequences of the ``fundamental calculation''.
The following result will allow us to apply the exponential construction
in Lemma~\ref{exponential} to many elements in $\cL$. (We regard $V$ as 
an algebra with trivial product $v\cdot v'=0$ for all $v,v'\in V$.)

\begin{lem} \label{wsdnilV} Assume that $V$ is $\fh$-diagonalizable. Let 
$\alpha\in \Phi$ and $0\neq e_\alpha\in \cL_\alpha$. Then the linear map 
$\rho(e_\alpha)\colon V\rightarrow V$ is nilpotent. Consequently,
$\exp(t\rho(e_\alpha))\in \GL(V)$ is defined for any $t\in \C$.
\end{lem}
  
\begin{proof} For any $\lambda\in P_\fh(V)$, we have $\rho(e_\alpha)
(V_\lambda) \subseteq V_{\lambda+\alpha}$; see Remark~\ref{fundcal}.
Hence, we have $\rho(e_\alpha)^m(V_\lambda) \subseteq V_{\lambda+m\alpha}$ 
for any integer $m>0$. Since $P_\fh(V)$ is a finite set, we must have $\rho
(e_\alpha)^m(V_\lambda)=\{0\}$ for some $m>0$ (depending on $\lambda$).
Since $P_\fh(V)$ is finite, there exists some $d>0$ such that $\rho
(e_\alpha)^d(V_\lambda)= \{0\}$ for all $\lambda \in P_\fh(V)$. Since $V=
\sum_{\lambda \in P_\fh(V)} V_\lambda$, it follows that $\rho(e_\alpha)^d=0$.
Then $t\rho(e_\alpha)$ is nilpotent for any $t\in \C$ and so we can apply
Lemma~\ref{exponential}.
\end{proof}

\begin{prop}[Integrality and $W$-invariance] \label{integralrep} Assume 
that $V$ is $\fh$-diagonalizable and let $\lambda\in P_\fh(V)$. Then 
$\lambda(h_i) \in \Z$ for all $i\in I$ and $w(\lambda)\in P_\fh(V)$ for 
all $w\in W$.
\end{prop}

\begin{proof} Let $i\in I$ and $\cL_i:=\langle e_i,f_i,h_i\rangle_{\C}
\subseteq \cL$, as above. We regard $V$ as a $\cL_i$-module (via 
restriction). Since $\cL_i\cong \slm_2(\C)$, the discussion in 
Remark~\ref{sl2highb} applies to the $\cL_i$-module~$V$, where $e_i,f_i,
h_i \in \cL_i$ play the role of $e,f,h \in \slm_2(\C)$, respectively. 
The following argument is very similar to that in Remark~\ref{astring}. 
Since $P_\fh(V)$ is a finite set, there are integers $p,q\geq 0$ such that 
$\lambda+m\alpha_i\in P_\fh(V)$ for $-q\leq m\leq p$ but $\lambda+(p+1) 
\alpha_i\not\in P_\fh(V)$ and $\lambda-(q+1)\alpha_i \not\in P_\fh(V)$. 
Then consider the subspace
\[ M:=V_{\lambda-q\alpha_i}\oplus \ldots \oplus V_{\lambda-\alpha_i}
\oplus V_{\lambda}\oplus V_{\lambda+\alpha_i}\oplus \ldots \oplus 
V_{\lambda+p\alpha_i} \subseteq V.\]
Using Remark~\ref{fundcal} one easily sees that $M$ is a $\cL_i$-submodule
of~$V$. Since $\alpha_i(h_i)=2$, the eigenvalues of $\rho_{h_i}\colon 
M\rightarrow M$ are given by $\lambda(h_i)+2m$ for $-q\leq m\leq p$ (not 
counting multiplicities). By Proposition~\ref{sl2modd}, the largest 
eigenvalue must be the negative of the smallest eigenvalue. Hence, 
$\lambda(h_i)+2p=-(\lambda(h_i)-2q)$ and so $\lambda(h_i)=q-p\in \Z$. 
Consequently, we have 
\[s_i(\lambda)= \lambda+(p-q) \alpha_i\in \{\lambda+m\alpha_i\mid -q\leq m
\leq p\},\]
and so $s_i(\lambda)\in P_\fh(V)$. Since $W=\langle s_i\mid i \in I\rangle$, 
this also implies that $w(\lambda)\in P_\fh(V)$ for all $w\in W$.
\end{proof}

%
 
The above result naturally leads to the definition of certain ``lattices''
in $E$. By definition, a subset $\Lambda \subseteq E$ is called a
\nm{lattice} if $\Lambda$ is a subgroup of $(E,+)$ and if there is a basis 
$\{v_i \mid i \in I\}$ of $E$ such that $\Lambda=\sum_{i\in I} \Z v_i$. For 
example, let $P:=\langle \Phi\rangle_\Z\subseteq E$, the set consisting of 
all $\Z$-linear combinations of roots $\alpha\in\Phi$. Then 
\[ P=\langle\Phi\rangle_\Z=\sum_{i\in I} \Z\alpha_i\subseteq E\quad
\mbox{is a lattice}\]
(since $\{\alpha_i\mid i \in I\}$ is an $\R$-basis of $E$). We call 
$P$ the \nm{root lattice} of $\cL$.

\begin{defn} \label{weightlat} Motivated by Proposition~\ref{integralrep},
we define 
\[\Omega=\Omega(\cL):=\{\lambda\in E\mid\lambda(h_j) \in\Z\mbox{ for 
all $j\in I$}\}\;\subseteq E.\]
This will be called the \nmi{weight lattice}{weight lattice (of $\cL$)} 
of $\cL$. By Remark~\ref{poshi}, we have $\beta(h_j)\in\Z$ for $j\in I$ and 
$\beta \in \Phi$; hence, we have $P=\langle \Phi\rangle_\Z\subseteq\Omega$. 

Clearly, $\Omega$ is a subgroup of $(E,+)$, but it is not entirely obvious 
that~$\Omega$ indeed is a lattice. In order to see this, let $\{\varpi_i 
\mid i \in I\}$ be the basis of $\fh^*$ that is dual to the basis $\{h_j 
\mid j\in I\}$ of $\fh$, that is, we have 
\[ \varpi_i(h_j)=\left\{\begin{array}{cl} 1 & \quad\mbox{if $i=j$},\\ 
0 &\quad\mbox{if $i\neq j$}.  \end{array}\right.\]
The $\varpi_i$ are called \nm{fundamental weights}. Note the following 
identity:
\begin{equation*}
\lambda=\sum_{i \in I} \lambda(h_i)\varpi_i \qquad\mbox{for any $\lambda
\in \fh^*$}.\tag{$*$}
\end{equation*}
(In order to prove this, just evaluate both sides on $h_j$ for $j\in I$.) 
In particular, we obtain that 
\begin{equation*}
\alpha_j=\sum_{i\in I} \alpha_j(h_i)\varpi_i=\sum_{i\in I} 
a_{ij} \varpi_i \qquad \mbox{for $j\in I$}.\tag{$*^\prime$}
\end{equation*}
Since $A=(a_{ij})$ is invertible and has entries in~$\Z$, the above 
equations show that $\varpi_i \in \langle \alpha_j \mid j \in I\rangle_\Q
\subseteq E$ for all $i\in I$. We now claim that 
\[\Omega= \langle \varpi_i\mid i \in I\rangle_\Z \quad \mbox{is a lattice
in $E$}.\]
Indeed, the inclusion ``$\subseteq$'' immediately follows using ($*$).
Conversely, let $\lambda \in \langle \varpi_j\mid j \in I
\rangle_\Z$ and write $\lambda=\sum_{j\in I} m_j \varpi_j$ with
$m_j \in\Z$ for all $j\in I$. Then $\lambda(h_i)=m_i\in\Z$ for all $i\in I$,
and so $\lambda \in \Omega$. Thus, the above equality is proved. 
\end{defn}

\addtocounter{table}{11}
\begin{table}[htbp] \caption{Fundamental groups for $A$ indecomposable} 
\label{fundtab} 
\begin{center}
{\small $\renewcommand{\arraystretch}{1.2}
\begin{array}{lc@{\hspace{20pt}}l} \hline \mbox{Type} & \det(A) & \Omega/
\langle \Phi\rangle_\Z\\ \hline A_n\,(n{\geq}1) & n+1 & \; \Z/(n+1)\Z\\ 
B_n,C_n\, (n{\geq}2) & 2 & \;\Z/2\Z\\ D_{n}\, (n{\geq}3)  & 4& 
\left\{\begin{array}{cl} \Z/2\Z \times\Z/2\Z & \mbox{ ($n$ even)}\\
\Z/4\Z & \mbox{ ($n$ odd)} \end{array}\right.\\ E_6  & 3 & \;\Z/3\Z\\
E_7  & 2 & \;\Z/2\Z \\ G_2,F_4,E_8  & 1 & \;\{0\}\\ \hline 
\end{array}$}
\end{center}
\end{table}

\begin{rem} \label{weightlb} Since $P=\langle \Phi\rangle_\Z\subseteq 
\Omega$, we can form the factor group $\Omega/P$ which is called the 
\nm{fundamental group} of $\Phi$. We claim that 
\[ \Omega/P\;\mbox{is a finite abelian group with} \;\;
|\Omega/P|=\det(A).\]
This is seen as follows. By the above relation~($*^\prime$), the entries
of the structure matrix $A$ describe the expressions of the basis 
elements $\alpha_j$ of $P$ in terms of the basis elements $\varpi_i$ 
of~$\Omega$. So the assertion follows from the general \nm{theory of 
elementary divisors}; see, e.g., Artin \cite[Chap.~12, \S 4]{Ar}). This 
also yields a concrete method for determining the order and even the 
structure of $\Omega/P$. For indecomposable $A$, the results are in 
Table~\ref{fundtab}. (We leave it as an exercise for the reader to
verify the entries of that table.)
\end{rem}

\begin{rem} \label{lemstem0a} Consider the Weyl group $W=\langle
s_i\mid i\in I\rangle\subseteq \mbox{GL}(E)$. For $i\in I$ and
$\lambda\in E$ we have $s_i(\lambda)=\lambda-\lambda(h_i)\alpha_i$.
If $j\in I$, then $\varpi_j(h_i)=\delta_{ij}$ (Kronecker delta) and so 
\[s_i(\varpi_j)=\varpi_j-\delta_{ij}\alpha_i \in\Omega.\]
It follows that $s_i(\Omega)\subseteq \Omega$ and so $w(\Omega)=\Omega$
for all $w\in W$. Thus, the action of $W$ on $E$ induces an action
of $W$ on $\Omega$. 
\end{rem}

\begin{defn} \label{deflambda} Let $V$ be an $\fh$-diagonaliable 
$\cL$-module. By Proposition~\ref{integralrep}, we have $P_\fh(V)\subseteq 
\Omega$. Then 
\[ \Lambda(V):=\langle P_\fh(V)\rangle_\Z \subseteq \Omega\]
will be called the \nm{weight lattice of $V$}. The following result shows 
that $\Lambda(V)$ indeed is a lattice in $E$ if $V$ is a \nm{faithful
$\cL$-module}, that is, the representation $\rho \colon \cL \rightarrow 
\gl(V)$ is an injective homomorphism.
\end{defn}

\begin{prop} \label{defomega} Let $V$ be an $\fh$-diagonalisable
$\cL$-module. Then~$V$ is a faithful $\cL$-module if and only if 
$\Phi\subseteq \Lambda(V) \subseteq \Omega$. In this case, $\Lambda(V)$ 
is a lattice in~$E$; in particular, $\langle P_\fh(V) \rangle_\R=E$.
\end{prop}

\begin{proof} Recall from Proposition~\ref{integralrep} that $\Lambda(V)
\subseteq \Omega$. Assume first that~$V$ is faithful. Let $\alpha\in 
\Phi$ and $0\neq e_\alpha \in \cL_\alpha$. Since $V$ is faithful and 
$e_\alpha\neq 0$, the map $\rho_{e_\alpha} \colon V \rightarrow V$ is 
non-zero. So there is some $0\neq v\in V$ such that $e_\alpha.v=
\rho_{e_\alpha}(v)\neq 0$. Since $V$ is $\fh$-diagonalisable, we have 
$V=\sum_{\lambda \in P_\fh(V)} V_\lambda$. So we can write $v=v_1+\ldots 
+v_r$ where $0\neq v_i \in V_{\lambda_i}$ and $\lambda_i\in P_\fh(V)\subseteq
\Lambda(V)$. Consequently, we must have $e_\alpha.v_i\neq 0$ for 
some~$i$. But then $0\neq e_\alpha.v_i\in V_{\lambda_i+\alpha}$ (see 
Remark~\ref{fundcal}) and so $\lambda_i+ \alpha\in \Lambda(V)$. Hence,
we also have $\alpha=(\lambda_i +\alpha)- \lambda_i\in \Lambda(V)$. This
shows that $\Phi\subseteq \Lambda(V)$, as desired.

Conversely, let $\Phi\subseteq \Lambda(V)$ and assume, if possible,
that $V$ is not faithful. Then $\ker(\rho)$ is a non-zero ideal of $\cL$ 
and so there is some $i\in I$ such that $h_i \in \ker(\rho)$; see 
Lemma~\ref{mainideal1}. Since $\cH^*=\langle \Phi\rangle_\C$ and 
$\Phi\subseteq \Lambda(V)$, we also have $\cH^*=\langle \Lambda(V)
\rangle_\C=\langle P_\fh(V)\rangle_\C$. So there exists some $\lambda 
\in P_\fh(V)$ such that $\lambda(h_i)\neq 0$. Hence, for $0\neq v \in 
V_\lambda$, we have $h_i.v=\lambda(h_i)v \neq 0$, contradiction to
$h_i\in \ker(\rho)$. 

Finally, assume that $\Phi\subseteq \Lambda(V)=\langle P_\fh(V)
\rangle_\Z$. Since $E=\langle \Phi\rangle_\R$, it follows that 
$\langle P_\fh(V)\rangle_\R=E$. Since $\Lambda(V)$ is a finitely generated
subgroup of a lattice in $E$, it is known that there are 
$\lambda_1, \ldots,\lambda_m\in \Lambda(V)$, where $m\leq \dim E$, such 
that $\Lambda(V)=\sum_{1\leq i \leq m} \Z \lambda_i$; see, e.g., Artin 
\cite[Chap.~12, \S 4]{Ar}. It follows that $E=\langle P_\fh(V)\rangle_\R=
\langle \lambda_1,\ldots,\lambda_m\rangle_\R$. Hence, we must have 
$m=\dim E$ and $\{\lambda_1,\ldots,\lambda_m\}$ is a basis of $E$.
\end{proof}

\begin{exmp} \label{defhwmL} We regard $V:=\cL$ as a $\cL$-module via the
adjoint representation; thus, $\cL$ acts on $V$ via $x.v=[x,v]$ for
$x\in \cL$ and $v \in V$. Since $(\cL,\fh)$ is of Cartan--Killing type, 
$V$ is $\fh$-diagonalisable where $P_\fh(V)=\{\underline{0}\}\cup \Phi$. So, 
in this case, $\langle \Phi\rangle_\Z=\Lambda(V) \subseteq \Omega$. Note 
that $V=\cL$ is a faithful $\cL$-module, since $\cL$ is semisimple and
so $Z(\cL)=\{0\}$ (see Remark~\ref{defsolv3a} and 
Proposition~\ref{mainideal}).
\end{exmp}

\begin{exmp} \label{reesl2} Let $\fg=\slm_2(\C)$, with $I=\{1\}$ and
standard basis $\{e_1,f_1,h_1\}$ such that $[e_1,f_1]=h_1$. Let
$\fh=\langle h_1\rangle_\C$ and $\varpi_1\in \fh^*$ be defined by
$\varpi_1(h_1)=1$. Then $\Omega=\langle \varpi_1\rangle_\Z$ and
$\Phi=\{\pm \alpha_1\}$ where $\alpha_1=2\varpi_1$. Let $V$ be a 
faithful, $\fh$-diagonalisable $\fg$-module. Then $\Phi \subseteq 
\Lambda(V)\subseteq \Omega$ and so there are only two possibilities: 
\begin{center}
Either $\quad \Lambda(V)=\langle \alpha_1\rangle_\Z\quad$ or
$\quad\Lambda(V)=\langle \varpi_1\rangle_\Z=\Omega$.
\end{center}
Some concrete examples:

If $V=\fg$ (adjoint representation), then $P_\fh(V)=\{\pm \alpha_1\}
\cup \{\underline{0}\}$ and so $\Lambda(V)=\langle \alpha_1\rangle_\Z$;
see Example~\ref{defhwmL}.

If $V$ is an irreducible $\fg$-module with $\dim V=m+1$ (where $m\geq 1$), 
then $\rho_{h_1}\colon V \rightarrow V$ has eigenvalues $m-2i$ for 
$0\leq i \leq m$; see Corollary~\ref{sl2modc}. Hence, $P_\fh(V)=\{m\varpi_1-
i\alpha_1\mid 0 \leq i \leq m\}$ and so $\Lambda(V)=\langle \alpha_1 
\rangle_\Z$ if $m$ is even, and $\Lambda(V)=\Omega$ if $m$ is odd.
\end{exmp}

\begin{exmp} \label{defhwm2} Let $\cL=\slm_n(\C)$ ($n\geq 2$) and $\fh
\subseteq \cL$ be the subalgebra consisting of diagonal matrices. Let $V=
\C^n$ be the natural $\cL$-module. Since $\cL$ is a simple Lie algebra,
$V$ is a faithful module. Now $(\cL,\fh)$ is of Cartan--Killing type where 
$I=\{1,\ldots,n-1\}$; the roots~$\Phi$, the weight spaces $\cL_\alpha$, and 
the elements $e_i$, $f_i$, $h_i$ ($i \in I$) are explicitly
described in Example~\ref{cartsln}. For $1\leq i \leq n$ let $\varepsilon_i
\in \fh^*$ be the map that sends a diagonal matrix to its $i$-th diagonal 
entry. Let $\{v_1,\ldots, v_n\}$ be the standard basis of~$V$. Then 
$h.v_i=\varepsilon_i(h)v_i$ for all $h\in \fh$ and so 
\[V=V_{\varepsilon_1}\oplus\ldots \oplus V_{\varepsilon_n}
\qquad\mbox{where}\qquad \dim V_{\varepsilon_i}=1 \;\mbox{for all $i$}.\]
In particular, $V$ is $\fh$-diagonalisable, where 
\[P_\fh(V)=\{\varepsilon_1,\ldots,\varepsilon_n\}\qquad \mbox{and}\qquad
\Lambda(V)=\langle \varepsilon_1,\ldots, \varepsilon_n\rangle_\Z.\]
Also note that, by Example~\ref{weylan}, the set $P_\fh(V)$ is a single 
$W$-orbit; we have $s_i(\varepsilon_i)=\varepsilon_{i+1}$ for $1\leq i 
\leq n-1$. The fundamental weights of~$\fg$ are given as follows:
\[\varpi_i=\varepsilon_1+ \ldots +\varepsilon_i \qquad\mbox{for $1\leq 
i\leq n-1$}.\]
Indeed, by Example~\ref{cartsln}, $h_j$ (for $1\leq j\leq n-1$) is the
diagonal matrix with entries $1,-1$ at positions $j,j+1$ (and $0$ 
otherwise). Hence, we have $\varepsilon_i(h_i)=1$, $\varepsilon_i(h_{i-1})
=-1$ (if $i\geq 2$) and $\varepsilon_i(h_j)=0$ if $j\not\in\{i,i-1\}$. 
Consequently, $\varepsilon_1=\varpi_1$ and 
\[ \varepsilon_i=\sum_{1\leq j\leq n-1} \varepsilon_i(h_j)\varpi_j=
\varepsilon_i(h_i)\varpi_i+\varepsilon_i(h_{i-1})\varpi_{i-1}=
\varpi_i-\varpi_{i-1}\]
for $2\leq i\leq n-1$. This yields the above expressions for $\varpi_i$.
Since $\varepsilon_1+\ldots +\varepsilon_n=\underline{0}$, it now also 
follows that 
\[ \Omega=\langle \varpi_1,\ldots,\varpi_{n-1} \rangle_\Z=\langle 
\varepsilon_1,\ldots,\varepsilon_n\rangle_\Z=\Lambda(V).\]
Finally, recall that the simple roots are given by $\alpha_i=\varepsilon_i-
\varepsilon_{i+1}$ for $1\leq i \leq n-1$. Consequently, we have
\[ \Omega=\langle \varepsilon_1,\ldots,\varepsilon_n\rangle_\Z=
\langle \varepsilon_1,\varepsilon_2-\varepsilon_1,\ldots,\varepsilon_n-
\varepsilon_{n-1}\rangle_\Z=\langle\Phi,\varpi_1\rangle_\Z;\]
thus, the image of $\varpi_1=\varepsilon_1$ in $\Omega/\langle \Phi
\rangle_\Z \cong \Z/n\Z$ generates that group.
\end{exmp}

%

The above example, in which $P_\fh(V)$ is a single $W$-orbit, is a special
case of a general construction that we will develop in the following two
sections. Of course, there are $\cL$-modules for which the set of weights 
is not a single $W$-orbit. For example, consider $\cL$ itself as a
$\cL$-module via the adjoint representation. Then $P_\fh(\cL)=\Phi\cup
\{\underline{0}\}$ and so there are at least two $W$-orbits. 
Furthermore, if $\cL$ is simple and not of simply laced type, then
there are long roots and short roots (see Remark~\ref{rellen}), and 
these cannot all lie in the same $W$-orbit.

\begin{rem} \label{factorA} Let $V$ be a faithful $\fh$-diagonalisable
$\fg$-module. By Proposition~\ref{defomega}, we have $\Phi\subseteq 
\Lambda(V)$ and $\Lambda(V)$ is a lattice in $E$. Let $\{\mu_i\mid i\in 
I\}$ be a basis of $E$ such that $\Lambda(V)=\sum_{i \in I} \Z \mu_i$; in 
particular, $\mu_i\in \Lambda(V) \subseteq \Omega$ for $i \in I$. Since 
$\Phi\subseteq \Lambda(V)$ we can write $\alpha_i=\sum_{j \in I} b_{ij}
\mu_j$ for each $i\in I$, where $b_{ij}\in \Z$. Thus, we obtain a matrix 
$B:=(b_{ij})_{i,j\in I}\in M_I (\Z)$. We also define the matrix 
$\breve{B}:=(\mu_j(h_i))_{i,j\in I}\in M_I(\Z)$. We deduce that
\[ a_{ij}=\alpha_j(h_i)=\sum_{l \in I} b_{jl}\mu_l(h_i)=(\breve{B}\cdot
B^{\text{tr}})_{ij} \quad\mbox{for $i,j\in I$}.\]
Thus, the lattice $\Lambda(V)\subseteq \Omega$ gives rise to a 
factorisation of the structure matrix $A=\breve{B} \cdot B^{\text{tr}}$. 
Conversely, if we are given any factorisation $A=\breve{B} \cdot 
B^{\text{tr}}$ where $B=(b_{ij})_{i,j \in I}$ and $\breve{B}=
(\breve{b}_{ij})_{i,j\in I}$ are matrices with entries in~$\Z$, then we set
\[\Lambda':=\sum_{j \in I}\Z\mu_j\subseteq E \quad \mbox{where}\quad\mu_j
:=\sum_{i\in I}\breve{b}_{ij}\varpi_i\in\Omega\quad\mbox{for $j\in I$}.\]
Since $\det(A)\neq 0$, we also have $\det(\breve{B})\neq 0$ and so 
$\Lambda'$ is a lattice in $E$. We have $\mu_j(h_i)=\breve{b}_{ij}$ for all
$i,j\in I$. Furthermore, for $i\in I$ we have
\begin{align*}
\sum_{j \in I} b_{ij}\mu_j&=\sum_{j,l\in I} b_{ij}\breve{b}_{lj}\varpi_l=
\sum_{l\in I} \Bigl(\sum_{j\in I} b_{ij}\breve{b}_{lj}\Bigr)\varpi_l\\
&= \sum_{l\in I} (\breve{B}\cdot B^{\text{tr}})_{li}\varpi_l=
\sum_{l \in I} a_{li}\varpi_l=\alpha_i.
\end{align*}
Thus, $\Lambda'\subseteq E$ is a lattice such that $\Phi\subseteq \Lambda'$.
This correspondence between lattices in $E$ containing $\Phi$ on the
one hand, and factorisations of~$A$ as a product of two (square) integer 
matrices on the other hand, is a special case of what happens for 
``{\itshape root data}'' in the theory of reductive algebraic groups; 
see, e.g., \cite[Remark~1.2.13]{GeMa}.
\end{rem}

We note the following extension of the above discussion, which will be 
useful later on in Section~\ref{secdiagII}. It also plays an important role 
in the discussion of ``Chevalley groups over~$\Z$''.

\begin{rem} \label{factorAa} For a faithful $\fh$-diagonalisable
$\cL$-module $V$, we set 
\[\fh_\Z:=\{x\in \fh\mid \mu(x) \in \Z\mbox{ for all $\mu \in 
P_\fh(V)$}\};\]
this is an additive subgroup of $\fh$. As above let $\{\mu_j \mid j 
\in I\}$ be a basis of $E$ such that $\Lambda(V)= \sum_{j \in I} \Z \mu_j$. 
Write $\mu_j= \sum_{j \in I} \breve{b}_{ij} \varpi_i$ for $j\in I$, where 
$\breve{b}_{ij} \in \Z$. Since the matrix $(\breve{b}_{ij})_{i,j\in I}$ has 
a non-zero determinant, there exist numbers $r_{ij}\in \Q$ such that, for any 
$i,j\in I$, we have
\[ \sum_{l \in I} r_{il}\breve{b}_{lj}=\delta_{ij}:=\left\{\begin{array}{cl} 
1 & \qquad \mbox{if $i=j$},\\ 0 & \qquad \mbox{if $i\neq j$}.\end{array}
\right.\] 
For $i\in I$ we set $h_i':=\sum_{j \in I} r_{ij}h_j\in \fh$. We claim that
\[ \fh_\Z=\sum_{i \in I} \Z h_i' \qquad \mbox{and}\qquad \mu_j(h_i')=
\delta_{ij} \quad \mbox{for $i,j\in I$}.\]
This is seen as follows. The identity $\mu_j(h_i')=\delta_{ij}$
immediately follows from the definition of $h_i'$ and the fact that
$\mu_j(h_i)=\breve{b}_{ij}$. Consequently, we have $\mu(h_i')\in \Z$ for 
all $\mu \in P_\fh(V)$. Thus, $h_i'\in \fh_\Z$ for $i\in I$. Conversely, 
let $x\in \fh_\Z$. Since the matrix $(r_{ij})$ has a non-zero determinant, 
the elements $\{h_i'\mid i \in I\}$ form a basis of~$\fh$. So we can
write $x=\sum_{i \in I} x_ih_i'$ where $x_i\in \C$. We must show that
$x_i\in \Z$ for all $i\in I$. Now $\mu_j(x)=\sum_{i \in I} x_i\mu_j(h_i')=
x_j$. Hence, since $\mu_j(x)\in \Z$ for all $j\in I$, we also have 
$x_j\in \Z$ for $j\in I$, as required. In particular, $h_i=\sum_{j \in I}
\breve{b}_{ij}h_j'$ for all $i\in I$.
\end{rem}

\begin{xca} \label{factorAxca} Let $\fg=\slm_n(\C)$ and regard $V=\fg$ 
as a $\fg$-module via the adjoint representation. Explicitly determine 
$\fh_\Z$ in this case, as a subset of $\fh=\{\mbox{diagonal matrices 
in $\fg$}\}$.
\end{xca}



\begin{xca} \label{hwmgom} Let $\cL=\gom_n(Q_n,\C)$ where $n\geq 4$; if
$Q_n=Q_n^{\text{tr}}$ and $n$ is even, also assume that $n\geq 6$. Let 
$\fh\subseteq\cL$ be the subalgebra consisting of diagonal matrices, as 
in Section~\ref{sec1a5}. Let $V=\C^n$ be the natural $\cL$-module; again,
this is a faithful module. Let $m\geq 2$ be such that $n=2m+1$ (if $n$ 
is odd) or $n=2m$ (if $n$ is even). Then $(\cL,\fh)$ is of Cartan--Killing 
type where $I=\{1,\ldots,m\}$ and the structure matrix~$A$ is of type
$B_m$, $C_m$ or $D_m$; see Proposition~\ref{CKbcd}. Show that~$V$ is 
$\fh$-diagonalisable and that 
\[ V=V_{\underline{0}}\oplus V_{\varepsilon_1}\oplus \ldots
\oplus V_{\varepsilon_m}\oplus V_{-\varepsilon_1}\oplus \ldots \oplus
V_{-\varepsilon_m},\]
where $V_{\underline{0}}=\{0\}$ if $n=2m$, and $V_{\underline{0}}\neq 
\{0\}$ if $n=2m+1$; here, $\varepsilon_i\in \fh^*$ is again the map that 
sends a diagonal matrix to its $i$-th diagonal entry. Thus, the weight
lattice of $V$ is given by 
\[ \Lambda(V)=\langle \varepsilon_1,\ldots,\varepsilon_m\rangle_\Z.\]
Now $\alpha_1, \ldots,\alpha_m$ are explicitly described as linear 
combinations of $\varepsilon_1,\ldots,\varepsilon_m$ in Remark~\ref{sec164}. 
Use this to show that 
\begin{alignat*}{2} 
\Lambda(V) &=\langle\Phi\rangle_\Z&&\qquad \mbox{if $n$ is odd 
(type $B_m$)},\\ \Lambda(V) /\langle\Phi\rangle_\Z&\cong \Z/2\Z 
&&\qquad \mbox{if $n$ is even (type $C_m$ or $D_m$)}.
\end{alignat*}
Comparing with Table~\ref{fundtab}, we see that $\langle\Phi\rangle_\Z
\subsetneqq \Lambda(V)=\Omega$ for type~$C_m$; furthermore, $\langle
\Phi\rangle_\Z=\Lambda(V) \subsetneqq \Omega$ for $B_m$ and $\langle\Phi
\rangle_\Z\subsetneqq \Lambda(V) \subsetneqq \Omega$ for~$D_m$. Show that 
the fundamental weights are given as follows.\\
(a) Assume that $n=2m+1$ and $Q_n^{\text{tr}}=Q_n$ (type $B_m$). Then 
\[ \textstyle \varpi_1=\frac{1}{2}(\varepsilon_1{+}\ldots{+}\varepsilon_m),
\quad \varpi_i = \varepsilon_1{+} \ldots {+} \varepsilon_{m+1-i}\;\;
\mbox{for $2\leq i\leq m$}.\] 
Show that $\Omega=\langle \Phi,\varpi_1\rangle_\Z$ in this case.\\
(b) Assume that $n=2m$ and $Q_n^{\text{tr}}=-Q_n$ (type $C_m$). Then 
\[ \varpi_{i}=\varepsilon_1{+}\ldots {+} \varepsilon_{m+1-i}\quad
\mbox{for $1\leq i\leq m$}.\]
Show that $\Omega=\langle \Phi,\varpi_m\rangle_\Z$ in this case.\\
(c) Assume that $n=2m$ and $Q_n^{\text{tr}}=Q_n$ (type $D_m$). Then 
\[\varpi_1=\textstyle\frac{1}{2}(\varepsilon_1{+}\ldots{+}\varepsilon_{m-1}{+}
\varepsilon_m), \quad \varpi_2=\frac{1}{2}(\varepsilon_1{+} \ldots {+}
\varepsilon_{m-1}{-}\varepsilon_m)\]
and $\;\varpi_i= \varepsilon_1{+} \ldots {+} \varepsilon_{m+1-i}\;$
for $3\leq i\leq m$. Show that 
\[ \Lambda(V)=\langle \Phi,\varpi_m\rangle_\Z\qquad\mbox{and}\qquad 
\Omega=\langle \Phi,\varpi_1,\varpi_2,\varpi_m\rangle_\Z.\]
Note that $\varpi_1-\varpi_2+\varpi_m=\varepsilon_1+\varepsilon_{m}\in
\Phi$; see Proposition~\ref{sec163}(a). Furthermore, if $m$ is even,
then the images of $\varpi_1$, $\varpi_2$ and $\varpi_m$ in $\Omega/
\langle \Phi\rangle_\Z$ are distinct and have order~$2$ each; if $m$ is 
odd, then the images of $\varpi_1$ and $-\varpi_2$ in $\Omega/\langle 
\Phi\rangle_\Z$ are equal and have order~$4$. \\
\noindent {\footnotesize [{\it Hint}. Argue as in Example~\ref{defhwm2}.
The elements $h_i\in\cH$ are described in the proof of 
Proposition~\ref{CKbcd}. See also \cite[Planche II, III, IV]{B} and
\cite[\S 3.6]{Sam}.]}
\end{xca}

\begin{rem} \label{singlegom} Let $\cL=\gom_n(Q_n,\C)$ and $V=\C^n$ be the
natural module, as in the above exercise. We note that, in each case, we have 
\[ \varpi_m=\varepsilon_1\in P_\fh(V) \qquad \mbox{(where $n=2m$ or 
$n=2m+1$)}.\]
Now assume that $n=2m\geq 4$; if $Q_n=Q_n^{\text{tr}}$, we also assume 
that $n\geq 6$. Then $P_\fh(V)= \{\pm \varepsilon_1, \ldots,\pm
\varepsilon_m\}$. Again, it follows from Proposition~\ref{sec168} and 
Remark~\ref{weylD} that $P_\fh(V)$ is a single $W$-orbit. (If $n=2m+1$, then 
$\underline{0} \in P_\fh(V)$ and so $P_\fh(V)$ is not a single $W$-orbit.)
\end{rem}


\begin{xca} \label{hwmg2} Assume that the structure matrix $A$ is of 
type~$G_2$. Consider the realisation of $\cL$ as a subalgebra of $\gl_7(\C)$, 
as in Exercise~\ref{otherG2}. Thus, $V=\C^7$ is a $\cL$-module; again,
$V$ is faithful. Describe the weights of $\fh$ on~$V$. Explicitly verify 
that $\Lambda(V)=\langle\Phi\rangle_\Z$ in this case. (This would also 
follow from Table~\ref{fundtab}.)
\end{xca}

\begin{xca} \label{dirsumfaith} Let $V$ be an $\fh$-diagonalizable 
$\cL$-module such that $V=V_1\oplus V_2$ where $V_1,V_2\subseteq V$
are $\cL$-submodules; note that $V_1,V_2$ are also $\fh$-diagonalizable. 
Show that $P_\fh(V)=P_\fh(V_1)\cup P_\fh(V_2)$ and, hence, that $\Lambda(V)= 
\Lambda(V_1)+\Lambda(V_2)$.
\end{xca}

\begin{xca} \label{tensfaith} Let $V_1$ and $V_2$ be (finite-dimensional) 
$\cL$-modules. By Remark~\ref{tenslie}, the tensor product $V:=V_1\otimes 
V_2$ is a $\cL$-module. Assume now that $V_1$ and $V_2$ are 
$\fh$-diagonalisable. Show that $V$ is $\fh$-diagonalisable and that
\[P_\fh(V)=\{\lambda+\mu\mid\lambda\in P_\fh(V_1) \mbox{ and } \mu\in 
P_\fh(V_2)\}.\]
{\footnotesize [{\it Hint}. Let $n=\dim V_1$ and $m=\dim V_2$. Let 
$\{v_1,\ldots, v_n\}$ be a basis of $V_1$ and $\{w_1,\ldots, w_m\}$ be a 
basis of $V_2$, where each $v_i$ is a weight vector (of weight 
$\lambda_i$, say) and each $w_j$ is a weight vector (of weight $\mu_j$,
say).]}
\end{xca}



In the following two sections, we further study the relation between 
$P=\langle\Phi\rangle_\Z$, $\Lambda(V)$ and $\Omega$, where the principal 
aim is to show: 
\begin{center}
\fbox{\begin{tabular}{l} For \textit{every} subgroup $\Lambda'\subseteq 
\Omega$ such that $\Phi\subseteq \Lambda'$, there exists a\\faithful,
$\fh$-diagonalisable $\cL$-module $V$ such that $\Lambda(V)=\Lambda'$.
\end{tabular}}
\end{center}
For example, if $\cL$ is of type $E_6$ or $E_7$ or $B_m$ ($m\geq 3$) or
$D_m$ ($m\geq 4$), then we have not yet seen a $\cL$-module $V$ with 
$\Lambda(V)=\Omega$. For $\cL$ of type $A_n$ ($n\geq 3$), we have not yet 
seen $\cL$-modules $V$ with $\langle\Phi\rangle_\Z\subsetneqq\Lambda(V)
\subsetneqq\Omega$.

\section{Minuscule weights} \label{sec5a1}

We keep the basic setting of the previous section. Our next aim is to 
find a natural set $\cM$ of representatives for the cosets of $\langle \Phi
\rangle_\Z$ in~$\Omega$; furthermore, for each $0\neq \lambda\in \cM$, we 
will construct a $\cL$-module whose weights are precisely the $W$-orbit
of~$\lambda$.  We will see that all this has a very elegant solution. Recall 
that $E=\langle \alpha_i\mid i \in I\rangle_\R \subseteq \cH^*$ and that 
\[ \Phi\subseteq \Omega=\langle \varpi_i \mid i \in I\rangle_\Z 
\subseteq E\qquad\mbox{where}\qquad |\Omega/\langle\Phi\rangle_\Z|<\infty.\] 
It will be convenient to fix a $W$-invariant positive-definite scalar 
product $\langle \; ,\; \rangle \colon E\times E \rightarrow\R$ 
(see Remark~\ref{decabst1a}). For $0 \neq v \in E$ we denote 
$v^\vee:=2v/\langle v,v\rangle\in E$. By Lemma~\ref{astring2} we have
\[\lambda(h_\alpha)=\langle \alpha^\vee,\lambda\rangle \qquad\mbox{for
$\alpha\in\Phi$ and $\lambda\in E$}.\]

\begin{defn} \label{weightorder} For $\lambda, \mu\in E$, we write
$\mu\preceq \lambda$ if $\lambda- \mu$ is a finite sum (possibly empty) 
of elements of $\Phi^+$. Thus, if $\mu\preceq\lambda$, then 
\[ \lambda-\mu=\sum_{\alpha\in \Phi^+} m_\alpha\alpha\qquad\mbox{where
$m_\alpha\in \Z_{\geq 0}$ for all $\alpha \in \Phi^+$}.\]
Consequently, we have: $\;\mu\preceq \lambda \Leftrightarrow \lambda-\mu=
\sum_{i\in I} n_i\alpha_i$, where $n_i\in \Z_{\geq 0}$ for 
all $i\in I$. We leave it as an exercise for the reader to check that 
$\preceq$ is a partial order on~$E$. We call $\preceq$ the \nm{weight 
order relation}. 
\end{defn}

\begin{rem} \label{cosetdom} Let $D$ be a coset of $\langle\Phi\rangle_\Z$
in $\Omega$. For $\lambda \in D$ and $i\in I$, we have $s_i(\lambda)=\lambda-
\lambda(h_i)\alpha_i\in D$, and so $w(D)=D$ for all $w \in W$. Thus, 
$D$ is a union of $W$-orbits of weights. Furthermore, note that for any 
$\lambda \in D$, we have $\{\mu \in \Omega\mid \mu \preceq \lambda\}
\subseteq D$ (by the definition of $\preceq$). Thus, the cosets of 
$\langle\Phi\rangle_\Z$ in $\Omega$ behave well with respect to the action 
of~$W$ and with respect to the weight order relation $\preceq$.
\end{rem}

\begin{defn} \label{defdomw} Let $\lambda\in\Omega$. Then we say that 
$\lambda \in \Omega$ is a \nm{dominant weight} if $\lambda(h_i)=\langle 
\alpha_i^\vee,\lambda\rangle\in \Z_{\geq 0}$ for all $i\in I$. Let 
\[\Omega^+:=\{\lambda \in \Omega\mid \lambda \mbox{ is dominant}\}.\] 
This set is non-empty; for example, we have $\underline{0}\in \Omega^+$ 
and $\varpi_i\in \Omega^+$ for all $i\in I$. Note that each $\lambda \in
\Omega^+$ has an expression 
\[ \lambda=\sum_{i\in I} m_i\varpi_i \qquad\mbox{where $m_i=
\lambda(h_i)\in\Z_{\geq 0}$ for all $i\in I$}.\]
By Lemma~\ref{astring3}, we then also have $\lambda(h_\alpha)\in 
\Z_{\geq 0}$ for all $\alpha\in \Phi^+$. 
\end{defn}

\begin{prop} \label{worbdom} Let $\lambda\in \Omega$. Then the following
hold.\\
{\rm (a)} If $\lambda \in \Omega^+$, then $w(\lambda)\preceq
\lambda$ for all $w\in W$.\\
{\rm (b)} The orbit $\{w(\lambda)\mid w\in W\}\subseteq \Omega$
contains a unique $\lambda_0 \in \Omega^+$.
\end{prop}

\begin{proof} (a) We show by induction on $\ell(w)$ that $w(\lambda)
\preceq \lambda$. If $\ell(w)=0$, then $w=1$ and there is nothing to 
prove. Now let $\ell(w)\geq 1$. By Corollary~\ref{basiclength1}, 
we can write $w=w's_i$ where $\ell(w')=\ell(w)-1$ and $i\in I$ is such 
that $w(\alpha_i)\in \Phi^-$. Now
\[w(\lambda)=w's_i(\lambda)=w'(\lambda)-\langle\alpha_i^\vee,
\lambda\rangle w'(\alpha_i),\]
where $\langle \alpha_i^\vee,\lambda \rangle\geq 0$ since $\lambda
\in\Omega^+$. Furthermore, $w'=ws_i$ and so $w'(\alpha_i)=-w(\alpha_i)\in 
\Phi^+$. Hence, $w(\lambda)\preceq w'(\lambda)$. By induction, we also
have $w'(\lambda)\preceq \lambda$ and so $w(\lambda)\preceq w'(\lambda) 
\preceq \lambda$, as required.

(b) In order to prove the existence of $\lambda_0$, we set $\psi:=
\sum_{i\in I} \varpi_i\in E$. Then $\langle \alpha_i^\vee, \psi\rangle=
\psi(h_i)=1$ and so $s_i(\psi)=\psi-\alpha_i$ for all $i\in I$. 
Now choose $w\in W$ such that $\langle w(\lambda),\psi\rangle\in\R$ is as 
large as possible. (This exists since $|W|<\infty$.) Let $i\in I$. Then 
$\langle w(\lambda), \psi\rangle\geq \langle (s_iw)(\lambda),\psi\rangle$. 
Hence, using the $W$-invariance of the scalar product, we obtain:
\begin{align*}
\langle w(\lambda),\psi \rangle& \geq \langle (s_iw)(\lambda),\psi\rangle
=\langle w(\lambda),s_i(\psi)\rangle=\langle w(\lambda),\psi-\alpha_i
\rangle\\ &=\langle w(\lambda),\psi\rangle- \langle w(\lambda),\alpha_i
\rangle
\end{align*}
and so $\langle w(\lambda),\alpha_i\rangle\geq 0$. But then also
$\langle \alpha_i^\vee,w(\lambda)\rangle=\langle w(\lambda),
\alpha_i^\vee\rangle\geq 0$ for all $i\in I$. So $\lambda_0:=w(\lambda)\in 
\Omega^+$. Now let $w_1,w_2\in W$ be such that $\lambda_1:=w_1(\lambda)$ 
and $\lambda_2:=w_2(\lambda)$ are dominant. By (a), we have $\lambda_2=
w_2w_1^{-1}(\lambda_1) \preceq \lambda_1$ and $\lambda_1=w_1w_2^{-1}
(\lambda_2)\preceq \lambda_2$; hence, $\lambda_1= \lambda_2$.
\end{proof}

\begin{lem} \label{lemstem0} Let $\lambda\in\Omega^+$. Then $\lambda=
\sum_{j\in I} u_j\alpha_j$, where $u_j\in \R_{\geq 0}$ for all $j \in I$. 
Furthermore, $\{\mu\in\Omega^+ \mid \mu\preceq\lambda\}$ is a finite set.
\end{lem}

\begin{proof} Since $\{\alpha_i\mid i\in I\}$ is a basis of~$E$, we 
can write $\lambda=\sum_{j\in I} u_j\alpha_j$ where $u_j\in\R$ for 
all~$j$. Since $\lambda$ is dominant, $0\leq\langle\alpha_i^\vee,\lambda
\rangle=\sum_{j \in I} a_{ij}u_j$ for all $i\in I$. Hence, if $u=
(u_i)_{i\in I}\in\R^I$, then $Au\geq 0$ (using the notation in
Section~\ref{sec3a1}). Now arrange the rows and columns of $A$ such that
$A$ is a block diagonal matrix, where all diagonal blocks are 
indecomposable. Since $W$ is finite, we deduce from Remark~\ref{decabst} 
and Lemma~\ref{decabst1} that each diagonal block of $A$ is of type (FIN). 
Hence, by Theorem~\ref{kac2}, we have $u\geq 0$, as claimed. 

Now let $\mu\in\Omega^+$ be such that $\mu\preceq \lambda$. By the same 
argument as above, we can write $\mu=\sum_{j \in I} v_j\alpha_j$ where 
$v_j\in \R_{\geq 0}$ for all $j\in I$. On the other hand, we have $\lambda-
\mu= \sum_{j\in I} n_j\alpha_j$, where $n_j\in\Z_{\geq 0}$ for all~$j\in I$.
Hence, we conclude that $0\leq n_j= u_j-v_j\leq u_j$ for all~$j$. So there 
are only finitely many possibilities for the $n_j$.
\end{proof} 

\begin{defn}[Humphreys \protect{\cite[Exc.~13.13]{H}}] \label{defcalm} 
We let $\cM$ denote the set of all minimal elements of $\Omega^+$, 
that is, the set of all $\lambda \in \Omega^+$ for which there exists 
no $\mu\in\Omega^+$ such that $\mu \preceq \lambda$, $\mu\neq \lambda$. 

For example, $\underline{0}\in\cM$. (Indeed, if $\mu\in\Omega^+$ and 
$\mu \preceq \underline{0}$, then $\mu=-\sum_{i \in I} n_i \alpha_i$ 
where $n_i\geq 0$ for all $i$, by the definition of $\preceq$. But, by
Remark~\ref{lemstem0}, we also have $n_i\leq 0$ for all $i$, and so 
$\mu=\underline{0}$.) 
\end{defn}

We will show below that $\cM$ is the desired set of coset representatives 
of $\langle \Phi\rangle_\Z$ in $\Omega$, and we will determine the set $\cM$
explicitly. 

\begin{lem}[Stembridge] \label{lemstem1} Let $\lambda,\mu\in \Omega^+$ be
such that $\lambda-\mu\in\langle\Phi\rangle_\Z$, that is, $\lambda$ and 
$\mu$ belong to the same coset of $\langle\Phi\rangle_\Z$ in $\Omega$. 
Then there exists some $\nu\in\Omega^+$ such that $\nu\preceq\lambda$ 
and $\nu\preceq \mu$.
\end{lem}

\begin{proof} Write $\lambda=\sum_{j\in I} a_j\alpha_j$ and $\mu=
\sum_{j\in I} b_j\alpha_j$ where $a_j,b_j\in \R$ for all~$j$. Since 
$\lambda-\mu\in\langle\Phi\rangle_\Z$, we have $a_j-b_j\in \Z$ for 
all~$j$. Now set $c_j:=\min\{a_j,b_j\}$ for all~$j\in I$ and define 
$\nu:=\sum_{j\in I} c_j\alpha_j\in E$. First note that $a_j-c_j\in 
\Z_{\geq 0}$ for all~$j$ and so $\lambda-\nu\in\langle\Phi\rangle_\Z$. 
Hence, $\nu\in\Omega$; furthermore, $\nu\preceq \lambda$ and 
$\nu\preceq\mu$. So it remains to show that $\nu\in\Omega^+$. Let 
$i\in I$. Now $\lambda\in\Omega^+$ and so $0\leq\langle\alpha_i^\vee,
\lambda\rangle=\sum_{j\in I} \langle \alpha_i^\vee, \alpha_j\rangle
a_j=\sum_{j\in I} a_{ij}a_j$. Since $a_{ii}=2$ and $a_{ij}\leq 0$ for 
$i\neq j$, we obtain that 
\[ 2a_i\geq \sum_{j\in I,j\neq i} (-a_{ij})a_j\geq \sum_{j\in I,j\neq i} 
(-a_{ij})c_j.\]
Similarly, since $\mu\in\Omega^+$, we obtain $2b_i \geq 
\sum_{j\in I,j\neq i} (-a_{ij})c_j$. But then 
\[ 2c_i=2\min\{a_i,b_i\}\geq \sum_{j\in I,j\neq i} (-a_{ij})c_j\]
and so $\langle\alpha_i^\vee,\nu\rangle=2c_i+\sum_{j\in I,j\neq i} 
a_{ij}c_j \geq 0$. Thus, $\nu$ is dominant.
\end{proof}

\begin{thm} \label{thmstem3} Let $D$ be a coset of $\langle\Phi\rangle_\Z$ 
in $\Omega$. Then $D$ contains a unique element of $\cM$. Consequently, we
have $|\cM|=|\Omega/\langle\Phi\rangle_\Z| <\infty$. Furthermore, if 
$\lambda\in\cM$, then $\langle \alpha^\vee, \lambda\rangle \in 
\{0,\pm 1\}$ for all $\alpha \in \Phi$.
\end{thm}

\begin{proof} We start with any $\lambda\in D$. By 
Proposition~\ref{worbdom}(b), there exists some $w\in W$ such that 
$\lambda':=w(\lambda) \in \Omega^+$. By Remark~\ref{cosetdom}, we have 
$\lambda'\in D$ and the set $\Omega':=\{\mu \in \Omega^+ \mid \mu 
\preceq \lambda'\}$ is contained in~$D$; furthermore, $\Omega'$ is 
finite by Remark~\ref{lemstem0}. So we can just pick an element 
$\lambda_0\in \Omega'$ that is minimal with respect to~$\preceq$; then 
$\lambda_0 \in D\cap\cM$ and so $D$ contains at least some element 
of~$\cM$. If we also have $\lambda_0'\in D\cap\cM$, then 
Lemma~\ref{lemstem1} shows that there is some $\nu\in \Omega^+$ such 
that $\nu \preceq \lambda_0$, $\nu \preceq \lambda_0'$. Since 
$\lambda_0, \lambda_0'$ are minimal, $\lambda_0=\nu= \lambda_0'$.

Now let $\lambda\in\cM$. Assume, if possible, that there exists some
$\alpha\in \Phi$ such that $\langle \alpha^\vee,\lambda\rangle \not\in 
\{0,\pm 1\}$. Replacing $\alpha$ by $-\alpha$ if necessary, we can assume
that $m:=\langle \alpha^\vee,\lambda\rangle >1$. We have $\lambda-
\alpha\in \Omega$. So, by Proposition~\ref{worbdom}(b), there exists 
some $w\in W$ such that $\lambda_1:=w(\lambda-\alpha)\in \Omega^+$.
Using the $W$-invariance of $\langle\;, \;\rangle$, we obtain:
\begin{align*}
\langle\lambda_1,\lambda_1\rangle&=\langle \lambda-\alpha,\lambda-\alpha
\rangle=\langle\lambda,\lambda \rangle-2\langle\alpha,\lambda\rangle+
\langle \alpha,\alpha\rangle\\ &=\langle \lambda,\lambda\rangle-m\langle
\alpha,\alpha  \rangle+ \langle \alpha,\alpha\rangle <
\langle \lambda,\lambda\rangle \quad\mbox{(since $m>1$)}.
\end{align*}
Now $\lambda$, $\lambda-\alpha$ and $\lambda_1$ all belong to the same 
coset of $\langle\Phi\rangle_\Z$ in $\Omega$ (see Remark~\ref{cosetdom}). 
Hence, by Lemma~\ref{lemstem1}, there exists some $\nu \in \Omega^+$ such 
that $\nu \preceq \lambda_1$ and $\nu\preceq \lambda$. But $\lambda\in\cM$
and so $\lambda=\nu\preceq \lambda_1$. Since $\lambda \in \Omega^+$ 
and $\lambda_1 \in \Omega^+$, we have $\lambda+\lambda_1\in\Omega^+$ 
and so $\lambda+\lambda_1=\sum_{j \in I} m_j\varpi_j$, where $m_j
\in \Z_{\geq 0}$ for all $j\in I$. Since $\lambda\preceq\lambda_1$,
we also have $\lambda_1-\lambda=\sum_{i \in I} n_i\alpha_i$, where
$n_i\in \Z_{\geq 0}$ for all $i \in I$. Finally, since  
\[\langle\alpha_i,\varpi_j\rangle=\textstyle{\frac{1}{2}}
\langle \alpha_i,\alpha_i\rangle \langle \alpha_i^\vee,\varpi_j
\rangle \geq 0 \qquad \mbox{for all $i,j\in I$},\]
we conclude that 
\[\langle \lambda_1,\lambda_1\rangle-\langle \lambda,\lambda\rangle=
\langle \lambda_1-\lambda, \lambda_1+\lambda\rangle=\sum_{i,j\in I} 
n_im_j\langle \alpha_i,\varpi_j\rangle\geq 0,\]
contradiction. Hence, we do have $\langle\alpha^\vee,\lambda\rangle 
\in\{0,\pm 1\}$ for all $\alpha\in\Phi$.
\end{proof}

The further condition on the elements of $\cM$ in 
Theorem~\ref{thmstem3} leads to the following definition.

\begin{defn}[Cf.\ Bourbaki \protect{\cite[Ch.~VI, \S 1, Exc.~24]{B}}]
\label{minusc} We say that $\lambda\in \Omega$ is a \nm{minuscule weight}
if $\langle\alpha^\vee,\lambda\rangle\in\{0,\pm 1\}$ for all $\alpha
\in \Phi$. Clearly, $\lambda=\underline{0}$ is minuscule.
\end{defn}

\begin{rem}\label{minusc0} Let $\lambda\in\Omega$. By the 
$W$-invariance of $\langle\;,\;\rangle$, we have  
\[ w(\lambda^\vee)=\frac{2w(\lambda)}{\langle \lambda,\lambda \rangle}
=\frac{2w(\lambda)}{\langle w(\lambda),w(\lambda) \rangle}
=w(\lambda)^\vee \qquad\mbox{for all $w\in W$}.\]
Hence, if $\lambda\in\Omega$ is minuscule, then $\langle \alpha^\vee,
w^{-1}(\lambda)\rangle=\langle w(\alpha^\vee),\lambda\rangle=\langle 
w(\alpha)^\vee,\lambda\rangle \in \{0,\pm 1\}$ for all $w\in W$ (since 
$w(\alpha)\in\Phi$). So all weights in the orbit $\{w(\lambda)\mid 
w\in W\}$ are minuscule and that orbit contains a unique dominant minuscule 
weight (by Proposition~\ref{worbdom}). 
\end{rem}

\begin{lem} \label{lemstem4} Let $\lambda\in \Omega$ be minuscule. If
$\mu \in \Omega^+$ is such that $\mu \preceq \lambda$, then $\mu=
\lambda$. Consequently, if $\lambda$ is dominant, then $\lambda\in \cM$.
\end{lem}

\begin{proof} Let $\mu\in \Omega^+$. For any $\lambda\in \Omega$ such that 
$\mu\preceq\lambda$, we write $\lambda- \mu=\sum_{i \in I} n_i\alpha_i$, 
where $n_i\in \Z_{\geq 0}$ for $i\in I$, and set $n(\lambda):=\sum_{i\in I} 
n_i\geq 0$. Now let $\lambda \in \Omega$ be minuscule such that $\mu
\preceq \lambda$. Assume, if possible, that $\mu\neq \lambda$. Then 
$n(\lambda)>0$. If $n(\lambda)=1$, then $\lambda-\mu=\alpha_i$ for 
some $i\in I$. Since $\lambda$ is minuscule, we have $\langle \alpha_i^\vee,
\lambda\rangle \in \{0,\pm 1\}$ and so 
\[ \langle \alpha_i^\vee,\mu\rangle=\langle \alpha_i^\vee,
\lambda\rangle -\langle \alpha_i^\vee, \alpha_i \rangle=
\langle \alpha_i^\vee, \lambda\rangle-2\leq -1,\]
contradiction to $\mu\in \Omega^+$. Now let $n(\lambda)>1$.
Since $\mu\neq \lambda$, we have 
\[\sum_{i \in I} n_i \langle \alpha_i,\lambda-\mu\rangle=\langle \lambda-
\mu,\lambda-\mu \rangle>0.\]
Hence, there is some $j\in I$ such that $n_j>0$ and $\langle \alpha_j, 
\lambda-\mu\rangle>0$. Then we also have $\langle \alpha_j^\vee, \lambda-
\mu\rangle>0$ and so $\langle \alpha_j^\vee, \lambda\rangle> \langle
\alpha_j^\vee, \mu \rangle \geq 0$, since $\mu$ is dominant. Since 
$\lambda$ is minuscule, we must have $\langle \alpha_j^\vee, \lambda
\rangle=1$. But then $\lambda':= \lambda-\alpha_j=s_j(\lambda)\in \Omega$ 
is also minuscule; see Remark~\ref{minusc0}. Furthermore, since $n_j>0$, 
we have $\mu \preceq \lambda'$ and $n(\lambda')=n(\lambda)-1\geq 1$. If 
we still have $n(\lambda')>1$, then we repeat the argument with $\lambda'$ 
instead of~$\lambda$. Hence, there is some $j'\in I$ such that $\lambda''
:=\lambda'-\alpha_{j'}=s_{j'}(\lambda')\in \Omega$ is minuscule, 
$\mu\preceq \lambda''$ and $n(\lambda'')=n(\lambda')-1$. After 
finitely many repetitions, we find some minuscule $\tilde{\lambda}\in 
\Omega$ such that $\mu \preceq \tilde{\lambda}$ and $n(\tilde{\lambda})=1$.
But then we obtain a contradicton as above.
\end{proof}

\begin{cor} \label{corminus} Let $\lambda \in \Omega^+$. Then $\lambda 
\in \cM$ (see Definition~\ref{defcalm}) if and only if $\lambda$ is
minuscule (see Definition~\ref{minusc}). Thus, every coset of $\langle
\Phi\rangle_\Z$ in $\Omega$ contains a unique $W$-orbit of minuscule 
weights.
\end{cor}

\begin{proof} By Theorem~\ref{thmstem3}, we have ``$\lambda \in 
\cM\Rightarrow \lambda$ minuscule''. The reverse implication holds
by Lemma~\ref{lemstem4}. Now let $D\subseteq \Omega$ be a coset of 
$\langle\Phi\rangle_\Z$. There is a unique $\lambda \in D
\cap \cM$, and $\lambda$ is minuscule; see once more Theorem~\ref{thmstem3}. 
By Remark~\ref{cosetdom}, the whole $W$-orbit of $\lambda$ is contained
in~$D$. 
\end{proof}

\begin{rem} \label{corminus1} Let $\Lambda'\subseteq \Omega$ be an
arbitrary subgroup with $\langle \Phi\rangle_\Z\subseteq \Lambda'$. 
Then $\langle\Phi \rangle_\Z$ has finite index in $\Lambda'$; let 
$D_0,D_1,\ldots,D_m$ be the cosets of $\langle\Phi\rangle_\Z$ in 
$\Lambda'$, where $D_0=\langle \Phi\rangle_\Z$ and $D_l\neq \langle
\Phi\rangle_\Z$ for $1\leq l\leq m$. (Here, $m=0$ if $\Lambda'=
\langle \Phi\rangle_\Z$.) If $l\geq 1$, then $D_l$ contains 
a unique dominant minuscule weight $\lambda_l^\circ \in \cM\setminus
\{\underline{0}\}$; see Theorem~\ref{thmstem3}. We claim that 
\begin{equation*}
\Lambda'=\langle Q\rangle_\Z\qquad \mbox{where} \qquad
Q:=\Phi \cup \{\lambda_l^\circ \mid 1\leq l \leq m\}.\tag{a}
\end{equation*}
(Note that $Q$ is a finite set.) Indeed, the inclusion ``$\supseteq$'' 
is clear, since $\Phi\subseteq D_0\subseteq \Lambda'$ and $\lambda_l^\circ
\in D_l \subseteq \Lambda'$ for $1\leq l\leq m$. Conversely, let $\lambda 
\in \Lambda'$. Then $\lambda \in D_l$ for a unique $l\in \{0,1,\ldots,m\}$. 
If $l=0$, then $\lambda \in D_0=\langle\Phi\rangle_\Z\subseteq \langle Q
\rangle_\Z$. If $l\geq 1$, then $\lambda=\lambda_l^\circ+\lambda'$ where 
$\lambda'\in \langle\Phi\rangle_\Z\subseteq \langle Q\rangle_\Z$. Hence, 
since $\lambda_l^\circ \in Q$, we have $\lambda\in \langle Q\rangle_\Z$. 
Thus, the above claim is proved. Furthermore, for $l\geq 1$ let $\Psi_l$ be
the $W$-orbit of $\lambda_l^\circ$. Then $\Psi_l\subseteq D_l$ (see 
Remark~\ref{cosetdom}) and so we also have:
\begin{equation*}
\Lambda'=\big\langle \Phi\cup \Psi_1\cup \ldots \cup \Psi_m
\big\rangle_\Z.\tag{b}
\end{equation*}
In the next section, we shall see that one can always construct a 
$\cL$-module $V$ with $P_\fh(V)=\{\underline{0}\}\cup \Phi\cup
\Psi_1\cup\ldots\cup\Psi_m$ and, hence, $\Lambda'=\Lambda(V)$. 
\end{rem}

\begin{table}[htbp]
\caption{Non-zero dominant minuscule weights (marked by ``$\circ$'')}
\label{Mtabminus}
{\small \begin{center}
\begin{picture}(298,125)
\put(  5, 25){$E_6$}
\put( 30, 25){\circle{6}}
\put( 28, 31){$1$}
\put( 33, 25){\line(1,0){20}}
\put( 50, 25){\circle*{5}}
\put( 48, 31){$3$}
\put( 50, 25){\line(1,0){20}}
\put( 70, 25){\circle*{5}}
\put( 68, 31){$4$}
\put( 70, 25){\line(0,-1){20}}
\put( 70,  5){\circle*{5}}
\put( 75,  3){$2$}
\put( 70, 25){\line(1,0){20}}
\put( 90, 25){\circle*{5}}
\put( 88, 31){$5$}
\put( 90, 25){\line(1,0){17}}
\put(110, 25){\circle{6}}
\put(108, 31){$6$}

\put(170, 25){$E_7$}
\put(195, 25){\circle*{5}}
\put(193, 31){$1$}
\put(195, 25){\line(1,0){20}}
\put(215, 25){\circle*{5}}
\put(213, 31){$3$}
\put(215, 25){\line(1,0){20}}
\put(235, 25){\circle*{5}}
\put(233, 31){$4$}
\put(235, 25){\line(0,-1){20}}
\put(235,  5){\circle*{5}}
\put(240,  3){$2$}
\put(235, 25){\line(1,0){20}}
\put(255, 25){\circle*{5}}
\put(253, 31){$5$}
\put(255, 25){\line(1,0){20}}
\put(275, 25){\circle*{5}}
\put(273, 31){$6$}
\put(275, 25){\line(1,0){17}}
\put(295, 25){\circle{6}}
\put(293, 31){$7$}

\put(  3,60){$C_n$}
\put(  2,52){$\scriptstyle{n \geq 2}$}
\put( 30,55){\circle*{5}}
\put( 28,61){$1$}
\put( 30,57){\line(1,0){20}}
\put( 30,53){\line(1,0){20}}
\put( 36,53){$>$}
\put( 50,55){\circle*{5}}
\put( 48,61){$2$}
\put( 50,55){\line(1,0){20}}
\put( 70,55){\circle*{5}}
\put( 68,61){$3$}
\put( 70,55){\line(1,0){10}}
\put( 90,55){\circle*{1}}
\put(100,55){\circle*{1}}
\put(110,55){\circle*{1}}
\put(120,55){\line(1,0){27}}
\put(130,55){\circle*{5}}
\put(121,61){$n{-}1$}
\put(138,55){\line(1,0){9}}
\put(150,55){\circle{6}}
\put(148,61){$n$}

\put(  3,90){$B_n$}
\put(  2,82){$\scriptstyle{n \geq 2}$}
\put( 30,85){\circle{6}}
\put( 28,91){$1$}
\put( 32,83){\line(1,0){17}}
\put( 32,87){\line(1,0){17}}
\put( 36,83){$<$}
\put( 50,85){\circle*{5}}
\put( 48,91){$2$}
\put( 50,85){\line(1,0){30}}
\put( 70,85){\circle*{5}}
\put( 68,91){$3$}
\put( 90,85){\circle*{1}}
\put(100,85){\circle*{1}}
\put(110,85){\circle*{1}}
\put(120,85){\line(1,0){10}}
\put(130,85){\circle*{5}}
\put(121,91){$n{-}1$}
\put(130,85){\line(1,0){20}}
\put(150,85){\circle*{5}}
\put(148,91){$n$}

\put(  3,120){$A_{n}$}
\put(  2,112){$\scriptstyle{n \geq 1}$}
\put( 30,115){\circle{6}}
\put( 28,121){$1$}
\put( 33,115){\line(1,0){14}}
\put( 50,115){\circle{6}}
\put( 48,121){$2$}
\put( 53,115){\line(1,0){14}}
\put( 70,115){\circle{6}}
\put( 67,121){$3$}
\put( 73,115){\line(1,0){11}}
\put( 90,115){\circle*{1}}
\put(100,115){\circle*{1}}
\put(110,115){\circle*{1}}
\put(118,115){\line(1,0){9}}
\put(130,115){\circle{6}}
\put(121,121){$n{-}1$}
\put(133,115){\line(1,0){14}}
\put(150,115){\circle{6}}
\put(148,121){$n$}

\put(180,87){$D_n$}
\put(180,79){$\scriptstyle{n \geq 3}$}
\put(195,105){\circle{6}}
\put(200,105){$1$}
\put(195,65){\circle{6}}
\put(201,60){$2$}
\put(197,103){\line(1,-1){18}}
\put(197,67){\line(1,1){18}}
\put(215,85){\circle*{5}}
\put(213,91){$3$}
\put(218,85){\line(1,0){10}}
\put(235,85){\circle*{1}}
\put(245,85){\circle*{1}}
\put(255,85){\circle*{1}}
\put(263,85){\line(1,0){10}}
\put(275,85){\circle*{5}}
\put(266,91){$n{-}1$}
\put(275,85){\line(1,0){17}}
\put(295,85){\circle{6}}
\put(293,91){$n$}
\end{picture}
(For the types $G_2$, $F_4$ and $E_8$, the only minuscule weight 
is $\underline{0}$.)
\end{center}}
\end{table}

\begin{prop} \label{minusc1} Assume that $A$ is indecomposable. Then 
the non-zero dominant minuscule weights are $\{\varpi_i \mid i \in 
I^\circ\}$, where $I^\circ \subseteq I$ is the set of indices with vertex 
marked by ``$\,\circ$'' in Table~\ref{Mtabminus} (p.~\pageref{Mtabminus}). 
\end{prop}

\begin{proof} By Corollary~\ref{corminus}, we have $\cM=\{\underline{0}\}
\cup \cM'$, where $\cM'$ denotes the set of all non-zero dominant minuscule 
weights. First we show that $|\cM'|\leq |I^\circ|$; more precisely,
\[ \cM'\subseteq \{\varpi_i \mid i \in I^\circ\}.\]
This is seen as follows. Let $\lambda\in\cM'$. Since $\lambda$
is dominant, we have $\lambda=\sum_{j\in I} m_j\varpi_j$ where 
$m_j\in \Z_{\geq 0}$ for all $j\in I$. Since $\lambda$ is minuscule, we 
have $m_i=\langle \alpha_i^\vee,\lambda\rangle\in \{0,\pm 1\}$ for all
$i \in I$. So we can already conclude that $m_i\in\{0,1\}$ for 
all $i\in I$. Thus, $\lambda=\sum_{j\in I'} \varpi_j$ for some 
subset $I'\subseteq I$. We have $I'\neq \varnothing$ since $\lambda
\neq \underline{0}$. 

If $A$ is simply laced, let $\alpha_0\in\Phi^+$ be the highest root 
as in Remark~\ref{highestr}. Writing $\alpha_0=\sum_{i\in I} n_i
\alpha_i$ with $n_i\in \Z_{\geq 0}$, we then also have $\alpha_0^\vee=
\sum_{i\in I} n_i\alpha_i^\vee$; see Exercise~\ref{xcastdbase2}. This 
yields $\sum_{j\in I'} n_j=\langle\alpha_0^\vee,\lambda\rangle\in
\{0,1\}$. The coefficients $(n_j)_{j\in I}$ are listed in 
Table~\ref{highroot} (p.~\pageref{highroot}). By inspection, we see 
that $n_j\geq 1$ for all $j\in I$. Hence, we must have $|I'|=1$ and 
$n_j=1$ for the unique index $j\in I'$. A further inspection shows that 
that index~$j$ is one of those marked by ``$\circ$'' in 
Table~\ref{Mtabminus}. For example, if $A$ is of type~$E_7$, then 
$(n_i)_{i\in I}=(2,2,3,4,3,2,1)$; there is only one coefficient equal 
to~$1$, and this corresponds to the vertex marked by ``$\circ$'' in 
Table~\ref{Mtabminus}. Similarly, if $A$ is of type $E_8$, then 
$(n_i)_{i\in I}=(2,3,4,6,5,4,3,2)$; there is no coefficient equal 
to~$1$ and so there is no minuscule weight at all. (Note that, for
the purposes of this argument, we do not need to know that the roots 
$\alpha_0$ in Table~\ref{highroot} are really the highest roots; we 
just need to know that each $\alpha_0$ in that table is a root at all.)

If $A$ is not simply laced, then let $\alpha_0'\in\Phi^+$ be the highest 
short root; see Exercise~\ref{highshort}. Explicit expressions for
$\alpha_0'^\vee$ are given in Table~\ref{highroot}. (Again, we just 
need to know that $\alpha_0'$ is a root at all.) Then the same 
reasoning as above yields that $|I'|=1$, and the unique index in $I'$ is 
one of those marked by ``$\circ$'' in Table~\ref{Mtabminus}.

Finally, Theorem~\ref{thmstem3} shows that $|\cM|=|\Omega/\langle\Phi
\rangle_\Z|$. Comparing Table~\ref{fundtab} and Table~\ref{Mtabminus}, we 
observe that $|\Omega/\langle\Phi\rangle_\Z|=|I^\circ|+1$. Hence, since 
$\cM=\{\underline{0}\} \cup \cM'$, we have $|I^\circ|=|\cM'|$. But we have
seen above that $\cM'\subseteq \{\varpi_i \mid i \in I^\circ\}$. So this
inclusion must be an equality. 
\end{proof}


For each minuscule weight $\varpi_{i_0}$ ($i_0\in I^\circ$) as in 
Proposition~\ref{minusc1}, the size of the corresponding $W$-orbit 
in $\Omega$ is shown in Table~\ref{Morbits}. See the examples below 
for further explanations.

\begin{table}[htbp]
\caption{Orbits of minuscule weights} \label{Morbits}
\begin{center}
$\renewcommand{\arraystretch}{1.2} \begin{array}{ccl} \hline \mbox{Type} 
& |\Omega/\langle\Phi\rangle_\Z| & \mbox{Size of orbit of minuscule 
$\varpi_{i_0}$}\\ \hline
A_{n-1} \; (n\geq 2)& n& \binom{n}{i_0} \;\; (1\leq i_0 \leq n-1)\\ 
B_n\; (n\geq 2) & 2 & 2^n\;\; (i_0=1) \\ C_n\; (n\geq 2) & 2 & 2n\;\; 
(i_0=n) \\ D_n \; (n\geq 3)& 4 & 2^{n-1}\;\; (i_0=1,2), \quad 2n\;\; 
(i_0=n) \\ E_6 & 3 & 27 \;\; (i_0=1,6) \\ E_7 & 2 & 56 \;\; 
(i_0=7) \end{array}$
\end{center}
\end{table}

\begin{exmp} \label{minusAn} Assume that $\fg$ is of type $A_{n-1}$,
$n\geq 2$. By Example~\ref{defhwm2}, the fundamental weights are given by 
$\varpi_r=\varepsilon_1+\ldots+ \varepsilon_r$ for $1\leq r \leq n-1$.
By Table~\ref{Mtabminus}, they are all minuscule. One easily sees that 
the $W$-orbit of $\varpi_r$ consists of all weights of the form 
\[ \varepsilon_{i_1}+ \ldots +\varepsilon_{i_r}\qquad\mbox{where}
\qquad 1\leq i_1<\ldots <i_r\leq n.\]
Hence, the size of that orbit is $\binom{n}{r}$; see Table~\ref{Morbits}. 
In particular, for $r=1$, we have $\varpi_1=\varepsilon_1$ and the
$W$-orbit is $\{\varepsilon_1,\ldots,\varepsilon_n\}$. For further 
details see Bourbaki \cite[Ch.~VIII, \S 13, no.~1]{B78}.
\end{exmp}

\begin{exmp} \label{minusBn} Assume that $\fg$ is of type $B_n$, $n
\geq 2$. By Table~\ref{Mtabminus}, $\varpi_1$ is the only fundamental
weight that is minuscule. By Exercise~\ref{hwmgom}, we have
$\varpi_1=\textstyle{\frac{1}{2}}(\varepsilon_1+\ldots +\varepsilon_n)$.
One easily sees that the $W$-orbit of $\varpi_1$ consists of all weights 
of the form $\frac{1}{2}(\pm \varepsilon_1\pm\ldots\pm\varepsilon_n)$,
for any choice of the signs. Hence, the size of that orbit is $2^n$; see
Table~\ref{Morbits}. For further details see Bourbaki \cite[Ch.~VIII, 
\S 13, no.~2]{B78}. 
\end{exmp}

%
%

\begin{exmp} \label{minusCn} Assume that $\fg$ is of type $C_n$, $n
\geq 2$. By Table~\ref{Mtabminus}, $\varpi_n$ is the only fundamental
weight that is minuscule. We already noted in Remark~\ref{singlegom}
that $\varpi_n=\varepsilon_1$ and that the $W$-orbit of $\varpi_n$
consists of all weights of the form $\pm \varepsilon_i$ for $1\leq i 
\leq n$.  Hence, the size of that orbit is $2n$; see Table~\ref{Morbits}. 
For further details see Bourbaki \cite[Ch.~VIII, \S 13, no.~3]{B78}.
\end{exmp}

\begin{exmp} \label{minusDn} Assume that $\fg$ is of type $D_n$, $n
\geq 3$. By Table~\ref{Mtabminus}, $\varpi_1$, $\varpi_2$ and $\varpi_n$ are
the only fundamental weights that are minuscule. We already noted in 
Remark~\ref{singlegom} that $\varpi_n=\varepsilon_1$ and that the 
$W$-orbit of $\varpi_n$ consists of all weights of the form 
$\pm \varepsilon_i$ for $1\leq i \leq n$. Hence, the size of that orbit 
is $2n$; see Table~\ref{Morbits}. Now consider the weights $\varpi_1$ 
and $\varpi_2$. By  Exercise~\ref{hwmgom}, we have
\[ \varpi_1=\textstyle{\frac{1}{2}}(\varepsilon_1+\ldots+
\varepsilon_{n-1}+\varepsilon_n).\]
One checks that the $W$-orbit of $\varpi_1$ consists of all weights of
the form 
\[ \textstyle{\frac{1}{2}}(\pm \varepsilon_1\pm\ldots\pm\varepsilon_n)\]
where the number of minus signs is even. Thus, the size of that orbit 
is $2^{n-1}$; see Table~\ref{Morbits}. Similarly, the $W$-orbit of 
\[ \varpi_2=\textstyle{\frac{1}{2}}(\varepsilon_1+\ldots+
\varepsilon_{n-1}-\varepsilon_n)\]
consists of all weights as above, but where the number of minus 
signs is odd. Hence, again, the size of that orbit is $2^{n-1}$.
For further details see Bourbaki \cite[Ch.~VIII, \S 13, no.~4]{B78}.
\end{exmp}

For comments on the minuscule weights in types $E_6$ and $E_7$, see
Example~\ref{minusE67} in the section below.

%

\begin{xca} \label{xcamindec} Assume that $A$ is decomposable. Then, as 
in Remark~\ref{cindec5}, we have a partition $I=\bigsqcup_{s\in S} I_s$ 
(disjoint union), where $S$ is a finite index set, $I_s\neq \varnothing$ 
for $s\in S$ and $A$ is block diagonal with indecomposable diagonal blocks 
$A_s:=(a_{ij})_{i,j\in I_s}$ for $s\in S$. Since 
$|W(A)|<\infty$, we have $|W(A_s)|<\infty$ for all $s\in S$ and so each
$A_s$ is of type (FIN) (see Lemma~\ref{decabst1}). We have a direct sum
decomposition $E=\bigoplus_{s\in S} E_s$, where $E_s:=\langle \alpha_i
\mid i \in I_s\rangle_\C$ and these subspaces are orthogonal to 
each other with respect to $\langle \;,\; \rangle$. Show that 
\[ \cM=\Bigl\{\sum_{s\in S} \lambda_s \,\big|\, \lambda_s\in \cM_s
\mbox{ for all $s\in S$}\Bigr\},\]
where $\cM_s:=\{\underline{0}\}\cup \{\varpi_i\mid i\in I_s^\circ\}$ 
and $I_s^\circ$ is the set of indices $i\in I_s$ with vertex marked by 
``$\circ$'' in the Dynkin diagram of $A_s$ (as in Table~\ref{Mtabminus}). 
\end{xca}

\section{Highest weight modules} \label{sechwm}

There is an important class of $\cL$-modules $V$ for which the set
of weights $P_\fh(V)$ behaves particularly well with respect to the
weight order relation $\preceq$ on $\fh^*$ and the cosets of $P=\langle
\Phi\rangle_\Z$ in the weight lattice $\Omega$. These are the ``highest 
weight modules'' which play a key role in the representation theory of 
semisimple Lie algebras. Here, we will just establish some very basic 
results about these modules. We will also explicitly construct highest 
weight modules corresponding to minuscule weights as discussed in the 
previous section.


\begin{defn}[Cf.\ \protect{Remark~\ref{sl2highb}}] \label{defprimitive} For 
$i\in I$ let $\{e_i,f_i,h_i\} \subseteq \cL$ be the corresponding 
\nms{$\slm_2$-triple}{sl$_2$-triple}, as above. Let $V$ be a $\cL$-module.
Let $0\neq v^+ \in V$ and $\lambda \in \fh^*$. Then $v^+$ is called a 
\nm{primitive vector} of~$V$ (of weight $\lambda$) if $e_i.v^+=0$ for 
all $i\in I$ and $h.v^+=\lambda(h)v^+$ for all $h\in \fh$. In particular, 
$\lambda \in P_\fh(V)$ and $v^+\in V_\lambda$.
\end{defn}

\begin{lem} \label{primitiv0} For every $\cL$-module $V$, there exists a 
primitive vector $0\neq v^+\in V$. If $v^+$ is of weight $\lambda$,
then $\lambda(h_\alpha)\in \Z_{\geq 0}$ for $\alpha\in \Phi^+$.
\end{lem}

\begin{proof} As in Remark~\ref{defTDa}, let $\fn^+:=\sum_{\alpha\in 
\Phi^+}\cL_\alpha\subseteq\cL$, where $\Phi^+\subseteq \Phi$ are the 
positive roots with respect to $\Delta=\{\alpha_i\mid i \in I\}$. Then 
consider the subalgebra $\fb^+:=\fh+\fn^+ \subseteq \cL$ in 
Proposition~\ref{borelsub}. We have $[\fb^+,\fb^+]=\fn^+$ and $\fb^+$ is 
solvable. By restriction, we can regard $V$ as a $\fb^+$-module. Let
$n=\dim V$. Since $\fb^+$ is solvable, Lie's Theorem~\ref{liethm} shows 
that there is a basis $\{v_1,\ldots,v_n\}$ of $V$ such that, for any 
$x\in \fb^+$, the corresponding matrix of $\rho_x\colon V\rightarrow V$ 
is upper triangular with $\lambda_1(x),\ldots,\lambda_n(x)$ on the 
diagonal, where $\lambda_j\in \fb^+$ are such that $\fn^+=[\fb^+,\fb^+]
\subseteq \ker(\lambda_j)$ for all~$j$. Let $v^+:=v_1$ and $\lambda\in
\fh^*$ be the restriction of $\lambda_1$ to~$\fh$. Then $v^+\in
\cL_\lambda$ since $b.v^+=\lambda(b)v^+$ for all $b\in \fb^+$. 
Furthermore, $e_i.v^+=0$ for all $i\in I$, since $e_i\in \fn^+\subseteq 
\ker(\lambda_1)$. So $v^+$ is a primitive vector of $V$, of weight 
$\lambda$.

Now let $\alpha\in \Phi^+$. As in Remark~\ref{propnrs0}, we have a 
subalgebra $\cL_\alpha= \langle e_\alpha,f_{\alpha},h_\alpha \rangle_\C
\subseteq \cL$ isomorphic to $\slm_2(\C)$. By restriction, we can regard
$V$ as $\cL_\alpha$-module. Then $v^+$ also is a primitive vector for 
the $\cL_\alpha$-module $V$. So, by Remark~\ref{sl2highb}(b), we 
have $\lambda(h_\alpha)\in \Z_{\geq 0}$.
\end{proof}

%
%
%
%



\begin{prop} \label{primitive1} Let $V$ be a $\cL$-module and $0\neq v^+
\in V$ be a primitive vector, of weight $\lambda \in \fh^*$. Let 
$V'\subseteq V$ be the subspace spanned by all $v\in V$ of the form
\begin{equation*}
v=f_{i_1}.(f_{i_2}.(\ldots (f_{i_l}.v^+)\ldots)) \quad
\mbox{where $l\geq 0$, $\;i_1,\ldots,i_l\in I$}. \tag{$*$}
\end{equation*} 
Then $V'$ is an $\fh$-diagonalisable $\cL$-submodule of $V$, where 
$V_\lambda'=\langle v^+\rangle_\C$ and $P_\fh(V')\subseteq \{\mu \in \fh^*
\mid \mu \preceq \lambda\}\subseteq \lambda+\langle\Phi\rangle_\Z$.
\end{prop}

\begin{proof} Since $h.v^+=\lambda(h)v^+$ for all $h\in \fh$, 
Remark~\ref{fundcal} and an easy induction on~$l$ show that, for every 
$v$ of the form ($*$), we have:
\begin{equation*}
h.v=\mu(h)v \quad \mbox{for all $h \in \fh$, where $\mu=
\lambda-\alpha_{i_1}- \ldots-\alpha_{i_l}$}.\tag{$*^\prime$}
\end{equation*}
Next we show that $V'\subseteq V$ is a $\cL$-submodule. By 
Proposition~\ref{genlie} and Exercise~\ref{xcagenerator1}, it is 
sufficient to show that $e_i.V'\subseteq V'$ and $f_i.V' \subseteq
V'$ for all $i \in I$. By the definition of $V'$, we certainly have 
$f_i.V' \subseteq V'$. It remains to show that $e_i.v\in V'$ for $i\in I$ 
and~$v$ of the form~($*$). Again, we use induction on $l$. If $l=0$, then 
$e_i.v^+=0$ since $v^+$ is a primitive vector. Now let $l\geq 1$ and set 
$v':= f_{i_2}.(f_{i_3}.(\ldots, (f_{i_l}.v^+)\ldots))$. Then $v=f_{i_1}.v'$ 
and so
\[ e_i.v=e_i.(f_{i_1}.v')=[e_i,f_{i_1}].v'+ f_{i_1}.(e_i.v').\]
By induction, $e_i.v'\in V'$ and so $f_{i_1}.(e_i.v')\in V'$.
Furthermore, assume that $0\neq [e_i,f_{i_1}]\in \cL_{\alpha_i-
\alpha_{i_1}}$. If $i=i_1$, then $[e_i,f_{i_1}]=h_i\in \fh$ and so 
$[e_i,f_{i_1}].v'=h_i.v'\in V'$ by ($*^\prime$). If $i \neq i_1$, then 
$[e_i,f_{i_1}]=0$. Hence, in all cases, $e_i.v\in V'$, as required.
Thus, $V'\subseteq V$ is a $\cL$-submodule. 

Now let $n=\dim V'$ and $\{v_1,\ldots,v_n\}$ be a basis of $V'$ where 
each $v_i$ is of the form ($*$) and the notation is such that $v_1=v^+$. 
Then ($*^\prime$) shows that $v_i\in V_{\mu_i}$ where $\mu_i \preceq 
\lambda$ for all~$i$; furthermore, $\mu_i\neq \lambda$ for $i\geq 2$.
This already implies that $V'$ is $\fh$-diagonalisable. Now let $\mu 
\in \fh^*$ be such that $V_\mu'\neq \{0\}$. Let $0\neq v \in V_\mu'$ 
and write $v=\sum_i c_iv_i\in \sum_i V_{\mu_i}'$, where $c_i\in \C$. By 
Exercise~\ref{xcaweights}, we must have $\mu=\mu_i$ for some $i$ and so 
$\mu\preceq \lambda$. Finally, assume that $0\neq u\in V_\lambda'$ and 
write $u=c_1v_1+\ldots +c_nv_n$ where $c_i\in \C$. Hence, $u-c_1v_1\in 
V_\lambda'$ and $u-c_1v_1\in \sum_{i\geq 2} V_{\mu_i}'$. Since $\mu_i\neq 
\lambda$ for all $i\geq 2$, Exercise~\ref{xcaweights} implies that 
$u-c_1v_1=0$ and so $u\in \langle v^+\rangle_\C$, as desired. 
\end{proof}

\begin{defn}  \label{defhwm} In the set-up of Proposition~\ref{primitive1}, 
let $0\neq v^+\in V$ be a primitive vector (of weight $\lambda\in \fh^*$) 
and assume that $V=V'$; in particular, $V$ is $\fh$-diagonalizable. Since 
$\mu \preceq \lambda$ for all $\mu\in P_\fh(V)$, the weight $\lambda$ is 
called the \nm{highest weight} of $V$, and $V$ itself is called a \nm{highest 
weight module}.
\end{defn}


\begin{exmp} \label{irredhwm} Assume that $V$ is an irreducible
$\cL$-module. Then~$V$ automatically is a highest weight module.
For, in this case, the submodule $V'\subseteq V$ in 
Proposition~\ref{primitive1} must be equal to $V$ (since it is non-zero). 
This also shows that $V$ is $\fh$-diagonalisable. Also note: If $0\neq 
v^+\in V$ is a primitive vector (of weight $\lambda\in \fh^*$), then any
other primitive vector in $V$ is a scalar multiple of~$v^+$. Indeed, let 
also $0\neq w^+\in V$ be a primitive vector, of weight $\mu\in \fh^*$. 
Then $\mu\in P_\fh(V)$ and so $\mu\preceq \lambda$ by 
Proposition~\ref{primitive1}. But we can also apply the construction in 
Proposition~\ref{primitive1} to $w^+$ and, hence, deduce that $\lambda
\preceq \mu$. Thus, we must have $\lambda=\mu$ and $w^+ \in V_\lambda=
\langle v^+\rangle_\C$.
\end{exmp}

\begin{exmp} \label{defhwm1} 
(a) Assume that $\cL$ is a simple Lie algebra. We regard $V:=\cL$ as a
$\cL$-module via the adjoint representation (as in Example~\ref{defhwmL}). 
Then a $\cL$-submodule of $V$ is the same thing as an ideal in $\cL$. So 
$V$ is irreducible and, hence, a highest weight module for $\cL$. What is
the highest weight? Now, by Remark~\ref{highestr}, there is a unique 
root $\alpha_0\in \Phi^+$ of maximal height. Let $0\neq v^+ \in 
V_{\alpha_0}$. Since $\alpha_0+ \alpha_i\not\in \Phi$ for $i\in I$, we 
have $e_i.v^+=[e_i,v^+] \subseteq \fg_{\alpha_0+\alpha_i}=\{0\}$ for 
$i\in I$. Hence, $v^+$ is a primitive vector, and $\alpha_0$ is the 
highest weight.

(b) Let $\cL=\slm_n(\C)$ or $\cL=\gom_n(Q_n,\C)$, as in 
Example~\ref{defhwm2} and Exercise~\ref{hwmgom}. Let $V=\C^n$ be the 
natural $\cL$-module. In each case, one easily sees that the first 
vector $v_1$ in the standard basis of $V$ is a primitive vector, with 
corresponding highest weight~$\varepsilon_1$. Note that $\varepsilon_1=
\varpi_1$ for $\cL=\slm_n(\C)$, and $\varepsilon_1=\varpi_m$ for
$\cL=\gom_n(Q_n,\C)$ (where, as usual, $n=2m$ or $n=2m+1$).
\end{exmp}


\begin{exmp} \label{trivmod} Let $V$ be a $1$-dimensional $\cL$-module;
then $V$ is irreducible. Furthermore, since $\cL=[\cL,\cL]$, we must have 
$x.v=0$ for $x\in \cL$ and $v\in V$ (see Example~\ref{slnss0a}). Thus,
the heighest weight of~$V$ is $\lambda=\underline{0}$. Conversely, 
assume that $V$ is an irreducible $\cL$-module with highest weight 
$\lambda=\underline{0}$. We claim that then $\dim V=1$. Indeed, let 
$0\neq v^+\in V_{\underline{0}}$ be a primitive vector. Assume that $v'
:=f_j.v^+\neq 0$ for some $j\in I$. Then, $v'\in V_{\underline{0}-
\alpha_j}=V_{-\alpha_j}$ by Remark~\ref{fundcal}, and so $-\alpha_j\in 
P_\fh(V)$. But then we also have $\alpha_j=s_j(-\alpha_j)\in P_\fh(V)$ by 
Proposition~\ref{integralrep} and so $\alpha_j\preceq \underline{0}$, 
contradiction. Thus, we must have $f_j.v^+=0$ for all $j\in I$ and, hence,
$V=\langle v^+\rangle_\C$ by Proposition~\ref{primitive1}. 
\end{exmp}

Our next aim is to construct highest weight modules such that the
highest weight is ``minuscule'' in the sense of Definition~\ref{minusc}.

\begin{defn}[Jantzen \protect{\cite[5A.1]{Ja}}\footnote{Note that Jantzen 
\cite{Ja} actually deals with the quantum group case, which gives rise to 
a number of technical complications which are not present in our setting 
here.}, see also \protect{\cite[\S 2]{G2}}] \label{defj} We consider
a non-empty subset $\Psi\subseteq \Omega$ that is a finite union of 
$W$-orbits of non-zero minuscule weights; thus, $\langle \alpha^\vee,
\mu\rangle \in\{0,\pm 1\}$ for $\alpha\in\Phi$ and $\mu\in\Psi$ (see 
Remark~\ref{minusc0}). Let~$M$ be a $\C$-vector space with a basis
$\{z_\mu \mid \mu \in \Psi\}$. For $i\in I$ we define linear maps 
$\tilde{e}_i\colon M\rightarrow M$ and $\tilde{f}_i\colon M\rightarrow M$ 
as follows, where $\mu\in \Psi$:
\begin{align*}
\tilde{e}_i(z_\mu)&:=\left\{\begin{array}{cl} z_{\mu+\alpha_i} & \quad
\mbox{if $\langle \alpha_i^\vee,\mu\rangle=-1$},\\0 &\quad\mbox{otherwise},
\end{array}\right.\\
\tilde{f}_i(z_\mu)&:=\left\{\begin{array}{cl} z_{\mu-\alpha_i} & \quad
\mbox{if $\langle \alpha_i^\vee,\mu\rangle=1$},\\0 &\quad\mbox{otherwise},
\end{array}\right.
\end{align*}
These maps are well-defined: if $\langle \alpha_i^\vee,\mu\rangle=-1$, 
then $\mu+\alpha_i=s_i(\mu)\in \Psi$; similarly, if $\langle 
\alpha_i^\vee,\mu\rangle=1$, then $\mu-\alpha_i=s_i(\mu)\in\Psi$. 
Note also that, for any $\mu\in \Psi$ there exists an $i\in I$ such
that $\langle \alpha_i^\vee,\mu\rangle=\pm 1$ (since $\{\alpha_i^\vee
\mid i \in I\}$ is a basis of $E$). In particular, there exists at least
some $i\in I$ such that $\tilde{e}_i\neq 0$ or $\tilde{f}_i\neq 0$.
\end{defn}



\begin{rem} \label{rem10} Let $|\Psi|=d\geq 1$ and let us choose an 
enumeration $\Psi=\{\mu_1,\ldots,\mu_d\}$ such that $i>j$ whenever $\mu_i
\preceq \mu_j$ and $\mu_i\neq \mu_j$. Then the above formulae show that 
each $\tilde{e}_i$ is represented by a strictly upper triangular matrix 
and each $\tilde{f}_i$ is represented by a strictly lower triangular 
matrix. In particular, the $\tilde{e}_i$ and $\tilde{f}_i$ are 
nilpotent linear maps; in fact, we have $\tilde{e}_i^2=\tilde{f}_i^2=0$ 
for $i\in I$. This is easily seen as follows. Let $\mu\in \Psi$ and 
assume that $\tilde{e}_i(z_\mu)\neq 0$. Then $\langle \alpha_i^\vee,
\mu\rangle=-1$ and $\tilde{e}_i(z_\mu)=z_{\mu+\alpha_i}$. But then 
$\langle \alpha_i^\vee, \mu+\alpha_i\rangle=-1+2=1$ and so 
$\tilde{e}_i(z_{\mu+\alpha_i})=0$, as required. The argument for 
$\tilde{f}_i$ is analogous.
\end{rem}

\begin{lem}[Cf.\ \protect{\cite[5A.1]{Ja}, \cite[\S 2]{G2}}] \label{lemj} 
Let us also define $\tilde{h}_i:=[\tilde{e}_i,\tilde{f}_i]\in \gl(M)$ 
for $i\in I$. Then the linear maps $\tilde{e}_i,\tilde{f}_i,\tilde{h}_i
\in \gl(M)$ satisfy the following \nm{Chevalley relations}, exactly as
in Section~\ref{sec3a3}:
\begin{align*}
[\tilde{h}_i,\tilde{h}_j]&=0, \quad [\tilde{h}_i,\tilde{e}_j]=a_{ij}
\tilde{e}_j, \quad [\tilde{h}_i,\tilde{f}_j]=-a_{ij}\tilde{f}_j
\quad \mbox{for all $i,j\in I$};\\
[\tilde{e}_i,\tilde{f}_j]&=0 \quad\mbox{for all $i,j\in I$ such 
that $i\neq j$}.
\end{align*}
Furthermore, we have $\tilde{h}_i(z_\mu)=\langle \alpha_i^\vee,\mu 
\rangle z_\mu=\mu(h_i)z_\mu$ for all $\mu\in \Psi$.
\end{lem}

\begin{proof} First we prove the formula for $\tilde{h}_i(z_\mu)$.
If $\langle \alpha_i^\vee,\mu \rangle=\mu(h_i)=-1$, then 
$\tilde{e}_i(z_\mu)=z_{\mu+\alpha_i}$, $\tilde{f}_{i}(z_\mu)=0$ and so 
\[ \tilde{h}_i(z_\mu)=[\tilde{e}_i, \tilde{f}_i](z_\mu)=-\tilde{f}_i
(\tilde{e}_i(z_\mu))=-\tilde{f}_i(z_{\mu+\alpha_i})=-z_\mu,\]
since $\langle \alpha_i^\vee,\mu+\alpha_i \rangle=1$. Thus, the
formula holds in this case. If $\langle \alpha_i^\vee,\mu\rangle=1$, 
the argument is analogous. If $\langle \alpha_i^\vee,\mu\rangle=0$, then 
$\tilde{e}_i(z_\mu)=\tilde{f}_i(z_\mu)=0$ and so $\tilde{h}_i(z_\mu)=
[\tilde{e}_i,\tilde{f}_i](z_\mu)=\tilde{e}_i(\tilde{f}_i(z_\mu))-
\tilde{f}_i(\tilde{e}_i(z_\mu))=0$. Hence, the formula holds in this case
as well. The formula for $\tilde{h}_i(z_\mu)$ shows, in particular,
that $\tilde{h}_i$ is represented by a diagonal matrix. So it is clear 
that $[\tilde{h}_i,\tilde{h}_j]=0$ for all $i,j$. 


Now consider the remaining relations. We begin with 
$[\tilde{h}_i,\tilde{e}_j]$. Let $\mu\in \Psi$. If $\langle 
\alpha_j^\vee,\mu\rangle=-1$, then $\tilde{e}_j(z_\mu)=
z_{\mu+\alpha_j}$ and so 
\begin{align*}
[\tilde{h}_i,\tilde{e}_j](z_\mu)&=\tilde{h}_i(z_{\mu+\alpha_j})-
\langle \alpha_i^\vee,\mu\rangle \tilde{e}_j(z_\mu)\\
&=\langle \alpha_i^\vee,\mu+\alpha_j\rangle z_{\mu+\alpha_j}-
\langle \alpha_i^\vee,\mu\rangle z_{\mu+\alpha_j}\\
&=\langle \alpha_i^\vee,\alpha_j\rangle \tilde{e}_j(z_\mu)=
a_{ij}\tilde{e}_j(z_\mu),
\end{align*} 
as required. If $\langle \alpha_j^\vee,\mu\rangle\neq -1$, then 
$\tilde{e}_j(z_\mu)=0$ and 
\[[\tilde{h}_i,\tilde{e}_j](z_\mu)=\tilde{h}_i(\tilde{e}_j(z_\mu))-
\langle \alpha_i^\vee,\mu\rangle \tilde{e}_j(z_\mu)=0.\]
Hence, we obtain again $[\tilde{h}_i,\tilde{e}_j](z_\mu)=a_{ij}
\tilde{e}_j(z_\mu)$. The argument for proving $[\tilde{h}_i,
\tilde{f}_j]=-a_{ij} \tilde{f}_j$ is completely analogous.

Finally, let $j\in I$, $i\neq j$. We must show that 
$[\tilde{e}_i,\tilde{f}_j]=0$, that is, $\tilde{e}_i(\tilde{f}_j
(z_\mu))=\tilde{f}_j(\tilde{e}_i(z_\mu))$. Now, we have
\begin{align*}
\tilde{f}_j(\tilde{e}_i(z_\mu))&=\left\{\begin{array}{cl}
z_{\mu+\alpha_i-\alpha_j} & \mbox{if $\langle\alpha_i^\vee,\mu
\rangle=-1$ and $\langle \alpha_j^\vee,\mu+\alpha_i\rangle=1$},\\ 
0 & \mbox{otherwise}; \end{array}\right.\\
\tilde{e}_i(\tilde{f}_j(z_\mu))&=\left\{\begin{array}{cl}
z_{\mu+\alpha_i-\alpha_j} & \mbox{if $\langle \alpha_j^\vee,\mu
\rangle=1$ and $\langle\alpha_i^\vee,\mu-\alpha_j\rangle=-1$},\\ 
0 & \mbox{otherwise}. \end{array}\right.
\end{align*}
So it remains to show that the conditions on the right hand side are
equivalent. Assume first that $\langle \alpha_i^\vee,\mu\rangle=-1$ and
$\langle\alpha_j^\vee, \mu+\alpha_i\rangle=1$. Since $i\neq j$, we have
$\langle \alpha_j^\vee,\alpha_i\rangle=a_{ji}\leq 0$ and so 
$\langle \alpha_j^\vee,\mu\rangle=1-\langle \alpha_j^\vee,\alpha_i
\rangle\geq 1$. Hence, we must have $\langle \alpha_j^\vee,\mu\rangle
=1$ and $a_{ji}=\langle \alpha_j^\vee,\alpha_i\rangle=0$. But then we 
also have $\langle \alpha_i^\vee,\alpha_j\rangle=a_{ij}=0$ and so 
$\langle \alpha_i^\vee,\mu-\alpha_j\rangle=-1$, as required. The 
reverse implication is proved similarly.
\end{proof}

\begin{prop} \label{corj} In the set-up of Definition~\ref{defj}, there 
is a unique homomorphism of Lie algebras $\rho\colon \cL\rightarrow 
\gl(M)$ such that $\rho(e_i)=\tilde{e}_i$ and $\rho(f_i)=\tilde{f}_i$ 
for $i\in I$. Thus, $M$ is a $\cL$-module. We have $M_\mu=\langle z_\mu
\rangle_\C$ for all $\mu\in \Psi$; so $M$ is $\fh$-diagonalisable with 
$P_\fh(M)=\Psi$. If $A$ is indecomposable, then $M$ is a faithful 
$\cL$-module.
\end{prop}


\begin{proof} We define two subsets $I_1,I_2\subseteq I$ as follows.
\begin{align*}
I_1 &:=\{i \in I \mid \langle \alpha_i^\vee,\mu \rangle=0 \mbox{ for 
all $\mu \in \Psi$}\},\\
I_2 &:=\{j \in I \mid \langle \alpha_j^\vee,\mu \rangle \neq 0 \mbox{ for 
some $\mu \in \Psi$}\}.
\end{align*}
Then $I=I_1\cup I_2$ and $I_1\cap I_2=\varnothing$. Since $\{\alpha_j^\vee
\mid j \in I\}$ is a basis of~$E$, and since $\mu \neq \underline{0}$ for 
all $\mu \in \Psi$, we certainly have $I_2\neq \varnothing$. (But it could
happen that we also have $I_1\neq \varnothing$.) Note that Lemma~\ref{lemj} 
shows that $\tilde{h}_j\neq 0$ for all $j\in I_2$. Since $\tilde{h}_j=
[\tilde{e}_j,\tilde{f}_j]$, this also implies that $\tilde{e}_j\neq 0$ and 
$\tilde{f}_j\neq 0$ for all $j\in I_2$. On the other hand, if $i\in I_1$,
then Lemma~\ref{lemj} shows that $\tilde{h}_i=0$. This implies that 
$0=[\tilde{h}_i,\tilde{e}_i]=2\tilde{e}_i$ and $0=[\tilde{h}_i,
\tilde{f}_i]=2\tilde{f}_i$; hence, $\tilde{e}_i=0$ and $\tilde{f}_i=0$. 
Thus, we have 
\[ \tilde{\cL}:=\langle \tilde{e}_i,\tilde{h}_i,\tilde{f}_i\mid i \in I
\rangle_{\text{alg}}=\langle \tilde{e}_j,\tilde{h}_j,\tilde{f}_j\mid j \in 
I_2 \rangle_{\text{alg}}\subseteq \gl(M).\] 
Let us now first deal with the case where $I_1=\varnothing$. Then
Lemma~\ref{lemj} shows that $\tilde{\cL}$ and the elements $\{\tilde{e}_i,
\tilde{h}_i, \tilde{f}_i\mid i \in I=I_2\}$ satisfy (Ch0), (Ch1), (Ch2)
in Section~\ref{sec3a3}. As noted above, we have $\tilde{e}_j\neq 0$ and
$\tilde{f}_j\neq 0$ for all $j\in I=I_2$. So all the assumptions of
Proposition~\ref{triang3b} are satisfied. We conclude that the Lie algebra 
$\tilde{\cL}$ is of Cartan--Killing type, with structure matrix~$A$. 
By the Isomorphism Theorem~\ref{isothm}, there is a unique isomorphism of 
Lie algebras $\cL\cong \tilde{\cL}$ such that $e_i\mapsto \tilde{e}_i$, 
$f_i\mapsto \tilde{f}_i$ for $i\in I$. This yields the desired Lie algebra
homomorphism $\rho \colon \cL\rightarrow \gl(M)$; since it is injective, 
the module $M$ is faithful in this case. The formula for the action of 
$\tilde{h}_i$ ($i\in I$) shows that $M$ is $\fh$-diagonalisable, with 
$P_\fh(M)=\Psi$. 

Now assume that $I_1\neq \varnothing$. We claim that then $a_{ij}=0$ for 
all $i\in I_1$ and all $j \in I_2$. (In particular, $A$ is not indecomposable
in this case.) Indeed, as noted above, we have $\tilde{h}_i=0$ for $i\in 
I_1$, and $\tilde{e}_j\neq 0$ for $j\in I_2$. Hence, $0=[\tilde{h}_i,
\tilde{e}_j]=a_{ij}\tilde{e}_j$ and so $a_{ij}=0$, as claimed. By 
Remark~\ref{mainideal3} and Proposition~\ref{levisub}, we have $\cL=
\cL_1\oplus \cL_2$ where $\cL_s$ (for $s=1,2$) is of Cartan--Killing type 
with structure matrix $A_s=(a_{ij})_{i,j\in I_s}$; furthermore, $[\cL_1,
\cL_2]=\{0\}$. As above, we see that $\tilde{\cL}$ is of Cartan--Killing
type with structure matrix $A_2$. Hence, again, the Isomorphism Theorem
yields an isomorphism of Lie algebras $\cL_2\cong \tilde{\cL}$ such that 
$e_i\mapsto \tilde{e}_i$, $f_i\mapsto \tilde{f}_i$ for $i\in I_2$. 
Composing this with the projection $\cL\rightarrow \cL_2$, we obtain
the desired homomorphism $\rho\colon \cL\rightarrow \gl(M)$. Note that
$\{0\}\neq \cL_1\subseteq \ker(\rho)$, so the module $M$ is not faithful in 
this case. The formula for the action of $\tilde{h}_i$ ($i\in I$) shows 
again that $M$ is $\fh$-diagonalisable, with $P_\fh(M)=\Psi$. 
\end{proof}

\begin{rem} \label{jhwm} In the set-up of Definition~\ref{defj}, assume 
that $\Psi$ is a single $W$-orbit of a non-zero dominant minuscule 
weight~$\lambda_0\in \Omega$. As remarked above, we have $z_{\lambda_0}
\in M_{\lambda_0}$. Since $\lambda_0$ is dominant, we have $\langle
\alpha_i^\vee,\lambda_0\rangle\geq 0$ and, hence, $\langle 
\alpha_i^\vee,\lambda_0\rangle\neq -1$ for $i\in I$. But then 
$e_i.z_{\lambda_0}=0$ for $i\in I$. So $z_{\lambda_0}$ is a primitive 
vector of weight~$\lambda_0$. We claim that:
\[\text{\textit{$M$ is irreducible with highest weight $\lambda_0$}}.\]
Indeed, let $U\subseteq M$ be an irreducible submodule. Then $U$ is 
$\fh$-diagonalisable and $\varnothing \neq P_\fh(U) \subseteq P_\fh(M)=\Psi$; 
see Proposition~\ref{wsdprop1}. Since $\Psi$ is a single $W$-orbit and since 
$P_\fh(U)$ is a union of $W$-orbits (see Proposition~\ref{integralrep}),
we must have $P_\fh(U)=\Psi$. Consequently, we have $U_\mu\neq \{0\}$ for 
all $\mu \in \Psi$ and so $\dim U \geq |\Psi|=\dim M$, that is, $U=M$. 
Hence, $M$ is irreducible and, consequently, a highest weight module 
(see Example~\ref{irredhwm}).
\end{rem}

\begin{exmp} \label{exmpj1} Let $\cL=\slm_n(\C)$, $n\geq 2$, and $\fh
\subseteq \cL$ be the subalgebra of diagonal matrices. As usual, for $1\leq 
i \leq n$ let $\varepsilon_i\in \fh^*$ be the linear map that sends a 
diagonal matrix to its $i$-th diagonal entry. Now consider the fundamental 
minuscule weight $\varpi_1=\varepsilon_1$. By Example~\ref{minusAn}, the 
$W$-orbit of $\varpi_1$ is $\{\varepsilon_1,\ldots,\varepsilon_n\}$. So, 
according to Definition~\ref{defj}, let $M$ be a $\C$-vector space with a 
basis $z_{\varepsilon_1},\ldots,z_{\varepsilon_n}$. By Example~\ref{cartsln}, 
the simple roots of $\cL$ are given by $\alpha_i=\varepsilon_i-
\varepsilon_{i+1}$ for $1\leq i \leq n-1$. Now note that 
\[ \renewcommand{\arraystretch}{0.9} \langle \alpha_i^\vee,\varepsilon_j
\rangle=\varepsilon_j(h_i)=\left\{\begin{array}{rl} 1 & \quad \mbox{if 
$i=j$},\\ -1 & \quad \mbox{if $j=i+1$},\\ 0 & \quad \mbox{otherwise}.
\end{array}\right.\]
Hence, the linear maps $\tilde{e}_i$ and $\tilde{f}_i$ in 
Definition~\ref{defj} act on $M$ by the following formulae.

If $j=i+1$, then $\langle\alpha_i^\vee,\varepsilon_j\rangle=-1$ and
$\varepsilon_j+\alpha_i=\varepsilon_{i+1}+(\varepsilon_i-
\varepsilon_{i+1})=\varepsilon_i$. Hence, 
$\tilde{e}_i.z_{\varepsilon_{i+1}}=z_{\varepsilon_{i}}$. Otherwise,
if $j\neq i+1$, we have $\tilde{e}_i.z_{\varepsilon_j}=0$. 

Similarly, if $j=i$, then $\langle\alpha_i^\vee,\varepsilon_j\rangle=1$
and $\varepsilon_j-\alpha_i=\varepsilon_{i}-(\varepsilon_i-
\varepsilon_{i+1})=\varepsilon_{i+1}$. Hence, 
$\tilde{f}_i.z_{\varepsilon_i}=z_{\varepsilon_{i+1}}$. Otherwise,
if $j\neq i$, we have $\tilde{f}_i.z_{\varepsilon_j}=0$. 

On the other hand, consider the standard module $V=\C^n$ and the 
Chevalley generators $e_i,f_i$ described in Example~\ref{cartsln}. 
We see that $e_i,f_i$ act in the same way on $\C^n$ as $\tilde{e}_i$, 
$\tilde{f}_i$ act on $M$. Thus, $M$ can be canonically 
identified with the standard module $V=\C^n$.
\end{exmp}

\begin{xca} \label{exmpj2} Let $\cL=\gom_n(Q_n,\C)$ where $n=2m\geq 4$ 
is even; if $Q_n=Q_n^{\text{tr}}$, also assume that $n=2m\geq 6$. Thus, 
$\fg$ is of type $C_m$ or $D_m$. By Examples~\ref{minusCn} 
and~\ref{minusDn}, the fundamental weight $\varpi_m=\varepsilon_1$ is 
minuscule, with $W$-orbit given by $\{\pm \varepsilon_1,\ldots, \pm 
\varepsilon_m\}$. In particular, the module $M$ in Definition~\ref{defj} 
has dimension $n=2m$. Show that $M$ can be canonically identified with 
the natural module $V=\C^n$. (Argue as in Example~\ref{exmpj1} and use 
the results in Section~\ref{sec1a5}, especially the description of the
Chevalley generators in Proposition~\ref{CKbcd}.)
\end{xca}

\begin{exmp} \label{exmpj3} Let $\cL=\gom_n(Q_n,\C)$ where $n\geq 4$ and
$Q_n=Q_n^{\text{tr}}$. Thus, either $n=2m+1$ and $\cL$ is of type $B_m$, 
or $n=2m$ and $\cL$ is of type $D_m$ where, in the latter case, we also 
assume that $n\geq 6$. 

Assume first that $\cL$ is of type $B_m$. By Example~\ref{minusBn}, the 
fundamental weight $\varpi_1$ is minuscule, with $W$-orbit of size $2^m$.
Following Bourbaki \cite[Ch.~VIII, \S 13, no.~2]{B78}, the 
corresponding $\cL$-module~$M$ in Definition~\ref{defj} is called
the \nm{spinor representation} of~$\fg$ (or {\itshape ``repr\'esentation
spinorielle}'' in French).

Now assume that $\cL$ is of type $D_m$. By Example~\ref{minusDn}, the 
fundamental weights $\varpi_1$ and $\varpi_2$ are minuscule, with 
$W$-orbits of size $2^{m-1}$. Following Bourbaki \cite[Ch.~VIII, \S 13, 
no.~4]{B78}, the corresponding $\cL$-modules~$M$ in Definition~\ref{defj} 
are called the \nm{semi-spinor representations}\footnote{Carter 
\cite[\S 13.5]{Ca3} and Fulton--Harris \cite[Chap.~20]{FH} use the terms 
\nm{half-spin representations} and \nm{spin representations}; we shall
also use these terms later on.} of~$\cL$ (or {\itshape ``repr\'esentations 
semi-spinorielles}'' in French). Their direct sum will be called the 
{\itshape spinor representation} of~$\cL$.

(Both in \cite{B78} and in \cite{FH}, these representations are constructed 
using the Clifford algebra of~$V=\C^n$; see also Carter~\cite[\S 13.5]{Ca3}.)
\end{exmp}

\begin{exmp} \label{minusE67} Let $\fg$ be of type $E_6$ or $E_7$.
There is a long history about the representations with a minuscule
highest weight in these cases. This includes numerous connections with 
geometry (keyword: ``the $27$ lines on a cubic surface'') and also with
finite group theory; we just mention Aschbacher \cite{AscE6}, Hulek
\cite[Chap.~5]{Hul}, Lurie \cite{lur}, Springer \cite{Spr2}, Vavilov 
\cite{vav1} and the further references there. 

As an exercise for the reader, show that the two modules corresponding to 
$\varpi_1$ and to~$\varpi_6$ (in type~$E_6$) are dual to each other in 
the sense of Exercise~\ref{dualspace}.
\end{exmp}

\begin{exmp} \label{expminus} In {\sf ChevLie} the minuscule
dominant weights are specified in the component {\tt .minuscule}
of {\tt LieAlg}. The following example shows how to create the
$W$-orbit of $\varpi_1$ for $\fg$ or type $E_6$.

\begin{verbatim}
  julia> l=LieAlg(:e,6)
  julia> l.minuscule     
  1  6                   # w_1 and w_6 are minuscule
  julia> v=zeros(Int8,6);v[1]=1;    # create W-orbit
  julia> println(weightorbit(l,v))  # of w_1
\end{verbatim}

{\small \begin{verbatim}
  [[1,0,0,0,0,0],[-1,0,1,0,0,0],[0,0,-1,1,0,0],
   [0,1,0,-1,1,0],[0,-1,0,0,1,0],[0,1,0,0,-1,1],
   [0,-1,0,1,-1,1],[0,1,0,0,0,-1],[0,0,1,-1,0,1],
   [0,-1,0,1,0,-1],[1,0,-1,0,0,1],[0,0,1,-1,1,-1],
   [-1,0,0,0,0,1],[1,0,-1,0,1,-1],[0,0,1,0,-1,0],
   [-1,0,0,0,1,-1],[1,0,-1,1,-1,0],[-1,0,0,1,-1,0],
   [1,1,0,-1,0,0],[-1,1,1,-1,0,0],[1,-1,0,0,0,0],
   [-1,-1,1,0,0,0],[0,1,-1,0,0,0],[0,-1,-1,1,0,0],
   [0,0,0,-1,1,0],[0,0,0,0,-1,1],[0,0,0,0,0,-1]]
\end{verbatim}}

\noindent (If $[n_1,\ldots,n_6]$ is one of the above $27$ vectors, 
then $n_1\varpi_1{+}\ldots{+}n_6\varpi_6\in \Omega$ is an element in
the $W$-orbit of $\varpi_1$.) The corresponding $27$-di\-men\-sion\-al 
representation of $\fg$ is created using the function
{\tt rep$\_$minuscule}.
\end{exmp}

We can now solve the problem stated at the end of Section~\ref{sechighw}, 
concerning the existence of $\cL$-modules with a prescribed weight lattice.

\begin{thm} \label{corj2} Let $\Lambda'\subseteq \Omega$ be an arbitrary 
subgroup with $\Phi \subseteq \Lambda'$. Then there exists a faithful,
$\fh$-diagonalizable $\cL$-module $V$ with $\Lambda(V)=\Lambda'$. In 
fact,~$V$ can be taken as the direct sum of $\cL$ (adjoint module) and,
possibly, various irreducible $\cL$-modules with a minuscule highest weight. 
(Note that $\cL$ also is a direct sum of irreducible $\cL$-submodules
by Example~\ref{liedirsum}.)
\end{thm}

\begin{proof} If $\Lambda'=\langle \Phi\rangle_\Z$, then we can just take
$V:=\fg$. Now assume that $\langle \Phi\rangle_\Z\subsetneqq \Lambda'$. As 
in Remark~\ref{corminus1}, there exist non-zero dominant minuscule weights 
$\lambda_1^\circ,\ldots,\lambda_m^\circ \in \cM$ (for some $m\geq 1$) 
such that 
\[ \Lambda'=\big\langle \Phi\cup\Psi_1\cup\ldots\cup\Psi_m\big\rangle_\Z,\]
where $\Psi_l$ is the $W$-orbit of $\lambda_l^\circ$ for $l=1,\ldots,m$. 
For $l=0$, let $M_0:=\cL$; then $P_\fh(M_0)=\{\underline{0}\} \cup \Phi$. 
For $l\geq 1$, consider the $\cL$-module $M_l$ of Proposition~\ref{corj}, 
constructed from the $W$-orbit $\Psi_l$; then $P_\fh(M_l)=\Psi_l$. Now 
take the direct sum $V:=M_0\oplus M_1\oplus \ldots \oplus M_m$. This is 
a $\cL$-module such that 
$P_\fh(V)=\{\underline{0}\} \cup \Phi\cup \Psi_1\cup \ldots \cup \Psi_m$
(see Exercise~\ref{dirsumfaith}). Hence, $\Lambda(V)=\langle P_\fh(V)
\rangle_\Z=\Lambda'$, as required. Note that, since $\Phi\subseteq 
P_\fh(V)$, the module $V$ is faithful. Furthermore, the modules 
$M_1,\ldots,M_m$ are irreducible by Remark~\ref{jhwm}.
\end{proof}

The construction of $V$ in the above proof is not optimised, in the sense
that $V$ does not necessarily have the smallest possible dimension such 
that $\Lambda(V)=\Lambda'$. In concrete situations, one may well be able
to find better models of $V$, especially when $\fg$ is simple.

\begin{exmp} \label{oprmod} Let $\fg$ be simple. Let $\Lambda' \subseteq 
\Omega$ be a subgroup with $\Phi\subseteq \Lambda'$. If $\Lambda'=\langle 
\Phi\rangle_\Z$, then we have $\Lambda(V)=\Lambda'$ for $V=\fg$ (adjoint 
module). For $\fg$ of type $G_2$, $F_4$ or $E_8$, there are no further 
cases to consider (see Table~\ref{fundtab}, p.~\pageref{fundtab}). Now 
let $\langle \Phi\rangle_\Z \subsetneqq \Lambda'$. We claim that:
\begin{center}
\begin{tabular}{l} {\itshape If $\Lambda'/\langle \Phi\rangle_\Z$ is 
cyclic, then $\Lambda'=\Lambda(V)$ where $V$ is an}\\
{\itshape irreducible $\fg$-module with a minuscule highest weight}.
\end{tabular}
\end{center}
Indeed, the assumption implies that $\Lambda'=\langle \Phi,\varpi_i
\rangle_\Z$ for some minuscule fundamental weight $\varpi_i$ ($i\in I$).
Let $V$ be the $\fg$-module of Proposition~\ref{corj}, constructed from 
the $W$-orbit of~$\varpi_i$. Since $\fg$ is simple, the module $V$ is 
faithful and, hence, $\Phi\subseteq \Lambda(V)$; we also have $\varpi_i
\in \Lambda(V)$. Since $P_{\fh}(V)$ just consists of the $W$-orbit of 
$\varpi_i$ (and this is contained in~$\Lambda'$), we conclude that 
$\Lambda'=\Lambda(V)$, as clained. 

By Table~\ref{fundtab} (p.~\pageref{fundtab}), this actually covers all 
cases except one. The exception occurs when $\Lambda'=\Omega$ and $\fg$ 
is of type $D_n$ with $n\geq 4$ even. In this case, $\Lambda'/\langle \Phi
\rangle_\Z=\Omega/\langle \Phi\rangle_\Z \cong \Z/2\Z\times \Z/2\Z$ and we 
have $\Lambda'=\Lambda(V)$, where $V$ is the direct sum of the two 
\nm{half-spin representations}.
\end{exmp}

\begin{xca} \label{expA1A1} Assume that $\fg=\fg_1\oplus \fg_2$
where $\fg_1,\fg_2$ are subalgebras such that $[\fg_1,\fg_2]=\{0\}$
and $\fg_s\cong \slm_2(\C)$ for $s=1,2$. Let $\{e_1,e_2,f_1,f_2,h_1,h_2\}$
be Chevalley generators such that $[e_i,f_i]=h_i$, $[h_i,e_i]=2e_i$,
$[h_i,f_i]=-2f_i$, $[e_i,e_j]=[f_i,f_j]=[e_i,f_j]=0$ for $i,j=1,2$ and 
$i\neq j$. We have $\Phi=\{\pm \alpha_1,\pm \alpha_2\}$; furthermore, 
$\Omega=\langle \varpi_1,\varpi_2\rangle_\Z$ where $\alpha_1=2\varpi_1$
and $\alpha_2=2\varpi_2$. Since $\Omega/\langle\Phi\rangle_\Z\cong\Z/2\Z 
\times \Z/2\Z$, there are five subgroups $\Lambda'\subseteq \Omega$ such 
that $\Phi \subseteq\Lambda'$. For $s=1,2$ let $V_s=\C^2$ be the natural 
$\cL_s$-module; we also regard $\fg_s$ itself as a $\fg_s$-module via the
adjoint representation. Show that, for each $\Lambda'$, a corresponding 
$\cL$-module $V$ such that $\Lambda(V)=\Lambda'$ can be constructed as 
follows.

\noindent (a) If $\Lambda'=\langle \Phi\rangle_\Z$, then we can take 
$V=\cL_1\oplus \cL_2$.

\noindent (b) If $\Lambda'=\Omega$, then we can take $V=V_1\oplus V_2$.

\noindent (c) If $\Lambda'=\langle\Phi,\varpi_1\rangle_\Z$, then 
we can take $V=V_1\oplus\cL_2$. Similarly, if $\Lambda'=\langle \Phi,
\varpi_2\rangle_\Z$, then we can take $V=\fg_1\oplus V_2$.

\noindent (d) Finally, let $\Lambda'=\langle\Phi,\varpi_1+\varpi_2
\rangle_\Z$.  Note that $\varpi_1+\varpi_2$ is minuscule by 
Proposition~\ref{minusc1} and Exercise~\ref{xcamindec}. Consider the 
$W$-orbit of $\varpi_1+\varpi_2$. By Remark~\ref{lemstem0a}, that 
$W$-orbit is given by
\[ \{\varpi_1+\varpi_2,\,\varpi_1+\varpi_2-\alpha_1,\,\varpi_1+\varpi_2-
\alpha_2,\,\varpi_1+\varpi_2-\alpha_1-\alpha_2\}.\]
Then we can take $V=M$ as in Proposition~\ref{corj}. Using 
Remark~\ref{tenslie1}, check that $V\cong V_1\otimes V_2$.
\end{xca}

\section{Admissible lattices in $\cL$-modules} \label{sec5a2}

We keep our basic assumptions. Let $\cL$ be a Lie algebra of Cartan--Killing
type, with Cartan subalgebra $\fh\subseteq \cL$ and root system $\Phi
\subseteq \fh^*$. Let $\{e_i,f_i\mid i \in I\}$ be Chevalley generators 
of $\cL$, as in Remark~\ref{astring0}. We now introduce the machinery 
needed to define Chevalley groups of non-adjoint type over an arbitrary 
field~$K$. This will be modeled on the procedure in Section~\ref{sec3a7}: 
first we define a group over $\C$, then we make sure that certain 
integrality conditions hold, which finally allow us to pass to a group 
over~$K$. Throughout, we tacitly assume that all $\cL$-modules that we 
consider have finite dimension. 
 
\begin{defn} \label{nonadj1} Let $V$ be a faithful, $\fh$-diagonalisable 
$\cL$-module, with corresponding representation $\rho\colon \cL\rightarrow 
\gl(V)$. By Lemma~\ref{wsdnilV}, for each $i\in I$, the linear maps 
$\rho(e_i) \colon V \rightarrow V$ and $\rho(f_i)\colon V\rightarrow V$ 
are nilpotent. So, for $i\in I$ and $t \in \C$, we can define 
\begin{equation*}
x_i(t; V):=\exp(t\rho(e_i)) \quad \mbox{and}\quad y_i(t;V)
:=\exp(t\rho(f_i))\quad \mbox{in $\GL(V)$}.
\end{equation*}
(Here we regard $V$ as an algebra with trivial product $v\cdot v'=0$ for 
$v,v'\in V$.) In analogy to the initial definition in Section~\ref{sec3a7},
we set 
\begin{equation*}
G_\C(V):=\big\langle x_i(t;V), y_i(t;V)\mid i\in I, t\in \C
\big\rangle \;\subseteq\; \GL(V).
\end{equation*}
For $i\in I$ and $t \in \C^\times$, we also have the elements
\begin{align*}
n_i(t;V)&:=x_i(t;V)y_i(-t^{-1};V)x_i(t;V)\in G_\C(V),\\
h_i(t;V)&:=n_i(t;V)n_i(-1;V)\in G_\C(V).
\end{align*}
Furthermore, let $\alpha\in \Phi$ and $0\neq \be_\alpha^+\in \cL_\alpha$,
an element of Lusztig's canonical basis (see Section~\ref{sec1a7}). Again, 
by Lemma~\ref{wsdnilV}, the linear map $\rho(\be_\alpha^+)\colon V 
\rightarrow V$ is nilpotent. So we can also define
\begin{equation*}
x_\alpha(t;V):=\exp(t\rho(\be_\alpha^+))\in \GL(V) \qquad
\mbox{for $\alpha\in \Phi$ and $t\in \C$}.
\end{equation*}
Since $\be_{\alpha_i}^+=\pm e_i$ and $\be_{-\alpha_i}^+=\pm f_i$, we 
have $x_{\alpha_i}(t;V)=x_i(\pm t;V)\in G_\C(V)$ and $x_{-\alpha_i}
(t;V)=y_i(\pm t;V)\in G_\C(V)$. But at this stage, it is not yet clear 
if $x_\alpha(t;V)$ belongs to $G_\C(V)$ for any $\alpha\in\Phi$.
\end{defn}

The following example shows that the above setting indeed is considerably 
wider than that for groups of ``adjoint type''.

\begin{exmp} \label{sl2reg} Let $\cL=\slm_2(\C)$ with $I=\{1\}$ and 
standard basis elements $e_1,f_1,h_1$ such that $[e_1,f_1]=h_1$. Let $V$ 
be an irreducible $\cL$-module with $\dim V\geq 2$, and let $\rho\colon 
\cL\rightarrow \gl(V)$ be the corresponding representation. Write $\dim V=
m+1$ with $m \geq 1$. Then, since $\fg$ is simple, we automatically have 
that $V$ is faithful. By Corollary~\ref{sl2modc}, there exists a basis 
$\{v_0,v_1,\ldots,v_m\}$ of $V$ such that 
\begin{center}
$\begin{array}{c@{\hspace{5pt}}l}
\rho_{e_1}\colon V\rightarrow V, & \qquad v_i\mapsto (m-i+1)v_{i-1}, \\
\rho_{f_1}\colon V\rightarrow V, & \qquad v_i\mapsto (i+1)v_{i+1}, \\
\rho_{h_1}\colon V\rightarrow V, & \qquad v_i\mapsto (m-2i)v_i,
\end{array}$
\end{center}
(where $v_{-1}=v_{m+1}=0$). Note that each basis vector $v_i$ is a weight 
vector for~$h_1$; so $V$ is $h$-diagonalizable. Now, we have 
\[ \rho_{e_1}^j(v_i)=(m{-}i{+}j)(m{-}i{+}(j{-}1))\cdots (m{-}i{+}1)v_{i-j}
\quad \mbox{for $0\leq j \leq i$},\]
and $\rho_{e_1}^j(v_i)=0$ for $j>i$. This yields that 
\[\frac{1}{j!}\rho_{e_1}^j(v_i)=\binom{m{-}i{+}j}{j}v_{i-j}
\qquad\mbox{for $0 \leq j \leq i$}.\] 
Consequently, for $t\in \C$, we obtain that 
\[ x_1(t;V)(v_i)=\sum_{0 \leq j \leq i} \binom{m{-}i{+}j}{j}t^j 
v_{i-j} =\sum_{0 \leq j \leq i} \binom{m{-}j}{m{-}i}t^{i-j}v_{j}.\]
Similarly, we have $\rho_{f_1}^j(v_i)=(i{+}j)(i{+}j{-}1)\cdots (i{+}1)
v_{i+j}$ for $0\leq j \leq m-i$, and $\rho_{f_1}^j(v_i)=0$ for $j>m-i$. 
This yields
\[ \frac{1}{j!}\rho_{f_1}^j(v_i)=\binom{i{+}j}{j}v_{i+j}=
\binom{i{+}j}{i}v_{i+j} \qquad \mbox{for $0 \leq j \leq m-i$}.\]
Consequently, we have 
\[ y_1(t;V)(v_i)=\sum_{0\leq j \leq m-i} \binom{i{+}j}{i}t^jv_{i+j} 
=\sum_{i \leq j \leq m} \binom{j}{i} t^{j-i}v_{j}.\]
For example, for $m=4$, the matrices of $x_1(t;V)$ and $y_1(t;V)$ with 
respect to the basis $\{v_0,v_1,v_2, v_3,v_4\}$ of $V$ are given as follows.
\[ x_1(t;V):\left(\begin{array}{@{\hspace{4pt}}c@{\hspace{4pt}}
c@{\hspace{4pt}}c@{\hspace{4pt}}c@{\hspace{4pt}}c@{\hspace{4pt}}}
        1&    4t&  6t^2&  4t^3&    t^4 \\
        0&      1&    3t&  3t^2&    t^3 \\
        0&      0&      1&    2t&    t^2 \\
        0&      0&      0&      1&      t \\
        0&      0&      0&      0&      1
\end{array}\right), \quad
y_1(t;V):\left(\begin{array}{@{\hspace{4pt}}c@{\hspace{4pt}}
c@{\hspace{4pt}}c@{\hspace{4pt}}c@{\hspace{4pt}}c@{\hspace{4pt}}}
        1&      0&      0&      0&      0 \\
        t&      1&      0&      0&      0 \\
      t^2&    2t&      1&      0&      0 \\
      t^3&  3t^2&    3t&      1&      0 \\
      t^4&  4t^3&  6t^2&    4t&      1 
\end{array}\right).\]
We also find that the matrices of $n_1(t;V)$ and $h_1(t;V)$ are given by:
\[ n_1(t;V):\left(\begin{array}{@{\hspace{1pt}}c@{\hspace{4pt}}
c@{\hspace{4pt}}c@{\hspace{4pt}}c@{\hspace{4pt}}c@{\hspace{4pt}}}
   0& 0& 0& 0& t^4 \\ 
   0& 0& 0& -t^2& 0 \\ 
   0& 0& 1& 0& 0 \\ 
   0& -t^{-2}& 0& 0& 0 \\ 
   t^{-4}& 0& 0& 0& 0 
\end{array}\right), \quad
h_1(t;V):\left(\begin{array}{@{\hspace{4pt}}c@{\hspace{4pt}}
c@{\hspace{4pt}}c@{\hspace{4pt}}c@{\hspace{2pt}}c@{\hspace{2pt}}}
   t^4& 0& 0& 0& 0 \\ 
   0& t^2& 0& 0& 0 \\ 
   0& 0& 1& 0& 0 \\ 
   0& 0& 0& t^{-2}& 0 \\ 
   0& 0& 0& 0& t^{-4} 
\end{array}\right).\]
Thus, $n_1(t;V)$ is still ``monomial'' and $h_1(t;V)$ is still ``diagonal'',
but this does not appear to be obvious at all from the above description of
$x_1(t;V)$ and $y_1(t;V)$. And what happens for faithful 
$\fh$-diagonalisable $\slm_2(\C)$-modules $V$ of arbitrary dimension? 
This, and much further information about the possible groups associated 
with $\fg=\slm_2(\C)$, will be discussed in Section~\ref{secA1} at the end 
of this chapter.
\end{exmp}

\begin{rem} \label{xcadet1} Let $n=\dim V$ and $B=\{v_1,\ldots,v_n\}$
be a basis of~$V$ such that each $v_i$ is a weight vector, of weight
$\mu_i\in \Omega$ say. Then $P_\fh(V)=\{\mu_1,\ldots,\mu_n\}$. Since the 
\nm{weight order relation} $\preceq$ is a partial order, we may choose 
the numbering such that, if $\mu_i\preceq \mu_j$ and $\mu_i\neq \mu_j$, 
then $i>j$. We claim:
\begin{itemize}
\item[(a)] For $\alpha\in \Phi^+$ and $t\in \C$, the matrix of 
$x_\alpha(t;V)$ with respect to~$B$ is upper triangular with~$1$ 
along the diagonal.
\item[(b)] For $\alpha\in \Phi^-$ and $t\in \C$, the matrix of 
$x_\alpha(t;V)$ with respect to~$B$ is lower triangular with~$1$ 
along the diagonal.
\end{itemize}
In particular, $\det(x_\alpha(t;V))=1$ for all $\alpha\in \Phi$ and
$t\in \C$. Let us prove (a). By the ``fundamental calculation'' in 
Remark~\ref{fundcal} we have $\be_\alpha^+.v_i\in \rho(\be_\alpha^+)
(V_{\mu_i})\subseteq V_{\mu_i+\alpha}$ for each $i$. So, if 
$\be_\alpha^+.v_i\neq 0$, then $\mu_i+\alpha=\mu_j$ for some $j$. But 
then, since $\alpha\in \Phi^+$, we have $\mu_i\preceq \mu_j$ and $\mu_i
\neq \mu_j$; so we must have $i>j$. Hence, $\be_\alpha^+.v_i$ is a 
linear combination of basis vectors $v_j$ where $i>j$. This means that 
the matrix of $\rho(\be_\alpha^+)$ with respect to $B$ is upper triangular 
with~$0$ along the diagonal. Then the same is also true for the matrix 
of $(t\rho(\be_\alpha^+))^m$, for any $m\geq 1$. Hence, finally, the
matrix of $x_\alpha(t;V)$ with respect to $B$ is the identity matrix 
plus a sum of upper triangular matrices with $0$ along the diagonal, 
as claimed. The proof of (b) is analogous.
\end{rem}

In order to convert the group $G_\C(V)$ in Definition~\ref{nonadj1}
into a group over the field~$K$, the following definition will be crucial.

\begin{defn}[Chevalley, Ree] \label{defadm} Let $V$ be a $\cL$-module, 
with corresponding representation $\rho\colon \cL\rightarrow \gl(V)$. 
Let $n:=\dim V<\infty$. Following Ree \cite[\S 1]{ree2}, a basis $\cB=
\{v_1,\ldots,v_n\}$ of $V$ is called a \nm{regular basis} if the 
following conditions hold:
\begin{itemize}
\item[(A1)] Every vector in $\cB$ is a \nm{weight vector} for $\cH$, that
is, for each~$i$ there exists some $\mu_i\in P_\fh(V)$ such that
$v_i\in V_{\mu_i}$.
\item[(A2)] The matrices of $\frac{1}{m!}\rho(e_i)^m$ and $\frac{1}{m!}
\rho(f_i)^m$ with respect to $\cB$ have entries in $\Z$, for all $m\geq 1$
and $i \in I$.
\end{itemize}
Assume that $\cB$ is such a basis. Note that (A1) shows that $V$ is 
$\fh$-diagonalisable, with $P_\fh(V)=\{\mu_1,\ldots,\mu_n\}$. For any 
$\mu\in P_\fh(V)$, the set $\cB\cap V_{\mu}$ is a basis of $V_{\mu}$. 
Following Chevalley \cite[\S 4]{Ch2}, we say that the $\Z$-module 
$\langle \cB \rangle_\Z \subseteq V$ is an \nm{admissible lattice} 
in~$V$. By (A2), that $\Z$-module is invariant under all endomorphisms 
$\frac{1}{m!} \rho(e_i)^m$ and $\frac{1}{m!}\rho(f_i)^m$ for $m\geq 1$ 
and $i\in I$.
\end{defn}

\begin{rem} \label{defadm1a} (a) As stated by Chevalley \cite[\S 4]{Ch2}, 
for every $\cL$-module~$V$ there exists a regular basis $\cB$ as above. 
The first published proof seems to be due to Ree \cite[Theorem~(1.6)]{ree2}.
See also Steinberg \cite[Chap.~2, Cor.~1]{St}, Borel \cite[Part~A, 
\S 2]{Borel}, Humphreys \cite[\S 27.1]{H} and Bourbaki \cite[Ch.~XIII, 
\S 12]{B78} (in chronological order). We will not need to use this general 
existence result here. For all our purposes, the statement in 
Theorem~\ref{existad} below will be sufficient.

(b) We will see in Proposition~\ref{admxalpha} below that, if the 
condition (A2) holds for the matrices of the Chevalley generators $e_i,
f_i$ ($i \in I$), then it will also hold for the matrices of all 
$\be_\alpha^+\in \cL_\alpha$ ($\alpha\in\Phi$). 
\end{rem}

\begin{rem} \label{defadm1b} As already mentioned above, the condition (A1) 
implies that $V$ is $\fh$-diagonalizable and, hence, that $\rho(e_i)$ and 
$\rho(f_i)$ are nilpotent for all $i\in I$. In a number of cases that we 
shall consider, we actually have $\rho(e_i)^2=\rho(f_i)^2 =0$ for all 
$i\in I$. In that case, (A2) reduces to the condition that the matrices 
of $\rho(e_i)$ and $\rho(f_i)$ with respect to~$\cB$ have entries in~$\Z$, 
for all $i\in I$.
\end{rem}

\begin{exmp} \label{adjadm} Let $V=\cL$ and consider the adjoint 
representation $\ad_\cL \colon \cL\rightarrow \gl(V)$. Let $\bB$
be Lusztig's canonical basis, as in Corollary~\ref{canbash}. It is 
implicit in the proof of Theorem~\ref{luform} that $\bB$ is a regular 
basis of $\cL$, but let us make it completely explicit here. Every 
element of $\bB=\{h_j^+\mid j\in I\} \cup \{\be_\alpha^+\mid \alpha\in
\Phi\}$ certainly is a weight vector for~$\fh$; we have $P_\fh(\fg)=
\Phi\cup\{\underline{0}\}$. In the proof of Lemma~\ref{weyl5} we saw 
that $\ad_\cL(e_i)(h)=-\alpha_i(h)e_i$, $\ad_\cL(f_i)(h)=\alpha_i(h)f_i$ 
and $\ad_\cL(e_i)^2(h)=\ad_\cL(f_i)^2=0$ for all $h\in \fh$. Since 
$h_j^+=-\epsilon(j) h_j$, we obtain
\[ \ad_\cL(e_i)(h_j^+)=-a_{ji}e_i \qquad\mbox{and}\qquad
\ad_\cL(f_i)(h_j^+)=a_{ji}f_i.\]
We certainly have $\ad_\cL(e_i)(e_i)=\ad_\cL(f_i)(f_i)=0$. Furthermore, 
\[ \ad_\cL(e_i)(f_i)=h_i=-\epsilon(i)h_i^+, \qquad \ad_\cL(e_i)^2(f_i)=
-2e_i\]
and $\ad_\cL(e_i)^3(f_i)=0$; similarly, 
\[ \ad_\cL(f_i)(e_i)=-h_i=\epsilon(i)h_i^+, \qquad \ad_\cL(f_i)^2(e_i)=
-2f_i\]
and $\ad_\cL(f_i)^3(e_i)=0$. Hence, we obtain
\[ \textstyle \frac{1}{2}\ad_\cL(e_i)^2(f_i)=-e_i\qquad \mbox{and}\qquad 
\frac{1}{2}\ad_\cL(f_i)^2(e_i)=-f_i.\]
Finally, let $\alpha\in\Phi$ be such that $\alpha\neq \pm \alpha_i$. As 
in the proof of Theorem~\ref{luform}, we see that 
\[ \textstyle \frac{1}{m!}\ad_\cL(e_i)^m(\be_\alpha^+)=\binom{q_{i,\alpha}
+m}{m}\be_{\alpha+m\alpha_i}^+\qquad\mbox{for $1\leq m\leq 
p_{i,\alpha}$},\]
and $\ad_\cL(e_i)^m(\be_\alpha^+)=0$ for $m>p_{i,\alpha}$. Similarly,
\[ \textstyle \frac{1}{m!}\ad_\cL(f_i)^m(\be_\alpha^+)=\binom{p_{i,\alpha}
+m}{m}\be_{\alpha-m\alpha_i}^+\qquad\mbox{for $1\leq m\leq 
q_{i,\alpha}$},\]
and $\ad_\cL(f_i)^m(\be_\alpha^+)=0$ for $m>q_{i,\alpha}$. The above 
formulae show that the matrices of $\frac{1}{m!}\ad_\cL(e_i)^m$ and 
$\frac{1}{m!}\ad_\cL(f_i)^m$ with respect to $\bB$ have entries in $\Z$
for all $m\geq 1$. Hence, $\bB$ is a regular basis.

We have seen in Section~\ref{sec1a6} that $0\leq p_{i,\alpha}+
q_{i,\alpha} \leq 3$. Hence, in any case, we conclude that 
$\ad_\cL(e_i)^4=\ad_\cL(f_i)^4=0$. 
\end{exmp}

\begin{exmp} \label{admAn} Let $n\geq 2$ and $\cL=\slm_n(\C)$. Let $1\leq 
r,s\leq n$, $r\neq s$. As in Section~\ref{sec05}, we denote by $E_{r,s}
\in M_n(\C)$ the elementary matrix with~$1$ at position $(r,s)$, 
and $0$ otherwise. Let 
\[ e_i:=E_{i,i+1}\quad \mbox{and}\quad f_i:=E_{i+1,i} \quad\mbox{for
$i\in I:=\{1,\ldots,n-1\}$}.\]
By Example~\ref{cartsln}, the elements $\{e_i,f_i\mid  i \in I\}$ are 
Chevalley generators for~$\cL$, where $h_i:=[e_i,f_i]$ is the diagonal 
matrix with entries $1,-1$ at positions $i,i+1$ (and $0$ otherwise). 
Let $V=\C^n$ (column vectors) and $\cB=\{b_1,\ldots,b_n\}$ be the standard 
basis of~$V$. The subalgebra $\fh=\langle h_i\mid i\in I\rangle_\C
\subseteq \cL$ consists of diagonal matrices and so it is clear that (A1) 
holds. Now we simply compute that $e_i^2=f_i^2=0$ for $i\in I$; obviously,
the entries of $e_i$ and $f_i$ are integers. Hence, (A2) also holds and 
so $\cB$ is a regular basis. 
\end{exmp}

\begin{exmp} \label{admminus} Let $\Psi\subseteq \Omega$ be a (non-empty)
union of $W$-orbits of non-zero minuscule weights. Let $M$ be $\C$-vector
space with a basis $\cB:=\{z_\mu\mid \mu \in \Psi\}$. Then $M$ becomes
an $\fh$-diagonalisable $\cL$-module where $e_i$ and $f_i$ act via the
formulae in Definition~\ref{defj}; see Proposition~\ref{corj}. Let
$\rho\colon \fg\rightarrow \gl(M)$ be the corresponding representation.
By Remark~\ref{rem10}, we have $\rho(e_i)^2=\rho(f_i)^2=0$ for all $i\in I$.
By Lemma~\ref{lemj}, each basis vector $z_\mu$ is a weight vector. 
Furthermore, the formulae in Definition~\ref{defj} show that the entries 
of the matrices of $\rho(e_i)$ and $\rho(f_i)$ with respect to $\cB$ are 
integers (in fact, only~$0$ and~$1$ occur). Thus, $\cB$ is a regular
basis of $M$.
\end{exmp}

\begin{xca} \label{regdirprod} Let $V$ be a $\fg$-module and $V_1,V_2
\subseteq V$ be submodules such that $V=V_1\oplus V_2$. Show that if $V_1$ 
and $V_2$ are faithful, then~$V$ is also faithful. Furthermore, show that 
if $\cB_1$ is a regular basis of $V_1$ and $\cB_2$ is a regular basis of 
$V_2$, then $\cB:=\cB_1\cup \cB_2$ is a regular basis of~$V$.
\end{xca}

\begin{thm}[Chevalley, Ree] \label{existad} Let $\Lambda'\subseteq 
\Omega$ be any subgroup with $\Phi \subseteq \Lambda'$. Then there 
exists a faithful, $\fh$-diagonalizable $\cL$-module $V$ such that 
$\Lambda(V)=\Lambda'$ and such that $V$ admits a regular basis. 
\end{thm}

\begin{proof} We have seen in the proof of Theorem~\ref{corj2} that there 
is a faithful, $\fh$-diagonalizable $\cL$-module~$V$ such that $\Lambda(V)
=\Lambda'$. Furthermore, $V$ can be written as the direct sum of submodules 
$M_0,M_1,\ldots,M_m$ where $M_0=\cL$ (adjoint module) and each $M_i$ ($i 
\geq 1$) is irreducible with a non-zero minuscule highest weight. Now 
$\cB_0:=\bB$ is a regular basis of~$M_0$; see Example~\ref{adjadm}. 
Furthermore, there is a regular basis $\cB_i$ for each $M_i$; see 
Example~\ref{admminus}. Then $\cB:=\cB_0\cup\cB_1\cup\ldots\cup\cB_m$ is 
a regular basis of~$V$; see Exercise~\ref{regdirprod}. 
\end{proof}

Let us now fix a faithful representation $\rho\colon \cL\rightarrow \gl(V)$ 
which satisfies the two conditions in Definition~\ref{defadm}. Let $K$ be any 
field. We would like to define a corresponding Chevalley group by following 
the general procedure in Section~\ref{sec3a7} and adapting it to the 
present setting. So let $\cB$ be a regular basis of~$V$ and let $V_\Z:=
\langle \cB\rangle_\Z \subseteq V$, a lattice in~$V$. We set 
\[\bar{V}:=K\otimes_{\Z} V_\Z\qquad\mbox{and}\qquad \bar{b}:=1\otimes b
\in \bar{V}\quad \mbox{for $b \in \cB$}.\]
Then $\bar{V}$ is a $K$-vector space and $\bar{\cB}:=\{\bar{b}\mid b 
\in \cB\}$ is a basis of~$\bar{V}$. Some more notation. If $\varphi
\in\End(V)$ is such that $\varphi(V_\Z)\subseteq V_\Z$, then there is an 
induced map $\bar{\varphi}\in \End(\bar{V})$. Let $M_\cB(\varphi)$ be the 
matrix of $\varphi$ with respect to~$\cB$, and $M_{\bar{\cB}}(\bar{\varphi})$ 
be the matrix of $\bar{\varphi}$ with respect to~$\bar{\cB}$. Then all
entries of $M_\cB(\varphi)$ are in~$\Z$ and we have $M_{\bar{\cB}}
(\bar{\varphi})=\overline{M_\cB(\varphi)}$ where, for a matrix $X$ with 
entries in $\Z$, we denote by $\bar{X}$ the matrix obtained by applying the 
canonical map $\Z\rightarrow K$ to the entries of~$X$. 

\begin{rem} \label{remadm2} Let $i\in I$ and $m\geq 0$. We set 
\begin{alignat*}{2}
e_i^{[m]}&:={\textstyle\frac{1}{m!}\rho(e_i)^m}\in \End(V), &\qquad 
\quad E_i^{[m]}&:=M_\cB\bigl(e_i^{[m]}\bigr),\\
f_i^{[m]}&:={\textstyle\frac{1}{m!}\rho(f_i)^m}\in \End(V), &\qquad 
\quad F_i^{[m]}&:=M_\cB\bigl(f_i^{[m]}\bigr).
\end{alignat*}
By (A2), we have $e_i^{[m]}(V_\Z)\subseteq V_\Z$ and so the matrix 
$E_i^{[m]}$ has entries in~$\Z$. Similarly, $f_i^{[m]} (V_\Z)\subseteq 
V_\Z$ and so $F_i^{[m]}$ has entries in $\Z$. Note that $e_i^{[m]}=f_i^{[m]}
=\underline{0}$ for $m\geq n:=\dim V$, since $\rho(e_i)$ and $\rho(f_i)$ 
are nilpotent. With this notation, we can now write
\[ x_i(t;V)=\sum_{m\geq 0} t^m e_i^{[m]} \qquad\mbox{and}\qquad
 y_i(t;V)=\sum_{m\geq 0} t^m f_i^{[m]}\]
for any $t\in \C$. Then 
\[ X_i(t):=\sum_{m\geq 0} t^m E_i^{[m]} \qquad \mbox{and}\qquad 
Y_i(t):=\sum_{m\geq 0} t^m F_i^{[m]}\]
are the matrices of $x_i(t;V)$ and $y_i(t;V)$ with respect to~$\cB$. 
We will usually assume that the elements of $\cB$ are arranged as in 
Remark~\ref{xcadet1}. We have $E_i^{[0]}=F_i^{[0]}=I_n$ (identity 
matrix). Furthermore, for $m\geq 1$, the argument in 
Remark~\ref{xcadet1} shows that
\begin{align*}
E_i^{[m]} &\mbox{ is upper triangular with $0$ along the diagonal},\\
F_i^{[m]} &\mbox{ is lower triangular with $0$ along the diagonal}.
\end{align*}
Now we can pass to $K$. For $\zeta \in K$, we define
\begin{align*}
\bar{x}_i(\zeta;V,\cB)&:=\sum_{m\geq 0} \zeta^m \bar{e}_i^{[m]}
\in \End(\bar{V}),\\
\bar{y}_i(\zeta;V,\cB)&:=\sum_{m\geq 0} \zeta^m \bar{f}_i^{[m]}
\in \End(\bar{V}),
\end{align*}
where $\bar{e}_i^{[m]}\in \End(\bar{V})$ and $\bar{f}_i^{[m]} \in 
\End(\bar{V})$ are the induced linear maps. Note that we also have
$\bar{e}_i^{[m]}=\bar{f}_i^{[m]}=\underline{0}$ for $m\geq n=\dim V$. 
Furthermore, the matrix of $\bar{x}_i(\zeta;V,\cB)$ with respect to 
$\bar{\cB}$ will be upper triangular with~$1$ along the diagonal, and 
the matrix of $\bar{y}_i(\zeta;V,\cB)$ will be lower triangular with~$1$ 
along the diagonal. In particular,
\[ \det\bigl(\bar{x}_i(\zeta;V,\cB)\bigr)=\det\bigl(\bar{y}_i(\zeta;V,\cB)
\bigr)=1 \qquad \mbox{for all $\zeta\in K$}.\]
\end{rem}


\begin{defn}[Chevalley, Ree] \label{defnonad} The subgroup
\[ G_K(V,\cB):=\big\langle \bar{x}_i(\zeta;V,\cB),\bar{y}_i
(\zeta;V,\cB)\mid i\in I,\zeta \in K\big\rangle \subseteq \GL(\bar{V})\]
is called the \nm{Chevalley group} associated with $(V,\cB)$ over $K$.\\
Whenever it is convenient, we will identify elements in $G_K(V,\cB)$ 
with their matrices with respect to $\bar{\cB}$.
\end{defn}

\begin{exmp} \label{defnonad1} (a) Let $K=\C$, $i\in I$ and $t\in \C$. Then,
clearly, $\bar{x}_i(t;V,\cB)$ and $\bar{y}_i(t;V,\cB)$ are the linear maps
$x_i(t;V)$ and $y_i(t;V)$, respectively, as in Definition~\ref{nonadj1}.
Thus, the construction in Definition~\ref{nonadj1} is a special case of
that in Definition~\ref{defnonad}.

(b) Let $V=\fg$ and $\rho$ be the adjoint representation $\ad_\cL \colon 
\cL \rightarrow \gl(V)$. First of all, since $\cL$ is semisimple, we 
have $\ker(\ad_\cL)=Z(\cL)=\{0\}$; so $V$ is faithful. Let $\cB=\bB$ be 
Lusztig's canonical basis of $\cL$. Then $\cB$ is a regular basis by 
Example~\ref{adjadm}. We see that $G_K(\fg,\bB)$ agrees with~$G_K(\fg)$ 
as in Definition~\ref{deflu}.

(c) Let $n\geq 2$ and $\cL=\slm_n(\C)$. Consider the natural representation 
$\rho\colon \cL\hookrightarrow \gl(V)$, where $V=\C^n$ with standard basis
$\cB=\{b_1,\ldots,b_n\}$, and where we identify $\gl(V)=\gl_n(\C)$. Clearly, 
the representation $\rho\colon \cL\hookrightarrow \gl(V)$ is faithful. By 
Example~\ref{admAn}, $\cB$ is a regular basis and 
\[X_i(t)=I_n+t E_{i,i+1}\qquad\mbox{and}\qquad Y_i(t)=I_n+t E_{i+1,i}\]
for $t\in \C$ and $1\leq i \leq n-1$. Consequently, for any $\zeta\in K$,
the matrix of $\bar{x}_i(\zeta;V,\cB)$ with respect to $\bar{\cB}$ is given 
by $I_n+\zeta E_{i,i+1}$, and that of $\bar{y}_i(\zeta;V,\cB)$ is given by 
$I_n+\zeta E_{i+1,i}$. So Proposition~\ref{gensln} shows that $G_K(V,\cB)
\cong \SL_n(K)$. In particular, we see that we do get new groups by the 
construction in Definition~\ref{defnonad}.
\end{exmp}

\begin{exmp} \label{admexp2} Assume that $\rho\colon \fg\rightarrow 
\gl(V)$ is faithful and that $\rho(e_i)^2=\rho(f_i)^2=0$ for all 
$i\in I$. Then $e_i^{[1]}=\rho(e_i)$, $f_i^{[1]}=\rho(f_i)$ and
$e_i^{[m]}=f_i^{[m]}=\underline{0}$ for all $m\geq 2$. It follows that 
\[ \bar{x}_i(\zeta;V,\cB)= \id_{\bar{V}}+\zeta\overline{\rho(e_i)}\quad
\mbox{and}\quad \bar{y}_i(\zeta;V,\cB)=\id_{\bar{V}}+\zeta 
\overline{\rho(f_i)}\]
for all $\zeta\in K$. In this case, the construction of $G_K(V,\cB)$ 
becomes particularly simple. An important class of examples is given
by the modules constructed from a union of orbits of non-zero minuscule
weights, as in Example~\ref{admminus}. 
\end{exmp}

\begin{table}[htbp] \caption{Matrix generators for a Chevalley group
of type $G_2$} \label{matLieG2a}
\begin{center}
{$\renewcommand{\arraystretch}{0.7} \renewcommand{\arraycolsep}{3pt}
\bar{x}_1(\zeta):\left(\begin{array}{ccccccc}
1&.&.&.&.&.&. \\  .&1&\zeta&.&.&.&. \\  .&.&1&.&.&.&. \\ .&.&.&1&.&.&. \\  
.&.&.&.&1&\zeta&. \\ .&.&.&.&.&1&. \\  .&.&.&.&.&.&1  \end{array}\right),\quad
\bar{x}_2(\zeta):\left(\begin{array}{ccccccc}
1&\zeta&.&.&.&.&. \\  .&1&.&.&.&.&. \\  .&.&1&2\zeta&\zeta^2&.&. \\ 
.&.&.&1&\zeta&.&. \\  .&.&.&.&1&.&. \\ .&.&.&.&.&1&\zeta \\  
.&.&.&.&.&.&1  \end{array}\right)$,}\\
{$\renewcommand{\arraystretch}{0.7} \renewcommand{\arraycolsep}{3pt}
\bar{y}_1(\zeta):\left(\begin{array}{ccccccc}
1&.&.&.&.&.&. \\  .&1&.&.&.&.&. \\  .&\zeta&1&.&.&.&. \\ .&.&.&1&.&.&. \\  
.&.&.&.&1&.&. \\ .&.&.&.&\zeta&1&. \\  .&.&.&.&.&.&1 \end{array}\right),\quad
\bar{y}_2(\zeta):\left(\begin{array}{ccccccc}
1&.&.&.&.&.&. \\  \zeta&1&.&.&.&.&. \\  .&.&1&.&.&.&. \\ .&.&\zeta&1&.&.&. 
\\  .&.&\zeta^2&2\zeta&1&.&. \\ .&.&.&.&.&1&. \\  .&.&.&.&.&\zeta&1 
\end{array}\right).$}
{\footnotesize (Here, a dot ``.'' stands for $0$.)}
\end{center}
\end{table}

\begin{exmp} \label{admG2} Let $\cL=\langle e_1,e_2,f_1,f_2
\rangle_{\text{alg}} \subseteq \gl_7(\C)$ be a simple Lie algebra of 
type $G_2$ as in Exercise~\ref{otherG2}, where $e_1,e_2,f_1,f_2$ are 
defined in Table~\ref{matLieG2} (p.~\pageref{matLieG2}). 
Let $V=\C^7$ (column vectors) and $\cB=\{b_1,\ldots,b_7\}$ be the standard 
basis of~$V$. Identifying $\gl(V)=\gl_7(\C)$, the inclusion 
$\cL\hookrightarrow \gl(V)$ is a faithful representation. By 
Exercise~\ref{otherG2}, $\fh=\langle h_1,h_2\rangle_\C$ consists of 
diagonal matrices; so it is clear that condition (A1) in 
Definition~\ref{defadm} holds. Now we simply compute that $e_1^2=
e_2^3=f_1^2=f_2^3=0_{7\times 7}$; furthermore, both $e_2^2$ and $f_2^2$ 
have precisely one non-zero entry, which is~$2$. Hence, $\frac{1}{2}e_2^2$ 
and $\frac{1}{2}f_2^2$ still have integer entries. So (A2) also holds 
and $\cB$ is a regular basis. If $K$ is any field, then the matrices of 
$\bar{x}_i(\zeta;V,\cB)$ and $\bar{y}_i(\zeta;V,\cB)$ are given in 
Table~\ref{matLieG2a}. Here, we do not get a new group (as we shall
see later), but we get a $7$-dimensional realization of the Chevalley 
group of type $G_2$, instead of the $14$-dimensional realization in 
terms of the adjoint representation.
\end{exmp}

\begin{rem} \label{remfieldext} Let $K_1\supseteq K$ be a field extension.
Let us denote $\bar{V}_1:=K_1\otimes_\Z V_\Z$ and $\bar{\cB}_1:=\{1\otimes 
b \mid b \in \cB\}\subseteq \bar{V}_1$. Then we also have the Chevalley 
group $G_{K_1}(V,\cB)\subseteq \GL(\bar{V}_1)$. For $i\in I$ and $\zeta 
\in K$, we can form the elements $\bar{x}_i(\zeta;V,\cB)$ and $\bar{y}_i
(\zeta;V,\cB)$ either within $\GL(\bar{V})$ or within $\GL(\bar{V}_1)$. 
But one immediately sees that the matrices of these elements with respect 
to $\bar{\cB}$ and~$\bar{\cB}_1$, respectively, are exactly the same. 
Consequently, we may naturally identify $G_K(V,\cB)$ with a subgroup of 
$G_{K_1}(V,\cB)$. In particular, this holds if we take for $K_1$ an algebraic
closure of~$K$.
\end{rem}

\begin{xca} \label{xcafieldauto} Let $\sigma\colon K \rightarrow K$ be a 
field automorphism. Show that there is a group automorphism $\tilde{\sigma}
\colon G_K(V,\cB)\rightarrow G_K(V,\cB)$ such that, for all $i\in I$ and 
$\zeta \in K$, we have 
\[\tilde{\sigma}\bigl(\bar{x}_i(\zeta;V,\cB)\bigr)=
\bar{x}_i(\sigma(\zeta);V,\cB),\quad \tilde{\sigma}\bigl(\bar{y}_i(\zeta;
V,\cB)\bigr)= \bar{y}_i(\sigma(\zeta);V,\cB).\]
Such an automorphism of $G_K(V,\cB)$ is called a \nm{field automorphism}.\\
\noindent {\footnotesize [{\it Hint.} By taking matrices with respect to
the basis $\bar{\cB}$ of $\bar{V}$, we may identify $G_K(V,\cB)$ with a 
subgroup of $\GL_n(K)$, where $n=\dim V$. We obtain a group automorphism 
$\GL_n(K)\rightarrow \GL_n(K)$ by applying $\sigma$ to the entries of a 
matrix in $\GL_n(K)$. Check that this automorphism preserves the subgroup
$G_K(V,\cB)\subseteq \GL_n(K)$.]}
\end{xca}

In Section~\ref{sec3a7} it was crucial that we could not only work
with the elements $x_i(t)$, $y_i(t)$ in $G_\C(V)$ and their matrices 
over~$\C$, but also argue at a ``polynomial level''. We now introduce
the required formalism in the present setting.

\begin{rem} \label{niinter0} Let $\Z[T]$ be the ring of polynomials in an 
indeterminate~$T$ over $\Z$. Let $i\in I$. Then $\rho(e_i)^n=\rho(f_i)^n=
\underline{0}$ where $n=\dim V$. Using the notation in Remark~\ref{remadm2}, 
we have 
\[ x_i(t;V)=\sum_{0\leq m\leq n} t^m e_i^{[m]} \qquad \mbox{and} \qquad 
y_i(t;V)=\sum_{0\leq m \leq n} t^m f_i^{[m]},\]
for $t\in \C$. Now we define the following matrices 
with entries in $\Z[T]$:
\[ X_i(T):=\sum_{0\leq m\leq n} T^m E_i^{[m]} \qquad \mbox{and} \qquad 
Y_i(T):=\sum_{0\leq m \leq n} T^m F_i^{[m]}.\]
If $t\in \C$ and we substitute $T\mapsto t$, then we obtain the matrices
\begin{equation*}
X_i(t)=M_\cB\bigl(x_i(t)\bigr)\qquad \mbox{and}\qquad Y_i(t)=M_\cB
\bigl(y_i(t)\bigr).\tag{a}
\end{equation*}
If $\zeta \in K$, then we have a canonical ring homomorphism $\Z[T]
\rightarrow K$ which sends $T$ to~$\zeta$. We denote by $\bar{X}_i
(\zeta)$ and $\bar{Y}_i(\zeta)$ the matrices obtained by applying that 
homomorphism to the entries of $X_i(T)$ and $Y_i(T)$, respectively. Then 
\begin{equation*}
\bar{X}_i(\zeta)=M_{\bar{\cB}}\bigl(\bar{x}_i(\zeta;V,\cB)\bigr)
\quad \mbox{and} \quad \bar{Y}_i(\zeta)=M_{\bar{\cB}}\bigl(\bar{y}_i
(\zeta;V,\cB)\bigr).\tag{b}
\end{equation*}
As in Section~\ref{sec3a7}, the above constructions will be most useful
in order to derive relations in $G_K(V,\cB)$ from similar ones
in $G_\C(V,\cB)$.
\end{rem}

Here is an example which shows how to obtain the above matrices in 
{\sf ChevLie}. We consider a Lie algebra $\cL$ of type $C_2$ and the
representation given by the unique non-zero minuscule highest weight.

\begin{verbatim}
  julia> l=LieAlg(:c,2)         # Lie algebra of type C_2
  julia> l.minuscule            # minuscule weights
  2                             # see Table 13 
  julia> r=rep_minuscule(l,2)   # e_i,f_i,h_i, I={1,2}
  julia> using Nemo
  julia> R,x=polynomial_ring(ZZ,"T")
  julia> x1=expliemat(r[1][1],x); x2=expliemat(r[1][2],x)
            1  0  0  0               1  T  0  0
            0  1  T  0               0  1  0  0 
            0  0  1  0               0  0  1  T
            0  0  0  1               0  0  0  1
  julia> y1=expliemat(r[2][1],x); y2=expliemat(r[2][2],x)
            1  0  0  0               1  0  0  0
            0  1  0  0               T  1  0  0 
            0  T  1  0               0  0  1  0
            0  0  0  1               0  0  T  1
\end{verbatim}

\section{The elements $\bar{x}_\alpha(\zeta;V,\cB)$ in 
$G_K(V,\cB)$} \label{relsGp}

Let us fix a faithful $\fg$-module $V$ and assume that $V$ admits
a regular basis $\cB$; see Definition~\ref{defadm}. Let $K$ be an 
arbitrary field. As in the previous section, let $\bar{V}=K \otimes_\Z
\langle \cB\rangle_\Z$ and $\bar{B}:=\{\bar{b}\mid b \in \cB\}$
where $\bar{b}=1\otimes b\in \bar{V}$. For $i\in I$ and $\zeta\in K$, we 
have the elements
\[ \bar{x}_i(\zeta;V,\cB)\in \GL(\bar{V}) \qquad\mbox{and}\qquad
\bar{y}_i(\zeta;V,\cB)\in \GL(\bar{V});\]
these generate the group $G_K(V,\cB)$; see Definition~\ref{defnonad}. 
For $\zeta\neq 0$ we also have the following elements of $G_K(V,\cB)$:
\begin{align*}
\bar{n}_i(\zeta;V,\cB)&=\bar{x}_i(\zeta;V,\cB)\bar{y}_i(-\zeta^{-1};V,\cB)
\bar{x_i}(\zeta;V,\cB),\\
\bar{h}_i(\zeta;V,\cB)&=\bar{n}_i(\zeta;V,\cB)\bar{n}_i(-1;V,\cB).
\end{align*}
Since $V,\cB$ will be fixed throughout this section, we omit the symbols
$V,\cB$ from the notation and simply denote the above elements by 
$\bar{x}_i(\zeta)$, $\bar{y}_i(\zeta)$, $\bar{n}_i(\zeta)$ and 
$\bar{h}_i(\zeta)$, respectively. If $K=\C$, then we just write 
$x_i(t)$, $y_i(t)$, $n_i(t)$ and $h_i(t)$ for these elements in
$G_\C(V)$. We would like to establish a number of relations among 
these elements, analogous to those in a group of adjoint type, as in 
Sections~\ref{sec3a7}--\ref{secweylG}.

The starting point is the following result, which heavily relies on the 
\nm{Transfer Lemma} (see Lemma~\ref{superexp}) and on Theorem~\ref{luform3}.

\begin{prop} \label{Cree0} Let $i\in I$ and $u\in \C^\times$.
For any $\alpha\in \Phi$, we have 
\[n_i(u)\circ \rho(\be_\alpha^+) \circ n_i(u)^{-1}=c_i(\alpha)u^{-\langle 
\alpha_i^\vee,\alpha\rangle} \rho(\be_{s_i(\alpha)}^+)\]
where $c_i(\alpha)\in \{\pm 1\}$ does not depend on~$u$; we have 
the relation $c_i(\alpha)c_i(-\alpha)=(-1)^{\langle \alpha_i^\vee,\alpha
\rangle}$. Furthermore, $c_i(\alpha_i)=1$ and 
\begin{align*}
n_i(u)\circ\rho(e_i)\circ n_i(u)^{-1}&=-u^{-2} \rho(f_i),\\
n_i(u)\circ\rho(f_i)\circ n_i(u)^{-1}&=-u^{2} \rho(e_i). 
\end{align*}
Finally, we have $n_i(u) \circ \rho(h_\alpha) \circ n_i(u)^{-1}=
\rho(h_{s_i(\alpha)})$ for $\alpha\in \Phi$.
\end{prop}
%

\begin{proof} In order to avoid any danger of confusion, we denote
\[ x_i^\fg(t):=\exp(t\, \ad_\cL(e_i)) \qquad \mbox{and} \qquad
y_i^\fg(t):=\exp(t\, \ad_\cL(f_i))\]
for $i \in I$ and $t \in \C$; these are elements of $G_\C(\fg)\subseteq 
\GL(\fg)$, as in Section~\ref{sec3a7}. We also set $n_i^\fg(t):=x_i^\fg(t)
y_i^\fg(-t^{-1})x_i^\fg(t)\in G_\C(\fg)$ for $t \neq 0$. Now consider the 
present set-up, where 
\[n_i(u)=x_i(u)y_i(-u^{-1})x_i(u) \in G_\C(V).\]
We will use the {\itshape Transfer Lemma} three times. First, we apply it 
with $x=ue_i$ and a given $y\in \fg$. With the above notation, this yields 
that
\[ \rho\bigl(x_i^\fg(u)(y)\bigr)=x_i(u)\circ \rho(y) \circ x_i(u)^{-1}.\]
Next, let $x=-u^{-1}f_i$ and $y'=x_i^\fg(u)(y) \in \fg$. Then we obtain
\begin{align*}
\rho\bigl(y_i^\fg&(-u^{-1})(y')\bigr)=y_i (-u^{-1})\circ \rho(y')
\circ y_i(-u^{-1})^{-1}\\
&=y_i(-u^{-1})\circ \bigl(x_i(u)\circ \rho(y)
\circ x_i(u)^{-1}\bigr)\circ y_i(-u^{-1})^{-1}.
\end{align*}
Finally, since $n_i^\fg(u)(y)=x_i^\fg(u)\bigl(y_i^\fg(-u^{-1})(y')\bigr)$, 
a third application with $x=ue_i$ and $y''= y_i^\fg(-u^{-1})y'\in \fg$ 
yields 
\begin{equation*}
\rho\bigl(n_i^\fg(u)(y)\bigr)=n_i(u)\circ\rho(y)\circ n_i(u)^{-1}.\tag{$*$}
\end{equation*}
Now let $y=\be_\alpha^+$ where $\alpha\in\Phi$. By 
Theorem~\ref{luform3}, we have 
\[n_i^\fg(u)(\be_\alpha^+)=c_i(\alpha)u^{-\langle \alpha_i^\vee,
\alpha\rangle} \be_{s_i(\alpha)}^+,\]
where $c_i(\alpha)\in \{\pm 1\}$ does not depend on~$u$. So ($*$) yields
the desired formula for $y=\be_\alpha^+$. Actually, we have $c_i(\alpha)=
(-1)^{q_{i,\alpha}}$ where $q_{i,\alpha}$ is defined via the
$\alpha_i$-string through~$\alpha$. Since $q_{i,-\alpha}=p_{i,\alpha}$ 
and $\langle \alpha_i^\vee,\alpha\rangle=q_{i,\alpha}-p_{i,\alpha}$, we 
conclude that $c_i(\alpha)c_i(-\alpha)=(-1)^{\langle \alpha_i^\vee,
\alpha\rangle}$, as claimed. 

Now take $y=e_i$. We have $e_i=\epsilon(i)\be_{\alpha_i}^+$ and 
$f_i=-\epsilon(i) \be_{-\alpha_i}^+$; furthermore, $s_i(\alpha_i)=
-\alpha_i$ and $q_{i,\alpha_i}=2$. Hence, $c_i(\alpha_i)=1$ and 
$n_i(u)(e_i)=-u^{-2}f_i$. So ($*$) yields the desired formula for 
$y=e_i$. Similarly, since $q_{i,-\alpha_i}=0$, we obtain the desired 
formula for $y=f_i$,

Finally, by Proposition~\ref{weyl5c}, we have $n_i^\fg(u)(h_\alpha)=
h_{s_i(\alpha)}$. So, again, ($*$) yields the desired formula for 
$y=h_\alpha$.
\end{proof}

There are a number of applications. The first one is a strengthening of 
the conditions concerning ``admissibility'' in Definition~\ref{defadm}. 

\begin{prop} \label{admxalpha} The condition {\rm (A2)} in 
Definition~\ref{defadm} also holds for $\rho(\be_\alpha^+)$, that is, 
the matrix of $\frac{1}{m!}\rho(\be_\alpha^+)^m$ with respect to $\cB$ 
has entries in $\Z$, for any $\alpha\in \Phi$ and any $m\geq 1$. 
\end{prop}

\begin{proof} There exists some $w\in W$ and $i\in I$ such that $\alpha=
w(\alpha_i)$. Write $w=s_{i_1}\cdots s_{i_r}$ where $r\geq 0$ and $i_1,
\ldots,i_r\in I$. Then set $\eta:=n_{i_1}(1) \cdots n_{i_r}(1) \in 
G_\C(V)$. By a repeated appliction of Proposition~\ref{Cree0}, we obtain 
\[\eta\circ \rho(e_i)\circ \eta^{-1}=\eta\circ \rho(\pm \be_{\alpha_i}^+)
\circ \eta^{-1}=\pm \rho(\be_{w(\alpha_i)}^+)=\pm \rho(\be_\alpha^+).\]
Hence, we also have 
$\textstyle \eta\circ \bigl(\frac{1}{m!}\rho(e_i)^m\bigr)\circ
\eta^{-1}=\pm \frac{1}{m!}\rho(\be_\alpha^+)^m$ for any $m\geq 1$.

Now, for any $j\in I$, the matrices of $x_j(\pm 1)$ and 
$y_j(\pm 1)$ with respect to $\cB$ have entries in $\Z$ and determinant 
equal~$1$. So an analogous statement also holds for the matrix of 
each $n_j(1)=x_j(1)y_j(-1)x_j(1)$ and, hence, also for the matrix 
of~$\eta$. Finally, since $\det(\eta)=1$, the matrix of $\eta^{-1}$ will 
also have entries in $\Z$. Hence, since (A2) holds for $\rho(e_i)$, it 
follows that (A2) also holds for $\rho(\be_\alpha^+)$.
\end{proof}

We can now extend the notation in Remark~\ref{remadm2} to all roots
in~$\Phi$. For any $\alpha\in \Phi$ and $m\in \Z_{\geq 0}$ we set 
\[\be_\alpha^{[m]}:={\textstyle\frac{1}{m!}\rho(\be_\alpha^+)^m}\in 
\End(V),\qquad \quad \bE_\alpha^{[m]}:=M_\cB\bigl(\be_\alpha^{[m]}\bigr).\]
Then $\be_\alpha^{[m]}(V_\Z) \subseteq V_\Z$ and $\bE_\alpha^{[m]}$ has 
entries in~$\Z$. Note again that $\be_\alpha^{[m]}=\underline{0}$ for 
$m\geq \dim V$ since $\rho(\be_\alpha^+)$ is nilpotent. Thus, we have
\[x_\alpha(t;V)=\exp\bigl(t\rho(\be_\alpha^+)\bigr)=\sum_{m\geq 0} t^m 
\be_\alpha^{[m]} \qquad \mbox{for any $t\in \C$}.\]
Now turn to the field $K$ and $\bar{V}=K\otimes V_\Z$. Let 
$\bar{\be}_\alpha^{[m]} \in \End(\bar{V})$ be the induced linear map 
and $\bar{\bE}_\alpha^{[m]}:=M_{\bar{\cB}}\bigl(\bar{\be}_\alpha^{[m]}
\bigr)$. Then we define 
\[ \bar{x}_\alpha(\zeta;V,\cB):=\sum_{m \geq 0} \zeta^m 
\bar{\be}_{\alpha}^{[m]}\in \End(\bar{V}) \qquad \mbox{for 
$\zeta \in K$}.\]
As before, we see that $\det(\bar{x}_\alpha(\zeta;V,\cB))=1$ and so 
$\bar{x}_\alpha(\zeta;V,\cB)\in \GL(\bar{V})$. In the following, if there 
is no danger of confusion, we simply write $x_\alpha(t)$ instead of 
$x_\alpha(t;V)$ and $\bar{x}_\alpha(\zeta)$ instead of $\bar{x}_\alpha
(\zeta;V,\cB)$. 

\begin{exmp} \label{admG2a} Consider the $7$-dimensional realisation
of a Chevalley group of type $G_2$ from Example~\ref{admG2}. Here,
we have 
\[ \Phi^+=\{\alpha_1,\alpha_2,\alpha_1+\alpha_2,
\alpha_1+2\alpha_2,\alpha_1+3\alpha_2,2\alpha_1+3\alpha_2\};\]
see Example~\ref{cartanG2}. Furthermore, let us fix the function
$\epsilon \colon I\rightarrow\{\pm1\}$ as in Table~\ref{Mdynkineps} 
(p.~\pageref{Mdynkineps}), that is, $\epsilon(1)=1$ and $\epsilon(2)=
-1$. The matrices of the elements $\bar{x}_\alpha(\zeta)$, for 
$\alpha\in \Phi^+$, are displayed in Table~\ref{matLieG2b}. 
\end{exmp}

\begin{table}[htbp] \caption{$\bar{x}_\alpha(\zeta)$, $\alpha
\in \Phi^+$, for a Chevalley group of type $G_2$} \label{matLieG2b}
\begin{center}
{\small $\renewcommand{\arraystretch}{0.7} 
\renewcommand{\arraycolsep}{2pt}
\begin{array}{ll}
\multicolumn{2}{l}{\bar{x}_{\alpha_1}(\zeta)=\bar{x}_1(\zeta), \qquad\quad
\bar{x}_{\alpha_2}(\zeta)=\bar{x}_2(-\zeta), \qquad\text{(see
Table~\ref{matLieG2a})}}\\\\
\bar{x}_{\alpha_1+\alpha_2}(\zeta)=& \quad \bar{x}_{\alpha_1+2\alpha_2}
(\zeta)= \\[5pt]
\renewcommand{\arraycolsep}{2pt}
\qquad \left(\begin{array}{ccccccc}
 1  &    . &   \zeta &     .  &    .&       . &    .\\
 .  &    1 &     . & -2\zeta  &   . & -\zeta^2  &    .\\
 .  &    . &     1 &     .  &    .&       . &    .\\
 .  &    . &     . &     1  &    .&     \zeta &    .\\
 .  &    . &     . &     .  &    1&       . & -\zeta\\
 .  &    . &     . &     .  &    .&       1 &    .\\
 .  &    . &     . &     .  &    .&       . &    1
\end{array}\right),
\renewcommand{\arraycolsep}{2pt}
&\quad\qquad \left(\begin{array}{ccccccc} 
 1  &    . &     . & -2\zeta  &   . &    .  & \zeta^2 \\
 .  &    1 &     . &     .  &  \zeta&     . &     . \\
 .  &    . &     1 &     .  &    .&   \zeta &     . \\
 .  &    . &     . &     1  &    .&     . &  -\zeta \\
 .  &    . &     . &     .  &    1&     . &     . \\
 .  &    . &     . &     .  &    .&     1 &     . \\
 .  &    . &     . &     .  &    .&     . &     1 
\end{array}\right),\\ \\
\bar{x}_{\alpha_1+3\alpha_2}(\zeta)=& \quad \bar{x}_{2\alpha_1+3\alpha_2}
(\zeta)= \\[5pt]
\renewcommand{\arraycolsep}{3pt}
\qquad \left(\begin{array}{ccccccc}
  1&    .&    .&    .&  \zeta&     .&    . \\
  .&    1&    .&    .&       .&    .&    . \\
  .&    .&    1&    .&       .&    .& -\zeta \\
  .&    .&    .&    1&       .&    .&    . \\
  .&    .&    .&    .&       1&    .&    . \\
  .&    .&    .&    .&       .&    1&    . \\
  .&    .&    .&    .&       .&    .&    1 
\end{array}\right),
\renewcommand{\arraycolsep}{5pt}
&\quad\qquad \left(\begin{array}{ccccccc} 
  1&    .&    .&    .&   .& -\zeta&    . \\
  .&    1&    .&    .&   .&      .& -\zeta \\
  .&    .&    1&    .&   .&      .&    . \\
  .&    .&    .&    1&   .&      .&    . \\
  .&    .&    .&    .&   1&      .&    . \\
  .&    .&    .&    .&   .&      1&    . \\
  .&    .&    .&    .&   .&      .&    1 
\end{array}\right).\\\\
\multicolumn{2}{c}{\mbox{\footnotesize (Here, a dot ``.'' stands 
for $0$.)}} \end{array}$}
\end{center}
\end{table}

The proofs of the following results will involve arguments at the
``polynomial level''. For $\alpha\in \Phi$ we define
\[ X_\alpha(T):=\sum_{m \geq 0} T^m \bE_\alpha^{[m]},\]
a matrix with entries in the polynomial ring $\Z[T]$ where $T$ is an 
indeterminate. Upon substituting $T\mapsto t$ for any $t\in\C$, we obtain
the matrix $X_\alpha(t):=M_\cB(x_\alpha(t))$. If $\zeta \in K$, then we 
can apply the canonical ring homomorphism $\Z[T] \rightarrow K$ which sends 
$T$ to~$\zeta$; this yields the matrix $\bar{X}_\alpha(\zeta):=M_{\bar{\cB}}
(\bar{x}_\alpha(\zeta))$.

\begin{prop} \label{action2a} {\rm (a)} Let $\alpha\in \Phi$. Then 
$\bar{x}_\alpha(\zeta+\xi)=\bar{x}_\alpha(\zeta)\bar{x}_\alpha(\xi)$ for 
all $\zeta,\xi \in K$. Furthermore, $\bar{x}_\alpha(0)=\id_{\bar{V}}$ and 
$\bar{x}_\alpha(\zeta)^{-1}=\bar{x}_\alpha(-\zeta)$.

\noindent {\rm (b)} We have $\bar{x}_i(\zeta+\xi)=\bar{x}_i(\zeta)
\bar{x}_i(\xi)$ and $\bar{y}_i(\zeta+\xi)=\bar{y}_i(\zeta)\bar{y}_i(\xi)$ 
for all $i\in I$ and $\zeta,\xi \in K$. Furthermore, $\bar{n}_i(\xi)^{-1}
= \bar{n}_i(-\xi)$ for $\xi\in K^\times$. 

\noindent {\rm (c)} Let $\alpha ,\beta\in \Phi$ be such that $\beta\neq 
-\alpha$ and $\alpha+\beta \not\in \Phi$. Then 
\[ \bar{x}_\alpha(\zeta)\bar{x}_\beta(\xi)=\bar{x}_\beta(\xi)
\bar{x}_\alpha(\zeta)\qquad \mbox{for all $\zeta,\xi\in K$}.\]
\end{prop}

\begin{proof} (a) First we work over $\C$. We clearly have $x_\alpha(0)=
\id_V$. Let $t,u\in \C$. By Lemma~\ref{exponential}, we have 
\[ x_\alpha(t)^{-1}=\exp(t\rho(\be_\alpha^+))^{-1}=\exp(-t\rho
(\be_\alpha^+))=x_\alpha(-t).\]
Furthermore, $(t+u)\rho(\be_\alpha^+)=t\rho(\be_\alpha^+)+u\rho
(\be_\alpha^+)$ and the two summands on the right hand side commute with 
each other. Hence, we obtain 
\[ \exp((t+u)\rho(\be_\alpha^+))=\exp(t\rho(\be_\alpha^+))\circ 
\exp(u\rho(\be_\alpha^+))\]
by Exercise~\ref{expocomm}, that is, $x_\alpha(t+u)=x_\alpha(t)x_\alpha
(u)$. Now we can pass to $K$. Once the above identity is established for
all $t,u\in \C$, we obtain an identity at the ``polynomial level'':
\[ X_\alpha(T+U)=X_\alpha(T){\cdot}X_\alpha(U),\]
where we work over the ring of polynomials $\Z[T, U]$ in two commuting 
indeterminates $T,U$. (See the proof of Lemma~\ref{fromctok1} for a 
similar argument.) Given $\zeta,\xi\in K$, it then remains to apply 
the canonical ring homomorphism $\Z[T, U]\rightarrow K$ which sends 
$T$ to~$\zeta$ and $U$ to~$\xi$.

(b) Note that $\be_{\alpha_i}^+=\epsilon(i)e_i$ and $\be_{-\alpha_i}^+=
-\epsilon(i)f_i$; hence, we have
\[ \bar{x}_i(\zeta)=\bar{x}_{\alpha_i}(\epsilon(i) \zeta)\qquad 
\mbox{and}\qquad \bar{y}_i(\zeta)=\bar{x}_{-\alpha_i}(-\epsilon(i)
\zeta).\]
So the statements concerning $\bar{x}_i(\zeta+\xi)$ and 
$\bar{y}_i(\zeta+\xi)$ are an immediate consequence of~(a).
Furthermore, if $\xi \neq 0$, then 
\begin{align*}
\bar{n}_i(\xi)^{-1}&=\bigl(\bar{x}_i(\xi)\bar{y}_i(-\xi^{-1})
\bar{x}_i(\xi)\bigr)^{-1}= \bar{x}_i(\xi)^{-1}\bar{y}_i
(-\xi^{-1})^{-1}\bar{x}_i(\xi)^{-1}\\&= \bar{x}_i(-\xi)\bar{y}_i
(\xi^{-1})\bar{x}_i(-\xi)= \bar{n}_i(-\xi)
\end{align*}
where the third equality holds again by (a).

(c) As in (a), it is enough to prove this over $\C$. Let $t,u\in \C$. Since 
$\alpha+\beta\not\in \Phi$ and $\beta \neq -\alpha$, we have $[\be_\alpha^+,
\be_\beta^+]=0$. Hence, since the map $\rho\colon\cL\rightarrow \gl(V)$ is 
a Lie algebra homomorphism, we also have $[\rho(\be_\alpha^+),
\rho(\be_\beta^+)]=0$, that is, $\rho(\be_\alpha^+)\colon V\rightarrow V$ 
and $\rho(\be_\beta^+) \colon V\rightarrow V$ commute with each other. 
Consequently, the maps $(t\rho(\be_\alpha^+))^m$ and 
$(u\rho(\be_\beta^+))^l$ commute with each other for all $m,l\in 
\Z_{\geq 0}$. So, finally,
\[ x_\alpha(t)=\sum_{m\geq 0} \frac{(t\, \ad_\cL(\be_\alpha^+))^m}{m!} 
\quad\mbox{and} \quad x_\beta(u)=\sum_{l\geq 0} \frac{(u\, 
\ad_\cL(\be_\beta^+))^l}{l!} \]
also commute with each other, as claimed.
\end{proof}

\begin{prop} \label{Cree0a} Let $i \in I$, $\alpha\in \Phi$, $\zeta\in K$
and $\xi \in K^\times$. Then 
\begin{align*}
\bar{n}_i(\xi)\bar{x}_\alpha(\zeta)\bar{n}_i(\xi)^{-1}&=\bar{x}_{s_i
(\alpha)}\bigl(c_i(\alpha)\zeta\xi^{-\langle \alpha_i^\vee,\alpha\rangle}
\bigr),\tag{a}\\ 
\bar{h}_i(\xi)\bar{x}_\alpha(\zeta)\bar{h}_i(\xi)^{-1}&=\bar{x}_\alpha
\bigl(\zeta\xi^{\langle \alpha_i^\vee,\alpha\rangle}\bigr), \tag{b}
\end{align*}
with $c_i(\alpha)\in \{\pm 1\}$ as in Proposition~\ref{Cree0}. 
\end{prop}

\begin{proof} First we work over $\C$. Let $t\in \C$ and $u\in \C^\times$.
We have 
\[n_i(u)\circ\rho(\be_\alpha^+)\circ n_i(u)^{-1}=c_i(\alpha)
u^{-\langle \alpha_i^\vee,\alpha\rangle}\rho(\be_{s_i(\alpha)}^+)\]
by Proposition~\ref{Cree0}. Hence, for any $m\geq 0$, we have 
\[n_i(u)\circ \rho(t\,\be_\alpha^+)^m \circ n_i(u)^{-1}=c_i(\alpha)^m
u^{-\langle \alpha_i^\vee, \alpha\rangle m}\rho(t\,\be_{s_i(\alpha)}^+)^m.\]
Multiplying each such term by $\frac{1}{m!}$ and then summing over all 
$m\geq 0$ yields the first of the above identities. Now consider the 
second identity. Using the formula $h_i(u)=n_i(u)n_i(-1)$, we obtain:
\begin{align*}
h_i(u)x_\alpha(t)h_i(u)^{-1}&=n_i(u)\bigl(x_{s_i(\alpha)}\bigl(c_i(\alpha)
(-1)^{\langle \alpha_i^\vee,\alpha\rangle}t\bigr)\bigr)n_i(u)^{-1}\\&=
x_{\alpha}\bigl(c\,tu^{-\langle \alpha_i^\vee,s_i(\alpha)\rangle}\bigr),
\end{align*}
where, by (a), $c:=c_i(\alpha)c_i(s_i(\alpha))(-1)^{\langle \alpha_i^\vee,
\alpha\rangle}\in \{\pm 1\}$. So it remains to show that $c=1$ and 
$\langle \alpha_i^\vee,s_i(\alpha) \rangle=-\langle \alpha_i^\vee,
\alpha\rangle$. Now, we do have 
\[ \langle \alpha_i^\vee,s_i(\alpha)\rangle= 
2\frac{\langle \alpha_i,s_i(\alpha)\rangle}{\langle \alpha_i,
\alpha_i\rangle}=2\frac{\langle s_i(\alpha_i),\alpha\rangle}{\langle 
\alpha_i, \alpha_i\rangle}=-2\frac{\langle \alpha_i,\alpha\rangle}{\langle 
\alpha_i, \alpha_i\rangle}=-\langle \alpha_i^\vee,\alpha\rangle\]
as required, where we used the fact that $\langle\;,\;\rangle$ is 
$W$-invariant. Furthermore, since $n_i(-u)=n_i(u)^{-1}$, 
Proposition~\ref{Cree0} also shows that 
\[n_i(u)^{-1}\circ\rho(\be_{s_i(\alpha)}^+)\circ n_i(u)=c_i(s_i(\alpha))
(-u)^{-\langle \alpha_i^\vee,s_i(\alpha)\rangle}\rho(\be_\alpha^+).\]
Combining this with the above identity for $n_i(u)\circ\rho(\be_\alpha^+)
\circ n_i(u)^{-1}$, we find that $c_i(s_i(\alpha))=c_i(\alpha)
(-1)^{\langle \alpha_i^\vee,s_i(\alpha)\rangle}$ and, hence, $c=1$.
 
Now we lift the above relations to a ``polynomial level''.  We consider 
the ring $\Z[T, U^{\pm 1}]$ in two commuting indeterminates~$T,U$. Using
the matrices $X_i(T)$, $Y_i(T)$ in Remark~\ref{niinter0}, we define 
\begin{align*}
N_i(U)&:=X_i(U)\cdot Y_i(-U^{-1})\cdot X_i(U),\\
H_i(U)&:=N_i(U)\cdot N_i(-1).
\end{align*}
Then, upon substituting $U\mapsto u$ for $u\in \C^\times$, we obtain the 
matrices 
\[N_i(u):=M_{\cB}(n_i(u))\qquad\mbox{and}\qquad H_i(u):=M_{\cB}(h_i(u)).\]
Hence, the relations that we proved over $\C$ now read as follows:
\begin{align*}
N_i(u)\cdot X_\alpha(t)\cdot N_i(u)^{-1}&=X_{s_i(\alpha)}\bigl(
c_i(\alpha) tu^{-\langle \alpha_i^\vee,\alpha\rangle}\bigr),\\
H_i(u)\cdot X_\alpha(t)\cdot H_i(u)^{-1}&=X_\alpha \bigl(tu^{\langle 
\alpha_i^\vee, \alpha\rangle}\bigr).
\end{align*}
Since these hold for all $t\in \C$ and all $u\in \C^\times$, we also have 
the analogous identities for matrices over $\Z[T,U^{\pm 1}]$:
\begin{align*}
N_i(U)\cdot X_\alpha(T)\cdot N_i(U)^{-1}&=X_{s_i(\alpha)}\bigl(
c_i(\alpha) TU^{-\langle \alpha_i^\vee,\alpha\rangle}\bigr),\\
H_i(U)\cdot X_\alpha(T)\cdot H_i(U)^{-1}&=X_\alpha \bigl(TU^{\langle 
\alpha_i^\vee, \alpha\rangle}\bigr),
\end{align*}
Now, given $\zeta\in K$ and $\xi\in K^\times$, we have a canonical ring 
homomorphism $\Z[T, U^{\pm 1}]\rightarrow K$ which sends $T$ to~$\zeta$ 
and $U$ to~$\xi$. Applying this to the above identities, we obtain analogous 
identities for the matrices of $\bar{n}_i(\xi)$, $\bar{h}_i(\xi)$ and 
$\bar{x}_\alpha(\zeta)$ with respect to~$\bar{\cB}$, as desired. 
\end{proof}

\begin{rem} \label{Cree0a0} We also note the following formulae, which are 
the exact analogue of Lemma~\ref{barni1}. Let $\zeta \in K$ and $\xi\in 
K^\times$. Then 
\begin{align*}
\bar{n}_i(\xi)\bar{x}_i(\zeta)\bar{n}_i(\xi)^{-1}&=\bar{y}_i(-\zeta
\xi^{-2}),\\ \bar{n}_i(\xi)\bar{y}_i(\zeta)\bar{n}_i(\xi)^{-1}&=
\bar{x}_i(-\zeta\xi^{2}).
\end{align*}
Indeed, as noted in the proof of Proposition~\ref{action2a}(b), we have 
$\bar{x}_i(\zeta)=\bar{x}_{\alpha_i}(\epsilon(i)\zeta)$ and 
$\bar{y}_i(\zeta)=\bar{x}_{-\alpha_i}(-\epsilon(i) \zeta)$. So the above 
identities simply follow from Proposition~\ref{Cree0a}; we have 
$c_i(\alpha_i)=1$ by Proposition~\ref{Cree0}. Similarly, we have
\begin{align*}
\bar{h}_i(\xi)\bar{x}_i(\zeta)\bar{h}_i(\xi)^{-1}&=\bar{x}_i
\bigl(\zeta\xi^2\bigr),\\
\bar{h}_i(\xi)\bar{y}_i(\zeta)\bar{h}_i(\xi)^{-1}&=\bar{y}_i
\bigl(\zeta\xi^{-2}\bigr),
\end{align*}
which are the exact analogue of Corollary~\ref{barni2}. (Since 
$\bar{h}_i(\xi)=\bar{n}_i(\xi)\bar{n}_i(-1)$, the latter two identities 
immediately follow from the previous two identities.) Consequently, by 
exactly the same argument as in Corollary~\ref{kchar4}, it follows that
$G_K(V,\cB)$ is equal to its own commutator subgroup if $|K|\geq 4$.
\end{rem}

Let us just write $\bar{G}=G_K(V,\cB)$ from now on. 

\begin{cor} \label{corCree0a} We have $\bar{x}_\alpha(\zeta)\in 
\bar{G}$ for $\alpha\in \Phi$ and $\zeta\in K$. In particular, $\bar{G}
=\langle \bar{x}_\alpha(\zeta)\mid \alpha\in \Phi,\zeta\in K \rangle
\subseteq \GL(\bar{V})$. 
\end{cor}

\begin{proof} Given $\alpha\in \Phi$, there exist $w\in W$ and $i\in I$ 
such that $\alpha=w(\alpha_i)$. Now write $w=s_{i_1}\cdots s_{i_r}$ 
where $i_1,\ldots,i_r\in I$. Then consider $\eta:=\bar{n}_{i_1}\cdots 
\bar{n}_{i_r}\in \bar{G}$, where we set $\bar{n}_i:=\bar{n}_i(1)$
for any $i\in I$. Let $\zeta \in K$. By Proposition~\ref{Cree0a}, we 
have $\bar{n}_i\bar{x}_\beta(\zeta) \bar{n}_i^{-1}=\bar{x}_{s_i(\beta)}
(\pm \zeta)$ for all $i \in I$ and $\beta \in \Phi$. Using this rule 
we obtain that 
\begin{align*}
\eta\,\bar{x}_{\alpha_i}(\zeta)\,\eta^{-1}&=
\bar{n}_{i_1}\cdots \bar{n}_{i_r}\bar{x}_{\alpha_i}(\zeta)
\bar{n}_{i_r}^{-1}\cdots \bar{n}_{i_1}^{-1}\\
&=\bar{n}_{i_1}\cdots \bar{n}_{i_{r-1}}\bar{x}_{s_{i_r}(\alpha_i)}
(\pm \zeta)\bar{n}_{i_{r-1}}^{-1} \cdots \bar{n}_{i_1}^{-1}.
\end{align*}
Repeating the argument with $\bar{n}_{i_{r-1}}$ and so on, eventually 
we find that 
\[ \eta\,\bar{x}_{\alpha_i}(\zeta)\,\eta^{-1}=\bar{x}_{s_{i_1}\cdots 
s_{i_r}(\alpha_i)}(\pm \zeta)=\bar{x}_{w(\alpha_i)}(\pm \zeta)=
\bar{x}_\alpha(\pm \zeta)\]
(where the sign may change at each step, but the signs do not depend
on~$\zeta$). Since $\eta\in \bar{G}$ and $\bar{x}_{\alpha_i}(\zeta)=
\bar{x}_i(\epsilon(i)\zeta)\in \bar{G}$, we conclude that $\bar{x}_\alpha
(\pm \zeta)=\eta\,\bar{x}_{\alpha_i}(\zeta)\,\eta^{-1}\in \bar{G}$. Since 
this holds for all $\zeta\in K$, we also have $\bar{x}_\alpha(\zeta)
\in \bar{G}$ for all $\zeta\in K$.
\end{proof}

We now come to further applications of Proposition~\ref{Cree0}.

\begin{prop} \label{nontrivrhoi} Let $\alpha\in \Phi$. Then the greatest 
common divisor of the entries of the matrix $\bE_\alpha^{[1]}=M_{\cB}
(\rho(\be_\alpha^+))$ is~$1$. Hence, the induced linear map
$\overline{\rho(\be_\alpha^+)} \in \End(\bar{V})$ is non-zero.
\end{prop}

\begin{proof} Let $i\in I$ and $\eta\in G_\C(V,\cB)$ be as in the above
proof. Then $\rho(e_i)=\pm \eta^{-1}\circ \rho(\be_\alpha^+)\circ \eta$. 
Hence, since the corresponding matrices have all their entries in $\Z$,
it is sufficient to prove the assertion about the greatest common 
divisor for the matrix $E_i^{[1]}=M_\cB(\rho(e_i))$. Assume the 
result is false. Then there is a prime number~$p$ which divides all entries 
of $E_i^{[1]}$. As in the proof of Proposition~\ref{admxalpha}, the matrices 
of $n_i(1)$ and $n_i(1)^{-1}=n_i(-1)$ with respect to~$\cB$ have integer 
entries. By Proposition~\ref{Cree0}, we have $n_i(1)\circ \rho(e_i) 
\circ n_i(1)^{-1}=-\rho(f_i)$. So we conclude that all entries of 
$F_i^{[1]}=M_\cB(\rho(f_i))$ must also be divisible by~$p$. But then 
$p^2$ will divide all entries of
\[ M_\cB(\rho(h_i))=M_\cB(\rho([e_i,f_i]))=E_i^{[1]}\circ F_i^{[1]}-
F_i^{[1]}\circ E_i^{[1]}.\]
On the other hand, $M_\cB(\rho(h_i))$ is a diagonal matrix. Thus, all 
eigenvalues of $\rho(h_i)$ are seen to be integers that are divisible 
by~$p^2$.

Now consider the subalgebra $\fg_i=\langle e_i,f_i,h_i\rangle_\C\subseteq
\cL$ and regard $V$ as a $\fg_i$-module (by restriction). There is a chain 
of $\fg_i$-modules 
\[\{0\}=V_0 \subsetneqq V_1 \subsetneqq V_2 \subsetneqq \ldots \subsetneqq 
V_r=V\]
such that the factor modules $V_l/V_{l-1}$ are irreducible for $1\leq l 
\leq r$. Denote by $\rho_l\colon \fg_i\rightarrow \gl(V_l/V_{l-1})$ the
representation corresponding to $V_l/V_{l-1}$. Since $\fg_i\cong\slm_2(\C)$, 
we can apply the results in Section~\ref{sec04a}. Arguing as in 
Proposition~\ref{sl2modd}, we see that $c\in \C$ is an eigenvalue of 
$\rho(h_i)$ if and only if $c$ is an eigenvalue of $\rho_l(h_i)$ for 
some~$l\in\{1,\ldots,r\}$. 

If $\dim V_l/\dim V_{l-1}=1$ for all~$l$, then $\rho_l(h_i)=\underline{0}$
for all~$l$; see Corollary~\ref{sl2modc}. Consequently, $0$ is the only 
eigenvalue of $\rho(h_i)$. Since $\rho(h_i)$ is diagonalizable, this 
would imply that $\rho(h_i)=\underline{0}$, contradiction to $\rho$ being 
injective. So there is some $l$ such that $\dim V_l/V_{l-1}\geq 2$. But 
then Corollary~\ref{sl2modc} shows that $1$ or $2$ is an eigenvalue of 
$\rho_l(h_i)$. Consequently, $1$ or $2$ will also be an eigenvalue of 
$\rho(h_i)$, contradiction to all eigenvalues of $\rho(h_i)$ being 
divisible by~$p^2$.
\end{proof}

To state the next result, we use the following notation. Given 
any $\varphi\in \End(V)$ and $b,b'\in \cB$, let us simply denote by 
$\varphi_{bb'}\in \C$ the $(b,b')$-entry of the matrix of $\varphi$ 
with respect to~$\cB$. Similarly, for any $\psi\in \End(\bar{V})$ 
and $b,b'\in \cB$, we denote by $\psi_{bb'}\in K$ the $(\bar{b},
\bar{b}')$-entry of the matrix of $\psi$ with respect to~$\bar{\cB}$. 

\begin{prop} \label{steinbergp} Let $r\geq 1$ and $\beta_1,\ldots, 
\beta_r\in \Phi^+$ (not necessarily distinct). Let $\alpha\in \Phi^+$ be 
such that $\hgt(\alpha) \leq \hgt(\beta_i)$ for $i=1,\ldots,r$. Let 
$b,b'\in \cB$ be such that $(\bar{\bE}_{\alpha}^{[1]})_{bb'} \neq 0$. 
(Such $b,b'$ exist by Proposition~\ref{nontrivrhoi}, and we have 
$b\neq b'$.) Then, for $\zeta_1,\ldots,\zeta_r\in K$, we have 
\[ \bigl(\bar{x}_{\beta_1}(\zeta_1) \cdots \bar{x}_{\beta_r}(\zeta_r)
\bigr)_{bb'}=\bigl(\sum_{\atop{1\leq i \leq r}{\beta_i=\alpha}}\zeta_i
\bigr)(\bar{\bE}_\alpha^{[1]})_{bb'}.\]
\end{prop}


\begin{proof} First we work over $\C$. Let $d\geq 1$ be such that 
$\rho(\be_{\beta_i}^+)^d=\underline{0}$ for all $i$. For 
any integers $l_1,\ldots,l_r\geq 0$ we set 
\[ \sigma_{l_1,\ldots,l_r}:=\frac{1}{l_1!\cdots l_r!}\rho
(\be_{\beta_1}^+)^{l_1} \circ \ldots\circ\rho(\be_{\beta_r}^+)^{l_r}
\in \End(V).\]
Let $t_1,\ldots,t_r\in \C$. Then, by the definition of $\exp$ and 
the elements $x_{\beta_i}(t_i)\in G_\C(V,\cB)$, we have 
\[ x_{\beta_1}(t_1) \cdots x_{\beta_r}(t_r)= \sum_{0\leq l_1,
\ldots,l_r\leq d} t_1^{l_1}\cdots t_r^{l_r} \sigma_{l_1,\ldots,l_r}.\]
Let $\mu \in P_\fh(V)$ be such that $b\in V_\mu$. Then a repeated 
application of the ``\nm{fundamental calculation}'' in 
Remark~\ref{fundcal} shows that 
\[\sigma_{l_1,\ldots,l_r}(b)\in V_{\mu+l_1\beta_1+ \ldots + l_r
\beta_r} \qquad \mbox{for any $l_1,\ldots,l_r\geq 0$}.\]
Here, we have $\sigma_{0,\ldots,0}(b)=b$. Now let $l_1,\ldots,l_r
\geq 0$ be such that $l_1\beta_1+\ldots + l_r\beta_r=\alpha$. Since 
$\hgt(\alpha)\leq \hgt(\beta_i)$ for all~$i$, this implies that 
$\beta_i=\alpha$ and $l_i=1$ for exactly one~$i$, and $l_j=0$ for 
$i\neq j$. Note that, in this case, $\sigma_{l_1,\ldots,l_r}=
\rho(\be_{\beta_i}^+)=\rho(\be_\alpha^+)$. Hence, we obtain 
\[ \bigl(x_{\beta_1}(t_1) \cdots x_{\beta_r}(t_r)\bigr)(b)
\in b+\bigl(\sum_{\atop{1\leq i \leq r}{\beta_i=\alpha}} t_i\bigr)
\rho(\be_{\alpha}^+)(b)+\sum_{\mu'} V_{\mu'},\]
where the sum runs over all $\mu'\in P_\fh(V)$ which can be written as 
$\mu'=\mu+l_1\beta_1+\ldots +l_r\beta_r$ for integers $l_1,\ldots,l_r$ 
such that $l_1\beta_1+\ldots +l_r\beta_r$ is non-zero and not equal to
$\alpha$. Now note that $\rho(\be_{\alpha}^+)(b)\in V_{\mu+\alpha}$. 
Hence, since $(\bar{\bE}_{\alpha}^{[1]})_{bb'}\neq 0$, we also
have $(\bE_\alpha^{[1]})_{bb'}\neq 0$ and so $b'\in V_{\mu+\alpha}$; 
in particular, $b\neq b'$ since $V$ is the direct sum of its weight 
spaces. Passing to matrices, this also implies that 
\[ \bigl(X_{\beta_1}(t_1) \cdots X_{\beta_r}(t_r)\bigr)_{bb'}=
\bigl(\sum_{\atop{1\leq i \leq r}{\beta_i=\alpha}}t_i\bigr)
(\bE_{\alpha}^{[1]})_{bb'}.\]
Now we can lift this to a ``polynomial level'' where we work over the 
ring $\Z[T_1,\ldots,T_r]$ in commuting indeterminates $T_1,\ldots,T_r$. 
Since the above identity holds for all $t_1,\ldots,t_r\in \C$, we have 
\[ \bigl(X_{\beta_1}(T_1) \cdots X_{\beta_r}(T_r)\bigr)_{bb'}=
\bigl(\sum_{\atop{1\leq i \leq r}{\beta_i=\alpha}}T_i\bigr)
(\bE_{\alpha}^{[1]})_{bb'}.\]
Given elements $\zeta_1,\ldots,\zeta_r\in K$ we have a canonical 
ring homomorphism $\Z[T_1,\ldots,T_r]\rightarrow K$ such that $T_i\mapsto 
\zeta_i$ for all~$i$. Applying that homomorphism to the above identity, 
we obtain an analogous identity over $K$, as required.
\end{proof}

\begin{cor} \label{cornontrivrhoi} Let $\alpha\in \Phi$ and $0\neq 
\zeta\in K$. Then $\bar{x}_\alpha(\zeta)\neq \id_{\bar{V}}$.
\end{cor}

\begin{proof} Let $\zeta\in K$. As in the proof of Corollary~\ref{corCree0a}, 
there is some $i\in I$ and $\eta\in \bar{G}$ such that $\bar{x}_\alpha
(\pm \zeta)=\eta \bar{x}_{\alpha_i}(\zeta)\eta^{-1}$. So it is sufficient
to consider the case where $\alpha=\alpha_i$ for some $i\in I$. But in
this case, we can use Proposition~\ref{steinbergp} with $r=1$ and $\alpha
=\beta_1=\alpha_i$. It follows that $\bar{x}_{\alpha_i}(\zeta)_{bb'}=
\zeta(\bar{\bE}_{\alpha_i}^{[1]})_{bb'}$ where $b\neq b'$ and 
$(\bar{\bE}_{\alpha_i}^{[1]})_{bb'}\neq 0$. Hence, if $\zeta\neq 0$,
then $\bar{x}_{\alpha_i}(\zeta)_{bb'}\neq 0$ and so $\bar{x}_{\alpha_i}
(\zeta)\neq \id_{\bar{V}}$, as required.
\end{proof}

\begin{xca} \label{steinbergpxca}  Show that, in the setting and with the 
assumptions of Proposition~\ref{steinbergp}, we also have
\[ \bigl(\zeta_1\bar{\bE}_{\beta_1}^{[1]}+\ldots +\zeta_r
\bar{\bE}_{\beta_r}^{[1]}\bigr)_{bb'}=\bigl(\sum_{\atop{1\leq i 
\leq r}{\beta_i=\alpha}}\zeta_i \bigr)(\bar{\bE}_\alpha^{[1]})_{bb'}.\]
Conclude that, if $|\Phi^+|=N$ and $\Phi^+=\{\beta_1,\ldots,\beta_N\}$, 
then the induced maps $\overline{\rho(\be_{\beta_1}^+)}$, $\ldots$, 
$\overline{\rho(\be_{\beta_N}^+)}$ are linearly independent in 
$\End(\bar{V})$. 
\end{xca}

In Exercise~\ref{xcafieldauto} we introduced \nmi{field automorphisms}{field
automorphism} of~$\bar{G}$. We now construct further automorphisms of 
$\bar{G}$. Assume we are given a collection of elements $\{\gamma_i \mid i 
\in I\} \subseteq K^\times$. We extend this to a collection $\{\gamma_\alpha
\mid \alpha\in \Phi\}\subseteq K^\times$ as follows. For $\alpha\in \Phi$
we write $\alpha=\sum_{i \in I} n_i\alpha_i$ with $n_i\in \Z$, and
set $\gamma_\alpha:=\prod_{i \in I} \gamma_i^{n_i}$. Note that, for 
all $\alpha,\beta\in \Phi$ such that $\alpha+\beta\in \Phi$, we have
\[ \gamma_{-\alpha}=\gamma_\alpha^{-1}\qquad\mbox{and}\qquad
\gamma_{\alpha+\beta}=\gamma_\alpha\gamma_\beta.\]

\begin{prop}[\nmi{Diagonal automorphisms}{diagonal automorphism}] 
\label{diagauto} In the above setting, there exists a unique group 
automorphism $\varphi\colon \bar{G}\rightarrow \bar{G}$ such that 
\[ \varphi(\bar{x}_\alpha(\zeta))=\bar{x}_\alpha(\gamma_\alpha\zeta) 
\qquad\mbox{for all $\alpha\in \Phi$ and $\zeta\in K$}.\]
If $K$ is algebraically closed, then $\varphi$ is an inner automorphism.
\end{prop}

\begin{proof} Since the matrix $A=(a_{ij})_{i,j\in I}$ has a non-zero
determinant, there exist numbers $c_{ij}\in \Q$ such that, for any
$i,j\in I$, we have 
\begin{equation*}
\renewcommand{\arraystretch}{0.8}
\sum_{l \in I} c_{il}a_{lj}=\left\{\begin{array}{cl} 1 & \quad
\mbox{if $i=j$},\\ 0 &\quad\mbox{if $i\neq j$}.\end{array}\right.\tag{$*$}
\end{equation*}
Let $n \in \Z_{>0}$ be such that $nc_{ij}\in \Z$ for all $i\in I$; in
fact, we may take $n=\det(A)$. Now let $K_1\supseteq K$ be a field 
extension such that there exist elements $\tilde{\gamma}_i\in K_1^\times$ 
with $\tilde{\gamma}_i^n=\gamma_i$ for all $i\in I$. By 
Remark~\ref{remfieldext}, we may naturally identify $\bar{G}=G_K(V,\cB)$ 
with a subgroup of $\bar{G}_1:=G_{K_1}(V,\cB)$. Under this identification, 
the element $\bar{x}_\alpha(\zeta)$ (for $\alpha\in \Phi$ and $\zeta\in K$) 
is exactly the same in $\bar{G}$ and in $\bar{G}_1$. Now set 
\[ \eta:=\prod_{i,j \in I} \bar{h}_j\bigl(\tilde{\gamma}_i^{nc_{ij}}
\bigr)\in\bar{G}_1 \qquad \mbox{(note that $nc_{ij}\in \Z$)}\]
and let $\varphi_1\colon \bar{G}_1\rightarrow \bar{G}_1$ be the inner 
automorphism given by conjugation with~$\eta$; thus, we have $\varphi_1(g)=
\eta g\eta^{-1}$ for all $g \in \bar{G}_1$. Let $\alpha\in\Phi$ and 
$\zeta\in K$. Write $\alpha=\sum_{l\in I} n_l \alpha_l$ with $n_l\in \Z$. 
Then, using a repeated application of Proposition~\ref{Cree0a}(b), we obtain:
\begin{align*} 
\varphi_1(\bar{x}_\alpha(\zeta))&=\bar{x}_\alpha\Bigl(\zeta \prod_{i,j 
\in I} \tilde{\gamma}_i^{nc_{ij}\langle \alpha_j^\vee,\alpha \rangle}
\Bigr)=\bar{x}_\alpha\Bigl(\zeta \prod_{i,j,l \in I} \tilde{\gamma}_i^{n_l
nc_{ij}\langle \alpha_j^\vee,\alpha_l \rangle}\Bigr)\\
&= \bar{x}_\alpha \Bigl(\zeta \prod_{i,j,l\in I} \tilde{\gamma}_i^{n_ln
c_{ij}a_{jl}} \Bigr)=\bar{x}_\alpha \Bigl(\zeta \prod_{i,l \in I} 
\tilde{\gamma}_i^{n_ln\sum_{j \in I} c_{ij}a_{jl}}\Bigr)\\ & 
\stackrel{(*)}{=}\bar{x}_\alpha\Bigl(\zeta \prod_{i \in I}
\tilde{\gamma}_i^{n_in} \Bigr)=\bar{x}_\alpha\Bigl(\zeta\prod_{i \in I} 
\gamma_i^{n_i}\Bigr) =\bar{x}_\alpha(\gamma_\alpha\zeta).
\end{align*}
Since $\bar{G}=\langle \bar{x}_\alpha(\zeta)\mid \alpha\in \Phi,
\zeta\in K\rangle$, it follows that $\varphi_1(\bar{G})=\bar{G}$. Hence, 
$\varphi_1$ restricts to an automorphism of $\bar{G}$ with the 
desired property.
\end{proof}

In the following remark and exercises, we return to the setting of 
Chapter~\ref{chap3} and consider the group $\bar{G}=G_K(\fg,\bB)\subseteq 
\GL(\bar{\fg})$ of adjoint type, as in Example~\ref{defnonad1}(b). Recall
that 
\[\bB=\{h_j^+ \mid j \in I\}\cup \{\be_\alpha^+\mid \alpha\in \Phi\}.\]
The new feature here (as compared to Chapter~\ref{chap3}) is that we 
also have the elements $\bar{x}_\alpha(\zeta):=\bar{x}_\alpha(\zeta;\fg,
\bB)\in \bar{G}$ for {\itshape all} $\alpha\in\Phi$ and $\zeta \in K$.

\begin{rem} \label{remfoot} Assume that $V=\fg$ is the adjoint module,
as above. We claim that, in this case, the element $\eta \in \bar{G}_1=
G_{K_1}(\fg,\bB)$ constructed in the above proof is explicitly given by 
\[ \eta(\bar{h}_j^+)=\bar{h}_j^+ \quad (j\in I) \qquad \mbox{and}\qquad
\eta(\bar{\be}_\alpha^+)=\gamma_\alpha \bar{\be}_\alpha^+ \quad (\alpha
\in \Phi).\]
Thus, $\eta$ is represented by a diagonal matrix with entries in $K$
(and not just in $K_1$); furthermore, $\eta$ only depends on $\{\gamma_i
\mid i \in I\}$ but not on the choice of the elements $\tilde{\gamma}_i
\in K_1^\times$ in the above proof. This is sees as follows. By 
Proposition~\ref{luform3c}, the elements $\bar{h}_i(\xi)\in
\bar{G}_1$ are represented by diagonal matrices; hence, so is the 
element~$\eta$. The explicit formulae in that proposition show 
that $\eta(\bar{h}_j^+)=\bar{h}_j^+$ for all $j\in I$; furthermore,
the scalar by which $\eta$ acts on $\bar{\be}_\alpha^+$ is given by 
the product (over all $i,j\in I$) of the elements $\tilde{\gamma}_i^{n
c_{ij}\langle \alpha_j^\vee,\alpha \rangle}\in K_1$, which evaluates to 
$\gamma_\alpha\in K$ by the same computation as in the above proof. 

The elements $\eta$ indeed are the ``additional'' diagonal elements 
mentioned in the footnote to Definition~\ref{deflu}. We will see in 
the following section that, when $V$ is not necessarily the adjoint 
module, then $\eta$ is still represented by a diagonal matrix, but 
possibly with diagonal entries in the larger field $K_1\supseteq K$.
\end{rem}

\begin{xca} \label{xcanit1} Let $\bar{\omega}\colon \bar{\cL}\rightarrow
\bar{\cL}$ be as in Example~\ref{xcaomega1}. Show that 
\[ \bar{\omega}\circ \bar{x}_\alpha (\zeta) \circ \bar{\omega}^{-1}=
\bar{x}_{-\alpha}(-\zeta)\qquad \mbox{for all $\alpha\in \Phi$ and 
$\zeta \in K$}.\]
Hence, conjugation by $\bar{\omega}$ defines an automorphism of 
$\bar{G}=G_K(\fg,\bB)$.\\
\noindent {\footnotesize [{\it Hint.} Apply the usual argument: first 
over $\C$, then pass to $K$. See Example~\ref{xcaomega1} for the special
case $\alpha=\alpha_i$ ($i \in I$), but note that $\omega(\be_\alpha^+)=-
\be_{-\alpha}^+$ for arbitrary~$\alpha$.]}
\end{xca}

\begin{xca} \label{xcagraphauto} Let $A=(a_{ij})_{i,j\in I}$ be 
indecomposable of simply laced type; let $i\mapsto i^\prime$ be a 
permutation of $I$ as in Exercise~\ref{xcafold}. Let $\tau\colon \Phi
\rightarrow \Phi$ be the induced permutation of $\Phi$. Let $\bar{G}=
G_K(\fg,\bB)$. Then show that there is a unique automorphism $\bar{\tau} 
\colon \bar{G} \rightarrow\bar{G}$ such that 
\[ \bar{\tau}\bigl(\bar{x}_\alpha(\zeta)\bigr)=\bar{x}_{\tau(\alpha)}
(\zeta)\qquad \mbox{for all $\alpha\in \Phi$ and $\zeta\in K$}.\]
The map $\bar{\tau}\colon \bar{G} \rightarrow \bar{G}$ is called
a \nm{graph automorphism} of $\bar{G}$.

Note that usually there are certain signs involved in the description
of such an automorphism; see Carter \cite[Prop.~12.2.3]{Ca1} or 
Steinberg \cite[Theorem~29 (p.~91)]{St} (and its corollary). Here, these 
signs disappear since we are working with Lusztig's canonical basis.\\
\noindent {\footnotesize [{\it Hints}. Argue as in Example~\ref{xcaomega1} 
and Exercise~\ref{xcanit1}.]}
\end{xca}

There are more general types of graph automorphisms, also for groups of
non-adjoint type; for this we refer to Carter \cite[Chap.~12]{Ca1} and 
Steinberg \cite[Chap.~10]{St}. If $K$ is algebraically closed, then
there is also the important \nmi{Isogeny Theorem}{isogeny theorem} which
describes all possible homomorphisms $\bar{G}\rightarrow \bar{G}$ with a 
finite kernel; see Steinberg \cite{St99} and the further discussion (with 
many examples) in \cite[\S 1.4]{GeMa}.

\section{The diagonal and monomial subgroups} \label{secdiagII}

We keep the notation from the previous two sections and define the 
following subgroups of $\bar{G}=G_K(V,\cB)$:
\begin{align*}
\bar{N}&:=\bigl\langle \bar{n}_i(\xi;V,\cB)\mid i \in I,\xi\in K^\times
\bigr\rangle \qquad \mbox{``\nm{monomial subgroup}''},\\
\bar{H}&:=\bigl\langle \bar{h}_i(\xi;V,\cB)\mid i \in I,\xi\in K^\times
\bigr\rangle \qquad \mbox{``\nm{diagonal subgroup}''}.
\end{align*}
By the definition of the elements $\bar{h}_i(\xi;V,\cB)$, it is clear
that $\bar{H}\subseteq \bar{N}$. One of our aims will be to show that 
$\bar{H}$ is normal in $\bar{N}$ and that the factor group 
$\bar{N}/\bar{H}$ is isomorphic to the Weyl group~$W$ of $\fg$. 

If $K$ is algebraically closed, we will see that $\bar{H}$ is isomorphic 
to a direct product of $|I|$ copies of the multplicative group $K^\times$, 
that is, $\bar{H}$ is a ``\nm{torus}'' in the sense of the theory of 
algebraic groups. (This will be further developed in the following chapter.)

First we need some information about the action of $n_i(t;V,\cB)$
and $h_i(t;V,\cB)$ on~$V$, where $t\in \C^\times$. Ideally, we would like
to have analogues of Theorem~\ref{luform3} and Proposition~\ref{luform3c} 
(that were proved for groups of adjoint type). However, explicit formulae 
analogous to those in Theorem~\ref{luform3} are simply not available for
general $V,\cB$. Still, we will get around that difficulty and obtain 
formulae for $h_i(t;V, \cB)$, which are a direct generalisation of those 
in Proposition~\ref{luform3c}.

In the discussion below, $V$ and $\cB$ will be fixed throughout. So we 
shall usually omit the symbols $V,\cB$ from the notation. Thus,
\begin{alignat*}{2}
n_i(t)&=n_i(t;V,\cB), \quad &h_i(t)=h_i(t;V,\cB) &\qquad\mbox{for $t\in 
\C^\times$},\\ \bar{n}_i(\xi)&=\bar{n}_i(\xi;V,\cB), \quad \;\;
&\bar{h}_i(\xi) = \bar{h}_i(\xi;V,\cB) &\qquad\mbox{for $\xi\in K^\times$}.
\end{alignat*} 
The following result will serve as a weak substitute of Theorem~\ref{luform3}.

\begin{lem} \label{Cree0aa} Let $i\in I$ and $t\in \C^\times$. Then 
$n_i(t)(V_\mu)=V_{s_i(\mu)}$ for all $\mu \in P_\fh(V)$. In particular, 
$\dim V_\mu=V_{s_i(\mu)}$.
\end{lem}

\begin{proof} Let $\alpha\in \Phi$ and consider the element $h_\alpha
\in\fh$. By Proposition~\ref{Cree0}, we have 
\[ n_i(t) \circ \rho(h_\alpha) \circ n_i(t)^{-1}=\rho(h_{s_i(\alpha)})
\qquad \mbox{for $t \in \C^\times$}.\] 
Since $n_i(t)^{-1}=n_i(-t)$, we also have $n_i(t)^{-1} \circ 
\rho(h_\alpha)\circ n_i(t)=\rho(h_{s_i(\alpha)})$ for $t \in \C^\times$. 
Now let $v\in V_\mu$. Then $h_{s_i(\alpha)}.v=\mu(h_{s_i(\alpha)})v$ and 
\begin{align*}
h_\alpha.(n_i(t)(v))&=(\rho(h_\alpha)\circ n_i(t))(v)=(n_i(t)\circ 
\rho(h_{s_i(\alpha)}))(v)\\&=n_i(t)(h_{s_i(\alpha)}.v)=
\mu(h_{s_i(\alpha)})n_i(t)(v).
\end{align*} 
Now, by Proposition~\ref{weyl5c}, we have $h_{s_i(\alpha)}=h_\alpha-
\alpha_i(h_\alpha)h_i$. Hence,
\[ \mu(h_{s_i(\alpha)})=\mu(h_\alpha)-\alpha_i(h_\alpha)\mu(h_i).\]
On the other hand, $s_i(\mu)=\mu-\mu(h_i)\alpha_i$. Hence, $s_i(\mu)
(h_\alpha)=\mu(h_{s_i(\alpha)})$ and so $n_i(t)(v)\in V_{s_i(\mu)}$.
Thus, we have $n_i(t)(V_\mu)\subseteq V_{s_i(\mu)}$. 
Finally, since $n_i(t)\colon V \rightarrow V$ is bijective, we have 
$\dim V_\mu \leq \dim V_{s_i(\mu)}$. But, since $s_i^2=\id_{\fh^*}$, 
we also have $\dim V_{s_i(\mu)} \leq \dim V_{s_i(s_i(\mu))}=V_\mu$. 
\end{proof}

%

\begin{prop} \label{Cree0ab} Let $i\in I$ and $t\in \C^\times$. Then
\[ h_i(t)(v)=t^{\langle \alpha_i^\vee,\mu\rangle}v\qquad\mbox{for
any $\mu \in P_\fh(V)$ and $v\in V_\mu$}.\]
Thus, the matrix of $h_i(t)$ with respect to $\cB$ is diagonal. 
Furthermore, we have $n_i(t)^2=h_i(-1)$.
\end{prop}

\begin{proof} Let $d\geq 1$ be such that $\rho(e_i)^d=\rho(f_i)^d=
\underline{0}$. Using the notation introduced in Remark~\ref{remadm2}, 
we set 
\[\sigma_{m,i}:=\sum_{\atop{0\leq l_1,l_2,l_3\leq d}{l_1-l_2+l_3=m}}
(-1)^{l_2}\,e_i^{[l_1]}\circ f_i^{[l_2]} \circ e_i^{[l_3]}\in \End(V),\]
for any $m\in\Z$. Now let $t\in \C^\times$. Then, by the definition of 
$x_i(t)$ and $y_i(-t^{-1})$, we have 
\[n_i(t)=x_i(t)y_i(-t^{-1})x_i(t)=\sum_{-d\leq m \leq 2d} t^m 
\sigma_{m,i}.\]
Now let us fix $v \in V_\mu$. Using the ``\nm{fundamental calculation}'' 
in Remark~\ref{fundcal}, we find that 
\begin{align*}
e_i^{[l_3]}(v)&\in V_{\mu+l_3\alpha_i},\\
\bigl(f_i^{[l_2]}\circ e_i^{[l_3]} \bigr)(v)&\in 
V_{\mu+l_3\alpha_i-l_2\alpha_i},\\
\bigl(e_i^{[l_1]}\circ f_i^{[l_2]}\circ e_i^{[l_3]}
\bigr)(v)&\in V_{\mu+l_3\alpha_i-l_2\alpha_i+l_1\alpha_i}
\end{align*}
for any $l_1,l_2,l_3\geq 0$. (Note that $e_i \in \fg_{\alpha_i}$ and 
$f_i\in \fg_{-\alpha_i}$.) It follows that $v_m:=\sigma_{m,i}(v)\in 
V_{\mu+m\alpha_i}$ for $-d\leq m \leq 2d$. So we conclude that 
\[ n_i(t)(v)=\sum_{-d\leq m \leq 2d} t^mv_m \qquad \mbox{for all
$t\in \C^\times$},\]
where the vectors $v_m$ do not depend on~$t$. On the other hand, by 
Lemma~\ref{Cree0aa}, we know that $n_i(t)(v) \in V_{s_i(\mu)}$. Now 
$s_i(\mu)=\mu-\mu(h_i)\alpha_i$ and $\mu(h_i)=\langle \alpha_i^\vee,\mu
\rangle$. Hence, we must have
\[ n_i(t)(v)=t^mv_m \qquad \mbox{for all $t\in \C^\times$, where
$m=-\langle \alpha_i^\vee, \mu\rangle$}.\]
Since this holds for all $t$, we also have $n_i(-t)(v)=
(-t)^m v_m$ and $n_i(-1)(v)=(-1)^mv_m$. Since $n_i(-t)=
n_i(t)^{-1}$, we deduce that $v=(-t)^m n_i(t)(v_m)$. It follows that 
\begin{align*}
h_i(t)(v)&=\bigl(n_i(t)\circ n_i(-1)\bigr)(v)=(-1)^m 
n_i(t)(v_m)\\&=(-1)^m(-t)^{-m}v=t^{-m}v=t^{\langle \alpha_i^\vee,
\mu\rangle}v,
\end{align*}
as desired. Furthermore, we have 
\[ n_i(t)^2(v)=t^mn_i(t)(v_m)=t^m(-t)^{-m}v=(-1)^mv.\]
Since we also have $h_i(-1)(v)=(-1)^{-m}=(-1)^m$, we conclude that 
$n_i(t)^2(v)=h_i(-1)(v)$. Since this holds for all $\mu\in P_\fh(V)$
and all $v\in V_\mu$, we conclude that $n_i(t)^2=h_i(-1)$.
\end{proof}

Next, we would like to transfer the above identities from $\C$ to~$K$.
So, as usual, we need to work at a ``polynomial level''. We consider 
the ring of Laurent polynomials $\Z[T,T^{-1}]$ in an indeterminate~$T$.
Already in the proof of Proposition~\ref{Cree0a}, we introduced
the matrices $N_i(T)$ and $H_i(T)$ with entries in $\Z[T,T^{-1}]$.
Upon substituting $T\mapsto t$ for $t\in \C^\times$, we obtain the
matrices 
\[ N_i(t):=M_\cB(n_i(t))\qquad\mbox{and}\qquad H_i(t):=M_\cB(h_i(t)).\]
On the other hand, if $\xi \in K^\times$, then we have a canonical ring 
homomorphism $\Z[T,T^{-1}]\rightarrow K$ such that $T\mapsto \xi$.
Applying that homomorphism to $N_i(T)$ and $H_i(T)$, we obtain the matrices 
\[ \bar{N}_i(\xi):=M_{\bar{\cB}}(\bar{n}_i(\xi))\qquad \mbox{and}
\qquad \bar{H}_i(\xi):=M_{\bar{\cB}}((h_i(\xi))).\]
With this notation, we can now state:

\begin{cor} \label{polylev2} Let $i\in I$ and $\xi\in K^\times$. Then 
the matrix $\bar{H}_i(\xi)$ is a diagonal matrix. If $\mu\in P_\fh(V)$
and $b\in \cB\cap V_\mu$, then the $(\bar{b},\bar{b})$-diagonal entry of 
$\bar{H}_i(\xi)$ is given by $\xi^{\langle \alpha_i^\vee,\mu\rangle}$.
\end{cor}

\begin{proof} First let $t\in \C^\times$ and consider the matrix
$H_i(t)$ of $h_i(t)$. For $b,b'\in \cB$ let $H_i(t)_{bb'}$ be the
$(b,b')$-entry of $H_i(t)$. By Proposition~\ref{Cree0ab}, we have 
\[ H_i(t)_{bb'}=\left\{\begin{array}{cl} t^{\langle \alpha_i^\vee,\mu
\rangle} & \qquad \mbox{if $b=b'\in V_\mu$},\\ 0 & \qquad
\mbox{if $b\neq b'$}. \end{array}\right.\]
Since this holds for all $t\in \C^\times$, we also have identities 
at the ``polynomial level'' (over $\Z[T,T^{-1}]$ as above):
\[ H_i(T)_{bb'}=\left\{\begin{array}{cl} T^{\langle \alpha_i^\vee,\mu
\rangle}&\qquad\mbox{if $b=b'\in V_\mu$},\\0&\qquad\mbox{if $b\neq b'$}.
\end{array}\right.\]
Now let $\xi\in K^\times$ and apply the ring homomorphism $\Z[T]
\rightarrow K$ such that $T\mapsto \xi$. This yields that 
$\bar{H}_i(\xi)$ also is diagonal, with diagonal entries as stated.
\end{proof}


\begin{prop} \label{Cree0ac} Let $i,j\in I$ and $\zeta,\xi\in K^\times$. 
Then we have:
\begin{itemize}
\item[{\rm (a)}] $\bar{h}_i(\zeta)\bar{h}_j(\xi)=\bar{h}_j(\xi)
\bar{h}_i(\zeta)$ and $\bar{h}_i(\zeta\xi)=\bar{h}_i(\zeta)\bar{h}_i(\xi)$.
\item[{\rm (b)}] $\bar{n}_j(\xi)\bar{h}_i(\zeta)\bar{n}_j(\xi)^{-1}=
\bar{h}_i (\zeta)\bar{h}_j(\zeta^{-a_{ij}})$.
\item[{\rm (c)}] $\bar{n}_i(\xi)^2=\bar{h}_i(-1)$ and $\bar{n}_i(\xi)^4=
\bar{h}_i(1)=\id_{\bar{V}}$.
\end{itemize}
\end{prop}

\begin{proof} First we work with the elements $h_i(t)$, $n_i(u)$ and 
their matrices $H_i(t)$, $N_j(u)$ over $\C$, where $t,u\in 
\C^\times$. 

(a) Let $n=\dim V$ and write $\cB=\{v_1,\ldots,v_n\}$. For each $l\in 
\{1,\ldots,n\}$ let $\mu_l\in P_\fh(V)$ be the weight of $v_l$ and set 
$m_l:=\langle \alpha_i^\vee,\mu_l\rangle$ for $l=1,\ldots,n$. Then 
Proposition~\ref{Cree0ab} shows that the matrix $H_i(t)$ is diagonal with 
diagonal entries $t^{m_l}$ for $l=1,\ldots,n$. So, clearly, $H_i(t)$
and $H_j(u)$ commute with each other. Furthermore, $H_i(t)H_i(u)$ is 
the diagonal matrix with diagonal entries $t^{m_l}u^{m_l}=(tu)^{m_l}$. 
Hence, $H_i(t)H_i(u)=H_i(tu)$.

(b) Let $\mu \in P_\fh(V)$ and $v\in V_\mu$. We set $v':=n_j(u)^{-1}(v)$.
Since $n_j(u)^{-1}=n_j(-u)$, we have $v'\in V_{s_j(\mu)}$
by Lemma~\ref{Cree0aa}. Using also Proposition~\ref{Cree0ab}, we obtain
\begin{align*}
\bigl(n_j(u) h_i(t) &n_j(u)^{-1}\bigr)(v)=\bigl(n_j(u) 
h_i(t)\bigr)(v')\\&=t^{\langle \alpha_i^\vee,s_i(\mu)\rangle}
n_j(u)(v')=t^{\langle \alpha_i^\vee,s_j(\mu)\rangle}v.
\end{align*}
Now $s_j(\mu)=\mu-\mu(h_j)\alpha_j$ and $\mu(h_j)=\langle \alpha_j^\vee,
\mu\rangle$. It follows that $\langle \alpha_i^\vee,s_j(\mu)\rangle=
\langle \alpha_i^\vee,\mu\rangle -\langle \alpha_j^\vee,\mu\rangle
\langle \alpha_i^\vee,\alpha_j\rangle$, where $\langle \alpha_i^\vee,
\alpha_j\rangle=a_{ij}$. Hence,
\[ t^{\langle \alpha_i^\vee,s_j(\mu)\rangle}v=t^{\langle \alpha_i^\vee,
\mu\rangle}t^{-a_{ij}\langle \alpha_j^\vee,\mu\rangle}v=
\bigl(h_i(t)h_j(t^{-a_{ij}})\bigr)(v).\]
So the two linear maps $n_j(u) h_i(t) n_j(u)^{-1}$ and $h_i(t)h_j
(t^{-a_{ij}})$ have the same effect on~$v$. Since this holds for 
all~$v$ in a basis of $V$, the two maps must be equal. 

(c) By Proposition~\ref{Cree0ab}, we already know that $n_i(u)^2=h_i(-1)$.
Using (a), we obtain that $n_i(u)^4=h_i(-1)^2=h_i((-1)(-1))=h_i(1)$. 
Finally, Proposition~\ref{Cree0ab} also shows that $h_i(1)=\id_V$.

Now we can pass to $K$. Once the above relations are established for
all $t,u\in \C^\times$, we obtain identities at the ``polynomial level'':
\begin{gather*}
H_i(T){\cdot}H_j(U)=H_j(U){\cdot}H_i(T) \;\mbox{ and }\; H_i(TU)=
H_i(T){\cdot}H_i(U),\\ N_i(U){\cdot}H_i(T){\cdot}N_i(U)^{-1} = 
H_i(T){\cdot}H_j(T^{-a_{ij}}),\\ N_i(U)^2=H_i(-1) \; \mbox{ and } 
\; N_i(U)^4=H_i(1)=\mbox{identity matrix},
\end{gather*}
where we work over the ring of Laurent polynomials $\Z[T^{\pm 1}, U^{\pm 1}]$
in two commuting indeterminates $T,U$. Given $\zeta,\xi\in K^\times$, it
then remains to apply the canonical ring homomorphism $\Z[T^{\pm 1},
U^{\pm 1}]\rightarrow K$ which sends $T$ to~$\zeta$ and $U$ to~$\xi$.
\end{proof}

\begin{rem} \label{sl2triple} Let us fix $i\in I$. In analogy to 
Lemma~\ref{lemkchar4}, we define the subgroup $\bar{G}_i:=\bigl\langle 
\bar{x}_i(\zeta),\bar{y}_i(\zeta)\mid \zeta \in K\bigr\rangle\subseteq 
\bar{G}$. By copying, almost verbatim, the proof of 
Proposition~\ref{sl2tripGi}, one finds that
\[ \bar{G}_i=\big\{\bar{x}_i(\zeta_1)\bar{y}_i(\zeta_2)\bar{h}_i(\xi)
\bar{x}_i(\zeta_3)\mid \zeta_1,\zeta_2,\zeta_3\in K,\xi \in K^\times
\big\}.\]
Indeed, that proof essentially relied on certain relations among the elements 
$\bar{x}_i(\zeta)$, $\bar{y}_i(\zeta)$, $\bar{n}_i(\xi)$, $\bar{h}_i(\xi)$
(for a Chevalley group of adjoint type). The exact analogues of those 
relations do hold in the present, more general setting by 
Proposition~\ref{action2a}(b), Remark~\ref{Cree0a0} and
Proposition~\ref{Cree0ac}. So one can literally just repeat the proof of 
Proposition~\ref{sl2tripGi}; we leave further details to the reader.
\end{rem}

\begin{cor} \label{diagII} The diagonal subgroup $\bar{H}\subseteq 
\bar{G}$ (as defined above) is abelian and we have 
\[ \bar{H}=\Bigl\{\prod_{i \in I}\bar{h}_i(\xi_i)\,\big|\, \xi_i
\in K^\times \mbox{ for all $i\in I$}\Bigr\}.\]
Furthermore, $\bar{H}$ is a normal subgroup of $\bar{N}$ and we have 
\[ \bar{N}=\bigl\{\bar{n}_{i_1}(1)\cdots \bar{n}_{i_r}(1)h\mid h\in
\bar{H}\;\mbox{ and }\;i_1,\ldots,i_r\in I\; (r\geq 0)\bigr\}.\]
\end{cor}

\begin{proof} The fact that $\bar{H}$ is abelian and is normal in
$\bar{N}$ immediately follows from Proposition~\ref{Cree0ac}. Now let
$g\in \bar{N}$. By Proposition~\ref{action2a}(b), we have $\bar{n}_i
(\xi)^{-1}=\bar{n}_i(-\xi)$ for $i \in I$ and $\xi\in K^\times$. 
Hence, we can write $g=\bar{n}_{i_1}(\xi_1)\cdots\bar{n}_{i_r}(\xi_r)$
where $i_1,\ldots,i_r\in I$ ($r\geq 0$) and $\xi_1,\ldots,\xi_r\in 
K^\times$. We have $\bar{n}_i(-1)=\bar{n}_i(1)^{-1}$ and $\bar{h}_i(\xi)
=\bar{n}_i(\xi)\bar{n}_i(-1)$ for $\xi\in K^\times$. Hence, $\bar{n}_i
(\xi)=\bar{h}_i(\xi)\bar{n}_i(1)$. Since $\bar{H}$ is normal in 
$\bar{N}$, it follows that $g\in\bar{n}_{i_1}(1) \cdots \bar{n}_{i_r}(1)
\bar{H}$, as claimed. 
\end{proof}

Let us see to what extent the above description of the elements of
$\bar{H}$ is unique.

\begin{lem} \label{uniqueH} Let $h=\prod_{i \in I} \bar{h}_i(\xi_i)
\in \bar{H}$, as above. Then 
\[ h=\id_{\bar{V}} \qquad \Leftrightarrow \qquad \prod_{i \in I} 
\xi_i^{\langle \alpha_i^\vee,\mu\rangle}=1 \quad \mbox{for all 
$\mu \in \Lambda(V)$},\]
where $\Lambda(V)=\langle P_\fh(V)\rangle_\Z$; see 
Definition~\ref{deflambda}.
\end{lem}

\begin{proof} Let $\mu \in P_\fh(V)$ and $b \in \cB\cap V_\mu$. By 
Corollary~\ref{polylev2}, the $(\bar{b},\bar{b})$-diagonal entry of each 
$\bar{H}_i(\xi_i)$ is given by $\xi_i^{\langle \alpha_i^\vee,\mu
\rangle}$. Hence, the $(\bar{b},\bar{b})$-diagonal entry of the matrix of
$h$ with respect to $\bar{\cB}$ is given by $\prod_{i \in I} \xi_i^{\langle 
\alpha_i^\vee,\mu\rangle}$. Consequently, we have the equivalence:
\begin{equation*}
h=\id_{\bar{V}} \qquad \Leftrightarrow \qquad \prod_{i \in I} 
\xi_i^{\langle \alpha_i^\vee,\mu\rangle}=1 \quad \mbox{for all $\mu 
\in P_\fh(V)$}.\tag{$*$}
\end{equation*}
Now let $\mu,\mu'\in \Lambda(V)$. If $\prod_{i \in I} \xi_i^{\langle 
\alpha_i^\vee,\mu\rangle}=1$ and $\prod_{i \in I} \xi_i^{\langle 
\alpha_i^\vee,\mu'\rangle}=1$, then we also have 
\[\prod_{i \in I} \xi_i^{\langle \alpha_i^\vee,-\mu\rangle}=1\qquad
\mbox{and}\qquad \prod_{i \in I} \xi_i^{\langle \alpha_i^\vee,\mu+\mu'
\rangle}=1.\]
Using this and ($*$), we obtain the implication 
\[ h=\id_{\bar{V}} \qquad \Rightarrow \qquad \prod_{i \in I} 
\xi_i^{\langle \alpha_i^\vee,\mu\rangle}=1 \quad \mbox{for all 
$\mu \in \Lambda(V)$}.\]
The reverse implication is clear by the above equivalence ($*$).
\end{proof}

\begin{exmp} \label{uniqueHa} Assume that $\Lambda(V)=\Omega$, the full
weight lattice; see Definition~\ref{weightlat}. Hence, in this case, 
$\Lambda(V)=\langle \varpi_j\mid j \in I\rangle_\Z$ where $\langle 
\alpha_i^\vee,\varpi_j \rangle=\delta_{ij}$ (Kronecker delta) for all 
$i,j\in I$. So, if $h=\prod_{i \in I} \bar{h}_i(\xi_i)$ as in 
Lemma~\ref{uniqueH}, then we have 
\[ \prod_{i \in I} \xi_i^{\langle \alpha_i^\vee,\varpi_j\rangle}=
\xi_j \qquad \mbox{for all $j\in I$}.\]
It follows that $h=\id_{\bar{V}}$ if and only if $\xi_j=1$ for all 
$j \in I$. Hence, in this case, the expression $h=\prod_{i \in I} 
\bar{h}_i(\xi_i)$ is unique. 
\end{exmp}

%

The above example shows that, if $\Lambda(V)=\Omega$, then the map 
\[ (K^\times)^I\rightarrow \bar{H}, \qquad (\xi_i)_{i \in I}
\mapsto \prod_{i \in I} \bar{h}_i(\xi_i),\]
is a group isomorphism. In general, if $\Lambda(V)\neq \Omega$, then the 
analogous map $(K^\times)^I \rightarrow \bar{H}$ need not be injective; 
see Example~\ref{expHNsl2}. In order to obtain an ``injectivity'' statement, 
we now introduce ``additional'' diagonal elements, analogous to those
in Remark~\ref{remfoot}. For this purpose, we use Remark~\ref{factorAa}. 
There we defined an additive subgroup $\fh_\Z \subseteq \fh$ (which depends 
on $\Lambda(V)\subseteq \Omega$) and showed that $\fh_\Z=\sum_{i \in I} 
\Z h_i'$ for certain elements $h_i'\in \fh$ such that $\mu(h_i')\in\Z$ 
for all $\mu \in P_\fh(V)$. For $i\in I$ and $\xi\in K^\times$ we now 
define $\bar{h}_i'(\xi)\in \GL(\bar{V})$ by 
\[ \bar{h}_i'(\xi)(\bar{b}):=\xi^{\mu(h_i')}\bar{b} \qquad 
\mbox{where $b\in \cB\cap V_\mu$ and $\mu\in P_\fh(V)$}.\]
In particular, $\bar{h}_i'(\xi)$ is represented by a diagonal matrix with
respect to the basis $\bar{\cB}$ of~$\bar{V}$. So we obtain an abelian
subgroup
\[ \bar{H}':=\langle \bar{h}_i'(\xi) \mid i \in I,\xi \in K^\times\rangle
\subseteq \GL(\bar{V}).\]
The following result shows that $\bar{H}'$ normalises $\bar{G}$. Thus,
we can form the (possibly) slightly larger group $\bar{G}':=\bar{G}{\cdot}
\bar{H}' \subseteq \GL(\bar{V})$. 

\begin{lem}[Ree] \label{xcasteinb35} Let $i\in I$ and $\xi \in K^\times$. 
Then 
\[\bar{h}_i'(\xi)\bar{x}_\alpha(\zeta)\bar{h}_i'(\xi)^{-1}=
\bar{x}_\alpha\bigl(\xi^{\alpha(h_i')}\zeta\bigr)\qquad \mbox{for 
$\alpha\in \Phi$ and $\zeta\in K$}.\]
So conjugation by $\bar{h}_i'(\xi)$ inside $\GL(\bar{V})$ defines
an automorphism of $\bar{G}$.
\end{lem}

\begin{proof} Let $\mu \in P_\fh(V)$ and $b \in \cB \cap V_\mu$. We must 
show that 
\[ \bigl(\bar{h}_i'(\xi)\circ \bar{x}_\alpha(\zeta)\bigr)(\bar{b})=
\bigl(\bar{x}_\alpha(\xi^{\alpha(h_i')}\zeta) \circ \bar{h}_i'(\xi)\bigr)
(\bar{b}).\]
This is seen as follows. We have $\bar{x}_\alpha(\zeta)(\bar{b})= 
\sum_{m\geq 0} \zeta^m \bar{\be}_\alpha^{[m]}(\bar{b})$. Now note that,
by arguments that we already used several times above, 
$\bar{\be}_\alpha^{[m]}(\bar{b})\in \bar{V}$ is a linear combination of basis 
vectors $\bar{b}'\in \bar{\cB}$ for various $b'\in \cB\cap V_{\mu+m\alpha}$.
Consequently, we obtain
\[ \bar{h}_i'(\xi)\bigl(\bar{\be}_\alpha^{[m]}(\bar{b})\bigr)=
\xi^{(\mu+m\alpha)(h_i')}\bar{\be}_\alpha^{[m]}(\bar{b}).\]
This yields that $\bigl(\bar{h}_i'(\xi)\circ \bar{x}_\alpha(\zeta)\bigr)
(\bar{b})=\sum_{m\geq 0} \xi^{(\mu+m\alpha)(h_i')}\zeta^m
\bar{\be}_\alpha^{[m]}(\bar{b})$, which is the same as 
$\bigl(\bar{x}_\alpha(\xi^{\alpha(h_i')}\zeta)\circ 
\bar{h}_i'(\xi)\bigr)(\bar{b})$.
\end{proof}

\begin{exmp} \label{sl2hp}  Let $\bar{G}=G_K(\fg,\bB)$ be of adjoint
type, as in Remark~\ref{remfoot}. Then $\Lambda(\fg)=\langle \Phi
\rangle_\Z$, with $\Z$-basis given by $\{\alpha_j\mid j \in I\}$. 
Hence, we have $\alpha_j(h_i')=\delta_{ij}$ for all $i,j\in I$ (see 
Remark~\ref{factorAa}). Now fix a collection $\{\gamma_i\mid i \in I\}
\subseteq K^\times$ and extend it to a collection $\{\gamma_\alpha\mid
\alpha\in \Phi\}$ as in Proposition~\ref{diagauto}. Setting $\eta:=
\prod_{i \in I} \bar{h}_i'(\gamma_i)\in \bar{H}'$, one readily checks
that $\eta(\bar{h}_j^+)=\bar{h}_j$ and $\eta(\bar{\be}_{\alpha_j}^+)=
\gamma_j\bar{\be}_{\alpha_j}^+$ for all $j\in I$. Consequently, we also
have $\eta(\bar{\be}_\alpha^+)=\gamma_\alpha \bar{\be}_\alpha$ for all
$\alpha\in \Phi$. Hence, for groups of adjoint type, $\eta=\prod_{i \in I}
\bar{h}_i'(\gamma_i)\in \bar{H}'$ is exactly the same as the ``diagonal'' 
element considered in Remark~\ref{remfoot}.
\end{exmp}


\begin{prop}[Steinberg] \label{steinb35} The map $\varphi \colon 
(K^\times)^I\rightarrow \bar{H}'$ defined by sending a tuple
$(\xi_i)_{i\in I}$ to $\prod_{i \in I} \bar{h}_i'(\xi_i) \in \bar{H}'$
is a group isomorphism. We have $\bar{H}\subseteq\bar{H}'$, with equality
when $K$ is algebraically closed.
\end{prop}

\begin{proof} It is clear that the map $K^\times \rightarrow \GL(\bar{V})$, 
$\xi\mapsto \bar{h}_i'(\xi)$, is a group homomorphism. Hence, $\varphi
\colon (K^\times)^I \rightarrow \bar{H}'$ is a surjective homomorphism. 
In order to show that $\varphi$ is injective, we construct an inverse 
map. For each $\mu \in P_\fh(V)$ we choose a basis vector $b_\mu \in 
\cB \cap V_\mu$. Then we define $\pi_\mu\colon \bar{H}' \rightarrow 
K^\times$ by sending $\eta\in \bar{H}'$ to the $(\bar{b}_\mu,
\bar{b}_\mu)$-diagonal entry of the matrix of $\eta$ with respect to 
$\bar{\cB}$. Thus, we have 
\[ \pi_\mu(\eta)=\prod_{i \in I} \xi_i^{\mu(h_i')} \qquad
\mbox{for} \quad \eta=\prod_{i \in I} \bar{h}_i'(\xi_i)\in \bar{H}'.\]
Now let $\{\mu_j\mid j \in I\}\subseteq \Lambda(V)$ be as in 
Remark~\ref{factorAa}; then $\mu_j(h_i')=\delta_{ij}$ for all $i,j\in I$. 
Since $\Lambda(V)=\langle P_\fh(V)\rangle_\Z$, we can write 
\[\mu_j=\sum_{\mu \in P_\fh(V)} n_{j,\mu}\mu \qquad \mbox{where}\qquad
n_{j,\mu} \in \Z. \]
Then we define a map $\psi_j\colon \bar{H}'\rightarrow K^\times$ by 
\[ \psi_j(\eta):=\prod_{\mu \in P_\fh(V)} \pi_\mu(\eta)^{n_{j,\mu}} 
\qquad (\eta\in \bar{H}').\]
Since each $\pi_\mu$ is a group homomorphism, it is clear that $\psi_j$
also is a group homomorphism. Now let $\eta\in \bar{H}'$ and write 
$\eta=\prod_{i \in I} \bar{h}_i'(\xi_i)$ where $\xi_i\in K^\times$.
Then, for $j\in I$, we compute:
\begin{align*}
\psi_j(\eta)&=\prod_{\mu \in P_\fh(V)} \pi_\mu(\eta)^{n_{j,\mu}}=
\prod_{\mu\in P_\fh(V)} \prod_{i \in I} \xi_i^{\mu(h_i')n_{j,\mu}}\\
&=\prod_{i \in I} \xi_i^{\sum_{\mu \in P_\fh(V)} \mu(h_i')n_{j,\mu}}
=\prod_{i \in I} \xi_i^{\mu_j(h_i')}=\xi_j.
\end{align*}
It follows that the maps $\varphi\colon (K^\times)^I\rightarrow\bar{H}'$ and 
\[ \psi\colon \bar{H}'\rightarrow (K^\times)^I,\qquad \eta \mapsto 
\bigl(\psi_j(\eta)\bigr)_{j \in I},\]
are inverse to each other. In particular, each $\eta\in \bar{H}'$ has a 
unique expression $\eta=\prod_{i \in I} \bar{h}_i'(\xi_i)$ with $\xi_i
\in K^\times$ for $i\in I$.

Next we show that $\bar{H}\subseteq \bar{H}'$. For $i,j\in I$ let 
$\breve{b}_{ij}\in \Z$ and $r_{ij}\in \Q$ be defined as in 
Remark~\ref{factorAa}. Then $h_i=\sum_{j \in I} \breve{b}_{ij}h_j'$ 
and so 
\[ \xi^{\mu(h_i)}=\prod_{j \in I} \xi^{\breve{b}_{ij}\mu(h_j')}
=\prod_{j \in I} \bigl(\xi^{\breve{b}_{ij}}\bigr)^{\mu(h_j')}=
\prod_{j \in I} \xi^{\breve{b}_{ij}\mu(h_j')}\]
for any $\mu \in P_\fh(V)$ and $\xi \in K^\times$. Now let $b\in 
\cB\cap V_\mu$. Then
\[ \Bigl(\prod_{j \in I} \bar{h}_j'(\xi^{\breve{b}_{ij}})\Bigr)(\bar{b})=
\Bigl(\prod_{j \in I} \xi^{\breve{b}_{ij}\mu(h_j')}\Bigr)\bar{b}=
\xi^{\mu(h_i)}\bar{b}=\bar{h}_i(\xi)(\bar{b})\]
where the last equality holds by Corollary~\ref{polylev2}. Thus, we have 
\[\bar{h}_i(\xi)=\prod_{j \in I} \bar{h}_j'(\xi^{\breve{b}_{ij}})\in
\hat{H}',\]
as desired. Finally, assume that $K$ is algebraically closed and that 
we are given a tuple $(\xi_i)_{i\in I} \in (K^\times)^I$. Let $0\neq n
\in \Z$ be such that $nr_{ij}\in \Z$ for all $i,j\in I$. We can find 
$\tilde{\xi}_i\in K^\times$ such that $\tilde{\xi}_i^{\,n}=\xi_i$ for 
all $i \in I$. (This is the place, and the only place, where we use 
that $K$ is algebraically closed.) Then set $\gamma_i:=\prod_{l \in I} 
\tilde{\xi}_l^{\,nr_{li}}\in K^\times$. With these definitions, for a 
fixed $j\in I$ we have
\[ \prod_{i\in I} \gamma_i^{\breve{b}_{ij}}=\prod_{i,l\in I}
\tilde{\xi}_l^{\,nr_{li}\breve{b}_{ij}}=\prod_{l\in I} 
\tilde{\xi}_l^{\,\sum_{i \in I} nr_{li}\breve{b}_{ij}}=\prod_{l \in I} 
\tilde{\xi}_l^{\,n\delta_{jl}}=\tilde{\xi}_j^{\,n}=\xi_j.\]
It follows that, for any $\mu \in P_\fh(V)$, we have 
\[ \prod_{j \in I} \xi_j^{\mu(h_j')}=\prod_{i,j\in I} 
\gamma_i^{\breve{b}_{ij}\mu(h_j')}=\prod_{i \in I} \gamma_i^{\sum_{j 
\in I} \breve{b}_{ij} \mu(h_j')}=\prod_{i \in I} \gamma_i^{\mu(h_i)},\]
where we use again that $h_i=\sum_{j \in I} \breve{b}_{ij}h_j'$. As
above, this implies that 
\[ \varphi\bigl((\xi_j)_{j\in I}\bigr)=\prod_{j\in I} \bar{h}_j'(\xi_j)=
\prod_{i\in I} \bar{h}_i(\gamma_i) \in \bar{H}.\]
Thus, the image of $\varphi$ is contained in $\bar{H}$. 
\end{proof}

Finally, we consider the relation between $\bar{H}$ and $\bar{N}$. For 
this purpose, it will be convenient to define ``monomial'' elements for 
all roots $\alpha\in \Phi$. Thus, for $\xi \in K^\times$ we set 
\[\bar{n}_\alpha(\xi):=\bar{x}_\alpha(\xi)\bar{y}_\alpha\bigl(-
(-1)^{\hgt(\alpha)}\xi^{-1}\bigr) \bar{x}_\alpha(\xi)\quad\in\bar{G}.\]
(The extra factor $(-1)^{\hgt(\alpha)}$ is required because, in our 
setting, we have $[\be_\alpha^+,\be_{-\alpha}^+]=(-1)^{\hgt(\alpha)}
h_\alpha$; see Corollary~\ref{canbash}.) In particular, for $i\in I$ we
have $\be_{\alpha_i}^+=\epsilon(i)e_i$ and $\be_{-\alpha_i}^+=-
\epsilon(i)f_i$. This yields that 
\[ \bar{n}_{\alpha_i}(\xi)=\bar{n}_i(\epsilon(i)\xi)\qquad
\mbox{for $i\in I$}.\]


\begin{lem} \label{generalna1} Let $\alpha\in \Phi$, $\xi\in K^\times$. 
Then $\bar{n}_\alpha(\xi)^{-1}=\bar{n}_\alpha (-\xi)$ and 
\[\bar{n}_i(1)\bar{n}_\alpha(\xi)\bar{n}_i(1)^{-1}=
\bar{n}_{s_i(\alpha)}\bigl(c_i(\alpha)\xi\bigr)\qquad 
\mbox{for $i\in I$}.\]
\end{lem}

\begin{proof} The identity $\bar{n}_\alpha(\xi)^{-1}=\bar{n}_\alpha
(-\xi)$ immediately follows from the fact that $\bar{x}_\alpha
(\xi)^{-1}=\bar{x}_\alpha(-\xi)$. By Proposition~\ref{Cree0a}, we have 
\begin{align*}
\bar{n}_i(1)\bar{x}_\alpha(\xi)&\bar{n}_i(1)^{-1}=\bar{x}_{s_i
(\alpha)}\bigl(c_i(\alpha)\xi\bigr),\\
\bar{n}_i(1)\bar{x}_{-\alpha}(-(-1)^{\hgt(\alpha)}\xi^{-1})&
\bar{n}_i(1)^{-1}=\bar{x}_{-s_i(\alpha)}\bigl(-c_i(-\alpha)
(-1)^{\hgt(\alpha)}\xi^{-1}\bigr).
\end{align*}
So the desired identity holds if we can show that 
\[c_i(-\alpha)(-1)^{\hgt(\alpha)}=c_i(\alpha)(-1)^{\hgt(s_i(\alpha))}.\]
But this follows from the formula $s_i(\alpha)=\alpha-\langle \alpha_i^\vee, 
\alpha\rangle\alpha_i$ and the fact that $c_i(\alpha) c_i(-\alpha)=
(-1)^{\langle \alpha_i^\vee, \alpha\rangle}$ (see Proposition~\ref{Cree0}). 
\end{proof}

\begin{lem} \label{generalna2} Let $\alpha,\beta\in \Phi$ and $w\in W$
be such that $w(\alpha)=\beta$. Write $w=s_{i_1}\cdots s_{i_r}$
where $i_1,\ldots,i_r\in I$ ($r\geq 0$). Then, setting $g:=\bar{n}_{i_1}(1)
\cdots \bar{n}_{i_r}(1)\in \bar{N}$, we have 
\[ g \bar{x}_\alpha(\zeta)g^{-1}=\bar{x}_\beta(\pm \zeta)\qquad
\mbox{and}\qquad g \bar{n}_\alpha(\xi)g^{-1}=\bar{n}_\beta(\pm \xi)\]
for all $\zeta\in K$ and $\xi\in K^\times$, where the signs do not depend
on $\zeta$ or $\xi$.
\end{lem}

\begin{proof} The first identity follows by a repeated application of
Propo\-si\-tion~\ref{Cree0a}, as in the proof of Corollary~\ref{corCree0a}.
Similarly, the second identity follows using a repeated application of 
Lemma~\ref{generalna1}.
\end{proof}

\begin{xca} \label{xcageneralna} Let $\alpha\in \Phi$ and $\xi\in 
K^\times$. Show that $\bar{n}_\alpha(\xi)\in \bar{N}$; furthermore, 
$\bar{h}_\alpha(\xi):=\bar{n}_\alpha(\xi) \bar{n}_\alpha(-1)\in\bar{H}$.\\
{\footnotesize [{\it Hint}. Write $\alpha=w(\alpha_j)$ where $w\in W$
and $j\in I$. Then use Lemma~\ref{generalna2}.]}
\end{xca}

\begin{thm}[Braid relations] \label{HHp2neu} \nmi{}{braid relations} 
Let $i,j\in I$, $i\neq j$, and $m_{ij} \in \{2,3,4,6\}$ be the order of
$s_is_j\in W$. (See Exercise~\ref{xcabraid}.) Then the following identity 
holds in $\bar{N}$:
\[\underbrace{\bar{n}_i(1)\cdot \bar{n}_j(1) \cdot \bar{n}_i(1) 
\cdots}_{\text{$m_{ij}$ terms}} =\underbrace{\bar{n}_j(1)\cdot 
\bar{n}_i(1)\cdot \bar{n}_j(1) \cdots}_{\text{$m_{ij}$ terms}}.\]
\end{thm}

\begin{proof} For any integer $r\geq 1$ we set 
\begin{align*}
\bar{n}_i(r)&:=\bar{n}_i(1)\cdot \bar{n}_j(1) \cdot \bar{n}_i(1) 
\cdots\in \bar{N},\\
\bar{n}_j(r)&:=\bar{n}_j(1)\cdot \bar{n}_i(1) \cdot \bar{n}_j(1) 
\cdots\in \bar{N},
\end{align*}
both with $r$ terms in the product. So we must show that $\bar{n}_i(m_{st})
=\bar{n}_j(m_{st})$. For this purpose, let $g:=\bar{n}_i(m_{st}-1)$. 
Thus, $\bar{n}_j(1)g$ is the right hand side of the desired identity, 
while the left hand side is $g\bar{n}_l(1)$ where $l=i$ if 
$m_{ij}$ is odd, and $l=j$ if $m_{ij}$ is even. 

Let $w:=s_is_js_i\cdots \in W$, again with $m_{ij}-1$ terms in the 
product. Then $s_jw=s_js_is_j \cdots$ and $ws_l=s_is_js_i\cdots$ (both 
with $m_{ij}$ factors). Hence, $s_jw=ws_l$ since $(s_is_j)^{m_{ij}}=
\id_W$; see also Exercise~\ref{xcabraid}. By Remark~\ref{genrefl}, the 
equality $s_j=ws_lw^{-1}$ implies that $w(\alpha_l)=\pm \alpha_j$. By 
Exercise~\ref{matsum1}, we have $\ell(ws_l)=\ell(s_jw)=m_{st}$ and so 
$\ell(ws_l)=\ell(w)+1$. Hence, Lemma~\ref{exccond}(b) shows that 
$w(\alpha_l)=\alpha_j$. So the assumptions of Lemma~\ref{generalna2} are 
satisfied for the roots $\alpha=\alpha_l$, $\beta=\alpha_j$ and the 
above elements $w,g$. 

Hence, $g\bar{n}_l(1)=\bar{n}_j(\delta)g$ where $\delta\in \{\pm 1\}$;
note also that $\bar{n}_l(1)=\bar{n}_{\alpha_l}(\pm 1)$ and
$\bar{n}_j(1)=\bar{n}_{\alpha_j}(\pm 1)$. Furthermore, 
$\bar{h}_j(\delta)=\bar{n}_j(\delta)\bar{n}_j(-1)\in \bar{H}$ and 
$\bar{n}_j(-1)=\bar{n}_j(1)^{-1}$. This yields the identity 
\[ \bar{n}_i(m_{st})=g\bar{n}_l(1)=\bar{n}_j(\delta)g=
\bar{h}_j(\delta)\bar{n}_j(1)g=\bar{h}_j(\delta)\bar{n}_j(m_{st})\]
which is almost the identity that we are trying to prove, but there is
an extra factor $\bar{h}_j(\delta)$. Now, we can reverse the roles of $i$ 
and~$j$. Then an analogous argument shows that there is a sign $\delta'\in 
\{\pm 1\}$ such that 
\[ \bar{n}_j(m_{st})=\bar{h}_i(\delta')\bar{n}_i(m_{st}).\]
Consequently, $\bar{h}_i(\delta')\bar{h}_j(\delta)=\id_{\bar{V}}$.
Since $\bar{h}_j(\delta)^2=\id_{\bar{V}}$, we conclude that $h:=
\bar{h}_i(\delta')=\bar{h}_j(\delta)$. We claim that this implies that 
$h=\id_{\bar{V}}$.

Indeed, let $\mu\in P_\fg(V)$ and $b\in \cB\cap V_{\mu}$. Using the formula 
in Corollary~\ref{polylev2}, we obtain that 
\[ \delta'^{\langle \alpha_i^\vee,\mu\rangle}=\bar{h}_i(\delta')(\bar{b})
=h(\bar{b})=\bar{h}_j(\delta)(\bar{b})=\delta^{\langle \alpha_j^\vee,
\mu\rangle}.\]
Since this holds for all $\mu \in P_\fh(V)$, we conclude that 
\[ \delta'^{\langle \alpha_i^\vee,\mu\rangle}=\delta^{\langle 
\alpha_j^\vee, \mu\rangle} \qquad\mbox{for all $\mu\in \Lambda(V)$}.\]
But $\Phi\subseteq \Lambda(V)$ and so we also have 
$\delta'^{\langle \alpha_i^\vee,\alpha_k\rangle}=\delta^{\langle 
\alpha_j^\vee, \alpha_k\rangle}$ for all $k\in I$. Thus, since $\langle 
\alpha_i^\vee,\alpha_k\rangle=a_{ik}$ and $\langle \alpha_j^\vee,
\alpha_k\rangle=a_{jk}$, we have 
\[ \delta'^{a_{ik}}=\delta^{a_{jk}}\qquad \mbox{for all $k\in I$}.\]
Now, if $a_{ij}\neq 0$, then $a_{ij}=-1$ or $a_{ji}=-1$. In the first 
case, we choose $k=j$ and obtain $\delta'=\delta^{a_{jj}}=\delta^2=1$; 
hence, $h=\bar{h}_i(\delta')=\bar{h}_i(1)=\id_{\bar{V}}$. Similarly, in 
the second case, we choose $k=i$ and obtain $h=\bar{h}_j(\delta)=
\bar{h}_j(1)=\id_{\bar{V}}$. Finally, if $a_{ij}=0$, then we are in the 
situation of Proposition~\ref{action2a}(c) (where $\alpha=\alpha_i$
and $\beta=\alpha_j$). We conclude that $\bar{x}_i(\zeta)$ commutes 
with $\bar{x}_j(\xi)$ and with $\bar{y}_j(\xi)$ for all $\zeta,\xi\in K$.
Hence, $\bar{n}_i(1)$ and $\bar{n}_j(1)$ also commute with each other. But 
this is exactly the desired identity in this case, since $m_{ij}=2$ if
$a_{ij}=0$; see Example~\ref{xcabraid}. 
\end{proof}

Recall from Corollary~\ref{diagII} that $\bar{H}$ is a normal 
subgroup of~$\bar{N}$. Now consider the factor group $\bar{N}/\bar{H}$. 
In order to avoid any danger of confusion, we denote the cosets in 
$\bar{N}/\bar{H}$ by $[g]:=g\bar{H}$ for $g\in \bar{N}$.

\begin{cor} \label{HHp3} There exists a group isomorphism 
$\psi\colon W \rightarrow \bar{N}/\bar{H}$ such that $\psi(s_i)=
[\bar{n}_i(1)]$ for all $i\in I$. 
\end{cor}

\begin{proof} For $i\in I$ we have $\bar{n}_i(1)^2=\bar{h}_i(-1)$; see 
Corollary~\ref{Cree0ac}(c). Hence, $[\bar{n}_i(1)]^2=[\id_{\bar{V}}]$. 
Since the ``braid relations'' in Theorem~\ref{HHp2neu} also hold, it 
follows by a general argument that there exists a group homomorphism 
$\psi\colon W \rightarrow \bar{N}/\bar{H}$ such that $\psi(s_i)=
[\bar{n}_i(1)]$ for all $i\in I$; see ``Matsumoto's Lemma'' in 
Appendix~\ref{asecmats}. Using Corollary~\ref{diagII}, we see that 
$\psi$ is surjective. 

To show injectivity, let $w\in W$ be such that $\psi(w)=[\id_{\bar{V}}]$. 
Write $w=s_{i_1}\cdots s_{i_r}$ where $r\geq 0$ and $i_1,\ldots,i_r\in I$. 
Let $g:=\bar{n}_{i_1}(1)\cdots \bar{n}_{i_r}(1) \in \bar{N}$. 
Assume, if possible, that $w\neq \id$. Then there exists some $i \in I$ 
such that $w(\alpha_i)\in \Phi^-$. A repeated application of 
Proposition~\ref{Cree0a}(a) shows that $g\bar{x}_{\alpha_i}(1)g^{-1}=
\bar{x}_{w(\alpha_i)}(\pm 1)$. Hence, $g\bar{x}_{\alpha_i}(1)g^{-1}$ is 
represented by a lower triangular matrix with~$1$ along the diagonal.
But we also have $[\id_{\bar{V}}]=\psi(w)=[g]$ and so $g\in \bar{H}$. 
Then a repeated application of Proposition~\ref{Cree0a}(b) shows that 
$g\bar{x}_{\alpha_i}(1)g^{-1}=\bar{x}_{\alpha_i}(\pm 1)$. Hence, 
$g\bar{x}_{\alpha_i}(1)g^{-1}$ is also represented by an upper triangular 
matrix with~$1$ along the diagonal. We conclude that $g\bar{x}_{\alpha_i}
(1)g^{-1}=\id_{\bar{V}}$ and, hence, $\bar{x}_{\alpha_i}(1)=\id_{\bar{V}}$, 
contradiction to Corollary~\ref{cornontrivrhoi}. Thus, we must have 
$w=\id$, as required.
\end{proof}

\begin{rem} \label{canrepnw} The fact that the above braid relations 
hold has the following consequence. Let $w\in W$ and choose any reduced 
expression $w=s_{i_1} \cdots s_{i_r}$ where $\ell(w)=r$ and $i_1,\ldots,
i_r\in I$. Then, by Proposition~\ref{matsum}, the element $\bar{n}_w:=
\bar{n}_{i_1}(1) \cdots \bar{n}_{i_r}(1)\in \bar{N}$ does not depend 
on the choice of the reduced expression for~$w$. Thus, for each $w\in W$, 
we have a \textit{canonical} representative $\bar{n}_w\in \bar{N}$. 
Consequently, for $w\in W$ and $i \in I$, we have the following formula
(which allows for a recursive computation of $n_w$, starting with 
$\bar{n}_\id =\id_{\bar{V}}$):
\[ \bar{n}_i(1)\bar{n}_w=\left\{\begin{array}{cl} \bar{n}_{s_iw} & 
\quad \mbox{if $\ell(s_iw)=\ell(w)+1$},\\ \bar{h}_i(-1)\bar{n}_{s_iw} & 
\quad \mbox{if $\ell(s_iw)= \ell(w)-1$}. \end{array}\right.\]
Indeed, if $\ell(s_iw)=\ell(w)+1$, then a reduced expression for $s_iw$
is obtained by taking any reduced expression for $w$ and multiplying
on the left by $s_i$; hence, $\bar{n}_{s_iw}=\bar{n}_i(1)\bar{n}_w$ in 
this case. If $\ell(s_iw)=\ell(w)-1$ then we take a reduced expression 
$s_iw=s_{i_1} \cdots s_{i_r}$ where $r=\ell(w)-1$. Consequently, 
$w=s_is_{i_1}\cdots s_{i_r}$ is a reduced expression for $w$ and so 
$\bar{n}_w=\bar{n}_i(1)\bar{n}_{i_1}(1) \cdots \bar{n}_{i_r}(1)=\bar{n}_i
(1)\bar{n}_{s_iw}$. Since $\bar{n}_i(1)^2=\bar{h}_i(-1)$, this yields the 
desired formula for $\bar{n}_{s_iw}$.
\end{rem}

\section{Chevalley groups of type $A_1$} \label{secA1}

The simplest examples of Chevalley groups are the groups associated with 
the Lie algebra $\fg=\slm_2(\C)$. This is also the first example in 
Carter \cite[\S 4.5]{Ca1}, where it is shown that the corresponding 
groups of adjoint type are isomorphic to $\SL_2(K)/\{\pm I_2\}$. (We 
have seen this here in Proposition~\ref{propree1}, even for $\fg=
\slm_n(\C)$, any $n\geq 2$.) One of the purposes of this section is to 
describe {\itshape all} possible Chevalley groups $\bar{G}=G_K(V,\cB)$ 
associated with $\fg=\slm_2(\C)$, where $V$ is {\itshape any} faithful 
$\fg$-module with a regular basis~$\cB$. These groups will also play a role
in the study of arbitrary Chevalley groups. Along the way, we introduce 
some general methods (due to Ree \cite{ree2}) that will be useful later on.

Until further notice, let $\fg$ be a Lie algebra and $\fh\subseteq\fg$ 
be a subalgebra such that $(\fg,\fh)$ is of Cartan--Killing type, with 
structure matrix $A=(a_{ij})_{i,j\in I}$. Let $V$ be a faithful 
$\fg$-module and $\cB$ be a regular basis of~$V$. Given a field $K$, let 
\[ G_K(V,\cB)=\bigl\langle \bar{x}_\alpha(\zeta;V,\cB)\mid \alpha
\in \Phi,\zeta\in K\bigr\rangle\subseteq \GL(\bar{V})\]
be the corresponding Chevalley group. We now have to carry $V,\cB$ along 
in our notation, because we will want to compare groups associated with 
different $V,\cB$. 

We say that a $\fg$-submodule $V'\subseteq V$ is a \nm{regular submodule} 
if $V'$ is also faithful and if $\cB':=\cB\cap V'$ is a regular basis of 
$V'$. In this case, we can also form the Chevalley group $G_K(V',\cB')$. 

\begin{lem}[Cf.\ Ree \protect{\cite[(4.13)]{ree2}}] \label{mree1} Let
$V'\subseteq V$ be a regular submodule, as above. Then there is a unique
surjective group homomorphism $\pi' \colon G_K(V,\cB)\rightarrow 
G_K(V',\cB')$ such that 
\[ \bar{x}_\alpha(\zeta;V,\cB)\mapsto \bar{x}_\alpha(\zeta;V',\cB')
\qquad \mbox{for $\alpha\in \Phi$ and $\zeta\in K$}.\]
\end{lem}

\begin{proof} As in Section~\ref{relsGp}, we set $\be_\alpha^{[m]}:=
\frac{1}{m!}\rho(\be_\alpha^+)^m \in \End(V)$ for any $\alpha\in \Phi$ 
and any integer $m\geq 0$. Since $V'\subseteq V$ is a submodule, we have 
a corresponding representation $\rho'\colon \fg\rightarrow \gl(V')$ such 
that $\rho'(x)(v)=\rho(x)(v)\in V'$ for all $v\in V'$. We also set 
\[{\be_\alpha'}^{\!\![m]}:=\textstyle{\frac{1}{m!}\rho'(\be_\alpha^+)^m \in 
\End(V')}\qquad \mbox{for $\alpha\in \Phi$ and $m\geq 0$}.\]
Then $\be_{\alpha}^{[m]}(b)={\be_\alpha'}^{\!\![m]}(b)$ for all $b\in 
\cB'$. Now consider $\bar{V}=K\otimes_\Z \langle \cB\rangle_\Z$
and $\bar{\cB}=\{1\otimes b \mid b \in \cB\}$. We may certainly identify 
$\bar{V}'=K\otimes_\Z \langle \cB'\rangle_\Z$ with a subspace of $\bar{V}$ 
and $\bar{\cB}'=\{1\otimes b\mid b\in\cB'\}$ with a subset of~$\bar{\cB}$.
Hence, if $\overline{\be}_\alpha^{[m]}\colon \bar{V} \rightarrow \bar{V}$ 
and ${\overline{\be}_\alpha'}^{\!\![m]}\colon \bar{V}' \rightarrow \bar{V}'$ 
are the induced linear maps, then we also have
\[\overline{\be}_{\alpha}^{[m]}(1\otimes b)=
{\overline{\be}_\alpha'}^{\!\![m]}(1\otimes b)\in \bar{V}'
\qquad\mbox{for $b\in\cB'$}.\]
Consequently, by the definitions of $\bar{x}_\alpha(\zeta;V,\cB)$ and 
$\bar{x}_\alpha(\zeta; V',\cB')$, the above identity implies that 
\begin{equation*}
\bar{x}_\alpha(\zeta;V,\cB)(1\otimes b)=\bar{x}_\alpha(\zeta;V',\cB')
(1\otimes b)\in \bar{V}'\tag{$*$}
\end{equation*}
for $b \in \cB'$ and $\zeta\in K$. Since $\bar{V}'=\langle \bar{\cB}'
\rangle_K$, we conclude that $g(\bar{V}')\subseteq \bar{V}'$ for all 
$g\in G_K(V',\cB')$. This means that each $g\in G_K(V,\cB)\subseteq
\GL(\bar{V})$ induces an invertible linear map $g' \colon \bar{V}'
\rightarrow \bar{V}'$, simply by restriction. So we obtain a map  
\[\pi'\colon G_K(V,\cB)\rightarrow\GL(\bar{V}'),\qquad g\mapsto g',\]
which clearly is a group homomorphism. Now ($*$) shows that $\pi'$ sends 
the element $\bar{x}_\alpha(\zeta;V,\cB)\in G_K(V,\cB)$ to the element 
$\bar{x}_\alpha(\zeta;V',\cB')\in G_K(V',\cB')$, for all $\alpha\in 
\Phi$ and $\zeta\in K$, as required. This also shows that $\pi'$ is 
unique and that $\pi'(G_K(V,\cB))=G_K(V',\cB')$. 
\end{proof}

Now the problem is that, at this stage, we can not say much about the kernel
of $\pi'$. However, we can at least describe the intersection of 
$\ker(\pi')$ with the diagonal subgroup 
\[H_K(V,\cB):=\Bigl\{\prod_{i \in I} \bar{h}_i(\xi_i;V,\cB)\,\big|\,
\xi_i\in K^\times \mbox{ for all $i\in I$}\Bigr\}\]
of $G_K(V,\cB)$. Let $P_\fh(V)$ be the set of weights of $\fh$ on $V$.
Recall from Definition~\ref{deflambda} that we have the \nm{weight lattice}
of $V$:
\[ \Lambda(V):=\langle P_\fh(V)\rangle_\Z \subseteq \Omega.\]
Similarly, $\Lambda(V')=\langle P_\fh(V')\rangle_\Z\subseteq 
\Omega$ is the weight lattice of $V'$. Note that, clearly, we have 
$P_\fh(V') \subseteq P_\fh(V)$ and $\Lambda(V') \subseteq \Lambda(V)$.

\begin{lem} \label{mree2} Let $h:=\prod_{i \in I} \bar{h}_i(\xi_i;V, 
\cB)\in H_K(V,\cB)$ where $\xi_i\in K^\times$ for all $i \in I$. Then 
\[\pi'(h)=\id_{\bar{V}'}\qquad \Leftrightarrow\qquad \prod_{i\in I} 
\xi_i^{\langle \alpha_i^\vee,\mu\rangle}=1 \quad \mbox{for all $\mu 
\in \Lambda(V')$}.\]
In particular, if $\Lambda(V')=\Lambda(V)$, then $\ker(\pi')\cap 
H_K(V,\cB)=\{\id_{\bar{V}}\}$.
\end{lem}

\begin{proof} Since $\pi'(\bar{x}_\alpha(\zeta;V,\cB))=\bar{x}_\alpha
(\zeta;V',\cB')$ for $\alpha\in \Phi$ and $\zeta\in K$, we also have 
\begin{align*}
\pi'(\bar{x}_i(\zeta;V,\cB))&=\bar{x}_i(\zeta;V',\cB'), \\
\pi'(\bar{y}_i(\zeta;V,\cB))&=\bar{y}_i(\zeta;V',\cB')
\end{align*}
for $i\in I$ and $\zeta\in K$. This implies that 
\begin{align*}
\pi'\bigl(\bar{n}_i(\xi;V,\cB)\bigr)&=\bar{n}_i(\xi;V',\cB'),\\
\pi'\bigl(\bar{h}_i(\xi;V,\cB)\bigr)&=\bar{h}_i(\xi;V',\cB')
\end{align*}
for $i\in I$ and $\xi\in K^\times$ (simply by the definitions of these 
elements). Hence, $\pi'(h)=\prod_{i \in I} \bar{h}_i(\xi_i;V',\cB')$
and so, by Lemma~\ref{uniqueH}, we obtain the above equivalence.
And if $\Lambda(V')=\Lambda(V)$, then Lemma~\ref{uniqueH} also 
shows that $\pi'(h)=\id_{\bar{V}'} \Leftrightarrow h=\id_{\bar{V}}$.
\end{proof}


%
%
%

Let us now apply the above results in the case where $\fg=\slm_2(\C)$,
with $I=\{1\}$ and standard basis $\{e_1,f_1,h_1\}$ such that $[e_1,f_1]
=h_1$. Let $V$ be a faithful $\fg$-module with a regular basis $\cB$. The 
corresponding Chevalley group is given by $G_K(V,\cB)=\bigl\langle 
\fX_1,\fY_1\bigr\rangle$ where 
\begin{align*}
\fX_1&:=\{\bar{x}_1(\zeta;V,\cB)\mid \zeta \in K\},\\
\fY_1&:=\{\bar{y}_1(\zeta;V,\cB)\mid \zeta \in K\}.
\end{align*} 
Let us also define $\fH_1:=\{\bar{h}_1(\xi;V,\cB)\mid \xi\in K^\times\}$. 
Then, by Remark~\ref{sl2triple}, we have the factorisation
\[G_K(V,\cB)=\fX_1\cdot \fY_1 \cdot \fH_1\cdot \fX_1,\] 
which plays a crucial role in the following proof.

\begin{lem} \label{mree4} Let $\fg=\slm_2(\C)$ and $V,\cB$ as above. 
Let $V'\subseteq V$ be a regular submodule and $\pi'\colon G_K(V,\cB)
\rightarrow G_K(V',\cB')$ be the homomorphism of Lemma~\ref{mree1}. 
Then $\ker(\pi')\subseteq \fH_1$.
\end{lem}

\begin{proof} We shall assume that the elements of $\cB$ are arranged 
as in Remark~\ref{xcadet1}. In the following, when we consider matrices, 
these are always taken with respect to $\bar{\cB}$ or to $\bar{\cB}'$. 
For example, the matrix of $\bar{x}_1(\zeta;V,\cB)$ will be upper 
triangular with~$1$ along the diagonal; furthermore, by the proof of 
Lemma~\ref{mree1}, the matrix of $\bar{x}_1(\zeta;V',\cB')$ is
obtained by taking the submatrix of the matrix of $\bar{x}_1(\zeta;
V,\cB)$ with rows and columns corresponding to basis vectors in $\cB'
\subseteq \cB$. In particular, the matrix of $\bar{x}_1(\zeta;V',\cB')$ 
will also be upper triangular with $1$ along the diagonal. Similarly,
$\bar{y}_1(\zeta;V,\cB)$ and $\bar{y}_1(\zeta;V',\cB')$ are represented
by lower triangular matrices with $1$ along the diagonal. Furthermore, 
$\bar{h}_1(\xi;V, \cB)$ and $\bar{h}_1(\xi;V',\cB')$ are represented 
by diagonal matrices. Now we can argue as follows.

Let $g\in \ker(\pi')$. Since $G_K(V,\cB)=\fX_1\cdot \fY_1 \cdot \fH_1\cdot 
\fX_1$, we can write $g=xyh\tilde{x}$ where $x,\tilde{x}\in \fX_1$, $y\in 
\fY_1$ and $h\in \fH_1$. Now $\id_{\bar{V}'}=\pi'(g)= \pi'(x) \pi'(y) 
\pi'(h) \pi'(\tilde{x})$ and so 
\[ \pi'(y)=\pi'(x)^{-1}\pi'(\tilde{x})^{-1}\pi'(h)^{-1}.\]
By Lemma~\ref{mree1} and the above discussion, $\pi'(y)$ is represented by 
a lower triangular matrix with $1$ along the diagonal. Similarly, $\pi'(x)$ 
and $\pi'(\tilde{x})$ are represented by upper triangular matrices with~$1$ 
along the diagonal. Since $\pi'(h)$ is represented by a diagonal matrix, 
the right hand side of the above identity is represented by 
an upper triangular matrix. Hence, we must have $\pi'(y)=\id_{\bar{V}'}$. 
But $y=\bar{y}_1(\zeta;V,\cB)$ for some $\zeta\in K$ and so $\pi'(y)=
\bar{y}_1(\zeta;V',\cB')$. Since this equals $\id_{\bar{V}'}$, it 
follows that $\zeta=0$ by Corollary~\ref{cornontrivrhoi}. But then we 
also have $y=\bar{y}_1(0;V,\cB)=\id_{\bar{V}}$; hence, $g=zh$ where 
$z:=xh\tilde{x}h^{-1}\in \fX_1$. (Note that $\fH_1$ normalises $\fX_1$
by Remark~\ref{Cree0a0}.) Thus,
\[\mbox{\itshape every $g\in \ker(\pi')$ is represented by an upper 
triangular matrix}.\]
Now, since $\ker(\pi')$ is a normal subgroup, we also have $\tilde{g}:=
ngn^{-1}\in \ker(\pi')$, where $n:=\bar{n}_1(1; V,\cB)$. By 
Remark~\ref{Cree0a0} and Proposition~\ref{Cree0ac}(b), we have $n\fX_1
n^{-1}=\fY_1$ and $n\fH_1n^{-1}=\fH_1$. Hence, $\tilde{g}=\tilde{z}
\tilde{h}$ where $\tilde{z}=nzn^{-1}\in \fY_1$ and $\tilde{h}=nhn^{-1} 
\in \fH_1$. So $\tilde{g}$ is represented by a lower triangular matrix. 
But we have just seen that every element in $\ker(\pi')$ is represented by 
an upper triangular matrix. We conclude that $\tilde{g}=\tilde{z}\tilde{h}$ 
must be a diagonal matrix. This forces that $\tilde{z}=\id_{\bar{V}}$ and so $z=\id_{\bar{V}}$. 
Hence, $g=h\in \fH_1$. 
\end{proof}

\begin{thm}[Ree] \label{mree5} Let $\fg=\slm_2(\C)$ and $V$ be any 
faithful $\fg$-module with a regular basis~$\cB$. Then there exists a 
unique surjective group homomorphism $\pi\colon \SL_2(K) \rightarrow 
G_K(V,\cB)$ such that 
\[ \renewcommand{\arraystretch}{0.9} \renewcommand{\arraycolsep}{3pt}
\begin{pmatrix} 1 & \zeta \\ 0 & 1 \end{pmatrix} \mapsto 
\bar{x}_1(\zeta;V,\cB) \quad\mbox{and}\quad \begin{pmatrix} 1 & 0 
\\ \zeta & 1 \end{pmatrix} \mapsto \bar{y}_1(\zeta;V,\cB) \]
for all $\zeta\in K$. We have $\ker(\pi)=\{I_2\}$ if $\Lambda(V)=
\langle \varpi_1\rangle_\Z$ and $\ker(\pi)=\{\pm I_2\}$ if $\Lambda(V)=
\langle \alpha_1\rangle_\Z$ (notation as in Example~\ref{reesl2}). Thus, 
\[G_K(V,\cB)\cong \SL_2(K)\qquad\mbox{or} \qquad G_K(V,\cB)\cong 
\SL_2(K)/\{I_2\}.\]
\end{thm}

\begin{proof} Let $V_2=\C^2$ be the natural $\fg$-module with standard
basis $\cB_2=\{b_1,b_2\}$. Then $V_2$ is faithful and $\cB_2$ is a
regular basis. We have a unique group isomorphism $\SL_2(K)
\stackrel{\sim}{\longrightarrow} G_K(V_2;\cB_2)$ such that 
\[ \renewcommand{\arraystretch}{0.9} \renewcommand{\arraycolsep}{3pt}
\begin{pmatrix} 1 & \zeta \\ 0 & 1 \end{pmatrix} \mapsto 
\bar{x}_1(\zeta;V_2,\cB_2) \quad\mbox{and}\quad \begin{pmatrix} 1 & 0 
\\ \zeta & 1 \end{pmatrix} \mapsto \bar{y}_1(\zeta;V_2,\cB_2) \]
for $\zeta\in K$; see Example~\ref{defnonad1}(b). Now consider
the direct product $\tilde{V}=V\times V_2$; this is a $\fg$-module in a 
natural way (see Example~\ref{dirprodmod}). We may canonically identify 
$V$ and $V_2$ with submodules of $\tilde{V}$ such that $\tilde{V}=
V\oplus V_2$. Since $V$ and $V_2$ are faithful, $\tilde{V}$ is also 
faithful. Since $\cB$ and $\cB_2$ are regular bases of $V$ and $V_2$, 
respectively, it easily follows that $\tilde{\cB}:=\cB\cup \cB_2$ is 
a regular basis of $\tilde{V}$ (see Exercise~\ref{regdirprod}). Hence, 
the Chevalley group $G_K(\tilde{V},\tilde{\cB})$ is defined. By 
Exercise~\ref{xcaweights}, we have $P_\fh(\tilde{V})=P_\fh(V)\cup 
P_\fh(V_2)$ and, hence, 
\[ \Lambda(\tilde{V})=\langle P_\fh(V),P_\fh(V_2)\rangle_\Z=
\Lambda(V)+\Lambda(V_2)\subseteq \Omega.\]
Now, by Example~\ref{defhwm2}, we have $\Lambda(V_2)=\Omega$. So we 
conclude that we also have $\Lambda(\tilde{V})=\Omega$. Since $V$ and 
$V_2$ are regular submodules of~$\tilde{V}$, we have two surjective group 
homomorphisms
\[ \pi_1\colon G_K(\tilde{V},\tilde{\cB})\rightarrow G_K(V,\cB)
\;\mbox{and}\;\pi_2\colon G_K(\tilde{V},\tilde{\cB})\rightarrow 
G_K(V_2,\cB_2)\]
as in Lemma~\ref{mree1}. By Lemma~\ref{mree4}, we have $\ker(\pi_2)
\subseteq\tilde{\fH}_1$ where $\tilde{\fH}_1$ is the diagonal subgroup 
of $G_K(\tilde{V},\tilde{B})$. Furthermore, we have $\Lambda(\tilde{V})=
\Omega=\Lambda(V_2)$ and so Lemma~\ref{mree2} shows that $\pi_2$ is 
an isomorphism. So we obtain a surjective group homomorphism 
\[ \tilde{\pi}:=\pi_1\circ \pi_2^{-1}\colon G_K(V_2,\cB_2)\rightarrow
G_K(V,\cB)\]
such that $\bar{x}_1(\zeta;V_2,\cB_2)\mapsto \bar{x}_1(\zeta;V,\cB)$
and $\bar{y}_1(\zeta;V_2,\cB_2)\mapsto \bar{y}_1(\zeta;V,\cB)$ for 
all $\zeta\in K$. It remains to determine $\ker(\tilde{\pi})=
\ker(\pi_1)$. By Lemma~\ref{mree4}, we have again $\ker(\pi_1) \subseteq 
\tilde{\fH}_1$. Furthermore, Lemma~\ref{mree2} shows that $g=\bar{h}_1
(\xi;\tilde{V},\tilde{\cB})\in \tilde{\fH}_1$ (for $\xi \in K^\times$) 
belongs to $\ker(\pi_1)$ if and only if $\xi^{\langle \alpha_1^\vee,\mu
\rangle}=1$ for all $\mu \in \Lambda(V)$. Hence, if $\Lambda(V)=\langle
\varpi_1\rangle_\Z$, then $\xi=1$ and $g=\id$. On the other hand, if 
$\Lambda(V)=\langle \alpha_1 \rangle_\Z$, then $\xi^2=1$ and so 
$g=\bar{h}_1(\pm 1;\tilde{V}, \tilde{\cB})$. Now $\pi_2(g)=\bar{h}_1
(\pm 1; V_2,\cB_2)$ corresponds to the matrix $\pm I_2$ under the above 
isomorphism $G_K(V_2,\cB_2)\cong \SL_2(K)$. Hence, $\pi_2(g)=\pm \id$. 
But $\pi_2$ is an isomorphism and so $g=\pm \id$. Thus, we have shown 
that $\ker(\tilde{\pi})=\{\id\}$ if $\Lambda(V)=\langle\varpi_1\rangle_Z$, 
and $\ker(\tilde{\pi})=\{\pm \id\}$ if $\Lambda(V)=\langle \alpha_1 
\rangle_\Z$. Finally, composing $\tilde{\pi}$ with the above isomorphism 
$\SL_2(K) \rightarrow G_K(V_2; \cB_2)$, we obtain the required homomorphism 
$\pi\colon \SL_2(K) \rightarrow G_K(V,\cB)$, with the same kernel.
\end{proof}

Let again $\fg$ be an arbitrary Lie algebra of Cartan--Killing type, with 
structure matrix $A=(a_{ij})_{i,j\in I}$. For $V=\fg$ (adjoint 
representation), the following result is contained in Chevalley 
\cite[\S II]{Ch}; see also the (slightly different) exposition by 
Carter \cite[\S 6.3]{Ca1}. The general case is due to Ree \cite[(3.2)]{ree2}.
The proof that we give here is different from those in \cite{Ch}, 
\cite{Ca1}, \cite{ree2}.

\begin{cor} \label{mree6} Let $V$ be a faithful $\fg$-module and $\cB$ be 
a regular basis of~$V$. Then, for any $i \in I$, there exists a unique
group homomorphism $\pi_i\colon \SL_2(K)\rightarrow G_K(V,\cB)$ such that 
\[ \renewcommand{\arraystretch}{0.9} \renewcommand{\arraycolsep}{3pt}
\begin{pmatrix} 1 & \zeta \\ 0 & 1 \end{pmatrix} \mapsto 
\bar{x}_i(\zeta;V,\cB) \quad\mbox{and}\quad \begin{pmatrix} 1 & 0 
\\ \zeta & 1 \end{pmatrix} \mapsto \bar{y}_i(\zeta;V,\cB) \]
for all $\zeta\in K$; we have $\ker(\pi_i)\subseteq \{\pm I_2\}$. 
\end{cor} 

\begin{proof} We have $\slm_2(\C)\cong \langle e_i,f_i,h_i\rangle_\C$.
Via this isomorphism, we may regard $V$ as a faithful $\slm_2(\C)$-module; 
then, of course, $\cB$ is still a regular basis for this 
$\slm_2(\C)$-module. Let 
\[ G_{K,i}(V,\cB):=\bigl\langle \bar{x}_i(\zeta; V,\cB),\bar{y}_i(\zeta;
V, \cB)\mid \zeta \in K\bigr\rangle\subseteq \GL(\bar{V}). \]
Then $G_{K,i}(V,\cB)$ is contained in $G_K(V,\cB)$ on the one hand, but 
$G_{K,i}(V,\cB)$ is also the Chevalley group associated with $\slm_2(\C)$ 
and the $\slm_2(\C)$-module $V$ (with its regular basis $\cB$); note that
the endomorphisms $\bar{x}_i(\zeta; V,\cB)$ and $\bar{y}_i(\zeta;V,\cB)$ have
exactly the same definition in both cases. So the required homomorphism 
$\pi_i\colon\SL_2(K)\rightarrow G_K(V,\cB)$ is obtained by composing the 
homomorphism $\SL_2(K)\rightarrow G_{K,i}(V,\cB)$ from Theorem~\ref{mree5}
with the inclusion $G_{K,i}(V,\cB) \subseteq G_K(V,\cB)$. 
\end{proof}

\begin{rem} \label{mree7} The above result provides an ``explanation''
for the definition of the elements $\bar{n}_i(\zeta;V,\cB)$ and 
$\bar{h}_i(\zeta;V,\cB)$ in $G_K(V,\cB)$. Just note the following 
computations with matrices in $\SL_2(K)$:
\begin{align*}
\renewcommand{\arraystretch}{0.9}\renewcommand{\arraycolsep}{3pt}
N_i(\zeta)&:=\begin{pmatrix} 1 & \zeta \\ 0 & 1 \end{pmatrix} \cdot 
\begin{pmatrix} 1 & 0 \\ -\zeta^{-1} & 1 \end{pmatrix} \cdot
\begin{pmatrix} 1 & \zeta \\ 0 & 1 \end{pmatrix}=
\begin{pmatrix} 0 & \zeta \\ -\zeta^{-1} & 0 \end{pmatrix},\\
H_i(\zeta)&:=N_i(\zeta)\cdot N_i(-1)=
\begin{pmatrix} 0 & \zeta \\ -\zeta^{-1} & 0 \end{pmatrix} \cdot
\begin{pmatrix} 0 & -1 \\ 1 & 0 \end{pmatrix} =
\begin{pmatrix} \zeta & 0 \\ 0 & \zeta^{-1} \end{pmatrix},
\end{align*}
which explain why the elements in $G_K(V,\cB)$ corresponding to $N_i
(\zeta)$ and $H_i(\zeta)$ are called ``monomial'' and ``diagonal'',
respectively.
\end{rem}

Finally, we state the following, most general version of the existence
of homomorphisms between Chevalley groups associated with the same~$\fg$
but different pairs $(V,\cB)$ as above.

\begin{thm} \label{reest} Let $V_1$ and $V_2$ be faithful $\fg$-modules.
Let $\cB_1$ be a regular basis of $V_1$ and $\cB_2$ be a regular basis
of $V_2$. If $\Lambda(V_2)\subseteq \Lambda(V_1)$, then there is a
unique surjective group homomorphism 
\[ \phi\colon G_K(V_1,\cB_1) \rightarrow G_K(V_2,\cB_2)\]
such that $\bar{x}_\alpha(\zeta;V_1,\cB_1)\mapsto \bar{x}_\alpha(\zeta;
V_2,\cB_2)$ for $\alpha\in \Phi$ and $\zeta\in K$. The kernel of $\phi$ 
is contained in the center of $G_K(V_1,\cB_1)$ and consists of elements 
in the diagonal subgroup $H_K(V_1,\cB_1)$. Furthermore, $\phi$ is an 
isomorphism if $\Lambda(V_1)=\Lambda(V_2)$.
\end{thm}

We will prove this later; see Ree \cite[(3.10), 
(3.11)]{ree2} and Steinberg \cite[Cor.~5 (p.~29]{St} for the original 
proofs. Taking $V_1=V_2$ we deduce, in particular, that $G_K(V,\cB)$ is 
uniquely determined (up to isomorphism) by~$K$ and the module~$V$, and $G_K
(V,\cB)$ does not depend on the choice of the regular basis~$\cB$ of~$V$. 

The proofs of Ree and Steinberg are quite different. Ree uses a 
generalisation of the above proof of Theorem~\ref{mree5}, while Steinberg 
works with a set of defining relations; see \cite[Cor.~3 (p.~28)]{St}.
Below we give some further comments about these proofs.

\begin{defn} \label{defadjsc} If $\Lambda(V)$ equals $\Omega$, the full
weight lattice, then $G_K(V,\cB)$ is called a \nm{universal Chevalley 
group} (or \nm{Chevalley group of simply connected type}); see 
Humphreys \cite[\S 17.4]{H} or Steinberg \cite[p.~30]{St}. At the other 
extreme, if $\Lambda(V)$ equals $\langle \Phi\rangle_\Z$, the root 
lattice, then $G_K(V,\cB)$ is called an \nm{adjoint Chevalley group}
(or \nm{Chevalley group of adjoint type}). The above Theorem~\ref{reest} 
shows that, in the general case where $\Lambda(V)$ lies somewhere 
between $\langle\Phi\rangle_\Z$ and $\Omega$, there always exists a 
surjective homomorphism from the universal group onto $G_K(V,\cB)$, 
and a surjective homomorphism from $G_K(V,\cB)$ onto the adjoint group.
\end{defn}

Let us now sketch Ree's proof of Theorem~\ref{reest} and explain why 
we can not carry it out here and now. Given $V_1$ and $V_2$ such that
$\Lambda(V_2)\subseteq \Lambda(V_1)$, we form $V=V_1\oplus V_2$; then 
$\cB=\cB_1\cup \cB_2$ is a regular basis of~$V$. We have two surjective 
group homomorphisms
\[ \pi_1\colon G_K(V,\cB)\rightarrow G_K(V_1,\cB_1)
\;\mbox{and}\;\pi_2\colon G_K(V,\cB)\rightarrow G_K(V_2,\cB_2)\]
as in Lemma~\ref{mree1}. Furthermore, 
\[ \Lambda(V)=\langle P_\fh(V_1),P_\fh(V_2)\rangle_\Z=
\Lambda(V_1)+\Lambda(V_2)=\Lambda(V_1).\]
As in the proof of Theorem~\ref{mree5}, we would like to conclude 
that~$\pi_1$ is an isomorphism, which would allow us to define 
\[ \phi:=\pi_2\circ \pi_1^{-1} \colon G_K(V_1,\cB_1)\rightarrow 
G_K(V_2,\cB_2).\]
In order to show that $\pi_1$ is an isomorphism, it would be sufficient 
to show that $\ker(\pi_1)\subseteq H_K(V,\cB)$ (see Lemma~\ref{mree4} 
for the case $\fg=\slm_2(\C)$), because then we could use again 
Lemma~\ref{mree2} and complete the argument as before. Now, looking at 
the above proof of Lemma~\ref{mree4}, we see that it would be sufficient
to generalise the factorisation $G_K(V,\cB)=\fX_1\cdot \fY_1\cdot \fH_1
\cdot \fX_1$ (for $\fg=\slm_2(\C)$) to arbitrary~$\fg$. But this is exactly 
what seems to be difficult to obtain with the tools that are available to 
us now. Later we shall deduce the required factorisation in the general 
case from \nm{Chevalley's commutator relations}, which will be 
proven in a later chapter below.

Alternatively, we could try to follow Steinberg's argument but, 
again, it would not be possible to carry it out here and now because 
Chevalley's commutator relations form a subset of Steinberg's defining 
relations for $G_K(V,\cB)$. So, in either case, the commutator relations 
seem to be a crucial ingredient in the proof of Theorem~\ref{reest}.

\section*{Notes on Chapter~\ref{chap4}}

The material in Section~\ref{sechighw} is standard. But, generally
speaking, we tend to give more details in basic examples and initial 
steps of the general theory than seems to be usual. For example, the 
weight lattice $\Lambda(V)$ and Proposition~\ref{defomega} just appear 
as an exercise in Humphreys \cite[Exc.~21.5]{H}. The proof of 
Proposition~\ref{integralrep} (based on $\slm_2$-representations) is taken 
from Samelson \cite[\S 3.2, Theorem~B]{Sam}; an alternative proof is given 
by Lemma~\ref{Cree0aa}. The definition of $\fh_\Z$ in 
Remark~\ref{factorAa} appears in Steinberg \cite[Cor.~2 (p.~16)]{St}. 
For the discussion of highest weight modules we follow Serre 
\cite[Chap.~VII]{S}, but we have avoided the universal enveloping algebra, 
as in Samelson \cite[\S 3.2]{Sam}. There is much more to be said about this 
topic; see, e.g., the relevant chapters in Bourbaki \cite{B78}, 
Fulton--Harris \cite{FH}, Humphreys \cite{H} or Kac \cite[Chap.~9]{K}. The 
statement that {\itshape every} subgroup $\Lambda'\subseteq \Omega$
containing $\Phi$ arises as $\Lambda(V)$ for some $\fg$-module~$V$, is also 
contained in \cite[Exc.~21.5]{H}; there it can be deduced from the general 
results on finite-dimensional highest weight modules in \cite[\S 21.2]{H}. 
Here, we obtain that statement in Theorem~\ref{corj2}, using the results on 
modules with a minuscule highest weight. 

The development of the basic results on minuscule weights is very 
much inspired by Stembridge \cite{Stem1}, which contains a detailed study 
of the partially ordered set $(\Omega^+,\preceq)$; Lemma~\ref{lemstem1} 
and Lemma~\ref{lemstem4} appear in \cite[\S 1.1]{Stem1}. For a further 
discussion (e.g., the relation to the affine Weyl group), see Bourbaki 
\cite[Ch.~VI, \S 2, no.~3]{B}, \cite[Ch.~VIII, \S 7, no.~3]{B78} and 
Humphreys \cite[\S 13]{H}. The construction of modules with a minuscule 
highest weight in Definition~\ref{defj} is based on Jantzen \cite{Ja}; see 
also Vavilov \cite[\S 3]{vav2}. Note that the proof that the formulae in 
Definition~\ref{defj} indeed define a $\fg$-module structure on~$M$ 
essentially relies on Proposition~\ref{triang3b}~---~which is a weak 
version of Serre's theorem mentioned in Remark~\ref{serre}. The idea that 
one can obtain {\itshape all} possible Chevalley groups associated with a 
given simple Lie algebra $\fg$ by considering only the adjoint representation 
of $\fg$ and various minuscule representations is explicitly worked out 
in \cite{G2}. Here, this is extended to the case where $\fg$ is of 
Cartan--Killing type but not necessarily simple.

Proposition~\ref{Cree0} appears in \cite[(3.6)]{ree2} and \cite[Lemma~19(a),
p.~22]{St}; the proof here is closer to that in \cite{ree2}, where we use 
Theorem~\ref{luform3} to get some control over the signs $c_i(\alpha)$ 
occurring there. 
The proof of Proposition~\ref{diagauto} (concerning diagonal automorphisms
in general) works out a suggestion of Steinberg; see the exercise just 
after \cite[Lemma~58 (p.~92)]{St}. Proposition~\ref{steinbergp} is adapted
from the proof of \cite[Lemma~17]{St}. The proof of Proposition~\ref{Cree0ab} 
is also due to Steinberg \cite[Lemma~19]{St}. The proof of 
Corollary~\ref{HHp3} essentially follows the argument of Carter 
\cite[Theorem~7.2.2]{Ca1} and Steinberg \cite[Lemma~22 (p.~24)]{St}. But 
Carter and Steinberg use a certain ``non-standard'' presentation of the 
Weyl group~$W$; see \cite[\S 2.4]{Ca1}. Here, we use a slight variation 
based on the ``braid relations'' in Theorem~\ref{HHp2neu}. Our proof of 
these relations is an adaptation of the argument in 
\cite[Lemma~56 (p.~87)]{St}. Lemma~\ref{xcasteinb35} appears in Ree 
\cite[(3.4)]{ree2}. The proof of Proposition~\ref{steinb35} follows 
\cite[Lemma~35 (p.~40)]{St}. This will be important when we place 
$\bar{G}$ in the context of the theory of algebraic groups.

The results in Section~\ref{secA1} on groups of type $A_1$ are an adaptation
of the arguments of Ree \cite[\S 4]{ree2} to the special case where 
$\fg=\slm_2(\C)$. Those arguments also apply to any~$\fg$ but then require 
further preparations, which we will only discuss in a later chapter. 

Finally, some words about the different approaches of Ree \cite{ree2} 
and Steinberg \cite{St}. Firstly, there is the crucial issue of the
existence of ``admissible lattices'', or ``regular bases'' in the language
of Ree. For this purpose, Steinberg uses results like the 
Poincar\'e--Birkhoff--Witt theorem and Kostant's $\Z$-form 
of the universal enveloping algebra~$\mathcal{U}$ of~$\fg$. Instead, Ree 
uses Cartan's classification of the irreducible representations of~$\fg$, 
and quite explicit properties of them for all types of~$\fg$. (Ree himself
remarks at the end of \cite[\S 1]{ree2} that it would be desirable to 
find a general proof.) See also the work of Smith \cite{smith} which, 
however, did not seem to have any visible resonance in the subsequent 
developments.

The second crucial issue are \nm{Chevalley's commutator relations}. 
Steinberg \cite[Chap.~3]{St} starts the whole discussion of Chevalley 
groups with a result that proves those relations by an argument involving 
a computation in the formal power series ring in two commuting 
variables over Kostant's $\Z$-form of~$\mathcal{U}$. This is quite 
short and elegant, given the material that has been prepared beforehand.
But it does not seem to be obvious (at least not to us) how to break down 
that argument to the elementary level that we wish to pursue here. On 
the other hand, Ree reduces the proof of the commutator relations 
to the adjoint case, where one can invoke Carter \cite{Ca1}. 

So our synthesis of Ree \cite{ree2} and Steinberg \cite{St} consists of
following Ree (and Carter) as far as the commutator relations are concerned;
we shall also follow Ree in establishing the all-important homomorphisms 
in Theorem~\ref{reest}. Otherwise, the development of the structure theory 
of Chevalley groups in this and in the following chapter mainly follows 
Steinberg. As far as the existence of ``admissible lattices'' is concerned, 
we follow our approach based on \cite{G2}; it seems that this is sufficient
for many purposes where Chevalley groups arise. 

\appendix

\chapter{Some complements and auxiliary results} \label{appaux}

\section{Generation of $\SL_n(K)$} \label{asecSLn}

Let $K$ be any field and $n\geq 1$. For $1\leq i,j\leq n$ let $E_{ij}$
be the $n \times n$-matrix with~$1$ as its $(i,j)$-entry and zeroes 
elsewhere.  We define the following $n\times n$-matrices over $K$:
\[ x_i^*(\zeta):=I_n+\zeta E_{i,i+1}\quad \mbox{and}\quad y_i^*
(\zeta):=I_n+\zeta E_{i+1,i}\]
for $1\leq i \leq n-1$, where $\zeta\in K$ and $I_n$ is the $n\times 
n$-identity matrix over $K$. Then $x^*_i(\zeta)$ is upper triangular with 
$1$ along the diagonal; $y^*_i(\zeta)$ is lower triangular with $1$ along 
the diagonal. In particular, $\det(x_i^*(\zeta))=\det(y_i^*(\zeta))=1$.

\begin{prop} \label{gensln} With the above notation, we have 
\[\SL_n(K)=\langle x^*_i(\zeta), y_i^*(\zeta)\mid 1\leq i \leq n-1,
\zeta \in K \rangle.\]
\end{prop}

\begin{proof} We proceed by induction on $n$, where we start the induction
with $n=1$. Note that the assertion does hold for $\SL_1(K)=\{\id\}$.
Now let $n\geq 2$ and assume that the assertion is already proved for
$\SL_{n-1}(K)$. Let $G_n\subseteq \SL_n(K)$ be the subgroup generated
by the specified generators; we must show that $G_n=\SL_n(K)$. We set
\[x_{ij}^*(\zeta):=I_n+\zeta E_{ij}\qquad\mbox{for 
any $\zeta\in K$ and $1\leq i,j\leq n$, $i\neq j$}; \]
in particular, $x_{i}^*(\zeta)=x_{i,i+1}^*(\zeta)$ and $y_{i}^*(\zeta)=
x_{i+1,i}^*(\zeta)$. First we show:
\[ x_{i1}^*(\zeta)\in G_n \qquad \mbox{and}\qquad x_{1i}^*(\zeta)\in G_n
\qquad \mbox{for $2\leq i\leq n$}.\]
This is seen as follows. If $n=2$, there is nothing to show. Now
let $n\geq 3$. Let $i,j,k\in\{1,\ldots,n\}$ be pairwise distinct; 
then the following commutation rule is easily checked 
by an explicit computation:
\[ x_{jk}^*(-\zeta')\cdot x_{ij}^*(-\zeta)\cdot x_{jk}^*(\zeta')\cdot 
x_{ij}^*(\zeta)= x_{ik}^*(-\zeta \zeta')\]
for all $\zeta,\zeta'\in K$. Setting $\zeta'=-1$, $i=3$, $j=2$ and $k=1$, 
we obtain:
\[ x_{21}^*(1)\cdot x_{32}^*(-\zeta)\cdot x_{21}^*(-1)\cdot x_{32}^*
(\zeta)=x_{31}^*(\zeta)\] 
for all $\zeta\in K$. Hence, since the left hand side belongs to $G_n$, we
also have $x_{31}^*(\zeta)\in G_n$ for all $\zeta\in K$. Next, if
$n\geq 4$, then we set $\zeta'=-1$, $i=4$, $j=3$ and $k=1$. This yields 
\[ x_{31}^*(1)\cdot x_{43}^*(-\zeta)\cdot x_{31}^*(-1)\cdot x_{43}^*
(\zeta)=x_{41}^*(\zeta).\] 
Since the left hand side is already known to belong to $G_n$,
we also have $x_{41}^*(\zeta)\in G_n$. Continuing in this way, we find 
that $x_{i1}^*(\zeta)\in G_n$ for all $\zeta\in K$ and $2\leq i\leq n$.
The argument for $x_{1i}^*(\zeta)$ is analogous. 

Now let $T=(t_{ij})\in \SL_n(K)$ be arbitrary. It will be useful to 
remember that, for $i\geq 2$, the matrix $x_{i1}^*(\zeta)\cdot T$ is 
obtained by adding the first row of $T$, multiplied by~$\zeta$, to the 
$i$-th row of~$T$. Similarly, the matrix $T\cdot x_{1i}^*(\zeta)$ is 
obtained by adding the first column of~$T$, multiplied by~$\zeta$, to 
the $i$-th column of~$T$. We claim that there is a finite sequence of 
operations of this kind that transforms $T$ into a new matrix $B=(b_{ij})$ 
such that 
\[\renewcommand{\arraystretch}{1.1} B=\left(\begin{array}{c|c} 
\;1\; & 0\\\hline 0 & B' \end{array}\right) \qquad \mbox{where}
\qquad B'\in \SL_{n-1}(K).\]
Indeed, since $\det(T)\neq 0$, the first column of $T$ is non-zero 
and so there exists some $i\in\{1,\ldots, n\}$ such that 
$t_{i1}\neq 0$. If $i>1$, then 
\[T':=x_{i1}^*\bigl(t_{i1}^{-1}(1-t_{11})\bigr)\cdot T\]
has entry~$1$ at position $(1,1)$. But then we can add suitable
multiples of the first row of $T'$ to the other rows and obtain a new
matrix $A''$ that has entry $1$ at position $(1,1)$ and entry $0$ at
positions $(i,1)$ for $i\geq 2$. Next we can add suitable 
multiplies of the first column of $T''$ to the other columns and 
achieve that all further entries in the first row become~$0$. Thus, we 
have transformed $T$ into a new matrix $B$ as required.  On the 
other hand, if there is no $i>1$ such that $t_{i1}\neq 0$, then
$t_{11}\neq 0$ and $t_{i1}=0$ for $i\geq 2$. In that case,
the matrix $x_{21}^*(1)\cdot T$ has a non-zero entry at position $(2,1)$
and we are in the previous case. 

Now consider $B$ as above. By induction, we have $\SL_{n-1}(K)=G_{n-1}$; 
so the submatrix $B'$ can be expressed as a product of the specified
generators of $\SL_{n-1}(K)$. Under the embedding 
\[\renewcommand{\arraystretch}{1.1} 
\SL_{n-1}(K)\hookrightarrow \SL_n(K),\qquad C\mapsto 
\left(\begin{array}{c|c} \;1\; & 0\\\hline 0 & C \end{array}\right),\]
the generators of $\SL_{n-1}(K)$ are sent to the generators $x_i^*(\zeta)
\in \SL_n(K)$ and $y_i^*(\zeta)\in \SL_n(K)$, where $\zeta\in K$
and $2\leq i\leq n-1$. Hence, any $B$ as above can be expressed as a 
product of generators $x_i^*(\zeta)$ and $y_{i}^*(\zeta)$ in $\SL_n(K)$,
for various $\zeta\in K$ and $2\leq i \leq n-1$. Since $B$ was obtained 
from $T$ by a sequence of multiplications with matrices $x_{1i}^*(\zeta) 
\in G_n$ or $x_{i1}^*(\zeta)\in G_n$, we conclude that $T\in G_n$ (and 
we even described an algorithm for expressing $T$ in terms of the 
specified generators). 
\end{proof}

In particular, for $n=2$, we have 
\[\renewcommand{\arraystretch}{0.8} \SL_2(K)=\Big\langle 
\left(\begin{array}{cc} 1 & t \\ 0 & 1 \end{array} \right), 
\left(\begin{array}{cc} 1 & 0 \\ t & 1 \end{array}\right)\; \big|\; 
t\in K \Big\rangle.\]

\section{Matsumoto's Lemma} \label{asecmats}

Let $W$ be a group and $S\subseteq W$ be a subset such that $W=\langle
S\rangle$ and such that each element $s\in S$ has order~$2$. As in
Definition~\ref{deflength}, we define a length function $\ell\colon W
\rightarrow \Z_{\geq 0}$ (with respect to $S$) and the notion of a 
reduced expression for an element $w\in W$. We have again $\ell(w)=
\ell(w^{-1})$ and $\ell(w)-1\leq \ell(sw)\leq \ell(w)+1$ for $w\in W$ 
and $s\in S$. We assume that the following ``\nm{Exchange Condition}'' 
holds (analogous to Lemma~\ref{exccond}):
\begin{equation*}
{\rm (E)}\quad\left\{\begin{array}{l} \text{Let $w\in W$ and $s\in S$ be 
such that $\ell(sw)\leq \ell(w)$, and}\\\text{let $w=s_1\cdots s_r$ where 
$r=\ell(w)\geq 1$ and $s_1,\ldots,s_r\in S$.}\\\text{Then $ss_{1} \cdots 
s_{j-1}=s_1 \cdots s_{j-1}s_{j}$ for some $j\in \{1, \ldots,r\}$.}
\end{array}\right.
\end{equation*}
For $s\neq t$ in $S$ we denote by $m_{st}\geq 2$ the order of the product
$st\in W$. (Here, $m_{st}=\infty$ is allowed.) 

\begin{xca} \label{matsum1} Let $s\neq t$ in $S$ be such that $m_{st}<
\infty$. Show that the subgroup $W':=\langle s,t\rangle \subseteq W$ has
order $2m_{st}$ and that $\ell(w)\leq m_{st}$ for all $w\in W'$. Furthermore, 
let $w_0:=sts \cdots=tst \cdots \in W'$, with $m_{st}$ terms in the products
on both sides. Then show that $\ell(w_0)=m_{st}$.
\end{xca}

Now let $\cM$ be a monoid, that is, a set with an associative 
multiplication ``$*$'' for which there is an identity element $1_\cM$. 
Assume that there is a map $f\colon S \rightarrow \cM$ such that
\begin{equation*}
\underbrace{f(s){*}f(t){*}f(s){*} \ldots}_{\text{$m_{st}$ terms}}=
\underbrace{f(t){*}f(s){*}f(t){*}\ldots}_{\text{$m_{st}$ terms}}\qquad
(\in \cM)\tag{M}
\end{equation*}
for any $s\neq t$ in $S$ with $m_{st}<\infty$.

\begin{prop} \label{matsum} In the above setting, there is a well-defined 
map $\hat{f}\colon W \rightarrow \cM$ such that, for any $w\in W$ and any 
\nm{reduced expression} $w=s_1 s_2 \cdots s_r$ (where $r\geq 0$ and 
$s_1,\ldots,s_r\in S$), we have $\hat{f}(w)=f(s_1)*f(s_2)* \ldots *f(s_r)$.
\end{prop}

\begin{proof} We set $\hat{f}(1):=1_\cM$ and $\hat{f}(s):=f(s)$ for 
$s\in S$. Now let $w\in W$ and $p:=\ell(w)\geq 1$. Assume we are given two
reduced expressions 
\[ w=s_1\cdots s_p=t_1\cdots t_p \qquad\mbox{where $s_i,t_j\in S$}.\]
Then we must show that $f(s_1){*}\ldots {*}f(s_p)=f(t_1){*}\ldots {*}
f(t_p)$. We proceed by induction on $\ell(w)=p$. For $p=1$, we have 
$w=s_1=t_1$ and so there is nothing to prove. Now let $p\geq 2$ and assume 
that the assertion is wrong, that is, we have 
\begin{equation*}
f(s_1){*}f(s_2){*}\ldots {*}f(s_p)\neq f(t_1){*}f(t_2){*}\ldots {*}f(t_p).
\tag{0}
\end{equation*}
In this case, we say that $w=s_1\cdots s_p=t_1\cdots t_p$ are two
``bad'' expressions for $w$. 

Now $t_1w=t_1s_1s_2\cdots s_p=t_2\cdots t_p$ and so $\ell(t_1w)
\leq p-1<\ell(w)$. Applying (E) to $s:=t_1$ and the expression $w=s_1\cdots 
s_p$, there exists some $j\in \{1,\ldots,p\}$ such that $t_1s_1\cdots 
s_{j-1}=s_1\cdots s_j$ and so 
\[ t_1w=t_1.(s_1\cdots s_{j}).(s_{j+1}\cdots s_p)=t_1.(t_1s_1\cdots 
s_{j-1}).(s_{j+1}\cdots s_p),\]
Consequently, we have $w=t_1s_1\cdots s_{j-1}s_{j+1}\cdots s_p$, which is 
a new reduced expression for $w$ since there are exactly~$p$ factors. We 
claim that $j=p$. Assume, if possible, that $j<p$. Then the last factor
in the expression $w=t_1s_1\cdots s_{j-1}s_{j+1}\cdots s_p$ equals $s_p$
and so 
\begin{equation*}
ws_p=t_1s_1\cdots s_{j-1}s_{j+1}\cdots s_{p-1}=s_1s_2\cdots s_{p-1},\tag{1}
\end{equation*}
where the second equality holds since $w=s_1\cdots s_{p-1}s_p$. These are 
two reduced expressions for $ws_p$. By induction, $\hat{f}(ws_p)\in\cM$ 
is already well-defined. On the other hand, we have 
\begin{equation*}
t_1w=s_1\cdots s_{j-1}s_{j+1}\cdots s_p=t_2\cdots t_p,\tag{2}
\end{equation*}
where the second equality holds since $w=t_1t_2\cdots t_p$. These are 
also two reduced expressions for $t_1w$. By induction, $\hat{f}(t_1w)
\in \cM$ is already well-defined. Consequently, we obtain:
\begin{align*}
f(s_1)&{*}\cdots {*}f(s_p)=\bigl(f(s_1){*}\ldots {*}f(s_{p-1})\bigr){*}
f(s_p)= \hat{f}(ws_p){*}f(s_p)\\
&=\bigl(f(t_1){*}f(s_1){*}\ldots {*}f(s_{j-1}){*}f(s_{j+1}){*}\ldots 
{*}f(s_{p-1})\bigr){*}f(s_p)\\
&=f(t_1){*}\bigl(f(s_1){*}\ldots {*}f(s_{j-1}){*}f(s_{j+1}){*}\ldots 
{*}f(s_{p-1}){*}f(s_p)\bigr)\\
&=f(t_1){*}\hat{f}(t_1w)=f(t_1){*}\bigl(f(t_2){*}\ldots {*}f(t_p)\bigr)\\
&=f(t_1){*}f(t_2){*}\ldots {*}f(t_p),
\end{align*}
where we used (1) for the second/third equality, and (2) for the fifth/sixth 
equality. But this contradicts our assumption (0). Hence, we must have
$j=p$ and so we have the new reduced expression $w=t_1s_1\cdots s_{p-1}$. 
Then $t_1w=s_1\cdots s_{p-1}=t_2\cdots t_p$. By induction, we obtain 
\begin{align*}
f(t_1)&{*}f(s_1){*}\ldots {*}f(s_{p-1})=f(t_1){*}\hat{f}(t_1w)\\
&= f(t_1){*}f(t_2){*}\ldots{*}f(t_p)\neq 
f(s_1){*}\ldots {*}f(s_p).
\end{align*}
Thus, starting from the two bad expressions $w=s_1\cdots s_p=t_1\cdots t_p$,
we produced a new reduced expression $w=t_1s_1\cdots s_{p-1}$ such that
$w=t_1s_1\cdots s_{p-1}=s_1\cdots s_p$ are two bad expressions for~$w$. We 
now repeat the whole argument with these two reduced expressions for~$w$.
So we obtain two new bad expressions $w=s_1t_1s_1\cdots s_{p-2}=t_1s_1
\cdots s_{p-1}$. If $p\geq 3$, we repeat again and obtain two new bad 
expressions $w=t_1 s_1t_1s_1\cdots s_{p-3}=s_1t_1s_1\cdots s_{p-2}$. After 
$p$ repetitions we eventually find two bad expressions $w=t_1s_1t_1\cdots=
s_1t_1s_1\cdots$, with $p$ factors on each side. It follows that $(s_1
t_1)^p=1$ and so $m_{s_1t_1}\leq p<\infty$. Now note that $w\in W':=
\langle s_1,t_1\rangle\subseteq W$. By Exercise~\ref{matsum1}, we 
have $|W'|=2m_{s_1t_1}$ and $\ell(w)=p\leq m_{s_1t_1}$. Hence, $p=m_{s_1t_1}$ 
and we obtain a final contradiction to the assumption that $f(t_1){*} f(s_1)
{*}f(t_1){*}\ldots= f(s_1){*} f(t_1){*}f(s_1){*} \ldots$ (with $p=m_{s_1t_1}$
factors on both sides).
\end{proof}

\begin{cor} \label{cmatsum} Assume that $\cM$ is a group such that
$f(s)^2=1_\cM$ for all $s\in S$. Then the map $\hat{f}\colon W
\rightarrow \cM$ is a group homomorphism. 
\end{cor}

\begin{proof} We must show that $\hat{f}(ww')=\hat{f}(w)*\hat{f}(w')$
for all $w,w'\in W$. Since $W=\langle S\rangle$, it is sufficient to show 
that $\hat{f}(sw)=\hat{f}(s)*\hat{f}(w)$ for all $s\in S$ and $w\in W$.
We do this by induction in $\ell(w)$. If $\ell(w)=0$ then $w=1$ and $\hat{f}
(1_W)=1_\cM$; so the assertion holds in this case. Now let $r:=\ell(w)\geq 
1$ and consider a reduced expression $w=s_1\cdots s_r$, where $s_1,\ldots,
s_r\in S$. Let $s\in S$ be arbitrary and set $w':=sw$; then $\ell(w)-1
\leq \ell(w')\leq \ell(w)+1$. Now there are two cases.

If $\ell(sw)>\ell(w)$, then the expression 
$w'=sw=ss_1 \cdots s_r$ is reduced. So we have  
\[ \hat{f}(sw)=\hat{f}(w')=\hat{f}(ss_1\cdots s_r)=\hat{f}(s)*\hat{f}
(s_1)*\ldots * \hat{f}(s_r).\]
Now $\hat{f}(w)=\hat{f}(s_1)*\ldots *\hat{f}(s_r)$ and so 
$\hat{f}(sw)=\hat{f}(s)*\hat{f}(w)$, as required. 

Now assume that $\ell(sw)\leq \ell(w)$. By (E), there exists some $j \in 
\{1,\ldots,r\}$ such that $ss_1\cdots s_{j-1}=s_1\cdots s_{j-1}s_j$. 
Hence, we have 
\begin{align*}
 w'&=sw=(ss_1\cdots s_{j-1})(s_js_{j+1}\cdots s_r)\\
&=(s_1\cdots s_{j-1}s_j)(s_js_{j+1}\cdots s_r)=s_1\cdots s_{j-1}s_{j+1}
\cdots s_r,
\end{align*}
which shows that $\ell(w')<\ell(w)$. Since $s^2=1$, we have $w=sw'$. By 
induction, we obtain 
\[ \hat{f}(w)=\hat{f}(sw')=\hat{f}(s)*\hat{f}(w')=\hat{f}(s)*\hat{f}(sw).\]
Since $\hat{f}(s)^2=1_\cM$, this implies $\hat{f}(sw)=\hat{f}(s)*
\hat{f}(w)$, as required.
\end{proof}

If $W$ is a finite Coxeter group, then Proposition~\ref{matsum} 
already appeared in Iwahori \cite[Theorem~2.6]{Iwa}.

\backmatter

\begin{theindex}

  \item abelian Lie algebra, \hyperpage{6}
  \item abstract root system, \hyperpage{126}
  \item adjoint Chevalley group, \hyperpage{262}
  \item adjoint representation, \hyperpage{22}
  \item admissible lattice, \hyperpage{222}
  \item affine type, \hyperpage{117}
  \item algebra, \hyperpage{1}
  \item algebra homomorphism, \hyperpage{2}
  \item algebra isomorphism, \hyperpage{2}
  \item $\alpha $-string through $\beta $, \hyperpage{92}
  \item $\alpha _i$-string through $\beta $, \hyperpage{60}, 
		\hyperpage{100}
  \item alternating bilinear form, \hyperpage{36}

  \indexspace

  \item based root system, \hyperpage{53}, \hyperpage{128}
  \item bracket, \hyperpage{5}
  \item braid relations, \hyperpage{133}, \hyperpage{253}

  \indexspace

  \item Cartan decomposition, \hyperpage{53}
  \item Cartan subalgebra, \hyperpage{53}
  \item Cartan's First Criterion, \hyperpage{32}, \hyperpage{42}
  \item Cartan--Killing classification, \hyperpage{41}
  \item Cartan--Killing type, \hyperpage{iv}, \hyperpage{53}, 
		\hyperpage{189}
  \item Cauchy--Schwarz inequality, \hyperpage{95}, \hyperpage{131}
  \item center, \hyperpage{10}
  \item character, \hyperpage{3}
  \item Chevalley generators, \hyperpage{59}, \hyperpage{75}
  \item Chevalley group, \hyperpage{15}, \hyperpage{162}, 
		\hyperpage{226}
  \item Chevalley group of adjoint type, \hyperpage{162}, 
		\hyperpage{262}
  \item Chevalley group of simply connected type, \hyperpage{262}
  \item Chevalley involution, \hyperpage{110}
  \item Chevalley relations, \hyperpage{136}, \hyperpage{142}, 
		\hyperpage{144}, \hyperpage{212}
  \item Chevalley's commutator relations, \hyperpage{263}, 
		\hyperpage{265}
  \item classical Lie algebra, \hyperpage{36}, \hyperpage{82}
  \item co-root, \hyperpage{92}
  \item coding theory, \hyperpage{130}
  \item combinatorial graph of $A$, \hyperpage{70}, \hyperpage{77}
  \item commutator subgroup, \hyperpage{173}
  \item computer algebra approach, \hyperpage{128}, \hyperpage{133}, 
		\hyperpage{145}
  \item cycle (graph theory), \hyperpage{71}

  \indexspace

  \item derivation, \hyperpage{10}
  \item diagonal automorphism, \hyperpage{240}
  \item diagonal subgroup, \hyperpage{180}, \hyperpage{243}
  \item diagonalisable, \hyperpage{44}
  \item direct product, \hyperpage{4}, \hyperpage{14}
  \item dominant weight, \hyperpage{200}
  \item Dynkin diagram, \hyperpage{82}, \hyperpage{120}

  \indexspace

  \item elementary Chevalley group, \hyperpage{185}
  \item $\epsilon $-canonical Chevalley system, \hyperpage{152}
  \item Exchange Condition, \hyperpage{154}, \hyperpage{269}
  \item Existence Theorem, \hyperpage{145}

  \indexspace

  \item faithful ${\mathfrak  {g}}$-module, \hyperpage{193}
  \item field automorphism, \hyperpage{228}, \hyperpage{240}
  \item finite type, \hyperpage{117}
  \item forest, \hyperpage{103}
  \item forest (graph theory), \hyperpage{71}
  \item free Lie algebra, \hyperpage{8}
  \item free magma, \hyperpage{7}
  \item Frobenius--Perron Theorem, \hyperpage{125}
  \item fundamental calculation, \hyperpage{190}, \hyperpage{239}, 
		\hyperpage{244}
  \item fundamental group, \hyperpage{192}
  \item fundamental weights, \hyperpage{191}

  \indexspace

  \item generalized binomial formula, \hyperpage{4}, \hyperpage{29}
  \item generalized Cartan matrix, \hyperpage{119}, \hyperpage{142}
  \item generalized eigenspace, \hyperpage{29}
  \item graph automorphism, \hyperpage{136}, \hyperpage{242}
  \item graph of $A$, \hyperpage{70}

  \indexspace

  \item half-spin representations, \hyperpage{216}, \hyperpage{218}
  \item height, \hyperpage{54}, \hyperpage{134}
  \item Heisenberg Lie algebra, \hyperpage{16}
  \item highest root, \hyperpage{134}
  \item highest short root, \hyperpage{135}
  \item highest weight, \hyperpage{210}
  \item highest weight module, \hyperpage{210}
  \item hyperoctahedral group, \hyperpage{88}

  \indexspace

  \item ideal, \hyperpage{2}
  \item indecomposable blocks, \hyperpage{81}
  \item indecomposable matrix, \hyperpage{77}
  \item indecomposable subsystems, \hyperpage{82}
  \item indefinite type, \hyperpage{117}
  \item integrable, \hyperpage{142}
  \item invariant subspace, \hyperpage{23}
  \item irreducible module, \hyperpage{23}
  \item isogeny theorem, \hyperpage{242}
  \item Isomorphism Theorem, \hyperpage{109}

  \indexspace

  \item Kac--Moody algebra, \hyperpage{142}
  \item Key Lemma, \hyperpage{64}
  \item Killing form, \hyperpage{18}

  \indexspace

  \item $L$-module, \hyperpage{21}
  \item $L$-module homomorphism, \hyperpage{25}
  \item lattice, \hyperpage{191}
  \item Leibniz rule, \hyperpage{11}
  \item length, \hyperpage{153}
  \item length of $\alpha $, \hyperpage{130}
  \item Linear independence of characters, \hyperpage{3}
  \item locally nilpotent, \hyperpage{10}, \hyperpage{142}
  \item long roots, \hyperpage{131}
  \item longest element, \hyperpage{156}

  \indexspace

  \item magma, \hyperpage{3}
  \item magma algebra, \hyperpage{3}
  \item matrix representation, \hyperpage{22}
  \item maximal torus, \hyperpage{180}
  \item minuscule weight, \hyperpage{203}
  \item monomial subgroup, \hyperpage{243}
  \item monomials, \hyperpage{3}

  \indexspace

  \item negative root, \hyperpage{53}
  \item non-degenerate bilinear form, \hyperpage{36}
  \item normalizer, \hyperpage{6}

  \indexspace

  \item orbit algorithm, \hyperpage{69}
  \item orthogonal reflection, \hyperpage{64}

  \indexspace

  \item positive root, \hyperpage{53}
  \item primitive vector, \hyperpage{33}, \hyperpage{61}, 
		\hyperpage{208}

  \indexspace

  \item radical, \hyperpage{17}
  \item reduced expression, \hyperpage{153}, \hyperpage{270}
  \item reflection, \hyperpage{153}
  \item reflexive bilinear form, \hyperpage{36}
  \item regular basis, \hyperpage{222}
  \item regular submodule, \hyperpage{256}
  \item representation, \hyperpage{22}
  \item root lattice, \hyperpage{191}
  \item roots, \hyperpage{53}

  \indexspace

  \item Schur's Lemma, \hyperpage{26}
  \item semi-spinor representations, \hyperpage{216}
  \item semidirect product, \hyperpage{14}, \hyperpage{23}
  \item semisimple algebraic group, \hyperpage{187}
  \item semisimple Lie algebra, \hyperpage{18}
  \item Serre relations, \hyperpage{142}
  \item short roots, \hyperpage{131}
  \item signed permutation, \hyperpage{88}
  \item simple algebra, \hyperpage{7}
  \item simple roots, \hyperpage{126}
  \item simply laced, \hyperpage{131}, \hyperpage{135}, \hyperpage{160}
  \item ${\mathfrak  {sl}}_2$-triple, \hyperpage{59}, \hyperpage{189}, 
		\hyperpage{208}
  \item solvable algebra, \hyperpage{15}
  \item spin representations, \hyperpage{216}
  \item spinor representation, \hyperpage{216}
  \item structure matrix, \hyperpage{57}
  \item subalgebra, \hyperpage{2}
  \item submodule, \hyperpage{23}
  \item symmetric bilinear form, \hyperpage{19}, \hyperpage{36}
  \item symplectic group, \hyperpage{172}

  \indexspace

  \item theory of elementary divisors, \hyperpage{192}
  \item torus, \hyperpage{243}
  \item Transfer Lemma, \hyperpage{168}, \hyperpage{230}
  \item triangular decomposition, \hyperpage{40}, \hyperpage{50}, 
		\hyperpage{53}

  \indexspace

  \item universal Chevalley group, \hyperpage{262}
  \item universal property of ${\mathfrak  {g}}(A)$, \hyperpage{146}
  \item universal property of the free Lie algebra, \hyperpage{8}
  \item universal property of the free magma, \hyperpage{8}

  \indexspace

  \item weight, \hyperpage{31}, \hyperpage{44}, \hyperpage{139}, 
		\hyperpage{189}
  \item weight lattice, \hyperpage{257}
  \item weight lattice (of ${\mathfrak  {g}}$), \hyperpage{191}
  \item weight lattice of $V$, \hyperpage{193}
  \item weight order relation, \hyperpage{199}, \hyperpage{221}
  \item weight space, \hyperpage{46}, \hyperpage{139}, \hyperpage{189}
  \item weight vector, \hyperpage{46}, \hyperpage{222}
  \item Weyl group, \hyperpage{63}, \hyperpage{126}
  \item Witt algebra, \hyperpage{13}, \hyperpage{25}, \hyperpage{41}

\end{theindex}
\end{document}